\documentstyle[12pt,amscd]{amsart}
\newsymbol\boxtimes 1202

\newcommand{\nc}{\newcommand}

\nc{\ad}{{\mbox{\bf{ad}}}}
\nc{\AJ}{{\operatorname{AJ}}}
\nc{\divv}{{\operatorname{div}}}
\nc{\Aut}{{\operatorname{Aut}}}
\nc{\BAO}{{\overset{\circ}{\Bbb A}}{}}
\nc{\Bls}{{{{\cal B}}ls}}
\nc{\Boxtimes}{{\fbox{$\times$}}}
\nc{\blt}{{\bullet}}
\nc{\bSp}{{\mbox{\bf{Sp}}}}
\nc{\bSt}{{\mbox{\bf{St}}}}
\nc{\card}{{\operatorname{card}}}
\nc{\Cch}{{\check{C}}}
\nc{\cd}{{\operatorname{cd}}}
\nc{\Ch}{{\operatorname{Ch}}}
\nc{\chara}{{\operatorname{char}}}
\nc{\CHom}{{{\cal{H}}om}}
\nc{\Coker}{{\operatorname{Coker}}}
\nc{\codim}{{\operatorname{codim}}}
\nc{\Cone}{{\operatorname{Cone}}}
\nc{\cSgn}{{{{\cal S}}ign}}
\nc{\depth}{{\operatorname{depth}}}
\nc{\dirlim}{{\underset{\rightarrow}{\operatorname{lim}}}}
\nc{\dotbox}{{\overset{\bullet}{\boxtimes}}}
\nc{\dotimes}{{\overset{\bullet}{\otimes}}}
\nc{\DT}{{\overset{\bullet}{T}}}
\nc{\Ed}{{\operatorname{Edge}}}
\nc{\emp}{{\emptyset}}
\nc{\Ext}{{\operatorname{Ext}}}
\nc{\Fac}{{{{\cal F}}ac}}
\nc{\Fun}{{\operatorname{F}}}
\nc{\FS}{{\cal{FS}}}
\nc{\Hom}{{\operatorname{Hom}}}
\nc{\End}{{\underline{\operatorname{End}}}}
\nc{\had}{{{\hat{\mbox{\bf{ad}}}}}}
\nc{\hgt}{{\operatorname{ht}}}
\nc{\Id}{{\operatorname{Id}}}
\nc{\id}{{\operatorname{id}}}
\nc{\Ima}{{\operatorname{Im}}}
\nc{\ind}{{\operatorname{ind}}}
\nc{\Ind}{{\operatorname{Ind}}}
\nc{\infi}{{\operatorname{inf}}}
\nc{\infh}{{\frac{\infty}{2}}}
\nc{\invlim}{{\underset{\leftarrow}{\operatorname{lim}}}}
\nc{\Jac}{{{{\cal J}}ac}}
\nc{\Ker}{{\operatorname{Ker}}}
\nc{\sk}{{\sf k}}
\nc{\lcm}{{\operatorname{lcm}}}
\nc{\Locsys}{{{{\cal L}}ocsys}}
\nc{\Map}{{{\cal M}ap}}
\nc{\modul}{{\operatorname{mod}}}
\nc{\Mor}{{\operatorname{Mor}}}
\nc{\MS}{{\cal{MS}}}
\nc{\Ob}{{\operatorname{Ob}}}
\nc{\oC}{{\operatorname{C}}}
\nc{\od}{{\operatorname{d}}}
\nc{\ob}{{\operatorname{b}}}
\nc{\og}{{\operatorname{g}}}
\nc{\opp}{{\operatorname{opp}}}
\nc{\Or}{{{{\cal O}}r}}
\nc{\Ord}{{{{\cal O}}rd}}
\nc{\Part}{{{{\cal P}}art}}
\nc{\PGL}{{\operatorname{PGL}}}
\nc{\Pic}{{\operatorname{Pic}}}
\nc{\pr}{{\operatorname{pr}}}
\nc{\Rep}{{{\cal{R}}ep}}
\nc{\rk}{{\operatorname{rk}}}
\nc{\Sets}{{{\cal{S}}ets}}
\nc{\Sew}{{{\cal{S}}ew}}
\nc{\sgn}{{\operatorname{sgn}}}
\nc{\Sh}{{{{\cal S}}h}}
\nc{\Sign}{{{{\cal S}}ign}}
\nc{\Spe}{{\mbox{\bf{Sp}}}}
\nc{\supr}{{\operatorname{sup}}}
\nc{\Supp}{{\operatorname{Supp}}}
\nc{\supp}{{\operatorname{supp}}}
\nc{\Teich}{{{\cal{T}}eich}}
\nc{\tFS}{{\widetilde{\cal{FS}}}}
\nc{\Tor}{{\operatorname{Tor}}}
\nc{\totimes}{{\tilde{\otimes}}}
\nc{\tr}{{\operatorname{tr}}}
\nc{\tRep}{{\widetilde{{\cal R}ep}}}
\nc{\tTeich}{{\widetilde{{\cal T}eich}}}
\nc{\Vect}{{{{\cal V}}ect}}
\nc{\Ve}{{\operatorname{Vert}}}
\nc{\wt}{\widetilde}

\nc{\bo}{{\mbox{\bf{0}}}}
\nc{\One}{{\bf{1}}}
\nc{\one}{{\bf{1}}}

\nc{\BA}{{\Bbb A}}
\nc{\bA}{{\overline{A}}}
\nc{\ba}{{\mbox{\bf{a}}}}
\nc{\baB}{{\overline{B}}}
\nc{\baeta}{{\bar{\eta}}}
\nc{\baJ}{{\bar{J}}}
\nc{\BB}{{\Bbb B}}
\nc{\bB}{{\mbox{\bf{B}}}}
\nc{\bc}{{\mbox{\bf{c}}}}
\nc{\bC}{{\overline{C}}}
\nc{\BC}{{\Bbb{C}}}
\nc{\bCC}{{\overline{\cal{C}}}}
\nc{\bCM}{{\overline{\cal{M}}}}
\nc{\bD}{{\bar{D}}}
\nc{\BD}{{\overline{D}}}
\nc{\bd}{{\mbox{\bf{d}}}}
\nc{\BE}{{\overline{E}}}
\nc{\BF}{{\overline{F}}}
\nc{\bF}{{{\bf{F}}}}
\nc{\bg}{{{\bf{g}}}}
\nc{\bGamma}{{\overline{\Gamma}}}
\nc{\bL}{{\bf{L}}}
\nc{\blambda}{{\bar{\lambda}}}
\nc{\bM}{{\bf{M}}}
\nc{\bmu}{{\vec{\mu}}}
\nc{\BN}{{\Bbb{N}}}
\nc{\bnu}{{\mbox{\boldmath{${\nu}$}}}}
\nc{\bof}{{\bf{f}}}
\nc{\BP}{{\Bbb P}}
\nc{\bP}{{\bf{P}}}
\nc{\BPO}{{\overset{\circ}{\BP}}}
\nc{\BQ}{{\Bbb Q}}
\nc{\bq}{{\bf{q}}}
\nc{\BR}{{\Bbb{R}}}
\nc{\bR}{{\bf{R}}}
\nc{\br}{{\bf{r}}}
\nc{\breta}{{\bar{\eta}}}
\nc{\bs}{{\bf{s}}}
\nc{\bS}{{\bf{S}}}
\nc{\bt}{{\bf{t}}}
\nc{\bU}{{\bf{U}}}
\nc{\bu}{{\bf{u}}}
\nc{\BUpsilon}{{\bar{\Upsilon}}}
\nc{\bw}{{\bf{w}}}
\nc{\bx}{{\bf{x}}}
\nc{\BZ}{{\Bbb{Z}}}
\nc{\bz}{{\bf{z}}}
\nc{\bzero}{\mbox{\boldmath{$0$}}}

\nc{\CA}{{\cal A}}
\nc{\CAD}{{\overset{\bullet}{\cal{A}}}{}}
\nc{\CAO}{{\overset{\circ}{\cal{A}}}{}}
\nc{\CB}{{\cal B}}
\nc{\CC}{{\cal C}}
\nc{\Co}{{{\overset{\circ}{C}}}{}}
\nc{\TCo}{{{\overset{\circ}{TC}}}{}}
\nc{\DD}{{\overset{\bullet}{D}}{}}
\nc{\Do}{{{\overset{\circ}{D}}}{}}
\nc{\TDo}{{{\overset{\circ}{TD}}}{}}
\nc{\CCD}{{{\cal D}}}
\nc{\CE}{{\cal E}}
\nc{\CF}{{\cal F}}
\nc{\CH}{{\cal H}}
\nc{\CI}{{\cal I}}
\nc{\CID}{{\overset{\bullet}{\cal{I}}}{}}
\nc{\CJ}{{\cal J}}
\nc{\CK}{{\cal K}}
\nc{\CL}{{\cal L}}
\nc{\CM}{{\cal M}}
\nc{\CN}{{\cal N}}
\nc{\CO}{{\cal O}}
\nc{\CP}{{\cal P}}
\nc{\CPO}{{\overset{\circ}{\cal{P}}}{}}
\nc{\CQ}{{\cal Q}}
\nc{\CR}{{\cal R}}
\nc{\CS}{{\cal S}}
\nc{\CT}{{\cal T}}
\nc{\CTD}{{\overset{\bullet}{\cal{T}}}{}}
\nc{\CTPO}{{\overset{\circ}{\cal{T}\cal{P}}}{}}
\nc{\CU}{{\cal{U}}}
\nc{\CV}{{\cal V}}
\nc{\CX}{{\cal X}}
\nc{\CY}{{\cal Y}}
\nc{\CZ}{{\cal Z}}

\nc{\dCL}{{\overset{\bullet}{\cal{L}}}{}}
\nc{\dd}{{\operatorname{d}}}
\nc{\ddelta}{{\overset{\bullet}{\delta}}{}}
\nc{\dfu}{{\overset{\bullet}{\frak{u}}}{}}
\nc{\dlambda}{{\overset{\bullet}{\lambda}}{}}
\nc{\DO}{{\overset{\circ}{D}}{}}
\nc{\dpar}{{\partial}}
\nc{\dS}{{\overset{\bullet}{S}}{}}
\nc{\dT}{{\overset{\bullet}{T}}{}}

\nc{\fa}{{\frak{a}}}
\nc{\fA}{{\frak{A}}}
\nc{\fC}{{\frak{C}}}
\nc{\fD}{{\frak{D}}}
\nc{\fE}{{\frak{E}}}
\nc{\fF}{{\frak{F}}}
\nc{\ff}{{\frak{f}}}
\nc{\fg}{{\frak{g}}}
\nc{\fH}{{\frak{H}}}
\nc{\fl}{{\frak{l}}}
\nc{\fL}{\frak{L}}
\nc{\fM}{\frak{M}}
\nc{\fp}{{\frak{p}}}
\nc{\fu}{{\frak{u}}}
\nc{\fm}{{{\frak{m}}}}
\nc{\fs}{{\frak{s}}}

\nc{\hCH}{{\hat{\cal{H}}}}
\nc{\hCI}{{\hat{\cal{I}}}}
\nc{\hfC}{{\hat{\frak{C}}}}
\nc{\hfg}{{\hat{\frak{g}}}}
\nc{\hL}{{\hat{L}}}
\nc{\HO}{{\overset{\circ}{H}}{}}
\nc{\hpsi}{{\hat{\psi}}}
\nc{\hx}{{\hat{x}}}

\nc{\jo}{{\overset{\circ}{j}}{}}

\nc{\phid}{{\overset{\bullet}{\phi}}{}}

\nc{\tA}{{\tilde{A}}}
\nc{\ta}{{\tilde{a}}}
\nc{\tB}{{\tilde{B}}}
\nc{\tb}{{\tilde{b}}}
\nc{\tBP}{{\tilde{\BP}}}
\nc{\tC}{{\tilde{C}}}
\nc{\tc}{{\tilde{c}}}
\nc{\tCA}{{\tilde{\cal{A}}}}
\nc{\tCC}{{\tilde{\cal{C}}}}
\nc{\tCH}{{\tilde{\cal{H}}}}
\nc{\tCI}{{\tilde{\cal{I}}}}
\nc{\tCO}{{\tilde{\cal{O}}}}
\nc{\tCP}{{\tilde{\cal{P}}}}
\nc{\tCT}{{\tilde{\cal{T}}}}
\nc{\tCX}{{\tilde{\cal{X}}}}
\nc{\tD}{{\tilde{D}}}
\nc{\tDelta}{{\tilde{\Delta}}}
\nc{\tE}{{\tilde E}}
\nc{\tF}{{\tilde F}}
\nc{\tfD}{{\tilde{\frak{D}}}}
\nc{\tfF}{{\tilde{\frak{F}}}}
\nc{\tff}{{\tilde{\frak{f}}}}
\nc{\tfu}{{\tilde{\frak{u}}}}
\nc{\tJ}{{\tilde{J}}}
\nc{\tj}{{\tilde{j}}}
\nc{\tK}{{\tilde K}}
\nc{\tL}{{\tilde{L}}}
\nc{\tM}{{\tilde{M}}}
\nc{\tP}{{\tilde{P}}}
\nc{\tPhi}{{\tilde{\Phi}}}
\nc{\tpi}{\tilde{\pi}}
\nc{\TPO}{{\overset{\circ}{T\BP}}{}}
\nc{\tR}{{\tilde{R}}}
\nc{\tS}{{\tilde S}}
\nc{\tT}{{\tilde{T}}}
\nc{\ttau}{{\tilde{\tau}}}
\nc{\ttheta}{{\tilde{\theta}}}
\nc{\tU}{{\tilde{U}}}
\nc{\tUpsilon}{{\tilde{\Upsilon}}}
\nc{\ty}{{\tilde y}}
\nc{\tY}{{\tilde Y}}
\nc{\txi}{{\tilde{\xi}}}

\nc{\UD}{{\overset{\bullet}{U}}{}}
\nc{\UO}{{\overset{\circ}{U}}{}}
\nc{\TO}{{{\overset{\circ}{T}}}{}}
\nc{\TB}{{{\overset{\bullet}{T}}}{}}

\nc{\vA}{{\vec{A}}}
\nc{\valpha}{{\vec{\alpha}}}
\nc{\vbeta}{{\vec{\beta}}}
\nc{\vc}{{\vec{c}}}
\nc{\vD}{{\vec{D}}}
\nc{\vd}{{\vec{d}}}
\nc{\vgamma}{{\vec{\gamma}}}
\nc{\vK}{{\vec{K}}}
\nc{\vlambda}{{\vec{\lambda}}}
\nc{\vmu}{{\vec{\mu}}}
\nc{\vnu}{{\vec{\nu}}}
\nc{\vo}{{\vec{0}}}
\nc{\vu}{{\vec{u}}}
\nc{\vx}{{\vec{x}}}
\nc{\vy}{\vec{y}}
\nc{\vzero}{\vec{0}}

\nc{\XO}{{\overset{\circ}{X}}{}}

\nc{\ya}{{\operatorname{aj}}}

\nc{\nen}{\newenvironment}
\nc{\ol}{\overline}
\nc{\ul}{\underline}
\nc{\ra}{\rightarrow}
\nc{\lra}{\longrightarrow}
\nc{\Lra}{\Longrightarrow}
\nc{\lla}{\longleftarrow}
\nc{\Llra}{\Longleftrightarrow}
\nc{\hra}{\hookrightarrow}
\nc{\iso}{\overset{\sim}{\lra}}
\nc{\rlh}{\rightleftharpoons}

\nc{\Thm}[1]{Theorem~\ref{#1}}
\nc{\Prop}[1]{Proposition~\ref{#1}}
\nc{\Lem}[1]{Lemma~\ref{#1}}
\nc{\Cor}[1]{Corollary~\ref{#1}}
\nc{\Conj}[1]{Conjecture~\ref{#1}}
\nc{\Claim}[1]{Claim~\ref{#1}}
\nc{\Defn}[1]{Definition~\ref{#1}}
\nc{\Exa}[1]{Example~\ref{#1}}
\nc{\Rem}[1]{Remark~\ref{#1}}
\nc{\Note}[1]{Note~\ref{#1}}

\nen{thm}[1]{\label{#1}{\bf Theorem.\ } \em}{}
\nen{prop}[1]{\label{#1}{\bf Proposition.\ } \em}{}
\nen{lem}[1]{\label{#1}{\bf Lemma.\ } \em}{}
\nen{cor}[1]{\label{#1}{\bf Corollary.\ } \em}{}
\nen{conj}[1]{\label{#1}{\bf Conjecture.\ } \em}{}
\nen{claim}[1]{\label{#1}{\bf Claim.\ } \em}{}

\nen{defn}[1]{\label{#1}{\bf Definition.\ } }{}
\nen{exa}[1]{\label{#1}{\bf Example.\ } }{}

\nen{rem}[1]{\label{#1}{\em Remark.\ } }{}
\nen{note}[1]{\label{#1}{\em Note.\ } }{}
\nen{exer}[1]{\label{#1}{\em Exercise.\ } }{}

\begin{document}

\bigskip

\bigskip

\centerline{\bf\Large FACTORIZABLE SHEAVES AND QUANTUM GROUPS}

\bigskip
\bigskip

\centerline{\large Roman Bezrukavnikov, Michael Finkelberg, Vadim Schechtman}


\bigskip

\bigskip

\pagenumbering{roman}
\newpage
empty
\newpage
empty
\newpage
empty
\newpage
\centerline{-}
\vspace{4cm}
\centerline{\em To Alexandre Grothendieck on his 70th
birthday, with admiration}
\newpage
\begin{center}
{\bf Table of contents}
\end{center}
\vspace{1cm}

Introduction 1
\vspace{1cm}

Acknowledgement 5

\vspace{1cm}

{\bf Part 0. OVERVIEW}

\vspace{.4cm}

1. Introduction 6

\vspace{.3cm}

{\bf Chapter 1. Local}

\vspace{.3cm}  

2. The category $\CC$ 11

3. Braiding local systems 15

4. Factorizable sheaves 19

5. Tensor product 20

6. Vanishing cycles 23

\vspace{.3cm}

{\bf Chapter 2. Global (genus $0$)}

\vspace{.3cm}

7. Cohesive local systems 27

8. Gluing 30

9. Semiinfinite cohomology 31

10. Conformal blocks (genus 0) 33

11. Integration 35

12. Regular representation 36

13. Regular sheaf 37

\vspace{.3cm}

{\bf Chapter 3. Modular}

\vspace{.3cm}

14. Heisenberg local system 40

15. Fusion structures on $\FS$ 46

16. Conformal blocks (higher genus) 48

\newpage


{\bf Part I. INTERSECTION COHOMOLOGY
 
OF REAL ARRANGEMENTS}

\vspace{1cm}

1. Introduction 50

2. Topological preliminaries 51

3. Vanishing cycles functors 55

4. Computations for standard sheaves 63

\vspace{1cm}

{\bf Part II. CONFIGURATION SPACES AND QUANTUM GROUPS}

\vspace{1cm}

1. Introduction 71

\vspace{.3cm}

{\bf Chapter 1. Algebraic discussion}

\vspace{.3cm}

2. Free algebras and bilinear forms 73

3. Hochschild complexes 82

4. Symmetrization 84

5. Quotient algebras 87

\vspace{.3cm}

{\bf Chapter 2. Geometric discussion}

\vspace{.3cm}

6. Diagonal stratification and related algebras 89

7. Principal stratification 94

8. Standard sheaves 100

\vspace{.3cm}

{\bf Chapter 3. Fusion}

\vspace{.3cm}

9. Additivity theorem 110

10. Fusion and tensor products 112

\vspace{.3cm}

{\bf Chapter 4. Category $\CC$}

\vspace{.3cm}

11. Simply laced case 116

12. Non-simply laced case 119

\vspace{1cm}

\newpage

{\bf Part III. TENSOR CATEGORIES ARISING FROM 

CONFIGURATION SPACES}

\vspace{1cm}

1. Introduction 122

\vspace{.3cm}

{\bf Chapter 1. Category $\FS$}

\vspace{.3cm}

2. Space $\CA$ 124

3. Braiding local system $\CI$ 126

4. Factorizable sheaves 128

5. Finite sheaves 130

6. Standard sheaves 132

\vspace{.3cm}

{\bf Chapter 2. Tensor structure}

\vspace{.3cm}

7. Marked disk operad 134

8. Cohesive local systems $^K\CI$ 138

9. Factorizable sheaves over $^K\CA$ 139

10. Gluing 142

11. Fusion 145

\vspace{.3cm}

{\bf Chapter 3. Functor $\Phi$}

\vspace{.3cm}

12. Functor $\Phi$ 150

13. Main properties of $\Phi$ 159

\vspace{.3cm}

{\bf Chapter 4. Equivalence}

\vspace{.3cm}

14. Truncation functors 162

15. Rigidity 165

16. Steinberg sheaf 167

17. Equivalence 169

18. The case of generic $\zeta$ 170

\vspace{1cm}

{\bf Part IV. LOCALIZATION OVER $\BP^1$}

\vspace{1cm}

1. Introduction 173

\vspace{.3cm}

{\bf Chapter 1. Gluing over $\BP^1$}

\vspace{.3cm}

2. Cohesive local system 174

3. Gluing 176

\vspace{.3cm}

{\bf Chapter 2. Semiinifinite cohomology}

\vspace{.3cm}

4. Semiinfinite functors $Ext$ and $Tor$ in $\CC$ 180

5. Some calculations 184

\vspace{.3cm}

{\bf Chapter 3. Global sections}

\vspace{.3cm}

6. Braiding and balance in $\CC$ and $\FS$ 188

7. Global sections over $\CA(K)$ 189

8. Global sections over $\CP$ 190

9. Application to conformal blocks 193

\vspace{1cm}

{\bf Part V. MODULAR STRUCTURE ON THE CATEGORY $\FS$}

\vspace{1cm}

1. Introduction 197

\vspace{.3cm}

{\bf Chapter 1. Heisenberg local system}

\vspace{.3cm}

2. Notations and statement of the main result 199

3. The scheme of construction 202

4. The universal line bundle 204

5. The universal local system 207

6. Factorization isomorphisms 216

\vspace{.3cm}


{\bf Chapter 2. The modular property of the Heisenberg system}

\vspace{.3cm}

7. Degeneration of curves: recollections and notations 219

8. Proof of Theorem 7.6(a) 221

9. Proof of Theorem 7.6(b) 224

\vspace{.3cm}

{\bf Chapter 3. Regular representation}

\vspace{.3cm}

10. A characterization of the regular bimodule 231

11. The adjoint representation 233

\vspace{.3cm}

{\bf Chapter 4. Quadratic degeneration in genus zero}

\vspace{.3cm}

12. $I$-sheaves 235

13. Degenerations of quadrics 236

14. The $I$-sheaf $\CR$ 237

15. Convolution 239

\vspace{.3cm}

{\bf Chapter 5. Modular functor}

\vspace{.3cm}

16. Gluing over $C$ 242

17. Degeneration of factorizable sheaves 246

18. Global sections over $C$ 248

\vspace{.3cm}

{\bf Chapter 6. Integral representations of conformal blocks}

\vspace{.3cm}

19. Conformal blocks in arbitrary genus 249

\vspace{1cm}

{\bf References} 252

{\bf Index of Notation} 254

{\bf Index of Terminology} 283

\newpage
\pagenumbering{arabic}
\centerline{\bf Introduction}

\bigskip

The aim of this work is to present a geometric construction of modules 
over Lusztig's small quantum groups. To give a better idea of a problem, 
we first make some historical remarks. 

{\bf I.1.} 
The story begins with {\bf Knizhnik-Zamolodchikov} system of differential equations 
$$
\frac{\dpar\phi(z)}{\dpar z_i}=\frac{1}{\kappa}\sum_{j\neq i}\ 
\frac{\Omega_{ij}\phi(z)}{z_i-z_j},
\eqno{(I.1)}
$$
$i=1,\ldots, n$. Here $\phi(z)=\phi(z_1,\ldots,z_n)$ is a function 
of complex variables with values in a tensor product $M=M_1\otimes M_2\otimes
\ldots M_n$ of $\fg$-modules, $\fg$ being a semisimple complex Lie algebra. 
One fixes an invariant symmetric inner product on $\fg$, given 
by a tensor $\Omega\in\fg\otimes\fg$. For $i\neq j$, $\Omega_{ij}$ is 
the endomorphism of $M$ induced by the multiplication by $\Omega$, acting 
on $M_i\otimes M_j$, and $\kappa$ is a non-zero complex parameter. The system 
(I.1) is of primary importance.  
It was discovered by physicists, [KZ], as a system of differential 
equations on correlation functions in Wess-Zumino-Witten models of 
two dimensional conformal field theory. 

As was discovered in the other physical paper [DF] (based on the fundamental 
work by Feigin and Fuchs, [FF]), the correlation functions 
in conformal field theories may be expressed in terms of the following 
integrals of hypergeometric type
$$
I(z_1,\ldots,z_n)=\int r(z_1,\ldots,z_n; t_1,\ldots,t_N)\prod\ (z_i-t_j)^{a_{ij}}\prod\ 
(t_p-t_q)^{b_{pq}} dt
\eqno{(I.2)}
$$
where $r(z;t)$ is a rational function, $a_{ij}, b_{pq}$ are complex numbers. 
So, we are dealing with integrals over homology cycles of some 
one-dimensional local systems over a {\it configuration space}  
$$
X^N(z_1,\ldots,z_n)=\BC^N\setminus\bigcup_{1\leq i\leq n;\ 1\leq j\leq N}\ 
\{z_i=t_j\}\setminus\bigcup_{1\leq p<q\leq N}\ \{t_p=t_q\}
\eqno{(I.3)}
$$
(the coordinates on $\BC^N$ being $t_i$), 
these integrals depending on $z_1,\ldots,z_n$ as on parameters. The integrals 
$I(z_1,\ldots,z_n)$, as functions 
of $z$, satisfy the {\it Gauss-Manin} differential equations. These results 
have been obtained by the powerful technique of {\it bosonization} 
in the representation theory of Virasoro and Kac-Moody Lie algebras.  

Inspired by these (and other) physical works, one wrote down in [SV1] 
the complete (for generic $\kappa$) set of solutions of KZ 
equations in terms of the integrals (I.2). In this paper, an interesting 
phenomenon has been found. It turned out that 
the the contragradient Verma modules over $\fg$ and the irreducible ones 
could be realized in some spaces of logarithmic differential forms 
over the spaces $X^N(z)$ (in the same manner, their tensor products 
appeared over the spaces $X^N(z_1,\ldots,z_N)$). The picture  
resembled very much the classical Beilinson-Bernstein theory  of 
localization of $\fg$-modules, [B1]. The role of the flag space $G/B$ is 
played by the above configuration spaces; in order to restore the whole 
module, one has to consider the spaces $X^N(z)$ for all $N$. The Verma modules 
are connected with the $!$-extensions of some standard local systems, 
their contragradient duals --- with $*$-extensions.         

The other part of the story is the important {\bf Kohno-Drinfeld} theorem, [K], 
[D], which says that the monodromy of the solutions of (I.1) is expressible 
in terms of the $R$-matrix of the quantized enveloping algebra $U_q\fg$. 
The quantization parameter $q$ is connected with the parameter $\kappa$ in 
(I.1) through the identity 
$$
q=\exp(2\pi i/\kappa)
\eqno{(I.4)}
$$
Another proof of this theorem follows from the results of [SV1], since   
the monodromy of the integrals (I.2) may be explicitely calculated. Moreover, 
it turned out that the results of [SV1] may be transferred into the 
topological world, with cycles replacing differential forms, and this way 
one got the realization of standard representations of the quantum group 
$U_q\fg$, cf. [SV2]. As B.Feigin put it at that time, "differential forms 
are responsible for the Lie algebra $\fg$, and cycles --- for the quantum 
group". 

The idea to realize standard representations of quantum groups in certain
cohomology of configuration spaces appeared also in [FW]. The important
role of the adjoint representation and certain higher-dimensional local
systems over the configuration spaces of curves of higher genus was first
revealed in [CFW].

One noted (cf. [S]) that in the topological world, the standard 
representations of a quantum group  
are realized in certain spaces of {\it vanishing cycles}, and the 
operators of the quantum group as the Lefschetz-Deligne {\it canonical} 
and {\it variation} maps between them. This implied the idea that 
some categories of representations of quantum groups could be realized geometrically as categories of compatible systems of perverse sheaves 
(smooth along some natural stratifications) 
over various configuration spaces. Here the main point of inspiration 
was Beilinson's {\it Gluing theorem}, allowing one to glue perverse sheaves 
from various spaces of vanishing cycles, [B2]. In the present work, 
we suggest a construction of such geometric category.

{\bf I.2.} Here is a brief guide through the text. It consists of six parts. 
Part 0 contains a detailed overview of Parts I --- V, which constitute the main 
body of the work. We recommend 
the reader to look through this part. (Practically) all the results of 
the work are precisely formulated, and the exposition may be in some respects 
more clear. We also recommend to look through the introductions to all Parts,  
where their contents is described.  

The plan of Parts I and II resembles that of [SV1]. In part I, we give some 
general results concerning perverse sheaves on the affine plane smooth 
along the stratification associated by an arrangement of hyperplanes. 
This is a nice subject in itself, maybe of independent interest. To each 
perverse sheaf as above, we assign a linear algebra datum, by means 
of which one can compute the cohomology of this sheaf. 

Part II contains some important results concerning the perverse 
sheaves on configuration spaces. Technically, it is a base of our main results. 
In particular, it contains the proof of all results announced in [SV2] 
(so, it may be instructive to look through {\it op. cit.} before 
reading Part II). 

Part III contains our main definition --- of the tensor category $\FS$ 
of {\bf factorizable sheaves}. A factorizable sheaf is certain compatible 
collection of perverse sheaves over configuration spaces, cf. Definition 4.2. 
We recommend to read Part 0, Introduction, 1.2, where this definition is explained in 
a slightly different language. This is the main hero of the present work.  

The main theorem of this work   
claims that this category is equivalent to 
the category $\CC$ connected with Lusztig's small quantum group, cf. Theorem 
17.1. See also Theorem 18.4 for a similar result, but with   
the generic quantization parameter $\zeta$. 

In Parts IV and V, we want to globalize our constructions. The category 
$\FS$ is a local object --- factorizable sheaves leave on a (formal) disk. 
As it is usual in conformal field theory, we can place these sheaves 
into points of an arbitrary algebraic curve $C$. As a result, we get 
some interesting perverse sheaves on symmetric powers of $C$. 

In Part IV, we study the case of a projective line. It turns out that 
the cohomology of the above "glued" perverse sheaves are expressible 
in terms of the Arkhipov's "semiinfinite" cohomology of the quantum 
group, cf. Theorem 8.4. 
In Part V, we study the case of an arbitrary curve $C$, and families of curves. 
As an application, we get the semisimplicity of the local systems of conformal 
blocks in Wess-Zumino-Witten models, cf. Corollary 19.9. 

{\bf I.3.} In this work, we live in the topological world. 
To go back to the world of 
differential equations, one should pass from perverse sheaves to 
$\CCD$-modules. If one translates our constructions to the  
language of $\CCD$-modules, one gets a geometric description of 
$\fg$-modules (at least for a generic $\kappa$), cf. [KS1], [KS2]. 
If one could say 
that the present work puts [SV2] in a natural framework, then [KS1], [KS2] 
play the similar role for [SV1]. Kohno-Drinfeld theorem then 
translates into Riemann-Hilbert correspondence.     

{\bf I.4.} We believe that the geometry studied in this work is 
an instance of the
{\em semiinfinite geometry of the configuration space of a curve}. We learnt
of this concept from A.~Beilinson in 1991. It is intimately, 
though mysteriously,
related to the geometry of {\em semiinfinite flag space} we learnt 
from B.~Feigin since 1989. For example, the definition of the category $\cal{PS}$ of
factorizable sheaves on the semiinfinite flag space just copies the definition
of factorizable sheaves in the present work, cf. [FM]. Also, the semiinfinite
cohomology of the small quantum group arise persistently in the computation
of IC-hypercohomology of both configuration spaces and Quasimaps' spaces
(finite-dimensional approximations of the semiinfinite flag space), cf.
[FK] and [FFKM]. The more striking parallels between the two theories are
accounted for in [FKM].

{\bf I.5.} The ideas of A.~Beilinson and B.~Feigin obviously penetrate this whole work. 
During the preparation of the text, the authors benefited  
from the 
ideas and suggestions of P.~Deligne, G.~Felder and
D.~Kazhdan. 
We are very obliged to
G.~Lusztig for the permission to use his unpublished results. The second 
author is happy to thank L.~Positselski whose friendly help came at the 
crucial moment.  
Our special gratitude goes to D.~Kazhdan and Yu.I.~Manin  
for the important incouragement at the initial (resp. final) stages of the 
work.

\bigskip

\centerline{\bf References}

\bigskip

[B1] A.~Beilinson, Localization of representations of reductive Lie 
algebras, in: Proceedings of the ICM Warsaw 1983, 699-710, PWN, 
Warsaw, 1984.

[B2] A.~Beilinson, How to glue perverse sheaves, in: K-theory, arithmetic and geometry (Moscow, 1984-1986), 42-51, {\it Lect. Notes in Math.} {\bf 1289}, 
Springer-Verlag, Berlin-Heidelberg-New York, 1987. 

[CFW] M.~Crivelli, G.~Felder, C.~Wieczerkowski, Generalized hypergeometric 
functions on the torus and the adjoint representation of $U_qsl_2$, 
{\it Comm. Math. Phys.} {\bf 154} (1993), 1-23.      

[DF] Vl.S.~Dotsenko, V.A.~Fateev, Conformal algebra and multipoint correlation 
functions in 2D statistical models, {\it Nucl. Phys.} {\bf B240} (1984), 
312-348.

[D] V.G.~Drinfeld, Quasi-Hopf algebras, Algebra i Analiz {\bf 1} 
(1989), 114-148 (russian). 

[FF] B.L.~Feigin, D.B.~Fuchs, Representations of the Virasoro algebra, in: 
Representations of Lie groups and related topics, {\it Adv. Stud. 
Contemp. Math.} {\bf 7}, Gordon and Breach, New York, 1990. 

[FFKM] B.Feigin, M.Finkelberg, A.Kuznetsov, I.Mirkovi\'c,
Semiinfinite Flags. II. Local and global Intersection Cohomology of Quasimaps
spaces, Preprint alg-geom/9711009.

[FW] G.~Felder, C.~Wieczerkowski, Topological representations of the quantum 
group $U_qsl_2$, {\it Comm. Math. Phys.} {\bf 138} (1991), 683-605. 

[FK] M.Finkelberg, A.Kuznetsov, Global Intersection Cohomology
of Quasimaps' spaces, {\em Int. Math. Res. Notices} {\bf 1997}, no. 7, 
301-328.

[FKM] M.Finkelberg, A.Kuznetsov, I.Mirkovi\'c, Frobenius realization of
$U(\check{\frak n})$, in preparation.

[FM] M.Finkelberg, I.Mirkovi\'c, Semiinfinite Flags. I. Case
of global curve $\BP^1$, Preprint alg-geom/9707010.

[KS1] S.~Khoroshkin, V.~Schechtman, Factorizable $\CCD$-modules, {\it Math. 
Res. Letters} {\bf 4} (1997), 239-257. 

[KS2] S.~Khoroshkin, V.~Schechtman, Non-resonance $\CCD$-modules 
over arrangements of hyperplanes, Preprint MPI 98-14; math.AG/9801134.  
   
[KZ] V.G.~Knizhnik, A.B.~Zamolodchikov, Current algebra and Wess-Zumino model 
in two dimensions, {\it Nucl Phys.} {\bf B247} (1984), 83-103.

[K] T.~Kohno, Monodromy representations of braid groups and Yang-Baxter 
equations, {\it Ann. Inst. Fourier} (Grenoble) {\bf 37} (1987), 139-160. 

[S] V.~Schechtman, Vanishing cycles and quantum groups. I, {\it Int. Math. 
Res. Notices} {\bf 1992}, no. 3, 39-49; II, {\it ibid.} {\bf 1992}, no. 10, 207-215.    

[SV1] V.~Schechtman, A.~Varchenko, Arrangements of hyperplanes and Lie 
algebra homology, {\it Inv. Math.} {\bf 106} (1991), 139-194. 

[SV2] V.~Schechtman, A.~Varchenko, Quantum groups and homology of local 
systems, in: Algebraic geometry and analytic geometry (Tokyo, 1990), 
182-197, ICM-90 Satell. Conf. Proc., Springer, Tokyo, 1991.

\newpage
\begin{center}
{\bf Acknowledgement}
\end{center}
\vspace{1cm}
During the work on this project, the first author visited Independent 
Moscow University, the second author visited Harvard University,
Stony Brook University and Tel-Aviv University, and the third
author visited Harvard University and Max-Planck-Institut f\"ur Mathematik. 
We thank these institutions 
for the hospitality.

The second author was partially supported by the
grants from AMS, INTAS94-4720, CRDF RM1-265.
The third author was partially supported by NSF.
\newpage 
\begin{center}
{\large \bf Part 0. OVERVIEW}
\end{center}

\section{Introduction}


\subsection{} Let $(I,\cdot)$ be a Cartan datum of finite type and  
$(Y,X,\ldots)$ the simply connected root datum of type $(I,\cdot)$,
cf. ~\cite{l1}.  

Let $l>1$ be an integer. Set $\ell=l/(l,2)$; to simplify the exposition,  
we will suppose in this Introduction that $d_i:=i\cdot i/2$ divides $\ell$ 
for all $i\in I$; we set $\ell_i:=\ell/d_i$. We suppose that all 
$\ell_i>3$.
Let $\rho\in X$ be the half-sum
of positive roots; let $\rho_{\ell}\in X$ be defined by $\langle i,\rho_{\ell}
\rangle =\ell_i-1\ (i\in I)$. We define a lattice $Y_{\ell}\subset Y$ by 
$Y_{\ell}=\{\mu\in X|\ \mbox{for all}\ \mu'\in X, \mu\cdot\mu'\in\ell\BZ\}$,
and set $\dd_{\ell}=\card(X/Y_{\ell})$.    

We fix a base field $\sk$ of characteristic not dividing $l$. We suppose that
$\sk$ contains a primitive $l$-th root of unity $\zeta$, and fix it.
Starting from these data, one defines certain category 
$\CC$. Its objects are finite dimensional $X$-graded $\sk$-vector spaces
equipped with an action of Lusztig's "small" quantum group $\fu$ 
(cf. ~\cite{l2}) such that the action of its Cartan subalgebra is compatible 
with the $X$-grading. Variant: one defines certain algebra $\dfu$ which is an
"$X$-graded" version of $\fu$ (see \ref{nul' dot}), and an object of
$\CC$ is a finite 
dimensional $\dfu$-module.  For the precise definition of $\CC$, see
~\ref{nul' def c}, \ref{nul' comp lu}.
For $l$ prime and $\chara \sk=0$, the category $\CC$ was studied in ~\cite{ajs},
for $\chara \sk>0$ and arbitrary $l$, $\CC$ was studied in \cite{aw}.
The category $\CC$ 
admits a remarkable structure of a {\em ribbon}\footnote{in other terminology, 
braided balanced rigid tensor} category (Lusztig).

\subsection{}
The main aim of this work is to introduce      
certain tensor category $\FS$ of geometric origin, which is  
equivalent to $\CC$. Objects of $\FS$ are called {\bf (finite) factorizable 
sheaves}. 

A notion of a factorizable sheaf is the main new concept of this work.
Let us give an informal idea, what kind of an object is it.
Let $D$ be the unit disk in the complex plane\footnote{One could also take   
a complex affine line or a formal disk; after a suitable modification 
of the definitions, the resulting categories of factorizable sheaves 
are canonically equivalent.}. Let
$\CCD$ denote the space of positive $Y$-valued divisors on $D$. Its points
are formal linear combinations $\sum\nu\cdot x\ (\nu\in Y^+:=\BN[I], x\in D)$.
This space is a disjoint union
$$
\CCD=\coprod_{\nu\in Y^+}\ D^{\nu}
$$
where $D^{\nu}$ is the subspace of divisors of degree $\nu$. Variant: $D^{\nu}$ 
is the configuration space of $\nu$ points running on $D$; its points are 
(isomorphism clases of) pairs of maps

($J\lra D$, $j\mapsto x_j$;\ $\pi: J\lra I$, $j\mapsto i_j$),
 
$J$ being a finite set.  
We say that we have a finite set $\{x_j\}$ of (possibly equal) 
points of $D$, a point $x_j$ being "coloured" by the colour $i_j$. The sum  
(in $Y$) of colours of all points should be equal to $\nu$. The space 
$D^{\nu}$ is a smooth affine analytic variety; it carries a canonical 
stratification defined by various intersections of hypersurfaces $\{x_j=0\}$ and 
$\{x_{j'}=x_{j''}\}$. The open stratum (the complement of all the hypersurfaces 
above) is denoted by $A^{\nu\circ}$.

One can imagine $\CCD$ as a $Y^+$-graded smooth stratified variety.
Let $\CA=\coprod A^{\nu}\subset\CCD$ be the open $Y^+$-graded subvariety
of positive $Y$-valued divisors on $D-\{0\}$. We have an open $Y^+$-graded 
subvariety $\CA^{\circ}=\coprod A^{\nu\circ}\subset\CA$.

Let us consider the $(Y^+)^2$-graded variety  
$$
\CA\times\CA=\coprod\ A^{\nu_1}\times A^{\nu_2};
$$
we define another $(Y^+)^2$-graded variety $\widetilde{\CA\times\CA}$,
together with two maps

(a) $\CA\times\CA\overset{p}{\lla}\widetilde
{\CA\times\CA}\overset{m}{\lra}\CA$

respecting the $Y^+$-gradings\footnote{a $(Y^+)^2$-graded space is considered 
as a $Y^+$-graded by means of the addition $(Y^+)^2\lra Y^+$.}; the map 
$p$ is a homotopy equivalence. One can imagine the   
diagram (a) above as a "homotopy multiplication"    
$$
m_{\CA}:\CA\times\CA\lra\CA;
$$
this "homotopy map" is "homotopy associative"; the meaning of this is explained
in the main text. We say that
$\CA$ is a ($Y^+$-graded) "homotopy monoid"; $\CA^{\circ}\subset\CA$ is a
"homotopy submonoid".

The space $\CCD$ is a "homotopy $\CA$-space": there is a
"homotopy map"   
$$
m_{\CCD}:\CCD\times\CA\lra\CCD
$$
which is, as above, a diagram of usual maps between $Y^+$-graded
varieties

(b) $\CCD\times \CA\overset{p}{\lla}\widetilde{\CCD\times\CA}
\overset{m}{\lra}\CCD$,

$p$ being a homotopy equivalence.  

For each $\mu\in X,\ \nu\in Y^+$, one defines a one-dimensional 
$\sk$-local system
$\CI_{\mu}^{\nu}$ over $A^{\nu\circ}$. Its monodromies are defined by variuos
scalar products of "colours" of running points and the colour $\mu$ of the 
origin $0\in D$.
These local systems have the following
compatibility. For each $\mu\in X, \nu_1,\nu_2\in Y^+$, a 
{\em factorization isomorphism}
$$
\phi_{\mu}(\nu_1,\nu_2): m^*\CI_{\mu}^{\nu_1+\nu_2}\iso p^*(\CI^{\nu_1}_{\mu}
\boxtimes\CI^{\nu_2}_{\mu-\nu_1})
$$
is given (where $m, p$ are the maps in the diagram (a) above), 
these isomorphisms satisfying certain {\em (co)associativity property}. The
collection of local systems $\{\CI^{\nu}_{\mu}\}$ is an $X$-graded
local system $\CI=\oplus\ \CI_{\mu}$ over $\CA^{\circ}$. One could imagine
the collection of factorization isomorphisms $\{\phi_{\mu}(\nu_1,\nu_2)\}$ 
as an isomorphism 
$$
\phi:\ m_{\CA^{\circ}}^*\CI\iso\CI\boxtimes\CI
$$
We call $\CI$ the {\em braiding local system}. 

Let $\CI^{\bullet}$ be the perverse sheaf on $\CA$ which is the
Goresky-MacPherson 
extension of the local system $\CI$. The isomorphism $\phi$ above induces 
a similar isomorphism
$$
\phi^{\bullet}:\ m_{\CA}^*\CI^{\bullet}\iso\CI^{\bullet}\boxtimes\CI^{\bullet}
$$
This sheaf, together with $\phi^{\bullet}$, looks like a "coalgebra"; it  
is an incarnation of the quantum group $\dfu$.  

A factorizable sheaf is a couple

($\mu\in X$ ("the highest weight");
a perverse sheaf $\CX$ over $\CCD$, smooth along the canonical stratification).

Thus, $\CX$ is a collection of sheaves $\CX^{\nu}$ over $\CCD^{\nu}$. These
sheaves should be connected by {\em factorization isomorphisms}
$$
\psi(\nu_1,\nu_2):\ m^*\CX^{\nu_1+\nu_2}\iso p^*(\CX^{\nu_1}
\boxtimes\CI^{\nu_2\bullet}_{\mu-\nu_1})
$$
satisfying an associativity property. Here $m, p$ are as in the diagram (b) 
above. One could imagine the whole collection $\{\psi(\nu_1,\nu_2)\}$ as  
an isomorphism 
$$
\psi: m_{\CCD}^*\CX\iso\CI^{\bullet}\boxtimes\CX
$$
satisfying a (co)associativity property. We impose also certain finiteness 
(of singularities) condition on $\CX$. So, this object looks like a
"comodule" over $\CI^{\bullet}$. It is an incarnation of an $\dfu$-module. 

We should mention one more important part of the structure of the space
$\CCD$. It comes from natural closed embeddings
$\iota_{\nu}:\ \CCD\hra\CCD[\nu]$ (were $[\nu]$ denotes the
shift of the grading);
these mappings define certain inductive system, and
a factorizable sheaf is a sheaf on its inductive limit.

The latter inductive 
limit is an example of a {\bf "semiinfinite space"}.     

For the precise definitions, see Sections \ref{nul' braiding}, \ref{nul' factoriz}
below. 
  
The category $\FS$ has a structure of a  
braided balanced tensor category coming from geometry. The  
tensor structure on $\FS$ is defined using the functors of nearby cycles. 
The tensor equivalence
$$
\Phi:\FS\iso\CC
$$
is defined using vanishing cycles functors. It respects the  
braidings and balances. This is the contents of Parts I--III.

\subsection{} 
\label{nul' glue intro} Factorizable sheaves are local objects. It turns out that
one can "glue" them along complex curves. More precisely,
given a finite family $\{\CX_a\}_{a\in A}$ of factorizable sheaves
and a smooth proper curve $C$ together
with a family $\{x_a,\tau_a\}_A$ of distinct points $x_a\in C$ with
non-zero tangent vectors $\tau_a$ at them (or, more generally, a family of such 
objects over a smooth base $S$), one can define a perverse 
sheaf, denoted by   
$$
\Boxtimes_A^{(C)}\ \CX_a
$$
on the (relative) configuration space $C^{\nu}$. Here $\nu\in Y^+$ is 
defined by

(a) $\nu=\sum_{a\in A}\ \mu_a+(2g-2)\rho_{\ell}$

where $g$ is the genus of $C$, $\mu_a$ is the highest weight of $\CX_a$. 
We {\em assume} that the right hand side of the equality belongs to $Y^+$.

One can imagine this sheaf as an "exterior tensor product"
of the family $\{\CX_a\}$ along $C$. It is obtained by "planting" the 
sheaves $\CX_a$ into the points $x_a$. To glue them together, one needs 
a "glue". 
This glue is called the {\bf Heisenberg local system} $\CH$; it is a sister
of the braiding local system $\CI$. Let us describe what it is. 

For a finite set $A$, let $\CM_A$ denote the moduli stack of
punctured curves $(C,\{x_a,\tau_a\}_A)$ as above; let $\eta: C_A\lra\CM_A$ be
the universal curve. Let $\delta_A=\det (R\eta_*\CO_{C_A})$ be the determinant
line bundle on $\CM_A$, and $\CM_{A;\delta}\lra\CM_A$ be the total space
of $\delta_A$ with the zero section removed. For $\nu\in Y^+$, let
$\eta^{\nu}: C^{\nu}_A\lra\CM_A$ be the relative configuration space
of $\nu$ points running on $C_A$; let
$C^{\nu\circ}_A\subset C^{\nu}_A$ be the
open stratum (where the running points are distinct and distinct from
the punctures $x_a$). 
The complementary subscript $()_{\delta}$ will denote the base
change under $\CM_{A;\delta}\lra\CM_A$. The complementary subscript $()_g$
will denote the base change under $\CM_{A,g}\hra\CM_A$, $\CM_{A,g}$ being
the substack of curves of genus $g$.

The Heisenberg local system is a collection of local systems
$\CH_{\vmu;A,g}$ over the stacks $C^{\nu\circ}_{A,g;\delta}$. Here
$\vmu$ is an $A$-tuple $\{\mu_a\}\in X^a$; $\nu=\nu(\vmu;g)$ is defined
by the equality (a) above. We assume that $\vmu$ is such that the right
hand side of this equality really belongs to $Y^+$. The dimension of 
$\CH_{\vmu;A,g}$ is equal to $\dd_{\ell}^g$; the monodromy around the
zero section of the determinant line bundle is equal to

(b) $c=(-1)^{\card(I)}\zeta^{-12\rho\cdot\rho}$.

These local systems have a remarkable compatibility ("fusion") property 
which we do not specify here, see Section \ref{nul' fusion}.

\subsection{}    
In Part IV we study the sheaves $\Boxtimes_A^{(C)}\ \CX_a$ for
$C=\BP^1$. Their cohomology, when expressed algebraically,
turn out to coincide with certain "semiinfinite" $\Tor$ spaces in
the category $\CC$ introduced by S.M.Arkhipov. Due to the results
of Arkhipov, this enables one to prove that
"spaces of conformal blocks" in WZW models are the natural subquotients of
such cohomology spaces. 

\subsection{} 
In Part V we study the sheaves $\Boxtimes_A^{(C)}\ \CX_a$ for arbitrary
smooth families of punctured curves. Let $\Boxtimes_{A,g}\ \CX_a$ denotes
the "universal exterior product"  
living on $C^{\nu}_{A,g;\delta}$ where $\nu=\nu(\vmu;g)$ is as
in \ref{nul' glue intro} (a) above,
$\vmu=\{\mu_a\},\ \mu_a$ being the highest weight of $\CX_a$. Let us integrate 
it: consider 
$$
\int_{C^{\nu}}\ \Boxtimes_{A,g}\ \CX_a:=R\eta^{\nu}_*(\Boxtimes_{A,g}\ \CX_a);
$$
it is a complex of sheaves over $\CM_{A,g;\delta}$ with smooth cohomology;
let us denote it $\langle\otimes_A\ \CX_a\rangle_g$. The cohomology sheaves
of such complexes     
define a {\em fusion structure}\footnote{actually 
a family of such structures (depending on the degree of cohomology)} on $\FS$ 
(and hence on $\CC$).

The classical WZW model fusion category 
is, in a certain sense, a subquotient of one of them. The number $c$,
\ref{nul' glue intro} (b),
coincides with the "multiplicative central charge" of the model (in the
sense of \cite{bfm}).  
 
As a consequence of this geometric description and of the Purity Theorem,
\cite{bbd}, 
the local systems of conformal blocks (in arbitrary genus)  
are semisimple. The Verdier duality induces a canonical non-degenerate
Hermitian form on them (if $\sk=\BC$).


\subsection{} We should mention a very interesting 
related work of G.Felder and collaborators.  
The idea of realizing quantum 
groups' modules in the cohomology of configuration spaces appeared 
independently in \cite{fw}. The Part V of the present work
arose from our attempts to understand \cite{cfw}.

\subsection{}
This work consists of 6 parts. Part 0 is an overview of the Parts I--V.
More precisely, Chapter 0.1 surveys the Parts I--III; Chapter 0.2 is an
exposition of Part IV, and Chapter 0.3 is an exposition of Part V.
The proofs are given in the Parts I--V. Part 0 contains no proofs, but we
hope that it clarifies the general picture to some extent.
Each part starts with a separate Introduction. The references within one Part
look like 1.1.1, and the references to other parts look like I.1.1.1.
The first author would like to stress that his contribution is limited to
chapters 1,2 of Part V. The remaining authors would like to note that these
chapters consitute the core of Part V.



\newpage
\begin{center}
{\bf Chapter 1. Local.}
\end{center} 

\section{The category $\CC$} 

\subsection{} 
\label{nul' root} We will follow Lusztig's terminology
and notations concerning root systems, cf. ~\cite{l1}.  

We fix an irreducible Cartan datum $(I,\cdot)$ 
of finite type. Thus, $I$ is a finite set together with a symmetric
$\BZ$-valued bilinear form $\nu_1,\nu_2\mapsto\nu_1\cdot\nu_2$
on the free abelian group $\BZ[I]$, such that for any $i\in I$,
$i\cdot i$ is even and positive, for any $i\neq j$ in $I$, 
$2i\cdot j/i\cdot i\leq 0$ and the $I\times I$-matrix $(i\cdot j)$ is 
positive definite. We set $d_i=i\cdot i/2$.

Let $d=\max_{i\in I}\ d_i$. This number is equal to the least 
common multiple of the numbers $d_i$, and belongs to the set $\{1,2,3\}$.
For a simply laced root datum $d=1$.  

We set $Y=\BZ[I], X=\Hom(Y,\BZ)$; $\langle,\rangle:\ Y\times X\lra
\BZ$ will denote the obvious pairing. The obvious
embedding $I\hra Y$ will be denoted by $i\mapsto i$. We will denote by 
$i\mapsto i'$  the embedding 
$I\hra X$ given by $\langle i,j'\rangle=2i\cdot j/i\cdot i$. 
(Thus, $(Y,X,\ldots)$ is the simply connected root datum of type $(I,\cdot)$, 
in the terminology of Lusztig.) 

The above embedding $I\subset X$ extends by additivity to the embedding 
$Y\subset X$. We will regard $Y$ as the sublattice of $X$ by means
of this embedding.  
For $\nu\in Y$, we will denote by the same letter $\nu$
its image in $X$. We set $Y^+=\BN[I]\subset Y$.

We will use the following partial order on $X$. For $\mu_1, \mu_2\in X$,
we write $\mu_1\leq \mu_2$ if $\mu_2-\mu_1$ belongs to $Y^+$.    

\subsection{} 
\label{nul' l} We fix a base field $\sk$, an integer $l>1$ and a primitive root
of unity $\zeta\in \sk$ as in the Introduction.
 
We set $\ell=l$ if $l$ is odd and $\ell=l/2$ 
if $l$ is even. For $i\in I$, we set $\ell_i=\ell/(\ell,d_i)$. Here 
$(a,b)$ stands for the greatest common divisor of $a,b$. We set 
$\zeta_i=\zeta^{d_i}$. We will assume that $\ell_i>1$ for any $i\in I$ and 
$\ell_i>-\langle i,j'\rangle+1$ for any $i\neq j$ in $I$. 

We denote by $\rho$ (resp. $\rho_{\ell}$) the element of $X$ such that 
$\langle i,\rho\rangle=1$ (resp.  
$\langle i,\rho_{\ell}\rangle=\ell_i-1$) for all $i\in I$. 

For a coroot $\beta\in Y$, there exists an element $w$ of the Weyl
group $W$ of our Cartan datum  
and $i\in I$ such that $w(i)=\beta$. We set $\ell_{\beta}:=\frac{\ell}
{(\ell,d_i)}$; this number does not depend on the choice of $w$ and $i$. 
 
We have $\rho=\frac{1}{2}\sum\ \alpha;\ \rho_{\ell}=\frac{1}{2}\sum\
(\ell_{\alpha}-1)\alpha$, the sums over all positive roots $\alpha\in X$.  

For $a\in\BZ, i\in I$, we set $[a]_i=1-\zeta_i^{-2a}$.

\subsection{}
\label{nul' lattice} We use the same notation $\mu_1,\mu_2\mapsto\mu_1\cdot\mu_2$
for the unique extension of the bilinear form on $Y$ to a $\Bbb{Q}$-valued
bilinear form on $Y\otimes_{\BZ}\Bbb{Q}=X\otimes_{\BZ}\Bbb{Q}$.

We define a lattice $Y_{\ell}=\{\lambda\in X|\mbox{ for all }\mu\in X,
\lambda\cdot\mu\in\ell\BZ\}$.

\subsection{}
Unless specified otherwise, 
a "vector space" will mean a vector space over $\sk$; $\otimes$ will
denote the tensor product over $\sk$. A "sheaf" (or a "local system")
will mean a sheaf (resp. local system) of $\sk$-vector spaces.

If $(T,\CS)$ is an open subspace of
the space of complex points of a separate scheme of 
finite type over $\BC$, with the usual topology,
together with an algebraic stratification $\CS$ satisfying the properties 
\cite{bbd} 2.1.13 b), c), we will denote by 
$\CM(T;\CS)$ the category of perverse sheaves over $T$ lisse along 
$\CS$, with respect to the middle perversity, cf. \cite{bbd} 2.1.13, 
2.1.16. 

\subsection{} Let $'\ff$ be the free associative $\sk$-algebra with $1$
with generators $\theta_i\ (i\in I)$. For $\nu=\sum \nu_ii\in\BN[I]$,
let $'\ff_{\nu}$ be the subspace of $'\ff$ spanned by the monomials 
$\theta_{i_1} 
\cdot\ldots\cdot\theta_{i_a}$ such that $\sum_{j} i_j=\nu$ in $\BN[I]$.

Let us regard $'\ff\otimes\ '\ff$
as a $\sk$-algebra with the product $(x_1\otimes x_2)(y_1\otimes y_2)=
\zeta^{\nu\cdot\mu}x_1y_1\otimes x_2y_2\ (x_2\in\ '\ff_{\nu}, y_1\in\ 
'\ff_{\mu}$). Let $r$ denote a unique homomorphism of $\sk$-algebras
$'\ff\lra\ '\ff\otimes\ '\ff$ carrying $\theta_i$ to $1\otimes\theta_i+
\theta_i\otimes 1\ (i\in I)$.

\subsection{Lemma-definition} {\em There exists a unique $\sk$-valued bilinear
form 
$(\cdot,\cdot)$ on $'\ff$ such that 

{\em (i)} $(1,1)=1;\ (\theta_i,\theta_j)=\delta_{ij}\ (i,j\in I)$;
{\em (ii)} $(x,yy')=(r(x),y\otimes y')$ for all $x,y,y'\in\ '\ff$.  

This bilinear form is symmetric.} $\Box$

In the right hand side of the equality (ii) we use the same notation 
$(\cdot,\cdot)$ for the bilinear form on $'\ff\otimes\ '\ff$ 
defined by $(x_1\otimes x_2,y_1\otimes y_2)=(x_1,y_1)(x_2,y_2)$. 

The radical of the form $'\ff$ is a two-sided ideal of $'\ff$. 

\subsection{} Let us consider the associative $\sk$-algebra $\fu$
(with $1$) defined by the generators $\epsilon_i, \theta_i\ (i\in I), 
K_{\nu}\ (\nu\in Y)$ and the relations (a) --- (e) below. 

(a) $K_0=1,\ K_{\nu}\cdot K_{\mu}=K_{\nu+\mu}\ (\nu,\mu\in Y)$; 

(b) $K_{\nu}\epsilon_i=\zeta^{\langle\nu,i'\rangle}\epsilon_iK_{\nu}\ 
(i\in I, \nu\in Y)$;

(c) $K_{\nu}\theta_i=\zeta^{-\langle\nu,i'\rangle}\theta_iK_{\nu}\ 
(i\in I, \nu\in Y)$; 

(d) $\epsilon_i\theta_j-\zeta^{i\cdot j}\theta_j\epsilon_i=\delta_{ij}
(1-\tK_i^{-2})\ (i,j\in I)$. 

Here we use the notation $\tK_{\nu}=\prod_i\ K_{d_i\nu_ii}\ (\nu=\sum\nu_ii)$.

(e) If $f(\theta_i)\in\ '\ff$ belongs to the radical of the form 
$(\cdot,\cdot)$ then $f(\theta_i)=f(\epsilon_i)=0$ in $\fu$. 

\subsection{} There is a unique $\sk$-algebra homomorphism $\Delta: \fu\lra
\fu\otimes\fu$ such that 
$$
\Delta(\epsilon_i)=\epsilon_i\otimes 1+\tK_i^{-1}\otimes\epsilon_i;\
\Delta(\theta_i)=\theta_i\otimes 1+\tK_i^{-1}\otimes\theta_i;\
\Delta(K_{\nu})=K_{\nu}\otimes K_{\nu}
$$
for any $i\in I, \nu\in Y$. Here $\fu\otimes\fu$ is regarded as an algebra
in the standard way. 

There is a unique $\sk$-algebra homomorphism $e:\fu\lra \sk$
such that $e(\epsilon_i)=e(\theta_i)=0, e(K_{\nu})=1\ (i\in I, \nu\in Y)$. 

\subsection{} There is a unique $\sk$-algebra homomorphism $s:\fu\lra\fu^{
\opp}$ such that 
$$
s(\epsilon_i)=-\epsilon_i\tK_i;\ s(\theta_i)=-\theta_i\tK_i;\
s(K_{\nu})=K_{-\nu}\ (i\in I, \nu\in Y).
$$
There is a unique $\sk$-algebra homomorphism $s':\fu\lra\fu^{\opp}$ such that
$$
s'(\epsilon_i)=-\tK_i\epsilon_i;\ s'(\theta_i)=-\tK_i\theta_i;\
s'(K_{\nu})=K_{-\nu}\ (i\in I, \nu\in Y).
$$

\subsection{} The algebra $\fu$ together with the additional structure given 
by the comultiplication $\Delta$, the counit $e$, the antipode $s$ and
the skew-antipode $s'$, is a Hopf algebra.

\subsection{}
\label{nul' def c} Let us define a category $\CC$ as follows. An object
of $\CC$ is a $\fu$-module $M$ which is finite dimensional over $\sk$,
with a given direct sum decomposition $M=\oplus_{\lambda\in X}\ M_{\lambda}$
(as a vector space) such that $K_{\nu}x=\zeta^{\langle\nu,\lambda\rangle}
x$ for any $\nu\in Y, \lambda\in X, x\in M_{\lambda}$. A morphism in $\CC$
is a $\fu$-linear map respecting the $X$-gradings. 

Alternatively, an object of $\CC$ may be defined as an $X$-graded finite
dimensional vector space $M=\oplus\ M_{\lambda}$ equipped with linear
operators
$$
\theta_i: M_{\lambda}\lra M_{\lambda-i'},\
\epsilon_i: M_{\lambda}\lra M_{\lambda+i'}\ (i\in I, \lambda\in X)
$$
such that 

(a) for any $i,j\in I, \lambda\in X$, the operator 
$\epsilon_i\theta_j-\zeta^{i\cdot j}\theta_j\epsilon_i$ acts as the 
multiplication by 
$\delta_{ij}[\langle i,\lambda\rangle]_i$ on $M_{\lambda}$.

Note that $[\langle i,\lambda\rangle]_i=
[\langle d_ii,\lambda\rangle]=[i'\cdot\lambda]$. 

(b) If $f(\theta_i)\in\ '\ff$ belongs to the radical of the form 
$(\cdot,\cdot)$ then the operators $f(\theta_i)$ and $f(\epsilon_i)$ act 
as zero on $M$. 

\subsection{} In ~\cite{l2} Lusztig defines an algebra $\bu_{\CB}$ over the
ring $\CB$ which is a quotient of $\BZ[v,v^{-1}]$ ($v$ being an
indeterminate) by the $l$-th cyclotomic polynomial. Let us consider 
the $\sk$-algebra $\bu_\sk$ obtained from $\bu_{\CB}$ by the base change
$\CB\lra \sk$ sending $v$ to $\zeta$. The algebra $\bu_\sk$ is generated
by certain elements $E_i, F_i, K_i\ (i\in I)$. Here $E_i=E_{\alpha_i}^{(1)},
F_i=F_{\alpha_i}^{(1)}$ in the notations of {\em loc.cit.} 

Given an object $M\in\CC$, let us introduce the operators $E_i, F_i, K_i$ on 
it by 
$$
E_i=\frac{\zeta_i}{1-\zeta_i^{-2}}\epsilon_i\tK_i,\ F_i=\theta_i, K_i=K_i. 
$$

\subsection{Theorem} 
\label{nul' comp lu} {\em The above formulas define the action
of the Lusztig's algebra {\em $\bu_\sk$} on an object $M$.

This rule defines an equivalence of $\CC$ with the category whose 
objects are $X$-graded finite dimensional {\em $\bu_\sk$}-modules
$M=\oplus M_{\lambda}$ such that $K_ix=\zeta^{\langle i,\lambda\rangle}x$ 
for any $i\in I, \lambda\in X, x\in M_{\lambda}$.} $\Box$

\subsection{} The structure of a Hopf algebra on $\fu$ defines canonically 
a {\em rigid tensor} structure on $\CC$ (cf. \cite{kl}IV, Appendix).

The Lusztig's algebra $\bu_\sk$ also has an additional structure of a Hopf
algebra. It induces the same rigid tensor structure on $\CC$.

We will denote the duality in $\CC$ by $M\mapsto M^*$. The unit 
object will be denoted by $\One$.  

\subsection{} Let $\fu^-$ (resp. $\fu^+$, $\fu^0$) denote the $\sk$-subalgebra
generated by the elements $\theta_i\ (i\in I)$ (resp. $\epsilon_i\ 
(i\in I)$, $K_{\nu}\ (\nu\in Y)$). We have the triangular decomposition 
$\fu=\fu^-\fu^0\fu^+=\fu^+\fu^0\fu^-$. 

We define the "Borel" subalgebras $\fu^{\leq 0}=\fu^-\fu^0,\ \fu^{\geq 0}=
\fu^+\fu^0$; they are the Hopf subalgebras of $\fu$.   

Let us introduce the $X$-grading 
$\fu=\oplus\fu_{\lambda}$ as a unique grading compatible with 
the structure of an algebra such that $\theta_i\in\fu_{-i'},\  
\epsilon_i\in\fu_{i'},\ K_{\nu}\in\fu_0$. We will use the induced 
gradings on the subalgebras of $\fu$. 

\subsection{} 
\label{nul' ind} Let $\CC^{\leq 0}$ (resp. $\CC^{\geq 0}$) be the category
whose objects are $X$-graded finite dimensional $\fu^{\leq 0}$-
(resp. $\fu^{\geq 0}$-) modules $M=\oplus M_{\lambda}$ such that
$K_{\nu}x=\zeta^{\langle\nu,\lambda\rangle}x$ for any $\nu\in Y, 
\lambda\in X, x\in M_{\lambda}$. Morphisms are $\fu^{\leq 0}$-
(resp. $\fu^{\geq 0}$-) linear maps compatible with the $X$-gradings. 

We have the obvious functors $\CC\lra\CC^{\leq 0}$ (resp. 
$\CC\lra\CC^{\geq 0}$). These functors admit the exact 
left adjoints $\ind_{\fu^{\leq 0}}^{\fu}: \CC^{\leq 0}\lra\CC$ 
(resp. $\ind_{\fu^{\geq 0}}^{\fu}: \CC^{\geq 0}\lra\CC$). 

For example, $\ind_{\fu^{\geq 0}}^{\fu}(M)=\fu\otimes_{\fu^{\geq 0}}M$.
The triangular decomposition induces an isomorphism of graded  
vector spaces $\ind_{\fu^{\geq 0}}^{\fu}(M)\cong \fu^-\otimes M$.

\subsection{} For $\lambda\in X$, let us consider on object 
$\sk^{\lambda}\in\CC^{\geq 0}$ defined as follows. As a graded vector space,
$\sk^{\lambda}=\sk^{\lambda}_{\lambda}=\sk$. The algebra $\fu^{\geq 0}$ acts on
$\sk^{\lambda}$ as follows: $\epsilon_ix=0, K_{\nu}x=\zeta^{\langle\nu,\lambda
\rangle}x\ (i\in I, \nu\in Y, x\in \sk^{\lambda})$.

The object $M(\lambda)=\ind_{\fu^{\geq 0}}^{\fu}(\sk^{\lambda})$ is called
a {\em (baby) Verma module}. Each $M(\lambda)$ has a unique 
irreducible quotient
object, to be denoted by $L(\lambda)$. The objects $L(\lambda)\ (\lambda\in X)$ 
are mutually non-isomorphic and every irreducible object in $\CC$ is 
isomorphic to one of them. Note that the category $\CC$ is {\em artinian}, 
i.e. each object of $\CC$ has a finite filtration with irreducible 
quotients.  

For example, $L(0)=\One$.  

\subsection{} Recall that a {\em braiding} on $\CC$ is a collection 
of isomorphisms
$$
R_{M,M'}: M\otimes M'\iso M'\otimes M\ (M,M'\in\CC)
$$
satisfying certain compatibility with the tensor structure (see \cite{kl}IV,
A.11).

A {\em balance} on $\CC$ is an automorphism of the identity functor
$b=\{ b_M: M\iso M\ (M\in\CC)\}$ such that for every $M,N\in\CC$, 
$b_{M\otimes N}\circ(b_M\otimes b_N)^{-1}=R_{N,M}\circ R_{M,N}$ 
(see {\em loc. cit.}). 

\subsection{}
\label{nul' k big} Let $\varpi$ denote the determinant of the $I\times I$-matrix
$(\langle i,j'\rangle)$. From now on we 
assume that $\chara(\sk)$ does not divide $2\varpi$, and
$\sk$ contains an element $\zeta'$ such
that $(\zeta')^{2\varpi}=\zeta$; we fix such an element $\zeta'$.
For a number $q\in\frac{1}{2\varpi}\BZ$, $\zeta^q$ will denote
$(\zeta')^{2\varpi q}$.

\subsection{Theorem} (G. Lusztig) {\em There exists a unique braided 
structure $\{R_{M,N}\}$ on the tensor category $\CC$ such that for any
$\lambda\in X$ and $M\in\CC$, if $\mu\in X$ is such that   
$M_{\mu'}\neq 0$ implies $\mu'\leq\mu$, then
$$
R_{L(\lambda),M}(x\otimes y)
=\zeta^{\lambda\cdot\mu}y\otimes x
$$
for any $x\in L(\lambda), y\in M_{\mu}$.} $\Box$

\subsection{} Let $n:X\lra\frac{1}{2\varpi}\BZ$ be the function defined
by 
$$
n(\lambda)=\frac{1}{2}\lambda\cdot\lambda-\lambda\cdot\rho_{\ell}.
$$

\subsection{Theorem} {\em There exists a unique balance $b$ on $\CC$ 
such that 
for any $\lambda\in X$, $b_{L(\lambda)}=\zeta^{n(\lambda)}$.} $\Box$

\subsection{} The rigid tensor category $\CC$, together with the additional 
structure given by the above braiding and balance, is a {\em ribbon 
category} in the sense of Turaev, cf. \cite{k} and references
therein.

\section{Braiding local systems}
\label{nul' braiding}

\subsection{} 
\label{nul' color} For a topological space $T$ and a finite set $J$, $T^J$
will denote the space of all maps $J\lra T$ (with the topology 
of the cartesian product). Its points are $J$-tuples 
$(t_j)$ of points of $T$. We denote by $T^{Jo}$ the subspace consisting
of all $(t_j)$ such that for any $j'\neq j''$ in $J$, $t_{j'}\neq t_{j''}$.   

Let $\nu=\sum \nu_ii\in Y^+$. Let us call an {\em unfolding} of 
$\nu$ a map of finite sets $\pi: J\lra I$ such that $\card(\pi^{-1}(i))=
\nu_i$ for all $i$. Let $\Sigma_{\pi}$ denote the group of all 
automorphisms $\sigma: J\iso J$ such that $\pi\circ\sigma=\pi$. 

The group $\Sigma_{\pi}$ acts on the space $T^J$ in the obvious way. 
We denote by $T^{\nu}$ the quotient space $T^J/\Sigma_{\nu}$. This space 
does not depend, up to a unique isomorphism, on the choice of an unfolding
$\pi$. The points of $T^{\nu}$ are collections $(t_j)$ of $I$-colored 
points of $T$, such that for any $I$, there are $\nu_i$ points of color $i$.
We have the canonical projection $T^J\lra T^{\nu}$, also to be denoted 
by $\pi$. We set $T^{\nu o}=\pi(T^{Jo})$. The map $\pi$ restricted to 
$T^{Jo}$ is an unramified Galois covering with the Galois group 
$\Sigma_{\pi}$. 

\subsection{} For a real $r>0$, let $D(r)$ denote the open disk on the 
complex plane $\{ t\in\BC|\ |t|<r\}$ and $\bar{D}(r)$ its closure.
For $r_1<r_2$, denote by $A(r_1,r_2)$ the open annulus
$D(r_2)-\bar{D}({r_1})$.  

Set $D=D(1)$. Let $\DD$ denote the punctured disk $D-\{0\}$.

\subsection{}   
For an integer $n\geq 1$, consider the space  
$$
E_n=\{ (r_0,\ldots,r_n)\in\BR^{n+1}|\ 0=r_0<r_1<\ldots<r_n=1\}.
$$ 
Obviously, the space $E_n$ is contractible. 

Let $J_1,\ldots,J_n$ be finite sets. Set $J=\coprod_{a=1}^n\ J_a$.
Note that $D^J=D^{J_1}\times\ldots\times D^{J_n}$.  
Let us untroduce the subspace
$D^{J_1,\ldots,J_n}\subset E_n\times D^J$. By definition it consists
of points $((r_a);(t^a_j)\in D^{J_a},\ a=1,\ldots, n)$ such that 
$t^a_j\in A(r_{a-1},r_a)$ for $a=1,\ldots, n$.
The canonical projection $E_n\times D^J\lra D^J$ induces 
the map
$$ 
m(J_1,\ldots,J_n): D^{J_1,\ldots,J_n}\lra D^J.
$$
The image of the above projection lands in the subspace 
$D^{J_1}\times \DD^{J_2}\times\ldots\times \DD^{J_n}\subset D^{J_1}\times
\ldots\times D^{J_n}=D^J$. The induced map 
$$
p(J_1,\ldots,J_n): D^{J_1,\ldots,J_n}\lra D^{J_1}\times \DD^{J_2}\times\ldots
\times \DD^{J_n}
$$
is homotopy equivalence. 

Now assume that we have maps $\pi_a: J_a\lra I$ which  
are unfoldings of the elements $\nu_a$. 
Then their sum $\pi: J\lra I$ is an unfolding of $\nu=\nu_1+\ldots+\nu_n$.
We define the space $D^{\nu_1,\ldots,\nu_n}\subset E_n\times D^{\nu}$ as the 
image of $D^{J_1,\ldots,J_n}$ under the projection
$\Id\times\pi: E_n\times D^J\lra E_n\times D^{\nu}$. The maps 
$m(J_1,\ldots,J_n)$ and $p(J_1,\ldots,J_n)$ induce 
the maps 
$$
m(\nu_1,\ldots,\nu_n):D^{\nu_1,\ldots,\nu_n}\lra D^{\nu_1+\ldots+\nu_n}
$$
and
$$
p(\nu_1,\ldots,\nu_n):D^{\nu_1,\ldots,\nu_n}\lra D^{\nu_1}\times \DD^{\nu_2}
\times\ldots\times \DD^{\nu_n}
$$
respectively, the last map being homotopy equivalence.  

\subsection{} We define the open subspaces  
$$
\DD^{\nu_1,\ldots,\nu_n}=D^{\nu_1,\ldots,\nu_n}\cap (E_n\times \DD^{\nu})
$$
and
$$
\DD^{\nu_1,\ldots,\nu_no}=D^{\nu_1,\ldots,\nu_n}\cap (E_n\times \DD^{\nu o}).
$$
We have the maps
$$
m_a(\nu_1,\ldots,\nu_n): D^{\nu_1,\ldots,\nu_n}\lra D^{\nu_1,\ldots,
\nu_{a-1},\nu_a+\nu_{a+1},\nu_{a+2},\ldots,\nu_n}
$$
and
$$
p_a(\nu_1,\ldots,\nu_n): D^{\nu_1,\ldots,\nu_n}\lra 
D^{\nu_1,\ldots,\nu_a}\times \DD^{\nu_{a+1},\ldots,\nu_n}
$$
$(a=1,\ldots,n-1)$
defined in an obvious manner. They induce the maps
$$
m_a(\nu_1,\ldots,\nu_n): \DD^{\nu_1,\ldots,\nu_n}\lra \DD^{\nu_1,\ldots,
\nu_{a-1},\nu_a+\nu_{a+1},\nu_{a+2},\ldots,\nu_n}
$$
and
$$
p_a(\nu_1,\ldots,\nu_n): \DD^{\nu_1,\ldots,\nu_n}\lra
\DD^{\nu_1,\ldots,\nu_a}\times \DD^{\nu_{a+1},\ldots,\nu_n}
$$
and similar maps between "$o$"-ed spaces. All the maps $p$
are homotopy equivalences.

All these maps satisfy some obvious compatibilities. We will need the following
particular case.

\subsection{} 
\label{nul' rhomb} The {\em rhomb} diagram below commutes.

\begin{center}
  \begin{picture}(14,8)
    \put(7,8){\makebox(0,0){$D^{\nu_1+\nu_2+\nu_3}$}}
    \put(5,6){\makebox(0,0){$D^{\nu_1,\nu_2+\nu_3}$}}
    \put(9,6){\makebox(0,0){$D^{\nu_1+\nu_2,\nu_3}$}}
    \put(3,4){\makebox(0,0){$D^{\nu_1}\times \DD^{\nu_1+\nu_2}$}}
    \put(7,4){\makebox(0,0){$D^{\nu_1,\nu_2,\nu_3}$}}
    \put(11,4){\makebox(0,0){$D^{\nu_1+\nu_2}\times \DD^{\nu_3}$}}
    \put(5,2){\makebox(0,0){$D^{\nu_1}\times \DD^{\nu_1,\nu_2}$}}
    \put(9,2){\makebox(0,0){$D^{\nu_1,\nu_2}\times \DD^{\nu_3}$}}
    \put(7,0){\makebox(0,0){$D^{\nu_1}\times \DD^{\nu_2}\times \DD^{\nu_3}$}}

    \put(5.5,6.5){\vector(1,1){1}}
    \put(8.5,6.5){\vector(-1,1){1}}
    \put(4.5,5.5){\vector(-1,-1){1}}
    \put(6.5,4.5){\vector(-1,1){1}}
    \put(7.5,4.5){\vector(1,1){1}}
    \put(9.5,5.5){\vector(1,-1){1}} 
    \put(4.5,2.5){\vector(-1,1){1}}
    \put(6.5,3.5){\vector(-1,-1){1}}
    \put(7.5,3.5){\vector(1,-1){1}}
    \put(9.5,2.5){\vector(1,1){1}}
    \put(5.5,1.5){\vector(1,-1){1}}
    \put(8.5,1.5){\vector(-1,-1){1}}

    \put(5.5,7){\makebox(0,0){$m$}}
    \put(8.5,7){\makebox(0,0){$m$}}
    \put(3.5,5){\makebox(0,0){$p$}}
    \put(5.5,5){\makebox(0,0){$m$}}
    \put(8.5,5){\makebox(0,0){$m$}}
    \put(10.5,5){\makebox(0,0){$p$}}
    \put(3.5,3){\makebox(0,0){$m$}}
    \put(5.5,3){\makebox(0,0){$p$}}
    \put(8.5,3){\makebox(0,0){$p$}}
    \put(10.5,3){\makebox(0,0){$m$}}
    \put(5.5,1){\makebox(0,0){$p$}}
    \put(8.5,1){\makebox(0,0){$p$}}
    
  \end{picture} 
\end{center}

\subsection{} We will denote by $\CA^o$ and call an {\em open
$I$-coloured configuration space} the collection of all spaces 
$\{ \DD^{\nu_1,\ldots,\nu_n o}\}$ together with the maps
$\{ m_a(\nu_1,\ldots,\nu_n), 
p_a(\nu_1,\ldots,\nu_n)\}$ between their various products. 

We will call a {\em local system} over $\CA^o$, or a {\em braiding local
system} a collection of data (a), (b) below satisfying the property 
(c) below. 

(a) A local system $\CI_{\mu}^{\nu}$ over $\DD^{\nu o}$ given for any
$\nu\in  Y^+, \mu\in X$. 

(b) An isomorphism
$\phi_{\mu}(\nu_1,\nu_2): m^*\CI_{\mu}^{\nu_1+\nu_2}\iso p^*(\CI_{\mu}^{\nu_1}
\boxtimes\CI_{\mu-\nu_1}^{\nu_2})$
given for any $\nu_1,\nu_2\in  Y^+, \mu\in X$. 

Here $p=p(\nu_1,\nu_2)$ and 
$m=m(\nu_1,\nu_2)$ are the arrows in the diagram 
$\DD^{\nu_1 o}\times \DD^{\nu_2 o}\overset{p}{\lla} \DD^{\nu_1,\nu_2 o}
\overset{m}{\lra} \DD^{\nu_1+\nu_2 o}$.

The isomorphisms $\phi_{\mu}(\nu_1,\nu_2)$ are called the {\em factorization 
isomorphisms}. 

(c) (The {\em associativity} of factorization isomorphisms.) For any  
$\nu_1, \nu_2, \nu_3\in  Y^+, \mu\in X$, the octagon below commutes.
Here the maps 
$m, p$ are the maps in the rhombic diagram above, with $D$, $\DD$ replaced
by $\DD^o$.

\begin{center}
  \begin{picture}(14,12)
    \put(4,12){\makebox(0,0){$m^*m_1^*\CI^{\nu_1+\nu_2+\nu_3}_{\mu}$}}
    \put(10,12){\makebox(0,0){$m^*m_2^*\CI^{\nu_1+\nu_2+\nu_3}_{\mu}$}}
    \put(0,8){\makebox(0,0){$m^*p^*(\CI^{\nu_1}_{\mu}\boxtimes\CI^{\nu_2+\nu_3}
    _{\mu-\nu_1})$}}
    \put(14,8){\makebox(0,0){$m^*p^*(\CI^{\nu_1+\nu_2}_{\mu}\boxtimes\CI^{\nu_3}
    _{\mu-\nu_1-\nu_2})$}}
    \put(0,4){\makebox(0,0){$p^*m^*(\CI^{\nu_1}_{\mu}\boxtimes\CI^{\nu_2+\nu_3}
    _{\mu-\nu_1})$}}
    \put(14,4){\makebox(0,0){$p^*m^*(\CI^{\nu_1+\nu_2}_{\mu}\boxtimes\CI^{\nu_3}
    _{\mu-\nu_1-\nu_2})$}}
    \put(3,0){\makebox(0,0){$p_1^*p^*(\CI_{\mu}^{\nu_1}\boxtimes\CI_{\mu-\nu_1}
    ^{\nu_2}\boxtimes\CI_{\mu-\nu_1-\nu_2}^{\nu_3})$}}
    \put(11,0){\makebox(0,0){$p_2^*p^*(\CI_{\mu}^{\nu_1}\boxtimes\CI_{\mu-\nu_1}
    ^{\nu_2}\boxtimes\CI_{\mu-\nu_1-\nu_2}^{\nu_3})$}}

    \put(6,12){\line(1,0){2}}
    \put(6,11.9){\line(1,0){2}}
    \put(3.5,11.5){\vector(-1,-1){3}}
    \put(10,11.5){\vector(1,-1){3}}
    \put(0,7.5){\line(0,-1){3}}
    \put(.1,7.5){\line(0,-1){3}}
    \put(14,7.5){\line(0,-1){3}}
    \put(13.9,7.5){\line(0,-1){3}}
    \put(.5,3.5){\vector(1,-1){3}}
    \put(13.5,3.5){\vector(-1,-1){3}}
    \put(6.4,0){\line(1,0){1}}
    \put(6.4,.1){\line(1,0){1}}

    \put(0,10){\makebox(0,0){$\phi_{\mu}(\nu_1,\nu_2+\nu_3)$}}
    \put(14,10){\makebox(0,0){$\phi_{\mu}(\nu_1+\nu_2,\nu_3)$}}
    \put(0,2){\makebox(0,0){$\phi_{\mu-\nu_1}(\nu_2,\nu_3)$}}
    \put(14,2){\makebox(0,0){$\phi_{\mu}(\nu_1,\nu_2)$}}

  \end{picture} 
\end{center}

Written more concisely, the axiom (c) reads as a "cocycle" condition 

(c)$'$  $\phi_{\mu}(\nu_1,\nu_2)\circ\phi_{\mu}(\nu_1+\nu_2,\nu_3)=
\phi_{\mu-\nu_1}(\nu_2,\nu_3)\circ\phi_{\mu}(\nu_1,\nu_2+\nu_3)$. 

\subsection{Correctional lemma} {\em Assume we are given the data
{\em (a), (b)} as above, with  
one-dimensional local systems $\CI_{\mu}^{\nu}$. 
Then there exist a collection of constants $c_{\mu}(\nu_1,\nu_2)\in \sk^*\
(\mu\in X,\nu_1,\nu_2\in  Y^+)$ such that the corrected isomorphisms 
$\phi'_{\mu}(\nu_1,\nu_2)=c_{\mu}(\nu_1,\nu_2)\phi_{\mu}(\nu_1,\nu_2)$ 
satisfy the associativity axiom {\em (c)}}. $\Box$

\subsection{} The notion of a morphism between two braiding local
systems is defined in an obvious way. This defines the category 
of braiding local systems, $\Bls$.

Suppose that $\CI=\{\CI_{\mu}^{\nu}; \phi_{\mu}(\nu_1,\nu_2)\}$ and 
$\CJ=\{\CJ_{\mu}^{\nu}; \psi_{\mu}(\nu_1,\nu_2)\}$ are two braiding local 
systems. Their 
{\em tensor product} $\CI\otimes\CJ$ is defined by 
$(\CI\otimes\CJ)_{\mu}^{\nu}=\CI_{\mu}^{\nu}\otimes\CJ_{\mu}^{\nu}$, 
the factorization isomorphisms being $\phi_{\mu}(\nu_1,\nu_2)\otimes 
\psi_{\mu}(\nu_1,\nu_2)$. This makes $\Bls$ a tensor category.
The subcategory of one-dimensional braiding local systems is a Picard  
category. 

\subsection{Example} 
\label{nul' sign} {\em Sign local system.} Let $L$ be a one-dimensional
vector space. Let $\nu\in Y^+$. Let $\pi: J\lra I$ be an unfolding of $\nu$, 
whence the canonical projection $\pi: \DD^{J\circ}\lra \DD^{\nu\circ}$.
Pick a point $\tilde{x}=(x_j)\in \DD^{J\circ}$ with all $x_j$ being {\em
real}; let $x=\pi(\tilde{x})$. The choice of base points defines the 
homomorphism $\pi_1(\DD^{\nu\circ};x)\lra\Sigma_{\pi}$. Consider its
composition with the sign map $\Sigma_{\pi}\lra\mu_2\hra \sk^*$. Here
$\mu_2=\{\pm 1\}$ is the group of square roots of $1$ in $\sk^*$.
We get a map $s: \pi_1(\DD^{\nu\circ};x)\lra \sk^*$. It defines a local system
$\Sign^{\nu}(L)$ over $\DD^{\nu\circ}$ whose stalk at $x$ is equal to $L$.
This local system does not depend (up to the 
unique isomorphism) on the choices made. 

Given a diagram
$$
\DD^{\nu_1\circ}\times \DD^{\nu_2\circ}\overset{p}{\lla} \DD^{\nu_1,\nu_2\circ}
\overset{m}{\lra}\DD^{\nu_2+\nu_2\circ}
$$
we can choose base point $x_i\in \DD^{\nu_i}$ and $x\in \DD^{\nu_1+\nu_2\circ}$
such that there exists $y\in \DD^{\nu_1,\nu_2}$ with $p(y)=(x_1,x_2)$ and
$m(y)=x$. Let 
$$
\phi^{L_1,L_2}(\nu_1,\nu_2): m^*\Sign^{\nu_1+\nu_2}(L_1\otimes L_2)\iso 
p^*(\Sign^{\nu_1}(L_1)\boxtimes\Sign^{\nu_2}(L_2))
$$
be the unique isomorphism equal to $\Id_{L_1\otimes L_2}$ over $y$. 
These isomorphisms satisfy the associativity condition (cf. (c)$'$ above). 

Let us define the brading local system $\Sign$ by $\Sign_{\mu}^{\nu}=
\Sign^{\nu}(\One)$, $\phi_{\mu}(\nu_1,\nu_2)=\phi^{\One,\One}(\nu_1,\nu_2)$.
Here $\One$ denotes the standard vector space $\sk$.

\subsection{Example} 
\label{nul' stand loc} {\em Standard local systems $\CJ$, $\CI$.} Let
$\nu\in Y^+$, let $\pi: J\lra I$ be an unfolding of $\nu$, $n=\card(J)$, 
$\pi: \DD^{J\circ}\lra \DD^{\nu\circ}$ the
corresponding projection. Let $L$ be a one-dimensional vector space. 
For each isomorphism $\sigma: J\iso \{1,\ldots, n\}$, choose a point $x(\sigma)= 
(x_j)\in \DD^{J\circ}$ such that all $x_j$ are real and positive and for each
$j', j''$ such that $\sigma(j')<\sigma(j'')$, we have $x_{j'}<x_{j''}$.
  
Let us define a local system $\CJ_{\mu}^{\pi}(L)$ over $\DD^{J\circ}$ as follows.
We set $\CJ_{\mu}^{\pi}(L)_{x(\sigma)}=L$ for all $\sigma$. We define the 
monodromies 
$$
T_{\gamma}: \CJ_{\mu}^{\pi}(L)_{x(\sigma)}\lra \CJ_{\mu}^{\pi}(L)_{x(\sigma')}
$$
along the homotopy classes of paths $\gamma: x_{\sigma}\lra x_{\sigma'}$ 
which generate the fundamental groupoid. Namely, let $x_{j'}$ and $x_{j''}$ 
be some neighbour points in $x(\sigma)$, with $x_{j'}<x_{j''}$.
Let $\gamma(j',j'')^+$ (resp. $\gamma(j',j'')^-$) be the paths corresponding 
to the movement of $x_{j''}$ in the upper (resp. lower)
hyperplane to the left from $x_{j'}$ position. We set 
$$
T_{\gamma(j',j'')^{\pm}}=\zeta^{\mp\pi(j')\cdot\pi(j'')}.
$$
Let $x_j$ be a point in $x(\sigma)$ closest to $0$. Let $\gamma(j)$ be the path 
corresponding to the counterclockwise travel of $x_j$ around $0$. We set 
$$
T_{\gamma(j)}=\zeta^{2\mu\cdot\pi(j)}.
$$
The point is that the above fromulas give a well defined morphism from
the fundamental groupoid to the groupoid of one-dimensional vector spaces. 

The local system $\CJ_{\mu}^{\pi}(L)$ admits an obvious 
$\Sigma_{\pi}$-equivariant structure. This defines the local system 
$\CJ_{\mu}^{\nu}(L)$ over $\DD^{\nu\circ}$. Given a diagram
$$
\DD^{\nu_1\circ}\times \DD^{\nu_2\circ}\overset{p}{\lla} \DD^{\nu_1+\nu_2\circ}
\overset{m}{\lra} \DD^{\nu_1+\nu_2\circ},
$$
let
$$
\phi_{\mu}^{L_1,L_2}(\nu_1,\nu_2): m^*\CJ^{\nu_1+\nu_2}_{\mu}(L_1\otimes L_2)
\iso p^*(\CJ^{\nu_1}_{\mu}(L_1)\boxtimes\CJ^{\nu_2}_{\mu-\nu_1}(L_2))
$$
be the unique isomorphism equal to $Id_{L_1\otimes L_2}$ over compatible 
base points. Here the compatibility is understood in the same sense as in the 
previous subsection. These isomorphisms satisfy the associativity 
condition.  

We define the braiding local system $\CJ$ by $\CJ_{\mu}^{\nu}=\CJ_{\mu}^{\nu}
(\One)$, $\phi_{\mu}(\nu_1,\nu_2)=\phi_{\mu}^{\One,\One}(\nu_1,\nu_2)$. 

We define the braiding local system $\CI$ by $\CI=\CJ\otimes\Sign$. The local 
systems $\CJ$, $\CI$ are called the standard local systems. In the sequel 
we will mostly need the local system $\CI$.

\section{Factorizable sheaves}
\label{nul' factoriz}

\subsection{} In the sequel for each $\nu\in Y^+$, we will denote by 
$\CS$ the stratification on the space $D^{\nu}$ the closures of whose 
strata are various intersections of hypersurfaces given by the equations 
$t_j=0,\ t_{j'}=t_{j''}$. The same letter will denote the induced
stratifications on its subspaces.  

\subsection{} Set  
$\CI_{\mu}^{\nu\bullet}=j_{!*}\CI_{\mu}^{\nu}[\dim \DD^{\nu\circ}]
\in\CM(\DD^{\nu};\CS)$
where $j: \DD^{\nu\circ}\hra \DD^{\nu}$ is the embedding.
   
Given a diagram 
$$
\DD^{\nu_1}\times \DD^{\nu_2}\overset{p}{\lla} \DD^{\nu_1,\nu_2}\overset{m}
{\lra} \DD^{\nu_1+\nu_2},
$$
the structure isomorphisms $\phi_{\mu}(\nu_1,\nu_2)$ of the local 
system $\CI$ induce the isomorphisms, to be denoted by the same 
letters, 
$$
\phi_{\mu}(\nu_1,\nu_2): m^*\CI^{\nu_1+\nu_2\bullet}\iso 
p^*(\CI^{\nu_1\bullet}\boxtimes\CI^{\nu_2\bullet}).
$$
Obviously, these isomorphisms satisfy 
the associativity axiom. 

\subsection{} Let us fix a coset $c\in X/Y$. We will regard $c$ as a 
subset of $X$.   
We will call a {\em factorizable sheaf} $\CM$ {\em supported at $c$} 
a collection 
of data (w), (a), (b) below satisfying the axiom (c) below. 

(w) An element $\lambda=\lambda(\CM)\in c$.

(a) A perverse sheaf $\CM^{\nu}\in\CM(D^{\nu};\CS)$ given for each $\nu\in Y^+$.

(b) An isomorphism $\psi(\nu_1,\nu_2): m^*\CM^{\nu_1+\nu_2}\iso
p^*(\CM^{\nu_1}\boxtimes\CI_{\lambda-\nu_1}^{\nu_2\bullet})$
given for any $\nu_1, \nu_2\in Y^+$. 

Here $p, m$ denote the arrows in the diagram $D^{\nu_1}\times \DD^{\nu_2}
\overset{p}{\lla} D^{\nu_1,\nu_2}\overset{m}{\lra} D^{\nu_1+\nu_2}$. 

The isomorphisms $\psi(\nu_1,\nu_2)$ are called the {\em factorization 
isomorphisms}. 

(c) For any $\nu_1, \nu_2, \nu_3\in Y^+$, the following  
{\em associativity condition} is fulfilled:
$$
\psi(\nu_1,\nu_2)\circ\psi(\nu_1+\nu_2,\nu_3)=
\phi_{\lambda-\nu_1}(\nu_2,\nu_3)\circ\psi(\nu_1,\nu_2+\nu_3).
$$

We leave to the reader to draw the whole octagon expressing this axiom. 

\subsection{} Let $\CM=\{\CM^{\nu};\ \psi(\nu_1,\nu_2)\}$ be a factorizable
sheaf supported at a coset $c\in X/Y$, $\lambda=\lambda(\CM)$. For each 
$\lambda'\geq\lambda, \nu\in Y^+$, define a sheaf $\CM_{\lambda'}^{\nu}\in
\CM(D^{\nu};\CS)$ by  
$$
\CM_{\lambda'}^{\nu}=\left\{ \begin{array}{ll}\iota(\lambda'-\lambda)_*
\CM^{\nu-\lambda'+\lambda}&\mbox{ if }\nu-\lambda'+\lambda\in Y^+\\
0&\mbox{ otherwise.}\end{array}\right.
$$
Here 
$$
\iota(\nu'): D^{\nu}\lra D^{\nu+\nu'}
$$
denotes the closed embedding adding $\nu'$ points sitting at the origin. 
The factorization isomorphisms $\psi(\nu_1,\nu_2)$ induce
similar isomorphisms 
$$
\psi_{\lambda'}(\nu_1,\nu_2): m^*\CM_{\lambda'}(\nu_1+\nu_2)\iso 
p^*(\CM_{\lambda'}^{\nu_1}\boxtimes\CI_{\lambda'-\nu_1}^{\nu_2})\ (\lambda'
\geq\lambda)
$$

\subsection{} Let $\CM, \CN$ be two factorizable sheaves supported 
at $c$. Let $\lambda\in X$ 
be such that $\lambda\geq\lambda(\CM)$ and $\lambda\geq\lambda(\CN)$.
For $\nu\geq\nu'$ in $Y^+$, consider the following composition
$$
\tau_{\lambda}(\nu,\nu'):\Hom(\CM_{\lambda}^{\nu},\CN_{\lambda}^{\nu})
\overset{m_*}{\lra}\Hom(m^*\CM_{\lambda}^{\nu},m^*\CN_{\lambda}^{\nu})
\overset{\psi(\nu',\nu-\nu')_*}{\iso}
$$
$$
\iso\Hom(p^*(\CM_{\lambda}^{\nu'}
\boxtimes\CI_{\lambda-\nu'}^{\nu-\nu'\bullet}),p^*(\CN_{\lambda}^{\nu'}
\boxtimes\CI_{\lambda-\nu'}^{\nu-\nu'\bullet}))=\Hom(\CM_{\lambda}^{\nu'},
\CN_{\lambda}^{\nu'}).
$$
Let us define the space of homomorphisms $\Hom(\CM,\CN)$ by 
$$
\Hom(\CM,\CN)=\dirlim_{\lambda}\invlim_{\nu}\Hom(\CM_{\lambda}^{\nu},
\CN_{\lambda}^{\nu})
$$
Here the inverse limit is taken over $\nu\in Y^+$, the transition maps being 
$\tau_{\lambda}(\nu,\nu')$ and the direct limit is taken over $\lambda\in X$ 
such that $\lambda\geq\lambda(\CM), \lambda\geq\lambda(\CN)$,  
the transition maps being induced by the obvious isomorphisms
$$
\Hom(\CM_{\lambda}^{\nu},\CN_{\lambda}^{\nu})=
\Hom(\CM_{\lambda+\nu'}^{\nu+\nu'},\CN_{\lambda+\nu'}^{\nu+\nu'})\
(\nu'\in Y^+). 
$$
With these spaces of homomorphisms and the obvious compositions, 
the factorizable sheaves supported at $c$ form the category, to be 
denoted by $\tFS_c$. By definition, the category $\tFS$ is the direct
product $\prod_{c\in X/Y}\ \tFS_c$. 

\subsection{Finite sheaves} Let us call a factorizable sheaf 
$\CM=\{\CM^{\nu}\}\in\tFS_c$ {\em finite} if there exists only a finite 
number of $\nu\in Y^+$ such that the conormal bundle of the origin
$O\in\CA^{\nu}$ is contained in the singular support of $\CM^{\nu}$.
Let $\FS_c\subset\tFS_c$ be the full subcategory of finite 
factorizable sheaves. We define the category $\FS$ by 
$\FS=\prod_{c\in X/Y}\FS_c$. One proves (using the lemma below) that 
$\FS$ is an abelian category.  

\subsection{Stabilization Lemma} {\em Let $\CM, \CN\in\FS_c,\ 
\mu\in X_c, \mu\geq\lambda(\CM), \mu\geq\lambda(\CN)$. There exists 
$\nu_0\in Y^+$ such that for all $\nu\geq\nu_0$ the canonical maps 
$$
\Hom(\CM,\CN)\lra\Hom(\CM_{\mu}^{\nu},\CN_{\mu}^{\nu})
$$
are isomorphisms.} $\Box$

\subsection{Standard sheaves} 
\label{nul' stand} Given $\mu\in X$, let us define the
"standard sheaves" $\CM(\mu), \CCD\CM(\mu)$ and $\CL(\mu)$
supported at the coset $\mu+Y$, by 
$\lambda(\CM(\mu))=\lambda(\CCD\CM(\mu))=\lambda(\CL(\mu))=\mu$;
$$
\CM(\mu)^{\nu}=j_!\CI^{\nu\bullet}_{\mu};\
\CCD\CM(\mu)^{\nu}=j_*\CI^{\nu\bullet}_{\mu};\
\CL(\mu)^{\nu}=j_{!*}\CI^{\nu\bullet}_{\mu},
$$
$j$ being the embedding $\DD^{\nu}\hra D^{\nu}$. The factorization maps
are defined by functoriality from the similar maps for $\CI^{\bullet}$. 

One proves that all these sheaves are finite.

\section{Tensor product}

In this section we will give (a sketch of) the construction of the
tensor structure on the category $\tFS$. We will make the assumption of 
~\ref{nul' k big}\footnote{Note that these assumptions are not necessary for the
construction of the tensor structure. They are essential, however, 
for the construction of braiding.}. 

\subsection{} For $z\in\BC$ and a real positive $r$, let $D(z;r)$
denote the open disk $\{ t\in\BC|\ |t-z|<r\}$ and $\bar{D}(z;r)$ its
closure. 

\subsection{}
\label{nul' spaces} For $\nu\in Y^+$, let us define the space $D^{\nu}(2)$
as the product $\DD\times D^{\nu}$. Its points will be denoted
$(z;(t_j))$ where $z\in \DD, (t_j)\in D^{\nu}$. Let us define the open
subspaces 
$$
\DD^{\nu}(2)=\{(z;(t_j))\in D^{\nu}(2)|\ t_j\neq 0, z\mbox{ for all }j\};\
\DD^{\nu}(2)^{\circ}=\DD^{\nu}(2)\cap(\DD\times \DD^{\nu\circ}).
$$
For $\nu,\nu'\in Y^+$, let us define the space 
$D^{\nu,\nu'}(2)$ as the subspace of $\BR_{>0}\times D^{\nu+\nu'}(2)$
consisting of all elements $(r;z;(t_j))$ such that $|z|<r<1$ and $\nu$ 
of the points $t_j$ live inside the disk $D(r)$ and $\nu'$ of them inside
the annulus $A(r,1)$.

We have a diagram

(a) $D^{\nu}(2)\times \DD^{\nu'}\overset{p}{\lla} D^{\nu,\nu'}(2)
\overset{m}{\lra} D^{\nu+\nu'}(2)$. 

Here $p((r;z;(t_j)))=((z;(t_{j'})), (t_{j''}))$ where $t_{j'}$ (resp. $t_{j''}$) 
being the points from the collection $(t_j)$ lying in $D(r)$ 
(resp. in $A(r,1)$);
$m((r;z;(t_j)))=(z;(t_j))$. The map $p$ is a homotopy equivalence. 

For $\nu_1,\nu_2,\nu\in Y^+$, let $D^{\nu_1;\nu_2;\nu}(2)$ be the 
subspace of $\BR_{>0}\times\BR_{>0}\times D^{\nu_1+\nu_2+\nu}(2)$
consisting of all elements $(r_1;r_2;z;(t_j))$ such that 
$\bar{D}(r_1)\cup\bar{D}(z;r_2)\subset D$; $\bar{D}(r_1)\cap\bar{D}(z;r_2)
=\emptyset$; $\nu_1$ of the points $(t_j)$ lie inside $D(r_1)$, 
$\nu_2$ of them lie inside $D(z;r_2)$ and $\nu$ of them lie 
inside $D-(\bar{D}(r_1)\cup\bar{D}(z;r_2))$. 

We have a diagram 

(b) $D^{\nu_1}\times D^{\nu_2}\times \DD^{\nu}(2)\overset{p}{\lla}
D^{\nu_1;\nu_2;\nu}(2)\overset{m}{\lra} D^{\nu_1+\nu_2+\nu}(2)$. 

Here $p((r_1;r_2;z;(t_j)))=((t_{j'}); (t_{j''}-z); (t_{j'''}))$ where 
$t_{j'}$ (resp. $t_{j''}, t_{j'''}$) are the points lying inside $D(r_1)$ 
(resp. $D(z;r_2), D-(\bar{D}(r_1)\cup\bar{D}(z;r_2)))$; 
$m((r_1;r_2;z;(t_j)))=(z;(t_j))$. The map $p$ is a homotopy equivalence. 

\subsection{}    
We set $\DD^{\nu,\nu'}(2)=D^{\nu,\nu'}(2)\cap(\BR_{>0}
\times \DD^{\nu+\nu'}(2))$;
$\DD^{\nu,\nu^{\prime}}(2)^{\circ}=D^{\nu,\nu'}(2)\cap(\BR_{>0}\times
\DD^{\nu+\nu'}(2)^{\circ})$; $\DD^{\nu_1;\nu_2;\nu}(2)=D^{\nu_1;\nu_2;\nu}(2)
\cap(\BR_{>0}\times\BR_{>0}\times \DD^{\mu_1+\nu_2+\nu}(2))$;
$\DD^{\nu_1;\nu_2;\nu}(2)^{\circ}=D^{\nu_1;\nu_2;\nu}(2)\cap(\BR_{>0}
\times\BR_{>0}\times \DD^{\nu_1+\nu_2+\nu}(2)^{\circ})$.

\subsection{} Given $\mu_1,\mu_2\in X, \nu\in Y^+$, choose an unfolding 
$\pi: J\lra I$ of the element $\nu$.  
In the same manner as in \ref{nul' stand loc}, we define the
one-dimensional local system $\CJ_{\mu_1,\mu_2}^{\nu}$ over 
$\DD^{\nu}(2)^{\circ}$ with the following monodromies: the monodromy
around a loop corresponding to the counterclockwise travel of the point 
$z$ around $0$ (resp. $t_j$ around $0$, $t_j$ around $z$, 
$t_{j'}$ around $t_{j''}$) is equal to the multiplication 
by $\zeta^{-2\mu_1\cdot\mu_2}$ (resp. $\zeta^{2\mu_1\cdot\pi(j)}$, 
$\zeta^{2\mu_2\cdot\pi(j)}$, $\zeta^{-2\pi(j')\cdot\pi(j'')}$).  

As in {\em loc. cit.}, one defines isomorphisms

(a) $\phi_{\mu_1,\mu_2}(\nu,\nu'): m^*\CJ_{\mu_1,\mu_2}^{\nu+\nu'}\iso 
p^*(\CJ_{\mu_1,\mu_2}^{\nu}\boxtimes\CJ_{\mu_1+\mu_2-\nu}^{\nu'})$

where $p$, $m$ are the morphisms in the diagram \ref{nul' spaces} (a)
(restricted to the $\DD^{\circ}$-spaces) and

(b) $\phi_{\mu_1,\mu_2}(\nu_1;\nu_2;\nu): m^*\CJ_{\mu_1,\mu_2}^{\nu_1+
\nu_2+\nu}\iso p^*(\CJ_{\mu_1}^{\nu_1}\boxtimes\CJ_{\mu_2}^{\nu_2}\boxtimes
\CJ_{\mu_1-\nu_1,\mu_2-\nu_2}^{\nu})$ 

where $p, m$ are the morphisms in the diagram \ref{nul' spaces} (b) (restricted
to the $\DD^{\circ}$-spaces),
which satisfy the cocycle conditions 

(c) $\phi_{\mu_1,\mu_2}(\nu,\nu')\circ\phi_{\mu_1,\mu_2}(\nu+\nu',\nu'')=
\phi_{\mu_1+\mu_2-\nu}(\nu',\nu'')\circ\phi_{\mu_1,\mu_2}(\nu,\nu'+\nu'')$ 

(d) $(\phi_{\mu_1}(\nu_1,\nu'_1)\boxtimes\phi_{\mu_2}(\nu_2,\nu'_2))
\circ\phi_{\mu_1,\mu_2}(\nu_1+\nu'_1;\nu_2+\nu'_2;\nu)=$

$=\phi_{\mu_1-\nu_1,\mu_2-\nu_2}(\nu'_1;\nu'_2;\nu)\circ\phi_{\mu_1,\mu_2}
(\nu_1;\nu_2;\nu+\nu'_1+\nu'_2)$  

(we leave to the reader the definition of the corresponding spaces).

\subsection{} Let us consider the sign local systems introduced in \ref{nul' sign}.
We will keep the same notation $\Sign^{\nu}$ for the inverse image
of the local system $\Sign^{\nu}$ under the forgetting of 
$z$ map $\DD^{\nu}(2)^{\circ}\lra \DD^{\nu\circ}$. We have the factorization
isomorphisms

(a) $\phi^{\Sign}(\nu,\nu'): m^*\Sign^{\nu+\nu'}\iso p^*(\Sign^{\nu}\boxtimes
\Sign^{\nu'})$; 

(b) $\phi^{\Sign}(\nu_1;\nu_2;\nu): m^*\Sign^{\nu_1+\nu_2+\nu}\iso p^*(
\Sign^{\nu_1}\boxtimes\Sign^{\nu_2}\boxtimes\Sign^{\nu})$

which satisfy the cocycle conditions similar to (c), (d) above.    

\subsection{}  
We define the local systems $\CI^{\nu}_{\mu_1,\mu_2}$ over the spaces 
$\DD^{\nu}(2)^{\circ}$ by

$\CI_{\mu_1,\mu_2}^{\nu}=\CJ_{\mu_1,\mu_2}^{\nu}\otimes\Sign^{\nu}$.  
  
The collection of local systems $\{\CI_{\mu_1,\mu_2}^{\nu}\}$ together
with the maps  
$\phi^{\CI}_{\mu_1,\mu_2}(\nu,\nu')=\phi^{\CJ}_{\mu_1,\mu_2}(\nu,\nu')\otimes
\phi^{\Sign}(\nu,\nu')$ and    
$\phi^{\CI}_{\mu_1,\mu_2}(\nu_1;\nu_2;\nu)=\phi^{\CJ}_{\mu_1,\mu_2}(\nu_1;\nu_2;
\nu)\otimes\phi^{\Sign}(\nu_1;\nu_2;\nu)$, forms an object $\CI(2)$ which 
we call 
a {\em standard braiding local system over the configuration space $\CA(2)^{\circ}
=\{ \DD^{\nu}(2)^{\circ}\}$}. It is unique up to a (non unique) isomorphism.
We fix such a local system.  

\subsection{} We set $\CI_{\mu_1,\mu_2}^{\nu\bullet}=
j_{!*}\CI_{\mu_1,\mu_2}^{\nu}
[\dim \DD^{\nu}(2)^{\circ}]$ where $j: \DD^{\nu}(2)^{\circ}\hra \DD^{\nu}(2)$
is the open embedding. It is an object of the category $\CM(\DD^{\nu}(2);\CS)$
where $\CS$ is the evident stratification. The factorization isomorphisms 
for the local system $\CI$ induce the analogous isomorphisms 
between these sheaves, to be denoted by the same letter. The collection 
of these sheaves and factorization isomorphisms will be denoted 
$\CI(2)^{\bullet}$.     

\subsection{} 
Suppose we are given two factorizable sheaves $\CM, \CN$. Let us call their
{\em gluing}, and denote by $\CM\boxtimes\CN$, the collection of 
perverse sheaves $(\CM\boxtimes\CN)^{\nu}$ over the spaces
$D^{\nu}(2)$ $(\nu\in Y^+)$ together with isomorphisms
$$
\psi(\nu_1;\nu_2;\nu): m^*(\CM\boxtimes\CN)^{\nu}\iso p^*(\CM^{\nu_1}
\boxtimes\CN^{\nu_2}\boxtimes\CI^{\nu\bullet}_{\lambda(\CM)-\nu_1,
\lambda(\CN)-\nu_2}),
$$
$p, m$ being the maps in the diagram \ref{nul' spaces} (b), which satisfy the
cocycle condition

$(\psi^{\CM}(\nu_1,\nu'_1)\boxtimes\psi^{\CN}(\nu_2,\nu'_2))\circ
\psi(\nu_1+\nu'_1;\nu_2+\nu'_2;\nu)=$

$=\phi_{\lambda(\CM)-\nu_1,\lambda(\CN)-\nu_2}(\nu'_1;\nu'_2;\nu)\circ
\psi(\nu_1;\nu_2;\nu+\nu_1'+\nu_2')$ 

for all $\nu_1,\nu_1',\nu_2,\nu_2',\nu\in Y^+$. 

Such a gluing exists and is unique, up to a unique isomorphism. The 
factorization isomorphisms $\phi_{\mu_1,\mu_2}(\nu_1;\nu_2;\nu)$ 
for $\CI(2)^{\bullet}$ and the ones for $\CM, \CN$,
induce the isomorphisms 
$$
\psi^{\CM\boxtimes\CN}(\nu,\nu'): m^*(\CM\boxtimes\CN)^{\nu+\nu'}\iso 
p^*((\CM\boxtimes\CN)^{\nu}\boxtimes\CI^{\nu'\bullet}_{\lambda(\CM)+
\lambda(\CN)-\nu}).
$$
satisfying the obvious cocycle condition. 

\subsection{} 
\label{nul' def tens} Now we can define the tensor product $\CM\otimes\CN\in\tFS$.
Namely, set $\lambda(\CM\otimes\CN)=\lambda(\CM)+\lambda(\CN)$. For 
each $\nu\in Y^+$, set
$$
(\CM\otimes\CN)^{\nu}=\Psi_{z\ra 0}((\CM\boxtimes\CN)^{\nu}).
$$
Here $\Psi_{z\ra 0}: \CM(D^{\nu}(2))\lra\CM(D^{\nu})$ denotes the functor 
of nearby cycles for the function $D^{\nu}(2)\lra D$ sending $(z;(t_j))$ to 
$z$. Note that 
$$
\Psi_{z\ra 0}(\CI_{\mu_1,\mu_2}^{\nu})=\CI_{\mu_1+\mu_2}^{\nu}.
$$
The factorization isomorphisms $\psi^{\CM\boxtimes\CN}$ induce the factorization
isomorphisms between the sheaves $(\CM\otimes\CN)^{\nu}$. 
This defines a factorizable sheaf $\CM\otimes\CN$.

One sees at once that this construction is functorial; thus it 
defines a functor of tensor product $\otimes: \tFS\times\tFS\lra\tFS$. 

The subcategory $\FS\subset\tFS$ is stable under the tensor product. 
The functor $\otimes:\ \FS\times\FS\lra\FS$ extends uniquely to a 
functor $\otimes:\ \FS\otimes\FS\lra\FS$ (for the discussion of the 
tensor product of abelian categories, see \cite{d2} 5).    

\subsection{} The half-circle travel of the point $z$ around $0$ from $1$ to
$-1$ in the upper halfplane defines the braiding isomorphisms 
$$
{R}_{\CM,\CN}: \CM\otimes\CN\iso\CN\otimes\CM.
$$
We will not describe here the precise definition of the associativity 
isomorphisms for the tensor product $\otimes$. We just mention that 
to define them one should introduce into the game certain configuration  
spaces $D^{\nu}(3)$ whose (more or less obvious) definition we leave 
to the reader. 

The unit of this tensor structure is the sheaf $\One=\CL(0)$ 
(cf. \ref{nul' stand}).
  
Equipped with these complementary structures, the category $\tFS$ becomes 
a {\em braided tensor category}.

\section{Vanishing cycles}

{\em GENERAL GEOMETRY}

\subsection{} 
\label{nul' facets} Let us fix a finite set $J$, and consider the space
$D^J$. Inside this space, let us consider the subspaces $D^J_{\BR}=
D^J\cap\BR^J$ and $D^{J+}=D^J\cap\BR^J_{\geq 0}$.

Let $\CH$ be the set ({\em arrangement}) 
of all real hyperplanes in $D^J_{\BR}$ of the form
$H_j:\ t_j=0$ or $H_{j'j''}:\ t_{j'}=t_{j''}$. An {\em edge} $L$ of the
arrangement $\CH$ is 
a subspace of $D^J_{\BR}$ which is a non-empty intersection $\bigcap H$
of some hyperplanes from $\CH$. We denote by $L^{\circ}$ the complement 
$L-\bigcup L'$, the union over all edges $L'\subset L$ of smaller 
dimension.   
A {\em facet} of $\CH$ is a connected
component $F$ of some $L^{\circ}$. We call a facet {\em positive} if
it lies entirely inside $D^{J+}$.

For example, we have a unique smallest facet $O$ --- the origin. For each 
$j\in J$, we have a positive one-dimensional facet $F_j$ given 
by the equations $t_{j'}=0\ (j'\neq j);\ t_j\geq 0$.

Let us choose a point $w_F$ on each positive facet $F$. We call a 
{\em flag} a sequence of embedded positive facets $\bF:\ F_0\subset\ldots F_p$; 
we say that $\bF$ {\em starts} from $F_0$. To such a flag we assign 
the simplex $\Delta_{\bF}$ --- the convex hull of the points
$w_{F_0},\ldots, w_{F_p}$. 

To each positive facet $F$ we assign the following two spaces: 
$D_F=\bigcup\Delta_{\bF}$, the union over all flags $\bF$ starting from $F$, 
and $S_F=\bigcup\Delta_{\bF'}$, the union over all flags $\bF'$ starting 
from a facet which properly contains $F$. Obviously, $S_F\subset D_F$. 

\subsection{} Given a complex $\CK\in\CCD(D^J;\CS)$ and a positive
facet $F$, we introduce a complex of vector spaces $\Phi^+_F(\CK)$ by
$$
\Phi^+_F(\CK)=R\Gamma(D_F,S_F;\CK)[-\dim F].
$$
This is a well defined object of the bounded derived category $\CCD(*)$ of
finite dimensional vector spaces,   
not depending on the choice of points $w_F$. It is called the {\em complex 
of vanishing cycles of $\CK$ across $F$}. 

\subsection{Theorem} 
\label{nul' verdier} {\em We have canonically $\Phi^+_F(D\CK)=D\Phi^+_F(\CK)$
where $D$ denotes the Verdier duality in the corresponding derived 
categories.} $\Box$

\subsection{Theorem} {\em If $\CM\in\CM(D^J;\CS)$ then  
$H^i(\Phi^+_F(\CM))=0$ for $i\neq 0$. Thus, $\Phi^+_F$ induces an exact functor
$$
\Phi^+_F: \CM(D^J;\CS)\lra\Vect
$$
to the category $\Vect$ of finite dimensional vector spaces.} $\Box$

\subsection{} 
\label{nul' dir sum} Given a positive facet $E$
and $\CK\in\CCD(D^J,\CS)$,
we have $S_E=\bigcup_{F\in\CF^1(E)} D_{F}$, the union over the set 
$\CF^1(E)$
of all positive facets $F\supset E$ with $\dim F=\dim F+1$, and 
$$
R\Gamma(S_E,\bigcup_{F\in\CF^1(E)}\ S_F;\CK)=\oplus_{F\in\CF^1(E)}\ 
R\Gamma(D_F,S_F;\CK).
$$

\subsection{} For two positive facets $E$ and $F\in\CF^1(E)$, and 
$\CK(D^J;\CS)$, let us define the natural map
$$
u=u^F_E(\CK): \Phi^+_F(\CK)\lra\Phi^+_E(\CK)
$$  
called {\em canonical}, as the composition
$$
R\Gamma(D_F,S_F;\CK)[-p]\lra R\Gamma(S_E,\bigcup_{F'\in\CF^1(E)} S_{F'};
\CK)[-p]\lra R\Gamma(S_E;\CK)[-p]\lra 
$$
$$
\lra R\Gamma(D_E,S_E)[-p+1]
$$
where $p=\dim F$, the first arrow being induced by the equality in 
\ref{nul' dir sum}, the last one being the coboundary map.

Define the natural map 
$$
v=v^E_F(\CK): \Phi^+_E(\CK)\lra\Phi^+_F(\CK)
$$
called {\em variation}, as the map dual to the composition
$$
D\Phi^+_F(\CK)=\Phi^+_F(D\CK)\overset{u(D\CK)}{\lra}\Phi^+_E(D\CK)=
D\Phi^+_E(\CK).
$$

{\em BACK TO FACTORIZABLE SHEAVES}

\subsection{} Let $\nu\in Y$. We are going to give two equivalent
definitions of an exact functor, called {\em vanishing cycles at the origin}
$$
\Phi: \CM(D^{\nu};\CS)\lra\Vect.
$$
{\em First definition.} Let $f: D^{\nu}\lra D$ be the function 
$f((t_j))=\sum t_j$. For an object $\CK\in\CCD(D^{\nu};\CS)$,
the Deligne's complex of vanishing cycles $\Phi_f(\CK)$
(cf. \cite{d3}) is concentrated at the origin of the hypersurface $f^{-1}(0)$. 
It is $t$-exact with respect to the middle $t$-structure. We set
by definition, $\Phi(\CM)= H^0(\Phi_f(\CM))\ (\CM\in\CM(D^{\nu};\CS))$.

{\em Second definition.} Choose an unfolding of $\nu$, $\pi: J\lra I$. 
Let us consider the canonical projection $\pi: D^J\lra D^{\nu}$. For 
$\CK\in\CCD(D^{\nu};\CS)$, the complex $\pi^*\CK$
is well defined as an element of the $\Sigma_{\pi}$-equivariant 
derived category, hence $\Phi^+_O(\pi^*(\CK))$ is a well defined object
of the $\Sigma_{\pi}$-equivariant derived category of vector spaces 
($O$ being the origin facet in $D^J$). 
Therefore, the complex  of $\Sigma_{\pi}$-invariants $\Phi^+_O(\pi^*\CK)
^{\Sigma_{\pi}}$ is a well defined object of $\CCD(*)$. If
$\CK\in\CM(D^{\nu};\CS)$ then all the cohomology of 
$\Phi^+_O(\pi^*\CK)^{\Sigma_{\pi}}$ in non-zero degree vanishes.
We set 
$$
\Phi(\CM)=H^0(\Phi^+_O(\pi^*\CM)^{\Sigma_{\pi}})\
$$
$(\CM\in\CM(D^{\nu};\CS))$. The equivalence of the two definitions follows 
without difficulty from the proper base change theorem. In computations 
the second definition is used.\footnote{its independence of the choice 
of an unfolding follows from its equivalence to the first definition.}

\subsection{} Let $\CM$ be a factorizable sheaf supported at $c\in X/Y$, 
$\lambda=\lambda(\CM)$.  
For $\nu\in Y^+$, define a vector space $\Phi(\CM)_{\lambda-\nu}$ by
$$
\Phi(\CM)_{\lambda-\nu}=\Phi(\CM^{\nu}).
$$
If $\mu\in X,\ \mu\not\leq\lambda$, set $\Phi(\CM)_{\mu}=0$. One sees easily
that this way we get an exact functor $\Phi$ from $\tFS_c$ to the category
of $X$-graded vector spaces with finite dimensional components. 
We extend it to the whole category $\tFS$ by additivity. 

\subsection{} Let $\CM\in\tFS_c$, $\lambda=\lambda(\CM)$, 
$\nu=\sum\nu_ii\in Y^+$. Let $i\in I$ be such that $\nu_i>0$. Pick an
unfolding of $\nu$, $\pi: J\lra I$. For each $j\in\pi^{-1}(i)$, 
the restriction of $\pi$, $\pi_j: J-\{j\}\lra I$, is an unfolding of 
$\nu-i$. 

For each $j\in\pi^{-1}(i)$, we have canonical and variation morphisms
$$
u_j: \Phi^+_{F_j}(\pi^*\CM^{\nu})\rlh\Phi^+_O(\pi^*\CM^{\nu}): v_j
$$
(the facet $F_j$ has been defined in \ref{nul' facets}).
Taking their sum over $\pi^{-1}(i)$, we get the maps 
$$
\sum u_j:\ \oplus_{j\in\pi^{-1}(i)}\Phi^+_{F_j}(\pi^*\CM^{\nu})\rlh
\Phi^+_O(\pi^*\CM):\sum v_j
$$
Note that the group $\Sigma_{\pi}$ acts on both sides and the maps 
respect this action. After passing to $\Sigma_{\pi}$-invariants, we get
the maps

(a) $\Phi^+_{F_j}(\pi^*\CM^{\nu})^{\Sigma_{\pi_j}}=(\oplus_{j'\in\pi^{-1}(i)}
\Phi^+_{F_{j'}}(\pi^*\CM^{\nu}))^{\Sigma_{\pi}}\rlh
\Phi^+_O(\pi^*\CM^{\nu})^{\Sigma_{\pi}}=\Phi(\CM)_{\nu}$.

Here $j\in\pi^{-1}(i)$ is an arbitrary element. Let us consider the space
$$
F_j^{\perp}=\{(t_{j'})\in D^J|t_j=r, t_{j'}\in D(r')\mbox{ for all }
j'\neq j\} 
$$
where $r, r'$ are some fixed real numbers such that $0<r'<r<1$. 
The space $F^{\perp}_j$ is transversal to $F_j$ and may be identified 
with $D^{J-\{j\}}$. The factorization isomorphism induces the isomorphism 
$$
\pi^*\CM^{\nu}|_{F^{\perp}_j}\cong\pi^*_j\CM^{\nu-i}\otimes(\CI_{\lambda-\nu+i}
^i)_{\{r\}}=\pi^*\CM^{\nu-i}
$$
which in turn induces the isomorphism
$$
\Phi^+_{F_j}(\pi^*\CM^{\nu})\cong\Phi^+_O(\pi^*_j\CM^{\nu-i}).
$$
which is $\Sigma_{\pi_j}$-equivariant. Taking $\Sigma_{\pi_j}$-invariants 
and composing with the maps (a), we get the maps 

(b) $\epsilon_i:\Phi(\CM)_{\nu-i}=\Phi^+_O(\pi^*_j\CM^{\nu-i})^{\Sigma_j}
\rlh \Phi(\CM)_{\nu}: \theta_i$

which do not depend on the choice of $j\in\pi^{-1}(i)$. 

\subsection{} For an arbitrary $\CM\in\tFS$,  
let us define the $X$-graded vector space $\Phi(\CM)$ as
$\Phi(\CM)=\oplus_{\lambda\in X}\Phi(\CM)_{\lambda}$.

\subsection{Theorem} {\em The operators $\epsilon_i, \theta_i\ (i\in I)$ 
acting on the $X$-graded vector space $\Phi(\CM)$ satisfy the relations
~\ref{nul' def c} {\em (a), (b)}.} $\Box$

\subsection{} A factorizable sheaf $\CM$ is finite
iff the space $\Phi(\CM)$ is finite dimensional.   
The previous theorem says that $\Phi$ defines an exact functor 
$$
\Phi: \FS\lra\CC.
$$
One proves that $\Phi$ is a tensor functor.  

\subsection{Example.} For every $\lambda\in X$, the factorizable sheaf 
$\CL(\lambda)$ (cf. \ref{nul' stand}) is finite. It is an irreducible object
of $\FS$, and every irreducible object in $\FS$ is isomorphic to some 
$\CL(\lambda)$. We have $\Phi(\CL(\lambda))=L(\lambda)$.

The next theorem is the main result of the present work. 

\subsection{Theorem} {\em The functor $\Phi$ is an equivalence of  
braided tensor categories.} $\Box$

\subsection{Remark} As a consequence, the category $\FS$ is rigid. We do not
know a geometric construction of the rigidity; it would be very interesting 
to find one.

\newpage
\begin{center}
{\bf Chapter 2. Global (genus $0$)}
\end{center}

\section{Cohesive local systems} 

\subsection{} From now on until the end of part 0 we make the assumptions
of \ref{nul' k big}. In the operadic notations below we partially follow
\cite{bd}. 

\subsection{Operad $\CCD$} For a nonempty finite
set $J$, let $D(J)$ denote the space 
whose points are $J$-tuples $\{x_j,\tau_j\}\ (j\in J)$ where $x_j\in D$ and 
$\tau_j$ is a non-zero tangent vector at $x_j$, such that all points $x_j$ 
are distinct. 

Let $\tD(J)$ be the space whose points are $J$-tuples 
$\{\phi_j\}$ of holomorphic maps $\phi_j: D\lra D\ (j\in J)$, 
each $\phi_j$ having the form $\phi_j(z)=x_j+\tau_jz\ (x_j\in D, \tau_j\in
\BC^*)$, such that $\phi_j(D)\cap\phi_{j'}(D)=\emptyset$ for
$j\neq j'$. We shall identify the $j$-tuple $\{\phi_j\}$ with the
$J$-tuple $\{x_j,\tau_j\}$, and consider $\tau_j$ as a non-zero tangent 
vector from $T_{x_j}D$, thus identifying $T_{x_j}D$ with $\BC$ using
the local coordinate $z-x_j$. So, $\tau_j$ is the image under $\phi_j$ of 
the unit tangent vector at $0$.   
We have an obvious map $p(J):\tD(J)\lra D(J)$ 
which is a homotopy equivalence.  

If $\rho: K\lra J$ is an epimorphic map of finite sets, the composition defines 
a holomorphic map
$$
m(\rho):\ \prod_J \tD(K_j)\times \tD(J)\lra \tD(K)
$$
where $K_j:=\rho^{-1}(j)$. 
If $L\overset{\sigma}{\lra}K\overset{\rho}{\lra}J$ are two epimorphisms of finite 
sets, the square
$$\begin{array}{ccc}
\prod_K \tD(L_k)\times\prod_J \tD(K_j)\times \tD(J)&\overset{m(\sigma)}{\lra}&
\prod_K \tD(L_k)\times\tD(K)\\
\prod m(\sigma_j)\downarrow&\ &\downarrow m(\sigma)\\
\prod_J \tD(L_j)\times\tD(J)&\overset{m(\rho\sigma)}{\lra}&\tD(L)
\end{array}$$
commutes. Here $\sigma_j: L_j\lra K_j$ are induced by $\sigma$.

Let $*$ denote the one element set. 
The space $\tD(*)$ has a marked point, also to be denoted by $*$, corresponding 
to the identity map $\phi: D\lra D$. 

If $\rho: J'\iso J$ is an isomorphism, it induces in the obvious way 
an isomorphism $\rho^*:\tD(J)\iso\tD(J')$ (resp. $D(J)\iso D(J')$). The first 
map coincides with $m(\rho)$ restricted to $(\prod_J *)\times\tD(J)$. 
In particular,
for each $J$, the group $\Sigma_J$ of automorphisms of the set $J$, acts 
on the spaces $\tD(J), D(J)$.   

The map $m(J\lra *)$ restricted to $*\times\tD(J)$ is the identity of $\tD(J)$. 

We will denote the collection of the spaces and maps $\{\tD(J), m(\rho)\}$ 
by $\CCD$, and call it the {\em operad of disks with tangent vectors}.

\subsection{Coloured local systems over $\CCD$} If $\rho: K\lra J$
is an epimorphism of finite sets and $\pi: K\lra X$ is a map of sets, 
we define the map $\rho_*\pi: J\lra X$ by $\rho_*\pi(j)=\sum_{K_j}\pi(k)$.
For $j\in J$, we denote $K_j:=\rho^{-1}(j)$ as above, and $\pi_j:\ K_j\lra X$ 
will denote the restriction of $\pi$.

Let us call an {\em $X$-coloured local system} $\CJ$ over $\CCD$
a collection of local systems $\CJ(\pi)$ over the spaces $\tD(J)$ given 
for every map $\pi:J\lra X$, $J$ being a non-empty finite set, together 
with {\em factorization isomorphisms}
$$
\phi(\rho):\ m(\rho)^*\CJ(\pi)\iso\Boxtimes_J\CJ(\pi_j)\boxtimes\CJ(\rho_*\pi)
$$
given for every epimorphism $\rho: K\lra J$ and $\pi: K\lra X$, which
satisfy the properties (a), (b) below.  

(a) {\em Associativity}. Given a map $\pi: L\lra X$ and a pair of epimorphisms
$L\overset{\sigma}{\lra}K\overset{\rho}{\lra}J$, the square below commutes.
$$\begin{array}{ccc}
m(\rho)^*m(\sigma)^*\CJ(\pi)&\overset{\phi(\sigma)}{\lra}&\Boxtimes_K\CJ(\pi_k)
\boxtimes m(\rho)^*\CJ(\sigma_*\pi)\\
\phi(\rho\sigma)\downarrow&\ &\downarrow\phi(\rho)\\
\Boxtimes_J m(\sigma_j)^*\CJ(\pi_j)\boxtimes\CJ(\rho_*\sigma_*\pi)&\overset
{\boxtimes\phi(\sigma_j)}{\lra}&\Boxtimes_K\CJ(\pi_k)\boxtimes\Boxtimes_J
\CJ((\sigma_*\pi)_j)\boxtimes\CJ(\rho_*\sigma_*\pi)
\end{array}$$
Note that $\pi_{j*}\sigma_j=(\sigma_*\pi)_j$.

For $\mu\in X$, let $\pi_{\mu}: *\lra X$ be defined by $\pi_{\mu}(*)=\mu$. 
The isomorphisms $\phi(\id_*)$ restricted to the marked points 
in $D(*)$, give the isomorphisms $\CJ(\pi_{\mu})_*\iso \sk$ (and imply that
the local systems $\CJ(\pi_{\mu})$ are one-dimensional). 

(b) For any $\pi: J\lra X$, the map $\phi(\id_J)$ restricted to 
$(\prod_J *)\times \tD(J)$, equals $\id_{\CJ(\pi)}$. 

The map $\phi(J\lra *)$ restricted to $*\times\tD(J)$, equals $\id_{\CJ(\pi)}$. 

\subsection{} The definition above implies that the local systems $\CJ(\pi)$ 
are functorial with respect to isomorphisms. In particular,
the action of the group $\Sigma_{\pi}$ on $\tD(J)$ lifts to $\CJ(\pi)$. 

\subsection{Standard local system over $\CCD$}
\label{nul' stand d} Let us define the "standard" local systems
$\CJ(\pi)$ by a version of the construction \ref{nul' stand loc}.

We will use
the notations of {\em loc. cit.} for the marked points 
and paths. We can identify $D(J)=D^{J\circ}\times(\BC^*)^J$ where
we have identified all tangent spaces $T_xD\ (x\in D)$ with $\BC$ using
the local coordinate $z-x$. 

We set $\CJ(\pi)_{x(\sigma)}=\sk$, the monodromies being
$T_{\gamma(j',j'')^{\pm}}=\zeta^{\mp\pi(j')\cdot\pi(j'')}$, and the monodromy 
$T_j$  
corresponding to the counterclockwise circle of a tangent vector $\tau_j$ 
is $\zeta^{-2n\pi(j)}$. The factorization isomorphisms 
are defined by the same condition as in {\em loc. cit}.

This defines 
the {\em standard local system $\CJ$} over $\CCD$. Below, the notation $\CJ$
will be reserved for this local system.

\subsection{} 
\label{nul' right} Set $P=\BP^1(\BC)$. We pick a point $\infty\in P$ and
choose a global coordinate $z:\BA^1(\BC)=P-\{\infty\}\iso\BC$.
This gives local coordinates: $z-x$ at $x\in\BA^1(\BC)$ and
$1/z$ at $\infty$. 

For a non-empty finite 
set $J$, let $P(J)$ denote the space of $J$-tuples $\{ x_j,\tau_j\}$ where
$x_j$ are distinct points on $P$, and $\tau_j$ is a non-zero tangent vector 
at $x_j$. Let $\tP(J)$ denote the space whose points are $J$-tuples of 
holomorphic embeddings $\phi_j: D\lra P$ with non-intersecting images 
such that each $\phi_j$ is a restriction of an algebraic morphism $P\lra P$. 
We have the $1$-jet projections $\tP(J)\lra P(J)$. We will use the notation 
$P(n)$ for $P(\{1,\ldots,n\})$, etc.   

An epimorphism $\rho: K\lra J$ induces the maps
$$
m_P(\rho):\ \prod_J \tD(K_j)\times \tP(J)\lra\tP(K)
$$
and
$$
\bar{m}_P(\rho):\ \prod_J D(K_j)\times\tP(J)\lra P(K).
$$
For a pair of epimorphisms $L\overset{\sigma}{\lra}K\overset{\rho}{\lra}J$, 
the square
$$\begin{array}{ccc}
\prod_K \tD(L_k)\times\prod_J \tD(K_j)\times \tP(J)&\overset{m_P(\sigma)}{\lra}&
\prod_K \tD(L_k)\times\tP(K)\\
\prod m(\sigma_j)\downarrow&\ &\downarrow m_P(\sigma)\\
\prod_J \tD(L_j)\times\tP(J)&\overset{m_P(\rho\sigma)}{\lra}&\tP(L)
\end{array}$$
commutes. 

If $\rho: J'\iso J$ is an isomorphism, it induces in the obvious way 
an isomorphism $\rho^*:\tP(J)\iso\tP(J')$ (resp. $P(J)\iso P(J')$). This 
last map 
coincides with $m_P(\rho)$ (resp. $\bar{m}_P(\rho)$)  
restricted to $(\prod_J *)\times\tP(J)$ (resp. $(\prod_J *)\times P(J)$). 

The collection of the spaces and maps $\{\tP(J), m_P(\rho)\}$ form an 
object $\tP$ called a {\em right module} over the operad $\CCD$.

\subsection{} 
\label{nul' cls def} Fix an element $\mu\in X$. Let us say that a map $\pi: J\lra X$
has {\em level $\mu$} if $\sum_J\pi(j)=\mu$.  

A {\em cohesive local
system $\CH$ of level $\mu$} over $P$  
is a collection of local systems $\CH(\pi)$ over 
$\tP(J)$ given for every $\pi: J\lra X$ of level $\mu$,
together with {\em factorization isomorphisms}
$$
\phi_P(\rho): m_P(\rho)^*\CH(\pi)\iso\Boxtimes_J\CJ(\pi_j)\boxtimes
\CH(\rho_*\pi)
$$
given for every epimorphism $\rho: K\lra J$ and $\pi: K\lra X$ of level 
$\mu$, which satisfy the properties (a), (b) below. 

(a) {\em Associativity}. Given a map $\pi: L\lra X$ and a pair of epimorphisms
$L\overset{\sigma}{\lra}K\overset{\rho}{\lra}J$, the square below commutes.
$$\begin{array}{ccc}
m_P(\rho)^*m_P(\sigma)^*\CH(\pi)&\overset{\phi_P(\sigma)}{\lra}&
\Boxtimes_K\CJ(\pi_k)
\boxtimes m_P(\rho)^*\CH(\sigma_*\pi)\\
\phi_P(\rho\sigma)\downarrow&\ &\downarrow\phi_P(\rho)\\
\Boxtimes_J m(\sigma_j)^*\CJ(\pi_j)\boxtimes\CH(\rho_*\sigma_*\pi)&\overset
{\boxtimes\phi(\sigma_j)}{\lra}&\Boxtimes_K\CJ(\pi_k)\boxtimes\Boxtimes_J
\CJ((\sigma_*\pi)_j)\boxtimes\CH(\rho_*\sigma_*\pi)
\end{array}$$

(b) For any $\pi: J\lra X$, the map $\phi_P(\id_J)$ restricted to 
$(\prod_J *)\times \tP(J)$, equals $\id_{\CH(\pi)}$.

\subsection{} The definition above implies that the local systems $\CH(\pi)$
are functorial with respect to isomorphisms. In particular,
the action of the group $\Sigma_{\pi}$ on $\tP(J)$, lifts to $\CH(\pi)$. 

\subsection{Theorem} {\em For each $\mu\in X$ such that 
$\mu\equiv\ 2\rho_{\ell}\ \modul Y_{\ell}$, there exists a unique 
up to an isomorphism one-dimensional cohesive local system $\CH=\CH^{(\mu)}$ 
of level $\mu$ over $P$}. $\Box$

The element $\rho_{\ell}\in X$ is defined in \ref{nul' l} and the lattice
$Y_{\ell}$ in \ref{nul' lattice}.

>From now on, let us fix such a local system $\CH^{(\mu)}$ for each $\mu$ as 
in the theorem. 

\subsection{} Note that the obvious maps $p(J): \tP(J)\lra P(J)$ are homotopy 
equivalences. Therefore the local systems $\CH(\pi)\ (\pi: J\lra X)$ 
descend to the unique local systems over $P(J)$, to be denoted by the same 
letter.

\section{Gluing} 

\subsection{} 
\label{nul' diag} Let us fix a finite set $K$. For $\nu\in Y^+$, pick
an unfolding of $\nu$, $\pi: J\lra I$. Consider the space $P^{\nu}$;
its points are formal linear combinations $\sum_J \pi(j) x_j\ 
(x_j\in P)$. We define the space $P^{\nu}(K)=P(K)\times P^{\nu}$; its points 
are tuples $\{y_k,\tau_k,\sum\pi(j)x_j\}\ (k\in K, y_k\in P, \tau_k\neq 0$ in
$T_{y_k}P)$, all $y_k$ being distinct. Let $P^{\nu}(K)^{\circ}$
(resp. $P^{\nu}(K)^{\bullet}$)  
be the subspace whose points are tuples as above with all $x_j$ distinct 
from $y_k$ and pairwise distinct (resp. all $x_j$ distinct from $y_k$).
We will use the notation $P^{\nu}(n)$ for $P^{\nu}(\{1,\ldots,n\})$, etc.  

Let $\nu'\in Y^+$ and $\vnu=\{\nu_k\}\in (Y^+)^K$ be such that
$\sum_K\nu_k+\nu'=\nu$. Define the space $P^{\vnu;\nu'}\subset\tP(K)\times
P^{\nu}$ consisting of tuples $\{\phi_k,\sum\pi(j)x_j\}$ such that
for each $k$, $\nu_k$ of the points $x_j$ lie inside $\phi_k(D)$ and
$\nu'$ of them lie outside all closures of these disks. 

Let $\vzero\in(Y^+)^K$ be the zero $K$-tuple. Define the space 
$\tP^{\nu}(K):=P^{\vzero;\nu}$.  
We have obvious maps 

(a) $\prod_K D^{\nu_k}\times\tP^{\nu'}(K)\overset{p(\vnu;\nu')}{\lla}
P^{\vnu;\nu'}\overset{m(\vnu;\nu')}{\lra}P^{\nu}(K)$. 

Let $\vnu^1, \vnu^2\in(Y^+)^K$ be such that $\vnu^1+\vnu^2=\vnu$. Define the 
space 

$P^{\vnu^1;\vnu^2;\nu}\subset\tP(K)\times\BR^K_{>0}\times\tP(K)
\times P^{\nu}$ consisting of all tuples $\{\phi_k; r_k; \phi'_k;
\sum\pi(j)x_j\}$ such that $r_k<1; \phi_k(z)=\phi'_k(r_kz)$,
and $\nu^1_k$ (resp. $\nu^2_k$, $\nu'$) from the points $x_j$ lie inside
$\phi_k(D)$ (resp. inside the annulus $\phi'_k(D)-\overline{\phi_k(D)}$,
inside $P-\bigcup\overline{\phi'_k(D)}$). We set $\tP^{\vnu;\nu'}:=
P^{\vzero;\vnu;\nu'}$, cf. \ref{nul' spaces}.

We have a commutative romb (cf. \ref{nul' rhomb}):

(b) \begin{center}
  \begin{picture}(14,8)
    \put(7,8){\makebox(0,0){$P^{\nu}(K)$}}
    \put(5,6){\makebox(0,0){$P^{\vnu^1;\nu^2+\nu'}$}}
    \put(9,6){\makebox(0,0){$P^{\vnu^1+\vnu^2;\nu'}$}}
    \put(3,4){\makebox(0,0){$\prod D^{\nu^1_k}\times\tP^{\nu^2+\nu'}(K)$}}
    \put(7,4){\makebox(0,0){$P^{\vnu^1;\vnu^2;\nu'}$}}
    \put(11,4){\makebox(0,0){$\prod D^{\nu^1_k+\nu^2_k}\times\tP^{\nu'}(K)$}}
    \put(5,2){\makebox(0,0){$\prod D^{\nu^1_k}\times\tP^{\vnu^2;\nu'}$}}
    \put(9,2){\makebox(0,0){$\prod D^{\nu^1_k,\nu^2_k}\times\tP^{\nu'}(K)$}}
    \put(7,0){\makebox(0,0){$\prod D^{\nu^1_k}\times\prod \DD^{\nu^2_k}
              \times\tP^{\nu'}(K)$}}

    \put(5.5,6.5){\vector(1,1){1}}
    \put(8.5,6.5){\vector(-1,1){1}}
    \put(4.5,5.5){\vector(-1,-1){1}}
    \put(6.5,4.5){\vector(-1,1){1}}
    \put(7.5,4.5){\vector(1,1){1}}
    \put(9.5,5.5){\vector(1,-1){1}} 
    \put(4.5,2.5){\vector(-1,1){1}}
    \put(6.5,3.5){\vector(-1,-1){1}}
    \put(7.5,3.5){\vector(1,-1){1}}
    \put(9.5,2.5){\vector(1,1){1}}
    \put(5.5,1.5){\vector(1,-1){1}}
    \put(8.5,1.5){\vector(-1,-1){1}}

    \put(5.5,7){\makebox(0,0){$m$}}
    \put(8.5,7){\makebox(0,0){$m$}}
    \put(3.5,5){\makebox(0,0){$p$}}
    \put(5.5,5){\makebox(0,0){$m$}}
    \put(8.5,5){\makebox(0,0){$m$}}
    \put(10.5,5){\makebox(0,0){$p$}}
    \put(3.5,3){\makebox(0,0){$m$}}
    \put(5.5,3){\makebox(0,0){$p$}}
    \put(8.5,3){\makebox(0,0){$p$}}
    \put(10.5,3){\makebox(0,0){$m$}}
    \put(5.5,1){\makebox(0,0){$p$}}
    \put(8.5,1){\makebox(0,0){$p$}}
    
  \end{picture} 
\end{center}

Here $\nu^2:=\sum\nu^2_k$.

\subsection{} Let $P(K;J)$ be the space consisting of tuples $\{y_k;\tau_k;
x_j\}\ (k\in K, j\in J, y_k, x_j\in P, \tau_k\neq 0$ in $T_{y_k}P$), where
all points $x_k, y_j$ are distinct. We have
$P^{\nu}(K)^{\circ}=P(K;J)/\Sigma_{\pi}$.  
We have an obvious projection $P(K\coprod J)\lra P(K;J)$.

\subsection{} 
\label{nul' sign p} Let $\{\CM_k\}$ be a $K$-tuple of factorizable sheaves
supported at some cosets in $X/Y$; let $\mu_k=\lambda(\CM_k)$.

Let $\tpi: K\coprod J\lra X$ be a map defined by $\tpi(k)=\mu_k$,
$\tpi(j)=-\pi(j)\in I\hra X$. The local system $\CH(\tpi)$ over $P(K\coprod J)$ 
descends to $P(K;J)$ since $\zeta^{2n(-i)}=1$ for all $i\in I$;  
this one in turn descends to the unique local system 
$\tCH_{\vmu}^{\nu}$ over $P^{\nu}(K)^{\circ}$, 
due to $\Sigma_{\pi}$-equivariance. Let us define the local system 
$\CH_{\vmu}^{\nu}:=\tCH_{\vmu}^{\nu}\otimes\Sign^{\nu}$.  
Here $\Sign^{\nu}$ denotes the inverse image of the sign local 
system on $P^{\nu\circ}$ (defined in the same manner as for the disk, cf. 
\ref{nul' sign}) under the forgetful map $P^{\nu}(K)^{\circ}\lra P^{\nu\circ}$.
    
Let $\CH^{\nu\bullet}_{\vmu}$
be the perverse sheaf over $P^{\nu}(K)^{\bullet}$ which is  
the middle extension of $\CH^{\nu}_{\vmu}[\dim P^{\nu}(K)]$.   
Let us denote by the same letter the inverse image of this perverse sheaf 
on the space $\tP^{\nu}(K)$ with respect to the evident projection 
$\tP^{\nu}(K)\lra P^{\nu}(K)^{\bullet}$. 

\subsection{} 
\label{nul' adm} Let us call an element $\nu\in Y^+$ {\em admissible}
(for a $K$-tuple $\{\mu_k\}$)
if $\sum\mu_k-\nu\equiv\ 2\rho_{\ell}\ \modul Y_{\ell}$, see \ref{nul' lattice}.

\subsection{Theorem - definition}
\label{nul' glu thm} {\em For each admissible $\nu$, there exists
a unique, up to a unique isomorphism, perverse sheaf, denoted by\  
$\Boxtimes^{(\nu)}_K\ \CM_k$, over $P^{\nu}(K)$, equipped with
isomorphisms 
$$
\psi(\vnu;\nu'): m^*\Boxtimes^{(\nu)}_K\ \CM_k
\iso p^*(\Boxtimes_K\ \CM_k^{\nu_k}
\boxtimes\CH^{\nu'\bullet}_{\vmu-\vnu})
$$
given for every diagram \ref{nul' diag} (a) such that for each rhomb
\ref{nul' diag} (b) the cocycle condition
$$
\phi(\vnu^2;\nu')\circ\psi(\vnu^1;\nu^2+\nu')=
(\Boxtimes_K\ \psi^{\CM_k}(\nu^1_k,\nu^2_k))\circ\psi(\vnu^1+\vnu^2;\nu')
$$
holds.} $\Box$

\subsection{} The sheaf $\Boxtimes^{(\nu)}_K\ \CM_k$ defines for each $K$-tuple
of $\vy=\{y_k,\tau_k\}$ of points of $P$ with non-zero tangent vectors,
the sheaf $\Boxtimes^{(\nu)}_{\vy}\ \CM_k$ over $P^{\nu}$, to be called the
{\em gluing of the factorizable sheaves $\CM_k$ into the points $(y_k,\tau_k)$}.

\subsection{Example} 
\label{nul' middle} The sheaf $\Boxtimes^{(\nu)}_K\ \CL(\mu_k)$ is equal
to the middle extension of the sheaf $\CH_{\vmu}^{\nu\bullet}$.

\section{Semiinfinite cohomology} 

In this section we review the theory of semiinfinite cohomology 
in the category $\CC$, due to S. M. Arkhipov, cf. \cite{a}.

\subsection{} Let $\CC_r$ be a category whose objects are {\em right} 
$\fu$-modules $N$, finite dimensional over $\sk$, with a given $X$-grading
$N=\oplus_{\lambda\in X}\ N_{\lambda}$ such that 
$xK_{\nu}=\zeta^{-\langle\nu,\lambda\rangle}x$ for any $\nu\in Y, 
\lambda\in X, x\in N_{\lambda}$. All definitions and results
concerning the
category $\CC$ given above and below, have the obvious versions 
for the category $\CC_r$. 

For $M\in\CC$, define $M^{\vee}\in\CC_r$ as follows: $(M^{\vee})_{\lambda}=
(M_{-\lambda})^*$ (the dual vector space); the action of the operators 
$\theta_i, \epsilon_i$ being the transpose of their action on $M$. This 
way we get an equivalence $^{\vee}:\CC^{\opp}\iso\CC_r$.         

\subsection{} Let us call an object $M\in\CC$ $\fu^-$- (resp. $\fu^+$-)
{\em good} if it admits a filtration whose successive quotients 
have the form $\ind_{\fu^{\geq 0}}^{\fu}(M')$ (resp.
$\ind_{\fu^{\leq 0}}^{\fu}(M'')$) for some $M'\in\CC^{\geq 0}$ 
(resp. $M''\in\CC^{\leq 0}$) (cf. \ref{nul' ind}). These classes of objects
are stable with respect to the tensor multiplication by an arbitrary object 
of $\CC$.  

If $M$ is $\fu^-$- (resp. $\fu^+$-) good then $M^*$ is $\fu^-$-
(resp. $\fu^+$-) good. 

If $M$ is $\fu^-$-good and $M'$ is $\fu^+$-good then $M\otimes M'$ is a 
projective object in $\CC$. 

\subsection{} Let us say that a complex $M^{\bullet}$ in $\CC$ is 
{\em concave} (resp. {\em convex}) if 

(a) there exists $\mu\in X$ such that all nonzero components 
$M^{\bullet}_{\lambda}$ have the weight $\lambda\geq\mu$ (resp. 
$\lambda\leq\mu$);  

(b) for any $\lambda\in X$, the complex $M^{\bullet}_{\lambda}$ is finite. 

\subsection{} For an object $M\in\CC$, we will call a {\em left} 
(resp. {\em right}) {\em $\fu^{\pm}$-good resolution} of $M$  an exact
complex 
$$
\ldots\lra P^{-1}\lra P^0\lra M\lra 0
$$
(resp.
$$ 
0\lra M\lra R^0\lra R^1\lra\ldots)
$$
such that all $P^i$ (resp. $R^i$) are $\fu^{\pm}$-good. 

\subsection{Lemma.} {\em Each object $M\in\CC$ admits a convex $\fu^-$-good
left resolution, a concave $\fu^+$-good left resolution, a concave 
$\fu^-$-good right resolution and a convex $\fu^+$-good right resolution.} 
$\Box$

\subsection{} 
\label{nul' pairing} For $N\in\CC_r, M\in\CC$, define a vector space $N\otimes_{\CC}M$
as the zero weight component of the tensor product $N\otimes_{\fu}M$ 
(which has an obvious $X$-grading). 

For $M, M'\in\CC$, we have an obvious perfect pairing 

(a) $\Hom_{\CC}(M,M')\otimes(M^{\prime\vee}\otimes_{\CC}M)\lra \sk$.

\subsection{} $M, M'\in\CC, N\in\CC_r$ and $n\in\BZ$, define the
{\em semiinifinite $\Ext$ and $\Tor$} spaces
$$
\Ext^{\infh+n}_{\CC}(M,M')=H^n(\Hom_{\CC}(R^{\bullet}_{\searrow}(M),
R^{\bullet}_{\nearrow}(M')))
$$
where $R^{\bullet}_{\searrow}(M)$ (resp. $R^{\bullet}_{\nearrow}(M')$) 
is an arbitrary $\fu^+$-good convex right resolution of $M$ 
(resp. $\fu^-$-good concave right resolution of $M'$),
$$
\Tor^{\CC}_{\infh+n}(N,M)=H^{-n}(P^{\bullet}_{\swarrow}(N)\otimes_{\CC}
R^{\bullet}_{\searrow}(M))
$$
where $P^{\bullet}_{\swarrow}(N)$ is an arbitrary $\fu^-$-good convex left 
resolution of $N$. 

This definition does not depend, up to a unique isomorphism, 
upon the choice of resolutions, and is functorial.
  
These spaces are finite dimensional and are non-zero only for finite
number of degrees $n$. 

The pairing \ref{nul' pairing} (a) induces perfect pairings
$$
\Ext^{\infh+n}_{\CC}(M,M')\otimes\Tor^{\CC}_{\infh+n}(M^{\prime\vee},M)\lra \sk\
(n\in\BZ).
$$

\section{Conformal blocks (genus $0$)}

In this section we suppose that $\sk=\BC$.

\subsection{} 
\label{nul' trivial} Let $M\in \CC$. We have a canonical embedding of vector
spaces 
$$
\Hom_{\CC}(\One,M)\hra M
$$
which identifies $\Hom_{\CC}(\One,M)$ with the maximal trivial
subobject of $M$. Here "trivial" means "isomorphic to a sum of a few 
copies of the object $\One$".   
Dually, we have a canonical epimorhism
$$
M\lra\Hom_{\CC}(M,\One)^*
$$
which identifies $\Hom_{\CC}(M,\One)^*$ with the maximal trivial quotient
of $M$. Let us denote by $\langle M\rangle$ the image
of the composition
$$
\Hom_{\CC}(\One,M)\lra M\lra\Hom_{\CC}(M,\One)^*
$$
Thus, $\langle M\rangle$ is canonically a subquotient of $M$.

One sees easily that if $N\subset M$ is a trivial direct summand of $M$
which is maximal, i.e. not contained in greater direct summand, then
we have a canonical isomorphism $\langle M\rangle\iso N$. For this reason,
we will call $\langle M\rangle$ {\em the maximal trivial
direct summand} of $M$.

\subsection{} Let $\gamma_0\in Y$ denote the highest coroot and $\beta_0\in Y$ 
denote the 
coroot dual to the highest root ($\gamma_0=\beta_0$ for a simply laced 
root datum).
 
Let us define the {\em first alcove} $\Delta_{\ell}\subset X$ by 
$$
\Delta_{\ell}=\{\lambda\in X|\ \langle i,\lambda+\rho\rangle>0\mbox{ for all }
i\in I;\ \langle\gamma_0,\lambda+\rho\rangle<\ell\}
$$
if $d$ does not divide $\ell$, i.e. if $\ell_i=\ell$ for all $i\in I$, 
and by
$$
\Delta_{\ell}=\{\lambda\in X|\ \langle i,\lambda+\rho\rangle >0\mbox{ for all }
i\in I;\ \langle\beta_0,\lambda+\rho\rangle <\ell_{\beta_0}\}
$$
otherwise, cf. \cite{ap} 3.19 ($d$ is defined in \ref{nul' root},
and $\ell_{\beta_0}$ in \ref{nul' l}). Note that $\ell_{\beta_0}=\ell/d$.

\subsection{} For $\lambda_1,\ldots,\lambda_n\in\Delta_{\ell}$, define the 
{\em space of conformal blocks} by
$$
\langle L(\lambda_1),\ldots,L(\lambda_n)\rangle :=
\langle L(\lambda_1)\otimes\ldots\otimes L(\lambda_n)\rangle.
$$
In fact, due to the ribbon structure on $\CC$, the right hand 
side is a {\em local system} over the space $P(n):=
P(\{1,\ldots,n\})$ (cf. \cite{d1}). It is more appropriate to 
consider the previous equality as the definition of the 
{\em local system of conformal blocks} over $P(n)$.

\subsection{Theorem} 
\label{nul' ar} (Arkhipov) {\em For each $\lambda_1,\ldots,\lambda_n
\in\Delta_{\ell}$, the space of conformal blocks  
$\langle L(\lambda_1),\ldots,L(\lambda_n)\rangle$ is naturally a 
subquotient of the space

$\Tor^{\CC}_{\infh+0}({\em\One}_r, L(\lambda_1)\otimes\ldots\otimes 
L(\lambda_n)\otimes L(2\rho_{\ell}))$.

More precisely, due to the ribbon structure on $\CC$, the latter space
is a stalk 
of a local system over $P(n+1)$, and inverse image of the local system
$\langle L(\lambda_1),\ldots,L(\lambda_n)\rangle$ under the projection 
onto the first coordinates $P(n+1)\lra P(n)$, is a natural 
subquotient of this local system.} 

Here $\One_r$ is the unit object in $\CC_r$.

Examples, also due to Arkhipov, show that the local systems of conformal blocks
are in general {\em proper} subquotients of the corresponding  
$\Tor$ local systems.

This theorem is an immediate consequence of the next lemma, 
which in turn follows from the geometric theorem \ref{nul' integral} below,
cf. \ref{nul' delta}.

\subsection{Lemma} 
\label{nul' lemma rho}{\em We have $\Tor^{\CC}_{\infh+n}({\em\One}_r,
L(2\rho_{\ell}))=\sk$ if $n=0$, and $0$ otherwise.} $\Box$

\subsection{} 
\label{nul' charge} Let $\fg$ be the simple Lie algebra (over $\sk$) associated
with our Cartan datum; let $\hfg$ be the corresponding affine 
Lie algebra.

Let $\MS$ denote the category of integrable $\hfg$-modules of central 
charge $\kappa-h$. Here $\kappa$ is a fixed positive integer, $h$ is the 
dual Coxeter number of our Cartan datum. $\MS$ is a semisimple abelian  
category whose irreducible objects are 
$\fL(\lambda),\ \lambda\in\Delta_{\ell}$ where  
$l=2d\kappa$, 
i.e. $\ell=d\kappa$ (we are grateful to Shurik Kirillov who pointed 
out the necessity of the factor $d$ here) and $\fL(\lambda)$ is the
highest weight module 
with a highest vector $v$ whose "top" part $\fg\cdot v$ is the irreducible
$\fg$-module of the highest weight $\lambda$. 

According to Conformal field theory, $\MS$ has a natural structure
of a ribbon category, cf. \cite{ms}, \cite{k}.

The usual local systems of conformal blocks in the WZW model may be defined as  
$$
\langle\fL(\lambda_1),\ldots,\fL(\lambda_n)\rangle=
\Hom_{\MS}(\One,\fL(\lambda_1)\otimes\ldots\otimes\fL(\lambda_n))
$$
the structure of a local system on the right hand side is due to the ribbon
structure on $\MS$.  

\subsection{}
\label{nul' ms} Let  $\zeta=\exp(\pi\sqrt{-1}/d\kappa)$.
We have an exact functor 
$$
\phi:\ \MS\lra\CC 
$$
sending $\fL(\lambda)$ to $L(\lambda)$. This functor identifies
$\MS$ with a full subcategory of $\CC$. 

The functor $\phi$ does not respect the tensor structures. It admits the 
left and right adjoints, $\phi^{\flat}, \phi^{\sharp}$. For $M\in\CC$, let 
$\langle M\rangle_{\MS}$ denotes the image of the composition 
$$
\phi\circ\phi^{\flat}(M)\lra M\lra\phi\circ\phi^{\sharp}(M).
$$
We have the following comparison theorem.

\subsection{Theorem}{\em We have naturally
$$
\phi(\fM\otimes\fM')=\langle\phi(\fM)\otimes\phi(\fM')\rangle_{\MS}.
$$}

This follows from the combination of the results of \cite{ap}, \cite{kl},
\cite{l3} and \cite{f}. 

\subsection{Corollary.} {\em For any $\lambda_1,\ldots,\lambda_n\in\Delta_{\ell}$, 
the functor $\phi$ induces an isomorphism of local systems 
$$
\langle\fL(\lambda_1),\ldots,\fL(\lambda_n)\rangle=
\langle L(\lambda_1),\ldots,L(\lambda_n)\rangle.\ \Box
$$}

\newpage
\section{Integration}

We keep the notations of the previous section.

\subsection{} Let $K$ be a finite set, $m=\card(K)$, $\{\CM_k\}$ a $K$-tuple
of finite factorizable sheaves, $\CM_k\in\FS_{c_k}$, $\mu_k:=\lambda
(\CM_k)$. Assume that $\nu:=\sum_K\ \mu_k-2\rho_{\ell}$ belongs to $Y^+$.

Let $\eta: P^{\nu}(K)\lra P(K)$ be the projection. 

\subsection{Theorem} 
\label{nul' integral} {\em We have canonical isomorphisms of local systems
over $P(K)$ 
$$
R^{a-2m}\eta_*(\Boxtimes^{(\nu)}_K\ \CM_k) =\Tor^{\CC}_{\infh-a}(\em{\One}_r,
\otimes_K\ \Phi(\CM_k))\ (a\in\BZ),
$$
the structure of a local system on the right hand side being induced 
by the ribbon structure on $\CC$.} $\Box$

\subsection{Proof of Lemma \ref{nul' lemma rho}}
\label{nul' delta} Apply the previous theorem
to the case when the $K$-tuple consists of one sheaf $\CL(2\rho_{\ell})$ and 
$\nu=0$. $\Box$

\subsection{} From now until the end of the section, $\sk=\BC$ and
$\zeta=\exp(\pi\sqrt{-1}/d\kappa)$.    
Let $\lambda_1,\ldots,\lambda_n\in\Delta_{\ell}$.
Let $\nu=\sum_{m=1}^n\ \lambda_m$; assume that $\nu\in Y^+$. Set   
$\vmu=\{\lambda_1,\ldots,\lambda_n,2\rho_{\ell}\}$.

Let $\eta$ be the projection   
$P^{\nu}(n+1):=P^{\nu}(\{1,\ldots,n+1\})\lra P(n+1)$ and 
$p: P(n+1)\lra P(n)$ be the projection on the first coordinates. 
 
Let $\CH^{\nu\sharp}_{\vmu}$ denote the middle extension of the sheaf 
$\CH^{\nu\bullet}_{\vmu}$. By the Example \ref{nul' middle},
$$
\CH^{\nu\sharp}_{\vmu}=
\Boxtimes_{1\leq a\leq n}^{(\nu)}\CL(\lambda_a)\boxtimes
\CL(2\rho_{\ell}).
$$   
\subsection{Theorem} 
\label{nul' subq} {\em The local system $p^*\langle L(\lambda_1),\ldots,
L(\lambda_n)\rangle$ is canonically a subquotient of the local 
system 
$$
R^{-2n-2}\eta_*\CH^{\nu\sharp}_{\vmu}.
$$} 

This theorem is an immediate corollary of the previous one and of the Theorem  
\ref{nul' ar}.

\subsection{Corollary} {\em In the notations of the previous theorem, 
the local system $\langle L(\lambda_1),\ldots,L(\lambda_n)\rangle$ 
is semisimple.} 

{\bf Proof.} The local system 
$p^*\langle L(\lambda_1),\ldots,L(\lambda_n)\rangle$ is a subquotient 
of the geometric local system $R^{-2n-2}\eta_*\CH^{\nu\sharp}_{\vmu}$, and 
hence is semisimple by
the Beilinson-Bernstein-Deligne-Gabber
Decomposition theorem, \cite{bbd}, Th\'{e}or\`{e}me 6.2.5. Therefore, 
the local system $\langle L(\lambda_1),\ldots,L(\lambda_n)\rangle$ is 
also semisimple, since the map $p$ induces the surjection on the fundamental 
groups. $\Box$

\subsection{} 
\label{nul' hermit} For a sheaf $\CF$, let
$\bar{\CF}$ denote the sheaf obtained from $\CF$ by the complex conjugation   
on the coefficients. 

If a perverse sheaf $\CF$ on $P^{\nu}$ is obtained by gluing some 
irreducible factorizable sheaves into some points of $P$
then its Verdier dual $D\CF$ is canonically isomorphic to $\bar{\CF}$. 
Therefore, the Poincar\'{e}-Verdier duality induces a perfect  
hermitian pairing on $R\Gamma(P^{\nu};\CF)$. 

Therefore,  
in notations of theorem \ref{nul' subq},
The Poincar\'{e}-Verdier duality 
induces a non-degenerate hermitian form on the local system 
$R^{-2n-2}\eta_*\CH^{\nu\sharp}_{\vmu}$.

By a little more elaborated argument using fusion, one can
introduce a canonical hermitian form on the 
systems of conformal blocks.    

Compare \cite{k}, where a certain hermitian form on the spaces
of conformal blocks (defined up to a positive constant) has been 
introduced. 

\subsection{} 
By the similar reasons, the Verdier duality defines a hermitian form
on all irreducible objects of $\CC$ (since the Verdier duality commutes 
with $\Phi$, cf. Theorem \ref{nul' verdier}).

\section{Regular representation}

\subsection{} 
\label{nul' modif} From now on we are going to modify slightly the definition
of the categories $\CC$ and $\FS$. Let $X_{\ell}$ be the lattice 
$$
X_{\ell}=\{\mu\in X\otimes\Bbb{Q}|\ \mu\cdot Y_{\ell}\in\ell\BZ\}
$$
We have obviously $X\subset X_{\ell}$, and $X=X_{\ell}$ if $d|\ell$. 

In this Chapter we will denote by $\CC$ a category of {\em $X_{\ell}$-graded}
(instead of $X$-graded) finite dimensional vector spaces $M=\oplus
_{\lambda\in X_{\ell}} M_{\lambda}$ equipped with linear operators
$\theta_i: M_{\lambda}\lra M_{\lambda-i'},\ \epsilon_i: M_{\lambda}
\lra M_{\lambda+i'}$ which satisfy the relations \ref{nul' def c} (a), (b).
This makes sense since $\langle d_ii,\lambda\rangle=i'\cdot\lambda\in\BZ$
for each $i\in I, \lambda\in X_{\ell}$. 

Also, in the definition of $\FS$ we replace $X$ by $X_{\ell}$. All the results
of the previous Chapters hold true {\em verbatim} with this modification.

We set $\dd_{\ell}=\card(X_{\ell}/Y_{\ell})$; this number is equal to 
the determinant of the form $\mu_1,\mu_2\mapsto\frac{1}{\ell}\mu_1
\cdot\mu_2$ on $Y_{\ell}$.

\subsection{} 
\label{nul' dot} Let $\tfu\subset\fu$ be the $\sk$-subalgebra generated by
$\tK_i, \epsilon_i, \theta_i\ (i\in I)$. Following the method of 
\cite{l1} 23.1, define a new algebra $\dfu$ (without unit) as follows. 

If $\mu', \mu''\in X_{\ell}$, we set 
$$
_{\mu'}\tfu_{\mu''}=\tfu/(\sum_{i\in I} (\tK_i-\zeta^{i\cdot\mu'})\tfu+
\sum_{i\in I} \tfu(\tK_i-\zeta^{i\cdot\mu''}));\
\dfu=\oplus_{\mu',\mu''\in X_{\ell}} (_{\mu'}\tfu_{\mu''}).
$$
Let $\pi_{\mu',\mu''}:\tfu\lra\ _{\mu'}\tfu_{\mu''}$ be the canonical
projection. We set $1_{\mu}=\pi_{\mu,\mu}(1)\in\dfu$. 
The structure of an algebra on $\dfu$ is defined as in {\em loc. cit.} 

As in {\em loc. cit.}, the category $\CC$ may be identified with the category
of finite dimensional (over $\sk$) (left) $\dfu$-modules $M$ which are
{\em unital}, i.e.  
 
(a) for every $x\in M$, $\sum_{\mu\in X_{\ell}}\ 1_{\mu}x=x$. 

If $M$ is 
such a module, the $X_{\ell}$-grading on $M$ is defined by $M_{\mu}=1_{\mu}M$. 

Let $\fu'$ denote  the quotient algebra of the algebra $\tfu$ by the 
relations $\tK_i^{l_i}=1\ (i\in I)$. Here $l_i:=\frac{l}{(l,d_i)}$. We have 
an isomorphism of vector spaces $\dfu=\fu'\otimes \sk[Y_{\ell}]$, cf. \ref{nul' adj}
below.   

\subsection{} Let $a:\CC\iso\CC_r$ be an equivalence defined 
by $aM=M\ (M\in\CC)$ as an $X_{\ell}$-vector space, $mx=s(x)m\
(x\in\fu, m\in M)$. Here $s: \fu\lra\fu^{\opp}$ is the antipode.
We will use the same notation $a$ for a similar equivalence $\CC_r\iso\CC$. 

Let us consider the category $\CC\otimes\CC$ (resp. $\CC\otimes\CC_r$) 
which may be identified with   
the category of finite dimensional $\dfu\otimes\dfu$- (resp.
$\dfu\otimes(\dfu)^{\opp}$-) modules satisfying 
a "unitality" condition similar to (a) above. Let us consider the 
algebra $\dfu$ itself as 
a regular $\dfu\otimes(\dfu)^{\opp}$-module. It is infinite dimensional, 
but is a union of finite 
dimensional modules, hence it may be considered as an
object of the category $\Ind(\CC\otimes\CC_r)=\Ind\CC\otimes\Ind\CC_r$ 
where $\Ind$ denotes the category of $\Ind$-objects, cf. \cite{d4} \S 4. 
Let us denote by $\bR$ the image of this object under the equivalence 
$\Id\otimes a:\Ind\CC\otimes\Ind\CC_r\iso\Ind\CC\otimes\Ind\CC$. 

Every object $\CO\in\CC\otimes\CC$ induces a functor $F_{\CO}:\CC\lra\CC$ 
defined by 
$$
F_{\CO}(M)=a(aM\otimes_{\CC}\CO).
$$
The same formula defines a functor $F_{\CO}:\Ind\CC\lra\Ind\CC$ for
$\CO\in\Ind(\CC\otimes\CC)$. 

We have $F_{\bR}=\Id_{\Ind\CC}$. 

We can consider a version of the above formalism using 
semiinfinite $\Tor$'s. An object $\CO\in\Ind(\CC\otimes\CC)$ defines functors
$F_{\CO;\infh+n}:\Ind\CC\lra\Ind\CC\ (n\in\BZ)$ defined by
$$
F_{\CO;\infh+n}(M)=a\Tor^{\CC}_{\infh+n}(aM,\CO).
$$
\subsection{Theorem} {\em {\em (i)} We have $F_{{\em\bR};\infh+n}=\Id_{\Ind\CC}$
if $n=0$, and $0$ otherwise. 

{\em (ii)} Conversely, suppose we have an object $Q\in\Ind(\CC\otimes\CC)$ 
together with an isomorphism of functors 
$\phi: F_{{\em\bR};\infh+\bullet}\iso F_{Q;\infh+\bullet}$. 
Then $\phi$ is induced 
by the unique isomorphism ${\em\bR}\iso Q$.} $\Box$

\subsection{Adjoint representation}
\label{nul' adj} For $\mu\in Y_{\ell}$, let $T(\mu)$
be a one-dimensional
$\dfu\otimes(\dfu)^{\opp}$-module equal to $L(\mu)$ (resp. to $aL(-\mu)$) 
as a $\dfu$- (resp. $(\dfu)^{\opp}$-) module. Let us consider the module
$T_{\mu}\bR=\bR\otimes T(\mu)\in\Ind(\CC\otimes\CC)$. This object 
represents the same functor $\Id_{\Ind\CC}$, hence we have a canonical 
isomorphism $t_{\mu}:\bR\iso T_{\mu}\bR$. 

Let us denote by $\had\in\Ind\CC$ the image of $\bR$ under the tensor 
product $\otimes:\Ind(\CC\otimes\CC)\lra\Ind\CC$. The isomorphisms $t_{\mu}$ 
above induce an action of the lattice $Y_{\ell}$ on $\had$. 
Set $\ad=\had/Y_{\ell}$. This is an object of $\CC\subset\Ind\CC$ which is 
equal to the algebra $\fu'$ considered as a $\dfu$-module by means
of the adjoint action. 

In the notations of \ref{nul' ms}, let us consider an object
$$
\ad_{\MS}:=\oplus_{\mu\in\Delta_{\ell}} \langle L(\mu)\otimes L(\mu)^*
\rangle_{\MS}\in\MS,
$$
cf. \cite{bfm} 4.5.3. 
  
\subsection{Theorem} {\em We have a canonical isomorphism   
$\langle{\em\ad}\rangle_{\MS}={\em\ad}_{\MS}$.} $\Box$

\section{Regular sheaf}
\label{nul' reg sheaf}

\subsection{Degeneration of quadrics} The construction below is taken
from \cite{kl}II 15.2. Let us consider the quadric $Q\subset\BP^1
\times\BP^1\times\BA^1$ given by the equation $uv=t$ where
$(u,v,t)$ are coordinates in the triple product. Let $f:Q\lra\BA^1$
be the projection to the third coordinate; for $t\in\BA^1$ denote
$Q_t:=f^{-1}(t)$. For $t\neq 0$, $Q_t$ is isomorphic to $\BP^1$;
the fiber $Q_0$ is a union of two projective lines clutched at a point:  
$Q_0=Q_u\cup Q_v$ where $Q_u$ (resp. $Q_v$) is an irreducible component  
given (in $Q_0$) by the equation $v=0$ (resp. $u=0$) and is isomorphic 
to $\BP^1$; their intersection being a point. We set
$'Q=f^{-1}(\BA^1-\{0\})$.

We have two sections $x_1, x_2: \BA^1\lra Q$ given by $x_1(t)=
(\infty,0,t),\ x_2(t)=(0,\infty,t)$. Consider two "coordinate charts" 
at these points: the maps $\phi_1,\phi_2:\BP^1\times\BA^1\lra Q$
given by
$$
\phi_1(z,t)=(\frac{tz}{z-1},\frac{z-1}{z},t);\
\phi_2(z,t)=(\frac{z-1}{z},\frac{tz}{z-1},t).
$$
This defines a map 

(a) $\phi: \BA^1-\{0\}\lra\tP(2)$,

in the notations of \ref{nul' right}.

\subsection{} For $\nu\in Y^+$, let us consider the corresponding 
(relative over $\BA^1$) configuration scheme
$f^{\nu}: Q^{\nu}_{/\BA^1}\lra\BA^1$.
For the brevity we will omit the subscript $_{/\BA^1}$ indicating that
we are dealing with the relative version of configuration spaces. We denote 
by $Q^{\nu\bullet}$ (resp. $Q^{\nu\circ}$) the subspace of configurations 
with the points distinct from $x_1, x_2$ (resp. also pairwise distinct).
We set $'Q^{\nu\circ}=Q^{\nu\circ}|_{\BA^1-\{0\}}$, etc.

The map $\phi$ above, composed with the canonical projection 
$\tP(2)\lra P(2)$, induces the maps 
$$
'Q^{\nu\circ}\lra P^{\nu}(2)\ (\nu\in Y^+)
$$
(in the notations of \ref{nul' diag}). For $\nu\in Y^+$ and
$\mu_1,\mu_2\in X_{\ell}$ such that $\mu_1+\mu_2-\nu= 2\rho_{\ell}$, 
let $\CJ_{\mu_1,\mu_2}^{\nu}$ denote the local
system over $'Q^{\nu\circ}$ which is the inverse image of the local
system $\CH^{\nu}_{\mu_1,\mu_2}$ over $P^{\nu}(2)^{\circ}$. 
Let $\CJ_{\mu_1,\mu_2}^{\nu\bullet}$ denote the perverse sheaf
over $'Q^{\nu\bullet}(\BC)$ which is the middle extension of
$\CJ_{\mu_1,\mu_2}^{\nu}[\dim Q^{\nu}]$.

\subsection{}
\label{nul' nearby}
Let us take the nearby cycles and get a perverse sheaf $\Psi_{f^{\nu}}
(\CJ^{\nu\bullet}_{\mu_1,\mu_2})$ over $Q_0^{\nu\bullet}(\BC)$. Let us consider the
space $Q_0^{\nu\bullet}$ more attentively. This is a reducible scheme 
which is a union
$$
Q^{\nu\bullet}_0=\bigcup_{\nu_1+\nu_2=\nu}\ \BA^{\nu_1}\times\BA^{\nu_2},
$$
the component $\BA^{\nu_1}\times\BA^{\nu_2}$ corresponding to configurations
where $\nu_1$ (resp. $\nu_2$) points are running on the affine line 
$Q_u-x_1(0)$ (resp. $Q_v-x_2(0)$). Here we identify these affine lines 
with a "standard" one using the coordinates $u$ and $v$ respectively.    
Using this decomposition we can define a closed embedding  
$$
i_{\nu}: Q^{\nu\bullet}_0\hra\BA^{\nu}\times\BA^{\nu}
$$
whose restriction to a component $\BA^{\nu_1}\times\BA^{\nu_2}$ sends
a configuration as above, to the configuration where all 
remaining points are equal to zero. 
Let us define a perverse sheaf 
$$
\CR_{\mu_1,\mu_2}^{\nu,\nu}=i_{\nu*}\Psi_{f^{\nu}}(\CJ^{\nu\bullet}
_{\mu_2,\mu_1})
\in\CM(\BA^{\nu}(\BC)\times\BA^{\nu}(\BC);\CS)
$$
Let us consider the collection of sheaves $\{\CR_{\mu_1,\mu_2}^{\nu,\nu}|\ 
\mu_1,\mu_2\in X_{\ell}, \nu\in Y^+, \mu_1+\mu_2-\nu=2\rho_{\ell}\}$.
One can complete this collection to an object $\CR$ of the category 
$\Ind(\FS\otimes\FS)$ where $\FS\otimes\FS$ is understood as
a category of finite factorizable sheaves corresponding to the 
{\em square} of our initial Cartan datum, i.e. $I\coprod I$, etc. 
For a precise construction, see Part V.

\subsection{Theorem} {\em We have $\Phi(\CR)={\em\bR}$.} $\Box$
 
\newpage
\begin{center}
{\bf Chapter 3. Modular}
\end{center}


\section{Heisenberg local systems}

In this section we sketch a construction of certain 
remarkable cohesive local systems on arbitrary smooth families of  
compact smooth curves,
to be called the {\em Heisenberg local systems}. 

In the definition and construction of local sustems below we will 
have to assume that our base field $\sk$ contains roots of unity
of sufficiently high degree; the characteristic of $\sk$ is assumed to be
prime to this degree.           

\subsection{} From now on until \ref{nul' last heis} we fix a smooth
proper morphism $f: C\lra S$ of relative dimension $1$,   
$S$ being a smooth connected scheme over $\BC$. For $s\in S$, we denote
$C_s:=f^{-1}(s)$. Let $g$ be the genus of fibres of $f$.

Let $S_{\delta}$ denote the total space of the determinant line bundle
$\delta_{C/S}=\det Rf_*\Omega^1_{C/S}$ without the zero section.
For any object (?) over $S$ (e.g., a scheme over $S$, a sheaf 
over a scheme over $S$, etc.), we will denote by (?)$_{\delta}$ its
base change under $S_{\delta}\lra S$.

Below, if we speak about a scheme as a topological (analytic)
space, we mean its set of $\BC$-points with the usual topology
(resp. analytic structure). 

\subsection{} 
\label{nul' rel conf}
We will use the relative versions of configuration spaces; to indicate 
this, we will use the subscript $_{/S}$.    
Thus, if $J$ is a finite set, $C^J_{/S}$ will denote the $J$-fold fibered 
product of $C$ with itself over $S$, etc. 

Let $C(J)_{/S}$ denote the subscheme of the $J$-fold cartesian power 
of the relative tangent bundle $T_{C/S}$ consisting of  
$J$-tuples $\{x_j, \tau_j\}$ where 
$x_j\in C$ and $\tau_j\neq 0$ in $T_{C/S,x}$, the points $x_j$ being pairwise
distinct. 
Let $\tC(J)_{/S}$ denote the 
space of $J$-tuples of holomorphic embeddings $\phi_j: D\times S\lra C$ 
over $S$ with disjoint images; we have the $1$-jet maps 
$\tC(J)_{/S}\lra C(J)_{/S}$. 

An epimorphism $\rho: K\lra J$ induces the maps
$$
m_{C/S}(\rho):\ \prod_J \tD(K_j)\times\tC(J)_{/S}\lra \tC(J)_{/S}
$$
and
$$
\bar{m}_{C/S}(\rho): \prod_J D(K_j)\times\tC(J)_{/S}\lra C(K)_{/S}
$$
which satisfy the compatibilities as in \ref{nul' right}.

\subsection{} We extend the function $n$ to $X_{\ell}$ (see \ref{nul' modif}) by
$n(\mu)=\frac{1}{2}\mu\cdot\mu-
\mu\cdot\rho_{\ell}\ (\mu\in X_{\ell})$. We will denote by $\CJ$ the 
$X_{\ell}$-coloured local system over the operad of disks $\CCD$ which
is defined exactly as in \ref{nul' stand d}, with $X$ replaced by $X_{\ell}$.

\subsection{} A {\em cohesive local system $\CH$ of level $\mu\in X_{\ell}$
over $C/S$} 
is a collection of local systems $\CH(\pi)$ over the spaces 
$C(J)_{/S;\delta}$
given for every map $\pi: J\lra X_{\ell}$ of level $\mu$ (note the base change  
to $S_{\delta}$!), together with the factorization isomorphisms
$$
\phi_C(\rho): m_{C/S}(\rho)^*\CH(\pi)\iso\Boxtimes_J \CJ(\pi_j)
\boxtimes \CH(\rho_*\pi).
$$
Here we have denoted by the same letter $\CH(\pi)$ the lifting of $\CH(\pi)$
to $\tC_{/S;\delta}$.
The factorization isomorphisms must satisfy the obvious analogs
of properties \ref{nul' cls def} (a), (b).

\subsection{} 
\label{nul' start heis} Now we will sketch a construction of certain cohesive
local system over $C/S$ of level $(2-2g)\rho_{\ell}$. For alternative 
beautiful constructions of $\CH$, see \cite{bp}.   

To simplify the exposition we will assume below that $g\geq 2$ 
(the construction for $g\leq 1$ needs some modification, and we omit it here,
see Part V).
Let us consider the group scheme $\Pic(C/S)\otimes X_{\ell}$ over $S$. 
Here $\Pic(C/S)$ is the relative Picard scheme. The group 
of connected components $\pi_0(\Pic(C/S)\otimes X_{\ell})$ is equal 
to $X_{\ell}$. Let us denote by $\Jac$  
the connected component corresponding to the element $(2-2g)\rho_{\ell}$; 
this is an
abelian scheme over $S$, due to the existence of the section 
$S\lra\Jac$ defined by $\Omega^1_{C/S}\otimes(-\rho_{\ell})$.  

For a scheme $S'$ over $S$, let $H_1(S'/S)$ denote the local system
of the first relative integral homology groups over $S$. We have 
$H_1(\Jac/S)=H_1(C/S)\otimes X_{\ell}$. We will denote by $\omega$ the
polarization of $\Jac$ (i.e. the  
skew symmetric form on the latter local system) equal to the tensor  
product of the standard form on $H_1(C/S)$ and the form 
$(\mu_1,\mu_2)\mapsto\frac{\dd^g_{\ell}}{\ell}\mu_1\cdot\mu_2$ on $X_{\ell}$. 
Note that the assumption $g\geq 2$ implies that 
$\frac{\dd^g_{\ell}}{\ell}\mu_1\cdot\mu_2\in\BZ$ for any $\mu_1,\mu_2\in
X_{\ell}$.   
Since the latter form is positive definite, $\omega$ is relatively ample 
(i.e. defines a relatively ample invertible sheaf on $\Jac$).

\subsection{} Let $\alpha=\sum n_{\mu}\cdot\mu\in\BN[X_{\ell}]$; set
$\Supp(\alpha)=\{\mu|\ n_{\mu}\neq 0\}$. Let us say that $\alpha$ 
is {\em admissible} if $\sum n_{\mu}\mu=(2-2g)\rho_{\ell}$. 
Let us denote by 
$$
\ya_{\alpha}: C^{\alpha}_{/S}\lra\Pic(C/S)\otimes X_{\ell}
$$
the Abel-Jacobi map sending $\sum \mu\cdot x_{\mu}$ to
$\sum x_{\mu}\otimes\mu$. If $\alpha$ is admissible then the map 
$\ya_{\alpha}$ lands in $\Jac$.

Let $D^{\alpha}$ denote the following relative divisor on $C^{\alpha}_{/S}$
$$
D^{\alpha}=\frac{\dd_{\ell}^g}{\ell}(\sum_{\mu\neq\nu}\mu\cdot\nu
\Delta_{\mu\nu}+\frac{1}{2}\sum_{\mu}\mu\cdot\mu\Delta_{\mu\mu}).
$$
Here $\Delta_{\mu\nu}\ (\mu,\nu\in \Supp(\alpha))$ denotes  
the corresponding diagonal in $C^{\alpha}_{/S}$. Note that all the 
multiplicities are integers.  

Let $\pi: J\lra X_{\ell}$ be an unfolding of $\alpha$. We will denote by
$D^{\pi}$ the pull-back of $D^{\alpha}$ to $C^J_{/S}$. Let us 
introduce the following line bundles  
$$
\CL(\pi)=\otimes_{j\in J}\CT_j^{\otimes\frac{\dd^g_{\ell}}{\ell}n(\pi(j))}
\otimes\CO(D^{\pi})
$$
on $C^J_{/S}$, and
$$
\CL_{\alpha}=\CL(\pi)/\Sigma_{\pi}
$$
on $C^{\alpha}_{/S}$ (the action of $\Sigma_{\pi}$ is an obvious one).
Here $\CT_j$ denotes the relative tangent line bundle on $C^J_{/S}$ in the 
direction $j$. Note that the numbers $\frac{\dd^g_{\ell}}{\ell}n(\mu)\
(\mu\in X_{\ell})$ are integers.  

\subsection{Proposition} {\em There exists a unique line bundle 
$\CL$ on $\Jac$ such that for each admissible $\alpha$, we have 
$\CL_{\alpha}=\ya_{\alpha}^*(\CL)$. The first Chern class 
$c_1(\CL)=-[\omega]$.} $\Box$

\subsection{} In the sequel if $\CL_0$ is a line bundle, let 
$\dCL_0$ denote the total space its with the zero section removed.
 
The next step is the construction of a certain local system $\fH$ over 
$\dCL_{\delta}$. Its dimension is equal to $\dd^g_{\ell}$ and
the monodromy 
around the zero section of $\CL$ (resp. of the determinant bundle) is equal to 
$\zeta^{-2\ell/\dd_{\ell}^g}$ (resp. $(-1)^{\rk(X)}\zeta^{-12\rho_{\ell}\cdot
\rho_{\ell}})$. The construction of $\fH$ is outlined below.

The previous construction assigns to a triple 

(a lattice $\Lambda$, a symmetric bilinear form $(\ ,\ ):\ \Lambda\times\Lambda
\lra\BZ$, $\nu\in\Lambda$)

an abelian scheme $\Jac_{\Lambda}:=(\Pic(C/S)\otimes\Lambda)_{(2g-2)\nu}$ 
over $S$, together with a line bundle $\CL_{\Lambda}$ on it (in the 
definition of $\CL_{\Lambda}$ one should use the function 
$n_{\nu}(\mu)=\frac{1}{2}\mu\cdot\mu+\mu\cdot\nu$). We considered the case 
$\Lambda=X_{\ell},\ (\mu_1,\mu_2)=\frac{\dd^g_{\ell}}{\ell}\mu_1\cdot\mu_2,\ 
\nu=-\rho_{\ell}$.  

Now let us apply this construction to the lattice 
$\Lambda=X_{\ell}\oplus Y_{\ell}$, the bilinear form
$((\mu_1,\nu_1),(\mu_2,\nu_2))=-\frac{1}{\ell}(\nu_1\cdot\mu_2+\nu_2\cdot\mu_1
+\nu_1\cdot\nu_2)$ and $\nu=(-\rho_{\ell},0)$. The first projection 
$\Lambda\lra X_{\ell}$ induces the morphism 
$$
p:\ \Jac_{\Lambda}\lra\Jac
$$
the fibers of $p$ are abelian varieties $\Jac(C_s)\otimes Y_{\ell}\ (s\in S)$.

\subsection{Theorem} {\em {\em (i)} The line bundle $\CL_{\Lambda}$ 
is relatively ample with respect to $p$. The direct image 
$\CE:=p_*\CL_{\Lambda}$ is  a locally free sheaf of rank $\dd_{\ell}^g$.

{\em (ii)} We have an isomorphism
$$
\det(\CE)=\CL\otimes\delta^{\dd^g_{\ell}(-\frac{1}{2}\rk(X_{\ell})+
6\frac{\rho_{\ell}\cdot\rho_{\ell}}{\ell})}.
$$} $\Box$

Here $\delta$ denotes the pull-back of the determinant bundle
$\delta_{C/S}$ to $\Jac$.

\subsection{} Let us assume for a moment that $\sk=\BC$ and
$\zeta=\exp(-\frac{\pi\sqrt{-1}}{\ell})$.   
By the result of Beilinson-Kazhdan, \cite{bk} 4.2,
the vector bundle $\CE$ carries a canonical flat projective connection. 
By {\em loc. cit.} 2.5, its lifting to $\det(\CE)^{\bullet}$ carries a flat
connection with the scalar monodromy around the zero section equal to 
$\exp(\frac{2\pi\sqrt{-1}}{\dd^g_{\ell}})$. We have an 
obvious map
$$
m: \dCL_{\delta}\lra\CL\otimes\delta^
{\dd^g_{\ell}(-\frac{1}{2}\rk(X_{\ell})+
6\frac{\rho_{\ell}\cdot\rho_{\ell}}{\ell})}.
$$
By definition, $\fH$ is the local system of horizontal sections
of the pull-back of $\CE$ to $\dCL_{\delta}$. The claim about its monodromies
follows from part (ii) of the previous theorem.

This completes the construction of $\fH$ for $\sk=\BC$ and $\zeta=\exp
(-\frac{\pi\sqrt{-1}}{\ell})$. The case of arbitrary $\sk$ (of sufficiently
large characteristic) and $\zeta$ follows from this one.       

\subsection{} 
\label{nul' last heis}
Let us consider an obvious map $q: C(J)_{/S;\delta}\lra C^J_{/S;\delta}$.
The pull-back  
$q^*\CL(\pi)$ has a canonical non-zero section $s$. Let $\tCH(\pi)$ be 
the pull-back of the local system $\fH$ to $q^*\CL(\pi)$. By definition,
we set $\CH(\pi)=s^*\tCH(\pi)$. For the construction of the factorization 
isomorphisms, see Part V.

\subsection{} 
\label{nul' level} Let $\fg$ be the simple Lie algebra connected with
our Cartan datum.   
Assume that $\zeta=\exp(\frac{\pi\sqrt{-1}}{d\kappa})$ for some
positive integer $\kappa$, cf. \ref{nul' charge} ($d$ is defined in \ref{nul' root}).

We have $12\rho_{\ell}\cdot\rho_{\ell}\equiv 12\rho\cdot\rho\ \modul\ l$, and 
$\rk(X)\equiv\dim\fg\ \modul\ 2$. By the strange formula of
Freudenthal-de Vries, we have $12\rho\cdot\rho=dh\dim\fg$ where $h$ 
is the dual Coxeter number of our Cartan datum.   
It follows that the monodromy of 
$\CH$ around the zero section of the determinant line bundle is equal to 
$\exp(\pi\sqrt{-1}\frac{(\kappa-h)\dim\fg}{\kappa})$. This number coincides 
with the multiplicative central charge of the conformal 
field theory associated with the affine Lie algebra $\hfg$ at level $\kappa$
(see \cite{bfm} 4.4.1, 6.1.1, 2.1.3, \cite{tuy} 1.2.2), cf 
\ref{nul' subquot} below.

{\em UNIVERSAL HEISENBERG SYSTEMS}

\subsection{} 
Let us define a category $\Sew$ as follows (cf. \cite{bfm} 4.3.2).
Its object $A$ is a finite set $\bA$ together with a collection $N_A=
\{n\}$ of non-intersecting two-element subsets $n\subset\bA$. Given
such an object, we set $A^1=\bigcup_{N_A}\ n,\ A^0=\bA-A^1$. A morphism
$f: A\lra B$ is an embedding $i_f: \baB\hra \bA$ and a collection $N_f$ of
non-intersecting two-element subsets of $\baB-\bA$ such that
$N_A=N_B\coprod N_f$. The composition of morphisms is obvious. ($\Sew$
coincides with the category $Sets^{\sharp}/\emp$, in the notations of 
\cite{bfm} 4.3.2.)  

For $A\in\Sew$, let us call an {\em $A$-curve} a data $(C, \{x_a,\tau_a\}
_{A^0})$ where $C$ is a smooth proper (possibly disconnected) complex  
curve, $\{x_a,\tau_a\}_{A^0}$ is an $A^0$-tuple of distinct points 
$x_a\in C$ together with non-zero tangent vectors $\tau_a$ at them.  
For such a curve, let $\bC_A$ denote the curve obtained from $C$ by clutching  
pairwise the points $x_{a'}$ with $x_{a''}\ (n=\{a',a''\})$
for all sets $n\in N_A$. Thus, the set $N_A$ is in the bijection with 
the set of nodes of the curve $\bC_A$ ($\bC_A=C$ if $N_A=\emp$). 

Let us call an {\em enhanced graph} a pair $\Gamma=(\bGamma,\bg)$ here $\bGamma$ 
is a non-oriented graph and $\bg=\{g_v\}_{v\in\Ve(\bGamma)}$ is a $\BN$-
valued $0$-chain of $\bGamma$. Here $\Ve(\bGamma)$ denotes the set of 
vertices of $\bGamma$. Let us assign to a curve $\bC_A$ an enhanced graph
$\Gamma(\bC_A)=(\bGamma(\bC_A), \bg(\bC_A))$. By definition,
$\bGamma(\bC_A)$ is a graph with $\Ve(\bGamma(\bC_A))=\pi_0(C)=\{$ 
the set of irreducible components of $\bC_A\}$ and the set of edges   
$\Ed(\bGamma(\bC_A))=N_A$, an edge $n=\{a',a''\}$ connecting the
vertices corresponding 
to the components of the points $x_{a'}, x_{a''}$. For $v\in\pi_0(C)$,
$g(\bC_A)_v$ is equal to the genus of the corresponding component 
$C_v\subset C$.   

\subsection{}   
Let $\CM_A$ denote the moduli stack of $A$-curves $(C,\ldots)$ such that 
the curve $\bC_A$ is stable in the sense of \cite{dm} (in particular
connected). The stack $\CM_A$ is smooth; we have   
$\CM_A=\coprod_{g\geq 0}\CM_{A,g}$ where $\CM_{A,g}$ is a  
substack of $A$-curves $C$ with $\bC_A$ having genus (i.e. 
$\dim H^1(\bC_A,\CO_{\bC_A})$) equal to $g$. In turn, we have the decomposition 
into connected components 
$$
\CM_{A,g}=\coprod_{\Gamma,\vA^0}\ \CM_{\vA^0,g,\Gamma}
$$ 
where $\CM_{\vA^0,g,\Gamma}$ is the stack of $A$-curves $(C,\{x_a,\tau_a\})$ as 
above, with $\Gamma(\bC_A)=\Gamma$, $\vA^0=\{A^0_v\}_{v\in\Ve(\bGamma)}$,
$A^0=\coprod\ A^0_v$, such that $x_a$ lives on the connected component 
$C_v$ for $a\in A^0_v$.   
  
We denote by $\eta: C_{A,g}\lra\CM_{A,g}$
(resp. $\baeta: \bC_{A,g}\lra\CM_{A,g})$ the universal smooth curve 
(resp. stable surve).
For $\nu\in Y^+$, we have the corresponding {\em relative} configuration 
spaces $C^{\nu}_{A,g}, C^{\nu\circ}_{A,g}, \bC^{\nu}_{A,g}$. For brevity we 
omit the relativeness subscript $/\CM_{A,g}$ from these notations. 
The notation $C^{\nu}_{\vA^0,g,\Gamma}$ etc., will mean the restriction 
of these 
configuration spaces to the component $\CM_{\vA^0,g,\Gamma}$. 

Let $\CM_g$ be the moduli stack of smooth connected curves of genus $g$, 
and $\bCM_g$ be its Grothendieck-Deligne-Mumford-Knudsen  
compactification, i.e. the moduli stack of stable curves of genus $g$. 
Let $\baeta: \bC_g\lra\bCM_g$ be the universal stable curve; let 
$\ol{\delta_g}=\det R\baeta_*(\omega_{\bC_g/\bCM_g})$ be the determinant line 
bundle; let $\bCM_{g;\delta}\lra\bCM_g$ be its total space with the zero
section removed.      

We have obvious maps $\CM_{A,g}\lra\bCM_g$.
Let the complementary subscript $(\cdot)_{\delta}$ denote the base change of
all the above objects under $\bCM_{g;\delta}\lra\bCM_g$.

\subsection{}
\label{nul' def h} Let us consider the configuration space $C^{\nu\circ}_{A,g}$;
it is the moduli stack of $\nu$ distinct points running on $A$-curves
$(C,\{x_a,\tau_a\})$ and not equal to the marked points $x_a$. This stack 
decomposes into connected components as follows: 
$$
C^{\nu\circ}_{A,g}=\coprod_{\Gamma,\vA^0,\vnu}\ C^{\vnu}_{\vA^0,g,\Gamma}
$$
where $\vnu=\{\nu_v\}_{v\in\Ve(\Gamma)}$ and $C^{\vnu}_{\vA^0,g,\Gamma}$ being 
the moduli stack of objects as above, with $\Gamma(\bC_A)=\Gamma$ and 
$\nu_v$ points running on the component $C_v$. The decomposition is taken 
over appropriate graphs $\Gamma$, decompositions $A^0=\coprod A^0_v$ and 
the tuples $\vnu$ with $\sum \nu_v=\nu$. 

Let us call an $A_0$-tuple $\vmu=\{\mu_a\}\in X_{\ell}^{A^0}$ 
{\em $(g,\nu)$-good} if 

(a) $\sum_{a\in A^0}\ \mu_a-\nu\equiv (2-2g)\rho_{\ell}\ \modul\ Y_{\ell}$.

Given such a tuple, we are going to define certain local system 
$\CH^{\nu}_{\vmu;A,g}$ over $C^{\nu\circ}_{A,g;\delta}$. Let us describe its
restriction $\CH^{\vnu}_{\vmu;\vA^0,g,\Gamma}$ to a connected 
component $C^{\vnu}_{\vA^0,g,\Gamma;\delta}$.

Let $\Gamma'$ be the first subdivision of $\bGamma$. We have
$\Ve(\Gamma')=\Ve(\bGamma)\coprod\Ed(\bGamma)=\pi_0(C)\coprod N_A$. 
The edges of $\Gamma'$ are 
indexed by the pairs $(n,a)$ where $n\in N_A, a\in n$, the corresponding 
edge $e_{n,a}$ having the ends $a$ and $n$. Let us define an orientation 
of $\Gamma'$ by the requierement that $a$ is the beginning of $e_{n,a}$. 
Consider the chain complex 
$$
C_1(\Gamma'; X_{\ell}/Y_{\ell})\overset{d}{\lra}C_0(\Gamma';X_{\ell}/Y_{\ell}). 
$$
Let us define a $0$-chain $c=c^{\vnu}_{\vmu}\in C_0(\Gamma';X_{\ell}/Y_{\ell})$ 
by 
$$
c(v)=\sum_{a\in A^0_v}\ \mu_a+(2g_v-2)\rho_{\ell}-\nu_v\ (v\in\pi_0(C));\
c(n)=2\rho_{\ell}\ (n\in N_A). 
$$
The goodness assumption (a) ensures that $c$ is a boundary.
By definition,
$$
\CH^{\vnu}_{\vmu;\vA^0,g,\Gamma}=\oplus_{\chi:\ d\chi=c}\ \CH_{\chi}
$$
Note that the set $\{\chi|\ d\chi=c\}$ is a torsor over the group 
$H_1(\Gamma';X_{\ell}/Y_{\ell})=H_1(\Gamma;X_{\ell}/Y_{\ell})$.
 
The local system $\CH_{\chi}$ is defined below, in \ref{nul' h chi}, after a little
notational preparation. 

\subsection{} Given two finite 
sets $J, K$, let $C(J;K)_g$ denote the moduli stack of objects 
$(C,\{x_j\},\{y_k,\tau_k\})$. Here $C$ is a smooth proper connected
curve of genus $g$, $\{x_j\}$ is a $J$-tuple of distinct points $x_j\in C$ 
and $\{y_k,\tau_k\}$ is a $K$-tuple of distinct points $y_k\in C$ together
with non-zero tangent vectors $\tau_k\in T_{y_k}C$. We suppose that
$y_k\neq x_j$ for all $k,j$.

We set $C(J)_g:=C(\emp;J)_g$. We have the forgetful maps 
$C(J\coprod K)_g\lra C(J;K)_g$. 

The construction of \ref{nul' start heis} --- \ref{nul' last heis} defines
the Heisenberg system $\CH(\pi)$ over the smooth stack  
$C(J)_{g;\delta}$ for each $\pi: J\lra X_{\ell}$.

Given $\nu\in Y^+$, choose an unfolding of $\nu$, $\pi: J\lra I$, and set
$C^{\nu}(K)^{\circ}_g:=C(J;K)/\Sigma_{\pi}$. Given a $K$-tuple $\vmu=
\{\mu_k\}\in X_{\ell}^K$, define a map $\tpi: J\coprod K\lra X_{\ell}$ by
$\tpi(j)=-\pi(j)\in - I\subset X_{\ell}\ (j\in J),\
\tpi(k)=\mu_k\ (k\in K)$. The local system $\CH(\tpi)$ over
$C(J\coprod K)_{g;\delta}$ descends to $C(J;K)_{g;\delta}$ since
$\zeta^{2n(-i)}=1$, and then to $C^{\nu}(K)_{g;\delta}^{\circ}$, by
$\Sigma_{\pi}$-equivariance; denote the latter local system by 
$\tCH_{\vmu}^{\nu}$, and set $\CH^{\nu}_{\vmu}=\tCH^{\nu}_{\vmu}\otimes 
\Sign^{\nu}$, cf. \ref{nul' sign p}.

\subsection{Lemma} {\em If  
$\vmu\equiv\vmu'\ \modul\ Y_{\ell}^K$ then we   
have canonical isomorphisms $\CH^{\nu}_{\vmu}=\CH^{\nu}_{\vmu'}$.} $\Box$

Therefore, it makes sense to speak about $\CH^{\nu}_{\vmu}$ for 
$\vmu\in (X_{\ell}/Y_{\ell})^K$. 

\subsection{}
\label{nul' h chi} Let us return to the situation at the end of \ref{nul' def h}.
We have $\Gamma=(\bGamma,\{g_v\}_{v\in\Ve(\bGamma)})$. Recall that 
$A^0=\coprod\ A^0_v$; we have also $A^1=\coprod A^1_v$ where $A^1_v:= 
\{a\in A^1|\ x_a\in C_v\}\ (v\in \Ve(\bGamma)$. Set $\bA_v=A^0_v\coprod A^1_v$,
so that $\bA=\coprod \bA_v$. We have an obvious map 

(a) $\CM_{\vA^0,g,\Gamma}\lra\prod_{\Ve(\bGamma)}\ C(\bA_v)_{g_v}$, 

and a map 

(b) $C^{\vnu}_{\vA^0,g,\Gamma}\lra \prod_{\Ve(\bGamma)}\ 
C^{\nu_v}(\bA_v)^{\circ}_{g_v}$ 

over (a). For each $v$, define an $\bA_v$-tuple $\vmu(\chi;v)$ equal to $\mu_a$ 
at $a\in A^0_v$ and $\chi(e_{a,n})$ at $a\in n\subset A^1$. By definition,
the local system $\CH_{\chi}$ over $C^{\vnu}_{\vA^0,g,\Gamma;\delta}$ is the
inverse image of the product 

$\Boxtimes_{\Ve(\bGamma)}\ \CH^{\nu_v}_{\vmu(\chi;v)}$

under the map (b) (pulled back to the determinant bundle).  

This completes the definition of the local systems $\CH^{\nu}_{\vmu;A,g}$. 
They have a remarkable compatibility property (when the object $A$ varies)  
which we are going to describe below, see theorem \ref{nul' cart}.

\subsection{} 
\label{nul' conf fib} Let $\tCT^{\nu}_{A,g}$ denote the fundamental groupoid
$\pi(C^{\nu\circ}_{A,g;\delta})$. We are going to show that these groupoids
form a cofibered category over $\Sew$. 

\subsection{} 
\label{nul' cofib} A morphism $f: A\lra B$ in $\Sew$ is called
a {\em sewing} (resp. {\em deleting}) if $A^0=B^0$ (resp.
$N_f=\emp$). A sewing $f$ with $\card(N_f)=1$ is called {\em simple}.  
Each morphism is a composition of a sewing and a deleting; each sewing 
is a composition of simple ones.  

(a) Let $f: A\lra B$ be a simple sewing. We have canonical morphisms
$$
\wp_f: \CM_{A,g;\delta}\lra\dT_{\dpar\bCM_{B,g;\delta}}\bCM_{B,g;\delta}
$$
and 
$$
\wp_f^{(\nu\circ)}:C^{\nu\circ}_{A,g;\delta}\lra
\dT_{\dpar\overline{C^{\nu\circ}_{B,g;\delta}}}\ \overline
{C^{\nu\circ}_{B,g;\delta}},
$$
over $\wp_f$, cf. \cite{bfm} 4.3.1. Here $\bCM_{B,g}$
(resp. $\overline{C^{\nu\circ}_{B,g}}$)
denotes the Grothendieck-Deligne-Mumford-Knudsen 
compactification of $\CM_{B,g}$ (resp. of $C^{\nu\circ}_{B,g}$)
and $\dpar\bCM_{B,g}$ (resp. $\dpar\overline{C^{\nu\circ}_{B,g}}$)
denotes the smooth locus of the boundary $\bCM_{B,g}-\CM_{B,g}$ (resp.
$\overline{C^{\nu\circ}_{B,g}}-\bC^{\nu}_{B,g}$).
The subscript $_{\delta}$ indicates the base change to the
determinant bundle, as before. 

Composing the specialization with the inverse image under $\wp_f^{(\nu\circ)}$, 
we get the canonical map $f_*:\tT^{\nu}_{A,g}\lra\tT^{\nu}_{B,g}$.

(b) Let $f:A\lra B$ be a deleting. It induces
the obvious morphisms (denoted by the same letter)
$$
f:\CM_{A,g;\delta}\lra\CM_{B,g;\delta}
$$
and
$$
f:C^{\nu\circ}_{A,g;\delta}\lra C^{\nu\circ}_{B,g;\delta}.
$$
The last map induces $f_*:\tT^{\nu}_{A,g}\lra\tT^{\nu}_{B,g}$.

Combining the  constructions (a) and (b) above, we get a category 
$\tT^{\nu}_g$  
cofibered in groupoids over $\Sew$, with fibers $\tT^{\nu}_{A,g}$. 

\subsection{}
\label{nul' rep}
Let $\Rep^{\nu}_{c;A,g}$ be the category of finite dimensional representations 
of $\tT^{\nu}_{A,g}$ (over $\sk$) with the monodromy $c\in \sk^*$ around
the zero section of the determinant bundle. The previous construction 
shows that these categories form a fibered category $\Rep^{\nu}_{c;g}$ over 
$\Sew$. 

\subsection{} For $A\in\Sew$, let us call an $A^0$-tuple 
$\vmu=\{\mu_a\}\in X^{A^0}_{\ell}$  
{\em good} if $\sum_{A^0}\ \mu_a\in Y$. 

If $f: B\lra A$ is a morphism, define a $B^0$-tuple $f^*\vmu=\{\mu'_b\}$ by
$\mu'_b=\mu_{i_f^{-1}(b)}$ if $b\in i_f(A^0)$, and $0$ otherwise. 
Obviously, $\sum_{B^0}\ \mu'_b=\sum_{A^0}\ \mu_a$.

Given a good $\vmu$, let us pick an element $\nu\in Y^+$ such that 
$\nu\equiv\sum\mu_a+(2g-2)\rho_{\ell}\ \modul\ Y_{\ell}$. We
can consider the local system $\CH^{\nu}_{\vmu;A,g}$   
as an object of $\Rep^{\nu}_{c;A,g}$ where 
$c=(-1)^{\card(I)}\zeta^{-12\rho\cdot\rho}$. 

\subsection{Theorem} 
\label{nul' cart} {\em For any morphism $f: B\lra A$ in $\Sew$ and
a $g$-good $\vmu\in X^{A^0}_{\ell}$, we have the canonical isomorphism 
$$
f^*\CH_{\vmu;A,g}^{\nu}=\CH_{f^*\vmu;B,g}^{\nu}.
$$
In other words, the local systems $\CH_{\vmu;B,g}^{\nu}$ define a
{\em{\bf cartesian section}} 
of the fibered category $\Rep^{\nu}_{c;g}$ over $\Sew/A$. Here 
$c=(-1)^{\card(I)}\zeta^{-12\rho\cdot\rho}$. $\Box$}

\section{Fusion structures on $\FS$}
\label{nul' fusion}

\subsection{} Below we will construct a family of "fusion structures"
on the category $\FS$ (and hence, due to the equivalence 
$\Phi$, on the category $\CC$) indexed by $m\in\BZ$. We should explain
what a fusion structure is. This is done in \ref{nul' def fus} below.
We will use a modification of the formalism from \cite{bfm} 4.5.4.  

\subsection{}  
Recall that we have a regular object $\CR\in\Ind
(\FS^{\otimes 2})$, cf. Section \ref{nul' reg sheaf}. We have the canonical
isomorphism $t(\CR)=\CR$ 
where $t: \Ind(\FS^{\otimes 2})\iso\Ind(\FS^{\otimes 2})$ is the permutation, 
hence an object $\CR_n\in\Ind(\FS^{\otimes 2})$ is well defined for 
any two-element set $n$. 

For an object $A\in\Sew$, we set $\tA=A^0\coprod N_A$. 
Let us call an {\em $A$-collection} of factorizable sheaves an 
$\tA$-tuple $\{\CX_{\ta}\}_{\ta\in\tA}$ where $\CX_{\ta}\in\FS_{c_{\ta}}$ if 
$\ta\in A^0$ and
$\CX_{\ta}=\CR_{\ta}$ if $\ta\in N_A$. We impose the condition that 
$\sum_{a\in A^0}\ c_a=0\in X_{\ell}/Y$.
We will denote such an object  
$\{\CX_a;\CR_n\}_A$. It defines an object 
$$
\otimes_A\ \{\CX_a; \CR_n\}\ :=(\otimes_{a\in A^0}\ \CX_a)\otimes(\otimes
_{n\in N_A}\ \CR_n)\in\Ind(\FS^{\otimes\bA}).
$$
If $f: B\lra A$ is a morphism in $\Sew$, we define a $B$-collection
$f^*\{\CX_a;\CR_n\}_B=\{\CY_{\tb}\}_{\tB}$ by $\CY_{\tb}=\CX_{i_f^{-1}(\tb)}$
for $\tb\in i_f(A^0)$, $\One$ if $\tb\in B^0-i_f(A^0)$ and $\CR_{\tb}$ if
$\tb\in N_B=N_A\cup N_f$.

\subsection{}    
Given an $A$-collection $\{\CX_a,\CR_n\}_A$ with $\lambda(\CX_a)=\mu_a$
such that 

(a) $\nu:=\sum_{A^0}\ \mu_a+(2g-2)\rho_{\ell}\in Y^+$,

one constructs (following the pattern of \ref{nul' glu thm}) a perverse sheaf
$\Boxtimes_{A,g}^{(\nu)}\ \{\CX_a; \CR_n\}$ over
$\bC^{\nu}_{A,g;\delta}$.
It is obtained by planting factorizable sheaves $\CX_a$ into the 
universal sections $x_a$ of the stable curve $\bC_{A,g}$, the regular
sheaves $\CR_{n}$ into the nodes $n$ of this curve and pasting them together  
into one sheaf by the Heisenberg system $\CH_{\vmu;A,g}^{\nu}$.  

\subsection{}  
Given $A\in\Sew$ and an $A$-collection $\{\CX_a;\CR_n\}_A$,   
choose elements $\mu_a\geq\lambda(\CX_a)$ in $X_{\ell}$ such that
(a) above holds (note that $2\rho_{\ell}\in Y$). Below $\nu$ will denote the
element as in (a) above.  

Let $\CX'_a$ denote the factorizable sheaf 
isomorphic to $\CX_a$ obtained from it by the change of the highest 
weight from $\lambda(\CX_a)$ to $\mu_a$.
For each $m\in\BZ$, define a local system
$\langle\otimes_A\ \{\CX_a; \CR_n\}\rangle^{(m)}_g$  
over $\CM_{A,g;\delta}$ as follows.

Let $\langle\otimes_A\ \{\CX_a; \CR_n\}\rangle^{(m)}
_{\vA^0,g,\Gamma}$ denote the restriction of the local system to be defined to 
the connected component $\CM_{\vA^0,g,\Gamma;\delta}$. By definition,
$$
\langle\otimes_{A}\ \{\CX_a; \CR_n\}\rangle^{(m)}_{\vA^0,g,\Gamma}:=
R^{m-\dim \CM_{\vA^0,g,\Gamma;\delta}}
\baeta_*^{\nu}(\Boxtimes_{\vA^0,g,\Gamma}^{(\nu)}\ \{\CX'_a;\CR_n\}).
$$
Here 
$\Boxtimes_{\vA^0,g,\Gamma}^{(\nu)}\ \{\CX'_a;\CR_n\}$ denotes the perverse
sheaf  
$\Boxtimes_{A,g}^{(\nu)}\ \{\CX'_a;\CR_n\}$ restricted to the subspace
$\bC^{\nu}_{\vA^0,g,\Gamma;\delta}$.
This definition does not depend on the choice of the elements $\mu_a$.  

\subsection{} Given a morphism $f: B\lra A$, we define, acting as in
\ref{nul' conf fib}, \ref{nul' rep}, a perverse sheaf
$f^*(\Boxtimes^{(\nu)}_{A,g}\ \{\CX_a;\CR_n\})$ over
$\bC^{\nu}_{B,g;\delta}$ and local systems
$f^*\langle\otimes_A\ \{\CX_a;\CR_n\}\rangle_g^{(m)}$ over 
$\CM_{B,g;\delta}$.

\subsection{Theorem} {\em In the above notations, we have 
canonical isomorphisms 
$$
f^*(\Boxtimes^{(\nu)}_{A,g}\ \{\CX_a; \CR_n\})=\Boxtimes_{B,g}^{(\nu)}\
f^*\{\CX_a; \CR_n\}.\ 
$$}

This is a consequence of Theorem \ref{nul' cart} above and the definition of
the regular sheaf $\CR$ as a sheaf of nearby cycles of the braiding local 
system, \ref{nul' nearby}.

\subsection{Corollary} 
\label{nul' fus cor} {\em We have canonical isomorphisms of local
systems  
$$
f^*\langle\otimes_A\ \{\CX_a; \CR_n\}\rangle^{(m)}_g=\langle\otimes_B\
f^*\{\CX_a; \CR_n\}
\rangle^{(m)}_g\ (m\in\BZ).\ \Box
$$}
\subsection{} 
\label{nul' def fus} The previous corollary may be expressed as follows. The
various $A$-collections of
factorizable sheaves (resp. categories $\tRep_{c;A}$ of finite dimensional 
representations of "Teichm\"{u}ller groupoids" 
$\tTeich_{A}=\pi(\CM_{A;\delta})$ having monodromy $c$ around the
zero section of the determinant bundle)
define a fibered category $\FS^{\sharp}$ (resp. $\tRep_{c}$) over $\Sew$.

For any $m\in \BZ$, the collection of local systems
$\langle\otimes_A\ \{\CX_a; \CR_n\}\rangle^{(m)}_g$ and 
the canonical isomorphisms of the previous theorem define a 
{\bf cartesian functor}
$$
\langle\ \rangle^{(m)}:\FS^{\sharp}\lra\tRep_{c}  
$$
where $c=(-1)^{\card(I)}\zeta^{-12\rho\cdot\rho}$. We call such a functor 
a {\em fusion structure of multiplicative central 
charge $c$} on the category $\FS$.   
The category $\FS$ with this fusion structure 
will be denoted by $\FS^{(m)}$. 

The difference from the definition of a fusion structure given in 
\cite{bfm} 4.5.4 is that  
our fibered categories live over $\Sew=\Sets^{\sharp}/\emp$ and not
over $\Sets^{\sharp}$, as in {\em op. cit.}   
 
\subsection{Example} For $A\in\Sew$, let us consider an $A$-curve 
$P=(\BP^1,\{x_a,\tau_a\})$; it defines a geometric point $Q$
of the stack $\CM_{A,g}$ where $g=\card(N_A)$ and hence a geometric point 
$P=(Q,1)$ of the stack $\CM_{A,g;\delta}$ since the determinant bundle is
canonically trivialized at $Q$. 

\subsubsection{} {\bf Theorem.} {\em For an $A$-collection 
$\{\CX_a; \CR_n\}$, the stalk of the local system 
$\langle\otimes_A\ \CX_a\rangle^{(m)}_g$ at a point $P$ is isomorphic to 
$$
\Tor_{\infh-m}^{\CC}({\em \One}_r,(\otimes_{A^0}\ \Phi(\CX_a))\otimes
{\em \ad}^{\otimes g}).
$$}
To prove this, one should apply theorem \ref{nul' integral} and the following
remark. Degenerating all nodes into cusps, one can include the nodal 
curve $\overline{P}_A$ into a  
one-parameter family whose special fiber is $\BP^1$, with
$\card(A^0)+g$ marked points. The nearby cycles of the sheaf 
$\Boxtimes_A\ \{\CX_a; \CR_n\}$ will be the sheaf obtained by the gluing of the
sheaves $\CX_a$ and $g$ copies of the sheaf $\Phi^{-1}(\ad)$ into these
marked points. 

\subsection{} Note that an arbitrary $A$-curve may be degenerated into 
a curve considered in the previous example. Due to \ref{nul' fus cor},
this determines the stalks of all our local systems (up to a non canonical  
isomorphism).

\section{Conformal blocks (higher genus)}

In this section we assume that $\sk=\BC$.

\subsection{} Let us make the assumptions of \ref{nul' charge}, \ref{nul' level}.
Consider the full subcategory $\MS\subset\CC$. Let us define the 
{\em regular object} $\bR_{\MS}$ by
$$
\bR_{\MS}=
\oplus_{\mu\in\Delta_{\ell}}\ (\fL(\mu)\otimes\fL(\mu)^*)\in\MS^{\otimes 2}. 
$$
As in the previous section, we have a notion of $A$-collection 
$\{L_a;\bR_{\MS n}\}_A\ (A\in\Sew, \{L_a\}\in\MS^{A^0})$. The classical
fusion structure on $\MS$, \cite{tuy}, defines for each $A$-collection 
as above, a local system 
$$
\langle\otimes_A\ \{L_a;\bR_{\MS n}\}\rangle_{\MS}
$$
on the moduli stack $\CM_{A;\delta}$.

We have $A=(\bA,N_A)$; let us define another object $A'=(\bA\cup\{*\},N_A)$.
We have an obvious deleting $f_A:A'\lra A$. Given an $A$-collection as
above, define an $A'$-collection in $\{L_a;L(2\rho_{\ell});\bR_n\}_{A'}$ in
the category $\CC$. Using the equivalence $\Phi$, we transfer to $\CC$ the 
fusion structures defined in the previous section on $\FS$; 
we denote them $\langle\ \rangle^{(m)}_{\CC}$.  

The following theorem generalizes Theorem \ref{nul' subq} to higher genus.

\subsection{Theorem} 
\label{nul' subquot} {\em For each $A\in\Sew$, the local system
$f_A^*\langle\otimes_A\ \{L_a;{\em\bR}_{\MS n}\}\rangle_{\MS}$ on 
$\CM_{A';\delta}$ is a canonical subquotient of the local system
$\langle\otimes_{A'}\ \{L_a;L(2\rho_{\ell});{\em\bR}_n\}\rangle^{(0)}_{\CC}$.}
$\Box$

\subsection{} Let us consider the special case $A=(A^0,\emp)$; 
for an $A$-collection $\{L_a\}_{A^0}$ in $\MS$, we have the classical  
local systems of conformal blocks $\langle\otimes_{A^0}\ \{L_a\}\rangle_{\MS}$ 
on $\CM_{A;\delta}$.

\subsection{Corollary} {\em The local systems 
$\langle\otimes_{A^0}\ \{L_a\}\rangle_{\MS}$ are semisimple. They carry a 
canonical non-degenerate Hermitian form.} 

In fact, the local system $f^*_A\langle\otimes_{A^0}\ \{L_a\}\rangle_{\MS}$ 
is semisimple by the previous theorem and by 
Beilinson-Bernstein-Deligne-Gabber, \cite{bbd} 6.2.5. The
map of fundamental groupoids
$f_{A*}:\tTeich_{A'}\lra\tTeich_A$ is surjective; therefore the initial
local system is semisimple.

The Hermitian form is defined in the same manner as in genus zero, cf.
\ref{nul' hermit}.

\newpage
\setcounter{section}{0}
\begin{center}{\large \bf Part I. INTERSECTION COHOMOLOGY}\end{center}

\begin{center}{\large\bf OF REAL ARRANGEMENTS}\end{center}

\section{Introduction}


\subsection{} Let $\BA$ be a complex affine space, $\CH$ a finite collection
of hyperplanes in $\BA$. We suppose that all $H\in\CH$ are given by real
equations, and denote by $H_{\BR}\subset H$ the subspace of real points.
An arrangement $\CH$ induces naturally a stratification of
$\BA$ denoted by $\CS_{\CH}$ (for precise definitions see section 2 below).
The main goal of this part is the study
of the category $\CM(\BA;\CS_{\CH})$ of perverse sheaves (of vector spaces
over a fixed ground field) over $\BA$ which are smooth along $\CS_{\CH}$.

In {\em section 1} we collect topological notations and known facts we will
need in the sequel.

{\em Section 2.}
To each object $\CM\in\CM(\BA;\CS_{\CH})$ we assign "a linear algebra data".
Namely, for each {\em facet} of $\CH$ --- i.e. a connected component
of an intersection of some $H_{\BR}$'s --- we define a vector space
$\Phi_F(\CM)$; if $F$ is contained and has codimension $1$ in the closure of
another facet $E$, we define two linear operators between $\Phi_F(\CM)$ and
$\Phi_E(\CM)$ acting in opposite directions. These data contain all
information about $\CM$ we need.

In fact, it is natural to expect that the sheaf $\CM$ may be {\em uniquely
reconstructed} from the linear algebra data above (one can check this in
the simplest cases).

The spaces $\Phi_F(\CM)$ are defined using certain relative cohomology.
They are similar to a construction by R. MacPherson (see ~\cite{b2}).
The spaces $\Phi_F(\CM)$ are analogous to functors of vanishing cycles.

The main technical properties of functors $\CM\mapsto\Phi_F(\CM)$ are
(a) {\em commutation with Verdier duality}, and (b) {\em exactness}.
In fact, (b) follows from (a) without difficulty; the proof of (a)
is the principal geometrical result of this part (see Theorem ~\ref{odin dual}).
Actually, we prove a more general statement concerning all complexes
smooth along $\CS_{\CH}$. This is a result of primary importance for us.

Using these linear algebra data, we define an exact functor from
$\CM(\BA;\CS_{\CH})$ to complexes of finite dimensional vector spaces
computing cohomology $R\Gamma(\BA;\bullet)$, see Theorem ~\ref{odin rgamma}.
This is the main result we will need below. The idea of constructing such
functorial complexes on categories of perverse sheaves
(quite analogous to the usual computation of singular cohomology using
cell decompositions) is due to A. Beilinson,
{}~\cite{b}. It was an important source of inspiration for us.

A similar problem was considered in ~\cite{c}.

In {\em section 3} we present explicit computations for standard extensions
of an arbitrary one-dimensional local system over an open stratum.
The main result is the computation of {\em intersection cohomology}, see
Theorem ~\ref{odin inters}. In our computations we use some simple geometric
ideas from M. Salvetti's work, ~\cite{sa}. However, we do not use a more
difficult  main theorem of this paper. On the contrary, maybe our
considerations shed some new light on it.

\subsection{} The idea of this work appeared several years ago,
in some discussions with H\'{e}l\`{e}ne Esnault and Eckart Viehweg.
We use this occasion to express to them our sincere gratitude. We are also
very grateful to Paul Bressler for a stimulating discussion.
We thank Kari Vilonen for pointing out an error in an earlier version of the
manuscript.


\section{Topological preliminaries}

In this section we introduce our notations and recall some basic
facts from topology. Main references are ~\cite{ks}, \cite{bbd}.

\subsection{}
\label{odin derived} Throughout this work, all our topological
spaces will be {\em locally compact}, in particular {\em Hausdorff}.
$\{ pt\}$ will denote a one-point space. For a topological space $X$,
$a_X:X\lra \{pt\}$ will be the unique map. If $Y\subset X$ is a subspace,
$\bar{Y}$ will denote the closure of $Y$.

Throughout this chapter we fix an arbitrary ground field $\sk$. {\em A vector
space} will mean a vector space over $\sk$. For a finite dimensional
vector space $V$, $V^*$ will denote a dual space.
$\Vect$ will denote a category
of vector spaces.

{\em A sheaf} will mean a sheaf of vector spaces. For a topological
space $X$, $\Sh(X)$ will denote a category of sheaves over $X$,
$\CCD^*(X)$ will denote the derived category of $\Sh(X)$,
$*=+,-,b$ or $\emp$ will have the usual meaning.

For $p\in\Bbb Z$, $\CCD^{\leq p}(X)\subset \CCD(X)$ will denote
a full subcategory of complexes $\CK$ with $H^i(\CK)=0$ for all
$i>p$.

If $\CA,\CB$ are abelian categories, we will say that a left exact functor
$F:\CCD^*(\CA)\lra\CCD(\CB)$ ($*=+,-,b$ or $\emp$)
has {\em cohomological dimension $\leq r$}, and write
$\cd (F)\leq r$, if $H^i(F(A))=0$ for any $A\in\CA$ and $i>r$. (Left exactness
here means that $H^i(F(A))=0$ for $i<0$).

\subsection{} Let $f:X\lra Y$ be a continuous map. We will use the following
notations for standard functors.

$f^*:\CCD(Y)\lra\CCD(X)$ --- the inverse image;
$f_*:\CCD^+(X)\lra\CCD^+(Y)$ --- (the right derived of) the direct image;
$f_!:\CCD^+(X)\lra\CCD^+(Y)$ --- (the right derived of) the direct image
with proper supports;
$f^!:\CCD^+(Y)\lra\CCD^+(X)$ --- the right adjoint to $f_!$ (defined
when $f_!$ has finite cohomological dimension), see ~\cite{ks}, 3.1.

We will denote the corresponding functors on sheaves as follows:
$f^*:\Sh(Y)\lra\Sh(X);\ R^0f_*:\Sh(X)\lra\Sh(Y);\
R^0f_!:\Sh(X)\lra\Sh(Y)$.

We will denote $R\Gamma(X;\cdot):=p_{X*};\  R\Gamma_c(X;\cdot):=p_{X!}$.
For $\CK\in\CCD^+(X),\ i\in\BZ$, we set
$H^i(X;\CK):=H^i(R\Gamma(X;\CK));\  H^i_c(X;\CK):=H^i(R\Gamma_c(X;\CK))$.

For $V\in\Vect=\Sh(\{ pt\})$ we denote by $V_X$ the constant sheaf $p^*_XV$.

If $Y$ is a subspace of $X$, we will use a notation $i_{Y,X}$,
or simply $i_Y$ if $X$ is understood, for the
embedding $Y\hra X$. If $Y$ is open in $X$, we will also write $j_Y$
instead of $i_Y$.

For $\CK\in\CCD(X)$, we will use notations $\CK|_Y:=i^*_{Y,X}\CK$;
$R\Gamma(Y;\CK):=R\Gamma(Y;\CK|_Y),\ H^i(Y;\CK):=H^i(Y;\CK_Y)$,
etc.

\subsection{} We have functors
\begin{equation}
R\CHom :\CCD^-(X)^{opp}\times\CCD^+(X)\lra\CCD^+(X),
\end{equation}
\begin{equation}
\otimes:\CCD^-(X)\times\CCD^*(X)\lra\CCD^*(X)
\end{equation}
where $*=,-,b$ or $\emp$;
\begin{equation}
R\Hom=R\Gamma\circ R\CHom:\CCD^-(X)^{opp}\times\CCD^+(X)\lra\CCD^+(\Vect)
\end{equation}
For $\CK,\CL\in\CCD^b(X)$ we have
\begin{equation}
\Hom_{\CCD(X)}(\CK,\CL)=H^0(R\Hom(\CK,\CL))
\end{equation}

We denote $\CK^*:=R\CHom(\CK,\sk_X)$.

\subsection{} Let $X$ be a topological space, $j:U\lra X$ an embedding
of an open subspace, $i: Y\lra X$ an embedding of the complement.
In this case $i_*=i_!$ and $j^!=j^*$. $R^0i_*$ and $R^0j_!$ are
exact, so we omit $R^0$ from their notation.
$j_!$ is the functor of extension by zero.
$i^!$ is the right derived of the functor $R^0i^!:\Sh(X)\lra \Sh(Y)$
of the subsheaf of sections with supports in $Y$ (in notations of ~\cite{ks},
$R^0i^!(F)=\Gamma_Y(F)$).

Let $\CK\in\CCD^+(X)$. We will use the following notations for relative
cohomology: $R\Gamma(X,Y;\CK):=R\Gamma(X;j_!\CK|_U);\
H^i(X,Y;\CK):=H^i(R\Gamma(X,Y;\CK))$.

If $Z\hra Y$ is a closed subspace,
we have a canonical exact triangle
\begin{equation}
i_{X-Y!}\CK|_{X-Y}\lra i_{X-Z!}\CK|_{X-Z}
\lra i_{Y-Z!}\CK|_{Y-Z} \lra i_{X-Y!}\CK|_{X-Y}[1]
\end{equation}
(of course, $i_{X-Y}=j$)
inducing a long cohomology sequence
\begin{equation}
\label{odin cohom}
\ldots\lra H^i(X,Y;\CK)\lra H^i(X,Z;\CK)\lra H^i(Y,Z;\CK)
\overset{\dpar}{\lra} H^{i+1}(X,Y;\CK)\lra\ldots
\end{equation}

\subsection{} Let $X$ be a topological space such that the functor
$R\Gamma_c(X;\cdot)$ has a finite cohomological dimension.
We define {\em a dualizing complex}:
\begin{equation}
\omega_X:=a_X^!\sk\in\CCD^b(X)
\end{equation}
For $\CK\in\CCD^b(X)$ we set $D_X\CK=R\CHom(\CK,\omega_X)\in\CCD^b(X)$.
If there is no risk of confusion we denote $D\CK$ instead of $D_X\CK$.
We get a functor
\begin{equation}
D:\CCD^b(X)^{opp}\lra\CCD^b(X)
\end{equation}
which comes together with a natural transformation
\begin{equation}
\Id_{\CCD^b(X)}\lra D\circ D
\end{equation}

\subsection{Orientations} (Cf. ~\cite{ks}, 3.3.) Let $X$ be an $n$-dimensional
$C^0$-manifold.
We define an {\em orientation sheaf} $\Or_X$ as the sheaf associated to
a presheaf $U\mapsto \Hom_{\BZ}(H^n_c(U;\BZ),\BZ)$. It is a locally constant
sheaf of abelian groups of rank $1$. Its tensor square is constant.
By definition, {\em orientation of $X$} is an isomorphism
$\Or_X\iso\BZ_X$.

We have a canonical isomorphism
\begin{equation}
\omega_X\iso\Or_X\otimes_{\BZ}\sk[n]
\end{equation}

\subsection{Homology} Sometimes it is quite convenient to use
homological notations. Let $Y\subset X$ be a closed subspace of
a topological space $X$ ($Y$ may be empty), $\CK\in\CCD^b(X)$.
We define homology groups as
$$
H_i(X,Y;\CF):=H^i(X,Y;\CF^*)^*.
$$
These groups behave covariantly with respect to continuous mappings.

Let $\sigma$ be {\em a relative singular $n$-cell}, i.e. a continuous
mapping
$$
\sigma: (D^n,\dpar D^n)\lra (X,Y)
$$
where, $D^n$ denotes a standard closed unit ball in $\BR^n$.
We supply $D^n$ with the standard orientation. Let $\DO^n$ denote the interior
of $D^n$.

Suppose that $\sigma^*\CK|_{\DO^n}$ is constant. Then by Poincar\'{e} duality
we have isomorphisms
$$
H_n(D^n,\dpar D^n;\sigma^*\CK)=H^n_c(\DO^n;\sigma^*\CK^*)^*\iso
H^0(\DO^n;\sigma^*\CK).
$$
(recall that $\DO^n$ is oriented).
Thus, given an element $s\in H^0(\DO^n;\CK)$, we can define a homology class
$$
cl(\sigma;s):=\sigma_*(s)\in H_n(X,Y;\CK).
$$
We will call the couple $(\sigma;s)$
{\em a relative singular $n$-cell for $\CK$}.

These classes enjoy the following properties.

\subsubsection{Homotopy} Let us call to cells $(\sigma_0;s_0)$
and $(\sigma_1,s_1)$ {\em homotopic} if there exists a map
$$
\sigma: (D^n\times I,\dpar D^n\times I)\lra (X,Y)
$$
(where $I$ denotes a unit interval) such that $\sigma^*\CK|_{\DO^n\times I}$
is constant,  and a section
$$
s\in H^0(\DO^n\times I;\sigma^*\CK)
$$
such that $(\sigma,s)$ restricted to $D^n\times\{ i\}$
is equal to $(\sigma_i;s_i),\ i=0,1$.

We have

(H) {\em if $(\sigma_0;s_0)$ is homotopic to $(\sigma_1;s_1)$ then
$cl(\sigma_0;s_0)=cl(\sigma_1;s_1)$.}

\subsubsection{Additivity} Suppose $D^n$ is represnted as a union of its
upper and lower half-balls $D^n=D^n_1\cup D^n_2$ where
$D^n_1=\{ (x_1,\ldots, x_n)\in D^n|\ x_n\geq 0\}$ and
$D^n_2=\{ (x_1,\ldots, x_n)\in D^n|\ x_n\leq 0\}$. Let us supply $D^n_i$
with the induced orientation.

Suppose we are given a relative $n$-cell $(\sigma;s)$ such that
$\sigma(D^n_1\cap D^n_2)\subset Y$.
Then its restriction to $D^n_i$ gives us two relative $n$-cells
$(\sigma_i;s_i)$,
$i=1,2$. We have

(A) $cl(\sigma;s)=cl(\sigma_1;s_1)+cl(\sigma_2;s_2)$.

\subsection{} We will call {\em a local system} a locally constant sheaf
of finite rank.

Let $X$ be a topological space, $i:Y\hra X$ a subspace. We will say that
$\CK\in\CCD(X)$ is {\em smooth} along $Y$ if all cohomology sheaves
$H^p(i^*\CK),\ p\in\BZ,$ are local systems.

We will call a {\em stratification} of $X$ a partition $X=\bigcup_{S\in\CS}S$
into a finite disjoint of locally closed subspaces ---
{\em strata} --- such that the closure of each stratum is the union
of strata.

We say that $\CK\in\CCD(X)$ is {\em smooth along $\CS$} if it is smooth
along each stratum. We will denote by $\CCD^*(X;\CS)$ the full subcategory
of $\CCD^*(X)$ ($*=+,-,b$ or $\emp$) consisting of complexes smooth along $\CS$.

\subsection{} Let us call a stratification $\CS$ of a topological space
$X$ {\em good} if the following conditions from ~\cite{bbd}, 2.1.13
hold.

(b) All strata are equidimensional topological varieties. If a stratum
$S$ lies in the closure of a stratum $T$, $\dim S<\dim T$.

(c) If $i:S\hra X$ is a stratum, the functor $i_*:\CCD^+(S)\lra \CCD^+(X)$
has a finite cohomological dimension. If $\CF\in\Sh(S)$ is a local system,
$i_*(\CF)$ is smooth along $\CS$.

Let $\CS$ be a good stratification such that all strata have even dimension.
We will denote by $\CM(X;\CS)\subset\CCD^b(X;\CS)$ the category of
smooth along $\CS$
perverse sheaves corresponding to the middle perversity $p(S)=-\dim\ S/2$, cf.
{\em loc. cit.}

\subsection{} Let $X$ be a complex algebraic variety. Let us call
its stratification $\CS$ {\em algebraic} if all strata are algebraic
subvarieties. Following ~\cite{bbd}, 2.2.1, let us call a sheaf
$\CF\in\Sh(X)$ {\em constructible} if it is smooth along some algebraic
stratification. (According to Verdier, every algebraic stratification
admits a good refinement.) We denote by $\CCD^b_c(X)$ a full subcategory
of $\CCD^b(X)$ consisting of complexes with constructible cohomology.

We denote by $\CM(X)\subset\CCD^b_c(X)$ the category of perverse sheaves
corresponding to the middle perversity. It is a filtered union of categories
$\CM(X;\CS)$ over all good algebraic stratifications $\CS$, cf. {\em loc.cit.}


\subsection{} Let $X$ be a complex manifold, $f:X\lra \BC$ a holomorphic map.
Set $Y=f^{-1}(0),\ U=X-Y$, let $i:Y\hra X,\ j:U\hra X$ denote the
embeddings.

We define a functor of {\em nearby cycles}
\begin{equation}
\label{odin near}
\Psi_f:\CCD^b(U)\lra\CCD^b(Y)
\end{equation}
as $\Psi_f(\CK)=\psi_f j_*(\CK)[-1]$ where $\psi_f$ is defined in ~\cite{ks},
8.6.1.

We define a functor of {\em vanishing cycles}
\begin{equation}
\label{odin vanish}
\Phi_f:\CCD^b(X)\lra\CCD^b(Y)
\end{equation}
as $\phi_f$ from {\em loc. cit.}

\subsection{}
\label{odin limit}
\begin{lem}{} Let $X$ be topological space, $Y\hra X$ a closed subspace,
$\CF\in\Sh(X)$. Then natural maps
$$
\dirlim\ H^i(U;\CF)\lra H^i(Y;\CF),\ i\in\BZ,
$$
where $U$ ranges through all the open neighbourhoods of $U$, are
isomorphisms.
\end{lem}

{\bf Proof.} See ~\cite{ks}, 2.5.1, 2.6.9. $\Box$

\subsection{Conic sheaves}
\label{odin conic} (Cf. ~\cite{ks}, 3.7.) Let $\BR^{*+}$ denote
the multiplicative group of positive real numbers. Let $X$ be a topological
space endowed with an $\BR^{*+}$-action.

Following {\em loc.cit.}, we will call a sheaf $\CF$ over $X$ {\em conic}
(with respect to the given $\BR^{*+}$-action)
if its restriction to every $\BR^{*+}$-orbit is constant.
We will call a complex $\CK\in\CCD(X)$ {\em conic} if all its cohomology
sheaves are conic.

We will denote by $\Sh_{\BR^{*+}}(X)\subset \Sh(X)$,
$\CCD^*_{\BR^{*+}}(X)\subset\CCD^*(X),\ *=b,+,-$ or $\emp$,
the full subcategories of conic objects.

\subsubsection{}
\label{odin conicl}
\begin{lem}{} Let $U\hra X$ be an open subset. Suppose that for every
$\BR^{*+}$-orbit $O\subset X$, $O\cap U$ is contractible (hence, non-empty).
Then for every conic $\CK\in\CCD^+(X)$ the restriction morphism
$$
R\Gamma(X;\CK)\lra R\Gamma(U;\CK)
$$
is an isomorphism.
\end{lem}

{\bf Proof.} See {\em loc. cit.}, 3.7.3. $\Box$

\section{Vanishing cycles functors}
\label{odin vanishcycl}

\subsection{Arrangements} Below we use some terminology from ~\cite{sa}.
Let $\BA_{\BR}$ be a
real affine space $\BA_{\BR}$ of dimension $N$, and $\CH_{\BR}=\{ H_{\BR}\}$
a finite set of distinct real affine hyperplanes in $\BA_{\BR}$. Such a set
is called a {\em real arrangement}. We pick once and for all a square root
$i=\sqrt{-1}\in\BC$.

Let $\BA=\BA_{\BR}\otimes_{\BR}\BC$ denote the complexification of
$\BA_{\BR}$, and $\CH=\{ H\}$ where $H:=H_{\BR}\otimes_{\BR}\BC$.
(A finite set of complex hyperplanes in a complex affine space
will be called a {\em complex arrangement}).

We will say that $\CH$ is {\em central} if $\bigcap_{H\in\CH}H$ consists
of one point.

For a subset $\CK\subset \CH$ denote
$$
H_{\CK}=\bigcap_{H\in \CK}H;\ _{\CK}H=\bigcup _{H\in \CK}H,
$$
and
$$
H_{\BR,\CK}=\bigcap_{H\in \CK}H_{\BR},
\ _{\CK}H_{\BR}=\bigcup_{H\in \CK}H_{\BR}.
$$
The nonempty subspaces $H_{\CK}$ and $H_{\CK,\BR}$ are called complex
and real {\em edges} respectively. Set
$$
\HO_{\CK}=H_{\CK}-\cup L,
$$
the union over all the complex edges $L\subset H_{\CK},\ L\neq H_{\CK}$. We set
$\HO_{\BR,\CK}=H_{\BR,\CK}\cap \BA_{\BR}$.
Connected components of $\HO_{\BR,\CK}$ are called {\em facets} of $\CH_{\BR}$.
Facets of codimension $0$ (resp., $1$) are called {\em chambers}
(resp., {\em faces}).

Let us denote by $\CS_{\CH}$ a stratification of $\BA$ whose strata are
all non-empty $\HO_{\CK}$. We will denote by $\BAO_{\CH}$ a unique open
stratum
$$
\BAO_{\CH}=\HO_{\emp}=\BA- _{\CH}H.
$$

In this Section we will study categories of sheaves $\CCD(\BA;\CS_{\CH})$
and $\CM(\BA;\CS_{\CH})$.

We will denote by $\CS_{\CH_{\BR}}$ a stratification of $\BA_{\BR}$ whose
strata are all facets. We set
$$
\BAO_{\CH,\BR}=\BA_{\BR}- _{\CH}H_{\BR}.
$$
It is a union of all chambers.

\subsection{Dual cells}
\label{odin salv}(cf. ~\cite{sa}).
Let us fix a point $^Fw$ on each facet $F$. We will call this set of points
$\bw =\{^Fw\}$ {\em marking} of our arrangement.

For two facets $F,\ E$ let us write $F<E$ if $F\subset\BE$ and
$\dim\ F<\dim\ E$.
We will say that $E$ is {\em adjacent} to $F$. We will denote by $\Ch (F)$
the set of all chambers adjacent to $F$.

Let us call a {\em flag} a sequence
of $q-p+1$ facets $\bF=(F_p<F_{p+1}<\ldots <F_q)$ with $\dim\ F_i=i$.
We say that $F_p$ is {\em the beginning} and $F_q$ {\em the end} of $\bF$.

Let us denote by $^{\bF}\Delta$ a closed $(q-p)$-symplex with vertices
$^{F_p}w,\ldots,^{F_q}w$. Evidently, $^{\bF}\Delta\subset\BF_q$.

Suppose we are given two facets $F_p<F_q$, $\dim F_i=i$. We will denote
$$
D_{F_p<F_q}=\bigcup\ ^{\bF}\Delta,
$$
the sum over all flags beginning at $F_p$ and ending at $F_q$. This is a
$(q-p)$-dimensional cell contained in $\BF_q$.

For a facet $F$ let us denote
$$
D_F=\bigcup_{C\in \Ch (F)}D_{F<C}
$$
We set
$$
S_F:=\bigcup_{E,C: F<E<C;}D_{E<C},
$$
the union over all facets $E$ and chambers $C$. The space $S_F$ is contained
in $D_F$ (in fact, in the defintion of $S_F$ it is enough to take the union
over all facets $E$ such that
$\dim\ E=\dim\ F+1$). If $q=\codim\ F$, then
$D_F$ is homeomorphic to a $q$-dimensional disc and $S_F$ --- to a $(q-1)$-
dimensional sphere, cf. ~\cite{sa}, Lemma 6. We denote $\DO_F:=D_F-S_F$.
We will call $D_F$ {\em a dual cell corresponding to $F$}.

We set $\DO_{F<C}:=D_{F<C}\cap \DO_F$.

\subsection{Generalized vanishing cycles} Let $\CK\in\CCD^b(\BA;\CS_{\CH})$,
$F$ a $p$-dimensional facet. Let us introduce a complex
\begin{equation}
\label{odin vc}
\Phi_F(\CK):=R\Gamma(D_F,S_F;\CK)[-p]\in\CCD^b(pt)
\end{equation}
This complex will be called {\em a complex of generalized vanishing
cycles of $\CK$ at a facet $F$}.

Formally, the definition of functor $\Phi_F$ depends
upon the choice of a marking $\bw$. However, functors defined using
two different markings are canonically isomorphic. This is evident.
Because of this, we omit markings from the notations.

\subsection{Transversal slices}
\label{odin slice} Let $F$ be a facet of dimension $p$ which
is a connected component of a real edge $M_{\BR}$ with the complexification
$M$.
Let us choose a real affine subspace $L_{\BR}$ of codimension $p$ transversal
to $F$ and passing through $^Fw$. Let $L$ be its complexification.

Let us consider a small disk $L_{\epsilon}\subset L$ with the centrum
at $^Fw=L\cap M$. We identify $L_{\epsilon}$ with an affine space
by dilatation. Our arrangement induces a central arrangement $\CH_L$
in $L_{\epsilon}$. Given  $\CK\in\CCD^b(\BA;\CS_{\CH})$,
define $\CK_L:=i_{L_{\epsilon}}^*\CK[-p]$.

\subsubsection{} {\bf Lemma.} {\em We have a natural isomorphism
\begin{equation}
\label{odin slicedu}
D(\CK_L)\iso (D\CK)_L
\end{equation}}

{\bf Proof.} Consider an embedding of smooth complex manifolds
$i_L: L\lra \BA$. Let us consider the following complexes:
$\omega_{L/\BA}:=i^!_L \sk_{\BA}$ (cf. ~\cite{ks}, 3.1.16 (i)) and
$\Or_{L/\BA}:=\Or_{L}\otimes_{\BZ}\Or_{\BA}$.
We have canonical isomorphism
$$
\omega_{L/\BA}\iso\Or_{L/\BA}[-2p]\otimes \sk.
$$
The chosen orientation of $\BC$ enables us to identify
$\Or_{\BA}$ and $\Or_L$, and hence $\Or_{L/\BA}$ with constant sheaves;
consequently, we get an isomorphism
$$
\omega_{L/\BA}\iso \sk[-2p].
$$
The canonical map
\begin{equation}
\label{odin closed}
i^!_L\CK\lra i^*\CK\otimes\omega_{L/\BA}
\end{equation}
(cf. ~\cite{ks}, (3.1.6)) is an isomorphism since singularities of $\CK$
are transversal to $L$ (at least in the neighbourhood of $^Fw$).
Consequently we get an isomorphism $i^!_L\CK\iso i^*\CK_L[-2p]$.

Now we can compute:
$$
D(\CK_L)=D(i_L^*\CK[-p])\iso D(i^*_L\CK)[2p]\iso i^!_LD\CK[p]
\iso i^*_LD\CK[-p]=(D\CK)_L,
$$
QED. $\Box$

\subsubsection{} {\bf Lemma.} {\em We have a natural
isomorphism
\begin{equation}
\label{odin sliceiso}
\Phi_{F}(\CK)\cong\Phi_{\{^Fw\} }(\CK_L)
\end{equation}}

This follows directly from the definition of functors $\Phi$.

This remark is often useful for reducing the study of functors $\Phi_F$ to
the case of a central arrangement.

\subsubsection{} {\bf Lemma.} {If $\CM\in\CM(\BA;\CS_{\CH})$
then $\CM_L\in\CM(L_{\epsilon};\CS_{\CH_{L_{\epsilon}}})$.}

{\bf Proof} follows from transversality of $L_{\epsilon}$ to singularities
of $\CM$. $\Box$

\subsection{Duality}
\label{odin dual}
\begin{thm}{} Functor $\Phi_F$ commutes with Verdier duality.
More precisely, for every $\CK\in\CCD^b(\BA,\CS_{\CH})$ there exists a natural
isomorphism
\begin{equation}
\label{odin duala}
\Phi_F(D\CK)\iso D\Phi_F(\CK)
\end{equation}
\end{thm}

This is the basic property of our functors.

\subsection{$1$-dimensional case}
\label{odin one-dim} To start with the proof, let us treat first the simplest case.
Consider an arrangement consisting of one point --- origin --- in
a one-dimensional
space $\BA$. It has one $0$-face $F$ and two $1$-faces $E_{\pm}$ --- real rays
$\BR_{>0}$ and $\BR_{<0}$. A marking consists of two points
$w_{\pm}\in E_{\pm}$ and $F$. We set $w_{\pm}=\pm 1$ (we pick a coordinate on
$\BA$).

For a positive
$r$ denote $\BA_{\leq r}:=\{ z\in\BA|\ |z|\leq r\},\ S_r:=\dpar\BA_{\leq r}$;
$\BA_{< r},\ \BA_{\geq r}$, etc. have an evident meaning. We also set
$\BA_{(r',r'')}:=\BA_{>r'}\cap \BA_{<r''}.$ A subscript $\BR$
will denote an intersection of these subsets with $\BA_{\BR}$.
Evidently, $D_F=\BA_{\leq 1,\BR},\ S_F=S_{1,\BR}$. Define $D_F^{opp}:=
\BA_{\geq 1,\BR},\ Y:=i\cdot D_F^{opp}$.

One sees easily (cf. {\em infra}, Lemma ~\ref{odin isom}) that one has isomorphisms
\begin{equation}
\label{odin isom-one}
\Phi_F(\CK)\cong R\Gamma(\BA,S_F;\CK)\iso R\Gamma(\BA,D_F^{opp};\CK).
\end{equation}

Let us choose real numbers $\epsilon, r', r''$ such that $0<\epsilon<r'<1<r''$.
Set $Y:=\epsilon i\cdot D_F^{opp}$. Denote $j:=j_{\BA-S_F}$.
We have natural isomorphisms
$$
D\Phi_{F}(\CK)\cong DR\Gamma(\BA,S_{F};\CK)
\cong R\Gamma_c(\BA; j_{*}j^*D\CK)
$$
(Poincar\'{e} duality)
$$
\cong R\Gamma(\BA,\BA_{\geq r''};j_{*}j^*D\CK)\cong
R\Gamma(\BA,Y\cup \BA_{\geq r''};j_{*}j^*D\CK)
$$
(homotopy).
Consider the restriction map
\begin{equation}
res: R\Gamma(\BA,Y\cup \BA_{\geq r''};j_*j^*D\CK)\lra
R\Gamma(\BA_{\leq r'},Y\cap \BA_{\leq r'};D\CK)
\end{equation}
We claim that $res$ is an isomorphism.
In fact, $\Cone (res)$ is isomorphic to
\begin{eqnarray}
R\Gamma(\BA,\BA_{\leq r'}\cup\BA_{\geq r''}\cup Y;j_*j^*D\CK)=
R\Gamma_c(\BA_{< r''},\BA_{\leq r'}\cup Y;j_*j^*D\CK)\cong\nonumber\\
\cong DR\Gamma(\BA_{< r''}-(\BA_{\leq r'}\cup Y);j_!j^*\CK)\nonumber
\end{eqnarray}
We have by definition
$$
R\Gamma(\BA_{< r''}-(\BA_{\leq r'}\cup Y);j_!j^*\CK)=
R\Gamma(\BA_{(r',r'')}-Y,S_F;\CK)
$$
On the other hand, evidently $\CK$ is smooth over $\BA_{(r',r'')}$, and we
have an evident retraction of $\BA_{(r',r'')}$ on $S_F$ (see Fig. 1).
Therefore, $R\Gamma(\BA_{(r',r'')}-Y,S_F;\CK)=0$ which proves the claim.

\begin{picture}(20,8)(-10,-4)

\put(0,0){\circle{0.2}}
\put(0,0){\oval(6,6)}
\put(-1.5,3.3){$S_{r''}$}

\put(-4,0){\line(1,0){8}}

\put(0,3){\line(0,-1){2}}
\put(0,-3){\line(0,1){2}}
\put(-0.4,1.8){$Y$}

\put(0,1){\circle*{0.2}}
\put(0,0.5){$i\epsilon w_+$}
\put(0,-1){\circle*{0.2}}

\put(2,0){\circle*{0.2}}
\put(2,-0.5){$w_+$}
\put(-2,0){\circle*{0.2}}
\put(-2,-0.5){$w_-$}

\put(0,0){\oval(3,3)}
\put(-1.5,1.6){$S_{r'}$}

\put(0,0){\oval(5,5)[tr]}
\put(0.5,2.5){\circle*{0.2}}
\put(0.5,2){$z$}
\put(0.5,2.5){\vector(1,0){1}}

\put(-0.5,-4){Fig. 1}

\end{picture}

A clockwise rotation by $\pi/2$ induces an isomorphism
$$
R\Gamma(\BA_{\leq r'},Y\cap \BA_{\leq r'};D\CK)\cong
R\Gamma(\BA_{\leq r'},\epsilon\cdot D_F^{opp};D\CK),
$$
and the last complex is isomorphic to $\Phi_F(D\CK)$ by dilatation and
{}~(\ref{odin isom-one}).
This proves the theorem for $\Phi_F$. The statement for  functors
$\Phi_{E_{\pm}}$ is evident.

Let us return to the case of an arbitrary arrangement.

\subsection{}
\label{odin isom}
\begin{lem}{} Suppose that $\CH$ is central. Let $F$ be the unique $0$-
dimensional facet. The evident restriction maps induce canonical isomorphisms
$$
\Phi_{F}(\CK)\overset{(1)}{\cong} R\Gamma(\BA_{\BR},S_{F};\CK)
\overset{(2)}{\cong} R\Gamma(\BA,S_{F};\CK)
\overset{(3)}{\cong} R\Gamma(\BA,D_{F}^{opp};\CK)
\overset{(4)}{\cong} R\Gamma(\BA_{\BR},D_{F}^{opp};\CK)
$$
where $D_{F}^{opp}=\BA_{\BR}-\DO_{F}$.
\end{lem}

{\bf Proof.} Let us fix a coordinate system with the origin at $F$, and
hence a metric on $\BA$.
For $\epsilon>0$ let $U_{\epsilon}\subset \BA_{\BR}$ denote
the set of points $x\in \BA_{\BR}$ having distance $<\epsilon$ from
$D_{F}$.

It follows from ~\ref{odin limit} that
$$
\Phi_{F}(\CK)=R\Gamma(D_{F},S_{F};\CK)\cong
\dirlim_{\epsilon}\ R\Gamma(U_{\epsilon},S_{F};\CK).
$$
On the other hand, from ~\ref{odin conicl} it follows that restriction maps
$$
R\Gamma(\BA_{\BR},S_{F};\CK)\lra R\Gamma(U_{\epsilon},S_{F};\CK)
$$
are isomorphisms. This establishes an isomorphism (1). To prove (3), one
remarks that its cone is acyclic.

The other isomorphisms are proven by the similar arguments. We leave
the proof to the reader. $\Box$

\subsection{Proof of ~\ref{odin dual}}
\label{odin mainconst} First let us suppose that $\CH$ is central
and $F$ is its $0$-dimensional facet. We fix a coordinate system in $\BA$ as
in the proof of the lemma above.
We have a decomposition $\BA=\BA_{\BR}\oplus i\cdot \BA_{\BR}$. We will denote
by $\Re,\Im:\BA\lra\BA_{\BR}$ the evident projections.

Let
$\BA_{\leq r}$, etc. have the meaning similar to the one-dimensional case
above.
Let us denote $j:=j_{\BA-S_F}$. We proceed as in one-dimensional case.

Let us choose positive numbers $r',r'',\ \epsilon,$ such that
$$
\epsilon D_F\subset \BA_{<r'}\subset \DO_F\subset D_F\subset
\BA_{< r''}
$$
Let us introduce a subspace
$$
Y=\epsilon i\cdot D_F^{opp}
$$
We have isomorphisms
$$
D\Phi_{F}(\CK)\cong DR\Gamma(\BA,S_{F};\CK)
$$
(by Lemma ~\ref{odin isom})
$$
\cong R\Gamma_c(\BA; j_{*}j^*D\CK)
\cong R\Gamma(\BA,\BA_{\geq r''};j_{*}j^*D\CK)\cong
R\Gamma(\BA,Y\cup \BA_{\geq r''};j_{*}j^*D\CK)
$$
(homotopy).
Consider the restriction map
\begin{equation}
res: R\Gamma(\BA,Y\cup \BA_{\geq r''};j_*j^*D\CK)\lra
R\Gamma(\BA_{\leq r'},Y\cap \BA_{\leq r'};D\CK)
\end{equation}
$\Cone (res)$ is isomorphic to
\begin{eqnarray}
R\Gamma(\BA,\BA_{\leq r'}\cup\BA_{\geq r''}\cup Y;j_*j^*D\CK)=
R\Gamma_c(\BA_{< r''},\BA_{\leq r'}\cup Y;j_*j^*D\CK)\cong\nonumber\\
\cong DR\Gamma(\BA_{< r''}-(\BA_{\leq r'}\cup Y);j_!j^*\CK)\nonumber
\end{eqnarray}
We have by definition
$$
R\Gamma(\BA_{< r''}-(\BA_{\leq r'}\cup Y);j_!j^*\CK)=
R\Gamma(\BA_{(r',r'')}-Y,S_F;\CK)
$$

\subsubsection{}
\label{odin acycl} {\bf Lemma.\ }{\em Set $\BA'_{\BR}:=\BA_{\BR}-\{ 0\}.$ We have
$R\Gamma(\BA-i\cdot \BA_{\BR},\BA'_{\BR};\CK)=0$.}

{\bf Proof.}
We have to prove that the restriction
map
\begin{equation}
\label{odin puncture}
R\Gamma(\BA-i\cdot \BA_{\BR};\CK)\lra R\Gamma(\BA'_{\BR};\CK)
\end{equation}
is an isomorphism. Consider projection to the real part
$$
\pi: \BA-i\cdot \BA_{\BR}\lra\BA'_{\BR}.
$$
Evidently,
$$
R\Gamma(\BA-i\cdot \BA_{\BR};\CK)=R\Gamma(\BA'_{\BR},\pi_*\CK).
$$
On the other hand, since $\CK$ is conic along fibers of $\pi$, we have
$$
\pi_*\CK\iso e^*\CK
$$
where $e:\BA'_{\BR}\lra \BA-i\cdot \BA_{\BR}$ denotes the embedding
(cf. \ref{odin conic}). This implies our lemma. $\Box$

It follows easily that $\Cone(res)$ is acyclic, therefore the map $res$ is an
isomorphism. In other words, we have constructed an isomorphism
$$
D\Phi_F(\CK)\iso R\Gamma(\BA_{\leq r'},Y\cap \BA_{\leq r'};D\CK).
$$
A clockwise rotation by $\pi/2$ induces an isomorphism
$$
R\Gamma(\BA_{\leq r'},Y\cap \BA_{\leq r'};D\CK)\cong
R\Gamma(\BA_{\leq r'},\epsilon\cdot D_F^{opp};D\CK),
$$
and the last complex is isomorphic to $\Phi_F(D\CK)$ by dilatation and
{}~\ref{odin isom} (3).

This proves the theorem for $\Phi_F$. Note that we have constructed an
explicit isomorphism.

\subsubsection{}
\label{odin transv} If $F$ is an arbitrary facet, we consider a transversal
slice $L$ as in ~\ref{odin slice}. By the results of {\em loc.cit.}, and
the above proven case, we have natural isomorphisms
$$
D\Phi_F(\CK)\iso D\Phi_{\{^Fw\} }(\CK_L)\iso \Phi_{\{^Fw\} }(D\CK_L)\iso
\Phi_{\{^Fw\}}((D\CK)_L)\iso \Phi_F(D\CK).
$$
This proves the theorem. $\Box$

\subsection{}
\label{odin exact}
\begin{thm}{} For every $\CM\in\CM(\BA;\CS_{\CH})$ and every facet $F$
we have $H^i(\Phi_{F}(\CM))=0$ for all $i\neq 0$.
\end{thm}

In other words, functors $\Phi_F$ are $t$-exact with respect to the middle
perversity.

{\bf Proof.} First let us suppose that $\CH$ is central and $F$ is its
$0$-dimensional facet. Let us prove that $\Phi_F$ is right exact,
that is, $H^i(\Phi_F(\CM))=0$ for $i>0$ and every $\CM$ as above.

In fact, we know that

\subsubsection{}
\label{odin estim} {\em if $S$ is any stratum of $\CS_{\CH}$ of dimension $p$
then $\CM|_S\in\CCD^{\leq -p}(S)$}

by the condition of perversity. In particular, $R\Gamma(\BA;\CM)=
i_F^*\CM\in\CCD^{\leq 0}(\{ pt\})$.

On the other hand, one deduces from ~\ref{odin estim} that
$R\Gamma(S_F;\CM)\in\CCD^{\leq -1}(\{ pt\})$. In fact, by definition
$S_F$ is a union of certain simplices $\Delta$, all of whose edges
lie in strata of positive dimension. This implies that
$R\Gamma(\Delta;\CM|_{\Delta})\in \CCD^{\leq -1}(\{ pt\})$, and one concludes
by Mayer-Vietoris argument, using similar estimates for intersections
of simplices.

Consequently we have
$$
\Phi_{F}(\CM)\cong R\Gamma(\BA,S_F;\CM)\in\CCD^{\leq 0}(\{ pt\}),
$$
as was claimed.

On the other hand by Duality theorem ~\ref{odin dual} we have an opposite
inequality, which proves that $\Phi_F$ is exact in our case.

The case of an arbitrary facet is reduced immediately to the central one by
noting that an operation of the restriction to a transversal slice
composed with a shift by its codimension is $t$-exact and using
{}~(\ref{odin sliceiso}). The theorem is proved. $\Box$

\subsection{} By the above theorem, the restriction of functors
$\Phi_F$ to the abelian subcategory $\CM(\BA,\CS_{\CH})$ lands
in subcategory $\Vect\subset\CCD(\{ pt\})$.

In other words, we get exact functors
\begin{equation}
\label{odin phi}
\Phi_F:\CM(\BA,\CS_{\CH})\lra\Vect
\end{equation}
These functors commute with Verdier duality.

We will also use the notation $\CM_F$ for $\Phi_F(\CM)$.

\subsection{Canonical and variation maps} Suppose we have a facet $E$.
Let us denote by $\Fac^1(E)$ the set of all facets $F$ such that
$E<F$, $\dim\ F=\dim\ E+1$. We have
\begin{equation}
\label{odin unican}
S_E=\bigcup_{F\in\Fac^1(E)}D_F
\end{equation}

Suppose we have $\CK\in\CCD(\BA,\CS_{\CH})$.

\subsubsection{}
\label{odin addcan} {\bf Lemma.} {\em We have a natural isomorphism
$$
R\Gamma(S_E,\bigcup_{F\in\Fac^1(E)}S_F;\CK)\cong
\oplus_{F\in\Fac^1{E}}R\Gamma(D_F,S_F;\CK)
$$}

{\bf Proof.} Note that $S_E-\bigcup_{F\in\Fac^1(E)}S_F=
\bigcup_{F\in\Fac^1(E)}\DO_F$ (disjoint union). The claim follows now from the
Poincar\'{e} duality. $\Box$

Therefore, for any $F\in\Fac^1(E)$ we get a natural inclusion map
\begin{equation}
\label{odin projcan}
i^F_E: R\Gamma(D_F,S_F;\CK)\hra R\Gamma(S_E,\bigcup_{F'\in\Fac^1(E)}S_{F'};\CK)
\end{equation}

Let us define a map
$$
u^F_E(\CK):\Phi_F(\CK)\lra\Phi_E(\CK)
$$
as a composition
$$
R\Gamma(D_F,S_F;\CK)[-p]\overset{i^F_E}{\lra}
R\Gamma(S_E,\bigcup_{F'\in\Fac^1(E)}S_{F'};\CK)[-p]\lra
R\Gamma(S_E;\CK)[-p]\lra
R\Gamma(D_E,S_E)[-p+1]
$$
where the last arrow is the coboundary map for the couple $(S_E,D_E)$,
and the second one is evident.

This way we get a natural transormation
\begin{equation}
\label{odin can}
u^F_E:\Phi_F\lra \Phi_E
\end{equation}
which will be called a {\em canonical map}.

We define a {\em variation map}
\begin{equation}
\label{odin var}
v^E_F:\Phi_E\lra\Phi_F
\end{equation}
as follows. By definition, $v^E_F(\CK)$ is the map dual to
the composition
$$
D\Phi_F(\CK)\iso\Phi_F(D\CK)\overset{u^F_E(D\CK)}{\lra}\Phi_E(D\CK)
\iso D\Phi_E(\CK).
$$

\subsection{Lemma}
\label{odin trans} {\em Suppose we have $4$ facets $A,B_1,B_2,C$ such that
$A<B_1<C,\ A<B_2<C$ and $\dim\ A=\dim\ B_i-1=\dim\ C-2$ (see Fig. 2). Then
$$
u^{B_1}_A\circ u^C_{B_1}=-u^{B_2}_A\circ u^C_{B_2}
$$
and
$$
v^{B_1}_C\circ v^A_{B_1}=-v^{B_2}_C\circ v^A_{B_2}.
$$}

For a proof, see below, ~\ref{odin prooftr}.

\begin{picture}(20,8)(-10,-4)

\put(0,-2){\line(1,1){4}}
\put(2.5,0){$B_2$}
\put(0,-2){\line(-1,1){4}}
\put(-3,0){$B_1$}

\put(0,-2){\circle*{0.2}}
\put(0,-2.7){$A$}

\put(0,3){$C$}

\put(-0.5,-4){Fig. 2}

\end{picture}

\subsection{Cochain complexes}
For each integer $p$, $0\leq p\leq N$, and $\CM\in\CM(\BA,\CS_{\CH})$
introduce vector spaces
\begin{equation}
C^{-p}_{\CH}(\BA;\CM)=\oplus_{F:\dim F=p}\ \CM_F
\end{equation}
For $i>0$ or $i<-N$ set $C^i_{\CH}(\BA;\CM)=0$.

Define operators
$$
d: C^{-p}_{\CH}(\BA;\CM)\lra C^{-p+1}_{\CH}(\BA;\CM)
$$
having components $u^F_E$.

\subsubsection{}
\label{odin nilp} {\bf Lemma.} {\em $d^2=0$}.

{\bf Proof.} Let us denote $X_p:=\bigcup_{F:\dim F=p}D_F$. We have
$X_p\supset X_{p+1}$. Evident embeddings
of couples $(D_F,S_F)\hra (X_p,X_{p+1})$ induce maps
$$
R\Gamma(X_p,X_{p+1}; \CM)\lra \oplus_{F:\dim F=p}\ R\Gamma(D_F,S_F;\CM)
$$
which are easily seen to be isomorphisms. Thus, we can identify
$C^{-p}(\BA;\CM)$ with $R\Gamma(X_p,X_{p+1};\CM)$. In these terms,
$d$ is a boundary homomorphism for the triple $X_p\supset X_{p+1}
\supset X_{p+2}$. After this description, the equality $d^2=0$ is a general
fact from homological algebra. $\Box$

\subsubsection{}
\label{odin prooftr} {\bf Proof of ~\ref{odin trans}.} The above lemma is equivalent
to the statement of ~\ref{odin trans} about maps $u$, which is thus proven.
The claim for variation maps follows by duality. $\Box$

This way we get a complex $C^{\bullet}_{\CH}(\BA;\CM)$ lying in
degrees from $-N$ to $0$. It will be called the {\em cochain
complex} of our arrangement $\CH_{\BR}$ with coefficients in $\CM$.

\subsection{Theorem}
\label{odin rgamma} {\em (i) A functor
$$
\CM\mapsto C^{\bullet}_{\CH}(\BA;\CM)
$$
is an exact functor from $\CM(\BA;\CS_{\CH})$ to the category of complexes
of vector spaces.

(ii) We have a canonical natural isomorphism in
$\CCD(\{ pt\})$
$$
C^{\bullet}_{\CH}(\BA;\CM)\iso R\Gamma(\BA;\CM)
$$}

{\bf Proof.} (i) is obvious from the exactness of functors $\Phi_F$,
cf. Thm. ~\ref{odin exact}.
To prove (ii), let us consider the filtration
$$
\BA\supset X_0\supset X_1\supset\ldots X_N\supset 0.
$$
It follows easily from homotopy argument (cf. ~\ref{odin limit}, ~\ref{odin conic})
that the restriction
$$
R\Gamma(\BA;\CM)\lra R\Gamma(X_0;\CM)
$$
is an isomorphism. On the other hand, a "Cousin" interpretation of
$C^{\bullet}_{\CH}(\BA;\CM)$ given in the proof of Lemma ~\ref{odin nilp},
shows that one has a canonical isomorphism
$R\Gamma(X_0;\CM)\iso C^{\bullet}_{\CH}(\BA;\CM)$. $\Box$

\section{Computations for standard sheaves}

\subsection{}
\label{odin poinc} Suppose we have a connected locally simply connected
topological space $X$ and a subspace $Y\subset X$ such that each connected
component of $Y$ is simply connected.
Recall that a {\em groupoid} is a category all of whose morphisms are
isomorphisms. Let us define a {\em Poincar\'{e} groupoid}
 $\pi_1(X;Y)$ as follows.

We set $\Ob \pi_1(X;Y)=\pi_0(Y)$. To define morphisms, let us choose
a point $y_i$ on each connected component $Y_i\subset Y$. By definition,
for two connected components $Y_i$ and $Y_j$, the set of homomorphisms
$\Hom_{\pi_1(X,Y)}(Y_i,Y_j)$ is the set of all homotopy classes of paths
in $X$ starting at $y_i$ and ending at $y_j$.

A different choice of points $y_i$ gives a canonically isomorphic
groupoid. If $Y$ is reduced to one point we come back to a usual definition of
the fundamental group.

Given a local system $\CL$ on $X$, we may assign to it a "fiber" functor
$$
F_{\CL}:\pi_1(X;Y)\lra\Vect,
$$
carrying $Y_i$ to the fiber $\CL_{y_i}$. This way we get an equivalence
of the category of local systems on $X$ and the category of functors
$\pi_1(X,Y)\lra\Vect$.

\subsection{} Return to the situation of the previous section. It is known
(cf. ~\cite{br}) that the homology group $H_1(\BAO_{\CH};\BZ)$ is a free
abelian group with a basis consisting of classes of small loops around
hyperplanes $H\in\CH$. Consequently, for each map
\begin{equation}
\label{odin monodr}
\bq: \CH\lra\BC^*
\end{equation}
there exists a one-dimensional local system $\CL(\bq)$ whose monodromy
around $H\in \CH$ is equal to $\bq(H)$. Such a
local system is unique up to a non-unique isomorphism.

Let us construct such a local system explicitely, using a language of the
previous subsection.

\subsection{} From now on we fix
a real equation for each $H\in\CH$, i.e. a linear function
$f_H:\BA_{\BR}\lra\BR$
such that $H_{\BR}=f^{-1}(0)$. We will denote also by $f_H$ the induced
function $\BA\lra\BC$.

The hyperplane $H_{\BR}$ divides $\BA_{\BR}$ into two halfspaces:
$\BA^+_{\BR,H}=\{ x\in\BA_{\BR}|f_H(x)>0\}$ and
{}~{$\BA^-_{\BR,H}=\{ x\in\BA_{\BR}|f_H(x)<0\}$}.

Let $F\subset H_{\BR}$ be a facet of dimension $N-1$.
We have two chambers $F_{\pm}$ adjacent to $F$, where
$F_{\pm}\subset\BA_{\BR,H}^{\pm}$. Pick a point $w\in F$.
Let us choose a real affine line $l_{\BR}\subset\BA_{\BR}$ transversal
to $H_{\BR}$ and passing through $w$. Let $l$ denote its complexification.

The function $f_H$ induces isomorphism $l\iso\BC$, and
$f_H^{-1}(\BR)\cap l=l_{\BR}$. Let us pick a real $\epsilon>0$ such that two
points
$f_H^{-1}(\pm\epsilon)\cap l_{\BR}$ lie in $F_{\pm}$ respectively. Denote
these points by $w_{\pm}$.

\begin{picture}(20,8)(-10,-4)

\put(-1,-3){\line(1,1){5}}
\put(-1.5,-3.5){$H_{\BR}$}
\put(0.3,-1.2){$F$}

\put(1.5,-0.5){\circle*{0.2}}
\put(1.2,-1){$w$}

\put(1.5,-0.5){\line(1,-1){1}}
\put(2,-1){\circle*{0.2}}
\put(2.3,-1){$w_+$}

\put(1.5,-0.5){\line(-1,1){1}}
\put(1,0){\circle*{0.2}}
\put(0.2,0){$w_-$}

\put(1,-3){\line(-1,1){5}}

\put(-1,0){$F_-$}
\put(3,0){$F_+$}

\put(-4,1){\line(1,0){8}}
\put(-3,1){\circle*{0.2}}
\put(3,1){\circle*{0.2}}

\put(0,-2){\circle*{0.2}}

\put(-0.5,-4){Fig. 3}

\end{picture}

Let us denote by  $\tau^+$ (resp., $\tau^-$) a counterclockwise
(resp., clockwise) path in the upper (resp., lower) halfplane connecting
$\epsilon$ with $-\epsilon$. Let us denote
$$
\tau^{\pm}_{F}=f_H^{-1}(\tau^{\pm})
$$
This way we get two well-defined homotopy classes of paths
connecting chambers $F_+$ and $F_-$. The argument $\arg\ f_H$
increases by $\mp\frac{\pi}{2}$ along $\tau^{\pm}_F$.

Note that if $H'$ is any other hyperplane of our arrangement then
$\arg f_{H'}$ gets no increase along $\tau^{\pm}_F$.

\subsection{} Now suppose we have $\bq$ as in ~(\ref{odin monodr}).
Note that all connected components of $\BAO_{\CH,\BR}$ --- chambers
of our arrangement --- are contractible.
Let us define a functor
$$
F(\bq^2):\pi_1(\BAO_{\CH},\BAO_{\CH,\BR})\lra\Vect_{\BC}
$$
as follows. For each chamber $C$ we set $F(\bq^2)(C)=\BC$.
For each facet $F$ of codimension $1$ which lies in a hyperplane $H$
we set
$$
F(\bq^2)(\tau^{\pm}_F)=\bq(H)^{\pm 1}
$$
It follows from the above remark on the structure of $H_1(\BAO_{\CH})$ that
we get a correctly defined functor.

The corresponding abelian local system over $\BA_{\CH}$ will be denoted
$\CL(\bq^2)$; it has a monodromy $\bq (H)^2$ around $H\in\CH$.
If all numbers $\bq(H)$ belong to some subfield $\sk\subset\BC$ then
the same construction gives a local system of $\sk$-vector spaces.

\subsection{} From now on until the end of this section we fix a map
$$
\bq:\CH\lra \sk^*,
$$
$\sk$ being a subfield of $\BC$, and denote by $\CL$ the local system
of $\sk$-vector spaces $\CL(\bq^2)$ constructed above. We denote by
$\CL^{-1}$ the dual local system $\CL(\bq^{-2})$.

Let $j:\BAO_{\CH}\lra\BA$ denote an open embedding. For $?$ equal to
one of the symbols $!,*$ or $!*$,
let us consider perverse sheaves $\CL_{?}:=j_?\CL[N]$. They belong to
$\CM(\BA;\CS_{\CH})$; these sheaves will be called {\em standard extensions}
of $\CL$, or simply {\em standard sheaves}. Note that
\begin{equation}
\label{odin dualstand}
D\CL_!\cong\CL^{-1}_*;\ D\CL_*\cong\CL^{-1}_!
\end{equation}
We have a canonical map
\begin{equation}
\label{odin map!*}
m: \CL_!\lra\CL_*,
\end{equation}
and by definition $\CL_{!*}$ coincides with its image.

Our aim in this section will be to compute explicitely the cochain
complexes of standard sheaves.

\subsection{Orientations}
\label{odin coor} Let $F$ be a facet which is a
connected component of a real edge $L_F$. Consider a linear space
$L^{\bot}_F=\BA_{\BR}/L_F$. Let us define
$$
\lambda_F:=H^0(L^{\bot}_F;\Or_{L^{\bot}_F}),
$$
it is a free abelian group of rank $1$. To choose an orientation
of $L^{\bot}_F$ (as a real vector space) is the same as to choose
a basis vector in $\lambda_F$.
We will call an orientation of $L^{\bot}_F$ a {\em coorientation}
of $F$.

We have an evident piecewise linear homeomorphism of
$D_F=\bigcup D_{F<C}$ onto a closed disk in $L^{\bot}_F$; thus, a
coorientation of $F$ is the same as an orientation of $D_F$ (as a $C^0$-
manifold); it defines orientations of all cells $D_{F<C}$.

\subsubsection{}
\label{odin sign} From now on until the end of the section, let us fix
coorientations of all facets.
Suppose we have a pair $E<F$, $\dim\ E=\dim\ F-1$. The cell $D_E$ is a part
of the boundary of $D_F$. Let us define the sign $\sgn(F,E)=\pm 1$ as follows.
Complete an orienting basis of $D_E$ by a vector directed outside
$D_F$; if we get the given orientation of $D_F$, set $\sgn(F,E)=1$,
otherwise set $\sgn(F,E)=-1$.

\subsection{Basis in $\Phi_F(\CL_!)^*$}

Let $F$ be a facet of dimension $p$. We have by definition
$$
\Phi_F(\CL_!)=H^{-p}(D_F,S_F;\CL_!)=H^{N-p}(D_F,S_F;j_!\CL)\cong
H^{N-p}(D_F,S_F\cup (_{\CH}H_{\BR}\cap D_F);j_!\CL).
$$
By Poincar\'{e} duality,
$$
H^{N-p}(D_F,S_F\cup (_{\CH}H_{\BR}\cap D_F);j_!\CL)^*\cong
H^0(D_F-(S_F\cup _{\CH}H_{\BR});\CL^{-1})
$$
(recall that we have fixed an orientation of $D_F$).
The space $D_F-(S_F\cup _{\CH}H_{\BR})$ is a disjoint union
$$
D_F-(S_F\cup _{\CH}H_{\BR})=\bigcup_{C\in \Ch (F)}\DO_{F<C}.
$$
Consequently,
$$
H^0(D_F-(S_F\cup _{\CH}H_{\BR});\CL^{-1})\cong
\oplus_{C\in \Ch (F)}H^0(\DO_{F<C};\CL^{-1}).
$$
By definition of $\CL$, we have canonical identifications
$H^0(\DO_{F<C};\CL^{-1})=\sk$. We will denote by $c(\CL_!)_{F<C}\in\Phi_F(\CL)^*$
the
image of $1\in H^0(\DO_{F<C};\CL^{-1})$ with respect to the embedding
$$
H^0(\DO_{F<C};\CL^{-1})\hra \Phi_F(\CL_!)^*
$$
following from the above.

Thus, classes $c(\CL_!)_{F<C},\ C\in \Ch (F)$, form a basis of
$\Phi_F(\CL_!)^*$.

\subsection{}
\label{odin can!} Let us describe canonical maps for
$\CL_!$. If $F<E,\ \dim\ E=\dim\ F+1$, let
$u^*:\Phi_F(\CL_!)^*\lra\Phi_E(\CL_!)^*$ denote the map dual to $u^E_F(\CL_!)$.
Let $C$ be a chamber adjacent to $F$.
Then
\begin{equation}
u^*(c(\CL_!)_{F<C})=\left\{ \begin{array}{ll}
                        \sgn(F,E)c(\CL_!)_{E<C}&\mbox{if $E<C$}\\
                          0&\mbox{otherwise}
                     \end{array}
             \right.
\end{equation}

\subsection{Basis in $\Phi_F(\CL_*)^*$} We have isomorphisms
\begin{equation}
\label{odin dualchains}
\Phi_F(\CL_*)\cong\Phi_F(D\CL_!^{-1})\cong\Phi_F(\CL_!^{-1})^*
\end{equation}
Hence, the defined above basis $\{ c(\CL^{-1}_!)_{F<C}\}_{C\in\Ch (F)}$
of $\Phi_F(\CL_!^{-1})^*$, gives a basis in  $\Phi_F(\CL_*)$.
We will denote by $\{ c(\CL_*)_{F<C}\}_{C\in\Ch (F)}$ the dual basis
of $\Phi_F(\CL_*)$.

\subsection{Example} Let us describe our chains explicitely in the simplest
one-dimensional case, in the setup ~\ref{odin one-dim}. We choose a natural
orientation on $\BA_{\BR}$.
A local system $\CL=\CL(q^2)$ is uniquely determined by one nonzero
complex number $q$. By definition, the upper (resp., lower) halfplane
halfmonodromy from  $w_+$ to $w_-$ is equal to $q$ (resp., $q^{-1}$).

\subsubsection{Basis in $\Phi_F(\CL_!)^*$} The space $\Phi_F(\CL_!)^*$
admits a basis consisting of two chains
$c_{\pm}=c(\CL_!)_{F<E_{\pm}}$ shown below, see Fig. 3(a).
By definition, a homology class is represented by a cell together with
a section of a local system $\CL^{-1}$ over it. The section of $\CL^{-1}$
over $c_+$ (resp., $c_-$) takes value $1$ over $w_+$ (resp., $w_-$).


\begin{picture}(20,8)(-10,-4)

\put(-5,0){\circle{0.2}}
\put(-5,0){\oval(6,6)}
\put(-6.5,3.3){$S_{r''}$}

\put(-3,0){\circle*{0.2}}
\put(-3,-0.5){$w_+$}
\put(-7,0){\circle*{0.2}}
\put(-7,-0.5){$w_-$}

\put(-7,0){\line(1,0){4}}
\put(-7,0){\vector(1,0){1.5}}
\put(-5,0){\vector(1,0){1.5}}

\put(-6.7,0.3){$c_-$}
\put(-3.5,0.3){$c_+$}


\put(-6,-3){\line(0,1){6}}
\put(-6,1){\vector(0,1){1}}
\put(-6,-3){\circle*{0.2}}
\put(-6,3){\circle*{0.2}}
\put(-5.8,2){$Dc_-$}

\put(-4,-3){\line(0,1){6}}
\put(-4,1){\vector(0,1){1}}
\put(-4,-3){\circle*{0.2}}
\put(-4,3){\circle*{0.2}}
\put(-3.8,2){$Dc_+$}

\put(-5.5,-3.5){(a)}


\put(5,0){\circle{0.2}}
\put(5,0){\oval(6,6)}
\put(3.5,3.3){$S_{r''}$}

\put(5,3){\line(0,-1){2}}
\put(5,-3){\line(0,1){2}}
\put(4.6,1.8){$Y$}

\put(5,1){\circle*{0.2}}
\put(5,-1){\circle*{0.2}}

\put(7,0){\circle*{0.2}}
\put(7,-0.5){$w_+$}
\put(3,0){\circle*{0.2}}
\put(3,-0.5){$w_-$}

\put(3,0){\line(1,0){4}}
\put(3,0){\vector(1,0){1.5}}
\put(5,0){\vector(1,0){1.5}}

\put(-6.7,0.3){$c_-$}
\put(-3.5,0.3){$c_+$}

\put(5,0){\oval(2,3)}
\put(6,0){\vector(0,1){0.5}}
\put(6.2,0.5){$\tilde{D}c_+$}
\put(4,0){\vector(0,1){0.5}}
\put(3,0.5){$\tilde{D}c_-$}
\put(5,1.5){\circle*{0.2}}
\put(5,-1.5){\circle*{0.2}}

\put(4.5,-3.5){(b)}

\put(-0.5,-4){Fig. 3}

\end{picture}

\subsubsection{Basis in $\Phi_F(\CL_*)^*$}
Let us adopt notations of ~\ref{odin mainconst}, with $\CK=\CL^{-1}_!$.
It is easy to find the basis $\{ Dc_+,Dc_-\}$ of the dual space
$$
\Phi_F(\CL_!^{-1})\cong H^0(\BA,\BA_{\geq r''};j^*j_*D\CL_!^{-1})^*=
H^1(\BA-\{ w_+,w_-\},\BA_{\geq r''};\CL)^*
$$
dual to $\{ c(\CL^{-1})_{F<E_+},c(\CL^{-1})_{F<E_-}\}$.
Namely, $Dc_{\pm}$ is represented by the relative $1$-chain
$$
\{ \pm\frac{1}{2}+y\cdot i|\ -\sqrt{(r'')^2-\frac{1}{4}}\leq y
\leq \sqrt{(r'')^2-\frac{1}{4}}\},
$$
with evident sections of $\CL^{-1}$ over them, see Fig. 3(a).

Next, one has to deform these chains to chains $\tilde{D}c_{\pm}$ with
their ends on $Y$, as in Fig. 3(b). Finally, one has to make a clockwise
rotation of the picture by $\pi/2$. As a result, we arrive at the
following two chains $c^*_{\pm}$ forming a basis of $\Phi_F(\CL_*)^*$:


\begin{picture}(20,8)(-10,-4)

\put(0,0){\circle{0.2}}

\put(2,0){\circle*{0.2}}
\put(2.2,-0.5){$w_+$}
\put(-2,0){\circle*{0.2}}
\put(-1.8,-0.5){$w_-$}

\put(0,0){\oval(4,4)}
\put(0,2){\vector(1,0){0.5}}
\put(0,2.3){$c^*_-$}
\put(0,-2){\vector(1,0){0.5}}
\put(0,-2.5){$c^*_+$}

\put(-2,0){\line(1,0){4}}
\put(-2,0){\vector(1,0){1.5}}
\put(0,0){\vector(1,0){1.5}}

\put(-1,0.3){$c_-$}
\put(0.8,0.3){$c_+$}

\put(-0.5,-4){Fig. 4}

\end{picture}

The section of $\CL^{-1}$ over $c^*_+$ (resp., $c_-$) has value $1$ at
$w_+$ (resp., $w_-$).

It follows from this description that the natural map
\begin{equation}
\label{odin m!*}
m:\Phi_F(\CL_*)^*\lra\Phi_F(\CL_!)^*
\end{equation}
is given by the formulas
\begin{equation}
\label{odin f!*}
m(c^*_+)=c_++qc_-;\ m(c^*_-)=qc_++c_-
\end{equation}

By definition, spaces $\Phi_{E_{\pm}}(\CL_?)^*$ may be identified
with fibers $\CL_{w_{\pm}}$ respectively, for both $?=!$ and $?=*$,
and hence with $\sk$. Let us denote by $c_{w_{\pm}}$ and
$c^*_{w_{\pm}}$ the generators corresponding to $1\in \sk$.

It follows from the above description that the canonical maps $u^*$ are given
by the formulas
\begin{equation}
\label{odin fu!}
u^*(c_+)=c_{w_+};\ u^*(c_-)=-c_{w_-};
\end{equation}
and
\begin{equation}
\label{odin fu*}
u^*(c_+^*)=c_{w_+}^*-qc_{w_-}^*;\
u^*(c_-^*)=qc_{w_+}^*-c_{w_-}^*;\
\end{equation}
Let us compute variation maps. To get them for the sheaf $\CL_!$, we should
by definition replace $q$ by $q^{-1}$ in ~(\ref{odin fu*}) and take the conjugate
map:
\begin{equation}
\label{odin fv!}
v^*(c_{w_+})=c_+ + q^{-1}c_-;\
v^*(c_{w_-})=-q^{-1}c_+-c_-
\end{equation}
To compute $v^*$ for $\CL_*$, note that the basis in
$$
H^0(\BA,\{ w_+,w_-\};\CL_!^{-1})^*=H_1(\BA,\{w_+,w_-,0\};\CL)
$$
dual to $\{ c^*_+,c^*_-\}$, is $\{ q^{-1}\tc_-,q^{-1}\tc_+\}$ {\em (sic!)}
where $\tc_{\pm}$ denote the chains defined in the same way as $c_{\pm}$,
with $\CL$ replaced by $\CL^{-1}$. From this remark it follows that
\begin{equation}
\label{odin fv*}
v^*(c_{w_+}^*)=q^{-1}c_-^*;\
v^*(c_{w_-}^*)=-q^{-1}c_+^*
\end{equation}

\subsection{} Let us return to the case of an arbitrary arrangement.
Let us say that a hyperplane $H\in\CH$ {\em separates} two chambers $C, C'$
if they lie in different halfspaces with respect to $H_{\BR}$.
Let us define numbers
\begin{equation}
\label{odin qsep}
\bq(C,C')=\prod\bq(H),
\end{equation}
the product over all hyperplanes $H\in \CH$ separating $C$ and $C'$.
In particular, $\bq(C,C)=1$.

\subsection{Lemma.} {\em Let $F$ be a face, $C\in\Ch(F)$. The canonical mapping
$$
m:\Phi_F(\CL_*)^*\lra\Phi_F(\CL_!)^*
$$
is given by the formula
\begin{equation}
\label{odin sform}
m(c(\CL_*)_{F<C})=\sum_{C'\in\Ch(F)}\bq(C,C')c(\CL_!)_{F<C'}
\end{equation}}

\subsubsection{} Since $\Phi_F(\CL_*)$ is dual to $\Phi_F(\CL_!)$, we may
view $m$ as  a {\em bilinear form} on $\Phi_F(\CL_!)$. By ~(\ref{odin sform})
it is {\em symmetric}.

{\bf Proof} of lemma. We generalize the argument of the previous example.
First consider the case of zero-dimensional $F$.
Given a chain $c(\CL_!^{-1})_{F<C}$, the corresponding dual chain may be
taken as
$$
Dc_{F<C}= \epsilon\cdot\ ^Cw\oplus i\cdot\BA_{\BR},
$$
were $\epsilon$ is a sufficiently small positive real.
Next, to get the dual chain $c(\CL_*)_{F<C}$, we should make a deformation
similar to the above one, and a rotation by $\frac{\pi}{2}$. It is convenient
to make the rotation first. After the rotation, we get a chain
$\BA_{\BR} -\epsilon i\cdot\ ^Cw$. The value of
$m(c(\CL_*)_{F<C})$ is given by the projection of this chain to $\BA_{\BR}$.

The coefficient at $c(\CL_!)_{F<C'}$ is given by the monodromy of $\CL^{-1}$
along the following path from $C$ to $C'$. First, go "down" from $^Cw$
to $-\epsilon i\cdot\ ^Cw$; next, travel in
$\BA_{\BR} -\epsilon i\cdot\ ^Cw$ along the straight line from
$-\epsilon i\cdot\ ^Cw$ to $-\epsilon i\cdot ^{C'}w$, and then go "up" to
$^{C'}w$. Each time we are passing under a hyperplane $H_{\BR}$
separating $C$ and $C'$, we gain a factor $\bq(H)$. This gives
desired coefficient for the case $\dim\ F=0$.

For an arbitrary $F$ we use the same argument by considering the intersection
of our picture with a transversal slice.   $\Box$

\subsection{} Let $E$ be a facet which is a component of a
real edge $L_{E,\BR}$; as usually $L$ will denote the complexification.
Let $\CH_L\subset \CH$ be a subset consisting of all hyperplanes
containg $L$. If we assign to a chamber $C\in\Ch(E)$ a unique chamber of the
subarrangement $\CH_L$ comtaining $C$, we get a bijection of $\Ch(E)$ with
the set of {\em all} chambers of $\CH_L$.

Let $F<E$ be another facet. Each chamber $C\in\Ch(F)$ is contained in a unique
chamber of $\CH_L$. Taking into account a previous bijection, we get a mapping
\begin{equation}
\label{odin pi}
\pi^F_E:\Ch(F)\lra\Ch(E)
\end{equation}

\subsection{Lemma.}{\em Let $F$ be a facet, $C\in\Ch(F)$. We have
\begin{equation}
\label{odin u*}
u^*(c(\CL_*)_{F<C})=\sum\sgn(F,E)\bq(C,\pi^F_E(C))c(\CL_*)_{E<\pi^F_E(C)},
\end{equation}
the summation over all facets $E$ such that $F<E$ and $\dim\ E=\dim\ F+1$.}

(Signs $\sgn(F,E)$ have been defined in ~\ref{odin sign}.)

{\bf Proof.} Again, the crucial case is $\dim\ F=0$ --- the case of arbitrary
dimension is treated using a transversal slice. So, let us suppose that
$F$ is zero-dimensional. In order to compute the coefficient of
$u^*(c(\CL_*)_{F<C})$ at $c(\CL^*)_{E<C'}$ where $E$ is a one-dimensional
facet adjacent to $F$ and $C'\in\Ch(E)$, we have to do the following.

Consider the intersection of a real affine subspace
$\BA_{\BR}-\epsilon\cdot i\cdot\ ^Cw$ (as in the proof of the previous lemma)
with a complex hyperplane $M_E$ passing through $^Ew$ and transversal
to $E$. The intersection will be homotopic to a certain chain
$c(\CL_*)_{E<C''}$ where $C''$ is easily seen to be equal to $\pi^F_E(C)$,
and the coefficient is obtained by the same rule as described
in the previous proof. The sign will appear in accordance with
compatibility of orientations of $D_F$ and $D_E$. $\Box$

\subsection{} Let us collect our results. Let us denote by
$\{ b(\CL_?)_{F<C}\}_{C\in\Ch(F)}$ the basis in $\Phi_F(\CL_?)$ dual
to $\{ c(\CL_?)_{F<C}\}$, where $?=!$ or $*$.

\subsection{Theorem.}{\em (i) The complex $C^{\bullet}_{\CH}(\BA;\CL_!)$
is described as follows. For each $p,\ 0\leq p\leq N$, the space
$C^{-p}_{\CH}(\BA;\CL_!)$ admits a basis consisting of all cochains
$b(\CL_!)_{F<C}$ where $F$ runs through all facets of $\CH_{\BR}$ of dimension
$p$, and $C$ through $\Ch(F)$. The differential
$$
d:C^{-p}_{\CH}(\BA;\CL_!)\lra C^{-p+1}_{\CH}(\BA;\CL_!)
$$
is given by the formula
\begin{equation}
\label{odin dl!}
d(b(L_!)_{F<C})=\sum_{E:E<F,\ \dim\ E=\dim\ F-1} \sgn(E,F)b(L_!)_{E<C}
\end{equation}

(ii) The complex $C^{\bullet}_{\CH}(\BA;\CL_*)$
is described as follows. For each $p,\ 0\leq p\leq N$, the space
$C^{-p}_{\CH}(\BA;\CL_*)$ admits a basis consisting of all cochains
$b(\CL_*)_{F<C}$ where $F$ runs through all facets of $\CH_{\BR}$ of dimension
$p$, and $C$ through $\Ch(F)$. The differential
$$
d:C^{-p}_{\CH}(\BA;\CL_*)\lra C^{-p+1}_{\CH}(\BA;\CL_*)
$$
is given by the formula
\begin{equation}
\label{odin dl*}
d(b(L_*)_{F<C})=\sum \sgn(E,F)\bq(C,C')b(L_*)_{E<C'},
\end{equation}
the summation over all facets $E<F$ such that $\dim\ E=\dim\ F-1$ and
all chambers $C'\in\Ch(E)$ such that $\pi^E_F(C')=C$.

(iii) The natural map of complexes
\begin{equation}
\label{odin mapc}
m:C^{\bullet}_{\CH}(\BA;\CL_!)\lra C^{\bullet}_{\CH}(\BA;\CL_*)
\end{equation}
induced by the canonical map $\CL_!\lra\CL_*$, is given by the formula
\begin{equation}
\label{odin formb}
m(b(\CL_!)_{F<C})=\sum_{C'\in\Ch(F)}\bq(C,C')b(\CL_*)_{F<C'}
\end{equation}}

All statements have already been proven.

\subsubsection{} {\bf Corollary.} {\em The complexes
$C^{\bullet}_{\CH}(\BA;\CL_!)$
and $C^{\bullet}_{\CH}(\BA;\CL_*)$ described explicitely in the above theorem,
compute the relative cohomology $H^{\bullet}(\BA,\ _{\CH}H;\CL)$ and the
cohomology of the open stratum $H^{\bullet}(\BAO,\CL)$ respectively,
and the map $m$ induces the canonical map in cohomology.}

{\bf Proof.} This follows immediately from ~\ref{odin rgamma}. $\Box$

This corollary was proven in ~\cite{v}, Sec. 2, by a different argument.

\subsection{Theorem.}
\label{odin inters} {\em The complex $C^{\bullet}_{\CH}(\BA;\CL_{!*})$
is canonically isomorphic to the image of ~(\ref{odin mapc}).}

{\bf Proof.} This follows from the previous theorem and the exactness
of the functor $\CM\mapsto C^{\bullet}_{\CH}(\BA;\CM)$,
cf. Thm. ~\ref{odin rgamma} (i). $\Box$

The above description of cohomology is analogous to ~\cite{sv1}, p. I,
whose results may be considered as a "quasiclassical" version
of the above computations.

\newpage
\setcounter{section}{0}
\begin{center}{\large \bf Part II. CONFIGURATION SPACES}\end{center}

\begin{center}{\large\bf AND QUANTUM GROUPS}\end{center}

\section{Introduction}

\subsection{} 
We are starting here the geometric study of the tensor category $\CC$
associated with a quantum group (corresponding to a
Cartan matrix of finite type) at a root of unity (see ~\cite{ajs}, 1.3 and
the present part, ~\ref{dva C} for the precise definitions).

The main results of this part are Theorems ~\ref{dva shaposym},
{}~\ref{dva stalks}, ~\ref{dva shsymnsl} and ~\ref{dva stalknsl} which

---  establish isomorphisms between homogeneous
components of irreducible objects in $\CC$ and spaces of vanishing
cycles at the origin of certain Goresky-MacPherson sheaves on
configuration spaces;

--- establish isomorphisms of the stalks at the origin
of the above GM sheaves with certain Hochschild complexes
(which compute the Hochschild homology of a certain "triangular"
subalgebra of our quantum group with coefficients in the coresponding
irreducible representation);

--- establish the analogous results for tensor products of irreducibles.
In geometry,
the tensor product of representations corresponds to a "fusion"
of sheaves on configuration spaces --- operation defined
using the functor of nearby cycles, see Section ~\ref{dva fus}.

We must mention that the assumption that we are dealing with a
Cartan matrix of finite type and a root of unity appears only at the very end
(see Chapter 4). We need these assumptions in order
to compare our representations with the conventional
definition of the category $\CC$. All previous results are valid in more
general assumptions. In particular a
Cartan matrix could be arbitrary and
a deformation parameter $\zeta$ not necessarily a root of unity.

\subsection{} Some of the results of this part constitute the description
of the cohomology of certain "standard" local systems over configuration
spaces in terms of quantum groups. These results, due to Varchenko
and one of us, were announced several years ago in ~\cite{sv2}.
The proofs may be found in ~\cite{v}.
Our proof of these results uses completely different approach.
Some close results were discussed in ~\cite{s}.

Certain results of a similar geometric spirit are discussed in
{}~\cite{fw}.

\subsection{} We are grateful to A.Shen who made our communication during
the writing of this part possible.

\subsection{}
\label{dva notations} {\em Notations.} We will use all the notations from
the part I. 
References to {\em loc. cit.} will look like I.1.1.
If $a,b$ are two integers, we will denote by $[a,b]$ the
set of all integers $c$ such that $a\leq c\leq b$; $[1,a]$ will be denoted
by $[a]$.
$\Bbb N$ will denote the set of non-negative integers.
For $r\in\BN$, $\Sigma_r$ will denote the group of all bijections $[r]\iso
[r]$.

We suppose that our ground field $\sk$ has characteristic $0$,
and fix an element $\zeta\in \sk$, $\zeta\neq 0$. For $a\in\BZ$ we will use
the notation
\begin{equation}
\label{dva azeta}
[a]_{\zeta}=1-\zeta^{-2a}
\end{equation}

The word "$t$-exact" will allways mean $t$-exactness with respect to
the middle perversity.

\newpage
\begin{center}
{\bf Chapter 1. Algebraic discussion}
\end{center}
\vspace{1cm}

\section{Free algebras and bilinear forms}

Most definitions of this section follow ~\cite{l1} and ~\cite{sv2}
(with slight modifications). We also add some new definitions and
computations important for the sequel. Cf. also ~\cite{v}, Section 4.

\subsection{} Until the end of this part, let us fix
a finite set $I$ and a symmetric $\BZ$-valued bilinear form $\nu,\nu'\mapsto
\nu\cdot\nu'$ on the free abelian group $\BZ[I]$ (cf. ~\cite{l1}, 1.1).
We will denote
by $X$ the dual abelian group $\Hom(\BZ[I],\BZ)$. Its elements will be
called {\em weights}. Given $\nu\in\BZ[I]$,
we will denote by $\lambda_{\nu}\in X$ the functional $i\mapsto
i\cdot \nu$. Thus we have
\begin{equation}
\label{dva lamnu}
\langle\lambda_{\nu},\mu\rangle=\nu\cdot\mu
\end{equation}
for all $\nu,\mu\in\BN[I]$.

\subsection{} Let ${\ '\ff}$ denote a free associative $\sk$-algebra with
$1$ with generators
$\theta_i,\ i\in I$. Let $\BN[I]$ be a submonoid of $\BZ[I]$ consisting
of all linear combinations of elements of $I$ with coefficients in $\BN$.
For $\nu=\sum\nu_ii\in\BN[I]$ we denote by ${\ '\ff}_{\nu}$ the $\sk$-subspace
of ${\ '\ff}$ spanned by all monomials $\theta_{i_1}\theta_{i_2}\cdot\ldots\cdot
\theta_{i_p}$ such that for any $i\in I$, the number of occurences of $i$
in the sequence $i_1,\ldots,i_p$ is equal to $\nu_i$.

We have a direct sum decomposition ${\ '\ff}=\oplus_{\nu\in\BN[I]}{\ '\ff}_{\nu}$,
all spaces ${\ '\ff}_{\nu}$ are finite dimensional, and we have ${\ '\ff}_0=\sk\cdot 1$,
${\ '\ff}_{\nu}\cdot{\ '\ff}_{\nu'}\subset{\ '\ff}_{\nu+\nu'}$.

Let $\epsilon:{\ '\ff}\lra \sk$ denote the augmentation ---
a unique $\sk$-algebra
map such that $\epsilon(1)=1$ and $\epsilon(\theta_i)=0$ for all $i$. Set
${\ '\ff}^+:=\Ker(\epsilon)$. We have ${\ '\ff}^+=\oplus_{\nu\neq 0}{\ '\ff}_{\nu}$.

An element $x\in{\ '\ff}$ is called {\em homogeneous} if it belongs to
${\ '\ff}_{\nu}$ for some $\nu$. We then set $|x|=\nu$.
We will use the notation $\depth(x)$ for the number $\sum_i\nu_i$
if $\nu=\sum_i\nu_ii$; it will be called {\em the depth} of $x$.

\subsection{}
\label{dva twist} Given a sequence
$\vec{K}=(i_1,\ldots,i_N),\ i_j\in I$, let us denote by $\theta_{\vec{K}}$
the monomial
$\theta_{i_1}\cdot\ldots\cdot\theta_{i_N}$. For an empty sequence we set
$\theta_{\emp}=1$.

For
$\tau\in\Sigma_N$
let us introduce the number
\begin{equation}
\label{dva zetaitau}
\zeta(\vec{K};\tau)=\prod\zeta^{i_a\cdot i_b},
\end{equation}
the product over all $a,b$ such that $1\leq a<b\leq N$ and $\tau(a)>\tau(b)$.

We will call this number {\em the twisting number of the sequence $\vec{K}$
with respect to the permutation $\tau$}.

We will use the notation
\begin{equation}
\label{dva permut}
\tau(\vec{K})=(i_{\tau(1)},i_{\tau(2)},\ldots,i_{\tau(N)})
\end{equation}

\subsection{}
\label{dva algff} Let us regard the tensor product ${\ '\ff}\otimes{\ '\ff}$ (in
the sequel $\otimes$ will mean $\otimes_\sk$ unless specified otherwise)
as a $\sk$-algebra with multiplication
\begin{equation}
\label{dva mult}
(x_1\otimes x_2)\cdot (x'_1\otimes x_2')=
\zeta^{|x_2|\cdot |x'_1|}x_1x'_1\otimes x_2x'_2
\end{equation}
for homogeneous $x_2,\ x'_1$. Let us define a map
\begin{equation}
\label{dva comult}
\Delta:{\ '\ff}\lra{\ '\ff}\otimes{\ '\ff}
\end{equation}
as a unique algebra homomorphism carrying $\theta_i$ to $\theta_i\otimes 1+
1\otimes\theta_i$.

\subsection{}
\label{dva coalgff} Let us define a coalgebra structure on
${\ '\ff}\otimes{\ '\ff}$ as
follows. Let us introduce the braiding isomorphism
\begin{equation}
\label{dva braid}
r:{\ '\ff}\otimes{\ '\ff}\iso{\ '\ff}\otimes{\ '\ff}
\end{equation}
by the rule
\begin{equation}
\label{dva braidform}
r(x\otimes y)=\zeta^{|x|\cdot |y|}y\otimes x
\end{equation}
for homogeneous $x,y$. By definition,
\begin{equation}
\label{dva deltaff}
\Delta_{{\ '\ff}\otimes{\ '\ff}}:{\ '\ff}\otimes
{\ '\ff}\lra({\ '\ff}\otimes{\ '\ff})\otimes({\ '\ff}\otimes{\ '\ff})
\end{equation}
coincides with the composition $(1_{{\ '\ff}}\otimes r\otimes 1_{{\ '\ff}})\circ
(\Delta_{{\ '\ff}}\otimes\Delta_{{\ '\ff}})$.

The multiplication
\begin{equation}
\label{dva multiplic}
{\ '\ff}\otimes{\ '\ff}\lra{\ '\ff}
\end{equation}
is a coalgebra morphism.

\subsection{}
\label{dva deltaform} Let us describe $\Delta$ more explicitely. Suppose a sequence
$\vec{K}=(i_1,\ldots,i_N),\ i_j\in I,$ is given. For a  subset
$A=\{ j_1,\ldots, j_a\}\subset [N],\ j_1<\ldots <j_a$, let
$A'=[N]-A=\{k_1,\ldots,k_{N-a}\},\ k_1<\ldots <k_{N-a}$. Define a permutation
$\tau_A$ by the formula
\begin{equation}
\label{dva taua}
(\tau(1),\ldots,\tau(N))=(j_1,j_2,\ldots,j_a,k_1,k_2,\ldots,k_{N-a})
\end{equation}
Set $\vec{K}_A:=(i_{j_1},i_{j_2},\ldots,i_{j_a}),\
\vec{K}_{A'}:=(i_{k_1},i_{k_2},\ldots,i_{k_{N-a}})$.

\subsubsection{}
\label{dva delta-expl} {\bf Lemma.} {\em
$$
\Delta(\theta_{\vec{K}})=\sum_{A\subset K}\zeta(\vec{K};\tau_A)
\theta_{\vec{K}_A}\otimes\theta_{\vec{K}_{A'}},
$$
the summation ranging over all subsets $A\subset [N]$.}

{\bf Proof} follows immediately from the definitions. $\Box$

\subsection{}
\label{dva powers} Let us denote by
\begin{equation}
\label{dva deltan}
\Delta^{(N)}:{\ '\ff}\lra{\ '\ff}^{\otimes N}
\end{equation}
iterated coproducts; by the coassociativity they are well defined.

Let us define a structure of an algebra on ${\ '\ff}^{\otimes N}$ as follows:
\begin{equation}
\label{dva product}
(x_1\otimes\ldots\otimes x_N)\cdot (y_1\otimes\ldots\otimes y_N)=
\zeta^{\sum_{j<i}|x_i|\cdot |y_j|} x_1y_1\otimes\ldots\otimes x_Ny_N
\end{equation}
for homogeneous $x_1,\ldots,x_N;y_1,\ldots,y_N$. The map $\Delta^{(N)}$ is an
algebra morphism.

\subsection{}
Suppose we have a sequence $\vec{K}=(i_1,\ldots,i_N)$.
Let us consider an element $\Delta^{(N)}(\theta_{\vec{K}})$;
let $\Delta^{(N)}(\theta_{\vec{K}})^+$ denote its projection to the
subspace ${\ '\ff}^{+\otimes N}$.

\subsubsection{}
\label{dva deltamon}
{\bf Lemma.} {\em
$$
\Delta^{(N)}(\theta_{\vec{K}})^+=\sum_{\tau\in\Sigma_N}\zeta(\vec{K};\tau)
\theta_{i_{\tau(1)}}\otimes\ldots\otimes\theta_{i_{\tau(N)}}
$$}

{\bf Proof} follows from ~\ref{dva delta-expl} by induction on $N$. $\Box$

\subsection{}
\label{dva duals}
For each component ${\ '\ff}_{\nu}$ consider the dual
$\sk$-space ${\ '\ff}^*_{\nu}$,
and set ${\ '\ff}^*:=\oplus{\ '\ff}^*_{\nu}$. Graded components
$\Delta_{\nu,\nu'}:{\ '\ff}_{\nu+\nu'}\lra{\ '\ff}_{\nu}\otimes
{\ '\ff}_{\nu'}$ define dual
maps ${\ '\ff}^*_{\nu}\otimes{\ '\ff}_{\nu'}^*\lra{\ '\ff}^*_{\nu+\nu'}$ which
give rise to a multiplication
\begin{equation}
\label{dva mult*}
{\ '\ff}^*\otimes{\ '\ff}^*\lra{\ '\ff}^*
\end{equation}
making ${\ '\ff}^*$ a graded associative algebra with $1$ (dual to the
augmentation of ${\ '\ff}$).
This follows from the coassociativity of $\Delta$, cf. ~\cite{l1}, 1.2.2.

Here and in the sequel, we will use identifications $(V\otimes W)^*=
V^*\otimes W^*$ (for finite dimensional spaces $V,W$) by the rule
$\langle\phi\otimes\psi,x\otimes y\rangle=\langle\phi,x\rangle\cdot
\langle\psi,y\rangle$.

The dual to ~(\ref{dva multiplic}) defines a comultiplication
\begin{equation}
\delta:{\ '\ff}^*\lra{\ '\ff}^*\otimes{\ '\ff}^*
\end{equation}
It makes ${\ '\ff}^*$ a graded coassociative coalgebra with a counit.

The constructions dual to ~\ref{dva algff} and ~\ref{dva coalgff}
equip ${\ '\ff}^*\otimes{\ '\ff}^*$ with a structure of a coalgebra and an algebra.
It follows from {\em loc. cit} that ~(\ref{dva mult*}) is a coalgebra
morphism, and $\delta$ is an algebra morphism.

By iterating $\delta$ we get maps
\begin{equation}
\label{dva deltait}
\delta^{(N)}:{\ '\ff}^*\lra{\ '\ff}^{*\otimes N}
\end{equation}
If we regard ${\ '\ff}^{*\otimes N}$
as an algebra by the same construction as in
{}~(\ref{dva product}), $\delta^{(N)}$ is an algebra morphism.

\subsection{Lemma.}
\label{dva form} {\em There exists a unique bilinear form
$$
S(\ ,\ ):{\ '\ff}\otimes{\ '\ff}\lra \sk
$$
such that

(a) $S(1,1)=1$ and $(\theta_i,\theta_j)=\delta_{i,j}$ for all $i,j\in I$;\\
(b) $S(x,y'y'')=S(\Delta(x),y'\otimes y'')$ for all $x,y',y''\in{\ '\ff}$;\\
(c) $S(xx',y'')=S(x\otimes x',\Delta(y''))$ for all $x,x',y''\in{\ '\ff}$.

(The bilinear form
$$
({\ '\ff}\otimes{\ '\ff})\otimes({\ '\ff}\otimes{\ '\ff})\lra \sk
$$
given by
$$
(x_1\otimes x_2)\otimes (y_1\otimes y_2)\mapsto S(x_1,y_1)S(x_2,y_2)
$$
is denoted again by $S(\ ,\ )$.)

The bilinear form $S(\ ,\ )$ on ${\ '\ff}$ is symmetric. The different
homogeneous components ${\ '\ff}_{\nu}$ are mutually orthogonal. }

{\bf Proof.} See ~\cite{l1}, 1.2.3. Cf. also ~\cite{sv2}, (1.8)-(1.11). $\Box$

\subsection{}
\label{dva epsilons} Following ~\cite{l1}, 1.2.13 and ~\cite{sv2}, (1.10)-(1.11),
let us introduce operators $\delta_i:{\ '\ff}\lra{\ '\ff},\ i\in I,$
as unique linear
mappings satisfying
\begin{equation}
\label{dva formeps}
\delta_i(1)=0;\ \delta_i(\theta_j)=\delta_{i,j},\ j\in I;\
\delta_i(xy)=\delta_i(x)y+\zeta^{|x|\cdot i}x\delta_i(y)
\end{equation}
for homogeneous $x$.

It follows from ~\ref{dva form} (c) that
\begin{equation}
\label{dva contrav}
S(\theta_ix,y)=S(x,\delta_i(y))
\end{equation}
for all $i\in I,\ x,y\in {\ '\ff}$, and obviously $S$ is determined uniquely
by this property, together with the requirement $S(1,1)=1$.

\subsection{}
\label{dva formulas}
{\bf Lemma.} {\em For any two sequences $\vec{K},\ \vec{K}'$ of $N$ elements
from $I$ we have
$$
S(\theta_{\vec{K}},\theta_{\vec{K}'})=\sum_{\tau\in\Sigma_N:\ \tau(\vec{K})=
\vec{K}'}\zeta(\vec{K};\tau).
$$}

{\bf Proof} follows from ~(\ref{dva contrav}) by induction on $N$, or else from
{}~\ref{dva deltamon}. $\Box$

\subsection{} Let us define elements $\theta^*_i\in{\ '\ff}_i^*$ by the rule
$<\theta^*_i,\theta_i>=1$. The form $S$ defines a homomomorphism of graded
algebras
\begin{equation}
\label{dva formap}
S:{\ '\ff}\lra{\ '\ff}^*
\end{equation}
carrying $\theta_i$ to $\theta^*_i$. $S$ is determined uniquely by
this property.

\subsection{}
\label{dva morcoalg} {\bf Lemma.} {\em The map $S$ is
a morphism of coalgebras.}

{\bf Proof.} This follows from the symmetry of $S$. $\Box$

\vspace{1cm}
{\em VERMA MODULES}
\vspace{1cm}

\subsection{}
\label{dva verma}
Let us pick a weight $\Lambda$. Our aim now will be to define
certain $X$-graded vector space $V(\Lambda)$ equipped with the following
structures.

(i) A structure of left ${\ '\ff}$-module
${\ '\ff}\otimes V(\Lambda)\lra V(\Lambda)$;

(ii) a structure of left ${\ '\ff}$-comodule $V(\Lambda)\lra{\ '\ff}\otimes
V(\Lambda)$;

(iii) a symmetric bilinear form $S_{\Lambda}$ on $V(\Lambda)$.

As a vector space, we set $V(\Lambda)={\ '\ff}$. We will define on $V(\Lambda)$
two gradings. The first one, $\BN[I]$-grading coincides with the grading
on ${\ '\ff}$. If $x\in V(\Lambda)$ is a homogeneous element, we will
denote by $\depth(x)$ its depth as an element of ${\ '\ff}$.

The second grading --- $X$-grading --- is defined
as follows. By definition, we set
$$
V(\Lambda)_{\lambda}=\oplus_{\nu\in\BN[I]|\Lambda-\lambda_{\nu}
=\lambda}{\ '\ff}_{\nu}
$$
for $\lambda\in X$. In particular, $V(\Lambda)_{\Lambda}={\ '\ff}_0=\sk\cdot 1$.
We will denote the element $1$ in $V(\Lambda)$ by $v_{\Lambda}$.

By definition, multiplication
\begin{equation}
\label{dva multv}
{\ '\ff}\otimes V(\Lambda)\lra V(\Lambda)
\end{equation}
coincides with the  multiplication in ${\ '\ff}$.

Let us define an $X$-grading in ${\ '\ff}$ by setting
$$
{\ '\ff}_{\lambda}=\oplus_{\nu\in\BN[I]|-\lambda_{\nu}=\lambda}{\ '\ff}_{\nu}
$$
for $\lambda\in X$. The map ~(\ref{dva multv})
is compatible with both $\BN[I]$ and $X$-gradings (we define gradings
on the tensor product as usually as a sum of gradings of factors).

\subsection{The form $S_{\Lambda}$}
\label{dva slambda} Let us define linear operators
$\epsilon_i:V(\Lambda)\lra V(\Lambda),\ i\in I,$ as unique
operators such that $\epsilon_i(v_{\Lambda})=0$
and
\begin{equation}
\label{dva epsilon}
\epsilon_i(\theta_jx)=[\langle\beta,i\rangle]_{\zeta}\delta_{i,j}x +
\zeta^{i\cdot j}\theta_j\epsilon_i(x)
\end{equation}
for $j\in I,\ x\in V(\Lambda)_{\beta}$.

We define $S_{\Lambda}: V(\Lambda)\otimes V(\Lambda)\lra \sk$ as a unique
linear map such that
$S_{\Lambda}(v_{\Lambda},v_{\Lambda})=1$, and
\begin{equation}
\label{dva vermacontra}
S_{\Lambda}(\theta_ix,y)=S_{\Lambda}(x,\epsilon_i(y))
\end{equation}
for all $x,y\in V(\Lambda),\ i\in I$.
Let us list elementary properties of $S_{\Lambda}$.

\subsubsection{}
\label{dva orthog} {\em Different graded components $V(\Lambda)_{\nu},\
\nu\in\BN[I]$, are orthogonal with respect to $S_{\Lambda}$}.

This follows directly from the definition.

\subsubsection{}
\label{dva ssymmetr} {\em The form $S_{\Lambda}$ is symmetric.}

This is an immediate corollary of the formula
\begin{equation}
\label{dva epsitheta}
S_{\Lambda}(\epsilon_i(y),x)=S_{\Lambda}(y,\theta_ix)
\end{equation}
which in turn is proved by an easy induction on $\depth(x)$.

\subsubsection{"Quasiclassical" limit}
\label{dva slclassic} Let us consider restriction
of our form to the homogeneous component $V(\Lambda)_{\lambda}$ of depth
$N$. If we divide our form by $(\zeta -1)^N$ and formally pass
to the limit $\zeta\lra 1$, we get the "Shapovalov" contravariant
form as defined in ~\cite{sv1}, 6.4.1.

The next lemma is similar to ~\ref{dva formulas}.

\subsection{}
\label{dva formulasv} {\bf Lemma.} {\em For any $\vec{K},\ \vec{K'}$ as in
{}~\ref{dva formulas} we have
$$
S_{\Lambda}(\theta_{\vec{K}}v_{\Lambda},\theta_{\vec{K'}}v_{\Lambda})=
\sum_{\tau\in\Sigma_N:\ \tau(\vec{K})=\vec{K'}}
\zeta(\vec{K};\tau)A(\vec{K},\Lambda;\tau)
$$
where
$$
A(\vec{K},\Lambda;\tau)=\prod_{a=1}^N[\langle\Lambda-\sum_{b:\ b<a,\tau(b)<
\tau(a)}\lambda_{i_b},i_a\rangle]_{\zeta}.
$$}

{\bf Proof.} Induction on $N$, using definition of $S_{\Lambda}$. $\Box$

\subsection{Coaction} Let us define a linear map
\begin{equation}
\label{dva coaction}
\Delta_{\Lambda}:V(\Lambda)\lra{\ '\ff}\otimes V(\Lambda)
\end{equation}
as follows. Let us introduce linear operators $t_i:{\ '\ff}^+\otimes V(\Lambda)
\lra{\ '\ff}^+\otimes V(\Lambda),\ i\in I,$ by the formula
\begin{equation}
\label{dva ti}
t_i(x\otimes y)=\theta_ix\otimes y-\zeta^{i\cdot\nu-2\langle\lambda,i\rangle}
\cdot x\theta_i\otimes y+\zeta^{i\cdot\nu}x\otimes\theta_iy
\end{equation}
for $x\in{\ '\ff}_{\nu}$ and $y\in V(\Lambda)_{\lambda}$.

By definition,
\begin{eqnarray}
\label{dva coactform}
\Delta_{\Lambda}(\theta_{i_N}\cdot\ldots\cdot\theta_{i_1}v_{\Lambda})=
1\otimes\theta_{i_N}\cdot\ldots\cdot\theta_{i_1}v_{\Lambda}\\ \nonumber
+[\langle\Lambda-\lambda_{i_{1}}-\ldots-\lambda_{i_{N-1}},i_N\rangle]_{\zeta}
\cdot\theta_{i_N}\otimes\theta_{i_{N-1}}\cdot\ldots\cdot\theta_{i_1}v_{\Lambda}
\\ \nonumber
+\sum_{j=1}^{N-1}
[\langle\Lambda-\lambda_{i_{1}}-\ldots-\lambda_{i_{j-1}},i_j\rangle]_{\zeta}
\cdot t_{i_N}\circ t_{i_{N-1}}\circ\ldots\circ t_{i_{j+1}}
(\theta_{i_j}\otimes\theta_{i_{j-1}}\cdot\ldots\cdot\theta_{i_1}v_{\Lambda})
\nonumber
\end{eqnarray}

\subsection{}
\label{dva quantcom} Let us define linear operators
\begin{equation}
\label{dva qadj}
\ad_{\theta_i,\lambda}:{\ '\ff}\lra{\ '\ff},\ i\in I,\ \lambda\in X
\end{equation}
by the formula
\begin{equation}
\label{dva quantadj}
\ad_{\theta_i,\lambda}(x)=\theta_ix-\zeta^{i\cdot\nu-2\langle\lambda,i\rangle}
\cdot x\theta_i
\end{equation}
for $x\in{\ '\ff}_{\nu}$.

Let us note the following relation
\begin{equation}
\label{dva adjdelta}
(\delta_i\circ\ad_{\theta_j,\lambda}-\zeta^{i\cdot
j}\cdot\ad_{\theta_j,\lambda}
\circ\delta_i)(x)=[\langle\lambda-\lambda_{\nu},i\rangle ]_{\zeta}\delta_{ij}x
\end{equation}
for $x\in{\ '\ff}_{\nu}$, where $\delta_i$ are operators defined in
{}~\ref{dva epsilons}, and $\delta_{ij}$ the Kronecker symbol.

\subsection{Formula for coaction}
\label{dva coactcom} Let us pick a sequence
$\vec{I}=(i_N,i_{N-1},\ldots,i_1)$. To shorten the notations, we set
\begin{equation}
\label{dva shortcom}
\ad_{j,\lambda}:=\ad_{\theta_{i_j},\lambda},\ j=1,\ldots, N
\end{equation}

\subsubsection{Quantum commutators}
\label{dva commut} For any non-empty subset $Q\subset [N]$, set
$\theta_{\vec{I},Q}:=\theta_{\vec{I}_Q}$ where $\vec{I}_Q$ denotes the
sequence obtained from $\vec{I}$ by omitting all entries $i_{j},\ j\in Q$.
We will denote ${\ '\ff}_Q={\ '\ff}_{\nu_Q}$ where $\nu_Q:=
\sum_{j\in Q}i_j$.

Let us define an element
$[\theta_{\vec{I},Q,\Lambda}]\in{\ '\ff}_{Q}$ as follows.
Set
\begin{equation}
\label{dva comone}
[\theta_{\vec{I},\{j\},\Lambda}]=\zeta^{i_j\cdot(\sum_{k>j}i_k)}
\theta_{i_j}
\end{equation}
for all $j\in [N]$.

Suppose now that $\card(Q)=l+1\geq 2$.
Let $Q=\{j_0,j_1,\ldots,j_l\},\ j_0<j_1<\ldots <j_l$.
Define the weights
$$
\lambda_a=\Lambda-\lambda_{\sum i_k},\ a=1,\ldots, l,
$$
where the summation is over $k$ from $1$ to $j_a-1$, $k\neq j_1,j_2,\ldots,
j_{a-1}$.

Let us define sequences $\vec{N}:=(N,N-1,\ldots,1),\
\vec{Q}=(j_l,j_{l-1},\ldots,j_0)$ and
$\vec{N}_Q$ obtained from $\vec{N}$ by omitting all entries $j\in Q$.
Define the permutation $\tau_Q\in \Sigma_N$ by the requirement
$$
\tau_Q(\vec{N})=\vec{Q}||\vec{N}_Q
$$
where $||$ denotes concatenation.

Set by definition
\begin{equation}
\label{dva comgen}
[\theta_{\vec{I},Q,\Lambda}]=\zeta(\vec{I},\tau_Q)\cdot
\ad_{j_l,\lambda_l}\circ\ad_{j_{l-1},\lambda_{l-1}}\circ\ldots\circ
\ad_{j_1,\lambda_1}(\theta_{i_{j_0}})
\end{equation}

\subsubsection{}
\label{dva formcoact} {\bf Lemma.} {\em We have
\begin{equation}
\label{dva sumcoact}
\Delta_{\Lambda}(\theta_{\vec{I}}v_{\Lambda})=
1\otimes \theta_{\vec{I}}v_{\Lambda}+
\sum_{Q}[\langle\Lambda-\lambda_{i_1}-\lambda_{i_2}-\ldots -
\lambda_{i_{j(Q)-1}},i_{j(Q)}\rangle]_{\zeta}\cdot
[\theta_{\vec{I},Q,\Lambda}]\otimes\theta_{\vec{I},Q}v_{\Lambda},
\end{equation}
the summation over all non-empty subsets $Q\subset [N]$,
$j(Q)$ denotes the minimal element of $Q$.}

{\bf Proof.} The statement of the lemma follows at once from the
inspection of definition ~(\ref{dva coactform}), after rearranging
the summands. $\Box$

Several remarks are in order.

\subsubsection{} Formula ~(\ref{dva sumcoact}) as similar to ~\cite{s}, 2.5.4.

\subsubsection{} If all elements $i_j$ are distinct then
the part of the sum  in the rhs of ~(\ref{dva sumcoact})
corresponding to one-element subsets $Q$ is
equal to
$\sum_{j=1}^N\theta_{i_j}\otimes\epsilon_{i_j}(\theta_{\vec{I}}v_{\Lambda})$.

\subsubsection{"Quasiclassical" limit}
It follows from the definition of quantum commutators that if
we divide the rhs of ~(\ref{dva sumcoact}) by $(\zeta -1)^N)$
and formally pass to the limit $\zeta\lra 1$, we get the expression for
the coaction
obtained in ~\cite{sv1}, 6.15.3.2.

\subsection{} Let us define the space $V(\Lambda)^*$ as the direct sum
$\oplus_{\nu}V(\Lambda)^*_{\nu}$. We define an $\BN[I]$-grading on it as
$(V(\Lambda)^*)_{\nu}=V(\Lambda)^*_{\nu}$, and an $X$-grading as
$V(\Lambda)^*_{\lambda}=\oplus_{\nu:\ \Lambda-\lambda_{\nu}=\lambda}
V(\Lambda)^*_{\nu}$.

The form $S_{\Lambda}$ induces the map
\begin{equation}
\label{dva maps}
S_{\Lambda}:V(\Lambda)\lra V(\Lambda)^*
\end{equation}
compatible with both gradings.

\subsection{Tensor products}
\label{dva tensprod} Suppose we are given $n$ weights
$\Lambda_0,\ldots,\Lambda_{n-1}$.

\subsubsection{} For every $m\in\BN$ we introduce a bilinear
form $S=S_{m;\Lambda_0,\ldots,\Lambda_{n-1}}$ on the tensor product
${\ '\ff}^{\otimes m}\otimes V(\Lambda_0)
\otimes\ldots\otimes V(\Lambda_{n-1})$ by the formula
$$
S(x_1\otimes\ldots\otimes x_m\otimes y_0\otimes\ldots
\otimes y_{n-1},
x'_1\otimes\ldots\otimes x'_m\otimes y'_0\otimes
\ldots\otimes y'_{n-1})=\prod_{i=1}^mS(x_i,x'_i)\prod_{j=0}^{n-1}
S_{\Lambda_j}(y_j,y'_j)
$$
(in the evident notations). This form defines mappings
\begin{equation}
\label{dva stens}
S:{\ '\ff}^{\otimes m}\otimes V(\Lambda_0)\otimes\ldots\otimes V(\Lambda_{n-1})
\lra
{\ '\ff}^{*\otimes m}\otimes V(\Lambda_0)^*\otimes\ldots\otimes V(\Lambda_{n-1})^*
\end{equation}

\subsubsection{} We will regard $V(\Lambda_0)
\otimes\ldots\otimes V(\Lambda_{n-1})$ as an ${\ '\ff}^{\otimes n}$-module with
an action
\begin{equation}
\label{dva multiter}
(u_0\otimes\ldots\otimes u_{n-1})\cdot (x_1\otimes\ldots\otimes x_{n-1})=
\zeta^{-\sum_{j<i}\langle\lambda_j,\nu_i\rangle}u_0x_0\otimes\ldots\otimes
u_{n-1}x_{n-1}
\end{equation}
for $u_i\in{\ '\ff}_{\nu_i},\ x_j\in V(\Lambda_j)_{\lambda_j}$, cf.
{}~(\ref{dva product}).
Here we regard ${\ '\ff}^{\otimes n}$ as an algebra according
to the rule of {\em loc. cit.}; one checks easily using ~(\ref{dva lamnu})
that we really get a module structure.

Using the iterated comultiplication $\Delta^{(n)}$, we get
a structure of an ${\ '\ff}$-module on
$V(\Lambda_0)\otimes\ldots\otimes V(\Lambda_{n-1})$.

\subsection{Theorem}
\label{dva coactshap} {\em We have an identity
\begin{equation}
\label{dva coactshapo}
S_{\Lambda}(xy,z)=S_{1;\Lambda}(x\otimes y,\Delta_{\Lambda}(z))
\end{equation}
for any $x\in{\ '\ff},\ y,z\in V(\Lambda)$ and any weight $\Lambda$.}

\subsection{Proof} We may suppose that $x,y$ and $z$ are monomials. Let
$z=\theta_{\vec{I}}v_{\Lambda}$ where $\vec{I}=(i_N,\ldots,i_1)$.

(a) Let us suppose first that all indices $i_j$ are distinct. We will
use the notations and computations from ~\ref{dva coactcom}. The sides of
{}~(\ref{dva coactshapo}) are non-zero only if $y$ is equal to
$z_Q:=\theta_{\vec{I},Q}v_{\Lambda}$ for some subset $Q\subset [N]$.

Therefore, it follows from Lemma ~\ref{dva formcoact} that it is enough to prove

\subsubsection{} {\bf Lemma.} {\em For every non-empty $Q\subset [N]$ and
$x\in{\ '\ff}_Q$ we have
\begin{equation}
\label{dva cominduct}
S_{\Lambda}(xz_Q,z)=
[\langle\Lambda,i_{j(Q)}\rangle-\mu_Q\cdot i_{j(Q)}]_{\zeta}
\cdot S(x,[\theta_{\vec{I},Q,\Lambda}])\cdot S_{\Lambda}(z_Q,z_Q)
\end{equation}
where $j(Q)$ denotes the minimal element of $Q$, and
$$
\mu_Q:=\sum_{a=1}^{j(Q)-1}i_a.
$$}

{\bf Proof.} If $\card(Q)=1$
the statement follows from the definiton ~(\ref{dva comone}). The proof will
proceed
by the simultaneous induction by $l$ and $N$.
Suppose that $x=\theta_{i_p}\cdot x'$, so $x'\in{\ '\ff}_{Q'}$ where
$Q'=Q-\{i_p\}$, $p=j_a$ for some $a\in [0,l]$. Let us set $\vec{I}'=
\vec{I}-\{i_p\},\ z'=z_{\{i_p\}}$, so that $z_Q=z'_{Q'}$.

We have
\begin{eqnarray}
S_{\Lambda}(\theta_{i_p}x'\cdot z_Q,z)=
S_{\Lambda}(x'\cdot z_Q,\epsilon_{i_p}(z))=\\ \nonumber
=[\langle\Lambda,i_p\rangle-(\sum_{k<p}i_k)\cdot i_p]_{\zeta}\cdot
\zeta^{(\sum_{k>p}i_k)\cdot i_p}\cdot
S(x'\cdot z_{Q},z')=\\ \nonumber
=[\langle\Lambda,i_p\rangle-(\sum_{k<p}i_k)\cdot i_p]_{\zeta}\cdot
[\langle\Lambda,i_{j(Q')}\rangle-\mu_{Q'}\cdot i_{j(Q')}]_{\zeta}\cdot
\zeta^{(\sum_{k>p}i_k)\cdot i_p}\cdot \\ \nonumber
\cdot S(x,[\theta_{\vec{I},Q,\Lambda}])\cdot S_{\Lambda}(z'_{Q'},z'_{Q'})
\nonumber
\end{eqnarray}
by induction hypothesis.
On the other hand,
$$
S(\theta_{i_p}\cdot x',[\theta_{\vec{I},Q,\Lambda}])=
S(x',\delta_{i_p}([\theta_{\vec{I},Q,\Lambda}])).
$$
Therefore, to complete the induction step it is enough to prove that
\begin{eqnarray}
\label{dva deltacommut}
[\langle\Lambda,i_{j(Q)}\rangle-\mu_Q\cdot i_{j(Q)}]_{\zeta}\cdot
\delta_{i_p}([\theta_{\vec{I},Q,\Lambda}])=\\ \nonumber
=[\langle\Lambda,i_p\rangle-(\sum_{k<p}i_k)\cdot i_p]_{\zeta}\cdot
[\langle\Lambda,i_{j(Q')}\rangle-\mu_{Q'}\cdot i_{j(Q')}]_{\zeta}\cdot
\zeta^{(\sum_{k>p}i_k)\cdot i_p}
[\theta_{\vec{I}',Q',\Lambda}] \nonumber
\end{eqnarray}
This formula follows directly from the definition of
quantum commutators ~(\ref{dva comgen}) and formula ~(\ref{dva adjdelta}).
One has to treat separately two cases: $a>0$, in which case $j(Q)=j(Q')=j_0$
and
$a=0$, in which case $j(Q)=j_0,\ j(Q')=j_1$. Lemma is proven.
$\Box$

This completes the proof of case (a).

(b) There are repeating indices in the sequence $\vec{I}$.
Suppose that $\theta_{\vec{I}}\in{\ '\ff}_{\nu}$. At this point we will use
symmetrization constructions (and simple facts) from Section ~\ref{dva symmetr}
below. The reader will readily see that there is no vicious circle.
So, this part of the proof must be read after {\em loc.cit.}

There exists a finite set $J$ and a map $\pi:J\lra I$ such that
$\nu=\nu_{\pi}$.
Using compatibility of the coaction and the forms $S$ with symmetrization ---
cf. Lemmata ~\ref{dva avercoact} and ~\ref{dva averterms} below ---
our claim is immediately
reduced to the analogous claim for the algebra $^{\pi}{\ '\ff}$,
the module $V(^{\pi}\Lambda)$ and homogeneous weight $\chi_J$
which does not contain multiple indices and therefore follows from (a) above.

This completes the proof of the theorem. $\Box$

\subsection{} Let us pick a weight $\Lambda$.
We can consider numbers $q_{ij}:=\zeta^{i\cdot j}$ and $r_i:=
\langle\Lambda,i\rangle$,  $i,j\in I$ as parameters of our bilinear
forms.

More precisely, for a given $\nu\in\BN[I]$ the matrix elements
of the form $S$ (resp., $S_{\Lambda}$) on ${\ '\ff}_{\nu}$ (resp.,
on $V(\Lambda)_{\Lambda-\lambda_{\nu}}$) in the standard bases of these spaces
are certain universal polynomials of $q_{ij}$ (resp., $q_{ij}$ and
$r_i$). Let us denote their determinants by $\det(S_{\nu})(\bq)$ and
$\det(S_{\Lambda,\nu})(\bq;\br)$ respectively. These determinants are
polynomials of corresponding variables with integer coefficients.

\subsubsection{}
\label{dva generic} {\bf Lemma.} {\em Polynomials $\det(S_{\nu})(\bq)$ and
$\det(S_{\Lambda,\nu})(\bq;\br)$ are not identically zero.

In other words, bilinear forms $S$ and $S_{\Lambda}$ are non-degenerate
for generic values of parameters  --- "Cartan matrix" $(q_{ij})$ and
"weight" $(r_i)$.}

{\bf Proof.} Let us consider the form $S_{\Lambda}$ first.
The specialization of the matrix of $S_{\Lambda,\nu}$ at $\zeta=1$ is the
identity matrix. It follows easily that $\det(S_{V,\nu})(\bq;\br)\neq 0$.

Similarly, the matrix of $S_{\nu}$ becomes identity at $\zeta=0$,
which implies the generic non-degeneracy.  $\Box$

\subsection{Theorem}
\label{dva coactas} {\em Coaction $\Delta_{\Lambda}$ is coassociative,
i.e.
\begin{equation}
\label{dva coasform}
(1_{{\ '\ff}}\otimes\Delta_{\Lambda})\circ\Delta_{\Lambda}=
(\Delta\otimes 1_{V(\Lambda)})\circ\Delta_{\Lambda}.
\end{equation}}

{\bf Proof.} The equality ~(\ref{dva coasform}) is a polynomial identity
depending on parameters $q_{ij}$ and $r_i$ of the preceding subsection.
For generic values of these parameters it is true due to associativity
of the action of ${\ '\ff}$ an $V(\Lambda)$, Theorem ~\ref{dva coactshap} and
Lemma ~\ref{dva generic}. Therefore it is true for all values
of parameters. $\Box$

\subsection{} The results Chapter 2 below provide a different,
geometric proof of Theorems ~\ref{dva coactshap} and ~\ref{dva coactas}.
Namely, the results of Section ~\ref{dva standsheaves} summarized in Theorem
{}~\ref{dva phi*} provide an isomorphism of our algebraic picture
with a geometric one, and in the geometrical language the above theorems
are obvious: they
are nothing but the naturality of the canonical morphism
between the extension by zero and the extension by star,
and the claim that a Cousin complex is a complex. Lemma ~\ref{dva generic}
also follows from geometric considerations: the extensions by zero and
by star coincide for generic values of monodromy.

\subsection{}
\label{dva comod} By Theorem ~\ref{dva coactas} the dual maps
\begin{equation}
\label{dva deltal*}
\Delta_{\Lambda}^*:{\ '\ff}\otimes V(\Lambda)^*\lra V(\Lambda)^*
\end{equation}
give rise to a structure
of a ${\ '\ff}^*$-module on $V(\Lambda)^*$.

More generally, suppose we are given $n$ modules $V(\Lambda_0),\ldots,
V(\Lambda_{n-1})$. We regard the tensor product
$V(\Lambda_0)^*\otimes\ldots\otimes V(\Lambda_{n-1})^*$ as a
${\ '\ff}^{*\otimes n}$-module according to the "sign" rule ~(\ref{dva multiter}).
Using iterated comultiplication ~(\ref{dva deltait}) we get a structure
of a ${\ '\ff}^*$-module on $V(\Lambda_0)^*\otimes\ldots\otimes
V(\Lambda_{n-1})^*$.

\subsubsection{}
\label{dva scomod} The square
$$\begin{array}{ccccc}
\;&{\ '\ff}\otimes V(\Lambda_0)\otimes\ldots\otimes
V(\Lambda_{n-1})&\lra&V(\Lambda_0)\otimes\ldots\otimes
V(\Lambda_{n-1})&\;\\ \nonumber
\;&S\downarrow&\;&\downarrow S&\; \\ \nonumber
\;&{\ '\ff}^*\otimes V(\Lambda_0)^*\otimes\ldots\otimes
V(\Lambda_{n-1})^*&\lra&V(\Lambda_0)^*\otimes\ldots\otimes
V(\Lambda_{n-1})^*&\;
\end{array}$$
commutes.

This follows from ~\ref{dva coactshap} and ~\ref{dva morcoalg}.

\section{Hochschild complexes}

\subsection{} If $A$ is an augmented $\sk$-algebra, $A^+$ --- the kernel
of the augmentation, $M$ an $A$-module, let $C_A^{\bullet}(M)$ denote
the following complex. By definition, $C_A^{\bullet}(M)$ is concentrated
in non-positive degrees. For $r\geq 0$
$$
C_A^{-r}(M)=A^{+\otimes r}\otimes M.
$$
We will use a notation $a_r|\ldots|a_1|m$ for $a_r\otimes\ldots a_1\otimes m$.

The differential $d:C_A^{-r}(M)\lra C_A^{-r+1}(M)$ acts as
$$
d(a_r|\ldots|a_1|m)=
\sum_{p=1}^{r-1}(-1)^{p}a_r|\ldots|a_{p+1}a_p|\ldots a_1|m+
a_r|\ldots a_{2}|a_1m.
$$
We have canonically $H^{-r}(C_A^{\bullet}(M))\cong \Tor^A_r(\sk,M)$ where
$\sk$ is considered as an $A$-module by means of the augmentation,
cf. ~\cite{mac}, Ch. X, \S 2.

We will be interested in the algebras ${\ '\ff}$ and ${\ '\ff}^*$. We define
the augmentation
${\ '\ff}\lra \sk$ as being zero on all
${\ '\ff}_{\nu},\ \nu\in\BN[I],\ \nu\neq 0$,
and identity on ${\ '\ff}_0$; in the same way it is defined on ${\ '\ff}^*$.

\subsection{}
\label{dva gradings} Let $M$ be a $\BN[I]$-graded ${\ '\ff}$-module.
Each term $C_{{\ '\ff}}^{-r}(M)$ is $\BN[I]$-graded
by the sum of gradings of tensor factors.
We will denote $_{\nu}C_{{\ '\ff}}^{-r}(M)$ the weight $\nu$ component.

For $\bnu=(\nu_0,\ldots,\nu_r)\in \BN[I]^{r+1}$ we set
$$
_{\bnu}C_{{\ '\ff}}^{-r}(M)=\ _{\nu_r,\ldots,\nu_0}C_{{\ '\ff}}^{-r}(M)=
{\ '\ff}_{\nu_r}\otimes\ldots\otimes{\ '\ff}_{\nu_1}\otimes M_{\nu_0}.
$$
Thus,
$$
_{\nu}C_{{\ '\ff}}^{-r}(M)=\oplus_{\nu_0+\ldots\nu_r=\nu}\
_{\nu_r,\ldots,\nu_0}C_{{\ '\ff}}^{-r}(M).
$$
Note that all $\nu_p$ must be $>0$ for $p>0$ since tensor factors lie in
${\ '\ff}^+$.

The differential $d$ clearly respects the  $\BN[I]$-grading; thus the whole
complex is $\BN[I]$-graded:
$$
C_{{\ '\ff}}^{\bullet}(M)=\oplus_{\nu\in\BN[I]}\
_{\nu}C_{{\ '\ff}}^{\bullet}(M).
$$
The same discussion applies to $\BN[I]$-graded ${\ '\ff}^*$-modules.

\subsection{}
\label{dva gradehoch} Let us fix weights
$\Lambda_0,\ldots,\Lambda_{n-1}$, $n\geq 1$. We will consider the Hochschild
complex
$C^{\bullet}_{{\ '\ff}}(V(\Lambda_0)\otimes\ldots\otimes V(\Lambda_{n-1}))$
where the structure of an ${\ '\ff}$-module on
$V(\Lambda_0)\otimes\ldots\otimes V(\Lambda_{n-1})$ has been introduced
in ~\ref{dva tensprod}.

\subsubsection{} In the sequel we will use the following notation.
If $K\subset I$ is a subset, we will denote by $\chi_K:=\sum_{i\in K}i
\in \BN[I]$.

\subsubsection{}
\label{dva gh} Suppose we have a map
\begin{equation}
\label{dva maprho}
\varrho: I\lra [-n+1, r]
\end{equation}
where $r$ is some non-negative integer. Let us introduce the elements
\begin{equation}
\label{dva nurho}
\nu_a(\varrho)=\chi_{\varrho^{-1}(a)},
\end{equation}
$a\in [-n+1,r]$.
Let us denote by $\CP_r(I;n)$ the set of all maps ~(\ref{dva maprho}) such that
$\varrho^{-1}(a)\neq\emp$ for all $a\in[r]$. It is easy to see that this set
is not empty iff $0\leq r\leq N$.

Let us assign to such a $\varrho$ the space
\begin{equation}
\label{dva homog}
_{\varrho}C^{-r}_{{\ '\ff}}(V(\Lambda_{0})\otimes\ldots\otimes V(\Lambda_{n-1})):=
{\ '\ff}_{\nu_r(\varrho)}\otimes\ldots\otimes{\ '\ff}_{\nu_1(\varrho)}\otimes
V(\Lambda_{0})_{\nu_0(\varrho)}\otimes\ldots\otimes
V(\Lambda_{n-1})_{\nu_{-n+1}(\varrho)}
\end{equation}
For each $\varrho\in\CP_r(I;n)$ this space is non-zero, and we have
\begin{equation}
\label{dva decompos}
_{\chi_I}C^{-r}_{{\ '\ff}}(V(\Lambda_{0})\otimes\ldots\otimes V(\Lambda_{n-1}))=
\oplus_{\varrho\in\CP_r(I;n)}\
_{\varrho}C^{-r}_{{\ '\ff}}(V(\Lambda_{0})\otimes\ldots\otimes V(\Lambda_{n-1}))
\end{equation}

\subsection{Bases}
\label{dva refine} Let us consider the set $\CP_N(I;n)$. Obviously, if
$\varrho\in\CP_N(I;n)$ then $\varrho(I)=[N]$, and the induced
map $I\lra [N]$ is a bijection; this way we get an isomorphism
between $\CP_N(I;n)$ and the set of all bijections $I\iso [N]$ or,
to put it differently, with the set of all total orders on $I$.

For an arbitrary $r$, let $\varrho\in\CP_r(I;n)$ and $\tau\in\CP_N(I;n)$.
Let us say that $\tau$ is {\em a refinement of $\varrho$}, and
write $\varrho\leq\tau$, if
$\varrho(i)<\varrho(j)$ implies $\tau(i)<\tau(j)$ for each $i,j\in I$.
The map $\tau$ induces total orders on all subsets $\varrho^{-1}(a)$.
We will denote by $\Ord(\varrho)$ the set of all refinements of a given
$\varrho$.

Given $\varrho\leq\tau$ as above, and $a\in [-n+1,r]$, suppose that
$\varrho^{-1}(a)=\{ i_1,\ldots,i_p\}$ and
$\tau(i_1)<\tau(i_2)<\ldots<\tau(i_p)$. Let us define a monomial
$$
\theta_{\varrho\leq\tau;a}=\theta_{i_p}\theta_{i_{p-1}}\cdot\ldots\cdot
\theta_{i_1}\in{\ '\ff}_{\nu_a(\varrho)}
$$
If $\varrho^{-1}(a)=\emp$, we set $\theta_{\varrho\leq\tau;a}=1$.
This defines a monomial
\begin{equation}
\label{dva monom}
\theta_{\varrho\leq\tau}=\theta_{\varrho\leq\tau;r}\otimes\ldots\otimes
\theta_{\varrho\leq\tau;1}\otimes\theta_{\varrho\leq\tau;0}v_{\Lambda_0}\otimes
\ldots\otimes\theta_{\varrho\leq\tau;-n+1}v_{\Lambda_{n-1}}\in\
_{\varrho}C^{-r}_{{\ '\ff}}(V(\Lambda_{0})\otimes\ldots\otimes V(\Lambda_{n-1}))
\end{equation}

\subsubsection{}
\label{dva basis} {\bf Lemma.} {\em The set  $\{\theta_{\varrho\leq\tau}|
\tau\in\Ord(\varrho)\}$ forms a basis of the space
{}~{$_{\varrho}C^{-r}_{{\ '\ff}}(V(\Lambda_{0})\otimes\ldots\otimes
V(\Lambda_{n-1}))$.}}

{\bf Proof} is obvious. $\Box$

\subsubsection{}
\label{dva basisj} {\bf Corollary.} {\em The set  $\{\theta_{\varrho\leq\tau}|
\varrho\in\CP_r(I;n),\ \tau\in\Ord(\varrho)\}$ forms a basis of the space
$_{\chi_I}C^{-r}_{{\ '\ff}}(V(\Lambda_{0})\otimes\ldots\otimes
V(\Lambda_{n-1}))$.} $\Box$

\subsection{}
\label{dva dualbases} We will also consider dual Hochschild complexes
$C^{\bullet}_{{\ '\ff}^*}(V(\Lambda_{0})^*\otimes\ldots\otimes
V(\Lambda_{n-1})^*)$ where $V(\Lambda_{0})^*\otimes\ldots\otimes
V(\Lambda_{n-1})^*$ is regarded as an ${\ '\ff}^*$-module as in ~\ref{dva scomod}.

We have obvious isomorphisms
$$
C^{-r}_{{\ '\ff}^*}(V(\Lambda_{0})^*\otimes\ldots\otimes
V(\Lambda_{n-1})^*)\cong
C^{-r}_{{\ '\ff}}(V(\Lambda_{0})\otimes\ldots\otimes
V(\Lambda_{n-1}))^*
$$
We define graded components
$$
_{\varrho}C^{-r}_{{\ '\ff}^*}(V(\Lambda_{0})^*\otimes\ldots\otimes
V(\Lambda_{n-1})^*),\
\varrho\in\CP_r(I;n),
$$
as duals to $_{\varrho}C^{-r}_{{\ '\ff}}(V(\Lambda_{0})\otimes\ldots\otimes
V(\Lambda_{n-1}))$.

We will denote by $\{\theta_{\varrho\leq\tau}^*|
\varrho\in\CP_r(I;n),\ \tau\in\Ord(\varrho)\}$ the basis of
$_{\chi_I}C^{-r}_{{\ '\ff}^*}(V(\Lambda_{0})^*\otimes\ldots\otimes
V(\Lambda_{n-1})^*)$ dual to the basis
$\{\theta_{\varrho\leq\tau}|
\varrho\in\CP_r(I;n),\ \tau\in\Ord(\varrho)\}$, ~\ref{dva basisj}.

\subsection{} The maps $S_{r;\Lambda_0,\ldots,\Lambda_{n-1}}$, cf.
{}~(\ref{dva stens}), for different $r$ are compatible with differentials
in Hochschild complexes, and therefore induce morphism
of complexes
\begin{equation}
\label{dva shapohoch}
S:C^{\bullet}_{{\ '\ff}}(V(\Lambda_{0})\otimes\ldots\otimes
V(\Lambda_{n-1})) \lra
C^{\bullet}_{{\ '\ff}^*}(V(\Lambda_{0})^*\otimes\ldots\otimes
V(\Lambda_{n-1})^*)
\end{equation}
This follows from ~\ref{dva scomod} and ~\ref{dva form} (b).

\section{Symmetrization}
\label{dva symmetr}

\subsection{}
\label{dva sym} Let us fix a finite set $J$ and a map $\pi:J\lra I$.
We set $\nu_{\pi}:=\sum_iN_ii\in\BZ[I]$ where
$N_i:=\card(\pi^{-1}(i))$.
The map $\pi$ induces a map $\BZ[J]\lra \BZ[I]$ also
to be denoted by $\pi$.
We will use the notation $\chi_K:=\sum_{j\in K}j\in\BN[J]$ for $K\subset J$.
Thus, $\pi(\chi_J)=\nu_{\pi}$.

We will denote also by $\mu,\mu'\mapsto
\mu\cdot\mu':=\pi(\mu)\cdot\pi(\mu')$ the bilinear form on $\BZ[J]$ induced
by the form on $\BZ[I]$.

We will denote by $\Sigma_{\pi}$ the group of all
bijections $\sigma:J\lra J$ preserving fibers of $\pi$.

Let $^{\pi}{\ '\ff}$ be a free associative $\sk$-algebra with $1$ with generators
$\ttheta_j,\ j\in J$. It is evidently $\BN[J]$-graded.
For $\nu\in\BN[J]$ the corresponding homogeneous component will be denoted
$^{\pi}{\ '\ff}_{\nu}$. The degree of a
homogeneous element $x\in\ ^{\pi}{\ '\ff}$ will be denoted by $|x|\in\BN[J]$.
The group $\Sigma_{\pi}$ acts on algebras
$^{\pi}{\ '\ff},\ ^{\pi}{\ '\ff}^*$ by
permutation of generators.

\subsection{}
\label{dva average} In the sequel, if $G$ is a group and $M$ is a $G$-module,
$M^G$ will denote the subset of $G$-invariants in $M$.

Let us define a $\sk$-linear "averaging" mapping
\begin{equation}
\label{dva avermap}
^{\pi}a:{\ '\ff}_{\nu_{\pi}}\lra (^{\pi}{\ '\ff}_{\chi_J})^{\Sigma_{\pi}}
\end{equation}
by the rule
\begin{equation}
\label{dva avermapform}
^{\pi}a(\theta_{i_1}\cdot\ldots\cdot\theta_{i_N})=
\sum \ttheta_{j_1}\cdot\ldots\cdot \ttheta_{j_N},
\end{equation}
the sum being taken over the set of all sequences $(j_1,\ldots,j_N)$ such
that $\pi(j_p)=i_p$ for any $p$. Note that this set is naturally a
$\Sigma_{\pi}$-torsor.
Alternatively, $^{\pi}a$ may be defined as follows. Pick some sequence
$(j_1,\ldots,j_N)$ as above, and consider an element
$$
\sum_{\sigma\in\Sigma_{\pi}}\sigma(\ttheta_{j_1}\cdot\ldots\cdot
\ttheta_{j_N});
$$
this element obviously lies in $(^{\pi}{\ '\ff}_{\chi_J})^{\Sigma_{\pi}}$ and
is equal to $^{\pi}a(\theta_{i_1}\cdot\ldots\cdot\theta_{i_N})$.

The map $\pi$ induces the map between homogeneous components
\begin{equation}
\label{dva ainv}
\pi:\ ^{\pi}{\ '\ff}_{\chi_J}\lra {\ '\ff}_{\nu_{\pi}}.
\end{equation}
It is clear that the composition $\pi\circ\ ^{\pi}a$ is equal to the
multiplication by $\card (\Sigma_{\pi})$, and $^{\pi}a\circ\pi$ ---
to the action of operator $\sum_{\sigma\in\Sigma_{\pi}}\sigma$.
As a consequence, we get

\subsubsection{}
{\bf Lemma,} ~\cite{sv1}, 5.11. {\em The map $^{\pi}a$ is an isomorphism.}
$\Box$

\subsection{}
\label{dva averdual} Let us consider the dual to the map ~(\ref{dva ainv}):
${\ '\ff}_{\nu_{\pi}}^*\lra\ ^{\pi}{\ '\ff}_{\chi_J}^*$; it is obvious that it
lands in the subspace of $\Sigma_{\pi}$-invariant functionals. Let us
consider the induced map
\begin{equation}
\label{dva avermap*}
^{\pi}a^*:{\ '\ff}_{\nu_{\pi}}^*\iso (^{\pi}{\ '\ff}_{\chi_J}^*)^{\Sigma_{\pi}}
\end{equation}
It follows from the above discussion that $^{\pi}a^*$ is an isomorphism.

\subsection{}
\label{dva averhoch} Given a weight $\Lambda\in X=\Hom(\BZ[I],\BZ)$, we will denote
by $^{\pi}\Lambda$ the composition $\BZ[J]\overset{\pi}{\lra}\BZ[I]
\overset{\Lambda}{\lra}\BZ$, and by $V(^{\pi}\Lambda)$ the corresponding
Verma module over $^{\pi}{\ '\ff}$.

Suppose we are given $n$ weights $\Lambda_0,\ldots,\Lambda_{n-1}$.
Let us consider the Hochschild complex
{}~{$C^{\bullet}_{^{\pi}{\ '\ff}}(V(^{\pi}\Lambda_0)\otimes\ldots\otimes
V(^{\pi}\Lambda_{n-1}))$}. By definition, its $(-r)$-th term coincides with
the tensor power $^{\pi}{\ '\ff}^{\otimes n+r}$. Therefore we can identify
the homogeneous component
$_{\chi_J}C^{-r}_{^{\pi}{\ '\ff}}(V(^{\pi}\Lambda_0)\otimes\ldots\otimes
V(^{\pi}\Lambda_{n-1}))$ with
 $(^{\pi}{\ '\ff}^{\otimes n+r})_{\chi_J}$ which in
turn is isomorphic to $^{\pi}{\ '\ff}_{\chi_J}$, by means of the multiplication
map $^{\pi}{\ '\ff}^{\otimes n+r}\lra\ ^{\pi}{\ '\ff}$. This defines a map
\begin{equation}
\label{dva hochcompj}
_{\chi_J}C^{-r}_{^{\pi}{\ '\ff}}(V(^{\pi}\Lambda_0)\otimes\ldots\otimes
V(^{\pi}\Lambda_{n-1}))\lra\ ^{\pi}{\ '\ff}_{\chi_J}
\end{equation}
which is an embedding when restricted to polygraded components.
The $\Sigma_{\pi}$-action on ${\ '\ff}$ induces the $\Sigma_{\pi}$-action on
$_{\chi_J}C^{-r}_{^{\pi}{\ '\ff}}(V(^{\pi}\Lambda_0)\otimes\ldots\otimes
V(^{\pi}\Lambda_{n-1}))$.

In the same manner we define a map
\begin{equation}
\label{dva hochcomp}
_{\nu_{\pi}}C^{-r}_{{\ '\ff}}(V(\Lambda_0)\otimes\ldots\otimes
V(\Lambda_{n-1}))\lra{\ '\ff}_{\nu_{\pi}}
\end{equation}

Let us define an averaging map
\begin{equation}
\label{dva avermaphoch}
^{\pi}a:\ _{\nu_{\pi}}C^{-r}_{{\ '\ff}}(V(\Lambda_0)\otimes\ldots\otimes
V(\Lambda_{n-1}))\lra\
_{\chi_J}C^{-r}_{^{\pi}{\ '\ff}}(V(^{\pi}\Lambda_0)\otimes\ldots\otimes
V(^{\pi}\Lambda_{n-1}))^{\Sigma_{\pi}}
\end{equation}
as the map induced by ~(\ref{dva avermap}). It follows at once that this map is
is an isomorphism.

These maps for different $r$ are by definition compatible with
differentials in Hochschild complexes. Therefore we get

\subsubsection{}
\label{dva avehoch} {\bf Lemma.} {\em The maps ~(\ref{dva avermaphoch})
induce isomorphism of complexes
\begin{equation}
\label{dva averhochiso}
^{\pi}a:\ _{\nu_{\pi}}C^{\bullet}_{{\ '\ff}}(V(\Lambda_0)\otimes\ldots\otimes
V(\Lambda_{n-1}))\iso\
_{\chi_J}C^{\bullet}_{^{\pi}{\ '\ff}}(V(^{\pi}\Lambda_0)\otimes\ldots\otimes
V(^{\pi}\Lambda_{n-1}))^{\Sigma_{\pi}}.\ \Box
\end{equation}}

\subsection{}
\label{dva avercoact} {\bf Lemma.} {\em The averaging
is compatible with coaction.
In other words, for any $\Lambda\in X$ the square
$$\begin{array}{ccccc}
\;&V(\Lambda)&\overset{\Delta_{\Lambda}}{\lra}&{\ '\ff}\otimes V(\Lambda)&\;\\
\;&^{\pi}a\downarrow&\;&\downarrow\ ^{\pi}a&\;\\
\;&V(^{\pi}\Lambda)&\overset{\Delta_{^{\pi}\Lambda}}{\lra}&^{\pi}{\ '\ff}
                                              \otimes V(^{\pi}\Lambda)&\;
\end{array}$$
commutes.}

{\bf Proof} follows at once by inspection of the definition ~(\ref{dva coactform}).
$\Box$

\subsection{}
\label{dva averhoch*} Consider the dual Hochschild complexes.
We have an obvious isomorphism
$$
_{\chi_J}C^{-r}_{^{\pi}{\ '\ff}^*}(V(^{\pi}\Lambda_0)^*\otimes\ldots\otimes
V(^{\pi}\Lambda_{n-1})^*)\cong\
_{\chi_J}C^{-r}_{^{\pi}{\ '\ff}}(V(^{\pi}\Lambda_0)\otimes\ldots\otimes
V(^{\pi}\Lambda_{n-1}))^*;
$$
using it, we define the isomorphism
$$
_{\chi_J}C^{-r}_{^{\pi}{\ '\ff}^*}(V(^{\pi}\Lambda_0)^*\otimes\ldots\otimes
V(^{\pi}\Lambda_{n-1})^*)\iso\ ^{\pi}{\ '\ff}^*_{\chi_J}
$$
as the dual to ~(\ref{dva hochcompj}). The $\Sigma_{\pi}$-action on the target
induces the action on
$_{\chi_J}C^{-r}_{^{\pi}{\ '\ff}^*}(V(^{\pi}\Lambda_0)^*\otimes\ldots\otimes
V(^{\pi}\Lambda_{n-1})^*)$.
Similarly, the isomorphism
$$
_{\nu_{\pi}}C^{-r}_{{\ '\ff}^*}(V(\Lambda_0)^*\otimes\ldots\otimes
V(\Lambda_{n-1})^*)\iso {\ '\ff}^*_{\nu_{\pi}}
$$
is defined. We define the averaging map
\begin{equation}
\label{dva avermaphoch*}
^{\pi}a^*: _{\nu_{\pi}}C^{-r}_{{\ '\ff}^*}(V(\Lambda_0)^*\otimes\ldots\otimes
V(\Lambda_{n-1})^*)\lra
_{\chi_J}C^{-r}_{^{\pi}{\ '\ff}^*}(V(^{\pi}\Lambda_0)^*\otimes\ldots\otimes
V(^{\pi}\Lambda_{n-1})^*)^{\Sigma_{\pi}}
\end{equation}
as the map which coincides with ~(\ref{dva avermap*}) modulo the above
identifications. Again, this map is an isomorphism.

Due to Lemma ~\ref{dva avercoact} these maps for different $r$ are compatible
with the differentials in Hochschild complexes. Therefore we get

\subsubsection{}
\label{dva avehoch*} {\bf Lemma.} {\em The maps ~(\ref{dva avermaphoch*})
induce isomorphism of complexes
\begin{equation}
\label{dva averhochiso*}
^{\pi}a^*: _{\nu_{\pi}}C^{\bullet}_{{\ '\ff}^*}(V(\Lambda_0)^*
\otimes\ldots\otimes
V(\Lambda_{n-1})^*)\iso\
_{\chi_J}C^{\bullet}_{^{\pi}{\ '\ff}^*}(V(^{\pi}\Lambda_0)^*
\otimes\ldots\otimes
V(^{\pi}\Lambda_{n-1})^*)^{\Sigma_{\pi}}.\ \Box
\end{equation}}

\vspace{1cm}
{\em BILINEAR FORMS}
\vspace{1cm}

\subsection{} Using the bilinear form on $\BZ[J]$ introduced above,
we define the symmetric bilinear form
$S(\ ,\ )$ on $^{\pi}{\ '\ff}$ exactly in the same way as the form $S$ on
${\ '\ff}$. Similarly, given $\Lambda\in X$, we define the bilinear
form $S_{^{\pi}\Lambda}$ on $V(^{\pi}\Lambda)$ as in ~\ref{dva slambda},
with $I$ replaced by $J$.

\subsubsection{}
\label{dva S-sym}
{\bf Lemma.} (i) {\em The square
$$\begin{array}{ccccc}
\;&{\ '\ff}_{\nu_{\pi}}&\overset{S}{\lra}&{\ '\ff}_{\nu_{\pi}}^*&\;\\
\;&^{\pi}a\downarrow&\;&\downarrow\ ^{\pi}a^*&\;\\
\;&^{\pi}{\ '\ff}_{\chi_J}&\overset{S}{\lra}&^{\pi}{\ '\ff}_{\chi_J}^*&\;
\end{array}$$
commutes.}

(ii) {\em For any $\Lambda\in X$ the square
$$\begin{array}{ccccc}
\;&V(\Lambda)_{\nu_{\pi}}&\overset{S_{\Lambda}}{\lra}&V(\Lambda)_{\nu_{\pi}}^*
&\;\\
\;&^{\pi}a\downarrow&\;&\downarrow\ ^{\pi}a^*&\;\\
\;&V(^{\pi}\Lambda)_{\chi_J}&\overset{S_{^{\pi}\Lambda}}{\lra}&
V(^{\pi}\Lambda)_{\chi_J}^*&\;
\end{array}$$
commutes.}

{\bf Proof.} (i) Let us consider an element
$\theta_{\vec{I}}=\theta_{i_1}\cdot\ldots\cdot\theta_{i_N}\in{\ '\ff}_{\nu_{\pi}}$
(we assume that $N=\card (J)$).
The functional $^{\pi}a^*\circ S(\theta_{\vec{I}})$ carries a monomial
$\ttheta_{j_1}\cdot\ldots\cdot\ttheta_{j_N}$ to
$$
S(\theta_{i_1}\cdot\ldots\cdot\theta_{i_N},
\theta_{\pi(j_1)}\cdot\ldots\cdot\theta_{\pi(j_N)}).
$$
On the other hand,
$$
S\circ\ ^{\pi}a(\theta_{\vec{I}})(\ttheta_{j_1}\cdot\ldots\cdot\ttheta_{j_N})=
\sum S(\ttheta_{k_1}\cdot\ldots\cdot\ttheta_{k_N},
\ttheta_{j_1}\cdot\ldots\cdot\ttheta_{j_N}),
$$
the summation ranging over all sequences $\vec{K}=(k_1,\ldots,k_N)$ such that
$\pi(\vec{K})=\vec{I}$. It follows from Lemma ~\ref{dva formulas} that
both expressions are equal.

(ii) The same argument as in (i), using Lemma ~\ref{dva formulasv} instead
of ~\ref{dva formulas}. $\Box$

More generally, we have

\subsection{Lemma}
\label{dva averterms} {\em For every $m\geq 0$ and weights $\Lambda_0,\ldots,
\Lambda_{n-1}\in X$ the square
$$\begin{array}{ccccc}
\;&({\ '\ff}^{\otimes m}\otimes V(\Lambda_0)\otimes\ldots\otimes
V(\Lambda_{n-1}))_{\nu_{\pi}}&
\overset{S_{m;\Lambda_0,\ldots,\Lambda_{n-1}}}{\lra}&
({\ '\ff}^{*\otimes m}\otimes V(\Lambda_0)^*\otimes\ldots\otimes
V(\Lambda_{n-1})^*)_{\nu_{\pi}}
&\;\\
\;&^{\pi}a\downarrow&\;&\downarrow\ ^{\pi}a^*&\;\\
\;&(^{\pi}{\ '\ff}^{\otimes m}\otimes V(^{\pi}\Lambda_0)\otimes\ldots\otimes
V(^{\pi}\Lambda_{n-1}))_{\chi_J}
&\overset{S_{m;^{\pi}\Lambda_0,\ldots,^{\pi}\Lambda_{n-1}}}{\lra}&
(^{\pi}{\ '\ff}^{*\otimes m}\otimes V(^{\pi}\Lambda_0)^*\otimes\ldots\otimes
V(^{\pi}\Lambda_{n-1})^*)_{\chi_J}&\;
\end{array}$$
commutes.}

{\bf Proof} is quite similar to the proof of the previous lemma.
We leave it to the reader. $\Box$

\section{Quotient algebras}

\subsection{} Let us consider the map ~(\ref{dva formap})
$S:{\ '\ff}\lra{\ '\ff}^*$. Let us consider its kernel $\Ker (S)$.
It follows at once from ~(\ref{dva contrav}) that $\Ker (S)$ is a left
ideal in ${\ '\ff}$. In the same manner, it is easy to see that
it is also a right ideal, cf. ~\cite{l1}, 1.2.4.

We will denote by $\ff$ the quotient algebra ${\ '\ff}/\Ker(S)$.
It inherits the $\BN[I]$-grading and the coalgebra structure from
${\ '\ff}$, cf. {\em loc.cit.} 1.2.5, 1.2.6.

\subsection{}
\label{dva llambda} In the same manner, given a weight $\Lambda$, consider
the kernel
of $S_{\Lambda}:V(\Lambda)\lra V(\Lambda)^*$. Let us denote
by $L(\Lambda)$ the quotient space $V(\Lambda)/\Ker(S_{\Lambda})$.
It inherits $\BN[I]$- and $X$-gradings from $V(\Lambda)$.
Due to Theorem ~\ref{dva coactshap} the structure of ${\ '\ff}$-module
on $V(\Lambda)$ induces the structure of $\ff$-module
on $L(\Lambda)$.

More generally, due to the structure of a coalgebra on $\ff$, all
tensor products $L(\Lambda_0)\otimes\ldots\otimes L(\Lambda_{n-1})$
become $\ff$-modules (one should take into account the "sign rule"
{}~(\ref{dva multiter})).

\subsection{} We can consider Hochschild complexes
$C^{\bullet}_{\ff}(L(\Lambda_0)\otimes\ldots\otimes L(\Lambda_{n-1}))$.

\subsubsection{} {\bf Lemma.} {\em The map
$$
S:C^{\bullet}_{{\ '\ff}}(V(\Lambda_0)\otimes\ldots\otimes V(\Lambda_{n-1}))\lra
C^{\bullet}_{{\ '\ff}^*}(V(\Lambda_0)^*\otimes\ldots\otimes V(\Lambda_{n-1})^*)
$$
factors through the isomorphism
\begin{equation}
\label{dva ims}
\Ima(S)\iso
C^{\bullet}_{\ff}(L(\Lambda_0)\otimes\ldots\otimes L(\Lambda_{n-1}))
\end{equation}}

{\bf Proof.} This follows at once from the definitions. $\Box$

\newpage
\begin{center}
{\bf Chapter 2. Geometric discussion}
\end{center}
\vspace{1cm}

\section{Diagonal stratification and related algebras}

\subsection{}
\label{dva diagsetup} Let us adopt notations of ~\ref{dva sym}. We set
$N:=\card(J)$.
Let $^{\pi}\BA_{\BR}$
denote a real affine space with coordinates $t_j,\ j\in J$, and
$^{\pi}\BA$ its complexification. Let us consider an arrangement $\CH_{\emp}$
consisting of all diagonals $\Delta_{ij},\ i,j\in J$. Let us denote
by $\CS_{\emp}$ the corresponding stratification;
$\CS_{\emp,\BR}$ will denote the corresponding real stratification of
$\BA_{\BR}$.

The stratification $\CS_{\emp}$ has a unique minimal stratum
\begin{equation}
\label{dva delta}
\Delta=\bigcap\ \Delta_{ij}
\end{equation}
--- main diagonal; it is one-dimensional.
We will denote by $^{\pi}\BAO_{\emp}$ (resp., $^{\pi}\BAO_{\emp,\BR}$)
the open stratum
of $\CS_{\emp}$ (resp., of $\CS_{\emp,\BR}$).

\subsection{}
\label{dva diagchamb} Let us describe the chambers of $\CS_{\emp,\BR}$.
If $C$ is a chamber and
$\bx=(x_j)\in C$, i.e. the embedding $J\hra\BR,\ j\mapsto x_j$,
it induces an obvious total order on $J$, i.e. a bijection
\begin{equation}
\label{dva tauc}
\tau_C:J\iso [N]
\end{equation}
Namely, $\tau_C$ is determined uniquely by the requirement
$\tau_C(i)<\tau_C(j)$ iff $x_i<x_j$; it does not depend on the choice
of $\bx$. This way we get  a one-to-one correspondence between
the set of chambers of $\CS_{\emp}$ and the set of all bijections
{}~(\ref{dva tauc}). We will denote by $C_{\tau}$ the chamber corresponding to
$\tau$.

Given $C$ and $\bx$ as above, suppose that we have $i,j\in J$ such that
$x_i<x_j$ and there is no $k\in J$ such that $x_i<x_k<x_j$. We will say that
$i,j$ are {\em neighbours} in $C$, more precisely that {\em $i$ is a left
neighbour of $j$}.

Let $\bx'=(x'_j)$ be a point with $x'_p=x_p$ for all $p\neq j$, and $x'_j$
equal to some number smaller than $x_i$ but greater than any $x_k$ such that
$x_k<x_i$.
Let $^{ji}C$ denote the chamber containing $\bx'$. Let us introduce a  homotopy
class of paths $^C\gamma_{ij}$ connecting $\bx$ and $\bx'$ as shown on
Fig. 1 below.

\begin{picture}(20,8)(-10,-4)

\put(0,0){\circle*{0.2}}
\put(0,-0.5){$i$}

\put(-4,0){\line(1,0){8}}

\put(1.5,0){\circle*{0.2}}
\put(1.5,-0.5){$j$}

\put(0,0){\oval(3,3)[t]}
\put(0,1.5){\vector(-1,0){0.5}}
\put(-0.5,2){$^C\gamma_{ij}$}

\put(-1.5,0){\circle*{0.2}}
\put(-1.5,-0.5){$j'$}

\put(3,0){\circle*{0.2}}
\put(-3,0){\circle*{0.2}}

\put(-0.5,-4){Fig. 1.}

\end{picture}

We can apply the discussion I.4.1
and consider the groupoid $\pi_1(^{\pi}\BAO_{\emp},^{\pi}\BAO_{\emp,\BR})$.
It has as the set of objects the set
of all chambers. The set of morphisms is generated by all morphisms
$^C\gamma_{ij}$
subject to certain evident braiding relations. We will need only the following
particular case.

To define a {\em one-dimensional} local system $\CL$ over $^{\pi}\BAO_{\emp}$
is
the same as to give a set of one-dimensional vector spaces
$\CL_C,\ C\in\pi_0(^{\pi}\BAO_{\emp,\BR})$,
together with arbitrary invertible linear operators
\begin{equation}
\label{dva halfmonodr}
^CT_{ij}:\CL_C\lra\CL_{^{ji}C}
\end{equation}
("half-monodromies")  defined for chambers having $i$ as a left neighbour of
$j$.

\subsection{}
\label{dva pici} We define a one-dimensional local system $^{\pi}\CI$ over
$^{\pi}\BAO_{\emp}$ as follows. Its fibers $^{\pi}\CI_{C}$
are one-dimensional
linear spaces with fixed basis vectors; they will be identified with $\sk$.

Half-monodromies are defined as
$$
^CT_{ij}=\zeta^{i\cdot j},\ i,j\in J
$$

\subsection{} Let $j:\ ^{\pi}\BAO_{\emp}\lra\ ^{\pi}\BA$ denote an open
embedding. We will study the
following objects of $\CM(^{\pi}\BA;\CS_{\emp})$:
$$
^{\pi}\CI_{?}=j_?^{\pi}\CI[N],
$$
where $?=!,*$. We have a canonical map
\begin{equation}
\label{dva map!*}
m:\ ^{\pi}\CI_!\lra\ ^{\pi}\CI_*
\end{equation}
and by definition $^{\pi}\CI_{!*}$ is its image, cf. I.4.5.

\subsection{}
\label{dva marking} For an integer $r$ let us denote by $\CP_r(J)$ the set of all
surjective mappings $J\lra [r]$. It is evident that $\CP_r(J)\neq\emp$ if and
only if $1\leq r\leq N$. To each $\rho\in\CP_r(J)$ let us assign a point
$w_{\rho}=(\rho(j))\in\ ^{\pi}\BA_{\BR}$. Let $F_{\rho}$ denote
the facet containing $w_{\rho}$. This way we get a bijection between
$\CP_r(J)$ and the set of $r$-dimensional facets. For $r=N$ we get the
bijection from ~\ref{dva diagchamb}.

At the same time we have defined a marking of $\CH_{\emp}$:
by definition, $^{F_{\rho}}w=w_{\rho}$. This defines cells $D_F,\ S_F$.

\subsection{}
The main diagonal $\Delta$ is a unique $1$-facet; it corresponds to the
unique element $\rho_0\in\CP_1(J)$.

We will denote by $\Ch$ the set of all chambers; it is the same as
$\Ch(\Delta)$ in notations of Part I.
Let $C_{\tau}$ be a chamber. The order $\tau$ identifies $C$ with
an open cone in the standard coordinate space $\BR^N$; we provide $C$
with the orientation induced from $\BR^N$.

\subsection{Basis in $\Phi_{\Delta}(^{\pi}\CI_*)$}
The construction
I.4.7 gives us the basis $\{c_{\Delta<C}\}$ in
$\Phi_{\Delta}(^{\pi}\CI_!)^*$ indexed
by $C\in\Ch$. We will use notation $c_{\tau,!}:=c_{\Delta<C_{\tau}}$.

A chain $c_{\tau,!}$ looks as follows.

\begin{picture}(20,6)(-10,-3)

\put(-4,0){\line(1,0){6}}

\put(-1,0){\vector(1,0){0.8}}
\put(-0.2,0){\circle*{0.15}}
\put(-0.3,-0.4){$j_1$}
\put(-1,0){\vector(1,0){1.6}}
\put(0.6,0){\circle*{0.15}}
\put(0.5,-0.4){$j_2$}

\put(2.2,0){$\ldots$}
\put(3,0){\line(1,0){2}}
\put(3,0){\vector(1,0){0.5}}
\put(3.5,0){\circle*{0.15}}
\put(3.5,-0.4){$j_N$}

\put(-1,-3){Fig. 2. A chain $c_{\tau,!}$.}

\end{picture}

Here $\tau(j_i)=i$. We will denote by $\{b_{\tau,!}\}$ the dual basis
in $\Phi_{\Delta}(\CI_!)$.

\subsection{Basis in $\Phi_{\Delta}(\CI_*)$}
Similarly, the definition I.4.9 gives us the basis $\{c_{\Delta<C}\},\
C\in\Ch$ in $\Phi_{\Delta}(\CI_*)^*$. We will use the notations
$c_{\tau,*}:=c_{\Delta<C_{\tau}}$.

If we specify
the definition I.4.9 and
its explanation I.4.12 to our arrangement, we get the following
picture for a dual chain $c_{\tau,*}$.

\begin{picture}(20,6)(-10,-3)

\put(-5,0){\line(1,0){9.5}}

\put(2,0){\circle*{0.15}}
\put(1.9,0.4){$j_1$}

\put(2,0){\oval(1.6,1.6)[b]}
\put(2,-0.8){\vector(1,0){0.3}}
\put(2.8,0){\circle*{0.15}}
\put(2.8,0.4){$j_2$}

\put(2,0){\oval(3.2,2.4)[b]}
\put(2,-1.2){\vector(1,0){0.8}}
\put(3.6,0){\circle*{0.15}}
\put(3.6,0.4){$j_3$}

\put(4.5,0){$\ldots$}
\put(5.5,0){\line(1,0){2}}

\put(2,0){\oval(10,4)[b]}
\put(2,-2){\vector(1,0){1.8}}
\put(7,0){\circle*{0.15}}
\put(7,0.4){$j_N$}

\put(-1,-3){Fig. 3. A chain $c_{\tau,*}$ .}

\end{picture}

This chain is represented by the section of a local system $\CI^{-1}$
over the cell in $^{\pi}\BAO_{\emp}$ shown above, which takes value $1$
at the point corresponding to the end of the travel (direction of travel
is shown by arrows).

To understand what is going on, it is instructive to treat the
case $N=2$ first, which essentially coincides with the Example I.4.10.

We will denote by $\{b_{\tau,*}\}$ the dual basis in $\Phi_{\Delta}(\CI_*)$.

\subsection{}
\label{dva signs} Obviously, all maps $\tau:J\lra [N]$ from $\CP_N(J)$ are
bijections. Given two such maps $\tau_1,\tau_2$, define the sign
$\sgn(\tau_1,\tau_2)=\pm 1$ as the sign of the permutation
$\tau_1 \tau_2^{-1}\in\Sigma_N$.

For any $\tau\in\CP_N(J)$ let us denote by $\vec{J}_{\tau}$ the sequence
$(\tau^{-1}(N),\tau^{-1}(N-1),\ldots,\tau^{-1}(1))$.

\subsection{} Let us pick $\eta\in\CP_N(J)$. Let us define
the following maps:
\begin{equation}
\label{dva phid!j}
^{\pi}\phi_{\Delta,!}^{(\eta)}:\Phi_{\Delta}(^{\pi}\CI_!)\lra\
^{\pi}{\ '\ff}_{\chi_J}
\end{equation}
which carries $b_{\tau,!}$ to $\sgn(\tau,\eta)\cdot\theta_{\vec{J}_{\tau}}$,
and
\begin{equation}
\label{dva phid*j}
^{\pi}\phi_{\Delta,*}^{(\eta)}:\Phi_{\Delta}(^{\pi}\CI_*)\lra\
^{\pi}{\ '\ff}^*_{\chi_J}
\end{equation}
which carries $b_{\tau,*}$ to
$\sgn(\tau,\eta)\cdot\theta^*_{\vec{J}_{\tau}}$.

\subsection{Theorem}
\label{dva diagiso} {\em (i) The maps $^{\pi}\phi_{\Delta,!}^{(\eta)}$ and
$^{\pi}\phi_{\Delta,*}^{(\eta)}$ are isomorphisms.
The square
$$\begin{array}{ccc}
\Phi_{\Delta}(^{\pi}\CI_!)&\overset{^{\pi}\phi_{\Delta,!}^{(\eta)}}{\iso}&
\ ^{\pi}{\ '\ff}_{\chi_J}\\
m\downarrow&\;&\downarrow S\\
\Phi_{\Delta}(^{\pi}\CI_*)&\overset{^{\pi}\phi_{\Delta,*}^{(\eta)}}{\iso}&
\ ^{\pi}{\ '\ff}^*_{\chi_J}
\end{array}$$
commutes.

(ii) The map $^{\pi}\phi_{\Delta,!}^{(\eta)}$ induces an isomorphism
\begin{equation}
\label{dva phi!*iso}
^{\pi}\phi_{\Delta,!*}^{(\eta)}:\Phi_{\Delta}(^{\pi}\CI_{!*})\iso\
^{\pi}\ff_{\chi_J}
\ \Box
\end{equation}}

{\bf Proof.} This theorem is particular case of I.14.16, I.4.17.
The claim about isomorphisms in (i) is clear. To prove the commutativity
of the square, we have to compute the action of the canonical map $m$
on our standard chains. The claim follows at once from their geometric
description given above. Note that here the sign in the definition
of morphisms $\phi$ is essential, due to orientations of our chains.
(ii) is a direct corollary of (i) $\Box$

\vspace{1cm}
{\em SYMMETRIZED CONFIGURATIONAL SPACES}
\vspace{1cm}

\subsection{Colored configuration spaces}
\label{dva color} Let us fix $\nu=\sum\nu_ii\in\BN[I]$, $\sum_i\nu_i=N$.
There exists a finite set $J$ and a morphism $\pi:J\lra I$ such that
$\card(\pi^{-1}(i))=\nu_i$ for all $i\in I$. Let us call such $\pi$
{\em an unfolding of $\nu$}. It is unique up to a non-unique isomorphism;
the automorphism group of $\pi$ is precisely $\Sigma_{\pi}$, and
$\nu=\nu_{\pi}$ in our previous notations.

Let us pick an unfolding $\pi$. As in the above discussion, we define
$^{\pi}\BA$ as a complex affine space with coordinates $t_j,\ j\in J$.
Thus, $\dim\ ^{\pi}\BA=N$.
The group
$\Sigma_{\pi}$ acts on the space $^{\pi}\BA$ by permutations
of coordinates.

Let us denote by $\CA_{\nu}$ the quotient manifold $^{\pi}\BA/\Sigma_{\pi}$.
As an
algebraic manifold, $\CA_{\nu}$ is also a complex $N$-dimensional affine
space. We have a canonical projection
\begin{equation}
\label{dva mappi}
\pi:\ ^{\pi}\BA\lra\CA_{\nu}
\end{equation}
The space $\CA_{\nu}$ does not depend on the choice of an unfolding $\pi$.
It will be called {\em the configuration space of $\nu$-colored points
on the affine line $\BA^1$}.

We will consider the stratification on $\CA_{\nu}$ whose strata are
$\pi(S),\ S\in\CS_{\emp}$; we will denote this stratification also by
$\CS_{\emp}$; this definition does not depend on the choice of $\pi$.
We will study the category $\CM(\CA_{\nu};\CS_{\emp})$.

We will denote by $\CAO_{\nu,\emp}$ the open stratum. It is clear that
$\pi^{-1}(\CAO_{\nu,\emp})=\ ^{\pi}\BAO_{\emp}$. The morphism $\pi$
is unramified over $\CAO_{\nu,\emp}$.

The action of $\Sigma_{\pi}$ on
$^{\pi}\BA$ may be extended in the evident
way to the local system $^{\pi}\CI$,
hence all our spaces of geometric origin ---
like $\Phi_{\Delta}(^{\pi}\CI_!)$, etc. ---
get an action of $\Sigma_{\pi}$.

\subsubsection{} If $M$ is an object with a $\Sigma_{\pi}$-action
(for example a vector
space or a sheaf), we will denote by $M^{\Sigma_{\pi},-}$ the subobject
$\{x\in M|\mbox{ for every }\sigma\in\Sigma_{\pi}\ \sigma x=\sgn(\sigma)x\}$
where $\sgn(\sigma)=\pm 1$ is the sign of a permutation.

A morphism $f:M\lra N$ between two objects with $\Sigma_{\pi}$-action will
be called {\em skew ($\Sigma_{\pi}$)-equivariant} if for any
$x\in M,\ \sigma\in\Sigma_{\pi},\ f(\sigma x)=\sgn(\sigma)\sigma f(x)$.

Let us define a local system over $\CA_{\nu}$
\begin{equation}
\label{dva cinu}
\CI_{\nu}=(\pi_*\ ^{\pi}\CI)^{\Sigma_{\pi},-}
\end{equation}

\subsection{}
\label{dva jcolor} Let $j:\ ^{\pi}\BAO_{\emp}\hra\ ^{\pi}\BA,\
j_{\CA}:\CAO_{\nu,\emp}\hra\CA_{\nu}$ be the open embeddings. Let us define
the following objects of $\CM(\CA_{\nu};\CS_{\emp})$:
\begin{equation}
\label{dva standca}
\CI_{\nu?}:=j_{\CA?}\CI_{\nu}[N]
\end{equation}
where $?=!,\ *$ or $!*$.
We have by definition
\begin{equation}
\label{dva inu!}
\CI_{\nu!}=(\pi_*\ ^{\pi}\CI_!)^{\Sigma_{\pi},-}
\end{equation}
The morphism $\pi$ is finite; consequently $\pi_*$ is $t$-exact
(see ~\cite{bbd}, 4.1.3) and commutes with the Verdier duality. Therefore,
\begin{equation}
\label{dva inu*}
\CI_{\nu*}=(\pi_*\ ^{\pi}\CI_*)^{\Sigma_{\pi},-};\
\CI_{\nu!*}=(\pi_*\ ^{\pi}\CI_{!*})^{\Sigma_{\pi},-}
\end{equation}

\subsection{} Let us define vector spaces
\begin{equation}
\label{dva phidelt}
\Phi_{\Delta}(\CI_{\nu?}):=(\Phi_{\Delta}(^{\pi}\CI_{?}))^{\Sigma_{\pi},-}
\end{equation}
where $?=!,*$ or $!*$.

Let us pick a $\Sigma_{\pi}$-equivariant marking of $\CH_{\emp}$, for
example the one from ~\ref{dva marking}; consider the corresponding
cells $D_{\Delta},\ S_{\Delta}$. It follows from ~(\ref{dva inu!}) and
{}~(\ref{dva inu*}) that
\begin{equation}
\label{dva relat}
\Phi_{\Delta}(\CI_{\nu?})=R\Gamma(\pi(D_{\Delta}),\pi(S_{\Delta});\
\CI_{\nu?})[-1]
\end{equation}
where $?=!,*$ or $!*$, cf. I.3.3.

\subsection{} The group $\Sigma_{\pi}$ is acting on
on $^{\pi}{\ '\ff}$. Let us pick $\eta\in\CP_N(J)$.
It follows from the definitions that
the isomorphisms $^{\pi}\phi_{\Delta,!}^{(\eta)},\
^{\pi}\phi_{\Delta,*}^{(\eta)}$ are skew
$\Sigma_{\pi}$-equivariant. Therefore,
passing to invariants in Theorem ~\ref{dva diagiso} we get

\subsection{Theorem}
\label{dva symdiag} {\em The maps $^{\pi}\phi_{\Delta,!}^{(\eta)},
^{\pi}\phi_{\Delta,*}^{(\eta)}$ induce isomorphisms included into a commutative
square
$$\begin{array}{ccc}
\Phi_{\Delta}(\CI_{\nu!})&
\overset{\phi_{\nu,!}^{(\eta)}}{\iso}&{\ '\ff}_{\nu}\\
m\downarrow&\;&\downarrow S\\
\Phi_{\Delta}(\CI_{\nu*})&\overset{\phi_{\nu,*}^{(\eta)}}{\iso}&{\ '\ff}^*_{\nu}
\end{array}$$
and
\begin{equation}
\label{dva phinu}
\phi_{\nu,!*}^{(\eta)}:\Phi_{\Delta}(\CI_{\nu!*})\iso\ff_{\nu}\ \ \ \Box
\end{equation}}

\section{Principal stratification}

The contents of this section is parallel to I, Section 3. However, we present
here certain modification of general constructions from {\em loc. cit.}

\subsection{}
\label{dva j} Let us fix a finite set $J$ of cardinality $N$. In this section
we will denote by $\BA_{\BR}$ a real affine space with fixed coordinates
$t_j:\BA_{\BR}\lra\BR,\ j\in J$, and by $\BA$ its complexification.
For $z\in \BC,\ i,j\in J$
denote by $H_j(z)\subset \BA$ a hyperplane $t_j=z$, and
by $\Delta_{ij}$ a hyperplane $t_i=t_j$.

Let us consider an arrangement
$\CH$ in $\Bbb A$ consisting of hyperplanes $H_i(0)$ and $\Delta_{ij}$,
$i,j\in J,\ i\neq j$. It is a complexification of an evident real arrangement
$\CH_{\BR}$ in $\BA_{\BR}$. As usual, the subscript $_{\BR}$ will denote
real points.

Denote by $\CS$ the corresponding stratification of $\BA$. To distinguish this
stratification from the diagonal stratification of the previous
section, we will call it {\em the principal stratification}.
To shorten the
notation, we will denote in this part by $\CCD(\CA,\CS)$ a category which would
be denoted $\CCD^b(\BA;\CS)$ in I. In this section
we will study the category $\CM(\BA;\CS)$.

\subsubsection{}
\label{dva posfacets} Let us consider a positive cone
$$
\BA^+_{\BR}=\{ (t_j)|\mbox{ all }t_j\geq 0\}\subset\BA_{\BR}
$$
A facet will be called {\em positive} if it lies inside
$\BA^+_{\BR}$.

A flag $\bF$ is called {\em positive} if all its facets
are positive.

\subsection{}
\label{dva poscells} Let us fix a marking $\bw=\{\ ^Fw\}$ of $\CH_{\BR}$
(cf. I.3.2).
For a positive facet $F$ define
$$
D^+_F=D_F\cap\BA_{\BR}^+;\ S^+_F=S_F\cap\BA_{\BR}^+;\
\DO^+_F=D_F^+-S_F^+.
$$
Note that $D^+_F$ coincides with the union of $^{\bF}\Delta$ over all
positive flags beginning at $F$, and $S^+_F$ coincides with the union
of $^{\bF}\Delta$ as above with $\dim\ ^{\bF}\Delta<\codim\ F$. It follows
that only marking points $^Fw$ for positive facets $F$ take part in the
definition of cells $D^+_F,\ S^+_F$.

\subsection{} Let $\CK$ be an object of $\CCD(\BA;\CS)$, $F$ a
positive facet of dimension $p$.
Let us introduce a notation
$$
\Phi_F^+(\CK)=\Gamma(D_F^+,S_F^+;\CK)[-p].
$$
This way we get a functor
\begin{equation}
\label{dva plus}
\Phi_F^+:\CCD(\BA;\ \CS)\lra\CCD^b(pt)
\end{equation}

\subsection{Theorem.}
\label{dva dual+} {\em Functors $\Phi_F^+$ commute
with Verdier duality. More precisely, we have canonical natural
isomorphisms
\begin{equation}
D\Phi_F^+(\CK)\iso\Phi_F^+(D\CK).
\end{equation}}

{\bf Proof} goes along the same lines as the proof of Theorem
I.3.5.

\subsection{}
\label{dva onedim} First let us consider the case $N=1$, cf. I.3.6. We will adopt
notations
from there and from I, Fig. 1. Our arrangement has one positive
$1$-dimensional facet $E=\BR_{>0}$, let $w\in E$ be a marking.

\begin{picture}(20,8)(-10,-4)

\put(0,0){\circle{0.2}}
\put(-0.5,-0.4){$F$}
\put(0,0){\oval(6,6)}
\put(-1.5,3.3){$S_{r''}$}

\put(0.1,0){\line(1,0){4}}
\put(0.5,0.3){$E$}

\put(0,3){\line(0,-1){2}}
\put(-0.4,1.8){$Y$}

\put(0,1){\circle*{0.2}}
\put(-0.5,0.5){$i\epsilon w$}

\put(2,0){\circle*{0.2}}
\put(2,-0.5){$w$}

\put(0,0){\oval(3,3)}
\put(-1.5,1.6){$S_{r'}$}

\put(-0.5,-4){Fig. 4}

\end{picture}

We have an isomorphism
\begin{equation}
\label{dva iso}
\Phi_{F}^+(\CK)\iso R\Gamma(\BA,\{ w\};\CK)
\end{equation}
Denote $j:=j_{\BA-\{ w\}}$. We have by Poincar\'{e} duality
\begin{equation}
D\Phi_F^+(\CK)\iso R\Gamma_c(\BA;j_*j^*D\CK)\iso
R\Gamma(\BA,\BA_{\geq r''};j_*j^*D\CK)
\end{equation}
Let us denote $D^{+opp}_F:=\BR_{\geq w}$, and $Y=\epsilon i\cdot D_F^{+opp}$.
We have
\begin{equation}
\label{dva homot}
R\Gamma(\BA,\BA_{\geq r''};j_*j^*D\CK)\iso
R\Gamma(\BA,Y\cup \BA_{\geq r''};j_*j^*D\CK)
\end{equation}
by homotopy. Consider the restriction map
\begin{equation}
\label{dva res}
res: R\Gamma(\BA,Y\cup \BA_{\geq r''};j_*j^*D\CK)\lra
R\Gamma(\BA_{\leq r'},Y\cap \BA_{\leq r'};D\CK)
\end{equation}
\subsubsection{}
\label{dva qiso} {\bf Claim.} {\em $res$ is an isomorphism.}

In fact, $\Cone (res)$ is isomorphic to
\begin{eqnarray}
R\Gamma(\BA,\BA_{\leq r'}\cup\BA_{\geq r''}\cup Y;j_*j^*D\CK)=
R\Gamma_c(\BA_{< r''},\BA_{\leq r'}\cup Y;j_*j^*D\CK)\cong\nonumber\\
\cong DR\Gamma(\BA_{< r''}-(\BA_{\leq r'}\cup Y);j_!j^*\CK)\nonumber
\end{eqnarray}
But
\begin{equation}
R\Gamma(\BA_{< r''}-(\BA_{\leq r'}\cup Y);j_!j^*\CK)=
R\Gamma(\BA_{< r''}-(\BA_{\leq r'}\cup Y),\{ w\};\CK)=0
\end{equation}
evidently. This proves the claim. $\Box$

A clockwise rotation by $\pi/2$ induces an isomorphism
$$
R\Gamma(\BA_{\leq r'},Y\cap \BA_{\leq r'};D\CK)\cong
R\Gamma(\BA_{\leq r'},\epsilon\cdot D_F^{opp};D\CK),
$$
and the last complex is isomorphic to $\Phi_F^+(D\CK)$. This proves the theorem
for $N=1$.

\subsection{}
\label{dva reduct} Now let us return to an arbitrary $J$. Let us prove
the theorem for $F$ equal to the unique
$0$-dimensional facet.

Let us introduce the following subspaces of
$\BA_{\BR}$ (as usually, a circle on the top will denote the interior).

$D_F^{opp}:=\BA_{\BR}-\DO_F$;
$D_F^{+opp}:=\BA^+_{\BR}-\DO_F^+$;
for each cell $D_{F<C},\ C\in \Ch(F)$ define $D_{F<C}^{opp}
:=C-\DO_{F<C}$.

It is easy to
see that the restriction induces isomorphism
$$
\Phi_F^+(\CK)\iso R\Gamma(\BA,D^{+opp}_F;\CK).
$$

We use the notations of I.3.8. Let us choose positive numbers
$r'<r'',\ \epsilon,$ such that
$$
\epsilon D_F\subset \BA_{<r'}\subset \DO_F\subset D_F\subset
\BA_{< r''}
$$
Define the subspace
$$
Y^+=\epsilon i\cdot D_F^{+opp};
$$
denote $j:=j_{\BA-S_F^+}$.
We have isomorphisms
\begin{eqnarray}
D\Phi_{F}^+(\CK)\cong DR\Gamma(\BA,S_{F}^+;\CK)
\cong R\Gamma_c(\BA; j_{*}j^*D\CK)\cong\\ \nonumber
\cong R\Gamma(\BA,\BA_{\geq r''};j_{*}j^*D\CK)\cong
R\Gamma(\BA,Y^+\cup \BA_{\geq r''};j_{*}j^*D\CK) \nonumber
\end{eqnarray}
Consider the restriction map
\begin{equation}
res: R\Gamma(\BA,Y^+\cup \BA_{\geq r''};j_*j^*D\CK)\lra
R\Gamma(\BA_{\leq r'},Y^+\cap \BA_{\leq r'};D\CK)
\end{equation}
$\Cone (res)$ is isomorphic to
\begin{eqnarray}
R\Gamma(\BA,\BA_{\leq r'}\cup\BA_{\geq r''}\cup Y^+;j_*j^*D\CK)=
R\Gamma_c(\BA_{< r''},\BA_{\leq r'}\cup Y^+;j_*j^*D\CK)\cong\nonumber\\
\cong DR\Gamma(\BA_{< r''}-(\BA_{\leq r'}\cup Y^+);j_!j^*\CK)=
R\Gamma(\BA_{(r',r'')}-Y^+,S_F^+;\CK)\nonumber
\end{eqnarray}

\subsubsection{}
\label{dva acycl+} {\bf Lemma.} (Cf. I.3.8.1.) {\em
$R\Gamma(\BA_{(r',r'')}-Y^+,S_F^+;\CK)=0$.}

{\bf Proof.} Let us define the following subspaces of $\BA$.

$A:=\{(t_j)|\mbox{ for all }j\ |t_j|<1;\ \mbox{ there exists }j: t_j\neq 0\}$;
$A_{\BR}^+:=A\cap\BA^+_{\BR}$. Note that $A_{\BR}^+\cap i\cdot A_{\BR}^+=
\emp$. Due to monodromicity, it is easy to see that
$$
R\Gamma(\BA_{(r',r'')}-Y^+,S_F^+;\CK)\cong
R\Gamma(A-i\cdot A_{\BR}^+,A_{\BR}^+;\CK).
$$
Therefore, it is enough to prove the following

\subsubsection{} {\bf Claim.} {\em The restriction map
\begin{equation}
\label{dva resmain}
R\Gamma(A-i\cdot A_{\BR}^+;\CK)\lra R\Gamma(A_{\BR}^+;\CK)
\end{equation}
is an isomorphism.}

{\bf Proof of the Claim.}
Let us introduce for each $k\in J$ open subspaces
$$
A_k=\{(t_j)\in A|t_k\not\in i\cdot\BR_{\geq 0}\}\subset A-i\cdot A^+_{\BR}
$$
and
$$
A'_k=\{(t_j)\in A_{\BR}^+|t_k>0\}\subset A_{\BR}^+
$$
Obviously $A'_k=A_k\cap A_{\BR}^+$.
For each subset $M\subset J$ set $A_M:=\bigcap_{k\in M}A_k;\
A'_M:=\bigcap_{k\in M}A'_k$.

For each non-empty $M$ define the spaces
$B_M:=\{(t_j)_{j\in M}|\mbox{ for all }j\ |t_j|<1,
t_j\not\in i\cdot\BR_{\geq 0}\}$ and
$B'_M:=\{(t_j)_{j\in M}|\mbox{ for all }j\ t_j\in\BR,\ 0<t_j<1\}$.
We have obvious projections
$f_M:A_M\lra B_M,\ f'_M:A'_M\lra B'_M$.

Let us look at fibers of $f_M$ and $f'_M$.
Given $b=(t_j)_{j\in M}\in B_M$, the fiber
$f^{-1}_M(b)$ is by definition $\{(t_k)_{k\in J-M}||t_k|<1\}$,
the possible singularities
of our sheaf $\CK$ are at the hyperplanes $t_k=t_j$ and $t_k=0$.
It is easy to see that $f_M$ is "lisse" with respect to $\CK$ which
means in particular that a stalk $(f_{M*}\CK)_b$ is equal to
$R\Gamma(f^{-1}_M(b);\CK)$. The same considerations apply to $f'_M$.
Moreover, it follows from I.2.12 that the restriction maps
$$
R\Gamma(f^{-1}_M(b);\CK)\lra R\Gamma((f'_M)^{-1}(b);\CK)
$$
are isomorphisms for every $b\in B'_M$. This implies that
$f'_{M*}\CK$ is equal to the restriction of $f_{M*}\CK$ to $B'_M$.

The sheaf $f_{M*}\CK$ is smooth along the diagonal stratification.
For a small $\delta>0$ let $U_{\delta}\subset B_M$ denote
an open subset $\{(t_j)\in B_M|\ |\arg(t_j)|<\delta\mbox{ for all }j\}$.
The restriction maps $R\Gamma(B_M;f_{M*}\CK)\lra R\Gamma(U_{\delta};
f_{M*}\CK)$
are isomorphisms. We have $B'_M=\bigcap_{\delta}U_{\delta}$, therefore
by I.2.12 the restriction
$R\Gamma(B_M;f_{M*}\CK)\lra R\Gamma(B_M';f_{M*}\CK)$ is an isomorphism.
This implies,
by Leray, that restriction maps
$$
R\Gamma(A_M;\CK)\lra R\Gamma(A'_M;\CK)
$$
are isomorphisms for every non-empty $M$.

Obviously $A-i\cdot A^+_{\BR}=\bigcup_{k\in J}A_k$ and
$A^+_{\BR}=\bigcup_{k\in J}A'_k$. Therefore,
by Mayer-Vietoris the map ~(\ref{dva resmain}) is an isomorphism.
This proves the claim, together with the lemma.
$\Box$

A clockwise rotation by $\pi/2$ induces an isomorphism
$$
R\Gamma(\BA_{\leq r'},Y^+\cap \BA_{\leq r'};D\CK)\cong
R\Gamma(\BA_{\leq r'},\epsilon\cdot D_F^{+opp};D\CK)\cong
R\Gamma(\BA_{\leq r'},\epsilon\cdot S^+_F;D\CK)
$$
and the last complex is isomorphic to $\Phi_F^+(D\CK)$ in view of
{}~\ref{dva reduct}. This proves the theorem for the case of the $0$-facet $F$.

\subsection{} Suppose that $F$ is an arbitrary positive facet.
{}From the description of positive facets (see {\em infra},
{}~\ref{dva listpos}) one sees that the cell $D^+_F$ is homeomorphic to a
cartesian product of the form
$$
D^+_{F_0}\times D^+_{F_1}\times\ldots\times D^+_{F_a}
$$
where $F_0$ is a $0$-facet of the principal
stratification in some affine space of smaller dimension, and
$D^+_{F_i}$ are the cells of {\em the diagonal} stratification
discussed in the previous section.

Using this remark, we apply a combination of the arguments of the previous
subsection (to the first factor) and of I.3.8 (to the remaining factors)
to get the required isomorphism. We leave details to the reader.

The theorem is proved. $\Box$

\subsection{} {\bf Theorem.} {\em All functors $\Phi^+_F$ are $t$-exact.
In other words,
for all positive facets $F$,
$$
\Phi^+_F(\CM(\BA;\CS))\subset\Vect\subset\CCD^b(pt).
$$}

{\bf Proof.} The same as that of I.3.9. $\Box$.

\subsection{} Thus we get exact functors
\begin{equation}
\label{dva phi+}
\Phi^+:\CM(\CA;\CS)\lra\Vect
\end{equation}
commuting with Verdier duality.

We will also denote vector spaces $\Phi^+_F(\CM)$ by $\CM^+_F$.

\subsection{Canonical and variation maps} Suppose we have a positive facet $E$.
Let us denote by $^+\Fac^1(E)$ the set of all positive facets $F$ such that
$E<F$, $\dim\ F=\dim\ E+1$. We have
\begin{equation}
\label{dva unican}
S^+_E=\bigcup_{F\in\ ^+\Fac^1(E)}D^+_F
\end{equation}

Suppose we have $\CK\in\CCD(\BA,\CS)$.

\subsubsection{}
\label{dva addcan} {\bf Lemma.} {\em We have a natural isomorphism
$$
R\Gamma(S^+_E,\bigcup_{F\in\ ^+\Fac^1(E)}S^+_F;\CK)\cong
\oplus_{F\in\ ^+\Fac^1(E)}\ R\Gamma(D^+_F,S^+_F;\CK)
$$}

{\bf Proof.} Note that $S^+_E-\bigcup_{F\in\ ^+\Fac^1(E)}S^+_F=
\bigcup_{F\in\ ^+\Fac^1(E)}\ \DO^+_F$ (disjoint union). The claim follows
now from the
Poincar\'{e} duality. $\Box$

Therefore, for any $F\in\ ^+\Fac^1(E)$ we get a natural inclusion map
\begin{equation}
\label{dva projcan}
i^F_E: R\Gamma(D^+_F,S^+_F;\CK)\hra
R\Gamma(S^+_E,\bigcup_{F'\in\Fac^1(E)}S^+_{F'};\CK)
\end{equation}

Let us define a map
$$
u^F_E(\CK):\Phi^+_F(\CK)\lra\Phi^+_E(\CK)
$$
as a composition
\begin{eqnarray}
R\Gamma(D^+_F,S^+_F;\CK)[-p]\overset{i^F_E}{\lra}
R\Gamma(S^+_E,\bigcup_{F'\in\ ^+\Fac^1(E)}S^+_{F'};\CK)[-p]\lra\\ \nonumber
R\Gamma(S^+_E;\CK)[-p]\lra
R\Gamma(D^+_E,S^+_E)[-p+1]\nonumber
\end{eqnarray}
where the last arrow is the coboundary map for the couple $(S^+_E,D^+_E)$,
and the second one is evident.

This way we get a natural transormation
\begin{equation}
\label{dva can}
^+u^F_E:\Phi^+_F\lra \Phi^+_E
\end{equation}
which will be called a {\em canonical map}.

We define a {\em variation morphism}
\begin{equation}
\label{dva var}
^+v^E_F:\Phi^+_E\lra\Phi^+_F
\end{equation}
as follows. By definition, $^+v^E_F(\CK)$ is the map dual to
the composition
$$
D\Phi^+_F(\CK)\iso\Phi^+_F(D\CK)\overset{^+u^F_E(D\CK)}{\lra}\Phi^+_E(D\CK)
\iso D\Phi^+_E(\CK).
$$

\subsection{Cochain complexes}
For each $r\in [0,N]$ and $\CM\in\CM(\BA,\CS)$
introduce vector spaces
\begin{equation}
^+C^{-r}(\BA;\CM)=\oplus_{F: F \mbox{ positive, }\dim\ F=r}\ \CM^+_F
\end{equation}
For $r>0$ or $r<-N$ set $^+C^r(\BA;\CM)=0$.

Define operators
$$
d:\ ^+C^{-r}(\BA;\CM)\lra\ ^+C^{-r+1}(\BA;\CM)
$$
having components $^+u^F_E$.

\subsubsection{}
\label{dva nilp+} {\bf Lemma.} {\em $d^2=0$}.

{\bf Proof.} The same as that of I.3.13.1. $\Box$

This way we get a complex $^+C^{\bullet}(\BA;\CM)$ lying in
degrees from $-N$ to $0$. It will be called the {\em
complex of positive cochains} of the sheaf $\CM$.

\subsection{Theorem}
\label{dva rgamma+} {\em (i) A functor
$$
\CM\mapsto \ ^+C^{\bullet}(\BA;\CM)
$$
is an exact functor from $\CM(\BA;\CS)$ to the category of complexes
of vector spaces.

(ii) We have a canonical natural isomorphism in
$\CCD(\{ pt\})$
$$
^+C^{\bullet}(\BA;\CM)\iso R\Gamma(\BA;\CM)
$$}

{\bf Proof.} One sees easily that restriction maps
$$
R\Gamma(\BA;\CM)\lra R\Gamma(\BA_{\BR}^+;\CM)
$$
are isomorphisms. The rest of the argument is the same as in I.3.14. $\Box$

\subsection{} Let us consider a function $\sum_jt_j:\BA\lra\BA^1$,
and the corresponding vanishing cycles functor
\begin{equation}
\label{dva vanish}
\Phi_{\Sigma\ t_j}:\CCD^b(\BA)\lra\CCD^b(\BA_{(0)})
\end{equation}
where $\BA_{(0)}=\{(t_j)|\sum_j t_j=0\}$, cf. ~\cite{ks}, 8.6.2.

If $\CK\in\CCD(\BA;\CS)$, it is easy to see that the complex
$\Phi_{\Sigma\ t_j}(\CK)$ has the support at the origin. Let us denote
by the same letter its stalk at the origin --- it is a complex
of vector spaces.

\subsubsection{} {\bf Lemma.} {\em We have a natural
isomorphism
\begin{equation}
\label{dva vanishin}
\Phi_{\Sigma\ t_j}(\CK)\iso\Phi^+_{\{0\}}(\CK)
\end{equation}}

{\bf Proof} is left to the reader. $\Box$

\subsection{}
\label{dva mainphi} Let us consider the setup of ~\ref{dva color}. Let us denote
by the same letter $\CS$ the stratification of $\CA_{\nu}$ whose strata
are subspaces $\pi(S)$, $S$ being a stratum of the stratification $\CS$
on $^{\pi}\BA$. This stratification will be called
{\em the principal stratification of $\CA_{\nu}$}.

The function $\Sigma\ t_j$ is obviously $\Sigma_{\pi}$-equivariant,
therefore it induces the function for which we will use the same notation,
\begin{equation}
\label{dva sum}
\Sigma\ t_j:\CA_{\nu}\lra\BA^1
\end{equation}
Again, it is easy to see that for $\CK\in\CCD^b(\CA_{\nu};\CS)$ the complex
$\Phi_{\Sigma t_j}(\CK)$ has the support at the origin.
Let us denote by $\Phi_{\nu}(\CK)$ its stalk at the origin.

It is known that the functor of vanishing cycles is $t$-exact
with respect to the middle perversity; whence we get an exact functor
\begin{equation}
\label{dva phi0}
\Phi_{\nu}:\CM(\CA_{\nu};\CS)\lra\Vect
\end{equation}

\subsubsection{} {\bf Lemma.} {\em Suppose that $\CN$ is a $\Sigma_{\pi}$-
equivariant complex of sheaves over $^{\pi}\BA$ which belongs to
to $\CCD^b(^{\pi}\BA;\CS)$ (after forgetting $\Sigma_{\pi}$-action).
We have a
natural isomorphism
\begin{equation}
\label{dva vanishinv}
\Phi_{\nu}((\pi_*\CN)^{\Sigma_{\pi},-})\iso(\Phi_{\Sigma\ t_j}(\CN))
^{\Sigma_{\pi},-}
\end{equation}}

{\bf Proof} follows from the proper base change for vanishing cycles
(see ~\cite{d3}, 2.1.7.1) and
the exactness of the functor $(\cdot)^{\Sigma_{\pi},-}$.
We leave details to the reader. $\Box$

\subsubsection{} {\bf Corollary.} {\em For a $\Sigma_{\pi}$-equivariant
sheaf $\CN\in\CM(^{\pi}\BA_{\nu};\CS)$ we have
a natural isomorphism of vector spaces
\begin{equation}
\label{dva vanishface}
\Phi_{\nu}(\CM)\iso(\Phi^+_{\{0\}}(\CN))^{\Sigma_{\pi},-}\
\end{equation}
where $\CM=(\pi_*\CN)^{\Sigma_{\pi},-}$. $\Box$}

\section{Standard sheaves}
\label{dva standsheaves}

Let us keep assumptions and notations of ~\ref{dva sym}.

\subsection{}
\label{dva group}
Let us denote by $^{\pi}\BA$ the complex affine space with coordinates $t_j,\
j\in J$. We will consider its principal stratification as in ~\ref{dva j}.

Suppose we are given a positive chamber $C$ and a point
$\bx=(x_j)_{j\in J}\in C$.
There extists a unique bijection
\begin{equation}
\label{dva bij}
\sigma_C:J\iso [N]
\end{equation}
such that for any $i,j\in J$, $\sigma_C(i)<\sigma_C(j)$ iff
$x_i<x_j$. This bijection does not depend upon the choice
of $\bx$.
This way we can identify the set of all positive chambers
with the set of all isomorphisms ~(\ref{dva bij}),
or, in other words, with the set of all total orderings of $J$.

Given $C$ and $\bx$ as above, suppose that we have $i,j\in J$ such that
$x_i<x_j$ and there is no $k\in J$ such that $x_i<x_k<x_j$. We will say that
$i,j$ are {\em neighbours} in $C$, more precisely that $i$ is a left neighbour
of $j$.

Let $\bx'=(x'_j)$ be a point with $x'_p=x_p$ for all $p\neq j$, and $x'_j$
equal to some number smaller than $x_i$ but greater than any $x_k$ such that
$x_k<x_i$.
Let $^{ji}C$ denote the chamber containing $\bx'$. Let us introduce a  homotopy
class of paths $^C\gamma_{ij}$ connecting $\bx$ and $\bx'$ as shown on
Fig. 5(a) below.

On the other hand, suppose that $i$ and $0$ are neighbours in $C$, there
is no $x_j$ between $0$ and $x_i$. Then we introduce the homotopy class
of paths from $\bx$ to itself as shown on Fig. 5 (b).

\begin{picture}(20,8)(-10,-4)

\put(-5,0){\circle*{0.2}}
\put(-5,-0.5){$i$}

\put(-9,0){\line(1,0){8}}

\put(-3.5,0){\circle*{0.2}}
\put(-3.5,-0.5){$j$}

\put(-5,0){\oval(3,3)[t]}
\put(-5,1.5){\vector(-1,0){0.5}}
\put(-5.5,2){$^C\gamma_{ij}$}

\put(-6.5,0){\circle*{0.2}}
\put(-6.5,-0.5){$j'$}

\put(-2,0){\circle*{0.2}}
\put(-8,0){\circle{0.2}}
\put(-8,-0.5){$0$}

\put(-5.5,-3){$(a)$}


\put(5,0){\circle*{0.2}}
\put(5.3,-0.5){$i$}

\put(1,0){\line(1,0){8}}

\put(3.5,0){\circle{0.2}}
\put(3.5,-0.5){$0$}
\put(3.5,0){\oval(3,3)}
\put(3.5,1.5){\vector(-1,0){0.5}}
\put(3.4,2){$^C\gamma_{i0}$}

\put(6.5,0){\circle*{0.2}}

\put(8,0){\circle*{0.2}}

\put(4.5,-3){$(b)$}


\put(-0.5,-4){Fig. 5}

\end{picture}

All chambers are contractible. Let us denote by $^{\pi}\BAO^+_{\BR}$ the union
of all positive chambers.
We can apply the discussion I.4.1
and consider the groupoid $\pi_1(^{\pi}\BAO,\ ^{\pi}\BAO^+_{\BR})$.
It has as the set of objects the set
of all positive chambers. The set of morphisms is generated by all morphisms
$^C\gamma_{ij}$ and $^C\gamma_{i0}$
subject to certain evident braiding relations. We will need only the following
particular case.

To define a {\em one-dimensional} local system $\CL$ over $^{\pi}\BAO$ is
the same as to give a set of one-dimensional vector spaces
$\CL_C,\ C\in\pi_0(^{\pi}\BAO^+_{\BR})$,
together with arbitrary invertible linear operators
\begin{equation}
\label{dva halfmon}
^CT_{ij}:\CL_C\lra\CL_{^{ji}C}
\end{equation}
("half-monodromies")  defined for chambers where $i$ is a left neighbour of
$j$ and
\begin{equation}
\label{dva mon}
^CT_{i0}:\CL_C\lra\CL_C
\end{equation}
defined for chambers with neighbouring $i$ and $0$.

\subsection{} Let us fix a weight $\Lambda\in X$.
We define a one-dimensional local system $\CI(^{\pi}\Lambda)$ over
$^{\pi}\BAO$ as follows. Its fibers $\CI(^{\pi}\Lambda)_{C}$ are
one-dimensional
linear spaces with fixed basis vectors; they will be identified with $\sk$.

Monodromies are defined as
$$
^CT_{ij}=\zeta^{i\cdot j},\
^CT_{i0}=\zeta^{-2\langle \Lambda,\pi(i)\rangle},
$$
for $i,j\in J$.

\subsection{} Let $j:\ ^{\pi}\BAO\lra\ ^{\pi}\BA$ denote an open embedding.
We will study the
following objects of $\CM(^{\pi}\BA;\CS)$:
$$
\CI(^{\pi}\Lambda)_{?}=j_?\CI(^{\pi}\Lambda)[N],
$$
where $?=!,*$. We have a canonical map
\begin{equation}
\label{dva !*}
m: \CI(^{\pi}\Lambda)_!\lra\CI(^{\pi}\Lambda)_*
\end{equation}
and by definition $\CI(^{\pi}\Lambda)_{!*}$ is its image, cf. I.4.5.

\vspace{1cm}
{\em COMPUTATIONS FOR $\CI(^{\pi}\Lambda)_!$}
\vspace{0.8cm}

\subsection{}
\label{dva listpos} We will use the notations ~\ref{dva gh} with $I=J$ and $n=1$.
For each $r\in [0,N]$, let us assign to a map $\varrho\in\CP_r(J;1)$ a point
$w_{\varrho}=(\varrho(j))_{j\in J}\in\ ^{\pi}\BA_{\BR}$. It is easy to see that
there exists a unique positive facet $F_{\varrho}$ containing $w_{\varrho}$,
and the rule
\begin{equation}
\label{dva rhofac}
\varrho\mapsto\ F_{\varrho}
\end{equation}
establishes an isomorphism between $\CP_r(J;1)$ and the set of all
positive facets of $\CS_{\BR}$. Note that $\CP_0(J;1)$ consists
of one element --- the unique map $J\lra [0]$; our stratification
has one zero-dimensional facet.

At the same time we have picked a point $^{F_{\varrho}}w:=w_{\varrho}$ on
each positive facet $F_{\varrho}$; this defines cells $D^+_F,\ S^+_F$
(cf. the last remark in ~\ref{dva poscells}).

\subsection{}
\label{dva listcouples} Given $\varrho\in\CP_r(J;1)$ and $\tau\in\CP_N(J;1)$,
it is easy to see that the chamber $C=F_{\tau}$ belongs to
$\Ch(F_{\varrho})$ iff $\tau$ is a refinement of $\varrho$ in the sense of
{}~\ref{dva refine}. This defines a bijection between the set of all refinements
of $\varrho$ and the set of all positive chambers containing $F_{\varrho}$.

We will denote the last set by $\Ch^+(F_{\varrho})$.

\subsubsection{Orientations} Let $F=F_{\varrho}$ be a positive facet and
$C=F_{\tau}\in\Ch^+(F)$. The map $\tau$ defines an isomorphism denoted
by the same letter
\begin{equation}
\label{dva tau}
\tau:J\iso [N]
\end{equation}
Using $\tau$, the natural order on $[N]$ induces a total order on $J$.
For each $a\in [r]$, let $m_a$ denote the minimal element
of $\varrho^{-1}(a)$, and set
$$
J'=J_{\varrho\leq \tau}:=J-\{ m_1,\ldots, m_r\}\subset J
$$
Let us consider the map
$$
(x_j)\in D_{F<C}\mapsto \{ x_j-m_{\varrho(j)}|j\in J'\}\in\BR^{J'}.
$$
It is easy to see that this mapping establishes an isomorphism
of the germ of the cell $D_{F<C}$ near the point $^Fw$ onto a germ
of the cone
$$
\{ 0\leq u_1\leq\ldots\leq u_{N-r}\}
$$
in $\BR^{J'}$ where we have denoted for a moment by $(u_i)$ the coordinates
in $\BR^{J'}$ ordered by the order induced from $J$.

This isomorphism together with the above order defines an orientation
on $D_{F<C}$.

\subsection{Basis in $\Phi_F^+(\CI(^{\pi}\Lambda)_!)^*$}
\label{dva basis!} We follow the pattern of
I.4.7. Let $F$ be a positive facet of dimension $r$. By definition,
$$
\Phi^+_F(\CI(^{\pi}\Lambda)_!)=H^{-r}(D^+_F,S^+_F;\CI(^{\pi}\Lambda)_!)=
H^{N-r}(D^+_F,S^+_F;j_!\CI(^{\pi}\Lambda))\cong
H^{N-r}(D^+_F,S^+_F\cup (_{\CH}H_{\BR}\cap D^+_F);j_!\CI(^{\pi}\Lambda)).
$$
Note that we have
$$
D^+_F-(S^+_F\cup\ _{\CH}H_{\BR})=\bigcup_{C\in\Ch^+(F)}\DO_{F<C}
$$
(disjoint union), therefore by additivity
$$
H^{N-r}(D^+_F,S^+_F\cup (_{\CH}H_{\BR}\cap D^+_F);j_!\CI(^{\pi}\Lambda))\cong
\oplus_{C\in\Ch^+(F)}
H^{N-r}(D_{F<C},S_{F<C};j_!\CI(^{\pi}\Lambda)).
$$
By Poincar\'{e} duality,
$$
H^{N-r}(D_{F<C},S_{F<C};j_!\CI(^{\pi}\Lambda))^*\cong
H^0(\DO_{F<C};\CI(^{\pi}\Lambda)^{-1})
$$
--- here we have used the orientations of cells $D_{F<C}$ introduced
above. By definition of the local system $\CI$, the last space is
canonically identified with $\sk$.

If $F=F_{\varrho},\ C=F_{\tau}$, we will denote by $c_{\varrho\leq \tau}
\in\Phi^+_F(\CI(^{\pi}\Lambda)_!)^*$
the image of $1\in H^0(\DO_{F<C};\CI(^{\pi}\Lambda)^{-1})$. Thus the
chains $c_{\varrho\leq\tau},\ \tau\in\Ord(\varrho)$, form
a basis of the space $\Phi^+_F(\CI(^{\pi}\Lambda)_!)^*$.

\subsection{Diagrams} It is convenient to use the following diagram
notations for chains $c_{\varrho\leq\tau}$.

Let us denote elements of $J$ by letters $a,b,c,\ldots$. An $r$-dimensional
chain $c_{\varrho\leq\tau}$ where $\varrho:J\lra [0,r]$, is represented
by a picture:

\begin{picture}(20,6)(-10,-3)

\put(-8,0){\line(1,0){10}}
\put(3,0){$\ldots$}
\put(4,0){\line(1,0){4}}

\put(-6,0){\circle{0.2}}
\put(-6.1,0.4){$0$}
\put(-6,0){\vector(1,0){0.8}}
\put(-5.2,0){\circle*{0.1}}
\put(-5.3,-0.4){$a$}

\put(-3,0){\circle{0.2}}
\put(-3.1,0.4){$1$}
\put(-3,0){\vector(1,0){0.8}}
\put(-2.2,0){\circle*{0.1}}
\put(-2.3,-0.4){$b$}
\put(-3,0){\vector(1,0){1.6}}
\put(-1.4,0){\circle*{0.1}}
\put(-1.5,-0.4){$c$}

\put(5,0){\circle{0.2}}
\put(4.9,0.4){$r$}
\put(5,0){\vector(1,0){0.8}}
\put(5.8,0){\circle*{0.1}}
\put(5.7,-0.4){$d$}
\put(5,0){\vector(1,0){1.6}}
\put(6.6,0){\circle*{0.1}}
\put(6.5,-0.4){$e$}

\put(-1,-3){Fig. 6. Chain $c_{\varrho\leq\tau}$.}

\end{picture}

A picture consists of $r+1$ fragments:

\begin{picture}(20,6)(-10,-3)

\put(-4,0){\line(1,0){8}}

\put(-1,0){\circle{0.2}}
\put(-1.1,0.4){$i$}
\put(-1,0){\vector(1,0){0.8}}
\put(-0.2,0){\circle*{0.1}}
\put(-0.3,-0.4){$a$}
\put(-1,0){\vector(1,0){1.6}}
\put(0.6,0){\circle*{0.1}}
\put(0.5,-0.4){$b$}
\put(-1,0){\vector(1,0){2.4}}
\put(1.4,0){\circle*{0.1}}
\put(1.3,-0.4){$c$}

\put(-1,-3){Fig. 7. Set $\varrho^{-1}(i)$.}

\end{picture}

where $i=0,\ldots, r$, the $i$-th fragment being a blank circle with
a number of small vectors going from it. These vectors are
in one-to-one correspondence with the set $\varrho^{-1}(i)$; their ends
are labeled by elements of this set. Their order
(from left to right) is determined by the order on $\varrho^{-1}(i)$
induced by $\tau$. The point $0$ may have no vectors (when $\varrho^{-1}(0)=
\emp$); all other points have at least one vector.

\subsection{} Suppose we have $\varrho\in\CP_r(J;1),\
\varrho'\in\CP_{r+1}(J;1)$. It is easy to see that $F_{\varrho}<
F_{\varrho'}$ if and only if there exists $i\in [0,r]$ such that
$\varrho=\delta_i\circ\varrho'$ where
$\delta_i:[0,r+1]\lra [0,r]$ carries $a$ to $a$ if $a\leq i$ and to $a-1$
if $a\geq i+1$. We will write in this case that $\varrho<\varrho'$.

Let us compute the dual to the canonical map
\begin{equation}
\label{dva can!}
^+u^*: \Phi^+_{F_{\varrho}}(\CI(^{\pi}\Lambda)_!)^*\lra
\Phi^+_{F_{\varrho'}}(\CI(^{\pi}\Lambda)_!)^*.
\end{equation}
Suppose we have $\tau\in\Ord(\varrho')$; set $C=F_{\tau}$, thus
$F_{\varrho}<F_{\varrho'}<C$. Let us define the sign
\begin{equation}
\label{dva signum}
\sgn (\varrho<\varrho'\leq\tau)=
(-1)^{\sum_{j=i+1}^{r+1}\card((\varrho')^{-1}(j)-1)}
\end{equation}
This sign has the following geometrical meaning. The cell
$D_{F_{\varrho'<C}}$ lies in the boundary of $D_{F_{\varrho<C}}$.
We have oriented these cells above. Let us define the compatibility
of these orientations as follows.
Complete an orienting basis of the smaller cell
by a vector directed outside the larger cell --- if we get the orientation
of the larger cell, we say that the orientations are compatible, cf. I.4.6.1.

It is easy to see from the definitions that the sign ~(\ref{dva signum}) is
equal $+1$ iff the orientations of $D_{F_{\varrho'<C}}$ and
$D_{F_{\varrho<C}}$ are compatible. As a consequence, we get

\subsubsection{} {\bf Lemma.}
\label{dva matrixu} {\em The map ~(\ref{dva can!}) has the following
matrix:
$$
^+u^*(c_{\varrho\leq\tau})=\sum\sgn(\varrho<\varrho'\leq\tau)c_{\varrho'\leq
\tau},
$$
the summation over all $\varrho'$ such that $F_{\varrho}<F_{\varrho'}<
F_{\tau}$ and $\dim F_{\varrho'}=\dim F_{\varrho}+1$.} $\Box$

\subsection{Isomorphisms $^{\pi}\phi_{!}$} We will use notations
of ~\ref{dva refine}
with $I$ replaced by $J$, ${\ '\ff}$ by $^{\pi}{\ '\ff}$, with $n=1$ and
$\Lambda_0=\ ^{\pi}\Lambda$.

Thus, for any $r\in [0,N]$ the set
$\{ b_{\varrho\leq\tau}|\varrho\in\CP_r(J;1),\
\tau\in\Ord(\varrho)\}$ is a basis of
$^+C^{-r}(^{\pi}\BA;\CI(^{\pi}\Lambda)_!)$.

Let us pick $\eta\in\CP_N(J)$. Any $\tau\in\CP_N(J;1)$ induces the bijection
$\tau':J\iso [N]$. We will denote by $\sgn(\tau,\eta)=\pm 1$ the
sign of the permutation $\tau'\eta^{-1}$.

Let us denote by
$\{ b_{\varrho\leq\tau}|\tau\in\Ord(\varrho)\}$
the basis in
$\Phi_{F_{\varrho}}^+(\CI(^{\pi}\Lambda)_!)$ dual to
$\{ c_{\varrho\leq\tau}\}$.

Let us define isomorphisms
\begin{equation}
\label{dva phirho}
^{\pi}\phi_{\varrho,!}^{(\eta)}:\Phi^+_{F_{\varrho}}(\CI(^{\pi}\Lambda)_!)\iso\
_{\varrho}C^{-r}_{^{\pi}{\ '\ff}}(V(^{\pi}\Lambda))
\end{equation}
by the formula
\begin{equation}
\label{dva phirhoform}
^{\pi}\phi_{\varrho,!}^{(\eta)}(b_{\varrho\leq\tau})=\sgn(\tau,\eta)
\sgn(\varrho)\theta_{\varrho\leq\tau}
\end{equation}
(see ~(\ref{dva monom})) where
\begin{equation}
\label{dva signrho}
\sgn(\varrho)=(-1)^{\sum_{i=1}^r(r-i+1)\cdot(\card(\varrho^{-1}(i))-1)}
\end{equation}
for $\varrho\in\CP_r(J;1)$. Taking the direct sum of $^{\pi}\phi_{\varrho,!}
^{(\eta)},\
\varrho\in\CP_r(J;1)$, we get isomorphisms
\begin{equation}
\label{dva phir}
^{\pi}\phi_{r,!}^{(\eta)}:\ ^+C^{-r}(^{\pi}\BA;\CI(^{\pi}\Lambda)_!)\iso\
_{\chi_J}C^{-r}_{^{\pi}{\ '\ff}}(V(^{\pi}\Lambda))
\end{equation}
A direct computation using ~\ref{dva matrixu} shows that the maps
$^{\pi}\phi_{r,!}^{(\eta)}$ are
compatible with differentials. Therefore, we arrive at

\subsection{}
\label{dva phi} {\bf Theorem.} {\em The maps $^{\pi}\phi_{r,!}^{(\eta)}$ induce
an isomorphism
of complexes
\begin{equation}
\label{dva phieq}
^{\pi}\phi_{!}^{(\eta)}:\ ^+C^{\bullet}(^{\pi}\BA;\CI(^{\pi}\Lambda)_!)\iso\
_{\chi_J}C^{\bullet}_{^{\pi}{\ '\ff}}(V(^{\pi}\Lambda))\ \Box
\end{equation}}

\vspace{1cm}
{\em COMPUTATIONS FOR $\CI(^{\pi}\Lambda)_*$}
\vspace{0.8cm}

\subsection{Bases} The Verdier dual to $\CI(^{\pi}\Lambda)_*$ is canonically
isomorphic to $\CI(^{\pi}\Lambda)^{-1}_!$. Therefore, by Theorem
{}~\ref{dva dual+} for
each positive facet $F$ we
have natural isomorphisms
\begin{equation}
\label{dva duality}
\Phi^+_F(\CI(^{\pi}\Lambda)_*)^*\cong \Phi^+_F(\CI(^{\pi}\Lambda)^{-1}_!)
\end{equation}
Let $\{ \tc_{\varrho\leq\tau}|\tau\in\Ord(\varrho)\}$ be
the basis in $\Phi^+_F(\CI(^{\pi}\Lambda)^{-1}_!)^*$ defined in ~\ref{dva basis!},
with $\CI(^{\pi}\Lambda)$ replaced by $\CI(^{\pi}\Lambda)^{-1}$.
We will denote by
$\{ c_{\varrho\leq\tau,*}|\tau\in\Ord(\varrho)\}$ the dual basis
in $\Phi^+_F(\CI(^{\pi}\Lambda)_*)^*$. Finally, we will denote by
$\{ b_{\varrho\leq\tau,*}|\tau\in\Ord(\varrho)\}$ the basis in
$\Phi^+_F(\CI(^{\pi}\Lambda)_*)$ dual to the previous one.

Our aim in the next subsections will be the description
of canonical morphisms $m:\Phi^+_F(\CI(^{\pi}\Lambda)_!)
\lra\Phi^+_F(\CI(^{\pi}\Lambda)_*)$ and
of the cochain complex $^+C^{\bullet}(^{\pi}\BA;\CI(^{\pi}\Lambda)_*)$ in
terms of our bases.

\subsection{Example.} Let us pick an element, $i\in I$ and let
$\pi:J:=\{ i\}\hra I$. Then we are in a one-dimensional situation,
cf. ~\ref{dva onedim}. The space $^{\pi}\BA$ has one coordinate $t_i$.
By definition, the local system $\CI(^{\pi}\Lambda)$ over
$^{\pi}\BAO=\ ^{\pi}\BA-\{ 0\}$ with the base point
$w: t_i=1$ has the fiber $\CI_w=\sk$ and monodromy equal to $\zeta^{-2\langle
\Lambda,i\rangle}$ along a counterclockwise loop.

The stratification $\CS_{\BR}$ has a unique $0$-dimensional
facet $F=F_{\varrho}$ corresponding to the unique element $\varrho
\in\CP_0(J;1)$,
and a unique $1$-dimensional positive facet $E=F_{\tau}$ corresponding
to the unique element $\tau\in\CP_1(J;1)$.

Let us construct the dual chain $c^*_1:=c_{\varrho\leq\tau,*}$. We adopt the
notations of ~\ref{dva onedim}, in particular of Fig.1.
The chain
$\tc_1:=\tc_{\varrho\leq\tau}\in
H^1(^{\pi}\BA,\{ 0,w\};j_!\CI(^{\pi}\Lambda)^{-1})^*$
is shown on Fig. 8(a) below. As a first step, we define the dual chain
$c^0_1\in H^1(^{\pi}\BA-\{ 0,w\},^{\pi}\BA_{\geq r''};\CI(^{\pi}\Lambda))^*$
--- it is also
shown on Fig. 8(a).

\begin{picture}(20,9)(-10,-4.5)

\put(-5,0){\circle{0.2}}
\put(-5.5,-0.4){$F$}
\put(-5,0){\oval(6,6)}
\put(-5.5,3.3){$S_{r''}$}

\put(-4.9,0){\line(1,0){4}}
\put(-1.5,0.3){$E$}

\put(-4.9,0){\vector(1,0){1}}
\put(-4,0){\circle*{0.1}}
\put(-4.7,-0.4){$\tc_1$}

\put(-4,-3){\line(0,1){6}}
\put(-4,0){\vector(0,1){1}}
\put(-3.8,2){$c^0_1$}

\put(-3,0){\circle*{0.2}}
\put(-3.1,-0.5){$w$}

\put(-5,0){\oval(3,3)}
\put(-6.5,1.6){$S_{r'}$}

\put(-5,-4){(a)}


\put(5,0){\circle{0.2}}
\put(5,0){\oval(6,6)}
\put(4.5,3.3){$S_{r'}$}

\put(5.1,0){\line(1,0){4}}

\put(5.1,0){\vector(1,0){1}}
\put(5.3,-0.4){$\epsilon\tc_1$}

\put(5,0){\oval(3,3)}
\put(5.5,1.5){\vector(-1,0){0.5}}
\put(5.5,1.7){$c^*_1$}

\put(6.5,0){\circle*{0.2}}
\put(6.6,-0.5){$\epsilon w$}

\put(5,-4){(b)}


\put(-0.5,-4.5){Fig. 8.}
\end{picture}

Next, we make a clockwise rotation of $c_1^0$ on $\pi/2$, and make a homotopy
inside the disk $^{\pi}\BA_{\leq r'}$ to a chain $c^*_1\in
H^1(^{\pi}\BA_{\leq r'}-\{ 0\},\epsilon\cdot w;\CI(^{\pi}\Lambda))^*$ beginning
and
ending at $\epsilon\cdot w$. This chain is shown on Fig. 8(b) (in a bigger
scale). Modulo "homothety" identification
$$
H^1(^{\pi}\BA_{\leq r'}-\{ 0\},\epsilon\cdot w;\CI(^{\pi}\Lambda))^*\cong
H^1(^{\pi}\BA-\{ 0\},w;\CI(^{\pi}\Lambda))^*=\Phi_F(\CI(^{\pi}\Lambda)_*)^*
$$
this chain coincides with $c_{\varrho\leq\tau,*}$.

The canonical map
$$
m^*: \Phi_F^+(\CI(^{\pi}\Lambda)_*)^*\lra\Phi_F^+(\CI(^{\pi}\Lambda)_!)^*
$$
carries $c_{\varrho\leq\tau,*}$ to
$[\langle\Lambda,i\rangle]_{\zeta}\cdot c_{\varrho\leq\tau}$. The
boundary map
$$
d^*: \Phi^+_F(\CI(^{\pi}\Lambda)_*)^*\lra\Phi^+_E(\CI(^{\pi}\Lambda)_*)^*
$$
carries $c_{\varrho\leq\tau,*}$ to
$[\langle\Lambda,i\rangle]_{\zeta}\cdot c_{\tau\leq\tau,*}$.

\subsection{Vanishing cycles at the origin}
\label{dva orig} Let us return
to the general situation. First let us consider an important case
of the unique $0$-dimensional facet --- the origin. Let $0$ denote the
unique element of $\CP_0(J;1)$. To shorten the notation,
let us denote $\Phi^+_{F_0}$ by $\Phi^+_0$.

The bases in $\Phi^+_0(\CI(^{\pi}\Lambda)_!)$,
etc. are
numbered by all bijections $\tau:J\iso [N]$.
Let us pick such a bijection. Let $\tau^{-1}(i)=\{j_i\},\ i=1,\ldots, N$.
The chain $c_{0\leq\tau}\in\Phi_{0}^+(\CI(^{\pi}\Lambda)_!)^*$
is depicted as follows:

\begin{picture}(20,6)(-10,-3)

\put(-4,0){\line(1,0){4}}

\put(-3,0){\circle{0.2}}
\put(-3.1,0.4){$0$}
\put(-3,0){\vector(1,0){0.8}}
\put(-2.2,0){\circle*{0.1}}
\put(-2.3,-0.4){$j_1$}
\put(-3,0){\vector(1,0){1.6}}
\put(-1.4,0){\circle*{0.1}}
\put(-1.5,-0.4){$j_2$}
\put(-3,0){\vector(1,0){2.4}}
\put(-0.6,0){\circle*{0.1}}
\put(-0.7,-0.4){$j_3$}

\put(0.5,0){$\ldots$}
\put(1.5,0){\line(1,0){2}}
\put(2,0){\vector(1,0){0.5}}
\put(2.5,0){\circle*{0.1}}
\put(2.4,-0.4){$j_N$}

\put(-1,-3){Fig. 9. $c_{0\leq\tau}$}

\end{picture}

Using considerations completely analogous to the one-dimensional case above,
we see that the dual chain
$c_{0\leq\tau,*}\in\Phi_{0}^+(\CI(^{\pi}\Lambda)_*)^*$ is portayed as
follows:

\begin{picture}(20,6)(-10,-3)

\put(-10,0){\line(1,0){10}}

\put(-3,0){\circle{0.2}}
\put(-3.1,0.2){$0$}

\put(-3,0){\oval(1.6,1.6)}
\put(-3,0.8){\vector(-1,0){0.3}}
\put(-2.2,0){\circle{0.2}}
\put(-2.1,-0.4){$j_1$}

\put(-3,0){\oval(3.2,2.4)}
\put(-3,1.2){\vector(-1,0){0.8}}
\put(-1.4,0){\circle{0.2}}
\put(-1.3,-0.4){$j_2$}

\put(-3,0){\oval(4.8,3.2)}
\put(-3,1.6){\vector(-1,0){1.3}}
\put(-0.6,0){\circle{0.2}}
\put(-0.5,-0.4){$j_3$}

\put(0.5,0){$\ldots$}
\put(1.5,0){\line(1,0){2}}

\put(-3,0){\oval(10,4)}
\put(-3,2){\vector(-1,0){1.8}}
\put(2,0){\circle{0.2}}
\put(2.1,-0.4){$j_N$}

\put(-1,-3){Fig. 10. $c_{0\leq\tau,*}$}

\end{picture}

The points $t_{j_i}$ are travelling independently along the corresponding
loops, in the indicated directions. The section of $\CI(^{\pi}\Lambda)^{-1}$
over
this cell is determined by the requierement to be equal to $1$ when
all points are equal to marking points, at the end of their the travel
(coming from below).

\subsection{Isomorphisms $\phi_{0,*}^{(\eta)}$} We will use notations of
{}~\ref{dva refine} and
{}~\ref{dva dualbases} with $I$ replaced by $J$, ${\ '\ff}$ ---
by $^{\pi}{\ '\ff}$,
with $n=1$ and
$\Lambda_0=\ ^{\pi}\Lambda$.
By definition, $C^0_{^{\pi}{\ '\ff}}(V(^{\pi}\Lambda))=V(^{\pi}\Lambda)$.
The space
$V(^{\pi}\Lambda)_{\chi_J}$
admits as a basis (of cardinality $N!$) the set consisting of all monomials
$$
\theta_{0\leq\tau}=\theta_{\tau^{-1}(N)}\theta_{\tau^{-1}(N-1)}\cdot
\ldots\cdot\theta_{\tau^{-1}(1)}v_{^{\pi}\Lambda},
$$
where $\tau$ ranges through the set of all bijections $J\iso [N]$.
By definition, $\{\theta^*_{0\leq\tau}\}$ is the dual basis of
$V(^{\pi}\Lambda)^*_{\chi_J}$.

Let us pick $\eta\in\CP_N(J)$.
Let us define an isomorphism
\begin{equation}
\label{dva piphio*}
^{\pi}\phi_{0,*}^{(\eta)}: \Phi^+_0(\CI(^{\pi}\Lambda)_*)\iso\
V(^{\pi}\Lambda)^*_{\chi_J}
\end{equation}
by the formula
\begin{equation}
\label{dva piphioform*}
^{\pi}\phi_{0,*}^{(\eta)}(b_{0\leq\tau})=\sgn(\tau,\eta)\theta^*_{0\leq\tau}
\end{equation}

\subsection{}
\label{dva shapov} {\bf Theorem.} {\em The square
$$\begin{array}{rcccl}
\;&\Phi^+_0(\CI(^{\pi}\Lambda)_!)&\overset{^{\pi}\phi_{0,!}^{(\eta)}}{\iso}&
         V(^{\pi}\Lambda)_{\chi_J}&\;\\
\;&m\downarrow&\;&\downarrow S_{^{\pi}\Lambda}&\;\\
\;&\Phi^+_0(\CI(^{\pi}\Lambda)_*)&\overset{^{\pi}\phi_{0,*}^{(\eta)}}{\iso}&
      V(^{\pi}\Lambda)^*_{\chi_J}&\;
\end{array}$$
commutes.}

{\bf Proof.} This follows directly from the discussion of ~\ref{dva orig}
and the definition of the form $S_{^{\pi}\Lambda}$. $\Box$

\subsection{}
\label{dva letus} Let us pass to the setup of ~\ref{dva color} and ~\ref{dva mainphi}.
Let $j_{\nu}:\CAO_{\nu}\hra\CA_{\nu}$ denote the embedding of the open
stratum of the principal stratification.

Set
\begin{equation}
\label{dva inulambda}
\CI_{\nu}(\Lambda)=(\pi_*\CI(^{\pi}\Lambda))^{\Sigma_{\pi},-}
\end{equation}
It is a local system over $\CAO_{\nu}$.

Let us define objects
\begin{equation}
\label{dva inula?}
\CI_{\nu}(\Lambda)_?:=j_{\nu?}\CI_{\nu}(\Lambda)[N]\in\CM(\CA_{\nu};\CS)
\end{equation}
where $?=!,*$ or $!*$. These objects will be called
{\em standard sheaves over $\CA_{\nu}$}.

The same reasoning as in ~\ref{dva jcolor} proves

\subsubsection{} {\bf Lemma.} {\em We have natural isomorphisms
\begin{equation}
\label{dva jstand}
\CI_{\nu}(\Lambda)_?\cong(\pi_*\CI(^{\pi}\Lambda)_?)^{\Sigma_{\pi},-}
\end{equation}
for $?=!,*$ or $!*$. $\Box$}

\subsection{} For a given $\eta\in\CP_N(J)$ the isomorphisms
$^{\pi}\phi_{0,*}^{(\eta)}$ and
$^{\pi}\phi_{0,!}^{(\eta)}$ are skew $\Sigma_{\pi}$-
equivariant. Therefore, after passing to invariants in Theorem ~\ref{dva shapov}
we get

\subsection{}
\label{dva shaposym} {\bf Theorem.} {\em The maps $^{\pi}\phi_{0,*}^{(\eta)}$  and
$^{\pi}\phi_{0,!}^{(\eta)}$ induce isomorphisms included into a commutative
square
$$\begin{array}{rcccl}
\;&\Phi_{\nu}(\CI_{\nu}(\Lambda)_!)&
\overset{\phi_{\nu,\Lambda,!}^{(\eta)}}{\iso}&
         V(\Lambda)_{\nu}&\;\\
\;&m\downarrow&\;&\downarrow S_{\Lambda}&\;\\
\;&\Phi_{\nu}(\CI_{\nu}(\Lambda)_*)&
\overset{\phi_{\nu,\Lambda,*}^{(\eta)}}{\iso}&
         V(\Lambda)^*_{\nu}&\;
\end{array}$$
and
\begin{equation}
\label{dva lnu}
\phi_{\nu,\Lambda,!*}^{(\eta)}:\Phi_{\nu}(\CI_{\nu}(\Lambda)_{!*})\iso
L(\Lambda)_{\nu}
\end{equation}}

{\bf Proof} follows from the previous theorem and Lemma ~\ref{dva S-sym}(ii).
$\Box$

\subsection{} Now suppose we are given an arbitrary $r$, $\varrho\in
\CP_r(J;1)$ and $\tau\in\Ord(\varrho)$. The picture of the dual chain
$c_{\varrho\leq\tau,*}$ is a combination of Figures 10 and 3. For example,
the chain dual to the one on Fig. 6 is portayed as follows:

\begin{picture}(20,6)(-10,-3)

\put(-8,0){\line(1,0){8}}
\put(0.5,0){$\ldots$}
\put(1.5,0){\line(1,0){6.5}}

\put(-6,0){\circle{0.2}}
\put(-6.1,0.3){$0$}
\put(-6,0){\oval(1.6,1.6)}
\put(-6,0.8){\vector(-1,0){0.5}}
\put(-5.2,0){\circle*{0.1}}
\put(-5.1,-0.4){$a$}

\put(-3,0){\circle{0.2}}
\put(-3.1,0.4){$1$}
\put(-3,0){\oval(1.6,1.6)[b]}
\put(-3,-0.8){\vector(1,0){0.5}}
\put(-2.2,0){\circle*{0.1}}
\put(-2.1,-0.4){$b$}

\put(-3,0){\oval(3.2,3.2)[b]}
\put(-3,-1.6){\vector(1,0){0.5}}
\put(-1.4,0){\circle*{0.1}}
\put(-1.2,-0.4){$c$}

\put(5,0){\circle{0.2}}
\put(4.9,0.4){$r$}
\put(5,0){\oval(1.6,1.6)[b]}
\put(5,-0.8){\vector(1,0){0.5}}
\put(5.8,0){\circle*{0.1}}
\put(5.9,-0.4){$d$}

\put(5,0){\oval(3.2,3.2)[b]}
\put(5,-1.6){\vector(1,0){0.5}}
\put(6.6,0){\circle*{0.1}}
\put(6.7,-0.4){$e$}

\put(-1,-3){Fig. 11. Chain $c_{\varrho\leq\tau,*}$.}

\end{picture}

\subsection{Isomorphisms $^{\pi}\Phi_*$} Let us pick $\eta\in\CP_N(J)$.
Let us define isomorphisms
\begin{equation}
\label{dva phirho*}
^{\pi}\phi_{\varrho,*}^{(\eta)}:\Phi^+_{F_{\varrho}}(\CI(^{\pi}\Lambda)_*)\iso\
_{\varrho}C^{-r}_{^{\pi}{\ '\ff}^*}(V(^{\pi}\Lambda)^*)
\end{equation}
by the formula
\begin{equation}
\label{dva phirhoform*}
^{\pi}\phi_{\varrho,*}^{(\eta)}(b_{\varrho\leq\tau,*})=\sgn(\tau,\eta)
\sgn(\varrho)\theta_{\varrho\leq\tau}^*
\end{equation}
where $\sgn(\varrho)$ is defined in ~(\ref{dva signrho}).

Taking the direct sum of $^{\pi}\phi_{\varrho,*}^{(\eta)},\
\varrho\in\CP_r(J;1)$, we get isomorphisms
\begin{equation}
\label{dva phir*}
^{\pi}\phi_{r,*}^{(\eta)}:\ ^+C^{-r}(^{\pi}\BA;\CI(^{\pi}\Lambda)_*)\iso\
_{\chi_J}C^{-r}_{^{\pi}{\ '\ff}^*}(V(^{\pi}\Lambda)^*)
\end{equation}

\subsection{Theorem}
\label{dva phi*} {\em The maps $^{\pi}\phi_{r,*}^{(\eta)}$ induce
an isomorphism
of complexes
\begin{equation}
\label{dva phieq*}
^{\pi}\phi_{*}^{(\eta)}:\ ^+C^{\bullet}(^{\pi}\BA;\CI(^{\pi}\Lambda)_*)\iso\
_{\chi_J}C^{\bullet}_{^{\pi}{\ '\ff}^*}(V(^{\pi}\Lambda)^*)
\end{equation}
which makes the square
$$\begin{array}{ccc}
^+C^{\bullet}(^{\pi}\BA;\CI(^{\pi}\Lambda)_!)&\overset{^{\pi}\phi_{!}^{(\eta)}}
\iso&\
_{\chi_J}C^{\bullet}_{^{\pi}{\ '\ff}}(V(^{\pi}\Lambda))\\
m\downarrow&\;&\downarrow S\\
^+C^{\bullet}(^{\pi}\BA;\CI(^{\pi}\Lambda)_*)&\overset{^{\pi}\phi_{*}^{(\eta)}}
\iso&\
_{\chi_J}C^{\bullet}_{^{\pi}{\ '\ff}^*}(V(^{\pi}\Lambda)^*)
\end{array}$$
commutative.}

{\bf Proof.} Compatibility with differentials is verified directly and
commutativity of the square are checked directly from the
geometric description of our chains (actually, it is sufficient to
check one of these claims --- the other one follows formally).

Note the geometric meaning of operators $t_i$ from ~(\ref{dva ti}) ---
they correspond to the deletion of the $i$-th loop on Fig. 10. $\Box$

\subsection{} Now let us pass to the situation ~\ref{dva letus}.
It follows from Theorem ~\ref{dva rgamma+} (after passing to skew $\Sigma_{\pi}$-
invariants) that the complexes $^+C^{\bullet}(^{\pi}\BA;\CI_{\nu?})
^{\Sigma_{\pi},-}$ compute the stalk of $\CI_{\nu?}$ at the origin.
Let us denote
this stalk by $\CI_{\nu?,0}$.

Therefore, passing to $\Sigma_{\pi}$-invariants in the previous theorem,
we get

\subsection{Theorem}
\label{dva stalks} {\em The isomorphisms $^{\pi}\phi^{(\eta)}_?$ where
$?=!,*$ or $!*$, induce isomorphisms in $\CCD^b(pt)$ included into a
commutative square
$$\begin{array}{ccc}
\CI_{\nu}(\Lambda)_{!,0}&\overset{_{\Lambda}\phi^{(\eta)}_{\nu,!,0}}
{\iso}&_{\nu}C^{\bullet}_{{\ '\ff}}(V(\Lambda))\\
m\downarrow&\;&\downarrow S\\
\CI_{\nu}(\Lambda)_{*,0}&\overset{_{\Lambda}\phi^{(\eta)}_{\nu,*,0}}
{\iso}&_{\nu}C^{\bullet}_{{\ '\ff}^*}(V(\Lambda)^*)
\end{array}$$
and
\begin{equation}
\label{dva stalkirr}
_{\Lambda}\phi^{(\eta)}_{\nu,!*,0}:
\CI_{\nu}(\Lambda)_{!*,0}\iso\ _{\nu}C^{\bullet}_{\ff}(L(\Lambda))\ \ \Box
\end{equation}}

\newpage
\begin{center}
{\bf Chapter 3. Fusion}
\end{center}
\vspace{1cm}

\section{Additivity theorem}

\subsection{} Let us start with the setup of ~\ref{dva j}. For a non-negative
integer $n$ let us denote by $(n)$ the set $[-n,0]$.
Let us introduce the following spaces.
$\BA^{(n)}$ - a complex affine space with a fixed system
of coordinates $(t_i),\ i\in (n)$.
Let $^nJ$ denote the disjoint union $(n)\cup J$;
$^n\BA$ --- a complex affine space with coordinates $t_j,\ j\in\ ^nJ$.

In general, for an affine space with a distinguished coordinate system of,
we will denote by $\CS_{\Delta}$ its diagonal stratification
as in ~\ref{dva diagsetup}.

Let $^n\BAO\subset\ ^n\BA$,
$\BAO\ ^{(n)}\subset\BA^{(n)}$ be the open strata of $\CS_{\Delta}$.

Let $^np:\ ^n\BA\lra\BA^{(n)}$ be the evident projection;
$^n\BB=\ ^np^{-1}(\BAO\ ^{(n)})$.
Given a point $\bz=(z_i)\in\BA^{(n)}$, let us denote by $^{\bz}\BA$ the
fiber $^np^{-1}(\bz)$ and by $^{\bz}\CS$ the stratification induced
by $\CS_{\Delta}$. We will consider $t_j,\ j\in J,$ as coordinates on
$^{\bz}\BA$.

The subscript $(.)_{\BR}$ will mean as usually "real points".

\subsection{}
Let us fix a point $\bz=(z_{-1},z_0)\in\BA^{(2)}_{\BR}$ such that
$z_{-1}<z_{0}$.
Let us concentrate on the fiber $^{\bz}\BA$. As an abstract variety
it is canonically isomorphic to $\BA$ --- a complex affine space
with coordinates $t_j,\ j\in J$; so we will suppress $\bz$ from its notation,
keeping it in the notation for the stratification $^{\bz}\CS$ where
the dependence on $\bz$ really takes place.

Let us fix a real $w>z_0$.
Let us pick two open non-intersecting disks $D_i\subset\BC$ with centra at
$z_i$ and not containing $w$. Let us pick two real numbers $w_{i}>z_i$
such that $w_i\in D_i$,
and paths $P_i$ connecting $w$ with $w_i$ as shown
on Fig. 12 below.

\begin{picture}(20,6)(-10,-3)

\put(-6,0){\circle{0.2}}
\put(-6.9,0){$z_{-1}$}


\put(-6,1){\oval(3,6)[b]}
\put(-6,0){\oval(2,2)}
\put(-6.4,-.6){$D_{-1}$}

\put(-.5,1){\oval(8,2)[t]}
\put(-.5,1){\oval(14,4)[t]}

\put(5,0){\circle*{0.2}}
\put(5.2,0){$w$}
\put(5,1){\oval(3,6)[b]}
\put(-6.3,1.5){$U_{-1}$}


\put(0.1,1){\oval(9.8,2.8)[t]}
\put(-5.5,1){\oval(1.4,2)[br]}
\put(5,0){\line(0,1){1}}
\put(-5.5,0){\circle*{0.2}}
\put(-5.9,.2){$w_{-1}$}
\put(.1,2.5){$P_{-1}$}


\put(0,0){\circle{0.2}}
\put(-.5,0.2){$z_{0}$}
\put(0,0){\oval(2,2)}
\put(-.2,-.6){$D_0$}


\put(.5,0){\circle*{.2}}
\put(0.2,.2){$w_0$}
\put(0.5,0){\line(1,0){4.5}}
\put(1.75,0.3){$P_0$}


\put(2.25,0){\oval(7.5,3)}
\put(1.75,-1){$U_0$}

\put(-1,-3){Fig. 12.}
\end{picture}

Let us denote by $\CQ_2(J)$ the set of all
maps $\rho:J\lra [-1,0]$. Given such a map, let us denote by $\BA_{\rho}
\subset\ \BA$ an open subvariety consisting of points
$(t_j)_{j\in J}$ such that $t_j\in D_{\rho(j)}$ for all $j$.
We will denote by the same letter $^{\bz}\CS$ the stratification of this
space induced by $^{\bz}\CS$.

Set $H_w:=\cup_j\ H_j(w)$; $P=P_{-1}\cup P_0$; $\tP=\{(t_j)\in\ \BA|
\mbox{ there exists $j$ such that }t_j\in P\}$.

Given $\CK\in\CCD(\BA;\ ^{\bz}\CS)$, the restriction map
$$
R\Gamma(\BA,H_w;\CK)\lra R\Gamma(\BA,\tP;\CK)
$$
is an isomorphism by homotopy. On the other hand, we have restriction maps
$$
R\Gamma(\BA,\tP;\CK)\lra R\Gamma(\BA_{\rho},
\tP_{\rho};\CK)
$$
where $\tP_{\rho}:=\tP\cap\BA_{\rho}$.
Therefore we have canonical maps
\begin{equation}
\label{dva resrho}
r_{\rho}:R\Gamma(\BA,H_w;\CK)\lra
R\Gamma(\BA_{\rho},\tP_{\rho};\CK)
\end{equation}

\subsection{Theorem}
\label{dva addtheor} {\em For every $\CK\in\CCD(\BA;\ ^{\bz}\CS)$ the
canonical map
\begin{equation}
\label{dva rest}
r=\sum r_{\rho}: R\Gamma(\BA,H_w;\CK)\lra\oplus_{\rho\in\CQ_2(J)}
R\Gamma(\BA_{\rho},\tP_{\rho};\CK)
\end{equation}
is an isomorphism.}

\subsection{Proof.} Let us pick two
open subsets $U_{-1},U_0\subset\BC$ as shown on Fig. 12.
Set $U=U_{-1}\cup U_0$, $\BA_U=\{(t_j)\in\BA|t_j\in U\mbox{ for all }
j\}$. It is clear that
the restriction morphism
$$
R\Gamma(\BA,\tP;\CK)\lra R\Gamma(\BA_U,\tP_U;\CK)
$$
where $\tP_U:=\tP\cap \BA_U$, is an isomorphism.

For each $\rho\in\CQ_2(J)$ set
\begin{equation}
\label{dva aurho}
\BA_{U,\rho}:=\{(t_j)\in \BA_U|t_j\in U_{\rho(j)}\mbox{ for all }j\};\
\tP_{U,\rho}:=\tP\cap\BA_{U,\rho}
\end{equation}
We have $\BA_U=\bigcup_{\rho}\BA_{U,\rho}$.

\subsubsection{} {\bf Lemma.} {\em For every $\CK\in\CCD(\BA;^{\bz}\CS)$
the sum of restriction maps
\begin{equation}
\label{dva mapq}
q:R\Gamma(\BA_U,\tP_U;\CK)\lra\sum_{\rho\in\CQ_2(J)}
R\Gamma(\BA_{U,\rho},\tP_{U,\rho};\CK)
\end{equation}
is an isomorphism.}

{\bf Proof.}
Suppose we have distinct $\rho_1,\ldots,\rho_m$ such that
$\BA_{U;\rho_1,\ldots,\rho_m}:=\BA_{U,\rho_1}\cap\ldots\cap \BA_{U,\rho_m}
\neq\emp$; set
$\tP_{U;\rho_1,\ldots,\rho_m}:=\tP\cap\BA_{U;\rho_1,\ldots,\rho_m}$.

Our lemma follows at once by Mayer-Vietoris argument from
the following

\subsubsection{} {\bf Claim.} {\em For every $m\geq 2$
$$
R\Gamma(\BA_{U;\rho_1,\ldots,\rho_m},
\tP_{U;\rho_1,\ldots,\rho_m};\CK)=0.
$$}

{\bf Proof of the claim.} It is convenient to use the following notations.
If $J=A\cup B$ is a disjoint union, we will denote by
$\rho_{A,B}$ the map $J\lra [-1,0]$ such that
$\rho^{-1}(-1)=A,\ \rho^{-1}(0)=B$, and
by $U_{A;B}$ the subspace $\BA_{U,\rho_{A,B}}$.

Let us prove the claim for the case $N=2$.
Let $J=\{i,j\}$.
In this case it is easy to see that the only non-trivial intersections are
$U^{(1)}=U_{ij;\emp}\cap U_{j;i}$ and $U^{(2)}=U_{ij;\emp}\cap U_{\emp;ij}$.

To prove our claim for $U^{(1)}$ we will use
{\em a shrinking neighbourhood argument}
based on Lemma I.2.12. Let $\CK'$ denote the sheaf on $U^{(1)}$ obtained
by extension by zero of $\CK|_{U^{(1)}-\tP}$. For each $\epsilon>0$
let us denote by $U_{-1,\epsilon}\subset\BC$ an open domain consisting
of points having distance $<\epsilon$ from $D_{-1}\cup P$. Set
$$
U^{(1)}_{\epsilon}=\{(t_i,t_j)|\ t_i\in U_{-1,\epsilon},\ t_j\in U_0\cap
U_{-1,\epsilon}\}.
$$
It is clear that restriction maps
$R\Gamma(U^{(1)};\CK')\lra R\Gamma(U^{(1)}_{\epsilon};\CK')$ are isomorphisms.
On the other hand, it follows from I.2.12 that
$$
\dirlim_{\epsilon}\ R\Gamma(U^{(1)}_{\epsilon};\CK')\cong
R\Gamma(\bigcap_{\epsilon}\ U^{(1)}_{\epsilon};\CK'),
$$
and the last complex is zero by the definition of $\CK'$
(the point $t_j$ is confined to $P$ in
$\bigcap_{\epsilon}\ U^{(1)}_{\epsilon}$).

The subspace $U^{(2)}$ consists of $(t_i,t_j)$ such that both $t_i$ and $t_j$
lie in $U_{-1}\cap U_{0}$. This case is even simpler. The picture is
homeomorphic to an affine plane with a sheaf smooth along
the diagonal stratification; and we are interested in its
cohomology modulo the coordinate cross. This is clearly equal to zero, i.e.
$R\Gamma(U^{(2)},\tP\cap U^{(2)};\CK)=0$.

This proves the claim
for $N=2$.
The case of an arbitrary $N$ is treated in a similar manner, and we leave
it to the reader. This completes the proof of the claim and of the lemma.
$\Box$

\subsubsection{} {\bf Lemma.} {\em For every $\rho\in\CQ_2(J)$ the restriction
map
$$
R\Gamma(\BA_{U,\rho},\tP_{U,\rho};\CK)\lra
R\Gamma(\BA_{\rho},\tP_{\rho};\CK)
$$
is an isomorphism.}

{\bf Proof.} Again let us consider the case $J=\{i,j\}$.
If $\rho=\rho_{ij;\emp}$ or $\rho_{\emp;ij}$ the statement is obvious.
Suppose $\rho=\rho_{i;j}$. Let us denote by $\BA'_{U,\rho}\subset\BA_{U,\rho}$
the subspace
$\{(t_i,t_j)|t_i\in D_{-1},\ t_j\in U_0\}$. It is clear that
the restriction map
$$
R\Gamma(\BA'_{U,\rho},\tP'_{U,\rho};\CK)\lra
R\Gamma(\BA_{\rho},\tP_{\rho};\CK)
$$
where $\tP'_{U,\rho}:=\tP\cap\BA'_{U,\rho}$,
is an isomorphism. Let us consider the restriction
\begin{equation}
\label{dva resaa}
R\Gamma(\BA_{U,\rho},\tP_{U,\rho};\CK)\lra
R\Gamma(\BA'_{U,\rho},\tP'_{U,\rho};\CK)
\end{equation}
The cone of this map is isomorphic to
$R\Gamma(\BA_{U,\rho},\tP_{U,\rho};\CM)$ where the sheaf $\CM$ has the same
singularities as $\CK$ and in addition is $0$ over the closure
$\bar{D}_{-1}$. Now, consider a system of shrinking neighbourhoods
of $P_{-1}\cup\bar{D}_{-1}$ as in the proof af the claim above,
we see that $R\Gamma(\BA_{U,\rho},\tP_{U,\rho};\CM)=0$, i.e.
{}~(\ref{dva resaa})
is an isomorphism. This implies the lemma for this case.

The case of arbitrary $J$ is treated exactly in the same way. $\Box$

Our theorem is an obvious consequence of two previous lemmas. $\Box$

\section{Fusion and tensor products}
\label{dva fus}

\subsection{Fusion functors} The constructions below were inspired by
{}~\cite{dr}.

For each integer $n\geq 1$ and $i\in [n]$, let
us define functors
\begin{equation}
\label{dva fusfunct}
^n\psi_i:\CCD(^n\BA;\CS_{\Delta})\lra\ \CCD(^{n-1}\BA;\CS_{\Delta})
\end{equation}
as follows. We have the $t$-exact nearby cycles functors
(see ~\cite{d3} or ~\cite{ks}, 8.6, but note the shift by [-1]!)
$$
\Psi_{t_{-i}-t_{-i+1}}[-1]:\CCD(^n\BA;\CS_{\Delta})\lra
\CCD(\BA';\CS_{\Delta})
$$
where $\BA'$ denotes (for a moment) an affine space with coordinates
$t_j,\ j\in ((n)\cup J)-\{-i\}$. We can identify the last space with
$^{n-1}\BA$ simply by renaming coordinates $t_j$ to $t_{j+1}$ for
$-n\leq j\leq -i-1$. By definition, $^n\psi_i$ is equal to
$\Psi_{t_{-i}-t_{-i+1}}[-1]$ followed by this identification.

\subsection{Lemma}
\label{dva assoc} {\em (i) For each $n\geq 2$, $i\in [n]$ have canonical
isomorphisms
\begin{equation}
\label{dva alphai}
^n\alpha_i:\ ^{n-1}\psi_i\circ\ ^n\psi_i\iso\ ^{n-1}\psi_i\circ\ ^n\psi_{i+1}
\end{equation}
and equalities
$$
^{n-1}\psi_{j}\circ\ ^n\psi_i=\ ^{n-1}\psi_i\circ\ ^n\psi_{j+1}
$$
for $j>i$,
such that

(ii) ({\em "Stasheff pentagon" identity})
the diagram below commutes:

\vspace{.8cm}

\begin{picture}(20,12)(-10,-6)

\put(0,6){$^{n-2}\psi_i\circ\ ^{n-1}\psi_i\circ\ ^n\psi_i$}

\put(5,2){$^{n-2}\psi_i\circ\ ^{n-1}\psi_i\circ\ ^n\psi_{i+1}$}
\put(5,-2){$^{n-2}\psi_i\circ\ ^{n-1}\psi_{i+1}\circ\ ^n\psi_{i+1}$}
\put(0,-6){$^{n-2}\psi_i\circ\ ^{n-1}\psi_{i+1}\circ\ ^n\psi_{i+2}$}

\put(-5,2){$^{n-2}\psi_i\circ\ ^{n-1}\psi_{i+1}\circ\ ^n\psi_i$}
\put(-5,-2){$^{n-2}\psi_i\circ\ ^{n-1}\psi_i\circ\ ^n\psi_{i+2}$}

\put(3.6,5.6){\vector(1,-1){2.6}}
\put(5,4.2){$^{n-2}\psi_i\circ\ ^n\alpha_i$}

\put(7,1){\vector(0,-1){2}}
\put(7.4,0){$^{n-1}\alpha_i\circ\ ^n\psi_{i+1}$}

\put(7,-2.3){\vector(-1,-1){3}}
\put(6,-4){$^{n-2}\psi_i\circ\ ^n\alpha_{i+1}$}

\put(0.4,5.6){\vector(-1,-1){2.4}}
\put(-3.8,4.6){$^{n-1}\alpha_i\circ\ ^n\psi_i$}

\put(-3.1,1){\line(0,-1){2}}
\put(-2.9,1){\line(0,-1){2}}

\put(-3,-2.3){\vector(1,-1){3}}
\put(-4.6,-4){$^{n-1}\alpha_i\circ\
^n\psi_{i+2}$}

\end{picture}
$\Box$}

\subsection{} Let us define a $t$-exact functor
\begin{equation}
\label{dva psi}
\psi:\CCD(^n\BA;\CS_{\Delta})\lra\CCD(\BA;\CS)
\end{equation}
as a composition $i_0^*[-1]\circ\ ^n\psi_1\circ\ ^{n-1}\psi_1\circ\ldots
\circ\ ^1\psi_1$, where
$$
i^*_0:\CCD(^0\BA;\CS_{\Delta})\lra\CCD(\BA;\CS)
$$
denotes the restriction to the subspace $t_0=0$. Note that $i^*_0[-1]$ is
a $t$-exact equivalence. (Recall that $\BA$ and $\CS$ denote the same as in
{}~\ref{dva j}).

\vspace{.5cm}
{\em STANDARD SHEAVES}
\vspace{.5cm}

The constructions and computations below generalize Section
{}~\ref{dva standsheaves}.

\subsection{} Let us make the following assumptions. Let us denote by $A$
the $I\times I$-matrix $(i\cdot j)$. Let us suppose that $\det\ A\neq 0$.
There exists a unique $\BZ[\frac{1}{\det\ A}]$-valued symmetric
bilinear form on $X$ (to be denoted by $\lambda,\mu\mapsto\lambda\cdot\mu$)
such that the map $\BZ[I]\lra X,\ \nu\mapsto\lambda_{\nu}$ respects
scalar products.

Let us suppose that our field $\sk$ contains an element $\zeta'$ such that
$(\zeta')^{\det\ A}=\zeta$, and fix such $\zeta'$. For $a=\frac{c}{\det\ A},\
c\in\BZ$, we set by definition $\zeta^a:=(\zeta')^c$.

\subsection{} Let us fix $\nu=\sum\nu_ii\in \BN[I]$ and its unfolding
$\pi:J\lra I$ as in ~\ref{dva color}, and an integer $n\geq 1$.
We will use the preceding notations
with this $J$.

Let us fix $n+1$ weights $\Lambda_0,
\Lambda_{-1},\ldots,\Lambda_{-n}\in X$.
We define a one-dimensional local system
$\CI(\Lambda_0,\ldots,\Lambda_{-n};\nu)$ over $^n\BAO$ exactly in the
same manner as in ~\ref{dva pici}, with half-monodromies equal to
$\zeta^{\pi(i)\cdot\pi(j)}$ if $i,j\in J$,
to $\zeta^{-\langle\Lambda_i,\pi(j)\rangle}$ if $i\in (n),j\in J$ and
to $\zeta^{\Lambda_i\cdot\Lambda_j}$ if $i,j\in (n)$.

Let
$j:\ ^n\BAO\lra\ ^n\BA$ be the embedding. Let us introduce the sheaves
\begin{equation}
\label{dva ilambdas}
\CI(\Lambda_0,\ldots,\Lambda_{-n};\nu)_?:=
j_?\CI(\Lambda_0,\ldots,\Lambda_{-n};\nu)[-n-N-1]\in\CM(^n\BA;\CS_{\Delta})
\end{equation}
where $?=!,*$ or $!*$.
Applying the functor $\psi$, we get the sheaves
$\psi\CI(\Lambda_0,\ldots,\Lambda_{-n};\nu)_?\in\CM(\BA;\CS)$.

All these objects are naturally $\Sigma_{\pi}$-equivariant. We define
the following sheaves on $\CA_{\nu}$:
\begin{equation}
\label{dva tensanu}
^{\psi}\CI_{\nu}(\Lambda_0,\ldots,\Lambda_{-n})_?:=
(\pi_*\psi\CI(\Lambda_0,\ldots,\Lambda_{-n};\nu)_?)^{\Sigma_{\pi},-}
\end{equation}
(cf. ~\ref{dva letus}).

The following theorem generalizes Theorem ~\ref{dva shaposym}.

\subsection{Theorem}
\label{dva shsymtens} {\em Given a bijection $\eta:J\iso [N]$, we have
natural isomorphisms included into a commutative square
$$\begin{array}{ccc}
\Phi_{\nu}(^{\psi}\CI_{\nu}(\Lambda_0,\ldots,\Lambda_{-n})_!)&
\overset{\phi_{!}^{(\eta)}}{\iso}&
         (V(\Lambda_0)\otimes\ldots\otimes V(\Lambda_{-n}))_{\nu}\\
m\downarrow&\;&\downarrow S_{\Lambda}\\
\Phi_{\nu}(^{\psi}\CI_{\nu}(\Lambda_0,\ldots,\Lambda_{-n})_*)&
\overset{\phi_{*}^{(\eta)}}{\iso}&
         (V(\Lambda_0)^*\otimes\ldots\otimes V(\Lambda_{-n})^*)_{\nu}
\end{array}$$
and
\begin{equation}
\label{dva lnucohtens}
\phi_{!*}^{(\eta)}:
\Phi_{\nu}(^{\psi}\CI_{\nu}(\Lambda_0,\ldots,\Lambda_{-n})_{!*})\iso
(L(\Lambda_0)\otimes\ldots\otimes L(\Lambda_{-n}))_{\nu}
\end{equation}
A change of $\eta$ multiplies these isomorphisms by the sign
of the corresponding permutation of $[N]$.}

{\bf Proof.} We may suppose that $\pi$ is injective, i.e. all
$\nu_i=0$ or $1$. The general case immediately follows from this one
after passing to $\Sigma_{\pi}$-skew invariants.

Let us consider the case $n=1$. In this case one sees easily
from the definitions that we have a canonical isomorphism
$$
\Phi_{\nu}(^{\psi}\CI_{\nu}(\Lambda_0,\Lambda_{-1})_?)\cong
R\Gamma(\BA,H_w;\CI(\Lambda_0,\Lambda_{-1};\nu)_?)
$$
in the notations of Additivity Theorem ~\ref{dva addtheor}. On the other hand,
the set $\CQ_2(J)$ is in one-to-one correspondence with the set of
all decompositions $\nu=\nu_0+\nu_{-1},\ \nu_i\in\BN[I]$, and if $\rho$
corresponds to such
a decomposition, we have a natural isomorphism
\begin{equation}
\label{dva factorstand}
R\Gamma(\BA_{\rho},\tP\cap\BA_{\rho};\CI_{\nu}(\Lambda_0,\Lambda_{-1})_?)\cong
\Phi_{\nu_0}(^{\psi}\CI(\Lambda_0;\nu_0)_?)\otimes
\Phi_{\nu_{-1}}(^{\psi}\CI(\Lambda_{-1};\nu_{-1})_?)
\end{equation}
by the K\"{u}nneth formula. Therefore,  Additivity Theorem implies
isomorphisms
$$
\Phi_{\nu}(^{\psi}\CI_{\nu}(\Lambda_0,\Lambda_{-1})_?)\cong
\oplus_{\nu_0+\nu_{-1}=\nu}
\Phi_{\nu_0}(^{\psi}\CI_{\nu_0}(\Lambda_0)_?)\otimes
\Phi_{\nu_{-1}}(^{\psi}\CI_{\nu_{-1}}(\Lambda_{-1})_?)
$$
which are the claim of our theorem.

The case $n>2$ is obtained similarly by the iterated use of the Additivity
Theorem. $\Box$

\subsection{} Next we will consider the stalks $(?)_0$ of our sheaves at $0$,
or what is the same (since they are $\BR^{+*}$-homogeneous), the complexes
$R\Gamma(\CA_{\nu};?)$.

The next theorem generalizes Theorem ~\ref{dva stalks}.

\subsection{Theorem}
\label{dva stalktens} {\em Given a bijection $\eta:J\iso [N]$, we have
natural isomorphisms included into a commutative square
$$\begin{array}{ccc}
^{\psi}\CI_{\nu}(\Lambda_0,\ldots,\Lambda_{-n})_{!0}&
\overset{\phi_{!,0}^{(\eta)}}{\iso}&\
_{\nu}C^{\bullet}_{{\ '\ff}}(V(\Lambda_0)
\otimes\ldots\otimes V(\Lambda_{-n}))\\
m\downarrow&\;&\downarrow S_{\Lambda}\\
^{\psi}\CI_{\nu}(\Lambda_0,\ldots,\Lambda_{-n})_{*0}&
\overset{\phi_{*,0}^{(\eta)}}{\iso}&\
_{\nu}C^{\bullet}_{{\ '\ff}^*}(V(\Lambda_0)^*
\otimes\ldots\otimes V(\Lambda_{-n})^*)
\end{array}$$
and
\begin{equation}
\label{dva lnustens}
\phi_{!*,0}^{(\eta)}:\
^{\psi}\CI_{\nu}(\Lambda_0,\ldots,\Lambda_{-n})_{!*0}\iso\
_{\nu}C^{\bullet}_{\ff}(L(\Lambda_0)\otimes\ldots\otimes L(\Lambda_{-n}))
\end{equation}
A change of $\eta$ multiplies these isomorphisms by the sign
of the corresponding permutation of $[N]$.}

{\bf Proof.} It is not hard to deduce from the previous theorem that
we have natural isomorphisms of complexes included into a
commutative square
$$\begin{array}{ccc}
^+C(\BA;\psi\CI(\Lambda_0,\ldots,\Lambda_{-n};\nu)_{!})&\iso&\
_{\chi_J}C^{\bullet}_{^{\pi}{\ '\ff}}(V(^{\pi}\Lambda_0)\otimes\ldots
\otimes V(^{\pi}\Lambda_{-n}))\\
m\downarrow&\;&\downarrow S_{\Lambda}\\
^+C(\BA;\psi\CI(\Lambda_0,\ldots,\Lambda_{-n};\nu)_{*})&\iso&\
_{\chi_J}C^{\bullet}_{^{\pi}{\ '\ff}^*}(V(^{\pi}\Lambda_0)^*\otimes\ldots
\otimes V(^{\pi}\Lambda_{-n})^*)
\end{array}$$
This implies our claim after passing to $\Sigma_{\pi}$-(skew) invariants.
$\Box$

\newpage
\begin{center}
{\bf Chapter 4. Category $\CC$}
\end{center}

\section{Simply laced case}
\label{dva compar}

\subsection{}
\label{dva zeta} From now on until the end of this part we will assume,
in addition to the assumptions
of ~\ref{dva notations}, that $\zeta$ is
a primitive $l$-th root of unity, where $l$ is a fixed integer $l>3$
prime to $2$ and $3$.

\subsection{}
\label{dva Lusztig}
We will use notations of ~\cite{l1}, Chapters 1, 2, which we briefly recall.

\subsubsection{}
Let $(I,\cdot)$ be {\em a simply laced Cartan datum of finite type}
(cf. {\em loc.cit.}, 1.1.1, 2.1.3),
that is, a finite set $I$ and a nondegenerate symmetric
bilinear form $\alpha,\beta\mapsto\alpha\cdot\beta$ on the free abelian
group $\BZ[I]$. This form satisfies conditions

(a) $i\cdot i=2$ for any $i\in I$;

(b) $i\cdot j\in\{0,-1\}$ for any $i\not=j$ in $I$.

\subsubsection{}
We will consider {\em the simply connected root datum of type $(I,\cdot)$},
that is (see {\em loc. cit.}, 2.2.2),
two free abelian groups $Y=\BZ[I]$ and $X=\Hom_{\BZ}(Y,\BZ)$ together with

(a) the canonical bilinear pairing
$\langle,\rangle:\;Y\times X\lra {\Bbb Z};$

(b) an obvious embedding $I\hookrightarrow Y\;(i\mapsto i)$ and an embedding
$I\hookrightarrow X\;(i\mapsto i')$, such that
$\langle i,j'\rangle=i\cdot j$ for any $i,j\in I$.

We will call $X$ {\em the lattice of weights}, and $Y$
{\em the lattice of coroots}.
An element of $X$ will be typically denoted by $\lambda,\mu,\nu,\ldots$;
and an element of $Y$ will be typically denoted by
$\alpha,\beta,\gamma,\ldots$.

\subsection{}
\label{dva C}
We consider the finite dimensional algebra $U$ over the field $\sk$
defined as in the section 1.3 of ~\cite{ajs}. We also consider the
category $\CC$ of finite dimensional $X$-graded $U$-modules
defined as in the section 2.3 of ~\cite{ajs}.

\subsubsection{}
\label{dva u}
The algebra $U$ is given by generators
$E_i,F_i,K^{\pm1}_i,\ i\in I$, subject to relations

(z) $K_i\cdot K_i^{-1}=1;\ K_iK_j=K_jK_i$;

(a) $K_jE_i=\zeta^{j\cdot i}E_iK_j$;

(b) $K_jF_i=\zeta^{-j\cdot i}F_iK_j$;

(c) $E_iF_j-F_jE_i=\delta_{ij}\frac{K_i-K_i^{-1}}{\zeta-\zeta^{-1}}$;

(d) $E_i^l=F_j^l=0$;

(e) $E_iE_j-E_jE_i=0$ if $i\cdot j=0$;
    $E_i^2E_j-(\zeta+\zeta^{-1})E_iE_jE_i+E_iE_j^2=0$ if $i\cdot j=-1$;

(f) $F_iF_j-F_jF_i=0$ if $i\cdot j=0$;
    $F_i^2F_j-(\zeta+\zeta^{-1})F_iF_jF_i+F_iF_j^2=0$ if $i\cdot j=-1$.

The algebra $U$ has a unique $\sk$-algebra $X$-grading $U=\oplus U_{\mu}$ for
which $|E_i|=i',\ |F_i|=-i',\ |K_i|=0$.

We define a comultiplication
\begin{equation}
\label{dva coprodu}
\Delta:U\lra U\otimes U
\end{equation}
as a unique $\sk$-algebra mapping such that
$$
\Delta(K^{\pm 1})=K^{\pm 1}\otimes K^{\pm 1};\
$$
$$
\Delta(E_i)=E_i\otimes K_i+1\otimes E_i;\
$$
$$
\Delta(F_i)=F_i\otimes 1+K_i^{-1}\otimes F_i.
$$
This makes $U$ a Hopf algebra (with obvious unit and counit).

\subsubsection{} The category $\CC$ is by definition a category
of finite dimensional $X$-graded $\sk$-vector spaces
$V=\oplus_{\mu\in X}V_{\mu}$,
equipped with a left action of $U$ such that the
$U$-action is compatible with the $X$-grading,
and
$$
K_ix=\zeta^{\langle i,\mu\rangle}x
$$
for $x\in V_{\mu},\ i\in I$.

Since $U$ is a Hopf algebra, $\CC$ has a canonical structure of a {\em tensor
category}.

\subsection{}
\label{dva fu}
We define an algebra $\fu$ having generators $\theta_i,\epsilon_i,
{K_i}^{\pm1},\ i\in I,$ subject to relations

(z) $K_i\cdot K_i^{-1}=1;\ K_iK_j=K_jK_i$;

(a) ${K_j}\epsilon_i=
{\zeta}^{j\cdot i}\epsilon_i{K_j}$;

(b) ${K_j}\theta_i=
{\zeta}^{-j\cdot i}\theta_i{K_j}$;

(c) $\epsilon_i\theta_j-\zeta^{i\cdot j}\theta_j\epsilon_i=
\delta_{ij}(1-K^{-2}_i)$

(d) if $f\in\Ker (S)\subset{\ '\ff}$ (see ~(\ref{dva formap})) then $f=0$;

(e) the same as (d) for the free algebra $\fE$ on the generators $\epsilon_i$.

\subsubsection{}
\label{dva comultu} Let us define the comultiplication
\begin{equation}
\label{dva comueq}
\Delta:\fu\lra\fu\otimes\fu
\end{equation}
by the formulas
$$
\Delta({K_i}^{\pm 1})={K_i}^{\pm 1}\otimes{K_i}^{\pm 1};\
$$
$$
\Delta(\theta_i)=\theta_i\otimes 1+K_i^{-1}\otimes\theta_i;\
$$
$$
\Delta(\epsilon_i)=\epsilon_i\otimes 1+K_i^{-1}\otimes\epsilon_i
$$
and the condition that $\Delta$ is a morphism of $\sk$-algebras.

This makes $\fu$ a Hopf algebra (with obvious unit and counit).

$\fu$ is an $X$-graded $\sk$-algebra, with an $X$-grading defined uniquely
by the conditions $|{K_i}^{\pm 1}|=0;\
|\theta_i|=-i';\ |\epsilon_i|=i'$.

\subsection{} We define $\tCC$ as a category of finite dimensional
$X$-graded vector spaces $V=\oplus V_{\lambda}$, equipped with a structure
of a left
$\fu$-module compatible with $X$-gradings and such that
$$
{K_i}x={\zeta}^{\langle i,\lambda\rangle}x
$$
for $x\in V_{\lambda},\ i\in I$.

Since $\fu$ is a Hopf algebra, $\tCC$ is a tensor category.

\subsection{}
Recall that for any $\Lambda\in X$ we have defined in ~\ref{dva llambda}
the $X$-graded $\ff$-module $L(\Lambda)$. It is a quotient-module of the
Verma module $V(\Lambda)$,
and it inherits its $X$-grading from the one of $V(\Lambda)$
(see ~\ref{dva verma}).
Thus $L(\Lambda)=\oplus L(\Lambda)_\lambda$, and we define the action of
generators ${K_i}$ on $L(\Lambda)_\lambda$ as multiplication by
${\zeta}^{\langle i,\lambda\rangle}$. Finally, we define the action of
generators $\epsilon_i$ on $V(\Lambda)$ as in ~\ref{dva slambda}. These operators
on $V(\Lambda)$ satisfy the relations (a) --- (c) above.

We
check immediately that this action descends to the quotient $L(\Lambda)$.
Moreover, it follows from Theorem ~\ref{dva coactshap}
that these operators acting on $L(\Lambda)$ satisfy the
relations (a) --- (e) above. So we have constructed the action of $\fu$
on $L(\Lambda)$, therefore we can regard it as an object of $\tCC$.

\subsubsection{} {\bf Lemma.} {\em $L(\Lambda)$ is an
irreducible object in $\tCC$.}

{\bf Proof.}
Let $I(\Lambda)$ be the maximal proper homogeneous (with respect to
$X$-grading)
submodule of $V(\Lambda)$ (the sum of all homogeneous submodules not containing
$v_\Lambda$). Then $V(\Lambda)/I(\Lambda)$ is irreducible, so it suffices to
prove that $I(\Lambda)=\ker(S_\Lambda)$. The inclusion $\ker(S_\Lambda)
\subset I(\Lambda)$ is obvious. Let us prove the opposite inclusion. Let
$y\in I(\Lambda)$. It is enough to check that $S_\Lambda(y,x)=0$ for any
$x\in V(\Lambda)$ of the form $\theta_{i_1}\ldots\theta_{i_n}v_\Lambda$.
By (22) we have $S_\Lambda(y,\theta_{i_1}\ldots\theta_{i_n}v_\Lambda)=
S_\Lambda(\epsilon_{i_n}\ldots\epsilon_{i_1}y,v_\Lambda)=0$ since
$\epsilon_{i_n}\ldots\epsilon_{i_1}y\in I(\Lambda)$. $\Box$

\subsection{}
Let us consider elements $E_i,F_i\in\fu$ given by the following formulas:

\begin{equation}
\label{dva formr}
E_i=\frac{\zeta^2}{\zeta-\zeta^{-1}}
\epsilon_i{K_i};\
F_i=\theta_i
\end{equation}
It is immediate to check that these elements satisfy the relations
{}~\ref{dva u} (a) --- (c).

Moreover, one checks without difficulty that
$$
\theta_i\theta_j-\theta_j\theta_i\in\Ker(S)\ \mbox{if}\ i\cdot j=0,
$$
and
$$
\theta_i^2\theta_j-(\zeta+\zeta^{-1})\theta_i\theta_j\theta_i+
\theta_j\theta_i^2\in\Ker(S)\ \mbox{if}\ i\cdot j=-1
$$
(cf. ~\cite{sv2}, 1.16).
Also, it is immediate that
$$
S(\theta_i^a,\theta_i^a)=\prod_{p=1}^a\frac{1-\zeta^{2p}}{1-\zeta^2}.
$$
It follows that $\theta_i^l\in\Ker(S)$ for all $i$.

It follows easily that the formulas ~(\ref{dva formr}) define a surjective
morphism of algebras
\begin{equation}
\label{dva mapr}
R:U\lra\fu
\end{equation}
Moreover, one checks at once that $R$ is a map of Hopf algebras.

Therefore, $R$ induces a tensor functor
\begin{equation}
\label{dva Q}
Q:\tCC\lra\CC
\end{equation}
which is an embedding of a full subcategory.

\subsection{Theorem}
\label{dva ceqcc}
{\em $Q$ is an equivalence.}

{\bf Proof.} It is enough to check that $\tilde{\CC}$ contains enough
projectives for $\CC$ (see e.g. Lemma A.15. of ~\cite{kl} IV).

First of all, we know from ~\cite{ajs}, section 4.1,
that the simple $\fu$-modules
$L(\Lambda),\Lambda\in X$ exhaust the list of simple objects of $\CC$.
Second, we know from ~\cite{apk}, Theorem 4.6 and Remark 4.7,
that the module $L(-\rho)$ is projective
where $-\rho\in X$ is characterized by the property $\langle i,-\rho\rangle=-1$
for any $i\in I$.
Finally we know, say, from ~\cite{apk}, Remark 4.7, and ~\cite{apw}, Lemma
9.11,
that the set of
modules $\{L(\Lambda)\otimes L(-\rho), \Lambda\in X\}$ is an ample system of
projectives for $\CC$.
$\Box$

\subsection{}
Denote by $\fu^0$ (resp., $U^0$) the subalgebra of $\fu$ (resp., of $U$)
generated by $K_i^{\pm 1},\ i\in I$. Obviously, both algebras
are isomorphic to the ring of Laurent polynomials in $K_i$, and the map
$R$ induces an identity isomorphism between them.

Denote by $\fu^{\leq 0}$ (resp., by $\fu^-$) the subalgebra of
$\fu$ generated by $\theta_i, {K_i}^{\pm 1}$
(resp., by $\theta_i$), $i\in I$. The last algebra may be
identified with $\ff$.
As a vector space, $\fu^{\leq 0}$
is isomorphic to $\ff\otimes\fu^0$.

Denote by $U^-\subset\ U$ the subalgebra generated by $F_i,i\in I$,
and by $U^{\leq 0}\subset\ U$ the subalgebra generated by $F_i,\
{K_i}^{\pm 1},i\in I$. As a vector space it is isomorphic to
$U^-\otimes U^0$.

\subsection{Theorem.} {\em (a) $R$ is an isomorphism;

(b) $R$ induces an isomorphism $U^-\iso\ff$.}

{\bf Proof.} Evidently it is enough to prove b). We know that $R$
is surjective, and
that $U^-$ is finite dimensional. So it suffices to prove that
$\dim\ U^-\leq\dim\ \ff$. We know by ~\cite{l1} 36.1.5. that
$\dim\ U^-=\dim L(-\rho)$. On the other hand, the map $\ff\lra L(-\rho),
f\mapsto f(v_{-\rho})$ is surjective by construction. $\Box$

\section{Non-simply laced case}

In the non-simply laced case all main results of this part hold true as well.
However, the definitions need some minor modifications.
In this section we will describe them.

\subsection{} We will use terminology and notations from ~\cite{l1},
especially from Chapters  1-3.
Let us fix a Cartan datum $(I,\cdot)$ of finite type, not necessarily
simply laced, cf. {\em loc. cit.}, 2.1.3.
Let $(Y=\BZ[I],X=\Hom (Y,\BZ),\langle,\rangle,I\hra X,I\hra Y)$
be the simply connected root datum associated with $(I,\cdot)$,
{\em loc.cit.}, 2.2.2.

We set $d_i:=\frac{i\cdot i}{2},\ i\in I$; these numbers are positive
integers. We set $\zeta_i:=\zeta^{d_i}$.

The embedding $I\hra X$ sends $i\in I$ to $i'$ such that
$\langle d_jj,i'\rangle =j\cdot i$ for all $i,j\in I$.

\subsection{} The category $\CC$ is defined in the same way as in the simply
laced case, where the definition of the Hopf algebra $U$ should be modified as
follows (cf. ~\cite{ajs}, 1.3).

By definition, $U$ has generators
$E_i,F_i,K^{\pm1}_i,\ i\in I$, subject to relations

(z) $K_i\cdot K_i^{-1}=1;\ K_iK_j=K_jK_i$;

(a) $K_jE_i=\zeta^{\langle j, i'\rangle}E_iK_j$;

(b) $K_jF_i=\zeta^{-\langle j, i'\rangle}F_iK_j$;

(c) $E_iF_j-F_jE_i=\delta_{ij}\frac{\tK_i-\tK_i^{-1}}{\zeta_i-\zeta_i^{-1}}$;

(d) $E_i^l=F_j^l=0$;

(e) $\sum_{p=0}^{1-\langle i,j'\rangle}(-1)^pE_i^{(p)}E_j
E_i^{(1-\langle i,j'\rangle-p)}=0$ for $i\neq j$;

(f) $\sum_{p=0}^{1-\langle i,j'\rangle}(-1)^pF_i^{(p)}F_j
F_i^{(1-\langle i,j'\rangle-p)}=0$ for $i\neq j$,

where we have used the notations:
$\tK_i:=K_i^{d_i}$;
$G_i^{(p)}:=G_i^p/[p]_i^!,\ G=E$ or $F$,
$$
[p]_i^!:=\prod_{a=1}^p\frac{\zeta_i^a-\zeta_i^{-a}}{\zeta_i-\zeta^{-1}_i}.
$$
The $X$-grading on $U$ is defined in the same way as in the simply laced case.

The comultiplication is defined as
$$
\Delta(K_i)=K_i\otimes K_i;
$$
$$
\Delta(E_i)=E_i\otimes\tK_i+1\otimes E_i;
$$
$$
\Delta(F_i)=F_i\otimes 1+\tK_i^{-1}\otimes F_i.
$$

\subsubsection{Remark} This algebra is very close to
(and presumably coincides with)
the algebra
$\bU$ from ~\cite{l1}, 3.1, specialized to $v=\zeta$. We use the opposite
comultiplication, though.

\subsection{} The definition of the Hopf algebra $\fu$ should be modified as
follows. It has generators
$\theta_i,\epsilon_i,
{K_i}^{\pm1},\ i\in I,$ subject to relations

(z) $K_i\cdot K_i^{-1}=1;\ K_iK_j=K_jK_i$;

(a) ${K_j}\epsilon_i=
{\zeta}^{\langle j, i'\rangle}\epsilon_i{K_j}$;

(b) ${K_j}\theta_i=
{\zeta}^{-\langle j, i'\rangle}\theta_i{K_j}$;

(c) $\epsilon_i\theta_j-\zeta^{i\cdot j}\theta_j\epsilon_i=
\delta_{ij}(1-\tK^{-2}_i)$

(d) if $f\in\Ker (S)\subset{\ '\ff}$ (see ~(\ref{dva formap})) then $f=0$;

(e) the same as (d) for the free algebra $\fE$ on the generators $\epsilon_i$.

The comultiplication is defined as
$$
\Delta({K_i}^{\pm 1})={K_i}^{\pm 1}\otimes{K_i}^{\pm 1};\
$$
$$
\Delta(\theta_i)=\theta_i\otimes 1+\tK_i^{-1}\otimes\theta_i;\
$$
$$
\Delta(\epsilon_i)=\epsilon_i\otimes 1+\tK_i^{-1}\otimes\epsilon_i
$$

The category $\tCC$ is defined in the same way as in the simply laced
case.

\subsection{} We define an $X$-grading on the free ${\ '\ff}$-module
$V(\Lambda)$ with generator $v_{\Lambda}$ by setting
$V(\Lambda)_{\Lambda}=\sk\cdot v_{\Lambda}$ and assuming that
operators $\theta_i$ decrease the grading by $i'$.

The definition of the form $S_{\Lambda}$ on $V(\Lambda)$ should be
modified as follows.
It is a unique bilinear from such that
$S_{\Lambda}(v_{\Lambda},v_{\Lambda})=1$ and
$S(\theta_ix,y)=S(x,\epsilon_iy)$ where the operators
$\epsilon_i:V(\Lambda)\lra V(\Lambda)$ are defined by the
requirements $\epsilon_i(v_{\Lambda})=0$,
$$
\epsilon_i(\theta_jx)=\zeta^{i\cdot j}\theta_j\epsilon_i(x)+
\delta_{ij}[\langle i,\lambda\rangle]_{\zeta_i}x
$$
for $x\in V(\Lambda)_{\lambda}$.

We define $L(\Lambda)$ as a quotient $V(\Lambda)/\Ker(S)$. As in the simply
laced case, $L(\Lambda)$ is naturally an object of $\CC$, and
the same argument proves that it is irreducible.

\subsection{} We define the morphism
\begin{equation}
\label{dva maprnsl}
R: U\lra\fu
\end{equation}
by the formulas
\begin{equation}
\label{dva formrnsl}
R(E_i)=\frac{\zeta_i^2}{\zeta_i-\zeta_i^{-1}}
\epsilon_i{\tK_i};\
R(F_i)=\theta_i;\ R(K_i)=K_i
\end{equation}
Using ~\cite{l1}, 1.4.3, one sees immediately that it is correctly
defined morphism of algebras. It follows at once from the definitions
that $R$ is a morphism of Hopf algebras.

Hence, we get a tensor functor
\begin{equation}
\label{dva Qnsl}
Q:\tCC\lra\CC
\end{equation}
and the same proof as in ~\ref{dva ceqcc} shows that $Q$ is an {\em equivalence
of categories}.

It is a result of primary importance for us.
It implies in particular that all irreducibles in $\CC$
(as well as their tensor products), come from $\tCC$.

\subsection{} Suppose we are given $\Lambda_0,\ldots,\Lambda_{-n}\in X$
and $\nu\in \BN[I]$. Let $\pi:J\lra I$, $\pi:\ ^{\pi}\BA\lra\CA_{\nu}$
denote the same as in ~\ref{dva color}. We will use the notations
for spaces and functors from Section ~\ref{dva fus}.

The definition of the local system $\CI(\Lambda_0,\ldots,\Lambda_{-n};
\nu)$ from {\em loc. cit.} should be modified: it should
have half-monodromies
$\zeta_j^{-\langle \pi(j),\Lambda_i\rangle}$ for
$i\in (n),\ j\in J$, the other formulas stay without change.

After that, the standard sheaves are defined as in {\em loc. cit.}
Now we are arriving at the main results of this part.
The proof is  the same as the proof of theorems ~\ref{dva shsymtens}
and ~\ref{dva stalktens}, taking into account the previous
algebraic remarks.

\subsection{Theorem}
\label{dva shsymnsl} {\em Let $L(\Lambda_0),\ldots,L(\Lambda_{-n})$ be
irreducibles of $\CC$, $\lambda\in X$,
$\lambda=\sum_{m=0}^n\Lambda_{-m}\ -\ \sum_i\nu_ii'$ for some
$\nu_i\in\BN$.  Set $\nu=\sum\nu_ii\in\BN[I]$.

Given a bijection $\eta:J\iso [N]$, we have
natural isomorphisms
\begin{equation}
\label{dva lnucohnsl}
\phi_{!*}^{(\eta)}:
\Phi_{\nu}(^{\psi}\CI_{\nu}(\Lambda_0,\ldots,\Lambda_{-n})_{!*})\iso
(L(\Lambda_0)\otimes\ldots\otimes L(\Lambda_{-n}))_{\lambda}
\end{equation}
A change of $\eta$ multiplies these isomorphisms by the sign
of the corresponding permutation of $[N]$. $\Box$}

\subsection{Theorem}
\label{dva stalknsl} {\em In the notations of the previous theorem we have
natural isomorphisms
\begin{equation}
\label{dva lnusnsl}
\phi_{!*,0}^{(\eta)}:\
^{\psi}\CI_{\nu}(\Lambda_0,\ldots,\Lambda_{-n})_{!*0}\iso\
C^{\bullet}_{\ff}(L(\Lambda_0)\otimes\ldots\otimes L(\Lambda_{-n}))_{\lambda}
\end{equation}
where we used the notation $C^{\bullet}(\ldots)_{\lambda}$ for
$_{\nu}C^{\bullet}(\ldots)$.
A change of $\eta$ multiplies these isomorphisms by the sign
of the corresponding permutation of $[N]$. $\Box$}

\newpage
\setcounter{section}{0}
\begin{center}{\large \bf Part III. TENSOR CATEGORIES}\end{center}

\begin{center}{\large \bf ARISING FROM CONFIGURATION SPACES}\end{center}

\section{Introduction}

\subsection{} 
In Chapter 1
we associate with every Cartan matrix of finite type and a non-zero
complex number $\zeta$ an
abelian artinian category $\FS$. We call
its objects {\em finite factorizable sheaves}.
They are certain infinite collections of perverse sheaves
on configuration spaces, subject to a compatibility
("factorization") and finiteness conditions.

In Chapter 2 the tensor structure on $\FS$ is defined
using functors of nearby cycles. It makes $\FS$ a braided tensor category.

In Chapter 3
we define, using vanishing cycles functors, an exact tensor functor
$$
\Phi:\FS\lra\CC
$$
to the category $\CC$ connected with the corresponding quantum group,
cf. ~\cite{ajs},
1.3 and ~\cite{fs} II, 11.3, 12.2.

In Chapter 4 we prove
\subsection{}
\label{tri thm} {\bf Theorem.} {\em $\Phi$ is an
equivalence of categories.}

One has to distinguish two cases.

(i) $\zeta$ is not root of unity. In this case it is wellknown that $\CC$
is semisimple. This case is easier to treat; ~\ref{tri thm} is
Theorem ~\ref{tri equiv thm gener}.

(ii) $\zeta$ is a root of unity. This is of course the most interesting
case; ~\ref{tri thm} is Theorem ~\ref{tri equiv thm}.

$\Phi$ may be regarded as a way of localizing
$\fu$-modules from category $\CC$ to the origin of the affine line ${\Bbb
A}^1$.
More generally, in order to construct tensor structure on $\FS$, we define
for each finite set $K$
certain categories $^K\FS$ along with the functor $g_K:\ \FS^K\lra\ ^K\FS$
(see the Sections ~\ref{tri sec kfs}, ~\ref{tri sec glu}) which
may be regarded as a way of localizing $K$-tuples of $\fu$-modules to
$K$-tuples of points of the affine line.
In the parts IV, V we will
show how to localize $\fu$-modules to the points of an arbitrary
smooth curve. For example, the case of a projective line
is already quite interesting, and is connected with
"semiinfinite" cohomology of quantum groups.

We must warn the reader that the proofs of some technical topological
facts are only sketched in this part.
The full details will appear later on.

\subsection{}
The construction of the space $\CA$ in Section 2 is inspired by the
idea of ``semiinfinite space of divisors on a curve'' one of us learnt from
A.Beilinson back in 1990. The construction of the braiding local system
$\CI$ in Section 3 is very close in spirit to P.Deligne's letter
{}~\cite{d1}. In terms of this letter, all the local systems $\CI^\alpha_\mu$
arise from the semisimple braided tensor category freely generated by the
irreducibles $i\in I$ with the square of $R$-matrix: $i\otimes j\lra
j\otimes i\lra i\otimes j$ given by the scalar matrix $\zeta^{i\cdot j}$.

We are very grateful to B.Feigin for the numerous inspiring discussions,
and to P.Deligne and L.Positselsky for the useful comments concerning the
definition of morphisms in $^K\FS$ in ~\ref{tri delpos}.

\subsection{Notations} We will use all the notations of parts I and II.
References to {\em loc.cit.} will look like Z.1.1 where Z=I or II.

During the whole part we fix a Cartan datum $(I,\cdot)$ of finite type
and denote by $(Y=\BZ[I];X=\Hom (Y,\BZ);I\hra Y,i\mapsto i;
I\hra X,i\mapsto i')$ the simply connected root datum associated with
$(I,\cdot)$, ~\cite{l1}, 2.2.2.
Given $\alpha=\sum a_ii\in Y$, we will denote by
$\alpha'$ the element $\sum a_ii'\in X$. This defines an embedding
\begin{equation}
\label{tri emb y x}
Y\hra X
\end{equation}

We will use the notation $d_i:=\frac{i\cdot i}{2}$.
We have $\langle j', d_ii\rangle=i\cdot j$. We will denote by $A$ the
$I\times I$-matrix $(\langle i,j'\rangle)$. We will denote by
$\lambda,\mu\mapsto\lambda\cdot\mu$ a unique $\BZ[\frac{1}{\det\ A}]$-
valued scalar product on $X$ such that ~(\ref{tri emb y x}) respects
scalar products. We have
\begin{equation}
\label{tri scal prod}
\lambda\cdot i'=\langle\lambda,d_ii\rangle
\end{equation}
for each $\lambda\in X,\ i\in I$.

We fix a non-zero complex number $\zeta'$ and suppose that our ground
field $\sk$ contains $\zeta'$. We set $\zeta:=(\zeta')^{\det\ A}$;
for $a=\frac{c}{\det\ A},\ c\in\BZ$, we will
use the notation $\zeta^a:=(\zeta')^c$.

We will use the following partial orders on $X$ and $Y$.
For $\alpha=\sum a_ii,\ \beta=\sum b_ii\in Y$ we write $\alpha\leq\beta$
if $a_i\leq b_i$ for all $i$. For $\lambda,\mu\in X$ we write
$\lambda\geq\mu$ if $\lambda-\mu=\alpha'$ for some
$\alpha\in Y,\ \alpha\geq 0$.

\subsection{} If $X_1,X_2$ are topological spaces, $\CK_i\in\CCD(X_i),
\ i=1,2$,
we will use the notation $\CK_1\boxtimes\CK_2$ for
$p_1^*\CK_1\otimes p_2^*\CK_2$ (where $p_i: X_1\times X_2\lra X_i$ are
projections). If $J$ is a finite set, $|J|$ will denote its cardinality.

For a constructible complex $\CK$, $SS(\CK)$ will denote
the singular support
of $\CK$ (micro-support in the terminology of ~\cite{ks}, cf. {\em loc.cit.},
ch. V).

\newpage
\begin{center}
{\bf Chapter 1. Category $\FS$}
\end{center}
\vspace{.8cm}

\section{Space $\CA$}

\subsection{} We will denote by $\BA^1$ the complex affine line with a fixed
coordinate $t$. Given real $c,c'$, $0\leq c<c'$, we will use the notations
$D(c)=\{t\in\BA^1|\ |t|<c\};\
\bD(c)=\{t\in\BA^1|\ |t|\leq c\};\
\BA^1_{>c}:=\{t\in\BA^1|\ |t|>c\}$;
$D(c,c')=\BA^1_{>c}\cap D(c')$.

Recall that we have introduced in II.6.12 configuration spaces
$\CA_{\alpha}$ for $\alpha\in\BN[I]$. If $\alpha=\sum a_ii$, the space
$\CA_{\alpha}$ parametrizes
configurations of $I$-colored points $\bt=(t_j)$ on $\BA^1$,
such that there are
precisely $a_i$ points of color $i$.

\subsection{} Let us introduce some open subspaces of $\CA_{\alpha}$.
Given a sequence
\begin{equation}
\label{tri seq pos roots}
\valpha=(\alpha_1,\ldots,\alpha_p)\in\BN[I]^p
\end{equation}
and a sequence of real numbers
\begin{equation}
\label{tri seq dist}
\vd=(d_1\ldots,d_{p-1})
\end{equation}
such that $0<d_{p-1}<d_{p-2}\ldots<d_{1},$ $p\geq 2$,
we define an open subspace
\begin{equation}
\label{tri ring conf}
\CA^{\valpha}(\vd)\subset\CA_{\alpha}
\end{equation}
which parametrizes configurations $\bt$ such that $\alpha_p$ of points $t_j$
lie inside the disk $D(d_{p-1})$, for $2\leq i\leq p-1$,
$\alpha_i$ of points lie inside the annulus $D(d_{i-1},d_i)$,
and $\alpha_1$ of points
lie inside the ring $\CA^1_{>d_1}$.

For $p=1$, we set $\CA^{\alpha}(\emp):=\CA_{\alpha}$.

By definition, a configuration space of empty collections of points
consists of one point. For example, so is $\CA^{0}(\emp)$.

\subsubsection{Cutting} Given $i\in [p-1]$ define subsequences
\begin{equation}
\label{tri subseq d}
\vd_{\leq i}=(d_1,\ldots, d_i);\ \vd_{\geq i}=(d_i,\ldots, d_{p-1})
\end{equation}
and
\begin{equation}
\label{tri subseq alpha}
\valpha_{\leq i}=(\alpha_1,\ldots,\alpha_i,0);\
\valpha_{\geq i}=(0,\alpha_{i+1},\ldots,\alpha_p)
\end{equation}
We have obvious {\em cutting isomorphisms}
\begin{equation}
\label{tri cut iso}
c_i:\CA^{\valpha}(\vd)\iso
\CA^{\valpha_{\leq i}}(\vd_{\leq i}) \times
\CA^{\valpha_{\geq i}}(\vd_{\geq i})
\end{equation}
satisfying the following compatibility:

for $i<j$ the square
\begin{center}
  \begin{picture}(14,6)
    \put(5,0){\makebox(4,2)
{$\CA^{\valpha_{\leq i}}(\vd_{\leq i})\times
\CA^{\valpha_{\geq i;\leq j}}(\vd_{\geq i;\leq j})
\times
\CA^{\valpha_{\geq j}}(\vd_{\geq j})$}}


    \put(5,4){\makebox(4,2){$\CA^{\valpha}(\vd)$}}


    \put(0,2){\makebox(4,2)
{$\CA^{\valpha_{\leq j}}(\vd_{\leq j})\times
\CA^{\valpha_{\geq j}}(\vd_{\geq j})$}}


    \put(10,2){\makebox(4,2)
{$\CA^{\valpha_{\leq i}}(\vd_{\leq i})\times
\CA^{\valpha_{\geq i}}(\vd_{\geq i})$}}


    \put(5.5,4.5){\vector(-2,-1){2}}
    \put(3.5,2.5){\vector(2,-1){2}}
    \put(10.5,2.5){\vector(-2,-1){2}}
    \put(8.5,4.5){\vector(2,-1){2}}

   \put(3.5,4){\makebox(1,0.5){$c_j$}}
   \put(3,1.5){\makebox(1,0.5){$c_i\times\id$}}
   \put(10.1,1.5){\makebox(1,0.5){$\id\times c_j$}}
   \put(9.7,4){\makebox(1,0.5){$c_i$}}

  \end{picture}
\end{center}
commutes.

\subsubsection{Dropping} For $i$ as above, let $\dpar_i\vd$ denote the
subsequence of $\vd$ obtained by dropping $d_i$, and set
\begin{equation}
\label{tri drop alpha}
\dpar_i\valpha=(\alpha_1,\ldots,\alpha_{i-1},
\alpha_{i}+\alpha_{i+1},\alpha_{i+2},\ldots,\alpha_p)
\end{equation}
We have obvious open embeddings
\begin{equation}
\label{tri drop emb}
\CA^{\dpar_i\valpha}(\dpar_i\vd)\hra
\CA^{\valpha}(\vd)
\end{equation}

\subsection{} Let us define $\CA_{\mu}^{\valpha}(\vd)$ as the space
$\CA^{\valpha}(\vd)$ equipped with an additional index $\mu\in X$.
One should understand $\mu$ as a weight assigned to the origin
in $\BA^1$. We will abbreviate the notation $\CA^{\alpha}_{\mu}(\emp)$ to
$\CA^{\alpha}_{\mu}$.

Given a triple $(\valpha,\vd,\mu)$ as above, let us define its $i$-cutting  ---
two triples
$(\valpha_{\leq i},\vd_{\leq i},\mu_{\leq i})$ and
$(\valpha_{\geq i},\vd_{\geq i},\mu)$, where
$$
\vmu_{\leq i}=\vmu-(\sum_{j=i+1}^p\alpha_j)'.
$$

We will also consider triples $(\dpar_i\valpha,\dpar_i\vd,\dpar_i\mu)$
where $\dpar_i\mu=\mu$ if $i<p-1$, and $\dpar_{p-1}\mu=\mu-\alpha_p'$.

The cutting isomorphisms ~(\ref{tri cut iso}) induce isomorphisms
\begin{equation}
\label{tri cut iso mu}
\CA^{\valpha}_{\mu}(\vd)\iso
\CA^{\valpha_{\leq i}}_{\mu_{\leq i}}(\vd_{\leq i}) \times
\CA^{\valpha_{\geq i}}_{\mu}(\vd_{\geq i})
\end{equation}

\subsection{}
\label{tri close emb} For each $\mu\in X,\ \alpha=\sum a_ii,
\beta=\sum b_ii\in\BN[I]$,
\begin{equation}
\label{tri sigma}
\sigma=\sigma^{\alpha,\beta}_{\mu}:\CA_{\mu}^{\alpha}\hra
\CA_{\mu+\beta'}^{\alpha+\beta}
\end{equation}
will denote a closed embedding which adds $b_i$ points of color $i$ equal
to $0$.

For $d>0$, $\CA^{(\alpha,\beta)}_{\mu+\beta'}(d)$ is an open
neighbourhood of $\sigma(\CA^{\alpha}_{\mu})$ in
$\CA^{\beta+\alpha}_{\mu+\beta'}$.

\subsection{} By definition, $\CA$ is a collection of all spaces
$\CA^{\valpha}_{\mu}(\vd)$ as above, together with the cutting isomorphisms
{}~(\ref{tri cut iso mu}) and the closed embeddings ~(\ref{tri sigma}).

\subsection{}
\label{tri conn comp} Given a coset $c\in X/Y$ (where we regard $Y$ as embedded
into $X$ by means of a map $i\mapsto i'$), we define
$\CA_c$ as a subset of $\CA$ consisting of $\CA_{\mu}^{\alpha}$ such that
$\mu\in c$. Note that the closed embeddings $\sigma$,
as well as cutting isomorphisms act inside $\CA_c$.
This subset will be called {\em a connected component} of $\CA$.
The set of connected components will be denoted $\pi_0(\CA)$.
Thus, we have canonically $\pi_0(\CA)\cong X/Y$.

\subsection{} We will be interested in two stratifications
of spaces $\CA_{\mu}^{\alpha}$. We will denote by
$\CAD^{\alpha}_{\mu}
\subset\CA^{\alpha}_{\mu}$ the complement
$$
A^{\alpha}_{\mu} - \bigcup_{\beta<\alpha}\ \sigma(\CA^{\beta}_
{\mu-\beta'+\alpha'}).
$$
We define a {\em toric stratification} of $\CA^{\alpha}_{\mu}$ as
$$
\CA_{\mu}^{\alpha}=\coprod\ \sigma(\CAD^{\beta}_
{\mu-\beta'+\alpha'}).
$$

Another stratification of $\CA_{\mu}^{\alpha}$ is {\em the principal
stratification} defined in II.7.14. Its open stratum will be denoted by
$\CAO^{\alpha}_{\mu}\subset\CA_{\mu}^{\alpha}$.
Unless specified otherwise, we will denote the prinicipal stratification
on spaces $\CA^{\alpha}_{\mu}$, as well as the induced stratifications
on its subspaces, by $\CS$.

The sign $\circ$ (resp., $\bullet$) over a subspace of $\CA^{\alpha}_{\mu}$
will denote the intersection of this subspace with $\CAO^{\alpha}_{\mu}$
(resp., with $\CAD^{\alpha}_{\mu}$).

\section{Braiding local system $\CI$}

\subsection{Local systems $\CI_{\mu}^{\alpha}$}
\label{tri def bls}
Let us recall some definitions from II. Let
$\alpha=\sum_ia_ii\in\BN[I]$
be given. Following II.6.12, let us choose an {\em unfolding} of $\alpha$,
i.e. a set $J$ together with a map $\pi: J\lra I$ such that $|\pi^{-1}(i)|=a_i$
for all $i$. We define the group
$\Sigma_{\pi}:=\{\sigma\in\Aut(J)|\ \sigma\circ\pi=\pi\}$.

We define $^{\pi}\BA$ as an affine space with coordinates $t_j,\ j\in J$;
it is equipped with the principal stratification defined by hyperplanes
$t_j=0$ and $t_i=t_j$, cf. II.7.1.
The group $\Sigma_{\pi}$ acts on $^{\pi}\BA$ by
permutations of coordinates, respecting the stratification. By definition,
$\CA_{\alpha}=\ ^{\pi}\BA/\Sigma_{\pi}$. We will denote by the same letter
$\pi$ the canonical projection $^{\pi}\BA\lra\CA_{\alpha}$.

If $^{\pi}\BAO\subset\ ^{\pi}\BA$ denotes the open stratum of the principal
stratification, $\pi(^{\pi}\BAO)=\CAO_{\pi}$, and the restriction of $\pi$
to $^{\pi}\BAO$ is unramified covering.

Suppose a weight $\mu\in X$ is given. Let us define a one dimensional
local system $^{\pi}\CI_{\mu}$ over $^{\pi}\BAO$ by the procedure II.8.1.
Its fiber over each positive chamber
$C\in\pi_0( ^{\pi}\BAO_{\BR})$ is identified with $\sk$; and monodromies
along standard paths shown on II, Fig. 5 (a), (b) are given by the formulas
\begin{equation}
\label{tri monodr}
^CT_{ij}=\zeta^{-\pi(i)\cdot \pi(j)},\ ^CT_{i0}=\zeta^{2\mu\cdot\pi(i)'}
\end{equation}
respectively (cf. ~(\ref{tri scal prod})). (Note that, by technical reasons,
this definition differs by the sign from that of II.8.2 and II.12.6).

We have a canonical $\Sigma_{\pi}$-equivariant structure on $^{\pi}\CI_{\mu}$,
i.e. a compatible system of isomorphisms
\begin{equation}
\label{tri equiv}
i_{\sigma}:\ ^{\pi}\CI_{\mu}\iso\sigma^*\ ^{\pi}\CI_{\mu},\
\sigma\in\Sigma_{\pi},
\end{equation}
defined uniquely by the condition that
$$
(i_{\sigma})_C=\id_\sk:\ (^{\pi}\CI_{\mu})_C=\sk\iso
(\sigma^*\ ^{\pi}\CI_{\mu})_{\sigma C}=\sk
$$
for all (or for some) chamber $C$. As a consequence, the group
$\Sigma_{\pi}$ acts on the local system $\pi_*\ ^{\pi}\CI_{\mu}$.

Let $\sgn:\Sigma_{\pi}\lra\{\pm 1\}$ denote the sign character. We define
a one-dimensional local system $\CI_{\mu}^{\alpha}$ over
$\CAO_{\mu}^{\alpha}=\CAO_{\alpha}$ as follows:
\begin{equation}
\label{tri def i}
\CI_{\mu}^{\alpha}:=(\pi_{*}\ ^{\pi}\CI_{\mu})^{\sgn}
\end{equation}
where the superscript $(\bullet)^{\sgn}$ denotes the subsheaf of sections $x$
such that $\sigma x=\sgn(\sigma)x$ for all $\sigma\in\Sigma_{\pi}$. Cf.
II.8.16.

Alternatively, we can define this local system as follows.
By the descent, there exists a unique
local system $\CI'_{\mu}$ over $\CA_{\alpha}$ such that $\pi^*\CI'_{\mu}$
is equal to
$^{\pi}\CI_{\mu}$ with the equivariant structure described above. In fact,
$$
\CI'_{\mu}=(\pi_*\ ^{\pi}\CI_{\mu})^{\Sigma_{\pi}}
$$
where the superscript $(\bullet)^{\Sigma_{\pi}}$ denotes invariants.
We have
\begin{equation}
\label{tri def2 i}
\CI_{\mu}^{\alpha}=\CI'_{\mu}\otimes\cSgn
\end{equation}
where $\cSgn$ denotes the one-dimensional local system over $\CA_{\alpha}$
associated with the sign representation $\pi_1(\CA_{\alpha})\lra\Sigma_{\pi}
\overset{\sgn}{\lra}\{\pm 1\}$.

This definition does not depend (up to a canonical isomorphism)
upon the choice of an unfolding.

\subsection{}
\label{tri def factoriz}
For each triple $(\valpha,\vd,\mu)$ as in the previous section,
let us denote by $\CI^{\valpha}_{\mu}(\vd)$ the restriction of
$\CI^{\alpha}_{\mu}$ to the subspace
$\CAO^{\valpha}_{\mu}(\vd)\subset\CAO^{\alpha}(\mu)$ where
$\alpha\in\BN[I]$ is the sum of components of $\valpha$.

Let us define {\em factorization isomorphisms}
\begin{equation}
\label{tri factor iso}
\phi_i=\phi^{\valpha}_{\mu;i}(\vd):
\CI^{\valpha}_{\mu}(\vd)
\iso
\CI^{\valpha_{\leq i}}_{\mu_{\leq i}}(\vd_{\leq i})
\boxtimes
\CI^{\valpha_{\geq i}}_{\mu}(\vd_{\geq i})
\end{equation}
(we are using identifications ~(\ref{tri cut iso mu})). By definition,
we have canonical
identifications of the stalks of all three local systems over a
point with real coordinates, with $\sk$. We define ~(\ref{tri factor iso}) as
a unique
isomorphism acting as identity when restricted to such a stalk. We will
omit irrelevant indices from the notation for $\phi$ if there is
no risk of confusion.

\subsection{Associativity}
\label{tri assoc} These isomorphisms have the following {\em associativity
property}.

For all $i<j$, diagrams
\begin{center}
  \begin{picture}(14,6)
    \put(5,0){\makebox(4,2)
{$\CI(\vd_{\leq i})\boxtimes
\CI(\vd_{\geq i;\leq j})
\boxtimes
\CI(\vd_{\geq j})$}}


    \put(5,4){\makebox(4,2){$\CI(\vd)$}}


    \put(0,2){\makebox(4,2)
{$\CI(\vd_{\leq j})\boxtimes \CI(\vd_{\geq j})$}}


    \put(10,2){\makebox(4,2)
{$\CI(\vd_{\leq i})\boxtimes \CI(\vd_{\geq i})$}}


    \put(5.5,4.5){\vector(-2,-1){2}}
    \put(3.5,2.5){\vector(2,-1){2}}
    \put(10.5,2.5){\vector(-2,-1){2}}
    \put(8.5,4.5){\vector(2,-1){2}}

   \put(3.5,4){\makebox(1,0.5){$\phi_j$}}
   \put(3,1.5){\makebox(1,0.5){$\phi_i\boxtimes\id$}}
   \put(10.1,1.5){\makebox(1,0.5){$\id\boxtimes\phi_j$}}
   \put(9.7,4){\makebox(1,0.5){$\phi_i$}}

  \end{picture}
\end{center}
are commutative.

In fact, it is enough to check the commutativity
restricted to some fiber $(\bullet)_C$, where it is obvious.

\subsection{}
\label{tri semiinf def} The collection of local systems
$\CI=\{\CI^{\alpha}_{\mu}\}$
together with factorization isomorphisms ~(\ref{tri factor iso}) will
be called {\em the braiding local system} (over $\CAO$).

The couple $(\CA,\CI)$ will be called {\em the semi-infinite configuration
space associated with the Cartan datum $(I,\cdot)$ and parameter $\zeta$}.

\subsection{} Let $j:\CAO^{\valpha}_{\mu}(\vd)\hra\CAD^{\valpha}_{\mu}(\vd)$
denote
an embedding; let us define a preverse sheaf
$$
\CID_{\mu}^{\valpha}(\vd):= j_{!*}\CI^{\valpha}_{\mu}(\vd)
[\dim\ \CA^{\valpha}_{\mu}]
\in\CM(\CAD^{\valpha}_{\mu}(\vd);\CS)
$$
By functoriality, factorization isomorphisms ~(\ref{tri factor iso})
induce analogous isomorphisms (denoted by the same letter)
\begin{equation}
\label{tri factor iso dot}
\phi_i=\phi^{\valpha}_{\mu;i}(\vd):
\CID^{\valpha}_{\mu}(\vd)
\iso
\CID^{\valpha_{\leq i}}_{\mu_{\leq i}}(\vd_{\leq i})
\boxtimes
\CID^{\valpha_{\geq i}}_{\mu}(\vd_{\geq i})
\end{equation}
satisfying the associativity property completely analogous to ~\ref{tri assoc}.

\section{Factorizable sheaves}

\subsection{} The aim of this section is to define certain $\sk$-linear
category $\tFS$. Its objects will be called {\em factorizable sheaves}
(over $(\CA,\CI)$). By definition, $\tFS$ is a direct product
of $\sk$-categories $\tFS_c$, where $c$ runs through $\pi_0(\CA)$
(see ~\ref{tri conn comp}). Objects of $\tFS_c$ will be called
{\em factorizable sheaves supported at $\CA_c$}.

In what follows we pick $c$, and denote by $X_c\subset X$ the corresponding
coset modulo $Y$.

\subsection{Definition}
\label{tri factor sh} {\em {\em A factorizable sheaf $\CX$ over
$(\CA,\CI)$ supported at $\CA_c$} is the following collection of data:

(a) a weight $\lambda\in X_c$;
it will be denoted by $\lambda(\CX)$;

(b) for each $\alpha\in\BN[I]$, a sheaf
$\CX^{\alpha}\in\CM(\CA^{\alpha}_{\lambda};\CS)$;

we will denote by $\CX^{\valpha}(\vd)$ perverse sheaves over
$\CA^{\valpha}_{\lambda}(\vd)$ obtained by taking the restrictions
with respect to the embeddings
$\CA^{\valpha}_{\lambda}(\vd)\hra \CA^{\alpha}_{\lambda}$;

(c) for each $\alpha,\beta\in \BN[I],\ d>0$, a {\em factorization isomorphism}
\begin{equation}
\label{tri factor iso sh}
\psi^{\alpha,\beta}(d):
\CX^{(\alpha,\beta)}(d)
\iso
\CID^{(\alpha,0)}_{\lambda-\beta'}(d)
\boxtimes
\CX^{(0,\beta)}(d)
\end{equation}

such that

{\em (associativity)}
for each $\alpha,\beta,\gamma\in\BN[I],\ 0<d_2<d_1$,
the square below must commute:

\begin{center}
  \begin{picture}(14,6)


    \put(5,4){\makebox(4,2){$\CX^{(\alpha,\beta,\gamma)}(d_1,d_2)$}}


    \put(0,2){\makebox(4,2)
{$\CID^{(\alpha,\beta,0)}_{\lambda-\gamma'}(d_1,d_2)
\boxtimes \CX^{(0,\gamma)}(d_2)$}}


    \put(10,2){\makebox(4,2)
{$\CID^{(\alpha,0)}_{\lambda-\beta'-\gamma'}(d_1)
\boxtimes \CX^{(0,\beta,\gamma)}(d_1,d_2)$}}

    \put(5,0){\makebox(4,2)
{$\CID^{(\alpha,0)}_{\lambda-\beta'-\gamma'}(d_1)\boxtimes
\CID^{(0,\beta,0)}_{\lambda-\gamma'}(d_1,d_2)
\boxtimes
\CX^{(0,\gamma)}(d_2)$}}


    \put(5.5,4.5){\vector(-2,-1){2}}
    \put(3.5,2.5){\vector(2,-1){2}}
    \put(10.5,2.5){\vector(-2,-1){2}}
    \put(8.5,4.5){\vector(2,-1){2}}

   \put(3.5,4){\makebox(1,0.5){$\psi$}}
   \put(3,1.5){\makebox(1,0.5){$\phi\boxtimes\id$}}
   \put(10.1,1.5){\makebox(1,0.5){$\id\boxtimes\psi$}}
   \put(9.7,4){\makebox(1,0.5){$\psi$}}

  \end{picture}
\end{center}}

\subsubsection{} Remark that with these definitions,
the braiding local system $\CI$ resembles a ``coalgebra'',
and a factorizable sheaf --- a ``comodule'' over it.

\subsection{Remark}
\label{tri fs gener} Note an immediate corollary of the factorization axiom.
We have isomorphisms
\begin{equation}
\label{tri res open}
\CX^{(\alpha,0)}(d)\cong\CX^0\otimes
\CID^{(\alpha,0)}_{\lambda}(d)
\end{equation}
(where $\CX^0$ is simply a vector space).

Our next aim is to define morphisms between factorizable sheaves.

\subsection{}
\label{tri trans maps}
Let $\CX$ be as above.
For each $\mu\geq\lambda$, $\mu=\lambda+\beta'$, and
$\alpha\in\BN[I]$,
let us define a sheaf $\CX^{\alpha}_{\mu}\in\CM(\CA^{\alpha}_{\mu};\CS)$ as
$\sigma_*\CX^{\alpha-\beta}$. For example,
$\CX^{\alpha}_{\lambda}=\CX^{\alpha}$. By taking restriction, the sheaves
$\CX^{\valpha}_{\mu}(\vd)\in\CM(\CA^{\valpha}_{\mu}(\vd);\CS)$ are defined.

Suppose $\CX,\CY$ are two factorizable sheaves supported at $\CA_c$,
$\lambda=\lambda(\CX),\ \nu=\lambda(\CY)$. Let $\mu\in X$,
$\mu\geq\lambda,\ \mu\geq\nu,\ \alpha,\beta\in\BN[I]$.
By definition we have canonical isomorphisms
\begin{equation}
\label{tri thetas}
\theta=\theta^{\beta,\alpha}_{\mu}:\
\Hom_{\CA^{\alpha}_{\mu}}(\CX^{\alpha}_{\mu},\CY^{\alpha}_{\mu})\iso
\Hom_{\CA^{\alpha+\beta}_{\mu+\beta'}}
(\CX^{\alpha+\beta}_{\mu+\beta'},\CY^{\alpha+\beta}_{\mu+\beta'})
\end{equation}

The maps ~(\ref{tri factor iso sh}) induce analogous isomorphisms
\begin{equation}
\label{tri factor iso sh mu}
\psi^{\alpha,\beta}_{\mu}(d):
\CX^{(\alpha,\beta)}_{\mu}(d)
\iso
\CID^{(\alpha,0)}_{\mu-\beta'}(d)
\boxtimes
\CX^{(0,\beta)}_{\mu}(d)
\end{equation}
which satisfy the same associativity property
as in ~\ref{tri factor sh}.

For $\alpha\geq\beta$ let us define maps
\begin{equation}
\label{tri taus}
\tau^{\alpha\beta}_{\mu}:
\Hom_{\CA^{\alpha}_{\mu}}(\CX^{\alpha}_{\mu},\CY^{\alpha}_{\mu})\lra
\Hom_{\CA^{\beta}_{\mu}}(\CX^{\beta}_{\mu},\CY^{\beta}_{\mu})
\end{equation}
as compositions
\begin{eqnarray}
\Hom_{\CA^{\alpha}_{\mu}}(\CX^{\alpha}_{\mu},\CY^{\alpha}_{\mu})
\overset{res}{\lra}
\Hom_{\CA^{(\alpha-\beta,\beta)}_{\mu}(d)}
(\CX^{(\alpha-\beta,\beta)}_{\mu}(d),\CY^{(\alpha-\beta,\beta)}_{\mu}(d))
\overset{\psi}{\iso}\\ \nonumber
\overset{\psi}{\iso}
\Hom_{\CA^{(\alpha-\beta,0)}_{\mu-\beta'}(d)}
(\CID^{(\alpha-\beta,0)}_{\mu-\beta'}(d),
 \CID^{(\alpha-\beta,0)}_{\mu-\beta'}(d))
\otimes
\Hom_{\CA^{(0,\beta)}_{\mu}(d)}(\CX^{(0,\beta)}_{\mu}(d),
\CY^{(0,\beta)}_{\mu}(d))
=\\ \nonumber
=\sk\otimes_\sk
\Hom_{\CA^{(0,\beta)}_{\mu}(d)}(\CX^{(0,\beta)}_{\mu}(d),
\CY^{(0,\beta)}_{\mu}(d))
=
\Hom_{\CA^{(0,\beta)}_{\mu}(d)}(\CX^{(0,\beta)}_{\mu}(d),
\CY^{(0,\beta)}_{\mu}(d))
\iso\\ \nonumber
\iso\Hom_{\CA^{\beta}_{\mu}}(\CX^{\beta}_{\mu},\CY^{\beta}_{\mu})\nonumber
\end{eqnarray}
where we have chosen some $d>0$,
the first map is the restriction, the second one is induced
by the factorization isomorphism, the last one is inverse to the
restriction. This definition does not depend on the choice
of $d$.

The associativity axiom implies that these maps satisfy an obvious
transitivity property. They are also compatible in the obvious way with the
isomorphisms $\theta$.

We define the space $\Hom_{\tFS_c}(\CX,\CY)$ as the following
inductive-projective limit
\begin{equation}
\label{tri ind proj}
\Hom_{\tFS_c}(\CX,\CY):=
\dirlim_{\mu}\invlim_{\beta}\Hom(\CX^{\beta}_{\mu},\CY^{\beta}_{\mu})
\end{equation}
where the inverse limit is over $\beta\in\BN[I]$, the transition maps
being $\tau^{\alpha\beta}_{\mu}$, $\mu$ being fixed, and the
direct limit over $\mu\in X$ such that
$\mu\geq\lambda,\ \mu\geq\nu$, the transition maps being
induced by ~(\ref{tri thetas}).

With these spaces of homomorphisms, factorizable sheaves supported at $\CA_c$
form a $\sk$-linear category to be denoted by $\tFS_c$ (the composition
of morphisms is obvious).

As we have already mentioned, the category
$\tFS$ is by definition a product $\prod_{c\in\pi_0(\CA)}\tFS_C$.
Thus, an object $\CX$ of $\tFS$ is a direct sum
$\oplus_{c\in\pi_0(\CA)}\CX_c$, where $\CX_c\in\tFS_c$.
If $\CX\in\tFS_c,\ \CY\in\tFS_{c'}$, then
$$
\Hom_{\tFS}(\CX,\CY)=\Hom_{\tFS_c}(\CX,\CY)
$$
if $c=c'$, and $0$ otherwise.

\subsection{}
\label{tri phi grad} Let $\Vect_f$ denote the category of finite dimensional
$\sk$-vector spaces.
Recall that in II.7.14 the functors
of "vanishing cycles at the origin"
$$
\Phi_{\alpha}:\CM(\CA^{\alpha}_{\mu};\CS)\lra\Vect_f
$$
have been defined.

Given $\CX\in\tFS_c$, let us define for each $\lambda\in X_c$ a vector space
\begin{equation}
\label{tri phi lambda}
\Phi_{\lambda}(\CX):=\Phi_{\alpha}(\CX^{\alpha}_{\lambda+\alpha'})
\end{equation}
where $\alpha\in\BN[I]$ is such that $\lambda+\alpha'\geq\lambda(\CX)$.
If $\lambda\in X-X_c$, we set $\Phi_{\lambda}(\CX)=0$.

Due to the definition of the sheaves $\CX^{\alpha}_{\mu}$,
{}~\ref{tri trans maps}, this vector space
does not depend on a choice of $\alpha$, up to a unique
isomorphism.

This way we get an exact functor
\begin{equation}
\label{tri phi x c}
\Phi:\tFS_c\lra\Vect_f^X
\end{equation}
to the category of $X$-graded vector spaces with finite dimensional components
which induces an exact functor
\begin{equation}
\label{tri phi x}
\Phi:\tFS\lra\Vect_f^X
\end{equation}

\subsection{Lemma}
\label{tri annih} {\em If $\Phi(\CX)=0$ then $\CX=0$.}

{\bf Proof.} We may suppose that $\CX\in\tFS_c$ for some $c$.
Let $\lambda=\lambda(\CX)$. Let us prove that for every
$\alpha=\sum a_ii\in\BN[I]$, $\CX^{\alpha}=0$. Let us do it by induction
on $|\alpha|:=\sum a_i$. We have $\CX^0=\Phi_{\lambda}(\CX)=0$ by assumption.

Given an arbitrary $\alpha$, it is easy to see from the factorizability and
induction hypothesis that $\CX^{\alpha}$ is supported at the origin
of $\CA_{\alpha}$. Since $\Phi_{\lambda-\alpha}(\CX)=0$, we conclude that
$\CX^{\alpha}=0$. $\Box$

\section{Finite sheaves}

\subsection{Definition} {\em A factorizable sheaf $\CX$ is called
{\em finite} if $\Phi(X)$ is finite dimensional.}

This is equivalent to saying that
there exists only finite number of $\alpha\in\BN[I]$
such that $\Phi_{\alpha}(\CX^{\alpha})\neq 0$ (or
$SS(\CX^{\alpha})$ contains the conormal
bundle to the origin $0\in\CA^{\alpha}_{\lambda}$,
where $\lambda:=\lambda(\CX)$).

\subsection{Definition} {\em The category of finite factorizable
sheaves (FFS for short) is a full subcategory $\FS\subset\tFS$
whose objects are finite factorizable sheaves.

We set $\FS_c:=\FS\cap\tFS_c$ for $c\in\pi_0(\CA)$.}

This category is our main character. It is clear that $\FS$ is a strictly
full subcategory of $\tFS$ closed with
respect to taking subobjects and quotients.

The next stabilization lemma is important.

\subsection{Lemma}
\label{tri stabilization} {\em Let $\CX,\CY$ be two FFS's supported at
the same connected component of $\CA$. For a fixed
$\mu\geq\lambda(\CX),\lambda(\CY)$ there exists $\alpha\in\BN[I]$ such that
for any $\beta\geq\alpha$ the transition map
$$
\tau^{\beta\alpha}_{\mu}:
\Hom_{\CA^{\beta}_{\mu}}(\CX^{\beta}_{\mu},\CY^{\beta}_{\mu})\lra
\Hom_{\CA^{\alpha}_{\mu}}(\CX^{\alpha}_{\mu},\CY^{\alpha}_{\mu})
$$
is an isomorphism.}

{\bf Proof.} Let us introduce a finite set
$$
N_{\mu}(\CY):=\{\alpha\in\BN[I]|\ \Phi_{\alpha}(\CY^{\alpha}_{\mu})\neq 0\}.
$$
Let us pick $\beta\in\BN[I]$. Consider a non-zero map
$f:\CX^{\beta}_{\mu}\lra\CY^{\beta}_{\mu}$. For each $\alpha\leq\beta$
we have a map $f^{\alpha}:=\tau_{\mu}^{\beta\alpha}(f):
\CX^{\alpha}_{\mu}\lra\CY^{\alpha}_{\mu}$. Let us consider subsheaves
$\CZ^{\alpha}:=\Ima(f^{\alpha})\subset\CY^{\alpha}_{\mu}$.
These subsheaves satisfy an obvious factorization property.

Let us consider the toric stratification of
$\CA^{\beta}_{\mu}$. For each $\alpha\leq\beta$ set
$\CA^{\alpha}:=\sigma(\CA^{\alpha}_{\mu+\alpha'-\beta'})
\subset\CA^{\beta}_{\mu}$;
$\CAD^{\alpha}:=\sigma(\CAD^{\alpha}_{\mu+\alpha'-\beta'})$.
Thus, the subspaces $\CAD^{\alpha}$ are strata of the toric stratification.
We have $\alpha_1\leq\alpha_2$ iff $\CA^{\alpha_1}\subset
\CA^{\alpha_2}$.

Let $\gamma$ denote a maximal element in the set
$\{\alpha|\ \CZ^{\beta}|_{\CAD^{\alpha}}\neq 0\}$. Then it is easy to see
that $\CZ^{\beta-\gamma}$ is a non-zero scyscraper on $\CA^{\beta-\gamma}$
supported at the origin. Therefore, $\Phi_{\beta-\gamma}
(\CZ^{\beta-\gamma})\neq 0$, whence
$\Phi_{\beta-\gamma}(\CY^{\beta-\gamma})\neq 0$, i.e.
$\beta-\gamma\in N_{\mu}(\CY)$.

Suppose that for some $\alpha\leq\beta$, $\tau_{\mu}^{\beta\alpha}(f)=0$.
Then
$$
\CZ^{\beta}|_{\CA^{(\beta-\alpha,\alpha)}(d)}=0
$$
for every $d>0$.
It follows that if $\CZ^{\beta}|_{\CAD^{\delta}}\neq 0$ then
$\delta\not\geq\beta-\alpha$.

Let us apply this remark to $\delta$ equal to $\gamma$ above.
Suppose that
$\gamma=\sum c_ii,\ \beta=\sum b_ii,\ \alpha=\sum a_ii$. There exists $i$
such that $c_i<b_i-a_i$. Recall that $\beta-\gamma\in N_{\mu}(\CY)$.
Consequently, we have

\subsubsection{} {\bf Corollary.} {\em Suppose that
$\alpha\geq \delta$ for all $\delta\in N_{\mu}(\CY)$. Then all the maps
$$
\tau^{\beta\alpha}_{\mu}:\Hom(\CX^{\beta}_{\mu},\CY^{\beta}_{\mu})\lra
\Hom(\CX^{\alpha}_{\mu},\CY^{\alpha}_{\mu}),
$$
$\beta\geq\alpha$, are injective.} $\Box$

Since all the spaces $\Hom(\CX^{\alpha}_{\mu},\CY^{\alpha}_{\mu})$ are
finite dimensional due to the constructibility of our sheaves,
there exists an $\alpha$ such that all $\tau^{\beta\alpha}_{\mu}$ are
isomorphisms. Lemma is proven. $\Box$

\subsection{}
\label{tri stabilization sec} For $\lambda\in X_c$ let us denote by
$\FS_{c;\leq\lambda}
\subset\FS_c$ the full subcategory whose objects are FFS's $\CX$
such that $\lambda(\CX_c)\leq\lambda$. Obviously
$\FS_c$ is a filtered union of these subcategories.

We have obvious functors
\begin{equation}
\label{tri proj bl}
p^{\beta}_{\lambda}:\FS_{c;\leq\lambda}\lra\CM(\CA^{\beta}_{\lambda};\CS),\
\CX\mapsto\CX^{\beta}_{\lambda}
\end{equation}

The previous lemma claims that for every $\CX,\CY\in\FS_{c;\leq\lambda}$
there exists $\alpha\in\BN[I]$ such that for every $\beta\geq\alpha$ the
map
$$
p^{\beta}_{\lambda}:
\Hom_{\FS}(\CX,\CY)\lra
\Hom_{\CA^{\beta}_{\lambda}}(\CX^{\beta}_{\lambda},\CY^{\beta}_{\lambda})
$$
is an isomorphism. (Obviously, a  similar claim holds true
for any finite number of FFS's.)

\subsection{Lemma}
\label{tri artin} {\em $\FS$ is an abelian artinian category.}

{\bf Proof.} $\FS$ is abelian by
Stabilization lemma. Each object has finite length by Lemma ~\ref{tri annih}.
$\Box$

\section{Standard sheaves}

\subsection{} For $\Lambda\in X$, let us define factorizable sheaves
$\CM(\Lambda), D\CM(\Lambda)_{\zeta^{-1}}$ and $\CL(\Lambda)$ as
follows. (The notation $D\CM(\Lambda)_{\zeta^{-1}}$ will be explained in
{}~\ref{tri dual fs} below).

Set
$$
\lambda(\CM(\Lambda))=\lambda(D\CM(\Lambda)_{\zeta^{-1}})=\lambda(\CL(\Lambda))=\Lambda.
$$
For $\alpha\in\BN[I]$ let $j$ denote the embedding
$\CAD_{\alpha}\hra\CA_{\alpha}$.
We define
$$
\CM(\Lambda)^{\alpha}=j_!\CID_{\Lambda}^{\alpha};\
D\CM(\Lambda)^{\alpha}_{\zeta^{-1}}=j_*\CID_{\Lambda}^{\alpha};\
\CL(\Lambda)^{\alpha}=j_{!*}\CID_{\Lambda}^{\alpha}.
$$
The factorization isomorphisms
are defined by functoriality from these isomorphisms for $\CID$.

Thus,
the collections $\{\CM(\Lambda)^{\alpha}\}_{\alpha}$, etc. form
factorizable sheaves to be denoted by
$\CM(\Lambda), D\CM(\Lambda)_{\zeta^{-1}}$ and $\CL(\Lambda)$ respectively.
Obviously, we have a canonical morphism
\begin{equation}
\label{tri mor!*}
m:\CM(\Lambda)\lra D\CM(\Lambda)_{\zeta^{-1}}
\end{equation}
and $\CL(\Lambda)$ is equal to its image.

\subsection{Theorem}
\label{tri all irreds} (i) {\em The factorizable sheaves
$\CL(\Lambda)$ are finite.}

(ii) {\em They are irreducible objects of $\FS$, non-isomorphic for
different $\Lambda$,
and they exhaust all irreducibles in $\FS$, up to isomorphism. }

{\bf Proof.} (i) follows from II.8.18.

(ii) Since the sheaves $\CL(\Lambda)^{\alpha}$ are
irreducible as objects of $\CM(\CA_{\alpha};\CS)$, the irreducibility
of $\CL(\Lambda)$ follows easily.
It is clear that they are non-isomorphic (consider the highest component).

Suppose $\CX$ is an irreducible FFS, $\lambda=\lambda(\CX)$. Let
$\alpha\in\BN[X]$ be a minimal among $\beta$ such that
$\Phi_{\lambda-\beta}(\CX)\neq 0$; set $\Lambda=\lambda-\alpha$.
By factorizability and the universal property of $!$-extension,
there exists a morphism if FS's $f:\CM(\Lambda)\lra\CX$ such that
$\Phi_{\Lambda}(f)\neq 0$ (hence is a monomorphism).
It follows from irreducibility of $\CL(\Lambda)$ that
the composition $\Ker(f)\lra\CM(\Lambda)\lra\CL(\Lambda)$ is equal
to zero, hence $f$ factors through a non-zero morphism
$\CL(\Lambda)\lra\CX$ which must be an isomorphism.
$\Box$

\subsection{}
\label{tri trivial} Let us look more attentively at the sheaf $\CL(0)$.

Let
$\tCA^{\alpha}\subset\CA^{\alpha}$ denote the open stratum
of the {\em diagonal} stratification, i.e. the complement to the diagonals.
Thus, $\CAO^{\alpha}\subset\tCA^{\alpha}$.
Let $\tCI^{\alpha}$ denote the local system over $\tCA^{\alpha}$ defined
in the same way as local systems $\CI^{\alpha}_{\mu}$, but using
only "diagonal" monodromies, cf. II.6.3.

One sees immediately that $\CL(0)^{\alpha}$ is equal to the middle
extenstion of $\tCI^{\alpha}$.

\newpage
\begin{center}
{\bf Chapter 2. Tensor structure}
\end{center}
\vspace{.8cm}

\section{Marked disk operad}

\subsection{} Let $K$ be a finite set. If $T$ is any set, we will denote
by $T^K$ the set of all mappings $K\lra T$; elements of $T^K$ will
be denoted typically by $\vx=(x_k)_{k\in K}$.

We will use the following partial orders on $X^K,\BN[I]^K$.
For $\vec{\lambda}=(\lambda_k),\ \vmu=
(\mu_k)\in X^K$, we write $\vlambda\geq\vmu$ iff
$\lambda_k\geq\mu_k$ for all $k$. An order on $\BN[I]^K$ is defined
in the same manner.

For $\valpha=(\alpha_k)\in \BN[I]^K$ we
will use the notation $\alpha$ for the sum of its components
$\sum_{k\in K}\alpha_k$; the same agreement will apply to $X^K$.

$\BA^K$ will denote the complex affine space with fixed
coordinates $u_k,\ k\in K$; $\BAO^K\subset\BA^K$ will denote
the open stratum of the diagonal stratification.

\subsection{Trees} We will call {\em a tree} a couple
\begin{equation}
\label{tri tree}
\tau=(\sigma,\vd)
\end{equation}
where $\sigma$, to be called {\em the shape} of $\tau$,
$$
\sigma=(K_p\overset{\rho_{p-1}}{\lra}
K_{p-1}\overset{\rho_{p-2}}{\lra}\ldots\overset{\rho_1}{\lra} K_1
\overset{\rho_0}{\lra} K_0)
$$
is a sequence of epimorhisms of finite sets, such that $\card(K_0)=1$,
$\vd=(d_0,d_1,\ldots,d_p)$, to be called
{\em the thickness} of $\tau$ --- a tuple of real numbers
such that $d_0=1>d_1>\ldots>d_p\geq 0$.

We will use a notation $\rho_{ab}$ for composition
$K_a\lra K_{a-1}\lra\ldots\lra K_b,\ a>b$.

A number $p\geq 0$ will be called
{\em the height} of $\tau$ and denoted $\hgt(\tau)$. Elements $k\in K_i$ will
be called {\em branches of height $i$}; $d_i$ will be called {\em
the thickness of $k$}. A unique branch of height $0$ will be called
{\em bole} and denoted by $*(\tau)$.

The set $K_p$ will be called {\em the base} of
$\tau$ and denoted $K_{\tau}$; we will also say that $\tau$ is
{\em $K_p$-based}; we will denote $d_p$ by $d_{\tau}$.
We will use notation $K(\tau)$ for the set
$\coprod_{i=0}^pK_i$ and $K'(\tau)$ for $\coprod_{i=0}^{p-1}K_i$.

A tree of height one will be called {\em elementary}. A tree $\tau$
whose branches of height $\hgt(\tau)$ have thickness $0$, will be called
{\em grown up}; otherwise it will be called {\em young}. We will assign to
every tree $\tau$ a grown up tree $\ttau$ by changing the thickness
of the thinnest branches to zero.

Thus, an elementary tree is essentially a finite set
and a real $0\leq d<1$; a grown up elementary tree is essentally
a finite set.

\subsubsection{Cutting}
\label{tri cutting} Suppose we have a tree $\tau$ as above,
and an integer $i$, $0<i<p$.
We define the operation of {\em cutting} $\tau$ at level $i$. It
produces from $\tau$ new trees $\tau_{\leq i}$ and
$\tau_{\geq k},\ k\in K_i$. Namely,
$$
\tau_{\leq i}=(\sigma_{\leq i},\vd_{\leq i}),
$$
where $\sigma_{\leq i}=(K_i\lra K_{i-1}\lra\ldots\lra K_0)$ and
$\vd_{\leq i}=(d_0,d_1,\ldots,d_i)$.

Second, for $k\in K_i$
$$
\tau_{\geq k}=(\sigma_{\geq k},\vd_{>i})
$$
where
$$
\sigma_{\geq k}=
(\rho^{-1}_{pi}(k)\lra\rho^{-1}_{p-1,i}(k)\lra\ldots\lra
\rho^{-1}_{i}(k)\lra\{k\}),
$$
$$
d_{>i}=(1,d_{i+1}d_i^{-1},\ldots,d_pd_i^{-1}).
$$

\subsubsection{}
\label{tri dpar} For $0<i\leq p$ we will denote by
$\dpar_i\tau$ a tree $(\dpar_i\sigma,\dpar_i\vd)$ where
$$
\dpar_i\sigma=(K_p\lra\ldots\lra\hat{K}_i\lra\ldots\lra K_0),
$$
and $\dpar_i\vd$ is obtained from $\vd$ by omitting $d_i$.

\subsection{Operad of disks}
\label{tri op disk}
For $r\in\BR_{\geq 0}\cup\{\infty\},\ z\in\BC$, we define an open disk
$D(z;r):=\{ u\in\BA^1|\ |u-z|<r\}$, and a closed disk
$\bD(z;r):=\{ u\in\BA^1|\ |u-z|\leq r\}$.

For a tree ~(\ref{tri tree})
we define a space
$$
\CO(\tau)=\CO(\sigma;\vd)
$$
parametrizing all collections $\vD=(\bD_k)_{k\in K_{\tau}}$ of closed
disks, such that $D_{*(\tau)}=D(0;1)$, for $k\in K_i$ the disk
$\bD_k$ has radius $d_i$,
for fixed $i\in [p]$ the disks $\bD_k,\ k\in K_i$, do not intersect,
and for each $i\in [0,p-1]$ and each $k\in K_{i+1}$ we have
$\bD_k\subset D_{\rho_{i}(k)}$.

Sometimes we will call such a collection
{\em a configuration of disks shaped by a tree $\tau$}.

\subsubsection{}
\label{tri int disks} Given such a configuration,
we will use the notation
\begin{equation}
\label{tri int disk}
\DO_k(\tau)=D_k-\bigcup_{l\in\rho_i^{-1}(k)}(\bD_l)
\end{equation}
if $k\in K_i$ and $i<p$, and we set
$\DO_k(\tau)=D_k$ if $i=p$.

If $\tau=(K\lra\{*\};d)$ is an elementary tree, we will use the notation
$\CO(K;d)$ for $\CO(\tau)$; if $d=0$, we will abrreviate the notation
to $\CO(K)$.

We have obvious embeddings
\begin{equation}
\label{tri emb elem}
\CO(K;d)\hra\CO(K)
\end{equation}
and
\begin{equation}
\label{tri emb k}
\CO(K)\hra\BAO^K
\end{equation}
this one is a homotopy equivalence.

We have open embeddings
\begin{equation}
\label{tri open}
\CO(\tau)\hra\CO(\ttau)
\end{equation}
obtained by changing the radius of smallest discs to zero.

\subsubsection{Substitution}
\label{tri substit} For each tree $\tau$ and $0<i<\hgt(\tau)$
we have the following {\em substitution isomorphisms}
\begin{equation}
\label{tri subst}
\CO(\tau)\cong
\CO(\tau_{\leq i})\times
\prod_{k\in K_i}\CO(\tau_{\geq k})
\end{equation}
In fact, a configuration
of disks shaped by a tree $\tau$ is the same as a configuration
shaped by $\tau_{\leq i}$, and for each $k\in K_i$ a configuration
shaped by $\tau_{>k}$ inside $D_k$ (playing the role of $D_0$; here
we have to make a dilation by $d_i^{-1}$).

These isomorphisms satisfy obvious quadratic relations connected
with pairs $0<i<j<\hgt(\tau)$. We leave their formulation
to the reader.

\subsection{Enhanced trees} We will call an {\em enhanced tree}
a couple $(\tau,\valpha)$ where $\tau$ is a tree and $\valpha\in
\BN[I]^{K'(\tau)}$. Vector $\valpha$ will be called {\em enhancement}
of $\tau$.

Let us define cutting for enhanced trees. Given $\tau$ and $i$ as in
{}~\ref{tri cutting}, let us note that $K'(\tau_{\leq i})$ and
$K'(\tau_{\geq k})$ are subsets of $K'(\tau)$. We define
$$
\valpha_{\leq i}\in \BN[I]^{K'(\tau_{\leq i})},\
\valpha_{\geq k}\in \BN[I]^{K'(\tau_{\geq k})}
$$
as the corresponding subsequences of $\valpha$.

Let us define operations $\dpar_i$ for enhanced trees. Namely, in the setup
of ~\ref{tri dpar}, we define
$\dpar_i\valpha=(\alpha'_k)\in\BN[I]^{K'(\dpar_i\tau)}$ as follows.
If $i=p$ then $K'(\dpar_p\tau)\subset K'(\tau)$, and we define
$\dpar_p\valpha$ as a corresponding subsequence. If $i<p$, we set
$\alpha'_k=\alpha_k$ if $k\in K_j,\ j>i$ or $j<i-1$. If $j=i-1$, we
set
$$
\alpha'_k=\alpha_k+\sum_{l\in\rho_{i-1}^{-1}(k)}\alpha_l.
$$

\subsection{Enhanced disk operad}
Given an enhanced tree $(\tau,\valpha)$,
let us define a configuration
space $\CA^{\valpha}(\tau)$ as follows. Its points are couples
$(\vD,\bt)$, where $\vD\in\CO(\tau)$ and
$\bt=(t_j)$ is an $\alpha$-colored configuration in $\BA^1$
(see II.6.12)
such that

{\em
for each $k\in K'(\tau)$ exactly $\alpha_k$ points lie inside $\DO_k(\tau)$
if $k\not\in K_{\tau}$ (resp., inside $D_k(\tau)$ if $k\in K_{\tau}$)
(see ~\ref{tri int disks}).}

In particular, all points lie inside $D_{*(\tau)}=D(0;1)$ and outside
$\bigcup_{k\in K_{\tau}}\bD_k$ if $\tau$ is young.
This space is an open subspace of the product
$\CO(\tau)\times\CA_{\alpha}$.

We will also use a notation
$$
\CA^{\alpha}(L;d):=\CA^{\alpha}(L\lra\{*\};d)
$$
for elementary trees and $\CA^{\alpha}(L)$ for $\CA^{\alpha}(L;0)$.

The isomorphisms ~(\ref{tri subst}) induce isomorphisms
\begin{equation}
\label{tri subst conf}
\CA^{\valpha}(\tau)\cong
\CA^{\valpha_{\leq i}}(\tau_{\leq i})\times
\prod_{k\in K_i}\CA^{\valpha_{\geq k}}(\tau_{\geq k})
\end{equation}

We have embeddings
\begin{equation}
\label{tri maps d}
d_i:\CA^{\valpha}(\tau)\lra
\CA^{\dpar_i\valpha}(\dpar_i\tau),\
0<i\leq p,
\end{equation}
--- dropping all disks $D_k,\ k\in K_i$.

We have obvious open embeddings
\begin{equation}
\label{tri tree base}
\CA^{\valpha}(\tau)\hra
\CA^{\alpha}(K_{\tau};d_{\tau})\hra
\CA^{\alpha}(K_{\tau})
\end{equation}

\subsection{Marked trees} We will call {\em a marked tree}
a triple $(\tau,\valpha,\vmu)$ where $(\tau,\valpha)$ is an
enhanced tree, and $\vmu\in X^{K_{\tau}}$. We will call
$\hgt(\tau)$ {\em the height} of this marked tree.

Let us define operations $\dpar_i$, $0<i\leq p=\hgt(\tau)$ for marked trees.
Namely, for $i<p$ we set $\dpar_i\vmu=\vmu$. For $i=p$ we define
$\dpar_p\vmu$ as $(\mu'_k)_{k\in K_{p-1}}$, where
$$
\mu_k'=\sum_{l\in\rho^{-1}_{p-1}(k)}\mu_l-\alpha'_l.
$$

Let us define cutting for marked trees. Namely, for
$1\leq i<p$ we define $\vmu_{\leq i}$ as
$\dpar_{i+1}\ldots\dpar_{p-1}\dpar_p\vmu$.

Next, for $k\in K_i$ we have $K_{\tau_{\geq k}}\subset K_{\tau}$, and
we define $\vmu_{\geq k}$ as a corresponding subsequence of $\vmu$.

\subsection{Marked disk operad} Now we can introduce our main objects.
For each
marked tree $(\tau,\valpha,\vmu)$ we define
$\CA_{\vmu}^{\valpha}(\tau)$ as a topological space
$\CA^{\valpha}(\tau)$ defined above, together with a marking
$\vmu$ of the tree $\tau$ considered as an additional index assigned to this
space.

We will regard
$\CA_{\vmu}^{\valpha}(\tau)$ as a space whose points are
configurations $(\vD,\bt)\in \CA^{\valpha}(\tau)$, together with
a marking of smallest disks $D_k,\ k\in K_{\tau}$, by weights $\mu_k$.

As above, we will use abbreviations $\CA_{\vmu}^{\alpha}(L;d)$
for $\CA_{\vmu}^{\alpha}(L\lra\{*\};d)$ (where $\vmu\in X^L$) and
$\CA_{\vmu}^{\alpha}(L)$ for $\CA_{\vmu}^{\alpha}(L;0)$.

We have natural open embeddings
\begin{equation}
\label{tri maps d mark}
d_i:\CA^{\valpha}_{\vmu}\lra
\CA^{\dpar_i\valpha}_{\dpar_i\vmu}(\dpar_i\tau),\
0<i\leq p,
\end{equation}
and
\begin{equation}
\label{tri tree base mark}
\CA^{\valpha}(\tau)\hra
\CA^{\alpha}(K_{\tau};d_{\tau})\hra
\CA^{\alpha}(K_{\tau})
\end{equation}
induced by the corresponding maps without marking.

The substitution isomorphisms ~(\ref{tri subst conf}) induce isomorphisms
\begin{equation}
\label{tri subst conf mark}
\CA^{\valpha}_{\vmu}(\tau)\cong
\CA^{\valpha_{\leq i}}_{\vmu_{\leq i}}(\tau_{\leq i})\times
\prod_{k\in K_i}\CA^{\valpha_{\geq k}}_{\vmu_{\geq k}}(\tau_{\geq k})
\end{equation}

\subsection{} We define closed embeddings
\begin{equation}
\label{tri sigma k}
\sigma=\sigma^{\alpha}_{\vmu;\vbeta}:\CA^{\alpha}_{\vmu}(K)\lra
\CA^{\alpha+\beta}_{\vmu+\vbeta'}(K)
\end{equation}
where $\vbeta=(\beta_k)_{k\in K},\ \beta_k=\sum_i\ b_k^i\cdot i$ and
$\beta=\sum_k\ \beta_k$. By definition,
$\sigma$ leaves points $u_k$ intact (changing their markings) and adds
$b_k^i$ copies of points of color $i$ equal to $u_k$.

\subsection{Stratifications} We set
$$
\CAD^{\alpha}_{\vmu}(K):=\CA^{\alpha}_{\vmu}(K)-
\bigcup_{\vgamma>0}\sigma(\CA^{\alpha-\gamma}_{\vmu-\vgamma}(K))
$$
We define a {\em toric stratification} of $\CA_{\vmu}^{\alpha}(K)$ as
$$
\CA_{\vmu}^{\alpha}(K)
=\coprod_{\vbeta<\valpha}\sigma(\CAD^{\beta}_{\vmu-\valpha'
+\vbeta'}(K))
$$

A {\em principal stratification} on $\CA^{\alpha}_{\vmu}(K)$ is defined
as follows.
The space $\CA^{\alpha}_{\vmu}(K)$ is a quotient of
$\BAO^K\times\BA^J$ where $\pi:J\lra I$ is an unfolding of $\alpha$
(cf. II.6.12). We define the principal stratification as the image of the
diagonal stratification on  $\BAO^K\times\BA^J$ under the canonical
projection $\BAO^K\times\BA^J\lra\CA^{\alpha}_{\vmu}(K)$.
We will denote by $\CAO^{\alpha}_{\vmu}(K)$
the open stratum of the principal stratification.

\section{Cohesive local systems $^K\CI$}

\subsection{}
\label{tri def i k} Let us fix a non-empty finite set $K$.
Suppose we are given $\vmu\in X^K$ and
$\alpha\in\BN[I]$. Let us pick an unfolding of $\alpha$,
$\pi:J\lra I$. Let
\begin{equation}
\label{tri pi alpha}
\pi^{\alpha}_{\vmu}:(D(0;1)^K\times D(0;1)^J)^{\circ}\lra
\CAO^{\alpha}_{\vmu}(K)
\end{equation}
denote the canonical projection (here $(D(0;1)^K\times D(0;1)^J)^{\circ}$
denotes
the open stratum of the diagonal stratification).

Let us define a one dimensional local system $^{\pi}\CI_{\vmu}$ by the
same procedure as in ~\ref{tri def bls}. Its fiber over each positive chamber
$C\in\pi_0((D(0;1)^K\times D(0;1)^J)^{\circ}_{\BR})$ is identified with $\sk$.
Monodromies
along the standard paths are given by the formulas
\begin{equation}
\label{tri monodr ki}
^CT_{ij}=\zeta^{-\pi(i)\cdot \pi(j)},\ ^CT_{ik}=\zeta^{2\mu_k\cdot\pi(i)'},\
^CT_{km}=\zeta^{-\mu_k\cdot\mu_m},
\end{equation}
$i,j\in J,\ i\neq j;\ k,m\in K,\ k\neq m$. Here $^CT_{ij}$ and
$^CT_{km}$ are half-circles, and $^CT_{ik}$ are full circles.
This definition essentially coincides with II.12.6, except for an overall
sign.

We define a one-dimensional local system $\CI_{\vmu}^{\alpha}(K)$
over $\CAO_{\vmu}^{\alpha}(K)$ as
\begin{equation}
\label{tri def ki}
\CI_{\vmu}^{\alpha}:=(\pi_{*}\ ^{\pi}\CI_{\vmu})^{\sgn}
\end{equation}
where the superscript $(\bullet)^{\sgn}$ has the same meaning as
in ~\ref{tri def bls}.

For each non-empty subset $L\subset K$ we can take a part
of weights $\vmu_L=(\mu_k)_{k\in L}$ and get a local system
$\CI_{\vmu}^{\alpha}(L)$ over $\CAO_{\vmu}^{\alpha}(L)$.

For each marked
tree $(\tau,\valpha,\vmu)$ with $K_{\tau}\subset K$, we define the
local system $\CI_{\vmu}^{\valpha}(\tau)$ as the restriction
of $\CI_{\vmu}^{\alpha}(K_{\tau})$ with respect to embedding
{}~(\ref{tri tree base mark}).

\subsection{Factorization}. The same construction as in
{}~\ref{tri def factoriz} defines
{\em factorization isomorphisms}
\begin{equation}
\label{tri factor i k}
\phi_i=\phi_{i;\vmu}^{\valpha}(\tau):
\CI_{\vmu}^{\valpha}(\tau)\cong
\CI_{\vmu_{\leq i}}^{\valpha_{\leq i}}(\tau_{\leq i})
\boxtimes
\fbox{$\times$}_{k\in K_i}\CI_{\vmu_{\geq k}}^{\valpha_{\geq k}}
(\tau_{\geq k})
\end{equation}

They satisfy the property of

\subsection{Associativity}
\label{tri assoc i} For all $0<i<j<p$ squares

\begin{center}
  \begin{picture}(14,6)

    \put(5,0){\makebox(4,2)
{$\CI(\tau_{\leq i})\boxtimes
\fbox{$\times$}_{k\in K_i}\ \CI(\tau_{\geq k;\leq j})
\boxtimes
\fbox{$\times$}_{l\in K_j}\
\CI(\tau_{\geq l})$}}


    \put(5,4){\makebox(4,2){$\CI(\tau)$}}


    \put(0,2){\makebox(4,2)
{$\CI(\tau_{\leq j})\boxtimes
\fbox{$\times$}_{l\in K_j}\ \CI(\tau_{\geq l})$}}


    \put(10,2){\makebox(4,2)
{$\CI(\tau_{\leq i})\boxtimes
\fbox{$\times$}_{k\in K_i}\ \CI(\tau_{\geq k})$}}


    \put(5.5,4.5){\vector(-2,-1){2}}
    \put(3.5,2.5){\vector(2,-1){2}}
    \put(10.5,2.5){\vector(-2,-1){2}}
    \put(8.5,4.5){\vector(2,-1){2}}

   \put(3.5,4){\makebox(1,0.5){$\phi_j$}}
   \put(3,1.5){\makebox(1,0.5){$\phi_i\boxtimes\id$}}
   \put(10.1,1.5){\makebox(1,0.5){$\id\boxtimes\phi_j$}}
   \put(9.7,4){\makebox(1,0.5){$\phi_i$}}

  \end{picture}
\end{center}

commute. (To unburden the notation we have omitted irrelevant indices ---
they are restored uniquely.)

\subsection{} The collection of local systems
$^K\CI=\{\CI^{\alpha}_{\vmu}(L),\ L\subset K\}$,
together with the factorization isomorphisms defined above,
will be called {\em the cohesive local system over $^K\CAO$}.

\subsection{} Let us define perverse sheaves
$$
\CID^{\valpha}_{\vmu}(\tau):=
j_{!*}\CI^{\valpha}_{\vmu}(\tau)
[\dim\ \CA^{\valpha}_{\vmu}(\tau)]\in
\CM(\CAD^{\valpha}_{\vmu}(\tau);\CS)
$$
where $j:\CAO^{\valpha}_{\vmu}(\tau)\hra
\CAD^{\valpha}_{\vmu}(\tau)$ denotes the embedding.
By functoriality, the factorization isomorphisms
{}~(\ref{tri factor i k}) induce isomorphisms
\begin{equation}
\label{tri factor idot k}
\phi_i=\phi_{i;\vmu}^{\valpha}(\tau):
\CID_{\vmu}^{\valpha}(\tau)\cong
\CID_{\vmu_{\leq i}}^{\valpha_{\leq i}}(\tau_{\leq i})
\boxtimes
\fbox{$\times$}_{k\in K_i}\CID_{\vmu_{\geq k}}
^{\valpha_{\geq k}}(\tau_{\geq k})
\end{equation}

These isomorphisms satisfy an associativity property completely
analogous to ~\ref{tri assoc i}; one should only replace
$\CI$ by $\CID$ in the diagrams.

\section{Factorizable sheaves over $^K\CA$}
\label{tri sec kfs}

We keep the assumptions of the previous section.

\subsection{} The first goal of this section is to define a $\sk$-linear
category $^K\tFS$ whose objects will be called {\em factorizable sheaves
(over $(^K\CA,\ ^K\CI)$)}. Similarly to $\tFS$, this category is by definition
a product of $\sk$-categories
\begin{equation}
\label{tri prod comp}
^K\tFS=\prod_{\vc\in\pi_0(\CA)^K}\ ^K\tFS_{\vc}.
\end{equation}
Objects of $^K\tFS_{\vc}$ will be called {\em factorizable sheaves
supported at $\vc$}.

\subsection{Definition}
\label{tri def kfactor sh} {\em A {\em factorizable sheaf $\CX$ over
$(^K\CA,\ ^K\CI)$ supported at $\vc=(c_k)\in\pi_0(\CA)^K$} is the
following collection of data:

(a) a $K$-tuple of weights $\vlambda=(\lambda_k)\in X^K$ such that
$\lambda_k\in X_{c_k}$, to be denoted
by $\vlambda(\CX)$;

(b) for each $\alpha\in\BN[I]$ a sheaf
$\CX^{\alpha}(K)\in\CM(\CA^{\alpha}_{\vlambda}(K);\CS)$.

Taking restrictions, as in ~\ref{tri def i k}, we get for each
$K$-based enhanced tree $(\tau,\valpha)$ sheaves
$\CX^{\valpha}(\tau)\in
\CM(\CA^{\valpha}_{\vlambda}(\tau);\CS)$.

(c) For each enhanced tree $(\tau,\valpha)$ of height $2$,
$\tau=(K\overset{\id}{\lra}K\lra \{*\};\ (1,d,0))$,
$\valpha=(\alpha,\vbeta)$ where
$\alpha\in\BN[I];\ \vbeta\in\BN[I]^K$,
a {\em factorization isomorphism}
\begin{equation}
\label{tri factor x k}
\psi(\tau):
\CX^{(\alpha,\vbeta)}(\tau)\cong
\CID_{\vlambda(\tau)_{\leq 1}}^{\alpha}(\tau_{\leq 1})
\boxtimes\CX^{(0,\vbeta)}(\tau)
\end{equation}

These isomorphisms should satisfy

{\em Associativity axiom.}

For all enhanced trees $(\tau,\valpha)$ of height $3$,
$\tau=(K\overset{\id}{\lra}K\overset{\id}{\lra}K\lra \{*\};
\ (1,d_1,d_2,0))$,
$\valpha=(\alpha,\vbeta,\vgamma)$ where $\alpha\in\BN[I];\ \vbeta,\vgamma
\in\BN[I]^K$,
the square

\begin{center}
  \begin{picture}(14,6)


    \put(5,4){\makebox(4,2){$\CX^{(\alpha,\vbeta,\vgamma)}(\tau)$}}


    \put(0,2){\makebox(4,2)
{$\CID^{(\alpha,\vbeta)}_{\vlambda(\tau)_{\leq 2}}(\tau_{\leq 2})\boxtimes
\CX^{(0,\vo,\vgamma)}(\tau)$}}


    \put(10,2){\makebox(4,2)
{$\CID^{\alpha}_{\vlambda(\tau)_{\leq 1}}(\tau_{\leq 1})\boxtimes
\CX^{(0,\vbeta,\vgamma)}(\tau)$}}

    \put(5,0){\makebox(4,2)
{$\CID^{\alpha}_{\vlambda(\tau)_{\leq 1}}(\tau_{\leq 1})
\boxtimes
\fbox{$\times$}_{k\in K}\
\CID^{\beta_k}_{\vlambda(\tau)_{\geq k;\leq 2}}
(\tau_{\geq k;\leq 2})
\boxtimes
\CX^{(0,\vo,\vgamma)}(\tau)$}}


    \put(5.5,4.5){\vector(-2,-1){2}}
    \put(3.5,2.5){\vector(2,-1){2}}
    \put(10.5,2.5){\vector(-2,-1){2}}
    \put(8.5,4.5){\vector(2,-1){2}}

   \put(3.5,4){\makebox(1,0.5){$\psi_2$}}
   \put(3,1.5){\makebox(1,0.5){$\phi_1\boxtimes\id$}}
   \put(10.1,1.5){\makebox(1,0.5){$\id\boxtimes\psi_2$}}
   \put(9.7,4){\makebox(1,0.5){$\psi_1$}}

  \end{picture}
\end{center}

commutes.}

\subsection{}
\label{tri trans maps k}
Let $\CX$ be as above.
For each $\vmu\in X^K$,
$\vmu\geq\vlambda$, so that $\vmu=\vlambda+\vbeta'$ for
some $\vbeta\in\BN[I]^K$, and
$\alpha\in\BN[I]$,
let us define a sheaf $\CX^{\alpha}_{\vmu}(K)\in
\CM(\CA^{\alpha}_{\vmu}(K);\CS)$ as
$\sigma_*\CX^{\alpha-\beta}(K)$. For example,
$\CX^{\alpha}_{\vlambda}(K)=\CX^{\alpha}(K)$.

Taking restrictions, the sheaves
$\CX^{\valpha}_{\vmu}(\tau)\in\CM(\CA^{\valpha}_{\vmu}(\tau);\CS)$
for all $K$-based trees $\tau$ are defined.

\subsection{}
Suppose $\CX,\CY$ are two factorizable sheaves supported at $\vc$,
$\vlambda=\vlambda(\CX),\ \vnu=\vlambda(\CY)$.
Let $\vmu\in X^K$, $\vmu\geq\vlambda,\ \vmu\geq\vnu$.
By definition, we have canonical isomorphisms
\begin{equation}
\label{tri thetas k}
\theta=\theta^{\alpha}_{\vmu;\vbeta}:\
\Hom_{\CA^{\alpha}_{\vmu}(K)}(\CX^{\alpha}_{\vmu}(K),\CY^{\alpha}_{\vmu}(K))
\iso
\Hom_{\CA^{\alpha+\beta}_{\vmu+\vbeta'}(K)}
(\CX^{\alpha+\beta}_{\vmu+\vbeta'}(K),\CY^{\alpha+\beta}_{\vmu+\vbeta'}(K))
\end{equation}
for each $\alpha\in\BN[I],\ \vbeta\in\BN[I]^K$.

\subsubsection{}
\label{tri taus kfsh} Suppose we are given
$\vbeta=(\beta_k)\in\BN[I]^K$. Let $\beta=\sum_k\ \beta_k$ as usually.
Choose a real $d$, $0<d<1$.

Consider a marked tree $(\tau_d,\ (0,\vbeta),\ \vmu)$ where
$$
\tau_d=(K\overset{\id}{\lra}K\lra \{*\};\ (0,d,1)).
$$
We have the restriction homomorphism
\begin{equation}
\label{tri def tau kfsh}
\xi_{\vmu;\vbeta;d}:
\Hom_{\CA^\beta_{\vmu}(K)}(\CX^\beta_{\vmu}(K),
\CY^\beta_{\vmu}(K))\lra
\Hom_{\CA^{(0,\vbeta)}_{\vmu}(\tau_d)}(\CX^{(0,\vbeta)}_{\vmu}(\tau_d),
\CY^{(0,\vbeta)}_{\vmu}(\tau_d))
\end{equation}

Suppose we are given $'\vbeta=('\beta_k)\in\BN[I]^K$ such that
$'\vbeta\leq\vbeta$. Let $'\beta=\sum_k\ '\beta_k$ as usually.
Choose a real $\varepsilon$, $0<\varepsilon<d$.

The restriction and the factorization isomorphisms $\psi$ induce the
map
\begin{equation}
\label{tri def theta kfsh}
\eta^{\vbeta;d}_{\vmu;\ '\vbeta;\varepsilon}:\
\Hom_{\CA^{(0,\vbeta)}_{\vmu}(\tau_d)}(\CX^{(0,\vbeta)}_{\vmu}(\tau_d),
\CY^{(0,\vbeta)}_{\vmu}(\tau_d))\lra
\Hom_{\CA^{(0,\ '\vbeta)}_{\vmu}(\tau_\varepsilon)}
(\CX^{(0,\ '\vbeta)}_{\vmu}(\tau_\varepsilon),
\CY^{(0,\ '\vbeta)}_{\vmu}(\tau_\varepsilon))
\end{equation}

The associativity axiom implies that these maps satisfy an obvious
transitivity property.

We define the space $\Hom_{^K\tFS}(\CX,\CY)$ as the following
inductive-projective limit
\begin{equation}
\label{tri ind proj k}
\Hom_{^K\tFS}(\CX,\CY):=
\dirlim\ \invlim\
\Hom_{\CA^{\alpha}_{\vmu}(K)}
(\CX^{\alpha}_{\vmu}(K),\CY^{\alpha}_{\vmu}(K))
\end{equation}
where the inverse limit is understood as follows. Its elements are collections
of maps
$$
\{ f^{\alpha}_K:
\CX^{\alpha}_{\vmu}(K)\lra\CY^{\alpha}_{\vmu}(K)\}
$$
given for all $\alpha\in\BN[I]$, $\vbeta\in\BN[I]^K$,
such that for every $\alpha,\ '\vbeta\leq\vbeta,\
0<\varepsilon<d<1$ as above, we have
$$
\eta^{\vbeta;d}_{\vmu;\ '\vbeta;\varepsilon} \xi_{\vmu;\vbeta;d}(f^\beta_K)
=\xi_{\vmu;\ '\vbeta;\varepsilon}(f^{'\beta}_K)
$$
$\vmu$ being fixed. The
direct limit is taken over $\vmu\in X^K$ such that
$\vmu\geq\vlambda,\ \vmu\geq\vnu$, the transition maps being
induced by ~(\ref{tri thetas k}).

With these spaces of homomorphisms, factorizable sheaves supported at $\vc$
form a $\sk$-linear category to be denoted by $^K\tFS_{\vc}$. As we have already
mentioned, the category of factorizable sheaves $^K\tFS$ is by
definition the product ~(\ref{tri prod comp}).

\vspace{.8cm}
{\em FINITE SHEAVES}
\vspace{.6cm}

\subsection{Definition} {\em A sheaf $\CX\in\ ^K\tFS_{\vc}$ is called
{\em finite}
if there exists only finitely many
$\vbeta\in\BN[I]^K$ such that the singular support of
$\CX^\alpha_{\vlambda}(K)$ contains the conormal bundle to
$\sigma^{\alpha-\beta}_{\vlambda-\vbeta';\vbeta}
(\CA^{\alpha-\beta}_{\vlambda-\vbeta'})$ (see (\ref{tri sigma k}))
for $\alpha\geq\beta=\sum_k \beta_k$.

A sheaf $\CX=\oplus_{\vc}\CX_{\vc}\in\ ^K\tFS,\ \CX_{\vc}\in\ ^K\tFS_{\vc}$
is called finite if all $\CX_{\vc}$ are finite.}

\subsection{}
\label{tri delpos}
Suppose we are given finite sheaves $\CX,\CY\in\ ^K\FS_{\vc}$; and
$\vmu\geq\vlambda(\CX),\vlambda(\CY)$.
As in the proof of the Lemma ~\ref{tri stabilization}, one can see that there
exists $'\vbeta\in\BN[I]^K$ such that for any $\vbeta\geq\ '\vbeta$
the map
\begin{equation}
\label{}
\eta^{'\vbeta;d}_{\vmu;\vbeta;\varepsilon}:\
\Hom_{\CA^{(0,\vbeta)}_{\vmu}(\tau_d)}(\CX^{(0,\vbeta)}_{\vmu}(\tau_d),
\CY^{(0,\vbeta)}_{\vmu}(\tau_d))\lra
\Hom_{\CA^{(0,\ '\vbeta)}_{\vmu}(\tau_\varepsilon)}
(\CX^{(0,\ '\vbeta)}_{\vmu}(\tau_\varepsilon),
\CY^{(0,\ '\vbeta)}_{\vmu}(\tau_\varepsilon))
\end{equation}

is an isomorphism. We will identify all the spaces
$\Hom_{\CA^{(0,\vbeta)}_{\vmu}(\tau_d)}(\CX^{(0,\vbeta)}_{\vmu}(\tau_d),
 \CY^{(0,\vbeta)}_{\vmu}(\tau_d))$ with the help of the above isomorphisms,
and we will denote this stabilized space by
${\overline{\Hom}}_{^K\FS}(\CX,\CY)$. Evidently, it does not depend
on a choice
of $'\vbeta$.

Quite similarly to the {\em loc.cit} one can see that for any
$\vbeta\geq\ '\vbeta$ the map
\begin{equation}
\label{}
\xi_{\vmu;\vbeta;d}:\
\Hom_{\CA^\beta_{\vmu}(K)}(\CX^\beta_{\vmu}(K),
\CY^\beta_{\vmu}(K))\lra
\Hom_{\CA^{(0,\vbeta)}_{\vmu}(\tau_d)}(\CX^{(0,\vbeta)}_{\vmu}(\tau_d),
\CY^{(0,\vbeta)}_{\vmu}(\tau_d))
\end{equation}
is an injection.

Thus we may view
$\Hom_{\CA^\beta_{\vmu}(K)}(\CX^\beta_{\vmu}(K),
 \CY^\beta_{\vmu}(K))$ as the subspace of
${\overline{\Hom}}_{^K\FS}(\CX,\CY)$.

We define $\Hom_{^K\FS}(\CX,\CY)\subset{\overline{\Hom}}_{^K\FS}(\CX,\CY)$
as the projective limit of the system of subspaces
$\Hom_{\CA^\beta_{\vmu}(K)}(\CX^\beta_{\vmu}(K),
 \CY^\beta_{\vmu}(K))$, $\vbeta\geq\ '\vbeta$.

With such definition of morphisms finite factorizable sheaves supported at
$\vc$ form an abelian category to be denoted by $^K\FS_{\vc}$. We set by
definition
\begin{equation}
\label{tri prod comp fin}
^K\FS=\prod_{\vc\in\pi_0(\CA)^K}\ ^K\FS_{\vc}
\end{equation}

\section{Gluing}
\label{tri sec glu}

\subsection{} Let
$$
\CA^{\alpha}_{\mu;1}\subset\CA^{\alpha}_{\mu}
$$
denote an open configuration subspace parametrizing configurations
lying entirely inside the unit disk $D(0;1)$. Due to monodromicity,
the restriction functors
$$
\CM(\CA^{\alpha}_{\mu};\CS)\lra\CM(\CA^{\alpha}_{\mu;1};\CS)
$$
are equivalences.

Let $\{*\}$ denote a one-element set. We have closed embeddings
$$
i:\CA_{\mu;1}^{\alpha}\hra\CA_{\mu}^{\alpha}(\{*\}),
$$
which identify the first space with the subspace of the second one
consisting of configurations with the small disk centered at $0$.
The inverse image functors
\begin{equation}
\label{tri res i}
i^*[-1]:\CM(\CA_{\mu}^{\alpha}(\{*\});\CS)\lra
\CM(\CA_{\mu}^{\alpha};\CS)
\end{equation}
are equivalences, again due to monodromicity.
Thus, we get equivalences
$$
\CM(\CA^{\alpha}_{\mu};\CS)\iso \CM(\CA^{\alpha}_{\mu}(\{*\};\CS)
$$
which induce canonical equivalences
\begin{equation}
\label{tri equiv pt tild}
\tFS\iso\tFS^{\{*\}}
\end{equation}
and
\begin{equation}
\label{tri equiv pt}
\FS\iso\FS^{\{*\}}
\end{equation}
Using these equivalences, we will sometimes identify these categories.

\subsection{Tensor product of categories}
Let $\CB_1,\CB_2$ be $\sk$-linear abelian categories. Their tensor product
category $\CB_1\otimes\CB_2$  is defined in \S5 of ~\cite{d2}.
It comes together with a canonical right biexact functor
$\CB_1\times\CB_2\lra\CB_1\otimes\CB_2$, and it is the initial object
among such categories.

\subsubsection{Basic Example} Let $M_i,\ i=1,2,$ be complex algebraic
varieties equipped with algebraic Whitney stratifications $\CS_i$.
Let $\CB_i=\CM(M_i;\CS_i)$. Then
$$
\CB_1\otimes\CB_2=\CM(M_1\times M_2;\CS_1\times\CS_2).
$$
The canonical functor
$\CB_1\times\CB_2\lra\CB_1\otimes\CB_2$ sends $(\CX_1,\CX_2)$ to
$\CX_1\boxtimes\CX_2$.

\subsubsection{} Recall the notations of ~\ref{tri taus kfsh}.
Let us consider the following category $\FS^{\otimes K}$.
Its objects are the collections of perverse sheaves
$\CX^{(0,\vbeta)}_{\vmu}(\tau_d)$ on the spaces
$\CA^{(0,\vbeta)}_{\vmu}(\tau_d)$ for sufficiently small $d$,
satisfying the usual factorization and finiteness conditions.
The morphisms are defined via the inductive-projective system with
connecting maps $\eta^{\vbeta;d}_{\vmu;\ '\vbeta;\varepsilon}$.
Using the above Basic Example, one can see easily that the category
$\FS^{\otimes K}$ is canonically equivalent to $\FS\otimes\ldots\otimes\FS$
($K$ times) which justifies its name.

By definition, the category $^K\FS$ comes together with the functor
$p_K:\ ^K\FS\lra\FS^{\otimes K}$ injective on morphisms. In effect,
$$
\Hom_{^K\FS}(\CX,\CY)\hra{\overline{\Hom}}_{^K\FS}(\CX,\CY)=
\Hom_{\FS^{\otimes K}}(p_K(\CX),p_K(\CY)).
$$
Let us construct a functor in the opposite direction.

\subsection{Gluing of factorizable sheaves}
For each $0<d<1$ let us consider a tree
$$
\tau_d=(K\overset{\id}{\lra}K\lra\{*\};(1,d,0)).
$$
Suppose we are given $\alpha\in\BN[I]$. Let $\CV(\alpha)$ denote the set of
all enhancements
$\valpha=(\alpha_*;(\alpha_k)_{k\in K})$ of $\tau$ such that
$\alpha_*+\sum_{k\in K}
\alpha_k=\alpha$. Obviously, the open subspaces
$\CA^{\valpha}(\tau_d)\subset\CA^{\alpha}(K)$, for varying $d$ and
$\valpha\in\CV(\alpha)$, form an open covering of $\CA^{\alpha}(K)$.

Suppose we are given a collection of factorizable sheaves
$\CX_k\in\FS_{c_k},\ k\in K$. Set $\vlambda=(\lambda(\CX_k))\in X^K$.
For each $d,\valpha$ as above
consider a sheaf
$$
\CX^{\valpha}(\tau_d):=
\CID^{\alpha_*}_{\vlambda_{\leq 1}}(\tau_{d;\leq 1})
\boxtimes\fbox{$\times$}_{k\in K}
\ \CX_k^{\alpha_k}
$$
over $\CA^{\valpha}_{\vlambda}(\tau_d)$.

Non-trivial pairwise intersections of the above open subspaces look as follows.
For $0<d_2<d_1<1$, consider a tree of height 3
$$
\varsigma=\varsigma_{d_1,d_2}=
(K\overset{\id}{\lra} K\overset{\id}{\lra}K\lra\{*\};
(1,d_1,d_2,0)).
$$
We have $\dpar_1\varsigma=\tau_{d_2},\ \dpar_2\varsigma=\tau_{d_1}$.
Let $\vbeta=(\beta_*,(\beta_{1;k})_{k\in K},(\beta_{2;k})_{k\in K})$
be an enhancement of $\varsigma$. Set $\valpha_1=\dpar_2\vbeta,\
\valpha_2=\dpar_1\vbeta$. Note that $\vbeta$ is defined uniquely by
$\valpha_1,\valpha_2$.
We have
$$
\CA^{\vbeta}_{\vlambda}(\varsigma)=
\CA^{\valpha_1}_{\vlambda}(\tau_{d_1})\cap
\CA^{\valpha_2}_{\vlambda}(\tau_{d_2}).
$$
Due to the factorization property for sheaves $\CID$ and $\CX_k$
we have isomorphisms
$$
\CX^{\valpha_1}(\tau_{d_1})|
_{\CA^{\vbeta}_{\vlambda}(\varsigma)}
\cong
\CID^{\dpar_3\vbeta}_{\dpar_3\vlambda}(\dpar_3\varsigma)\boxtimes
\fbox{$\times$}_{k\in K}(\CX_k)^{\beta_{2;k}}_{\lambda_k},
$$
and
$$
\CX^{\valpha_2}(\tau_{d_2})|
_{\CA^{\vbeta}_{\vlambda}(\varsigma)}
\cong
\CID^{\dpar_3\vbeta}_{\dpar_3\vlambda}(\dpar_3\varsigma)\boxtimes
\fbox{$\times$}_{k\in K}(\CX_k)^{\beta_{2;k}}_{\lambda_k}
$$
Taking composition, we get isomorphisms
\begin{equation}
\label{tri glue iso}
\phi_{d_1,d_2}^{\valpha_1,\valpha_2}:
\CX^{\valpha_1}(\tau_{d_1})|
_{\CA^{\valpha_1}_{\vlambda}(\tau_{d_1})\cap
  \CA^{\valpha_2}_{\vlambda}(\tau_{d_2})}\iso
\CX^{\valpha_2}(\tau_{d_2})|
_{\CA^{\valpha_1}_{\vlambda}(\tau_{d_1})\cap
  \CA^{\valpha_2}_{\vlambda}(\tau_{d_2})}
\end{equation}
{}From the associativity of the factorization for the sheaves $\CID$ and
$\CX_k$ it follows
that the isomorphisms ~(\ref{tri glue iso}) satisfy the cocycle condition;
hence they define a sheaf $\CX^{\alpha}(K)$ over
$\CA^{\alpha}_{\vlambda}(K)$.

Thus, we have defined a collection of sheaves $\{\CX^{\alpha}(K)\}$.
Using the corresponding data for the sheaves $\CX_k$, one defines easily
factorization isomorphisms ~\ref{tri def kfactor sh} (d)  and check that
they satisfy the associativity property. One also sees immediately that
the collection of sheaves $\{\CX^{\alpha}(K)\}$ is finite.
We leave this verification to the reader.

This way we get maps
$$
\prod_{k}\Ob(\FS_{c_k})\lra\Ob(^K\FS_{\vc}),\ \vc=(c_k)
$$
which extend by additivity to the map
\begin{equation}
\label{tri map glue}
g_K: \Ob(\FS^K)\lra \Ob(^K\FS)
\end{equation}

To construct the functor
\begin{equation}
\label{tri glu fun}
g_K: \FS^K\lra\ ^K\FS
\end{equation}

it remains to define $g_K$ on morphisms.

Given two collections of finite factorizable sheaves
$\CX_k,\CY_k\in\FS_{c_k},\ k\in K$,
let us choose $\vlambda=(\lambda_k)_{k\in K}$ such that $\lambda_k\geq
\lambda(\CX_k),\lambda(\CY_k)$ for all $k\in K$.
Suppose we have a collection of morphisms $f_k:\ \CX_k\lra\CY_k,\ k\in K$;
that is the maps $f_k^{\alpha_k}:\ (\CX_k)^{\alpha_k}_{\lambda_k}\lra
(\CY_k)^{\alpha_k}_{\lambda_k}$ given for any $\alpha_k\in\BN[I]$ compatible
with factorizations.

Given $\alpha\in\BN[I]$ and an enhancement $\valpha\in\CV(\alpha)$
as above we define the morphism
$f^{\valpha}(\tau_d):\ \CX^{\valpha}(\tau_d)\lra\CY^{\valpha}(\tau_d)$
over $\CA^{\valpha}_{\vlambda}(\tau_d)$ as follows:
it is the tensor product of the identity on
$\CID^{\alpha_*}_{\vlambda_{\leq 1}}(\tau_{d;\leq 1})$
with the morphisms $f_k^{\alpha_k}:\
(\CX_k)^{\alpha_k}_{\lambda_k}\lra
(\CY_k)^{\alpha_k}_{\lambda_k}$.

One sees easily as above that the morphisms $f^{\valpha}(\tau_d)$ glue
together to give a morphism $f^\alpha(K):\ \CX^\alpha(K)\lra\CY^\alpha(K)$;
as $\alpha$ varies they provide a morphism
$f(K):\ g_K((\CX_k))\lra g_K((\CY_k))$. Thus we have defined
the desired functor $g_K$.
Obviously it is $K$-exact, so by universal property it defines the same named
functor

\begin{equation}
\label{tri funct glue}
g_K:\ \FS^{\otimes K}\lra\ ^K\FS
\end{equation}

By the construction, the composition $p_K\circ g_K:\ \FS^{\otimes K}\lra
\FS^{\otimes K}$ is isomorphic to the identity functor.
Recalling that $p_K$ is injective on morphisms we see that $g_K$ and $p_K$
are quasiinverse. Thus we get

\subsection{Theorem} {\em The functors $p_K$ and $g_K$ establish a
canonical equivalence
$$
^K\FS\iso\FS^{\otimes K}\ \Box
$$}

\newpage
\section{Fusion}

\vspace{.5cm}
{\em BRAIDED TENSOR CATEGORIES}
\vspace{.8cm}

In this part we review the definition of a braided tensor category
following Deligne, ~\cite{d1}.

\subsection{} Let $\CC$ be a category, $Y$ a locally connected
locally simply connected topological space. By a {\em local system}
over $Y$ with values in $\CC$ we will mean a locally constant sheaf
over $Y$ with values in $\CC$. They form a category to be denoted
by $\Locsys(Y;\CC)$.

\subsubsection{} We will use the following basic example.
If $X$ is a complex algebraic variety
with a Whitney stratification $\CS$ then the category
$\CM(X\times Y;\CS\times\CS_{Y;tr})$ is equivalent to
$\Locsys(Y;\CM(X;\CS))$. Here $S_{Y;tr}$ denotes the trivial
stratification of $Y$, i.e. the first category consists
of sheaves smooth along $Y$.

\subsection{} Let $\pi:K\lra L$ be an epimorphism of non-empty
finite sets. We will use the notations of ~\ref{tri op disk}.
For real $\epsilon, \delta$ such that
$1>\epsilon>\delta> 0$, consider a tree
$$
\tau_{\pi;\epsilon,\delta}=(K\overset{\pi}{\lra}L\lra\{*\};\
(1,\epsilon,\delta))
$$
We have an isomorphism which is a particular case of ~(\ref{tri subst}):
$$
\CO(\tau_{\pi;\epsilon,\delta})\cong
\CO(L;\epsilon)\times\prod_{l\in L}\ \CO(K_l;\delta\epsilon^{-1})
$$
where $K_l:=\pi^{-1}(l)$.

\subsubsection{} {\bf Lemma} {\em There exists essentially unique functor
$$
r_{\pi}:\Locsys(\CO(K);\CC)\lra
\Locsys(\CO(L)\times\prod_l\ \CO(K_l);\CC)
$$
such that for each $\epsilon,\delta$ as above the square

\begin{center}
\begin{picture}(25,6)


\put(5,4){$\Locsys(\CO(K);\CC)$}

\put(9.3,4.1){$\vector(1,0){2.5}$}
\put(10,4.3){$r_{\pi}$}

\put(12,4){$\Locsys(\CO(L)\times\prod_l\ \CO(K_l);\CC)$}


\put(5.5,1){$\Locsys(\tau_{\pi;\epsilon,\delta};\CC)$}

\put(9.5,1.1){$\vector(1,0){2}$}
\put(10,1.3){$\sim$}

\put(12,1)
{$\Locsys(\CO(L;\epsilon)\times\prod_l\ \CO(K_l;\delta\epsilon^{-1});\CC)$}


\put(7.5,3.6){\vector(0,-1){2}}
\put(15.5,3.6){\vector(0,-1){2}}

\end{picture}
\end{center}

commutes.}

{\bf Proof} follows from the remark that $\CO(L)$ is a union of its open
subspaces
$$
\CO(L)=\bigcup_{\epsilon>0}\ \CO(L;\epsilon).\ \Box
$$

\subsection{}
\label{tri braided str} Let $\CC$ be a category. A {\em braided
tensor structure} on $\CC$ is the following collection of data.

(i) For each non-empty finite set $K$ a functor
$$
\otimes_K:\ \CC^K\lra\Locsys(\CO(K);\CC),\ \{X_k\}\mapsto\otimes_K\ X_k
$$
from the $K$-th power of $\CC$ to the category of local systems
(locally constant sheaves) over the space $\CO(K)$ with values
in $\CC$ (we are using the notations of ~\ref{tri op disk}).

We suppose
that $\otimes_{\{*\}}\ X$ is the constant local system with the fiber $X$.

(ii) For each $\pi:K\lra L$ as above a natural isomorphism
$$
\phi_{\pi}:(\otimes_K\ X_k)|_{\CO(L)\times\prod\CO(K_l)}\iso
\otimes_L\ (\otimes_{K_l}\ X_k).
$$
To simplify the notation,
we will write this isomorphism in the form
$$
\phi_{\pi}:
\otimes_K\ X_k\iso\otimes_L\ (\otimes_{K_l}\ X_k),
$$
implying that in the left hand side we must take restriction.

These isomorphisms must satisfy the following

{\em Associativity axiom.} For each pair of epimorphisms
$K\overset{\pi}{\lra} L\overset{\rho}{\lra} M$ the square

\begin{center}
  \begin{picture}(14,6)


    \put(5,4){\makebox(4,2){$\otimes_K\ X_k$}}


    \put(0,2){\makebox(4,2)
{$\otimes_M(\otimes_{K_m}\ X_k)$}}


    \put(10,2){\makebox(4,2)
{$\otimes_L(\otimes_{K_l}\ X_k)$}}

    \put(5,0){\makebox(4,2)
{$\otimes_M(\otimes_{L_m}(\otimes_{K_l}\ X_k))$}}


    \put(5.5,4.5){\vector(-2,-1){2}}
    \put(3.5,2.5){\vector(2,-1){2}}
    \put(10.5,2.5){\vector(-2,-1){2}}
    \put(8.5,4.5){\vector(2,-1){2}}

   \put(3.1,4){\makebox(1,0.5){$\phi_{\rho\pi}$}}
   \put(9.7,4){\makebox(1,0.5){$\phi_{\pi}$}}
   \put(2.1,1.5){\makebox(1,0.5)
{$\otimes_M\ \phi_{\pi|_{K_m}}$}}
   \put(10.1,1.5){\makebox(1,0.5){$\phi_{\rho}$}}

  \end{picture}
\end{center}

where $K_m:=(\rho\pi)^{-1}(m),\ L_m:=\rho^{-1}(m)$, commutes.

\subsection{} The connection with the conventional definition is as
follows. Given two objects $X_1,X_2\in\Ob\ \CC$, define an object
$X_1\dotimes X_2$ as the fiber of $\otimes_{\{1,2\}} X_k$ at the point
$(1/3,1/2)$.

We have natural isomorphisms
$$
A_{X_1,X_2,X_3}:
X_1\dotimes (X_2\dotimes X_3)\iso (X_1\dotimes X_2)\dotimes X_3
$$
coming from isomorphisms $\phi$ associated with two possible
order preserving epimorphic maps $\{1,2,3\}\lra\{1,2\}$, and
$$
R_{X_1,X_2}: X_1\dotimes X_2\iso X_2\dotimes X_1
$$
coming from the standard half-circle monodromy.
Associativity axiom for $\phi$ is equivalent to the
the usual compatibilities for these maps.

\subsection{} Now suppose that the data ~\ref{tri braided str} is given
for {\em all} (possibly empty) tuples and all (not necessarily
epimorphic) maps. The space $\CO(\emp)$ is by definition
one point, and a local system $\otimes_{\emp}$ over it
is simply an object of $\CC$; let us denote it $\One$ and
call a {\em unit} of our tensor structure. In this case
we will say that $\CC$ is a braided tensor category with unit.

In the conventional language, we have natural isomorphisms
$$
\One\dotimes X\iso X
$$
(they correspond to $\{2\}\hra\{1,2\}$)
satisfying the usual compatibilities with $A$ and $R$.

\vspace{.5cm}
{\em FUSION FUNCTORS}
\vspace{.8cm}

\subsection{} Let
$$
\CA_{\alpha;1}\subset\CA_{\alpha}
$$
denote the open subspace parametrizing configurations lying
inside the unit disk $D(0;1)$.

Let $K$ be a non-empty finite set.
Obviously, $\CA^{\alpha}(K)=\CA_{\alpha;1}\times\CO(K)$,  and we have
the projection
$$
\CA^{\alpha}(K)\lra\CO(K).
$$
Note also that we have an evident open embedding $\CO(K)\hra D(0;1)^K$.

Our aim in this part is to define certain {\em fusion functors}
$$
\Psi_K:\CCD(\CA^{\alpha}(K))\lra\CCD^{mon}(\CA^{\alpha}(\{*\})\times\CO(K))
$$
where $(\bullet)^{mon}$ denotes the full subcategory of complexes
smooth along $\CO(K)$. The construction follows the classical
definition of nearby cycles functor, ~\cite{d3}.

\subsection{Poincar\'{e} groupoid} We start with a topological notation.
Let $X$ be a connected locally simply connected topological space.
Let us denote by
$\widetilde{X\times X}$ the space
whose points are triples $(x,y,p)$, where $x,y\in X$; $p$ is a homotopy
class of paths in $X$ connecting $x$ with $y$.
Let
\begin{equation}
\label{tri free cov}
c_X:\widetilde{X\times X}\lra X\times X
\end{equation}
be the evident projection. Note that for a
fixed $x\in X$, the restriction of $c_X$ to $c_X^{-1}(X\times\{ x\})$
is a universal covering of $X$ with a group
$\pi_1(X;x)$.

\subsection{} Consider the diagram with cartesian squares

\begin{center}
  \begin{picture}(25,6)


\put(0,4){$\CA^{\alpha}(\{*\})\times\CO(K)$}
\put(6.5,4){$\widehat{\CA_{\alpha}(K)}\times\CO(K)$}
\put(12,4){$\CA^{\alpha}(K)\times\CO(K)$}
\put(18,4){$\widetilde{\CA^{\alpha}(K)\times\CO(K)}$}


\put(0,1){$D(0,1)\times\CO(K)$}
\put(6,1){$D(0;1)^K\times\CO(K)$}
\put(12,1){$\CO(K)\times\CO(K)$}
\put(18,1){$\widetilde{\CO(K)\times\CO(K)}$}


\put(2,3.6){\vector(0,-1){2}}
\put(8,3.6){\vector(0,-1){2}}
\put(14,3.6){\vector(0,-1){2}}
\put(20,3.6){\vector(0,-1){2}}


\put(4.1,4.1){$\vector(1,0){2}$}
\put(5,4.3){$\tDelta$}
\put(4,1.1){$\vector(1,0){1.8}$}
\put(5,1.3){$\Delta$}

\put(11.5,4.1){$\vector(-1,0){1.3}$}
\put(10.5,4.3){$\tj$}
\put(11.5,1.1){$\vector(-1,0){1.5}$}
\put(10.5,1.3){$j$}

\put(17.5,4.1){$\vector(-1,0){1.8}$}
\put(16.5,4.3){$\tc$}
\put(17.5,1.1){$\vector(-1,0){2}$}
\put(16.5,1.3){$c$}

  \end{picture}
\end{center}

where we have denoted
$\widehat{\CA^{\alpha}(K)}:=\CA_{\alpha;1}\times D(0;1)^K.$
Here $\Delta$ is induced by the diagonal embedding
$D(0;1)\hra D(0;1)^K$, $j$ --- by the open embedding
$\CO(K)\hra D(0;1)^K$, $c$ is the map ~(\ref{tri free cov}).
The upper horizontal arrows are defined by pull-back.

We define $\Psi_K$ as a composition
$$
\Psi_K=\tDelta^*\tj_*\tc_*\tc^*p^*[1]
$$
where $p: \CA^{\alpha}(K)\times\CO(K)\lra\CA^{\alpha}(K)$
is the projection.

This functor is $t$-exact and induces an exact functor
\begin{equation}
\label{tri def fus fun}
\Psi_K:\CM(\CA^{\alpha}(K);\CS)\lra
\CM(\CA^{\alpha}(\{*\})\times\CO(K);\CS\times\CS_{tr})
\end{equation}
where $\CS_{tr}$ denotes the trivial stratification of $\CO(K)$.

\subsection{} Set
$$
\CA^{\alpha}(K)_d:=\CA_{\alpha;1}\times\CO(K;d)
$$
The category $\CM(\CA^{\alpha}(K);\CS)$ is equivalent to the
"inverse limit" $"\invlim"\CM(\CA^{\alpha}(K)_d;\CS)$.

Let $\pi:K\lra L$ be an epimorphism. Consider a configuration space
$$
\CA^{\alpha}(\tau_{\pi;d}):=\CA_{\alpha;1}\times\CO(\tau_{\pi;d})
$$
where $\tau_{\pi;d}:=\tau_{\pi;d,0}$.
An easy generalization of the definition of $\Psi_K$ yields
a functor
$$
\Psi_{\pi;d}:\CM(\CA^{\alpha}(\tau_{\pi;d}))\lra
\CM(\CA^{\alpha}(L)_d\times\prod_{l\in L}\CO(K_l))
$$
(In what follows we will omit for brevity stratifications from the notations
of abelian categories $\CM(\bullet)$,
implying that we use the principal stratification on all configuration
spaces $\CA^{\alpha}(\bullet)$ and the trivial stratification
on spaces $\CO(\bullet)$, i.e. our sheaves are smooth along
these spaces.) Passing to the limit over $d>0$ we conclude that there
exists essentially unique functor
\begin{equation}
\label{tri psi kl}
\Psi_{K\lra L}:\CM(\CA^{\alpha}(K))\lra
\CM(\CA^{\alpha}(L)\times\prod_l\ \CO(K_l))
\end{equation}
such that all squares

\begin{center}
\begin{picture}(25,6)


\put(6.5,4){$\CM(\CA^{\alpha}(K))$}

\put(9.3,4.1){$\vector(1,0){2.5}$}
\put(10,4.3){$\Psi_{K\lra L}$}

\put(12,4){$\CM(\CA^{\alpha}(L)\times\prod_l\ \CO(K_l))$}


\put(6.5,1){$\CM(\CA^{\alpha}(\tau_{\pi;d}))$}

\put(9.5,1.1){$\vector(1,0){2}$}
\put(10,1.3){$\Psi_{\pi;d}$}

\put(12,1){$\CM(\CA^{\alpha}(L)_d\times\prod_l\ \CO(K_l))$}


\put(7.5,3.6){\vector(0,-1){2}}
\put(14.5,3.6){\vector(0,-1){2}}

\end{picture}
\end{center}

commute (the vertical arrows being restrictions). If $L=\{*\}$,
we return to $\Psi_K$.

\subsection{Lemma.} {\em All squares

\begin{center}
\begin{picture}(25,6)


\put(6.5,4){$\CM(\CA^{\alpha}(K))$}

\put(9.3,4.1){$\vector(1,0){2.5}$}
\put(10,4.3){$\Psi_{K}$}

\put(12,4){$\CM(\CA^{\alpha}(\{*\})\times\CO(K))$}


\put(4.5,1){$\CM(\CA^{\alpha}(L)\times\prod_l\ \CO(K_l))$}

\put(9.8,1.1){$\vector(1,0){1.7}$}
\put(10.2,1.3){$\Psi_{L}$}

\put(12,1){$\CM(\CA^{\alpha}(\{*\})\times\CO(L)\times\prod_l\ \CO(K_l))$}


\put(7.5,3.6){\vector(0,-1){2}}
\put(5.5,2.6){$\Psi_{K\lra L}$}
\put(14.5,3.6){\vector(0,-1){2}}
\put(15,2.6){$r_{\pi}$}

\end{picture}
\end{center}

$2$-commute. More precisely, there exist natural isomorphisms
$$
\phi_{K\lra L}:r_{\pi}\circ\Psi_K\iso
\Psi_L\circ\Psi_{K\lra L.}
$$
These isomorphisms satisfy a natural cocycle condition
(associated with pairs of epimorphisms $K\lra L\lra M$). }

\subsection{} Applying the functors $\Psi_K$ componentwise, we get functors
$$
\Psi_K:\ ^K\FS\lra\Locsys(\CO(K);\FS);
$$
taking composition with the gluing functor $g_K$, ~(\ref{tri glu fun}),
we get functors
\begin{equation}
\label{tri tens prods fs}
\otimes_K:\ \FS^K\lra\Locsys(\CO(K);\FS)
\end{equation}
It follows from the previous lemma that these functors define
a braided tensor structure on $\FS$.

\subsection{} Let us define a unit object in $\FS$ as
$\One_{\FS}=\CL(0)$ (cf. \ref{tri trivial}). One can show that it is
a unit for the braided tensor structure defined above.

\newpage
\begin{center}
{\bf Chapter 3. Functor $\Phi$}
\end{center}
\vspace{.8cm}

\section{Functor $\Phi$}

\subsection{}
Recall the category $\CC$ defined in II.11.3.2 and II.12.2.

Our main goal in this section will be the construction of a tensor functor
$\Phi:\ \FS\lra \CC$.

\subsection{}
\label{tri maps} Recall that we have already defined in ~\ref{tri phi grad}
a functor
$$
\Phi:\FS\lra\Vect_f^X.
$$
Now we will construct natural transformations
$$
\epsilon_i:\ \Phi_{\lambda}(\CX)\lra\Phi_{\lambda+i'}(\CX)
$$
and
$$
\theta_i:\ \Phi_{\lambda+i'}(\CX)\lra\Phi_{\lambda}(\CX).
$$
We may, and will, assume that $\CX\in\FS_c$ for some $c$.
If $\lambda\not\in X_c$ then there is nothing to do.

Suppose that $\lambda\in X_c$; pick $\alpha\in\BN[I]$ such that
$\lambda+\alpha'\geq\lambda(\CX)$. By definition.
$$
\Phi_{\lambda}(\CX)=\Phi_{\alpha}(\CX_{\lambda+\alpha'}^{\alpha})
$$
where $\Phi_{\alpha}$ is defined in II.7.14 (the definition will be
recalled below).

\subsection{}
\label{tri unfoldings} Pick an unfolding $\pi:\ J\lra I$ of $\alpha$, II.6.12;
we will use the same
notation for the canonical projection
$$
\pi:\ ^{\pi}\BA\lra\CA^{\alpha}_{\lambda+\alpha'}=\CA_{\alpha}.
$$
Let $N$ be the dimension of $\CA_{\alpha}$.

\subsection{}
\label{tri facets} Recall some notations from II.8.4. For each
$r\in [0,N]$ we have denoted by $\CP_r(J;1)$ the set of all
maps
$$
\varrho: J\lra [0,r]
$$
such that $\varrho(J)$ contains $[r]$. Let us assign to such
$\varrho$ the real point $w_{\varrho}=(\varrho(j))_{j\in J}\in\ ^{\pi}\BA$.

There exists a unique positive facet of $\CS_{\BR}$, $F_{\varrho}$
containing $w_{\varrho}$. This establishes a bijection
between $\CP_r(J;1)$ and the set $\Fac_r$ of $r$-dimensional
positive facets. At the same time we have fixed on each
$F_{\varrho}$ a point $w_{\varrho}$. This defines cells
$D^+_{\varrho}:=D^+_{F_{\varrho}},\ S^+_{\varrho}:=S^+_{F_{\varrho}}$, cf.
II.7.2.

Note that this "marking" of positive facets is $\Sigma_{\pi}$-invariant.
In particular, the group $\Sigma_{\pi}$ permutes the above mentioned cells.

We will denote by $\{0\}$ the unique zero-dimensional facet.

\subsection{} Given a complex $\CK$ from the bounded derived category
$\CCD^b(\CA_{\alpha})$, its inverse image $\pi^*\CK$ is correctly
defined as an element of the {\em equivariant} derived category
$\CCD^b(^{\pi}\BA,\Sigma_{\pi})$ obtained by localizing the
category of $\Sigma_{\pi}$-equivariant complexes on $^{\pi}\BA$.
The direct image $\pi_*$ acts between equivariant derived categories
$$
\pi_*:\CCD^b(^{\pi}\BA,\Sigma_{\pi})\lra\CCD^b(\CA_{\alpha},\Sigma_{\pi})
$$
(the action of $\Sigma_{\pi}$ on $\CA_{\pi}$ being trivial).

We have the functor of $\Sigma_{\pi}$-invariants
\begin{equation}
\label{tri inv}
(\bullet)^{\Sigma_{\pi}}:\CCD^b(\CA_{\alpha},\Sigma_{\pi})\lra
\CCD^b(\CA_{\alpha})
\end{equation}

\subsubsection{} {\bf Lemma.} {\em For every $\CK\in\CCD^b(\CA_{\alpha})$ the
canonical morphism
$$
\CK\lra (\pi_*\pi^*\CK)^{\Sigma_{\pi}}
$$
is an isomorphism.}

{\bf Proof.} The claim may be checked fiberwise. Taking of a fiber
commutes with taking $\Sigma_{\pi}$-invariants since our group
$\Sigma_{\pi}$ is finite and we are living over the field of characteristic
zero, hence $(\bullet)^{\Sigma_{\pi}}$ is exact. After that, the claim
follows from definitions. $\Box$

\subsubsection{}
\label{tri inv rg} {\bf Corollary.} {\em For every $\CK\in\CCD^b(\CA_{\alpha})$
we have canonical isomorphism
\begin{equation}
\label{tri inv rgamma}
R\Gamma(\CA_{\alpha};\CK)\iso R\Gamma(^{\pi}\BA;\pi^*\CK)^{\Sigma_{\pi}}\ \Box
\end{equation}}

\subsection{} Following II.7.13, consider the sum of coordinates
function
$$
\sum t_j:\ ^{\pi}\BA\lra\BA^1;
$$
and for $\CL\in\CCD(^{\pi}\BA;\CS)$ let
$\Phi_{\sum t_j}(\CL)$ denote the fiber at the origin
of the corresponding vanishing cycles functor.
If $H\subset\ ^{\pi}\BA$ denotes the inverse image $(\sum t_j)^{-1}(\{ 1\})$
then we have canonical isomorphisms
\begin{equation}
\label{tri rel coh}
\Phi_{\sum t_j}(\CL)\cong R\Gamma(^{\pi}\BA,H;\CL)\cong
\Phi^+_{\{0\}}(\CL)
\end{equation}
The first one follows from the definition of vanishing cycles and the second
one from homotopy argument.

Note that the if $\CL=\pi^*\CK$ for some $\CK\in\CCD(\CA_{\alpha};\CS)$
then the group $\Sigma_{\pi}$ operates canonically on all terms
of ~(\ref{tri rel coh}), and the isomorphisms are $\Sigma_{\pi}$-equivariant.

Let us use the same notation
$$
\sum t_j:\CA_{\alpha}\lra\BA^1
$$
for the descended function, and for $\CK\in\CCD(\CA_{\alpha};\CS)$
let $\Phi_{\sum t_j}(\CK)$ denote the fiber at the origin of the
corresponding vanishing cycles functor.
If $\CK$ belongs to $\CM(\CA_{\alpha};\CS)$ then $\Phi_{\sum t_j}(\CK)$
reduces to a single vector space and this is what we call
$\Phi_{\alpha}(\CK)$.

If
$\bar{H}$ denotes $\pi(H)=(\sum t_j)^{-1}(\{1\})\subset\CA_{\alpha}$ then
we have canonical isomorphism
$$
\Phi_{\sum t_j}(\CK)\cong R\Gamma(\CA_{\alpha},\bar{H};\CK)
$$

\subsection{Corollary} (i) {\em For every $\CK\in\CCD(\CA_{\alpha};\CS)$
we have a canonical isomorphism
\begin{equation}
\label{tri phi equiv}
i_{\pi}: \Phi_{\sum t_j}(\CK)\iso
\Phi^+_{\{0\}}(\pi^*\CK)^{\Sigma_{\pi}}.
\end{equation}

(ii) This isomorphism does not depend on the choice of an unfolding
$\pi: J\lra I$.}

Let us explain what (ii) means. Suppose $\pi':J'\lra I$ be another
unfolding of $\alpha$. There exists (a non unique) isomorphism
$$
\rho:J\iso J'
$$
such that $\pi'\circ\rho=\pi$. It induces isomorphisms
$$
\rho^*:\Sigma_{\pi'}\iso\Sigma_{\pi}
$$
(conjugation by $\rho$), and
$$
\rho^*:\Phi^+_{\{0\}}((\pi')^*\CK)\iso
\Phi^+_{\{0\}}(\pi^*\CK)
$$
such that
$$
\rho^*(\sigma x)=\rho^*(\sigma)\rho^*(x),
$$
$\sigma\in\Sigma_{\pi'}$, $x\in \Phi^+_{\{0\}}((\pi')^*\CK)$.
Passing to invariants, we get an isomorphism
$$
\rho^*:\Phi^+_{\{0\}}((\pi')^*\CK)^{\Sigma_{\pi'}}\iso
\Phi^+_{\{0\}}(\pi^*\CK)^{\Sigma_{\pi}}.
$$
Now (ii) means that $i_{\pi}\circ\rho^*=i_{\pi'}$. As a consequence, the last
isomorphism does not depend on the choice of $\rho$.

{\bf Proof.} Part (i) follows from the preceding discussion and
{}~\ref{tri inv rg}.

As for (ii), it suffices to prove that any automorphism
$\rho: J\iso J$ respecting $\pi$
induces the identity automorphism of the space of invariants
$\Phi^+_{\{0\}}(\pi^*\CK)^{\Sigma_{\pi}}$. But the action
of $\rho$ on the space $\Phi^+_{\{0\}}(\pi^*\CK)$ comes
from the action of $\Sigma_{\pi}$ on this space, and our claim
is obvious.
$\Box$

In computations the right hand side of ~(\ref{tri phi equiv}) will
be used as a {\em de facto} definition of $\Phi_{\alpha}$.

\subsection{}
\label{tri disj} Suppose that $\alpha=\sum a_ii$; pick an $i$ such that
$a_i>0$.

Let us introduce the following notation. For a partition
$J=J_1\coprod J_2$ and a positive $d$ let
$\BA^{J_1,J_2}(d)$ denote an open suspace of $^{\pi}\BA$ consisting
of all points $\bt=(t_j)$ such that $|t_j|>d$ (resp., $|t_j|<d$)
if $j\in J_1$ (resp., $j\in J_2$).

We have obviously
\begin{equation}
\label{tri disjoint}
\pi^{-1}(\CA^{i,\alpha-i}(d))=
\coprod_{j\in \pi^{-1}(i)} \BA^{\{j\},J-\{j\}}(d)
\end{equation}

\subsection{}
\label{tri subgroups}
For $j\in \pi^{-1}(i)$ let $\pi_j: J-\{j\}\lra I$ denote the restriction
of $\pi$; it is an unfolding of $\alpha-i$. The group $\Sigma_{\pi_j}$
may be identified with the subgroup of $\Sigma_{\pi}$ consisting
of automorphisms stabilizing $j$. For $j', j''\in J$
let $(j'j'')$ denotes their transposition. We have
\begin{equation}
\label{tri conj}
\Sigma_{\pi_{j''}}=(j'j'')\Sigma_{\pi_{j'}}(j'j'')^{-1}
\end{equation}
For a fixed $j\in\pi^{-1}(i)$ we have a partition into cosets
\begin{equation}
\label{tri cosets}
\Sigma_{\pi}=\coprod_{j'\in\pi^{-1}(i)}\ \Sigma_{\pi_j}(jj')
\end{equation}

\subsection{}
\label{tri one dim}
For $j\in J$ let $F_j$ denote a one-dimensional facet
corresponding to the map $\varrho_j: J\lra [0,1]$ sending all elements to $0$
except for $j$ being sent to $1$ (cf. ~\ref{tri facets}).

For $\CK\in\CCD(\CA_{\alpha};\CS)$ we have canonical and variation maps
$$
v_j: \Phi^+_{\{0\}}(\pi^*\CK)\overset{\lra}{\lla}\Phi^+_{F_j}(\pi^*\CK): u_j
$$
defined in II, (89), (90). Taking their sum, we get maps
\begin{equation}
\label{tri sum uv}
v_{i}: \Phi^+_{\{0\}}(\pi^*\CK)\overset{\lra}{\lla}
\oplus_{j\in\pi^{-1}(i)}\ \Phi^+_{F_j}(\pi^*\CK): u_{i}
\end{equation}
Note that the group $\Sigma_{\pi}$ operates naturally on
both spaces and both maps $v_{i}$ and $u_{i}$ respect this action.

Let us consider more attentively the action of $\Sigma_{\pi}$ on
$\oplus_{j\in\pi^{-1}(i)}\ \Phi^+_{F_j}(\pi^*\CK)$. A subgroup
$\Sigma_{\pi_j}$ respects the subspace $\Phi^+_{F_j}(\pi^*\CK)$. A
transposition
$(j'j'')$ maps $\Phi^+_{F_{j'}}(\pi^*\CK)$ isomorphically
onto $\Phi^+_{F_{j''}}(\pi^*\CK)$.

Let us consider the space of invariants
$$
(\oplus_{j\in\pi^{-1}(i)}\ \Phi^+_{F_j}(\pi^*\CK))^{\Sigma_{\pi}}
$$
For every $k\in \pi^{-1}(i)$ the obvious projection induces isomorphism
$$
(\oplus_{j\in\pi^{-1}(i)}\ \Phi^+_{F_j}(\pi^*\CK))^{\Sigma_{\pi}}
\iso
(\Phi^+_{F_k}(\pi^*\CK))^{\Sigma_{\pi_k}};
$$
therefore for two different $k,k'\in\pi^{-1}(i)$ we get an
isomorphism
\begin{equation}
\label{tri ikk}
i_{kk'}: (\Phi^+_{F_k}(\pi^*\CK))^{\Sigma_{\pi_k}}
\iso
(\Phi^+_{F_{k'}}(\pi^*\CK))^{\Sigma_{\pi_{k'}}}
\end{equation}
Obviously, it is induced by transposition $(kk')$.

\subsection{}
\label{tri def theta eps} Let us return to the situation ~\ref{tri maps} and
apply the preceding discussion to $\CK=\CX^{\alpha}_{\lambda+\alpha'}$.
We have by definition
\begin{equation}
\label{tri phi l}
\Phi_{\lambda}(\CX)=
\Phi_{\sum t_j}(\CX^{\alpha}_{\lambda+\alpha'})\cong
\Phi^+_{\{0\}}(\pi^*\CX^{\alpha}_{\lambda+\alpha'})^{\Sigma_{\pi}}
\end{equation}
On the other hand, let us pick some $k\in\pi^{-1}(i)$ and a real $d$,
$0<d<1$. The subspace
\begin{equation}
\label{tri perp}
F_k^{\perp}(d)\subset\BA^{\{j\},J-\{j\}}(d)
\end{equation}
consisting of points $(t_j)$ with $t_k=1$, is a transversal slice to
the face $F_k$. Consequently,
the factorization isomorphism for $\CX_{\lambda+\alpha'}^{\alpha}$
lifted to $\BA^{\{j\},J-\{j\}}(d)$ induces isomorphism
$$
\Phi^+_{F_k}(\pi^*\CX_{\lambda+\alpha'}^{\alpha})\cong
\Phi_{\{0\}}^+(\pi^*_k\CX_{\lambda+\alpha'}^{\alpha-i})\otimes
(\CI^i_{\lambda+i'})_{\{1\}}=
\Phi_{\{0\}}^+(\pi^*_k\CX_{\lambda+\alpha'}^{\alpha-i})
$$
Therefore we get isomorphisms
\begin{eqnarray}
\label{tri phi li}
\Phi_{\lambda+i'}(\CX)=
\Phi_{\{0\}}^+(\pi^*_k\CX_{\lambda+\alpha'}^{\alpha-i})^{\Sigma_{\pi_k}}
\cong
\Phi^+_{F_k}(\pi^*\CX_{\lambda+\alpha'}^{\alpha})^{\Sigma_{\pi_k}}
\cong \\ \nonumber
\cong (\oplus_{j\in\pi^{-1}(i)}\
\Phi^+_{F_j}(\pi^*\CX_{\lambda+\alpha'}^{\alpha}))^{\Sigma_{\pi}}
\end{eqnarray}
It follows from the previous discussion that this isomorphism
does not depend on the intermediate choice of $k\in\pi^{-1}(i)$.

Now we are able to define the operators $\theta_i$,
$\epsilon_i$:
$$
\epsilon_i:\Phi_{\lambda}(\CX)\overset{\lra}{\lla}
\Phi_{\lambda+i'}(\CX) :\theta_i
$$
By definition, they are induced by the maps $u_i,\ v_i$ (cf. ~(\ref{tri sum uv}))
respectively
(for $\CK=\CX_{\lambda+\alpha'}^{\alpha}$) after
passing to $\Sigma_{\pi}$-invariants and
taking into account isomorphisms ~(\ref{tri phi l}) and ~(\ref{tri phi li}).

\subsection{Theorem.}
\label{tri relations} {\em The operators $\epsilon_i$ and $\theta_i$
satisfy the relations {\em II.12.3}, i.e. the previous construction defines
functor
$$
\tPhi:\FS\lra\tCC
$$
where the category $\tCC$ is defined as in {\em loc. cit}}.

{\bf Proof} will occupy the rest of the section.

\subsection{} Let $\fu^+$ (resp., $\fu^-$) denote the subalgebra of $\fu$
generated by all $\epsilon_i$ (resp., $\theta_i$). For $\beta\in\BN[I]$
let $\fu^{\pm}_{\beta}\subset\fu^{\pm}$ denote the corresponding homogeneous
component.

The proof will go as follows. First, relations II.12.3 (z), (a), (b) are
obvious. We will do the rest in three steps.

{\em Step 1.} Check of (d). This is equivalent to showing that
the action of operators $\theta_i$ correctly defines maps
$$
\fu^-_{\beta}\otimes\Phi_{\lambda+\beta'}(\CX)\lra\Phi_{\lambda}(\CX)
$$
for all $\beta\in\BN[I]$.

{\em Step 2.} Check of (e). This is equivalent to showing that
the action of operators $\epsilon_i$ correctly defines maps
$$
\fu^+_{\beta}\otimes\Phi_{\lambda}(\CX)\lra\Phi_{\lambda+\beta'}(\CX)
$$
for all $\beta\in\BN[I]$.

{\em Step 3.} Check of (c).

\subsection{} Let us pick an arbitrary $\beta=\sum b_ii\in\BN[I]$ and
$\alpha\in\BN[I]$ such that $\lambda+\alpha'\geq\lambda(\CX)$ and
$\alpha\geq\beta$. We pick the data from ~\ref{tri unfoldings}. In what
follows we will generalize the considerations of ~\ref{tri disj} ---
{}~\ref{tri def theta eps}.

Let $U(\beta)$ denote the set of all subsets $J'\subset J$ such that
$|J'\cap\pi^{-1}(i)|=b_i$ for all $i$. Thus, for such $J'$,
$\pi_{J'}:=\pi|_{J'}:J'\lra I$ is an unfolding of $\beta$ and
$\pi_{J-J'}$ --- an unfolding of $\alpha-\beta$. We have a disjoint
sum decomposition
\begin{equation}
\label{tri disjoint beta}
\pi^{-1}(\CA^{\beta,\alpha-\beta}(d))=
\coprod_{J'\in U(\beta)}\ \BA^{J',J-J'}(d)
\end{equation}
(cf. ~(\ref{tri disjoint})).

\subsection{}
For $J'\in U(\beta)$ let $F_{J'}$ denote a one-dimensional facet
corresponding to the map $\varrho_{J'}: J\lra [0,1]$ sending $j\not\in J'$ to
$0$
$j\in J'$ --- to $1$ (cf. ~\ref{tri facets}).

For $\CK\in\CCD(\CA_{\alpha};\CS)$ we have canonical and variation maps
$$
v_{J'}: \Phi^+_{\{0\}}(\pi^*\CK)\overset{\lra}{\lla}
\Phi^+_{F_{J'}}(\pi^*\CK): u_{J'}
$$
Taking their sum, we get maps
\begin{equation}
\label{tri sum uv beta}
v_{\beta}: \Phi^+_{\{0\}}(\pi^*\CK)\overset{\lra}{\lla}
\oplus_{J'\in U(\beta)}\ \Phi^+_{F_{J'}}(\pi^*\CK): u_{\beta}
\end{equation}
The group $\Sigma_{\pi}$ operates naturally on
both spaces and both maps $v_{\beta}$ and $u_{\beta}$ respect this action.

A subgroup $\Sigma_{J'}$ respects the subspace
$\Phi^+_{F_{J'}}(\pi^*\CK)$. The projection induces an isomorphism
$$
(\oplus_{J'\in U(\beta)}\ \Phi^+_{F_{J'}}(\pi^*\CK))^{\Sigma_{\pi}}\iso
\Phi^+_{F_{J'}}(\pi^*\CK)^{\Sigma_{J'}}.
$$
We have the crucial

\subsection{Lemma} {\em Factorization isomorphism for $\CX$ induces
canonical isomorphism
\begin{equation}
\label{tri phi lbeta}
\fu^-_{\beta}\otimes\Phi_{\lambda+\beta'}(\CX)\cong
(\oplus_{J'\in U(\beta)}\ \Phi^+_{F_{J'}}(\pi^*\CX^{\alpha}_{\lambda+\alpha'}))
^{\Sigma_{\pi}}
\end{equation}}

{\bf Proof.} The argument is the same as in ~\ref{tri def theta eps}, using
II, Thm. 6.16. $\Box$

\subsection{} As a consequence, passing to $\Sigma_{\pi}$-invariants
in ~(\ref{tri sum uv beta}) (for $\CK=\CX^{\alpha}_{\lambda+\alpha'}$)
we get the maps
$$
\epsilon_{\beta}:\Phi_{\lambda}(\CX)
\overset{\lra}{\lla}
\fu^-_{\beta}\otimes\Phi_{\lambda+\beta'}(\CX):\theta_{\beta}
$$

\subsection{Lemma}
\label{tri assoc theta} {\em
The maps $\theta_{\beta}$ provide $\Phi(\CX)$ with a structure of a left
module over the negative subalgebra $\fu^-$.}

{\bf Proof.} We must prove the associativity. It follows from the
associativity of factorization isomorphisms. $\Box$

This lemma proves relations II.12.3 (d)
for operators $\theta_i$, completing Step 1 of the proof of our theorem.

\subsection{} Now let us consider operators
$$
\epsilon_{\beta}:\Phi_{\lambda}(\CX)\lra
\fu_{\beta}^-\otimes\Phi_{\lambda+\beta'}(\CX).
$$
By adjointness, they induce operators
$$
\fu^{-*}_{\beta}\otimes\Phi_{\lambda}(\CX)\lra
\Phi_{\lambda+\beta'}(\CX)
$$
The bilinear form $S$, II.2.10, induces isomorhisms
$$
S:\fu_{\beta}^{-*}\iso\fu_{\beta}^{-};
$$
let us take their composition with the isomorphism
of algebras
$$
\fu^-\iso\fu^+
$$
sending $\theta_i$ to $\epsilon_i$. We get isomorphisms
$$
S':\fu_{\beta}^{-*}\iso\fu_{\beta}^+
$$
Using $S'$, we get from $\epsilon_{\beta}$ operators
$$
\fu_{\beta}^+\otimes\Phi_{\lambda}(\CX)\lra\Phi_{\lambda+\beta'}(\CX)
$$

\subsubsection{} {\bf Lemma.} {\em The above operators
provide $\Phi(\CX)$ with a structure of a left $\fu^+$-module.

For $\beta=i$ they coincide with operators $\epsilon_i$ defined above.}

This lemma completes Step 2, proving relations II.12.3 (e)
for generators $\epsilon_i$.

\subsection{} Now we will perform the last Step 3 of the proof,
i.e. prove the relations II.12.3 (c) between operators $\epsilon_i,\theta_j$.
Consider a square

\begin{center}
  \begin{picture}(14,6)


    \put(5,4){\makebox(4,2){$\Phi_{\lambda+i'}(\CX)$}}


    \put(0,2){\makebox(4,2){$\Phi_{\lambda}(\CX)$}}


    \put(10,2){\makebox(4,2){$\Phi_{\lambda-j'+i'}(\CX)$}}


    \put(5,0){\makebox(4,2){$\Phi_{\lambda-j'}(\CX)$}}

   \put(3.5,3.5){\vector(2,1){2}}
   \put(3.5,4){\makebox(1,0.5){$\epsilon_i$}}
   \put(3.5,2.5){\vector(2,-1){2}}
   \put(3,1.5){\makebox(1,0.5){$\theta_j$}}
   \put(8.5,1.5){\vector(2,1){2}}
   \put(10.1,1.5){\makebox(1,0.5){$\epsilon_i$}}
   \put(8.5,4.5){\vector(2,-1){2}}
   \put(9.7,4){\makebox(1,0.5){$\theta_j$}}

  \end{picture}
\end{center}

We have to prove that
\begin{equation}
\label{tri eps delta}
\epsilon_i\theta_j-\zeta^{i\cdot j}\theta_j\epsilon_i=
\delta_{ij}(1-\zeta^{-2\lambda\cdot i'})
\end{equation}
(cf. ~(\ref{tri scal prod})).

As before, we may and will suppose that $\CX\in\FS_c$ and
$\lambda\in X_c$ for some $c$. Choose $\alpha\in \BN[I]$ such that
$\alpha\geq i,\ \alpha\geq j$ and
$\lambda-j'+\alpha'\geq\lambda(\CX)$. The above square may be identified
with the square

\begin{center}
  \begin{picture}(14,6)


    \put(5,4){\makebox(4,2)
     {$\Phi_{\alpha-i-j}(\CX_{\lambda-j'+\alpha'}^{\alpha-i-j})$}}


    \put(0,2){\makebox(4,2)
     {$\Phi_{\alpha-j}(\CX_{\lambda-j'+\alpha'}^{\alpha-j})$}}


    \put(10,2){\makebox(4,2)
     {$\Phi_{\alpha-i}(\CX_{\lambda-j'+\alpha'}^{\alpha-i})$}}


    \put(5,0){\makebox(4,2)
     {$\Phi_{\alpha}(\CX_{\lambda-j'+\alpha'}^{\alpha})$}}

   \put(3.5,3.5){\vector(2,1){2}}
   \put(3.5,4){\makebox(1,0.5){$\epsilon_i$}}
   \put(3.5,2.5){\vector(2,-1){2}}
   \put(3,1.5){\makebox(1,0.5){$\theta_j$}}
   \put(8.5,1.5){\vector(2,1){2}}
   \put(10.1,1.5){\makebox(1,0.5){$\epsilon_i$}}
   \put(8.5,4.5){\vector(2,-1){2}}
   \put(9.7,4){\makebox(1,0.5){$\theta_j$}}

  \end{picture}
\end{center}

\subsection{} Choose an unfolding $\pi:J\lra I$ of $\alpha$; let
$\Sigma=\Sigma_{\pi}$ be its automorphism group, and
$\pi:\BA=\ ^{\pi}\BA\lra\CA_{\alpha}$ denote the corresponding
covering. We will reduce our proof to certain statements about
(vanishing cycles of) sheaves on $\BA$.

Let us introduce a vector space
$$
V=H^0\Phi^+_{\{0\}}(\pi^*\CX^{\alpha}_{\lambda-j'+\alpha'});
$$
the group $\Sigma$ operates on it, and we have
\begin{equation}
\label{tri phi not}
\Phi_{\alpha}(\CX^{\alpha}_{\lambda-j'+\alpha'})\cong V^{\Sigma}
\end{equation}
For each $k\in J$ we have a positive one-dimensional facet
$F_k\subset\BA_{\BR}$
defined as in ~\ref{tri one dim}. Denote
$$
V_k=H^0\Phi_{F_k}^+(\pi^*\CX^{\alpha}_{\lambda-j'+\alpha'});
$$
we have canonically
\begin{equation}
\label{tri phi one}
\Phi_{\alpha-p}(\CX^{\alpha-p}_{\lambda-j'+\alpha'})\cong
(\oplus_{k\in\pi^{-1}(p)} V_k)^{\Sigma}
\end{equation}
for each $p\in I,\ p\leq \alpha$, cf. ~\ref{tri def theta eps}.

\subsubsection{} We have to extend considerations of ~\ref{tri def theta eps}
to two-dimensional facets.
For each pair of different $k,l\in J$ such that $\pi(k)=i,\pi(l)=j$,
let $F_{kl}$
denote a two-dimensional positive facet corresponding
to the map $\varrho:J\lra [0,2]$ sending $k$ to $1$, $l$ --- to $2$ and
all other elements --- to zero (cf. ~\ref{tri facets}).
Set
$$
V_{kl}=H^0\Phi_{F_{kl}}^+(\pi^*\CX^{\alpha}_{\lambda-j'+\alpha'}).
$$
Again, due to equivariance of our sheaf, the group $\Sigma$ operates
on $\oplus\ V_{kl}$ in such a way that
$$
\sigma(V_{kl})=V_{\sigma(k)\sigma(l)}.
$$
Let
$$
\pi_{kl}:J-\{k,l\}\lra I
$$
be the restriction of $\pi$. It is an unfolding of $\alpha-i-j$;
let $\Sigma_{kl}$ denote its automorphism group.
Pick $d_1,d_2$ such that $0<d_2<1<d_1<2$. The subspace
\begin{equation}
\label{tri perp two}
F^{\perp}_{kl}\subset\BA^{\{l\},\{k\},J-\{k,l\}}(d_1,d_2)
\end{equation}
consisting of all points $(t_j)$ such that $t_k=1$ and $t_l=2$ is a transversal
slice to $f_{kl}$. Consequently, factorization axiom implies canonical
isomorphism
$$
\Phi^+_{F_{kl}}(\pi^*\CX^{\alpha}_{\lambda-j'+\alpha'})\cong
\Phi^+_{\{0\}}(\pi^*_{kl}\CX^{\alpha-i-j}_{\lambda-j'+\alpha'})\otimes
(\CI^i_{\lambda+i'})_{\{1\}}\otimes(\CI^j_{\lambda})_{\{2\}}=
\Phi^+_{\{0\}}(\pi^*_{kl}\CX^{\alpha-i-j}_{\lambda-j'+\alpha'})
$$

{\em Symmetry.} Interchanging $k$ and $l$, we get isomorphisms
\begin{equation}
\label{tri transpos}
t:V_{kl}\iso V_{lk}
\end{equation}

Passing to $\Sigma$-invariants, we get isomorphisms
\begin{eqnarray}
\label{tri phi two}
\Phi_{\alpha-i-j}(\CX^{\alpha-i-j}_{\lambda-j'+\alpha'})=
\Phi^+_{\{0\}}(\pi^*_{kl}\CX^{\alpha-i-j}_{\lambda-j'+\alpha'})^{\Sigma_{kl}}
\cong\\ \nonumber
\cong
\Phi^+_{F_{kl}}(\pi^*\CX^{\alpha}_{\lambda-j'+\alpha'})^{\Sigma_{kl}}
\cong
(\oplus_{k\in\pi^{-1}(i),l\in\pi^{-1}(j),k\neq l}\ V_{kl})^{\Sigma}
\end{eqnarray}
(cf. ~(\ref{tri phi li})).

\subsection{}
\label{tri can var}
The canonical and variation maps induce linear operators
$$
v_k:V\overset{\lra}{\lla} V_k: u^k
$$
and
$$
v^k_{lk}:V_k\overset{\lra}{\lla} V_{lk}: u^{lk}_k
$$
which are $\Sigma$-equivariant in the obvious sense. Taking their sum,
we get operators
$$
V\overset{\lra}{\lla}\oplus_{k\in\pi^{-1}(p)} V_k
\overset{\lra}{\lla} \oplus_{l\in\pi^{-1}(q),k\in\pi^{-1}(p)} V_{lk}
$$
which induce, after passing to $\Sigma$-invariants,
operators $\epsilon_p,\epsilon_q$ and $\theta_p,\theta_q$.

Our square takes a form

\begin{center}
  \begin{picture}(14,6)


    \put(0,4){\makebox(4,2)
     {$(\oplus_{k\in\pi^{-1}(j),l\in\pi^{-1}(i),k\neq l} V_{lk})^{\Sigma}$}}


\put(10,4){\makebox(4,2)
     {$(\oplus_{k\in\pi^{-1}(j),l\in\pi^{-1}(i),k\neq l} V_{kl})^{\Sigma}$}}


    \put(0,2){\makebox(4,2)
     {$(\oplus_{k\in\pi^{-1}(j)} V_k)^{\Sigma}$}}


    \put(10,2){\makebox(4,2)
     {$(\oplus_{l\in\pi^{-1}(i)} V_l)^{\Sigma}$}}


    \put(5,0){\makebox(4,2)
     {$V^{\Sigma}$}}


   \put(6.6,4.8){$\overset{t}{\iso}$}

   \put(2,3.5){\vector(0,1){1}}
   \put(.8,3.7){\makebox(1,0.5){$\sum v^k_{lk}$}}
   \put(3.5,2.5){\vector(2,-1){2}}
   \put(3,1.5){\makebox(1,0.5){$\sum u^k$}}
   \put(8.5,1.5){\vector(2,1){2}}
   \put(10.1,1.5){\makebox(1,0.5){$\sum v_l$}}
   \put(12,4.5){\vector(0,-1){1}}
   \put(12.5,3.7){\makebox(1,0.5){$\sum u_l^{kl}$}}

  \end{picture}
\end{center}

Now we will formulate two relations between $u$ and $v$ which imply
the necessary relations between $\epsilon$ and $\theta$.

\subsection{Lemma}
\label{tri i eq j} {\em Suppose that $i=j$ and $\pi(k)=i$. Consider operators
$$
u^k:V_k\overset{\lra}{\lla} V:v_k.
$$
The composition $v_ku^k$ is equal to the multiplication by
$1-\zeta^{-2\lambda\cdot i'}$.}

{\bf Proof.} Consider the transversal slice $F_k^{\perp}(d)$ to the face
$F_k$ as in ~\ref{tri def theta eps}. It follows from the definition
of the canonical and variation maps, II.7.10, that composition
$v_ku^k$ is equal to $1-T^{-1}$ where $T$ is the monodromy acquired by
of $\Phi^+_{\{0\}}(\pi^*_k\CX^{\alpha-i}_{\lambda-i'+\alpha'})$ when
the point $t_k$ moves around the disc of radius $d$ where all other
points are living. By factorization, $T=\zeta^{2\lambda\cdot i'}$. $\Box$

\subsection{Lemma.}
\label{tri i noteq j} {\em For $k\in\pi^{-1}(j),\ l\in\pi^{-1}(i),\
k\neq l,$
consider the pentagon

\begin{center}
  \begin{picture}(14,6)


    \put(0,4){\makebox(4,2)
     {$V_{lk}$}}


\put(10,4){\makebox(4,2)
     {$V_{kl}$}}


    \put(0,2){\makebox(4,2)
     {$V_k$}}


    \put(10,2){\makebox(4,2)
     {$V_l$}}


    \put(5,0){\makebox(4,2)
     {$V$}}


   \put(6.6,4.8){$\overset{t}{\iso}$}

   \put(2,3.5){\vector(0,1){1}}
   \put(.8,3.7){\makebox(1,0.5){$v^k_{lk}$}}
   \put(3.5,2.5){\vector(2,-1){2}}
   \put(3,1.5){\makebox(1,0.5){$u^k$}}
   \put(8.5,1.5){\vector(2,1){2}}
   \put(10.1,1.5){\makebox(1,0.5){$v_l$}}
   \put(12,4.5){\vector(0,-1){1}}
   \put(12.5,3.7){\makebox(1,0.5){$u_l^{kl}$}}

  \end{picture}
\end{center}
We have
$$
v_lu^k=\zeta^{i\cdot j}u_k^{kl}\circ t\circ v_{lk}^k.
$$}

This lemma is a consequence of the following more general statement.

\subsection{} Let $\CK\in\CCD(\CA;\CS)$ be arbitrary. We have naturally
$$
\Phi_{F_{kl}}^+(\CK)\cong\Phi_{\{0\}}^+(\CK|_{F^{\perp}_{kl}(d_1,d_2)}[-2])
$$
Let
$$
t:\Phi_{F_{kl}}^+(\CK)\iso\Phi_{F_{lk}}^+(\CK)
$$
be the monodromy isomorphism induced by the travel of the point $t_l$ in the
upper half plane to the position to the left of $t_k$
(outside the disk of radius $d_1$).

\subsubsection{} {\bf Lemma.} {\em The composition
$$
v_{F_l}^{\{0\}}\circ u_{\{0\}}^{F_k}:
\Phi^+_{F_k}(\CK)\lra \Phi^+_{F_l}(\CK)
$$
is equal to $u_{F_l}^{F_{kl}}\circ t\circ v_{F_{lk}}^{F_k}$.}

\subsection{} It remains to note that due to ~\ref{tri can var}
the desired relation ~(\ref{tri eps delta}) is a formal consequence
of lemmas ~\ref{tri i eq j} and ~\ref{tri i noteq j}. This completes the proof
of Theorem ~\ref{tri relations}. $\Box$

\subsection{} Taking composition of $\tPhi$ with an inverse to the equivalence
$Q$, II (143), we get a functor
\begin{equation}
\label{tri final phi}
\Phi:\FS\lra\CC
\end{equation}

\section{Main properties of $\Phi$}

{\em TENSOR PRODUCTS}

\subsection{Theorem}
\label{tri tens fun} {\em $\Phi$ is a tensor functor, i.e. we have
natural isomorphisms
$$
\Phi(\CX\dotimes\CY)\iso\Phi(\CX)\otimes\Phi(\CY)
$$
satisfying all necessary compatibilities.}

{\bf Proof} follows from the Additivity theorem, II.9.3. $\Box$

\vspace{.5cm}
{\em DUALITY}

\subsection{}
\label{tri dual c} Let $\tCC_{\zeta^{-1}}$ denote the category $\tCC$ with
the value of parameter $\zeta$ changed to $\zeta^{-1}$. The notations
$\FS_{\zeta^{-1}}$, etc. will have the similar meaning.

Let us define a functor
\begin{equation}
\label{tri d tc}
D:\tCC^{opp}\lra\tCC_{\zeta^{-1}}
\end{equation}
as follows. By definition, for $V=\oplus V_{\lambda}\in Ob\ \tCC^{opp}=
Ob\ \tCC$,
$$
(DV)_{\lambda}=V^*_{\lambda},
$$
and operators
$$
\theta_{i,DV}:(DV)_{\lambda}\overset{\lra}{\lla}(DV)_{\lambda-i'}:
\epsilon_{i,DV}
$$
are defined as
$$
\epsilon_{i,DV}=\theta^*_{i,V};\
\theta_{i,DV}=-\zeta^{2\lambda\cdot i'}\epsilon^*_{i,V}.
$$
On morphisms $D$ is defined in the obvious way.
One checks directly that $D$ is well defined,
respects tensor structures, and is an equivalence.

\subsection{}
\label{tri dual fs} Let us define a functor
\begin{equation}
\label{tri d fs}
D:\FS^{opp}\lra\FS_{\zeta^{-1}}
\end{equation}
as follows. For $\CX\in Ob\ \FS^{opp}=Ob\ \FS$ we set
$\lambda(D\CX)=\lambda(\CX)$; $(D\CX)^{\alpha}=D(\CX^{\alpha})$ where
$D$ in the right hand side is Verdier duality. Factorization
isomorphisms for $D\CX$ are induced in the obvious way
from factorization isomorphisms for $\CX$. The value of $D$
on morphisms is defined in the obvious way.

$D$ is a tensor equivalence.

\subsection{Theorem.} {\em Functor $\tPhi$ commutes with $D$.}

{\bf Proof.} Our claim is a consequence of the following
topological remarks.

\subsection{}
\label{tri affine} Consider a standard affine space $\BA=\BA^J$ with a
principal stratification $\CS$ as in II.7. Let $\CK\in\CCD(\BA;\CS)$,
let $F_j$ be the one-dimensional facet corresponding to an element
$j\in J$ as in ~\ref{tri one dim}. Consider a transversal slice
$F_j^{\perp}(d)$ as in ~\ref{tri def theta eps}. We have canonically
$$
\Phi_{F_j}^+(\CK)\cong\Phi_{\{0\}}(\CK|_{F_j^{\perp}(d)}[-1]);
$$
when the point $t_j$ moves counterclockwise around the disk
of radius $d$, $\Phi_j^+(\CK)$ acquires monodromy
$$
T_j:\Phi_j^+(\CK)\iso\Phi_j^+(\CK).
$$

\subsubsection{} {\bf Lemma.} {\em Consider canonical and variation
maps
$$
v_{D\CK}:\Phi^+_{\{0\}}(D\CK)\overset{\lra}{\lla}
\Phi^+_{F_j}(D\CK):u_{D\CK}
$$
Let us identify $\Phi^+_{\{0\}}(D\CK),\Phi^+_{F_j}(D\CK)$ with
$\Phi^+_{\{0\}}(\CK)^*$ and $\Phi^+_{F_j}(\CK)^*$ respectively.
Then the maps become
$$
v_{D\CK}=u^*_{\CK};\ u_{D\CK}=-v^*_{\CK}\circ T_j^*
$$}

The theorem follows from this lemma. $\Box$

\vspace{.5cm}
{\em STANDARD OBJECTS}

\subsection{Theorem}
\label{tri irred} {\em We have naturally
$$
\Phi(\CL(\Lambda))\cong L(\Lambda)
$$
for all $\Lambda\in X$.}

Combining this with Theorem ~\ref{tri all irreds}, we get

\subsection{Theorem} {\em $\Phi$ induces bijection between sets
of isomorphism classes of irreducibles.} $\Box$

\subsection{Verma modules} Let $\fu^{\geq 0}\subset\fu$ denote the subalgebra
generated by $\epsilon_i, K_i,K^{-1},\ i\in I$.
For $\Lambda\in X$, let $\chi_{\Lambda}$ denote a one-dimensional
representation of $\fu^{\geq 0}$ generated by a vector $v_{\Lambda}$,
with the action
$$
\epsilon_iv_{\Lambda}=0,\
K_iv_{\Lambda}=\zeta^{\langle\Lambda,i\rangle}v_{\Lambda}.
$$
Let $M(\Lambda)$ denote the induced $\fu$-module
$$
M(\Lambda)=\fu\otimes_{\fu^{\geq 0}}\chi_{\Lambda}.
$$
Equipped with an obvious $X$-grading, $M(\Lambda)$ becomes
an object of $\tCC$. We will also use the same notation for the
corresponding object of $\CC$.

\subsection{Theorem} {\em The factorizable sheaves $\CM(\Lambda)$ are
finite. We have naturally
$$
\Phi(\CM(\Lambda))\cong M(\Lambda).
$$}

{\bf Proof} is given in the next two subsections.

\subsection{} Let us consider the space $\BA$ as in ~\ref{tri affine}.
Suppose that $\CK\in\CCD(\BA;\CS)$ has the form
$\CK=j_!j^*\CK$ where
$$
j:\BA-\bigcup_{j\in J}\{t_j=0\}\hra\BA.
$$
Let $F_{\Delta}$ be the positive facet whose closure is the main diagonal.

\subsubsection{} {\bf Lemma.} {\em The canonical map
$$
u:\Phi^+_{F_{\Delta}}(\CK)\lra\Phi^+_{\{0\}}(\CK)
$$
is an isomorphism.}

{\bf Proof.} Pick $j_0\in J$, and consider a subspace $Y=\{t_{j_0}=0\}
\cup\{t_{j_0}=1\}\subset\BA$. Set $\CK'=k_!k^*\CK$ where
$$
k:\BA-\bigcup_{j\in J}\{t_j=0\}-\bigcup_{j\in J}\{t_j=1\}\hra\BA.
$$
We have
$$
\Phi^+_{\{0\}}(\CK)=R\Gamma(\BA;\CK')
$$
On the other hand by homotopy we have
$$
\Phi^+_{F_{\Delta}}(\CK)\cong
R\Gamma(\{t_{j_0}=c\};\CK')[-1]
$$
where $c$ is any real between $0$ and $1$. Let us compute
$R\Gamma(\BA;\CK')$ using the Leray spectral sequence of
a projection
$$
p:\BA\lra\BA^1,\ (t_j)\mapsto t_{j_0}.
$$
The complex $p_*\CK'$ is equal to zero at the points $\{0\}$ and $\{1\}$,
and is constant with the fiber $R\Gamma(\{t_{j_0}=c\};\CK')$ over $c$.
It follows that
$$
R\Gamma(\BA;\CK')\cong R\Gamma(\BA^1;p_*\CK')
\cong R\Gamma(\{t_{j_0}=c\};\CK')[-1],
$$
and the inverse to this isomorphism may be identified with $u$. $\Box$

\subsection{} Suppose we have $\alpha\in\BN[I]$, let
$\pi:J\lra I;\ \pi:\BA\lra\CA_{\alpha}$ be the corresponding unfolding.
Let us apply the previous lemma to $\CK=\pi^*\CM(\Lambda)^{\alpha}$.
Note that after passing to $\Sigma_{\pi}$-invariants,
the map $u$ becomes
$$
\fu_{\alpha}^-\lra\tPhi_{\alpha}(\CM(\Lambda))
$$
by Theorem II.6.16. This identifies homogeneous components
of $\tPhi_{\alpha}(\CM(\Lambda))$ with the components of the Verma module.
After that, the action of $\epsilon_i$ and $\theta_i$ is identified with
the action of $\fu$ on it. The theorem is proven. $\Box$

\newpage
\begin{center}
{\bf Chapter 4. Equivalence}
\end{center}
\vspace{.8cm}

\section{Truncation functors}

\subsection{} Recall the notations of ~\ref{tri stabilization sec}.
We fix a coset $X_c\subset X$, and we denote the subcategory
$\FS_{\leq \lambda; c}\subset\FS$ by $\FS_{\leq\lambda}$ for simplicity until
further notice.

Given $\lambda\in X$, we will denote by $\CC_{\leq \lambda} \subset\CC$ the
full subcategory of all $\fu$-modules $V$ such that $V_{\mu}\neq 0$
implies $\mu\leq \lambda$. We denote by $q_{\lambda}$ the
embedding functor $\CC_{\lambda}\hra\CC$. Obviously,
$\Phi(\FS_{\leq \lambda})\subset\CC_{\lambda}$.

In this section we will construct functors
$$
\sigma^!_{\lambda},\sigma^*_{\lambda}:\FS\lra\FS_{\leq\lambda}
$$
and
$$
q^!_{\lambda},q^*_{\lambda}:\CC\lra\CC_{\leq\lambda},
$$
such that $\sigma^!_{\lambda}$ (resp., $\sigma^*_{\lambda}$)
is right (resp., left) adjoint to $\sigma_{\lambda}$ and
$q^!_{\lambda}$ (resp., $q^*_{\lambda}$)
is right (resp., left) adjoint to $q_{\lambda}$.

\subsection{} First we describe $\sigma^*_{\lambda},\sigma^!_{\lambda}$.
Given a  factorizable sheaf $\CX=\{\CX^{\alpha}\}$ with $\lambda(\CX)=
\mu\geq\lambda$ we define FS's $\CY:=\sigma^*_{\lambda}\CX$ and
$\CZ:=\sigma^!_{\lambda}\CX$ as follows.

We set $\lambda(\CY)=\lambda(\CZ)=\lambda$. For $\alpha\in\BN[I]$ we set
$$
\CY^{\alpha}=L^0\sigma^*\CX^{\alpha+\mu-\lambda}
$$
if $\alpha+\mu-\lambda\in\BN[I]$ and $0$ otherwise, and
$$
\CZ^{\alpha}=R^0\sigma^!\CX^{\alpha+\mu-\lambda}
$$
if $\alpha+\mu-\lambda\in\BN[I]$ and $0$ otherwise. Here $\sigma$ denotes
the canonical closed embedding (cf. ~\ref{tri close emb})
$$
\sigma:\CA^{\alpha}_{\lambda}\hra\CA^{\alpha+\mu-\lambda}_{\mu},
$$
and $L^0\sigma^*$ (resp., $R^0\sigma^!$) denotes the zeroth perverse
cohomology of $\sigma^*$ (resp., of $\sigma^!$).

The factorization isomorphisms for $\CY$ and $\CZ$ are induced from those
for $\CX$; associativity is obvious.

\subsection{Lemma.}
\label{tri adj fs} {\em Let $\CM\in\FS_{\leq \lambda}$,
$\CX\in\FS$. Then
$$
\Hom_{\FS}(\CX,\CM)=\Hom_{\FS_{\leq\lambda}}(\sigma^*_{\lambda}\CX,\CM)
$$
and
$$
\Hom_{\FS}(\CM,\CX)=\Hom_{\FS_{\leq\lambda}}(\CM,\sigma^!_{\lambda}\CX)
$$}

{\bf Proof.} Let $\mu=\lambda(\CX)$. We have
$$
\Hom_{\CA^{\alpha}_{\mu}}(\CX^{\alpha}_{\mu},\CM^{\alpha}_{\mu})=
\Hom_{\CA^{\alpha-\mu+\lambda}_{\mu}}
(\sigma^*_{\lambda}\CX_{\lambda}^{\alpha-\mu+\lambda},
\CM^{\alpha-\mu+\lambda}_{\lambda})
$$
by the usual adjointness. Passing to projective limit in $\alpha$,
we get the desired result for $\sigma^*$. The proof for $\sigma^!$ is similar.
$\Box$

\subsection{} Given $\lambda\leq\mu\in X_c$, we denote by
$\sigma_{\lambda\leq\mu}$ the embedding of the full subcategory
$$
\sigma_{\lambda\leq\mu}:\FS_{\leq\lambda}\hra\FS_{\leq\mu}.
$$
Obviously, the functor
$$
\sigma^*_{\lambda\leq\mu}:=\sigma^*_{\lambda}\circ\sigma_{\mu}:
\FS_{\leq\mu}\lra\FS_{\leq\lambda}
$$
is left adjoint to $\sigma_{\lambda\leq\mu}$.
Similarly, $\sigma^!_{\lambda\leq\mu}:=\sigma^!_{\lambda}\circ\sigma_{\mu}$
is the right adjoint to $\sigma_{\lambda\leq\mu}$.

For $\lambda\leq\mu\leq\nu$ we have obvious transitivities
$$
\sigma^*_{\lambda\leq\mu}\sigma^*_{\mu\leq\nu}=\sigma^*_{\lambda\leq\nu};\
\sigma^!_{\lambda\leq\mu}\sigma^!_{\mu\leq\nu}=\sigma^!_{\lambda\leq\nu}.
$$

\subsection{} For each $\alpha\in\BN[I]$ and $i\in I$ such that
$\alpha\geq i$ let $j_{\nu-i'\leq\nu}^{\alpha}$ denote the open
embedding
$$
j_{\nu-i'\leq\nu}^{\alpha}:\CA^{\alpha}_{\nu}-
\sigma(\CA^{\alpha-i}_{\nu-i'})\hra\CA^{\alpha}_{\nu}.
$$
Note that the complement of this open subspace is a divisor, so the
corresponding extension by zero and by $*$ functors are $t$-exact,
cf. ~\cite{bbd}, 4.1.10 (i).
Let us define functors
$$
j_{\nu-i'\leq\nu!},j_{\nu-i'\leq\nu*}:\FS_{\leq\nu}\lra\FS_{\leq\nu}
$$
as follows. For a factorizable sheaf
$\CX=\{\CX^{\alpha}_{\nu}\}\in\FS_{\leq\nu}$ we set
$$
(j_{\nu-i'\leq\nu!}\CX)^{\alpha}_{\nu}=
j_{\nu-i'\leq\nu!}^{\alpha}j_{\nu-i'\leq\nu}^{\alpha*}\CX^{\alpha}_{\nu}
$$
and
$$
(j_{\nu-i'\leq\nu*}\CX)^{\alpha}_{\nu}=
j_{\nu-i'\leq\nu*}^{\alpha}j_{\nu-i'\leq\nu}^{\alpha*}\CX^{\alpha}_{\nu},
$$
the factorization isomorphisms being induced from those for $\CX$.

\subsection{Lemma.}
\label{tri ex seq} {\em We have natural in $\CX\in\FS_{\leq\nu}$ exact
sequences
$$
j_{\nu-i',\nu!}\CX\lra\CX\overset{a}{\lra}\sigma^*_{\nu-i',\nu}\CX\lra 0
$$
and
$$
0\lra \sigma^!_{\nu-i',\nu}\CX\overset{a'}{\lra}\CX\lra j_{\nu-i',\nu*}\CX
$$
where the maps $a$ and $a'$ are the adjunction morphisms.}

{\bf Proof.} Evidently follows from the same claim at each finite level,
which is ~\cite{bbd}, 4.1.10 (ii).
$\Box$

\subsection{} Recall (see ~\ref{tri dual fs}) that we have the Verdier
duality functor
$$
D:\FS^{opp}\lra\FS_{\zeta^{-1}}.
$$
By definition, $D(\FS_{\leq\lambda}^{opp})\subset\FS_{\zeta^{-1};\leq\lambda}$
for all $\lambda$.

We have functorial isomorphisms
$$
D\circ\sigma^*_{\lambda}\cong\sigma^!_{\lambda}\circ D;\
D\circ\sigma^*_{\nu-i',\nu}\cong\sigma^!_{\nu-i',\nu}\circ D
$$
and
$$
D\circ j_{\nu-i',\nu*}\cong j_{\nu-i',\nu!}\circ D
$$
After applying $D$, one of the exact sequences in ~\ref{tri ex seq} becomes
another one.

\subsection{} Let us turn to the category $\CC$. Below we will identify
$\CC$ with $\tCC$ using the equivalence $Q$, cf. II.12.5. In other
words, we will regard objects of $\CC$ as $\fu$-modules.

For $\lambda\in X_c$
functors $q^!_{\lambda}$ and $q^*_{\lambda}$ are defined as follows.
For $V\in\CC$, $q^!_{\lambda}V$ (resp., $q^*_{\lambda}V$) is the maximal
subobject (resp., quotient) of $V$ belonging to the subcategory
$\CC_{\lambda}$.

For $\lambda\leq\mu\in X_c$ let $q_{\lambda\leq\mu}$ denotes an
embedding of a full subcategory
$$
q_{\lambda\leq\mu}:\CC_{\leq\lambda}\hra\CC_{\leq\mu}
$$
Define $q^!_{\lambda\leq\mu}:=q^!_{\lambda}\circ q_{\mu};\
q^*_{\lambda\leq\mu}:=q^!_{\lambda}\circ q_{\mu}$. Obviously,
the first functor is right adjoint, and the second one is left adjoint
to $q_{\lambda\leq\mu}$. They have an obvious transitivity property.

\subsection{} Recall that in ~\ref{tri dual c} the weight preserving
duality equivalence
$$
D:\CC^{opp}\lra\CC_{\zeta^{-1}}
$$
is defined.
By definition, $D(\CC_{\leq\lambda}^{opp})\subset\CC_{\zeta^{-1};\leq\lambda}$
for all $\lambda$.

We have functorial isomorphisms
$$
D\circ q^*_{\lambda}\cong q^!_{\lambda}\circ D;\
D\circ q^*_{\nu-i',\nu}\cong q^!_{\nu-i',\nu}\circ D.
$$

\subsection{} Given $i\in I$, let us introduce a "Levi" subalgebra
$\fl_i\subset\fu$ generated by $\theta_j,\ \epsilon_j,\ j\neq i$, and
$K_i^{\pm}$. Let $\fp_i\subset\fu$ denote the "parabolic"
subalgebra generated by $\fl_i$ and $\epsilon_i$.

The subalgebra $\fl_i$ projects isomorphically to $\fp_i/(\epsilon_i)$
where $(\epsilon_i)$ is a two-sided ideal generated by $\epsilon_i$.
Given an $\fl_i$-module $V$, we can consider it as a $\fp_i$-module
by restriction of scalars for the projection $\fp_i\lra\fp_i/(\epsilon_i)\cong
\fl_i$, and form the induced $\fu$-module $\Ind_{\fp_i}^{\fu}V$ ---
"generalized Verma".

\subsection{} Given an $\fu$-module $V\in\CC_{\leq\nu}$, let us consider a
subspace
$$
_iV=\oplus_{\alpha\in\BN[I-\{i\}]}V_{\nu-\alpha'}\subset V.
$$
It is an $X$-graded $\fp_i$-submodule of $V$. Consequently, we have
a canonical element
$$
\pi\in\Hom_{\CC}(\Ind_{\pi_i}^{\fu}\ _iV, V)=
\Hom_{\fp_i}(_iV,V)
$$
corresponding to the embedding $_iV\hra V$.

We will also consider the dual functor
$$
V\mapsto D^{-1}\Ind_{\fp_i}^{\fu}\ _i(DV).
$$
By duality, we have a natural morphism
in $\CC$, $V\lra D^{-1}\Ind_{\fp_i}^{\fu}\ _i(DV)$.

\subsection{Lemma.}
\label{tri ex seq c} {\em We have natural in $V\in\CC_{\leq\nu}$ exact
sequences
$$
\Ind_{\fp_i}^{\fu}\ _iV\overset{\pi}{\lra} V\lra q^*_{\nu-i'\leq\nu} V\lra 0
$$
and
$$
0\lra q^!_{\nu-i'\leq\nu} V\lra V\lra D^{-1}\Ind_{\fp_i}^{\fu}\ _i(DV).
$$
where the arrows $V\lra q^*_{\nu-i'\leq\nu} V$ and
$q^!_{\nu-i'\leq\nu} V\lra V$ are adjunction morphisms.}

{\bf Proof.} Let us show the exactness of the first sequence.
By definition, $q^*_{\nu-i'\leq\nu} V$ is the maximal
quotient of $V$ lying in the subcategory $\CC_{\nu-i'\leq\nu}\subset\CC$.
Obviously, $\Coker\ \pi\in\CC_{\nu-i',\nu}$. It remains to show that for
any morphism $h:V\lra W$ with $W\in\CC_{\nu-i'}$, the composition
$h\circ\pi:\Ind_{\fp_i}^{\fu}\ _iV\lra W$ is zero. But
$\Hom_{\CC}(\Ind_{\fp_i}^{\fu}\ _iV,W)=\Hom_{\fp_i}(_iV,W)=0$ by
weight reasons.

The second exact sequence is the dual to the first one. $\Box$

\subsection{Lemma.}
\label{tri phi j} {\em We have natural in $\CX\in\FS_{\nu}$ isomorphisms
$$
\Phi j_{\nu-i'\leq\nu!}\CX\iso\Ind_{\fp_i}^{\fu}\ _i(\Phi\CX)
$$
and
$$
\Phi j_{\nu-i'\leq\nu*}\CX\iso\ D^{-1}\Ind_{\fp_i}^{\fu}\ _i(D\Phi\CX)
$$
such that the diagram

\begin{center}
  \begin{picture}(25,6)


\put(1,4){$\Phi j_{\nu-i',\nu!}\CX$}
\put(7.8,4){$\Phi\CX$}
\put(13,4){$\Phi j_{\nu-i',\nu*}\CX$}


\put(1,1){$\Ind_{\fp_i}^{\fu}\ _i(\Phi\CX)$}
\put(7.8,1){$\Phi\CX$}
\put(12,1){$D^{-1}\Ind_{\fp_i}^{\fu}\ _i(D\Phi\CX)$}


\put(2,3.6){\vector(0,-1){2}}
\put(8.1,3.6){\line(0,-1){2}}
\put(8.2,3.6){\line(0,-1){2}}
\put(14,3.6){\vector(0,-1){2}}


\put(4,4.1){$\vector(1,0){3}$}
\put(4,1.1){$\vector(1,0){3}$}

\put(9.5,4.1){$\vector(1,0){3}$}
\put(9.5,1.1){$\vector(1,0){2}$}

\end{picture}
\end{center}

commutes.}

\subsection{Lemma.}
\label{tri phi sigma} {\em Let $\lambda\in X_c$. We have natural in
$\CX\in\FS$ isomorphisms
$$
\Phi\sigma^*_{\lambda}\CX\cong q_{\lambda}^*\Phi\CX
$$
and
$$
\Phi\sigma^!_{\lambda}\CX\cong q_{\lambda}^!\Phi\CX
$$}

{\bf Proof} follows at once from lemmas ~\ref{tri phi j},
{}~\ref{tri ex seq} and ~\ref{tri ex seq c}. $\Box$

\section{Rigidity}

\subsection{Lemma.} {\em Let $\CX\in\FS_{\leq 0}$. Then the natural map
$$
a:\Hom_{\FS}(\CL(0),\CX)\lra\Hom_{\CC}(L(0),\Phi(\CX))
$$
is an isomorphism.}

{\bf Proof.} We know already that $a$ is injective, so we have
to prove its surjectivity. Let $\CK(0)$ (resp.,
$K(0)$) denote the kernel of the projection
$\CM(0)\lra \CL(0)$ (resp., $M(0)\lra L(0)$). Consider a diagram with
exact rows:

\begin{center}
  \begin{picture}(25,6)


\put(3,4){$0$}
\put(6,4){$\Hom(\CL(0),\CX)$}
\put(12,4){$\Hom(\CM(0),\CX)$}
\put(18,4){$\Hom(\CK(0),\CX)$}


\put(3,1){$0$}
\put(6,1){$\Hom(L(0),\Phi(\CX))$}
\put(12,1){$\Hom(M(0),\Phi(\CX))$}
\put(18,1){$\Hom(K(0),\Phi(\CX))$}


\put(8,3.6){\vector(0,-1){2}}
\put(7.5,2.6){$a$}
\put(14,3.6){\vector(0,-1){2}}
\put(13.5,2.6){$b$}
\put(20,3.6){\vector(0,-1){2}}


\put(3.5,4.1){$\vector(1,0){2}$}
\put(3.5,1.1){$\vector(1,0){2}$}

\put(9.5,4.1){$\vector(1,0){2.1}$}
\put(10.1,1.1){$\vector(1,0){1.5}$}

\put(16,4.1){$\vector(1,0){1.7}$}
\put(16.4,1.1){$\vector(1,0){1.2}$}

  \end{picture}
\end{center}

All vertical rows are injective. On the other hand, $\Hom(M(0),\Phi(\CX))=
\Phi(\CX)_0$. The last space is isomorphic to a generic stalk of
$\CX^{\alpha}_0$ for each $\alpha\in\BN[I]$, which in turn is isomorphic
to $\Hom_{\FS}(\CM(0),\CX)$ by the universal property of the shriek extension.
Consequently, $b$ is isomorphism by the equality of dimensions. By diagram
chase, we conclude that $a$ is isomorphism. $\Box$

\subsection{Lemma.}
\label{tri iso irr} {\em For every $\CX\in\FS$ the natural maps
$$
\Hom_{\FS}(\CL(0),\CX)\lra\Hom_{\CC}(L(0),\Phi(\CX))
$$
and
$$
\Hom_{\FS}(\CX,\CL(0))\lra\Hom_{\CC}(\Phi(\CX),L(0))
$$
are isomorphisms.}

{\bf Proof.} We have
$$
\Hom_{\FS}(\CL(0),\CX)=\Hom_{\FS}(\CL(0),\sigma^!_0\CX)
$$
(by lemma ~\ref{tri adj fs})
$$
=\Hom_{\CC}(L(0),\Phi(\sigma^!_0\CX))
$$
(by the previous lemma)
$$
=\Hom_{\CC}(L(0),q_0^!\Phi(\CX))
$$
(by lemma ~\ref{tri phi sigma})
$$
=\Hom_{\CC}(L(0),\Phi(\CX)).
$$
This proves the first isomorphism. The second one follows
by duality. $\Box$

\subsection{} Recall (see ~\cite{kl}IV, Def. A.5) that an object $X$
of a tensor category is called {\em rigid} if there exists
another object $X^*$ together with morphisms
$$
i_X:\One\lra X\otimes X^*
$$
and
$$
e_X: X^*\otimes X\lra \One
$$
such that the compositions
$$
X=\One\otimes X\overset{i_X\otimes \id}{\lra}
X\otimes X^*\otimes X\overset{\id\otimes e_X}{\lra} X
$$
and
$$
X^*=X^*\otimes \One\overset{\id\otimes i_X}{\lra}
X^*\otimes X\otimes X^*\overset{e_X\otimes\id}{\lra} X^*
$$
are equal to $\id_X$ and $\id_{X^*}$ respectively. A tensor
category is called rigid if all its objects are rigid.

\subsection{Theorem.} {\em All irreducible objects $\CL(\lambda),\
\lambda\in X$, are rigid in $\FS$.}

{\bf Proof.}
We know (cf. ~\cite{ajs}, 7.3) that $\CC$ is rigid. Moreover, there exists an
involution $\lambda\mapsto\blambda$ on $X$ such that $L(\lambda)^*=
L(\blambda)$. Let us define $\CL(\lambda)^*:=\CL(\blambda)$;
$i_{\CL(\lambda)}$ corresponds to $i_{L(\Lambda)}$ under identification
$$
\Hom_{\FS}(\CL(0),\CL(\lambda)\otimes\CL(\blambda))=
\Hom_{\CC}(L(0),L(\lambda)\otimes L(\blambda))
$$
and $e_{\CL(\lambda)}$ corresponds to
$e_{L(\lambda)}$ under identification
$$
\Hom_{\FS}(\CL(\blambda)\otimes\CL(\lambda),\CL(0))=
\Hom_{\CC}(L(\blambda)\otimes L(\lambda),L(0)),
$$
cf. ~\ref{tri iso irr}. $\Box$

\section{Steinberg sheaf}
\label{tri steinb sh}

In this section we assume that $l$ is a positive number prime to $2,3$ and
that $\zeta'$ is a primitive $(l\cdot\det A)$-th root of $1$
(recall that $\zeta=(\zeta')^{\det A})$.

We fix a weight $\lambda_0\in X$ such that $\langle i,\lambda_0\rangle=
-1$(mod $l$) for any $i\in I$. Our goal in this section is the proof of the
following

\subsection{}
\label{tri steinberg}
\begin{thm}{}
The FFS $\CL(\lambda_0)$ is a projective object of the category $\FS$.
\end{thm}

{\bf Proof.} We have to check that $\Ext^1(\CL(\lambda_0),\CX)=0$ for any
FFS $\CX$. By induction on the length of $\CX$ it is enough to prove that
$\Ext^1(\CL(\lambda_0),\CL)=0$ for any simple FFS $\CL$.

\subsection{}
\label{tri principle}
To prove vanishing of $\Ext^1$ in $\FS$ we will use the following
principle. Suppose $\Ext^1(\CX,\CY)\not=0$, and let
$$
0\lra\CY\lra\CZ\lra\CX\lra 0
$$
be the corresponding nonsplit extension.
Let us choose a weight $\lambda$ which is bigger than
of $\lambda(\CX),\lambda(\CY),\lambda(\CZ)$. Then for any
$\alpha\in\BN[I]$ the sequence
$$
0\lra\CY^\alpha_\lambda\lra\CZ^\alpha_\lambda\lra\CX^\alpha_\lambda\lra 0
$$
is also exact, and for $\alpha\gg 0$ it is also nonsplit (see lemma
{}~\ref{tri stabilization}). That is, for $\alpha\gg 0$ we have
$\Ext^1(\CX^\alpha_\lambda,\CY^\alpha_\lambda)\not=0$ {\em in the category
of all perverse sheaves on the space} $\CA^\alpha_\lambda$. This latter
$\Ext$ can be calculated purely topologically. So its vanishing gives a
criterion of $\Ext^1$-vanishing in the category $\FS$.

\subsection{}
In calculating $\Ext^1(\CL(\lambda_0),\CL(\mu))$ we will distinguish between
the following three cases:

a) $\mu\not\geq\lambda_0$;

b) $\mu=\lambda_0$;

c) $\mu>\lambda_0$.

\subsection{} Let us treat the case a).

\subsubsection{}
\label{tri shriek}
{\bf Lemma.} {\em
For any $\alpha\in\BN[I]$ the sheaf $\CL(\lambda_0)^\alpha_{\lambda_0}$ is
the shriek-extension from the open stratum of toric stratification of
$\CA^\alpha_{\lambda_0}$.}

{\bf Proof.}
Due to the factorization property it is enough to check that the stalk
of $\CL(\lambda_0)^\alpha_{\lambda_0}$ at the point $\{0\}\in
\CA^\alpha_{\lambda_0}$ vanishes for any $\alpha\in\BN[I], \alpha\not=0$.
By the Theorem II.8.23, we have $(\CL(\lambda_0)^\alpha_{\lambda_0})_0=
\;_\alpha C^\bullet_{\ff}(L(\lambda_0))\simeq 0$ since $L(\lambda_0)$ is a free
$\ff$-module by ~\cite{l1} 36.1.5 and
Theorem II.11.10(b) and, consequently,
$C^\bullet_{\ff}(L(\lambda_0))\simeq H^0_{\ff}(L(\lambda_0))=\sk$
and has weight zero. $\Box$

Returning to the case a), let us choose $\nu\in X$, $\nu\geq\lambda_0,
\nu\geq\mu$. For any $\alpha$, the sheaf $\CL':=\CL(\mu)_{\nu}^{\alpha}$
is supported on the subspace
$$
\CA':=\sigma(\CA_{\mu}^{\alpha+\mu-\nu})\subset\CA_{\nu}^{\alpha}
$$
and $\CL'':=\CL(\lambda_0)_{\nu}^{\alpha}$ --- on the subspace
$$
\CA'':=\sigma(\CA_{\lambda_0}^{\alpha+\lambda_0-\nu})\subset\CA^{\alpha}_{\nu}.
$$
Let $i$ denote a closed embedding
$$
i: \CA''\hra\CA^{\alpha}_{\nu}
$$
and $j$ an open embedding
$$
j: \CAO'':=\CA''-\CA''\cap\CA'\hra\CA''.
$$
We have by adjointness
$$
R\Hom_{\CA^{\alpha}_{\nu}}(\CL'',\CL')=
R\Hom_{\CA''}(\CL'',i^!\CL')=
$$
(by the previous lemma)
$$
=R\Hom_{\CAO''}(j^*\CL'',j^*i^!\CL')=0
$$
since obviously $j^*i^!\CL'=0$. This proves the vanishing in the case a).

\subsection{}
\label{tri b}
In case (b), suppose
$$
0\lra\CL(\lambda_0)\lra\CX\lra\CL(\lambda_0)\lra 0
$$
is a nonsplit extension in $\FS$. Then for $\alpha\gg 0$ the restriction
of $\CX^\alpha_{\lambda_0}$ to the open toric stratum of
$\CA^\alpha_{\lambda_0}$ is a nonsplit extension
$$
0\lra\CID^\alpha_{\lambda_0}\lra{\overset{\bullet}{\CX}}^\alpha_{\lambda_0}
\lra\CID^\alpha_{\lambda_0}\lra 0
$$
(in the category of all perverse sheaves
on $\CAD^\alpha_{\lambda_0}$) (we can restrict to the open toric stratum
because of lemma ~\ref{tri shriek}). This contradicts to the
factorization property of FFS $\CX$. This contradiction
completes the consideration of case (b).

\subsection{}
\label{tri c}
In case (c), suppose $\Ext^1(\CL(\lambda_0),\CL(\mu))\not=0$ whence
$\Ext^1(\CL(\lambda_0)^\alpha_\mu, \CL(\mu)^\alpha_\mu)\not=0$ for some
$\alpha\in\BN[I]$ by the principle ~\ref{tri principle}. Here the latter $\Ext$
is taken in the category of all perverse sheaves on $\CA^\alpha_\mu$.
We have $\mu-\lambda_0=\beta'$ for some $\beta\in\BN[I],\beta\not=0$.

Let us consider the closed embedding
$$
\sigma:\CA':=\CA_{\lambda_0}^{\alpha-\beta}\hra\CA_{\mu}^{\alpha};
$$
let us denote by $j$ an embedding of the open toric stratum
$$
j:\CAD':= \CAD_{\lambda_0}^{\alpha-\beta}\hra\CA'.
$$
As in the previous case, we have
$$
R\Hom_{\CA^{\alpha}_{\mu}}(\CL(\lambda_0)^{\alpha}_{\mu},
\CL(\mu)^{\alpha}_{\mu})=
R\Hom_{\CA'}(\CL(\lambda_0)^{\alpha}_{\mu},\sigma^!\CL(\mu)^{\alpha}_{\mu})=
R\Hom_{\CAD'}(j^*\CL(\lambda_0)^{\alpha}_{\mu},
j^*\sigma^!\CL(\mu)^{\alpha}_{\mu}).
$$

We claim that $j^*\sigma^!\CL(\mu)^{\alpha}_{\mu}=0$.
Since the sheaf $\CL(\mu)^{\alpha}_{\mu}$ is Verdier auto-dual up to replacing
$\zeta$ by $\zeta^{-1}$, it suffices to check that
$j^*\sigma^*\CL(\mu)^{\alpha}_{\mu}=0$.

To prove this vanishing, by factorization property of $\CL(\mu)$, it is
enough to check that the stalk of the sheaf $\CL(\mu)^\beta_\mu$ at the
point $\{0\}\in \CA^\beta_\mu$ vanishes.

By the Theorem II.8.23, we have $(\CL(\mu)^\beta_\mu)_0=\;_\beta
C^\bullet_{\ff}(L(\mu))$. By the Theorem II.11.10 and Shapiro Lemma,
we have $_\beta C^\bullet_{\ff}(L(\mu))\simeq C^\bullet_U(M(\lambda_0)
\otimes L(\mu))$.

By the Theorem 36.1.5. of ~\cite{l1}, the canonical projection
$M(\lambda_0)\lra L(\lambda_0)$ is an isomorphism. By the autoduality of
$L(\lambda_0)$ we have $C^\bullet_U(L(\lambda_0)\otimes L(\mu))\simeq
R\Hom^\bullet_U(L(\lambda_0),L(\mu))\simeq 0$ since $L(\lambda_0)$ is a
projective $U$-module, and $\mu\not=\lambda_0$.

This completes the case c) and the proof of the theorem. $\Box$

\section{Equivalence}

We keep the assumptions of the previous section.

\subsection{Theorem.}
\label{tri equiv thm} {\em Functor $\Phi:\FS\lra\CC$ is an equivalence.}

\subsection{Lemma.} {\em For any $\lambda\in X$ the FFS
$\CL(\lambda_0)\dotimes\CL(\lambda)$ is projective.

As $\lambda$ ranges through $X$, these sheaves form an ample system
of projectives in $\FS$.}

{\bf Proof.} We have
$$
\Hom(\CL(\lambda_0)\dotimes\CL(\lambda),?)=
\Hom(\CL(\lambda_0),\CL(\lambda)^*\ \dotimes ?)
$$
by the rigidity, and the last functor is exact since
$\dotimes$ is a biexact functor in $\FS$, and $\CL(\lambda_0)$ is
projective. Therefore, $\CL(\lambda_0)\dotimes\CL(\lambda)$ is
projective.

To prove that these sheaves form an ample system of projectives,
it is enough to show that for each $\mu\in X$ there exists
$\lambda$ such that
$\Hom(\CL(\lambda_0)\dotimes\CL(\lambda),\CL(\mu))\neq 0$.
We have
$$
\Hom(\CL(\lambda_0)\dotimes\CL(\lambda),\CL(\mu))=
\Hom(\CL(\lambda),\CL(\lambda_0)^*\dotimes\CL(\mu)).
$$
Since the sheaves $\CL(\lambda)$ exhaust irreducibles in $\FS$,
there exists $\lambda$ such that
$\CL(\lambda)$ embeds into $\CL(\lambda_0)^*\dotimes\CL(\mu)$, hence
the last group is non-zero. $\Box$

\subsection{Proof of ~\ref{tri equiv thm}} As $\lambda$ ranges through $X$,
the modules $\Phi(\CL(\lambda_0)\dotimes
\CL(\lambda))=L(\lambda_0)\otimes L(\lambda)$ form an ample system of
projectives in $\CC$. By the Lemma A.15 of ~\cite{kl}IV we only have to
show that
$$
\Phi:\ \Hom_{\FS}(\CL(\lambda_0)\dotimes\CL(\lambda),\CL(\lambda_0)\dotimes
\CL(\mu))\lra
\Hom_{\CC}(L(\lambda_0)\otimes L(\lambda),L(\lambda_0)\otimes L(\mu))
$$
is an isomorphism for any $\lambda,\mu\in X$. We already know that
it is an injection. Therefore, it remains to compare the dimensions
of the spaces in question. We have
$$
\dim\Hom_{\FS}(\CL(\lambda_0)\dotimes\CL(\lambda),\CL(\lambda_0)\dotimes
\CL(\mu))=
\dim\Hom_{\FS}(\CL(\lambda_0),\CL(\lambda_0)\dotimes \CL(\mu)\dotimes
\CL(\lambda)^*)
$$
by rigidity,
$$
=[\CL(\lambda_0)\dotimes\CL(\mu)\dotimes \CL(\lambda)^*\ :\CL(\lambda_0)]
$$
because $\CL(\lambda_0)$ is its own indecomposable
projective cover in $\FS$,
$$
=[L(\lambda_0)\otimes L(\mu)\otimes L(\lambda)^*:L(\lambda_0)]
$$
since $\Phi$ induces an isomorphism of $K$-rings of the categories
$\FS$ and $\CC$,
$$
=\dim\Hom_{\CC}(L(\lambda_0)\otimes L(\lambda),L(\lambda_0)\otimes L(\mu))
$$
by the same argument applied to $\CC$. The theorem is proven. $\Box$

\section{The case of generic $\zeta$}
\label{tri generic zeta}

In this section we suppose that $\zeta$ is not a root of unity.

\subsection{}
\label{tri U versus u}
Recall the notations of II.11,12. We have the algebra $U$ defined in
II.12.2, the algebra $\fu$ defined in II.12.3, and the homomorphism
$R:\ U\lra\fu$ defined in II.12.5.

\subsubsection{}
{\bf Lemma.} {\em
The map $R:\ U\lra\fu$ is an isomorphism.}

{\bf Proof} follows from ~\cite{r}, no. 3, Corollaire. $\Box$

\subsection{}
We call $\Lambda\in X$ dominant if $\langle i,\Lambda\rangle\geq 0$ for any
$i\in I$. An irreducible $U$-module $L(\Lambda)$ is finite dimensional
if only if $\Lambda$ is dominant.
Therefore we will need a larger category $\CO$ containing
all irreducibles $L(\Lambda)$.

Define $\CO$ as a category consisting of all $X$-graded
$U$-modules $V=\oplus_{\mu\in X}V_\mu$ such that

a) $V_\mu$ is finite dimensional for any $\mu\in X$;

b) there exists $\lambda=\lambda(V)$ such that $V_\mu=0$ if
$\mu\not\geq\lambda(V)$.

\subsubsection{} {\bf Lemma}
\label{tri CO} {\em
The category $\CO$ is equivalent to the usual category $\CO_{\fg}$ over the
classical finite dimensional Lie algebra $\fg$.}

{\bf Proof.}
See ~\cite{f1}. $\Box$

\subsection{}
\label{tri Ext over U}
Let $W$ denote the Weyl group of our root datum. For $w\in W,\ \lambda\in X$
let $w\cdot\lambda$ denote the usual action of $W$ on $X$ centered at $-\rho$.

Finally, for $\Lambda\in X$ let $M(\Lambda)\in\CO$ denote the $U^-$-free
Verma module with highest weight $\Lambda$.

\subsubsection{} {\bf Corollary.}
\label{tri ext=0} {\em
Let $\mu,\nu\in X$ be such that $W\cdot\mu\not=W\cdot\nu$. Then
$\Ext^\bullet(M(\nu),L(\mu))=0$.}

{\bf Proof.} $M(\nu)$ and $L(\mu)$ have
different central characters. $\Box$

\subsection{}
\label{tri equiv thm gener}
{\bf Theorem.} {\em
Functor $\Phi:\ \FS\lra\CC$ is an equivalence.}

{\bf Proof.} We know that $\Phi(\CL(\Lambda))\simeq L(\Lambda)$ for any
$\Lambda\in X$. So $\Phi(\CX)$ is finite dimensional iff all the irreducible
subquotients of $\CX$ are of the form $\CL(\lambda),\ \lambda$ dominant.
By virtue of Lemma ~\ref{tri CO} above the category $\CC$ is semisimple:
it is equivalent to the category of finite dimensional $\fg$-modules. It
consists of finite direct sums of modules $L(\lambda),\ \lambda$ dominant.
So to prove the Theorem it suffices to check semisimplicity of $\FS$.
Thus the Theorem follows from

\subsection{Lemma}
\label{tri vanish ext} {\em
Let $\mu,\nu\in X$ be the dominant weights. Then
$\Ext^1(\CL(\mu),\CL(\nu))=0$.}

{\bf Proof.} We will distinguish between the following two cases:

(a) $\mu=\nu$;

(b) $\mu\not=\nu$.

In calculating $\Ext^1$ we will use the principle ~\ref{tri principle}.
The argument in case (a) is absolutely similar to the one in section
{}~\ref{tri b}, and we leave it to the reader.

In case (b) either $\mu-\nu\not\in Y\subset X$ --- and then the sheaves
$\CL(\nu)$ and $\CL(\mu)$ are supported on the different connected components
of $\CA$, whence $\Ext^1$ obviously vanishes, --- or there exists
$\lambda\in X$ such that $\lambda\geq\mu,\nu$. Let us fix such $\lambda$.
Suppose $\Ext^1(\CL(\mu),\CL(\nu))\not=0$. Then according to the principle
{}~\ref{tri principle} there exists $\alpha\in\BN[I]$ such that
$\Ext^1(\CL(\mu)^\alpha_\lambda,\CL(\nu)^\alpha_\lambda)\not=0$. The latter
$\Ext$ is taken in the category of all perverse sheaves on
$\CA^\alpha_\lambda$.

We have $\Ext^1(\CL(\mu)^\alpha_\lambda,\CL(\nu)^\alpha_\lambda)=
R^1\Gamma(\CA^\alpha_\lambda, D(\CL(\mu)^\alpha_\lambda\otimes
D\CL(\nu)^\alpha_\lambda))$ where $D$ stands for Verdier duality,
and $\otimes$ denotes the usual tensor product of constructible comlexes.

We will prove that
\begin{equation}
\label{tri zero}
\CL(\mu)^\alpha_\lambda\otimes D\CL(\nu)^\alpha_\lambda=0
\end{equation}
 and hence we will arrive at the contradiction.
Equality (\ref{tri zero}) is an immediate corollary of the lemma we presently
formulate.

For $\beta\leq\alpha$ let us consider the canonical embedding
$$
\sigma: \CAD^{\alpha-\beta}\hra\CA^\alpha
$$
and denote its image by $\CAD'$ (we omit the lower case indices).

\subsubsection{} {\bf Lemma.}
\label{tri vanish stalk} {\em
(i) If $\sigma^*\CL(\mu)^\alpha_\lambda\not=0$ then
$\lambda-\beta\in W\cdot\mu$.

(ii) If $\sigma^!\CL(\mu)^\alpha_\lambda\not=0$ then
$\lambda-\beta\in W\cdot\mu$.}

To deduce Lemma ~\ref{tri vanish ext} from this lemma we notice first that
the sheaf
$\CL(\mu)^\alpha_\lambda$ is Verdier autodual up to replacing $\zeta$ by
$\zeta^{-1}$. Second, since the $W$-orbits of $\nu$ and $\mu$ are disjoint,
we see that over any toric stratum $\CAD'\subset\CA^\alpha$ at least one
of the factors of (\ref{tri zero}) vanishes.

It remains to prove Lemma ~\ref{tri vanish stalk}.
We will prove (i), while (ii) is just dual.
Let us denote $\beta+\mu-\lambda$ by $\gamma$. If $\gamma\not\in\BN[I]$ then
(i) is evident. Otherwise,
by the factorizability condition it is enough to check that
the stalk of $\CL(\mu)^\gamma_\mu$ at the origin in $\CA^\gamma_\mu$ vanishes
if $\mu-\gamma\not\in W\cdot\mu$.
Let us denote $\mu-\gamma$ by $\chi$.

By the Theorem II.8.23, we have $(\CL(\mu)^\gamma_\mu)_0=\
_\gamma C^\bullet_{U^-}(L(\mu))\simeq C^\bullet_U(M(\chi)\otimes L(\mu))$
which is dual to $\Ext^\bullet_U(M(\chi),L(\mu))$. But the latter $\Ext$
vanishes by the Corollary ~\ref{tri ext=0} since $W\cdot\chi\not=W\cdot\mu$. $\Box$

This completes the proof of Lemma ~\ref{tri vanish ext} together with
Theorem ~\ref{tri equiv thm gener}.

\newpage
\setcounter{section}{0}
\centerline{\large \bf Part IV. LOCALIZATION ON $\BP^1$}

\section{Introduction}

\subsection{} 
Given a collection of $m$ finite factorizable sheaves $\{\CX_k\}$,
we construct here some perverse sheaves over configuration
spaces of points on a projective line $\BP^1$ with $m$ additional
marked points.

We announce here (with sketch proof) the computation of the
cohomology spaces of these sheaves. They turn out
to coincide with certain "semiinfinite" $\Tor$ spaces of the corresponding
$\fu$-modules introduced by S.Arkhipov. For a precise formulation see
Theorem ~\ref{chetyre global thm}.

This result is strikingly similar to the following hoped-for
picture of affine Lie algebra representation theory, explained to us by
A.Beilinson. Let $M_1,M_2$ be two modules over an affine Lie algebra $\hfg$
on the critical level. One hopes that there is a localization functor
which associates to these modules perverse sheaves $\Delta(M_1),
\Delta(M_2)$ over the semiinfinite flag space $\hat{G}/\hat{B}^0$. Suppose that
$\Delta(M_1)$ and $\Delta(M_2)$ are equivariant with respect to the opposite
Borel subgroups of $\hat{G}$. Then the intersection $S$ of their supports
is finite dimensional, and one hopes that
$$
R^{\blt}\Gamma(S,\Delta(M_1)\otimes\Delta(M_2))=
\Tor_{\infh-\blt}^{\hfg}(M_1,M_2)
$$
where in the right hand side we have
the Feigin (Lie algebra) semiinfinite homology.

As a corollary of Theorem ~\ref{chetyre global thm} we get a description of
local systems of conformal blocks in WZW models in genus zero
(cf. ~\cite{ms}) as natural subquotients of some semisimple local systems
of geometric origin. In particular, these local systems are
semisimple themselves.


\subsection{} We are grateful to G.Lusztig for the permission
to use his unpublished results. Namely, Theorem ~\ref{chetyre braid thm}
about braiding in the category $\CC$ is due to G.Lusztig.
Chapter 2 (semiinfinite homological algebra in $\CC$) is an exposition
of the results due to S.Arkhipov (see ~\cite{a}).

We are also grateful to A.Kirillov, Jr. who explained to us how
to handle the conformal blocks of non simply laced Lie algebras.

\subsection{}
Unless specified otherwise, we will keep the notations of parts I,II,III.
For $\alpha=\sum_i\ a_ii\in\BN[I]$ we will use the notation
$|\alpha|:=\sum_i\ a_i$.

References to {\em loc.cit.} will look like Z.1.1 where Z=I, II or III.

We will keep assumptions of III.1.4 and III.16. In particular,
a "quantization" parameter $\zeta$ will be a primitive $l$-th root
of unity where $l$ is a fixed positive number prime to $2,3$.

\newpage
\begin{center}
{\bf Chapter 1. Gluing over $\BP^1$}
\end{center}
\vspace{.8cm}

\section{Cohesive local system}

\subsection{Notations} Let $\alpha\in\BN[X]$,
$\alpha=\sum a_\mu\mu$.
We denote by $\supp\alpha$ the subset of $X$ consisting of all $\mu$ such that
$a_\mu\not=0$.
Let  $\pi:\ J\lra X$ be an {\em unfolding} of $\alpha$, that is
a map of sets such that
$|\pi^{-1}(\mu)|=a_\mu$ for any $\mu\in X$.
As always, $\Sigma_\pi$ denotes the group of
automorphisms of $J$ preserving the fibers of $\pi$.

$\BP^1$ will denote a complex projective line. The $J$-th cartesian power
$\BP^{1J}$ will be denoted by $\CP^J$. The group $\Sigma_\pi$ acts
naturally on $\CP^J$, and the quotient
space $\CP^J/\Sigma_\pi$ will be denoted by $\CP^{\alpha}$.

$\CP^{oJ}$ (resp., $\CP^{o\alpha}$) stands for the complement
to diagonals in $\CP^J$ (resp., in $\CP^{\alpha}$).

$T\CP^J$ stands for the total space of the tangent bundle to $\CP^J$;
its points are couples $((x_j),(\tau_j))$ where $(x_j)\in\CP^J$ and
$\tau_j$ is a tangent vector at $x_j$. An open subspace
$$
\DT\CP^J\subset T\CP^J
$$
consists of couples with $\tau_j\neq 0$ for all $j$.
So, $\DT\CP^J\lra\CP^J$ is a $(\BC^*)^J$-torsor.
We denote by $T\CP^{oJ}$ its restriction to $\CP^{oJ}$.
The group $\Sigma_\pi$ acts freely on $T\CP^{oJ}$, and we denote the quotient
$T\CP^{oJ}/\Sigma_\pi$ by $T\CP^{o\alpha}$.

The natural projection
$$
T\CP^{oJ}\lra T\CP^{o\alpha}
$$
will be denoted by $\pi$, or sometimes by $\pi_J$.


\subsection{}
Let $\BP^1_{st}$ ($st$ for "standard") denote "the"
projective line with fixed coordinate $z$; $D_{\epsilon}\subset\BP^1_{st}$ ---
the open disk of radius $\epsilon$ centered at $z=0$; $D:=D_1$.
We will also use the notation $D_{(\epsilon,1)}$ for the open
annulus $D-\overline{D}_{\epsilon}$ (bar means the closure).

The definitions of $D^J, D^{\alpha}, D^{oJ}, D^{o\alpha}, \DT D^J, TD^{oJ},
TD^{o\alpha}$,
etc., copy the above definitions, with $D$ replacing $\BP^1$.

\subsection{} Given a finite set $K$, let
$\tCP^K$ denote the space of $K$-tuples $(u_k)_{k\in K}$ of algebraic
isomorphisms $\BP^1_{st}\iso\BP^1$ such that the images $u_k(D)$
do not intersect.

Given a $K$-tuple $\valpha=(\alpha_k)\in\BN[X]^K$ such that
$\alpha=\sum_k\ \alpha_k$, define a space
$$
\CP^{\valpha}:=\tCP^K\times\prod_{k\in K}\ \DT D^{\alpha_k}
$$
and an open subspace
$$
\CP^{o\valpha}:=\tCP^K\times\prod_{k\in K}\ TD^{o\alpha_k}.
$$
We have an evident "substitution" map
$$
q_{\valpha}:\CP^{\valpha}\lra\DT\CP^{\alpha}
$$
which restricts to $q_{\valpha}:\CP^{o\valpha}\lra T\CP^{o\alpha}$.

\subsubsection{} In the same way we define the spaces
$TD^{\valpha}$, $TD^{o\valpha}$.

\subsubsection{}
\label{chetyre assoc space} Suppose that we have an epimorphism
$\xi:L\lra K$, denote $L_k:=\xi^{-1}(k)$. Assume that
each $\alpha_k$ is in turn decomposed
as
$$
\alpha_k=\sum_{l\in L_k} \ \alpha_{l},\ \alpha_{l}\in\BN[X];
$$
set $\valpha_k=(\alpha_{l})\in\BN[X]^{L_k}$. Set
$\valpha_L=(\alpha_l)\in\BN[X]^L$.

Let us define spaces
$$
\CP^{\valpha_L;\xi}=\tCP^{K}\times\prod_{k\in K}\ \DT D^{\valpha_k}
$$
and
$$
\CP^{o\valpha_L;\xi}=\tCP^{K}\times\prod_{k\in K}\ TD^{o\valpha_k}.
$$
We have canonical substitution maps
$$
q_{\valpha_L;\xi}^1:\CP^{\valpha_L;\xi}\lra\CP^{\valpha}
$$
and
$$
q_{\valpha_L;\xi}^2:\CP^{\valpha_L;\xi}\lra\CP^{\valpha_L}.
$$
Obviously,
$$
q_{\valpha}\circ q_{\valpha_L;\xi}^1=q_{\valpha_L}\circ q_{\valpha_L;\xi}^2.
$$

\subsection{Balance function}
\label{chetyre fun n}
Consider a function $n:\ X\lra\BZ[\frac{1}{2\det A}]$ such that $$
n(\mu+\nu)=n(\mu)+n(\nu)+\mu\cdot\nu
$$
It is easy to see that $n$ can be written in the following form:
$$
n(\mu)=\frac{1}{2}\mu\cdot\mu+\mu\cdot\nu_0
$$
for some $\nu_0\in X$.
{}From now on we fix such a function $n$ and hence the corresponding
$\nu_0$.

\subsection{} For an arbitrary $\alpha\in\BN[X]$,
let us define a one-dimensional local system
$\CI^{\alpha}_D$ on $TD^{o\alpha}$. We will proceed in the same way as in
III.3.1.

Pick an unfolding of $\alpha$, $\pi:J\lra X$. Define a local system
$\CI^J_D$ on $TD^{oJ}$ as follows: its stalk at each point
$((\tau_j),(x_j))$ where all $x_j$ are real, and all the tangent vectors
$\tau_j$ are real and directed to the right, is $\sk$. Monodromies are:

--- $x_i$ moves counterclockwise around $x_j$: monodromy is
    $\zeta^{-2\pi(i)\cdot\pi(j)}$;

--- $\tau_j$ makes a counterclockwise circle: monodromy is
    $\zeta^{-2n(\pi(j))}$.

This local system has an evident $\Sigma_{\pi}$-equivariant structure,
and we define a local system $\CI^{\alpha}_D$ as
$$
\CI^{\alpha}_D=(\pi_*\CI^J_D)^{\sgn}
$$
where $\pi:\ TD^{oJ}\lra TD^{o\alpha}$ is the canonical projection,
and $(\bullet)^{\sgn}$ denotes the subsheaf of skew $\Sigma_{\pi}$-invariants.

\subsection{} We will denote the unique homomorphism
$$
\BN[X]\lra X
$$
identical on $X$, as $\alpha\mapsto\alpha^{\sim}$.

\subsubsection{}
\label{chetyre admis}
{\bf Definition.} {\em An element $\alpha\in\BN[X]$
is called {\em admissible} if $\alpha^{\sim}\equiv -2\nu_0\bmod lY$.}

\subsection{} We have a canonical "$1$-jet at $0$" map
$$
p_K:\tCP^K\lra T\CP^{oK}
$$

\subsection{Definition} {\em A {\em cohesive local system} (CLS)
(over $\BP^1$) is the following collection of data:

(i) for each admissible $\alpha\in\BN[X]$ a one-dimensional
local system $\CI^{\alpha}$ over $T\CP^{o\alpha}$;

(ii) for each decomposition $\alpha=\sum_{k\in K}\ \alpha_k$, $\alpha_k\in
\BN[X]$, a
factorization isomorphism
$$
\phi_{\valpha}:q^*_{\valpha}\CI^{\alpha}\iso
p_K^*\pi_K^*\CI^{\alpha_K}\boxtimes\Boxtimes_k\ \CI^{\alpha_k^{\sim}}_D
$$
Here $\alpha_K:=\sum_k\ \alpha_k^{\sim}\in\BN[X]$ (note that
$\alpha_K$ is obviously admissible); $\pi_K: T\CP^{oK}\lra T\CP^{o\alpha}$
is the symmetrization map.

These isomorphisms must satisfy the following

{\em Associativity axiom.} In the assumptions of ~\ref{chetyre assoc space} the
equality
$$
\phi_{\valpha_L;\xi}\circ
q_{\valpha_L;\xi}^{1*}(\phi_{\valpha})=q_{\valpha_L;\xi}^{2*}(\phi_{\valpha_L})
$$
should hold. Here $\phi_{\valpha_L;\xi}$ is induced by the evident
factorization isomorphisms for local systems on the disk $\CI^{\valpha_k}_D$.}

Morphisms between CLS's are defined in the obvious way.

\subsection{Theorem.} {\em Cohesive local systems over $\BP^1$
exist. Every two CLS's are isomorphic.
The group of automorphisms of a CLS is $\sk^*$.}

This theorem is a particular case of a more general theorem, valid
for curves of arbitrary genus, to be proved in Part V. We leave the proof
in the case of $\BP^1$ to the interested reader.

\section{Gluing}

\subsection{}
\label{chetyre bal fun} Let us define an element $\rho\in X$ by the condition
$\langle \rho,i\rangle=1$ for all $i\in I$. From now on we choose a balance
function $n$, cf. ~\ref{chetyre fun n}, in the form
$$
n(\mu)=\frac{1}{2}\mu\cdot\mu + \mu\cdot\rho.
$$
It has the property that $n(-i')=0$ for all $i\in I$. Thus, in the notations
of {\em loc. cit.} we set
$$
\nu_0=\rho.
$$

We pick a corresponding CLS $\CI=\{\CI^{\beta},\ \beta\in\BN[X]\}$.

Given $\alpha=\sum a_ii\in \BN[I]$ and $\vmu=(\mu_k)\in X^K$, we define an
element
$$
\alpha_{\vmu}=\sum a_i\cdot (-i')+\sum_k\mu_k\in\BN[X]
$$
where the sum in the right hand side is a formal one. We say that a pair
$(\vmu,\alpha)$ is {\em admissible} if $\alpha_{\vmu}$ is admissible in
the sense of the previous section, i.e.
$$
\sum_k\ \mu_k-\alpha\equiv -2\rho\bmod lY.
$$
Note that given $\vmu$, there exists $\alpha\in\BN[I]$ such that
$(\vmu,\alpha)$
is admissible if and only if $\sum_k\ \mu_k\in Y$; if this holds true,
such elements $\alpha$ form an obvious countable set.

We will denote by
$$
e:\BN[I]\lra \BN[X]
$$
a unique homomorphism sending $i\in I$ to $-i'\in X$.

\subsection{} Let us consider the space $\DT\CP^K\times\DT\CP^{e(\alpha)}$;
its points are quadruples $((z_k),(\tau_k),(x_j),(\omega_j))$ where
$(z_k)\in\CP^K,\ \tau_k$ --- a non-zero tangent vector to $\BP^1$ at
$z_k$, $(x_j)\in\CP^{e(\alpha)}$, $\omega_j$ --- a non-zero tangent
vector at $x_j$. To a point $z_k$ is assigned a weight $\mu_k$,
and to $x_j$ --- a weight $-\pi(j)'$. Here $\pi:J\lra I$ is an unfolding
of $\alpha$ (implicit in the notation $(x_j)=(x_j)_{j\in J}$).

We will be interested in some open subspaces:
$$
\DT\CP^{\alpha}_{\vmu}:=T\CP^{oK}\times\DT\CP^{e(\alpha)}\subset
\DT\CP^{\alpha_{\vmu}}
$$
and
$$
T\CP^{o\alpha}_{\vmu}\subset
\DT\CP^{\alpha}_{\vmu}
$$
whose points are quadruples $((z_k),(\tau_k),(x_j),(\omega_j))\in
\DT\CP^{\alpha}_{\vmu}$  with all $z_k\neq x_j$. We have an obvious
symmetrization projection
$$
p^{\alpha}_{\vmu}: T\CP^{o\alpha}_{\vmu}\lra T\CP^{o\alpha_{\vmu}}.
$$
Define a space
$$
\CP^{\alpha}_{\vmu}=T\CP^{oK}\times\CP^{e(\alpha)};
$$
its points are triples $((z_k),(\tau_k),(x_j))$ where
$(z_k)$, $(\tau_k)$ and $(x_j)$ are as above; and to $z_k$ and $x_j$
the weights as above are assigned. We have the canonical projection
$$
\DT\CP^{\alpha}_{\vmu}\lra \CP^{\alpha}_{\vmu}.
$$
We define the open subspaces
$$
\CP^{o\alpha}_{\vmu}\subset
\CP^{\bullet\alpha}_{\vmu}\subset
\CP^{\alpha}_{\vmu}.
$$
Here the $\bullet$-subspace (resp., $o$-subspace)
consists of all $((z_k),(\tau_k),(x_j))$ with
$z_k\neq x_j$ for all $k,j$ (resp., with all $z_k$ and $x_j$ distinct).

We define the {\em principal stratification} $\CS$ of $\CP^{\alpha}_{\vmu}$
as the stratification generated by subspaces $z_k=x_j$ and $x_j=x_{j'}$
with $\pi(j)\neq\pi(j')$. Thus, $\CP^{o\alpha}_{\vmu}$ is the open
stratum of $\CS$. As usually, we will denote by the same letter the
induced stratifications on subspaces.

The above projection restricts to
$$
T\CP^{o\alpha}_{\vmu}\lra \CP^{o\alpha}_{\vmu}.
$$

\subsection{Factorization structure}
\subsubsection{} Suppose we are given $\valpha\in\BN[I]^K,\ \beta\in\BN[I]$;
set $\alpha:=\sum_k\ \alpha_k$. Define a space
$$
\CP^{\valpha,\beta}_{\vmu}\subset
\tCP^K\times\prod_k\ D^{\alpha_k}\times\CP^{e(\beta)}
$$
consisting of all collections $((u_k),((x^{(k)}_j)_k),(y_j))$
where $(u_k)\in\tCP^K,\ (x^{(k)}_j)_k\in D^{\alpha_k},\
(y_j)\in\CP^{e(\beta)}$, such that
$$
y_j\in\BP^1-\bigcup_{k\in K}\ \overline{u_k(D)}
$$
for all $j$ (the bar means closure).

We have canonical maps
$$
q_{\valpha,\beta}:\CP^{\valpha,\beta}_{\vmu}\lra\CP^{\alpha+\beta}_{\vmu},
$$
assigning to $((u_k),((x^{(k)}_j)_k),(y_j))$ a configuration
$(u_k(0)),(\overset{\bullet}{u}_k(\tau)),(u_k(x_j^{(k)})),(y_j))$,
where $\tau$ is the unit tangent vector to $D$ at $0$,
and
$$
p_{\valpha,\beta}:\CP^{\valpha,\beta}_{\vmu}\lra
\prod_k\ D^{\alpha_k}\times\CP^{\bullet\beta}_{\vmu-\valpha}
$$
sending $((u_k),((x^{(k)}_j)_k),(y_j))$ to $((u_k(0)),
(\overset{\bullet}{u}_k(\tau)),(y_j))$.

\subsubsection{} Suppose we are given $\valpha,\vbeta\in\BN[I]^K,\
\gamma\in\BN[I]$; set $\alpha:=\sum_k\ \alpha_k,\ \beta:=\sum_k\ \beta_k$.

Define a space $D^{\alpha,\beta}$ consisting of couples
$(D_{\epsilon},(x_j))$ where $D_{\epsilon}\subset D$ is some smaller disk
($0<\epsilon<1$), and $(x_j)\in D^{\alpha+\beta}$ is a configuration
such that $\alpha$ points dwell inside $D_{\epsilon}$, and $\beta$ points ---
outside $\overline{D_{\epsilon}}$.
We have an evident map
$$
q_{\alpha,\beta}:D^{\alpha,\beta}\lra D^{\alpha+\beta}.
$$

Let us define a space
$$
\CP^{\valpha,\vbeta,\gamma}_{\vmu}\subset
\tCP^K\times\prod_{k\in K}\ D^{\alpha_k,\beta_k}\times\CP^{e(\gamma)}
$$
consisting of all triples $((u_k),\bx,(y_j))$ where
$(u_k)\in\tCP^K,\ \bx\in \prod_{k}\ D^{\alpha_k,\beta_k},\
(y_j)\in\CP^{e(\gamma)}$ such that
$$
y_j\in\BP^1-\bigcup_k\ \overline{u_k(D)}.
$$
We have obvious projections
$$
q^1_{\valpha,\vbeta,\gamma}: \CP^{\valpha,\vbeta,\gamma}_{\vmu}\lra
\CP^{\valpha+\vbeta,\gamma}_{\vmu}
$$
and
$$
q^2_{\valpha,\vbeta,\gamma}: \CP^{\valpha,\vbeta,\gamma}_{\vmu}\lra
\CP^{\valpha,\beta+\gamma}_{\vmu}
$$
such that
$$
q_{\valpha+\vbeta,\gamma}\circ q^1_{\valpha,\vbeta,\gamma}=
q_{\valpha,\beta+\gamma}\circ q^2_{\valpha,\vbeta,\gamma}.
$$
We will denote the last composition by $q_{\valpha,\vbeta,\gamma}$.

We have a natural projection
$$
p_{\valpha,\vbeta,\gamma}:
\CP^{\valpha,\vbeta,\gamma}_{\vmu}\lra
\prod_k\ D_{\frac{1}{2}}^{\alpha_k}
\times\prod_k\ D_{(\frac{1}{2},1)}^{\beta_k}\times
\CP^{\bullet\gamma}_{\vmu-\valpha-\vbeta}.
$$

\subsection{}
\label{chetyre change} Let us consider a local system
$p^{\alpha*}_{\vmu}\CI^{\alpha_{\vmu}}$
over $T\CP^{o\alpha}_{\vmu}$.
By our choice of the balance function $n$, its monodromies with
respect to the rotating of tangent vectors $\omega_j$ at points $x_j$
corresponding to negative simple roots, are trivial. Therefore
it descends to a unique local system over
$\CP^{o\alpha}_{\vmu}$, to be denoted by $\CI^{\alpha}_{\mu}$.

We define a perverse sheaf
$$
\CI^{\bullet\alpha}_{\vmu}:= j_{!*}\CI^{\alpha}_{\vmu}
[\dim \CP^{\alpha}_{\vmu}]\in
\CM(\CP^{\bullet\alpha}_{\vmu};\CS).
$$

\subsection{Factorizable sheaves over $\BP^1$}
\label{chetyre fact p} Suppose we are given
a $K$-tuple of FFS's $\{\CX_k\},\ \CX_k\in\FS_{c_k},\ k\in K,\ c_k\in X/Y$,
where $\sum_k\ c_k=0$.
Let us pick
$\vmu=(\vmu_k)\geq (\lambda(\CX_k))$.

Let us call a {\em factorizable sheaf over $\BP^1$ obtained by
gluing the sheaves $\CX_k$} the following collection of data which
we will denote by $g(\{\CX_k\})$.

(i) For each $\alpha\in\BN[I]$ such that $(\vmu,\alpha)$ is admissible,
a sheaf $\CX^{\alpha}_{\vmu}\in\CM(\CP^{\alpha}_{\vmu};\CS)$.

(ii) For each $\valpha=(\alpha_k)\in\BN[I]^K,\ \beta\in \BN[I]$ such that
$(\mu,\alpha+\beta)$ is admissible (where $\alpha=\sum\ \alpha_k$),
a {\em factorization isomorphism}
$$
\phi_{\valpha,\beta}: q^*_{\valpha,\beta}\CX^{\alpha+\beta}_{\vmu}\iso
p^*_{\valpha,\beta}((\Boxtimes_{k\in K}\ \CX^{\alpha_k}_{\mu_k})
\boxtimes\CI^{\bullet\beta}_{\vmu-\valpha}).
$$

These isomorphisms should satisfy

{\em Associativity property.} The following two isomorphisms
$$
q^*_{\valpha,\vbeta,\gamma}\CX^{\alpha+\beta+\gamma}_{\vmu}\iso
p^*_{\valpha,\vbeta,\gamma}(
(\Boxtimes_k\ \CX^{\alpha_k}_{\mu_k})\boxtimes
(\Boxtimes_k\ \CI^{\bullet\beta_k}_{\mu_k-\alpha_k})\boxtimes
\CI^{\bullet\gamma}_{\vmu-\valpha-\vbeta})
$$
are equal:
$$
\psi_{\vbeta,\gamma}\circ
q^{2*}_{\valpha,\vbeta,\gamma}(\phi_{\valpha,\beta+\gamma})=
\phi_{\valpha,\vbeta}\circ
q^{1*}_{\valpha,\vbeta,\gamma}(\phi_{\valpha+\vbeta,\gamma}).
$$
Here $\psi_{\vbeta,\gamma}$ is the factorization isomorphism
for $\CI^{\bullet}$, and $\phi_{\valpha,\vbeta}$ is the tensor
product
of factorization isomorphisms for the sheaves $\CX_k$.

\subsection{Theorem}
\label{chetyre glu thm} {\em There exists a unique up to a canonical
isomorphism factorizable sheaf over $\BP^1$ obtained by gluing
the sheaves $\{\CX_k\}$.}

{\bf Proof} is similar to III.10.3. $\Box$

\newpage
\begin{center}
{\bf Chapter 2. Semiinfinite cohomology}
\end{center}
\vspace{.5cm}

In this chapter we discuss, following essentially ~\cite{a},
the "Semiinfinite homological algebra" in the category $\CC$.

\section{Semiinfinite functors $Ext$ and $Tor$ in $\CC$}

\subsection{} Let us call an $\fu$-module {\em $\fu^-$-induced}
(resp., {\em $\fu^+$-induced}) if it is induced from some
$\fu^{\geq 0}$ (resp., $\fu^{\leq 0}$)-module.

\subsubsection{}
\label{chetyre tens ind} {\bf Lemma.} {\em If $M$ is a $\fu^-$-induced, and
$N$ is $\fu^+$-induced then $M\otimes_\sk N$ is $\fu$-projective.}

{\bf Proof.} An induced module has a filtration whose factors are
corresponding Verma modules. For Verma modules the claim is easy. $\Box$

\subsection{Definition}
\label{chetyre downup} Let $M^{\bullet}=\oplus_{\lambda\in X}\
M^{\bullet}_{\lambda}$ be a complex (possibly unbounded) in $\CC$.
We say that $M^{\bullet}$ is {\em concave} (resp. {\em convex})
if it satisfies the properties (a) and (b) below.

(a) There exists $\lambda_0\in X$ such that for any $\lambda\in X$,
if $M^{\bullet}_{\lambda}\neq 0$ then $\lambda\geq\lambda_0$
(resp. $\lambda\leq\lambda_0$).

(b) For any $\mu\in X$ the subcomplex
$\oplus_{\lambda\leq\mu}\ M^{\bullet}_{\lambda}$
(resp. $\oplus_{\lambda\geq\mu}\ M^{\bullet}_{\lambda}$)
is finite.

We will denote the category of concave (resp., convex)
complexes by $\CC^{\uparrow}$ (resp., $\CC^{\downarrow}$).

\subsection{} Let $V\in\CC$. We will say that a surjection
$\phi: P\lra V$ is {\em good} if it satisfies the following
properties:

(a) $P$ is $\fu^-$-induced;


(b) Let $\mu\in\supp P$ be an extremal point, that is, there is no
$\lambda\in\supp P$ such that $\lambda>\mu$. Then
$\mu\not\in\supp(\ker\phi)$.

For any $V$ there exists a good surjection as above. Indeed, denote by $p$
the projection $p:M(0)\lra L(0)$,
and take for $\phi$ the map $p\otimes\id_V$.

\subsection{}
\label{chetyre rel hoch}
Iterating, we can construct a $\fu^-$-induced convex left resolution of
$\sk=L(0)$. Let us pick such a resolution and denote it by
$P^{\blt}_{\swarrow}$:
$$
\ldots\lra P^{-1}_{\swarrow}\lra
P^{0}_{\swarrow}\lra L(0)\lra 0
$$

We will denote by
$$
^*:\CC\lra\CC^{opp}
$$
the rigidity in $\CC$
(see e.g. ~\cite{ajs}, 7.3). We denote by $P^{\bullet}_{\nearrow}$
the complex $(P^{\bullet}_{\swarrow})^*$. It is a $\fu^-$-induced
concave right resolution of $\sk$. The fact that
$P^{\bullet}_{\nearrow}$ is $\fu^-$-induced follows since
$\fu^-$ is Frobenius (see e.g. ~\cite{x}).

\subsection{} In a similar manner, we can construct a
$\fu^+$-induced concave left resolution of $\sk$. Let us pick
such a resolution and denote it $P^{\blt}_{\nwarrow}$:
$$
\ldots\lra P^{-1}_{\nwarrow}\lra
P^{0}_{\nwarrow}\lra L(0)\lra 0
$$

We denote by $P^{\bullet}_{\searrow}$
the complex $(P^{\bullet}_{\nwarrow})^*$. It is a $\fu^+$-induced
convex right resolution of $\sk$.

\subsubsection{} For $M\in\CC$ we denote by
$P^{\bullet}_{\swarrow}(M)$ (resp., $P^{\bullet}_{\nearrow}(M)$,
$P^{\bullet}_{\nwarrow}(M)$, $P^{\bullet}_{\searrow}(M)$) the resolution
$P^{\bullet}_{\swarrow}\otimes_\sk M$ (resp.,
$P^{\bullet}_{\nearrow}\otimes_\sk M$, $P^{\bullet}_{\nwarrow}\otimes_\sk M$,
$P^{\bullet}_{\searrow}\otimes_\sk M$) of $M$.

\subsection{} We denote by $\CC_r$ the category of $X$-graded
right $\fu$-modules $V=\oplus_{\lambda\in X} V_{\lambda}$ such that
$$
K_i|_{V_{\lambda}}=\zeta^{-\langle i,\lambda\rangle}
$$
(note the change of a sign!), the operators $E_i$, $F_i$ acting as
$E_i:V_{\lambda}\lra V_{\lambda+i'},\ F_i:V_{\lambda}\lra V_{\lambda-i'}$.

\subsubsection{} Given $M\in\CC$, we define $M^{\vee}\in\CC_r$ as follows:
$(M^{\vee})_{\lambda}=(M_{-\lambda})^*$,
$E_i:(M^{\vee})_{\lambda}\lra (M^{\vee})_{\lambda+i'}$ is the
transpose of $E_i:M_{-\lambda-i'}\lra M_{-\lambda}$, similarly,
$F_i$ on $M^{\vee}$ is the transpose of $F_i$ on $M$.

This way we get an equivalence
$$
^{\vee}:\CC^{opp}\iso\CC_r.
$$
Similarly, one defines an equivalence $^{\vee}:\CC_r^{opp}\lra\CC$,
and we have an obvious isomorphism $^{\vee}\circ\ ^{\vee}\cong\Id$.

\subsubsection{} Given $M\in\CC$, we define $sM\in\CC_r$ as follows:
$(sM)_{\lambda}=M_{\lambda};\ xg=(sg)x$ for $x\in M, g\in\fu$ where
$$
s:\fu\lra\fu^{opp}
$$
is the antipode defined in ~\cite{ajs}, 7.2. This way we get an
equivalence
$$
s:\CC\iso\CC_r
$$
One defines an equivalence $s:\CC_r\iso\CC$ in a similar manner.
The isomorphism of functors $s\circ s\cong \Id$ is constructed
in {\em loc. cit.}, 7.3.

Note that the rigidity $*$ is just the composition
$$
^*=s\circ\ ^{\vee}.
$$

\subsubsection{} We define the categories $\CC^{\uparrow}_r$ and
$\CC^{\downarrow}_r$ in the same way as in ~\ref{chetyre downup}.

For $V\in\CC_r$ we define $P^{\bullet}_{\swarrow}(V)$ as
$$
P^{\bullet}_{\swarrow}(V)=sP^{\bullet}_{\swarrow}(sV);
$$
and $P^{\bullet}_{\nearrow}(V), P^{\bullet}_{\searrow}(V),
P^{\bullet}_{\nwarrow}(V)$ in a similar way.

\subsection{Definition} (i) Let $M,N\in\CC$. We define
$$
\Ext^{\infh+\bullet}_{\CC}(M,N):=
H^{\bullet}(\Hom_{\CC}(P^{\bullet}_{\searrow}(M),P^{\bullet}_{\nearrow}(N))).
$$

(ii) Let $V\in\CC_r,\ N\in\CC$. We define
$$
\Tor_{\infh+\bullet}^{\CC}(V,N):=
H^{-\bullet}(P^{\bullet}_{\swarrow}(V)\otimes_{\CC}P^{\bullet}_{\searrow}(N)).\
\Box
$$

Here we understand
$\Hom_{\CC}(P^{\bullet}_{\searrow}(M),P^{\bullet}_{\nearrow}(N)$ and
$P^{\bullet}_{\swarrow}(V)\otimes_{\CC}P^{\bullet}_{\searrow}(N)$ as simple
complexes associated with the corresponding double complexes.
Note that due to our boundedness properties of weights
of our resolutions, these double complexes are bounded.
Therefore all $\Ext^{\infh+i}$ and $\Tor_{\infh+i}$ spaces are finite
dimensional, and are non-zero only for finite number of $i\in\BZ$.

\subsection{Lemma.}
\label{chetyre ext tor} {\em For $M,N\in\CC$ there exist canonical
nondegenerate pairings
$$
\Ext^{\infh+n}_{\CC}(M,N)\otimes
\Tor_{\infh+n}^{\CC}(N^{\vee},M)\lra \sk.
$$}

{\bf Proof.} There is an evident non-degenenerate pairing
$$
\Hom_{\CC}(M,N)\otimes (N^{\vee}\otimes_{\CC}M)\lra \sk.
$$
It follows that the complexes computing
$\Ext$ and $\Tor$ are also canonically dual. $\Box$

\subsection{Theorem}
\label{chetyre semi thm} {\em
{\em (i)} Let $M,N\in\CC$. Let $R^{\bullet}_{\searrow}(M)$ be a
$\fu^+$-induced convex right resolution of $M$,
and $R^{\bullet}_{\nearrow}(N)$ --- a $\fu^-$-induced concave right resolution
of $N$. Then there is a canonical isomorphism
$$
\Ext^{\infh+\bullet}_{\CC}(M,N)\cong
H^{\bullet}(\Hom_{\CC}(R^{\bullet}_{\searrow}(M),R^{\bullet}_{\nearrow}(N))).
$$

{\em (ii)} Let $V\in\CC_r,N\in\CC$. Let $R^{\bullet}_{\swarrow}(V)$ be a
$\fu^-$-induced convex left resolution of $V$,
and $R^{\bullet}_{\searrow}(N)$ --- a $\fu^+$-induced convex right resolution
of $N$ lying in $\CC^{\downarrow}$.
Then there is a canonical isomorphism
$$
\Tor_{\infh+\bullet}^{\CC}(V,N)\cong
H^{-\bullet}(R^{\bullet}_{\swarrow}(V)\otimes_{\CC}R^{\bullet}_{\searrow}(N)).
$$}

{\bf Proof} will occupy the rest of the section.

\subsection{Lemma}
\label{chetyre filtr} {\em Let $V\in\CC_r$; let $R_i^{\bullet}$,
$i=1,2,$ be two $\fu^-$-induced convex left resolutions of $V$.
There exists a third $\fu^-$-induced
convex left resolution $R^{\bullet}$ of $V$, together
with two termwise surjective maps
$$
R^{\bullet}\lra R^{\bullet}_i,\ i=1,2,
$$
inducing identity on $V$.}

{\bf Proof.} We will construct $R^{\bullet}$ inductively,
from right to left. Let
$$
R^{\bullet}_i:\ \ldots\overset{d_i^{-2}}{\lra} R^{-1}_i
\overset{d_i^{-1}}{\lra} R^0_i\overset{\epsilon_i}{\lra} V.
$$
First, define $L_0:=R^0_1\times_V R^0_2$.
We denote by $\delta$ the canonical map $L^0\lra V$, and by
$q_i^0$ the projections $L^0\lra R^0_i$. Choose a good surjection
$\phi_0: R^0\lra L^0$ and define
$p^0_i:=q^0_i\circ\phi_0: R^0\lra R^0_i;\
\epsilon:=\delta\circ\phi_0:R^0\lra V$.

Set $K_i^{-1}:=\ker\epsilon_i; K:=\ker\epsilon$. The projections $p^0_i$
induce surjections $p^0_i:K^{-1}\lra K^{-1}_i$. Let us define
$$
L^{-1}:=\ker ((d_1^{-1}-p^0_1,d_2^{-1}-p_2^0):
R_1^{-1}\oplus K^{-1}\oplus R_2^{-1}\lra
K^{-1}_1\oplus K_2^{-1}).
$$
We have canonical projections
$q_i^{-1}:L^{-1}\lra R_i^{-1},\
\delta^{-1}:L^{-1}\lra K^{-1}$. Choose a good surjection
$\phi_{-1}:R^{-1}\lra L^{-1}$ and write $d^{-1}:R^{-1}\lra R^0$ for
$\delta^{-1}\circ\phi_{-1}$ composed with the inclusion
$K^{-1}\hra R^0$. We define $p_i^{-1}:=q_i^{-1}\circ\phi_{-1}$.

We have just described an induction step, and we can proceed in the same
manner. One sees directly that the left $\fu^-$-induced resolution
$R^{\bullet}$ obtained this way actually lies in $\CC^{\downarrow}_r$. $\Box$

\subsection{} Let $N\in\CC$, and let $R^{\bullet}_{\searrow}$ be
a $\fu^+$-induced convex right resolution of $N$.
For $n\geq 0$ let $b_{\geq n}(R^{\bullet}_{\searrow})$ denote the
stupid truncation:
$$
 0\lra R^0_{\searrow}\lra\ldots\lra R^n_{\searrow}\lra 0\lra\ldots
$$
For $m\geq n$ we have evident truncation maps
$b_{\geq m}(R^{\bullet}_{\searrow})\lra b_{\geq n}(R^{\bullet}_{\searrow})$.

\subsubsection{} {\bf Lemma.}
\label{chetyre stabil semi} {\em Let $R^{\bullet}_{\swarrow}$ be a
$\fu^-$-induced left resolution of a module $V\in\CC_r$. We have
$$
H^{\bullet}(R^{\bullet}_{\swarrow}\otimes_{\CC}R^{\bullet}_{\searrow})=
\invlim_n\
H^{\bullet}(R^{\bullet}_{\swarrow}\otimes_{\CC}b_{\leq n}
R^{\bullet}_{\searrow}).
$$

For every $i\in\BZ$ the inverse system
$$
\{H^i(R^{\bullet}_{\swarrow}\otimes_{\CC}b_{\leq n}R^{\bullet}_{\searrow})\}
$$
stabilizes.}

{\bf Proof.} All spaces
$H^i(R^{\bullet}_{\swarrow}\otimes_{\CC}b_{\leq n}R^{\bullet}_{\searrow})$
and only finitely many weight components of $R^{\bullet}_{\swarrow}$ and
$R^{\bullet}_{\searrow}$ contribute to $H^i$. $\Box$

\subsection{Proof of Theorem ~\ref{chetyre semi thm}} Let us consider case (ii),
and prove that $H^{\bullet}(R^{\bullet}_{\swarrow}(V)\otimes_{\CC}
R^{\bullet}_{\searrow}(N))$ does not depend, up to a canonical
isomorphism, on the choice of a resolution $R^{\bullet}_{\swarrow}(V)$.
The other independences are proved exactly in the same way.

Let $R^{\bullet}_i,\ i=1,2,$ be two left $\fu^-$-induced left convex
resolutions of $V$. According to Lemma ~\ref{chetyre filtr}, there exists a third
one, $R^{\bullet}$, projecting onto $R^{\bullet}_i$. Let us prove that
the projections induce isomorphisms
$$
H^{\bullet}(R^{\bullet}\otimes_{\CC}R^{\bullet}_{\searrow}(N))\iso
H^{\bullet}(R^{\bullet}_i\otimes_{\CC}R^{\bullet}_{\searrow}(N)).
$$
By Lemma ~\ref{chetyre stabil semi}, it suffices to prove that
$$
H^{\bullet}(R^{\bullet}\otimes_{\CC}b_{\leq n}R^{\bullet}_{\searrow}(N))\iso
H^{\bullet}(R^{\bullet}_i\otimes_{\CC}b_{\leq n}R^{\bullet}_{\searrow}(N)).
$$
for all $n$. Let $Q^{\bullet}_i$ be a cone of $R^{\bullet}\lra R^{\bullet}_i$.
It is an exact $\fu^-$-induced convex complex bounded from the right.
It is enough to check that
$H^{\bullet}(Q^{\bullet}_i\otimes_{\CC}b_{\leq n}R^{\bullet}_{\searrow}(N))=0$.

Note that for $W\in\CC_r,M\in\CC$ we have canonically
$$
W\otimes_{\CC}M=(W\otimes sM)\otimes_{\CC}\sk.
$$
Thus
$$
H^{\bullet}(Q^{\bullet}_i\otimes_{\CC}b_{\leq n}R^{\bullet}_{\searrow}(N))=
H^{\bullet}((Q^{\bullet}_i\otimes sb_{\leq n}R^{\bullet}_{\searrow}(N))
\otimes_{\CC}\sk)=0,
$$
since
$(Q^{\bullet}_i\otimes sb_{\leq n}R^{\bullet}_{\searrow}(N)$ is an exact
bounded from the right complex, consisting of modules
which are tensor products of
$\fu^+$-induced and $\fu^-$-induced, hence
$\fu$-projective modules (see Lemma ~\ref{chetyre tens ind}).

\subsubsection{} It remains to show that if $p'$ and $p''$ are two maps
between $\fu^-$-induced convex resolutions of $V$,
$R^{\bullet}_1\lra R^{\bullet}_2$,
inducing identity on $V$, then the isomorphisms
$$
H^{\bullet}(R^{\bullet}_1\otimes_{\CC}R^{\bullet}_{\searrow}(N))\iso
H^{\bullet}(R^{\bullet}_2\otimes_{\CC}R^{\bullet}_{\searrow}(N))
$$
induced by $p'$ and $p''$, coincide. Arguing as above, we see that it is enough
to prove this with $R^{\bullet}_{\searrow}(N)$ replaced by
$b_{\leq n}R^{\bullet}_{\searrow}(N)$. This in turn is equivalent to
showing that two isomorphisms
$$
H^{\bullet}((R^{\bullet}_1\otimes sb_{\leq n}R^{\bullet}_{\searrow}(N))
\otimes_{\CC}\sk)\iso
H^{\bullet}((R^{\bullet}_2\otimes sb_{\leq n}R^{\bullet}_{\searrow}(N))
\otimes_{\CC}\sk)
$$
coincide. But $R^{\bullet}_i\otimes sb_{\leq n}R^{\bullet}_{\searrow}(N)$
are complexes of projective $\fu$-modules, and the morphisms
$p'\otimes\id$ and $p''\otimes\id$ induce the same map
on cohomology, hence they are homotopic; therefore they
induce homotopic maps after tensor multiplication by $\sk$.

This completes the proof of the theorem. $\Box$

\section{Some calculations}

We will give a recipe for calculation of
$\Tor^{\CC}_{\infh+\blt}$ which will prove useful for the next chapter.

\subsection{} Recall that in III.13.2 the duality functor
$$
D:\CC_{\zeta^{-1}}\lra\CC^{opp}
$$
has been defined (we identify $\CC$ with $\tCC$ as usually). We will
denote objects of $\CC_{\zeta^{-1}}$ by letters with the subscript
$(\blt)_{\zeta^{-1}}$.

Note that $DL(0)_{\zeta^{-1}}=L(0)$.

Let us describe duals to Verma modules.
For $\lambda\in X$ let us denote by $M^+(\lambda)$ the
Verma module with respect to the subalgebra $\fu^+$ with the lowest
weight $\lambda$, that is
$$
M^+(\lambda):=\Ind^{\fu}_{\fu^{\leq 0}}\chi_{\lambda}
$$
where $\chi_{\lambda}$ is an evident one-dimensional representaion
of $\fu^{\leq 0}$ corresponding to the character $\lambda$.

\subsubsection{}
\label{chetyre dm} {\bf Lemma.} {\em We have
$$
DM(\lambda)_{\zeta^{-1}}=M^+(\lambda-2(l-1)\rho).
$$}

{\bf Proof} follows from ~\cite{ajs}, Lemma 4.10. $\Box$

\subsection{}
Let us denote by $K^{\blt}$ a two term complex in $\CC$
$$
L(0)\lra DM(0)_{\zeta^{-1}}
$$
concentrated in degrees $0$ and $1$, the morphism being dual to the
canonical projection $M(0)_{\zeta^{-1}}\lra L(0)_{\zeta^{-1}}$.

For $n\geq 1$ define a complex
$$
K^{\blt}_n:=b_{\geq 0}(K^{\blt\otimes n}[1]);
$$
it is concentrated in degrees from $0$ to $n-1$. For example,
$K^{\blt}_1=DM(0)_{\zeta^{-1}}$.

For $n\geq 1$ we will denote by
$$
\xi:K^{\blt}_n\lra K^{\blt}_{n+1}
$$
the map induced by the embedding $L(0)\hra DM(0)_{\zeta^{-1}}$.

We will need the following evident properties of the system
$\{K^{\blt}_n,\xi_n\}$:

(a) $K^{\blt}_n$ is $\fu^+$-induced;

(b) $K^{\blt}_n$ is exact off degrees $0$ and $n-1$;
$H^0(K^{\blt}_n)=\sk$. $\xi_n$ induces identity map between
$H^0(K^{\blt}_n)$ and $H^0(K^{\blt}_{n+1})$.

(c) For a fixed $\mu\in X$ there exists $m\in\BN$ such that for any
$n$ we have $(b_{\geq m}K^{\blt}_n)_{\geq \mu}=0$.
Here for $V=\oplus_{\lambda\in X}\in\CC$ we set
$$
V_{\geq\mu}:=\oplus_{\lambda\geq\mu}\ V_{\lambda}.
$$

\subsection{} Let $V\in\CC_r$; let $R^{\blt}_{\swarrow}(V)$ be a
$\fu^-$-induced convex left resolution of $V$. Let $N\in\CC$.

\subsubsection{}
\label{chetyre half tor}{\bf Lemma.} (i) {\em For a fixed $k\in\BZ$ the direct system
$\{H^k(R^{\blt}_{\swarrow}(V)\otimes_{\CC}(K^{\blt}_n\otimes N)),\xi_n\}$
stabilizes.}

(ii) {\em We have a canonical isomorphism
$$
\Tor^{\CC}_{\infh+\blt}(V,N)\cong
\dirlim_n\ H^{-\blt}(R^{\blt}_{\swarrow}(V)\otimes_{\CC}(K^{\blt}_n\otimes N)).
$$}

{\bf Proof.} (i) is similar to Lemma ~\ref{chetyre stabil semi}.
(ii) By Theorem ~\ref{chetyre semi thm} we can use any $\fu^+$-induced
right convex resolution of $N$ to compute
$\Tor^{\CC}_{\infh+\blt}(V,N)$. Now extend
$K^{\blt}_n\otimes N$ to a $\fu^+$-induced convex resolution of $N$
and argue like in the proof of Lemma ~\ref{chetyre stabil semi} again. $\Box$

\subsection{} Recall the notations of ~\ref{chetyre rel hoch}
and take for $R^{\blt}_{\swarrow}(V)$ the resolution
$P^{\blt}_{\swarrow}(V)=P^{\blt}_{\swarrow}\otimes V$. Then
$$
H^{\blt}(P^{\blt}_{\swarrow}(V)\otimes_{\CC}(K^{\blt}_n\otimes N))=
H^{\blt}(V\otimes_{\CC}(P^{\blt}_{\swarrow}\otimes K^{\blt}_n
\otimes N)).
$$
Note that $P^{\blt}_{\swarrow}\otimes K^{\blt}_n\otimes N$ is a right
bounded complex quasi-isomorphic to $K^{\blt}_n\otimes N$.
The terms of $P^{\blt}_{\swarrow}\otimes K^{\blt}_n$ are
$\fu$-projective by Lemma ~\ref{chetyre tens ind}, hence the terms of
$P^{\blt}_{\swarrow}\otimes K^{\blt}_n\otimes N$ are projective by rigidity
of $\CC$. Therefore,
$$
H^{-\blt}(V\otimes_{\CC}(P^{\blt}_{\swarrow}\otimes K^{\blt}_n
\otimes N))=
\Tor^{\CC}_{\blt}(V,K^{\blt}_n\otimes N).
$$
Here $\Tor^{\CC}_{\blt}(*,*)$ stands for the zeroth
weight component of  $\Tor^{\fu}_{\blt}(*,*)$.

Putting all the above together, we get

\subsection{Corollary}
{\em For a fixed $k\in\BZ$ the direct system
$\{\Tor_k^{\CC}(V,K^{\blt}_n\otimes N)\}$ stabilizes.

We have
$$
\Tor^{\CC}_{\infh+\blt}(V,N)=
\dirlim_n\ \Tor^{\CC}_{\blt}(V,K^{\blt}_n\otimes N).\ \Box
$$}

\subsection{} Dually, consider complexes $DK^{\blt}_{n,\zeta^{-1}}$. They form
a projective system
$$
\{\ldots\lra DK^{\blt}_{n+1,\zeta^{-1}}\lra DK^{\blt}_{n,\zeta^{-1}}\lra
\ldots \}
$$

These complexes enjoy properties dual to (a) --- (c) above.

\subsection{Theorem}
\label{chetyre two side} {\em For every $k\in\BZ$ we have canonical
isomorphisms
$$
\Tor^{\CC}_{\infh+k}(V,N)\cong
\invlim_m\dirlim_n\ H^{-k}((V\otimes sDK^{\blt}_{m,\zeta^{-1}})\otimes_{\CC}
(K^{\blt}_n\otimes N)).
$$
Both the inverse and the direct systems actually stabilize.}

{\bf Proof} follows from Lemma ~\ref{chetyre half tor}. We leave details
to the reader. $\Box$

Here is an example of calculation of $\Tor^{\CC}_{\infh+\blt}$.

\subsection{Lemma}
\label{chetyre tor rho} {\em
$\Tor^{\CC}_{\infh+\blt}(\sk,L(2(l-1)\rho))=\sk\ \mbox{in degree $0$.}$}

{\bf Proof.} According to Lemma ~\ref{chetyre ext tor} it suffices to prove
that
$\Ext_{\CC}^{\infh+\blt}(L(2(l-1)\rho),L(0))=\sk$. Choose a $\fu^+$-induced
right convex resolution
$$
L(2(l-1)\rho)\overset{\epsilon}{\lra} R^{\blt}_{\searrow}
$$
such that
$$
R^0_{\searrow}= DM(2(l-1)\rho)_{\zeta^{-1}}=M^+(0),
$$
and all the weights in $R^{\geq 1}_{\searrow}$ are $<2(l-1)\rho$.

Similarly, choose a $\fu^-$-free right concave resolution
$$
L(0)\overset{\epsilon}{\lra} R^{\blt}_{\nearrow}
$$
such that
$$
R^0_{\nearrow}=M(2(l-1)\rho)=DM^+(0)_{\zeta^{-1}},
$$
and all the weights of $R^{\geq 1}_{\nearrow}$ are $>0$. By
Theorem ~\ref{chetyre semi thm} we have
$$
\Ext_{\CC}^{\infh+\blt}(L(2(l-1)\rho),L(0))=
H^{\blt}(\Hom_{\CC}^{\blt}(R^{\blt}_{\searrow},R^{\blt}_{\nearrow})).
$$
Therefore it is enough to prove that

(a) $\Hom_{\CC}(R^0_{\searrow},R^0_{\nearrow})=\sk$;

(b) $\Hom_{\CC}(R^m_{\searrow},R^n_{\nearrow})=0$ for $(m,n)\neq (0,0)$.

(a) is evident. Let us prove (b) for, say, $n>0$. $R^m_{\searrow}$ has
a filtration with successive quotients of type
$M^+(\lambda)$, $\lambda\leq 0$; similarly,
$R^n_{\searrow}$ has
a filtration with successive quotients of type
$DM^+(\mu)_{\zeta^{-1}}$, $\mu > 0$. We have
$\Hom_{\CC}(M^+(\lambda),DM^+(\mu)_{\zeta^{-1}})=0$, therefore
$\Hom_{\CC}(R^m_{\searrow},R^n_{\nearrow})=0$. The proof for $m>0$
is similar. Lemma is proven. $\Box$

\vspace{.5cm}
{\em CONFORMAL BLOCKS AND $Tor^{\CC}_{\infh+\blt}$}
\vspace{.5cm}

\subsection{}
\label{chetyre triv sub} Let $M\in \CC$. We have a canonical embedding
$$
\Hom_{\CC}(\sk,M)\hra M
$$
which identifies $\Hom_{\CC}(\sk,M)$ with the maximal trivial
subobject of $M$. Dually, we have a canonical
epimorhism
$$
M\lra\Hom_{\CC}(M,\sk)^*
$$
which identifies $\Hom_{\CC}(M,\sk)^*$ with the maximal trivial quotient
of $M$. Let us denote by $\langle M\rangle$ the image
of the composition
$$
\Hom_{\CC}(\sk,M)\lra M\lra\Hom_{\CC}(M,\sk)^*
$$
Thus, $\langle M\rangle$ is canonically a subquotient of $M$.

One sees easily that if $N\subset M$ is a trivial direct summand of $M$
which is maximal, i.e. not contained in greater direct summand, then
we have a canonical isomorphism $\langle M\rangle\iso N$. By this reason,
we will call $\langle M\rangle$ {\em the maximal trivial
direct summand} of $M$.

\subsection{}
\label{chetyre alc} Let
$$
\Delta_l=\{\lambda\in X|\ \langle i,\lambda+\rho\rangle>0,\
\mbox{for all }i\in I;\ \langle\gamma,\lambda+\rho\rangle<l\}
$$
denote the first alcove. Here $\gamma\in \CR\subset Y$ is the
highest coroot.

For $\lambda_1,\ldots,\lambda_n\in\Delta_l$, {\em the space of conformal
blocks} is defined as
$$
\langle L(\lambda_1),\ldots,L(\lambda_n)\rangle :=
\langle L(\lambda_1)\otimes \ldots\otimes L(\lambda_n)\rangle
$$
(see e.g. ~\cite{an} and Lemma ~\ref{chetyre compar conf} below).

\subsection{Corollary}
\label{chetyre conf subq} {\em The space of conformal blocks
$\langle L(\lambda_1),\ldots,L(\lambda_n)\rangle$
is canonically a subquotient of
$\Tor^{\CC}_{\infh+0}(\sk,L(\lambda_1)\otimes \ldots\otimes
L(\lambda_n)\otimes L(2(l-1)\rho)).$}

{\bf Proof} follows easily from the definition of $\langle\blt\rangle$
and Lemma ~\ref{chetyre tor rho}. $\Box$

\subsection{}
\label{chetyre ex conf} Let us consider an example showing that
$\langle L(\lambda_1),\ldots,L(\lambda_n)\rangle$ is in general a
{\em proper} subquotient of
$\Tor^{\CC}_{\infh+0}(\sk,L(\lambda_1)\otimes\ldots\otimes
L(\lambda_n)\otimes L(2(l-1)\rho))$.

We leave the following to the reader.

\subsubsection{} {\bf Exercise.} {\em Let $P(0)$ be the indecomposable
projective cover of $L(0)$. We have
$\Tor^{\CC}_{\infh+0}(\sk,P(0))=\sk.\ \Box$}

We will construct an example featuring $P(0)$ as a direct summand of
$L(\lambda_1)\otimes \ldots\otimes L(\lambda_n)$.

Let us take a root datum of type $sl(2)$; take $l=5,\ n=4,\
\lambda_1=\lambda_2=2,\ \lambda_3=\lambda_4=3$
(we have identified $X$ with $\BZ$).

In our case $\rho=1$, so $2(l-1)\rho=8$. Note that $P(0)$ has highest
weight $8$, and it is a unique indecomposable projective with
the highest weight $8$. So, if we are able to find a
projective summand of highest weight $0$ in $V=L(2)\otimes L(2)\otimes
L(3)\otimes L(3)$ then $V\otimes L(8)$ will
contain a projective summand of highest weight $8$, i.e. $P(0)$.

Let $U_\sk$ denote the quantum group with divided powers over $\sk$
(see ~\cite{l2}, 8.1). The algebra $\fu$ lies inside $U_\sk$. It is wellknown
that all irreducibles $L(\lambda)$, $\lambda\in\Delta_l$, lift to simple
$U_\sk$-modules $\wt{L(\lambda)}$ and for
$\lambda_1,\ldots,\lambda_n\in\Delta$ the $U_\sk$-module
$\wt{L(\lambda_1)}\otimes\ldots\otimes\wt{L(\lambda_n)}$ is a direct sum
of irreducibles $\wt{L(\lambda)},\ \lambda\in\Delta_l,$
and indecomposable projectives $\wt{P(\lambda)},\ \lambda\geq 0$
(see, e.g. ~\cite{an}).

Thus $\wt{L(2)}\otimes\wt{L(2)}\otimes\wt{L(3)}\otimes\wt{L(3)}$
contains an indecomposable projective summand with the highest weight
$10$, i.e. $\wt{P(8)}$. One can check easily
that when restricted to $\fu$, $\wt{P(8)}$ remains projective
and contains a summand $P(-2)$. But the highest weight of $P(-2)$ is zero.

We conclude that
$L(2)\otimes L(2)\otimes L(3)\otimes L(3)\otimes L(8)$ contains
a projective summand $P(0)$, whence
$$
\langle L(2),L(2),L(3),L(3)\rangle\neq
\Tor^{\CC}_{\infh+0}(\sk,L(2)\otimes L(2)\otimes L(3)\otimes L(3)\otimes L(8)).
$$

\newpage
\begin{center}
{\bf Chapter 3. Global sections}
\end{center}
\vspace{.5cm}

\section{Braiding and balance in $\CC$ and $\FS$}

\subsection{}
\label{chetyre br} Let $U_\sk$ be the quantum group with divided powers,
cf ~\cite{l2}, 8.1. Let $_R\CC$ be the category of finite dimensional
integrable $U_\sk$-modules defined in ~\cite{kl}IV, \S 37.
It is a rigid braided tensor category. The braiding, i.e.
family of isomorphisms
$$
\tR_{V,W}:V\otimes W\iso W\otimes V,\ V,W\in\ _R\CC,
$$
satisfying the usual constraints, has been defined in
{}~\cite{l1}, Ch. 32.


\subsection{}
\label{chetyre ups} As $\fu$ is a subalgebra of $U_\sk$, we have the restriction
functor preserving $X$-grading
$$
\Upsilon:\ _R\CC\lra\CC.
$$
The following theorem is due to G.Lusztig (private communication).

\subsubsection{}
\label{chetyre braid thm} {\bf Theorem.} {\em
{\em (a)} There is a unique braided structure $(R_{V,W},\theta_V)$
on $\CC$ such that the restriction functor $\Upsilon$ commutes
with braiding.

{\em (b)} Let $V=L(\lambda)$, and let $\mu$ be the highest weight of
$W\in\CC$, i.e. $W_{\mu}\neq 0$ and $W_{\nu}\neq 0$ implies
$\nu\leq\mu$. Let $x\in V,\ y\in W_{\mu}$. Then
$$
R_{V,W}(x\otimes y)=\zeta^{\lambda\cdot\mu}y\otimes x;
$$

{\em (c)} Any braided structure on $\CC$
enjoying the property {\em (b)} above coincides with that defined in
{\em (a)}.
\ $\Box$}

\subsection{} Recall that an automorphism $\ttheta=\{\ttheta_V: V\iso V\}$
of the identity functor of $_R\CC$ is called {\em balance}
if for any $V,W\in\ _R\CC$ we have
$$
\tR_{W,V}\circ\tR_{V,W}=\ttheta_{V\otimes W}\circ (\ttheta_V\otimes
\ttheta_W)^{-1}.
$$

The following proposition is an easy application of the results
of ~\cite{l1}, Chapter 32.

\subsection{Proposition} {\em The category $_R\CC$ admits a unique balance
$\ttheta$ such that

--- if $\tL(\lambda)$ is an irreducible in $_R\CC$ with the
highest weight $\lambda$, then $\ttheta$ acts on $\tL(\lambda)$ as
multiplication by $\zeta^{n(\lambda)}$.\ $\Box$}

Here $n(\lambda)$ denotes the function introduced in ~\ref{chetyre bal fun}.

Similarly to ~\ref{chetyre braid thm}, one can prove

\subsection{Theorem} {\em {\em (a)} There is a unique balance
$\theta$ on $\CC$ such that $\Upsilon$ commutes with balance;

{\em (b)} $\theta_{L(\lambda)}=\zeta^{n(\lambda)}$;

{\em (c)} If $\theta'$ is a balance in $\CC$ having property {\em (b)},
then $\theta'=\theta$.\ $\Box$}

\subsection{}
\label{chetyre phi braid} According to Deligne's ideology, ~\cite{d1}, the gluing
construction of ~\ref{chetyre glu thm} provides the category $\FS$ with the
balance $\theta^{\FS}$. Recall that the braiding $R^{\FS}$ has
been defined in III.11.11 (see also 11.4). It follows easily from the
definitions that $(\Phi(R^{\FS}),\Phi(\theta^{\FS}))$ satisfy the
properties (b)(i) and (ii) above. Therefore, we have
$(\Phi(R^{\FS}),\Phi(\theta^{\FS}))=(R,\theta)$, i.e.
$\Phi$ is an equivalence of braided balanced categories.

\section{Global sections over $\CA(K)$}
\label{chetyre global}

\subsection{} Let $K$ be a finite non-empty set, $|K|=n$, and let
$\{\CX_k\}$ be a $K$-tuple of finite gactorizable sheaves. Let $\lambda_k:=
\lambda(\CX_k)$ and $\lambda=\sum_k\ \lambda_k$; let $\alpha\in\BN[I]$.
Consider the sheaf $\CX^{\alpha}(K)$ over $\CA^{\alpha}_{\vlambda}(K)$
obtained by gluing $\{\CX_k\}$, cf. III.10.3. Thus
$$
\CX^{\alpha}(K)=g_K(\{\CX_k\})^{\alpha}_{\vlambda}
$$
in the notations of {\em loc.cit.}

We will denote by $\eta$, or sometimes by $\eta^{\alpha}$,
or $\eta^{\alpha}_K$ the projection
$\CA^{\alpha}(K)\lra\CO(K)$. We are going to describe
$R\eta_*\CX^{\alpha}(K)[-n]$. Note that it is an element of
$\CCD(\CO(K))$  which is smooth, i.e. its cohomology sheaves are local
systems.

\subsection{} Let $V_1,\ldots,V_n\in\CC$. Recall (see II.3) that
$C^{\blt}_{\fu^-}(V_1\otimes\ldots\otimes V_n)$ denotes the
Hochschild complex of the $\fu^-$-module $V_1\otimes\ldots\otimes V_n$. It
is naturally $X$-graded, and its $\lambda$-component is denoted
by the subscript $(\blt)_{\lambda}$ as usually.

Let us consider a homotopy point $\bz=(z_1,\ldots,z_n)\in\CO(K)$
where all $z_i$ are real, $z_1<\ldots<z_n$. Choose a bijection
$K\iso [n]$. We want to describe a stalk
$R\eta_*\CX^{\alpha}(K)_{\bz}[-n]$. The following theorem generalizes
Theorem II.8.23. The proof is similar to {\em loc. cit}, cf. III.12.16,
and will appear later.

\subsection{Theorem}
\label{chetyre coh tens} {\em There is a canonical isomorphism,
natural in $\CX_i$,
$$
R\eta_*\CX^{\alpha}(K)_{\bz}[-n]\cong
C^{\blt}_{\fu^-}(\Phi(\CX_1)\otimes\ldots\otimes\Phi(\CX_n))_{\lambda-\alpha}.
\ \Box
$$}

\subsection{} The group $\pi_1(\CO(K);\bz)$ is generated
by counterclockwise loops of $z_{k+1}$ around $z_k$, $\sigma_i$,
$k=1,\ldots,n-1$. Let $\sigma_k$ act on
$\Phi(\CX_1)\otimes\ldots\otimes\Phi(\CX_n)$ as
$$
\id\otimes\ldots\otimes R_{\Phi(\CX_{k+1}),\Phi(\CX_k)}\circ
R_{\Phi(\CX_{k}),\Phi(\CX_{k+1})}\otimes\ldots\otimes\id.
$$
This defines an action of $\pi_1(\CO(K);\bz)$ on
$\Phi(\CX_1)\otimes\ldots\otimes\Phi(\CX_n)$,
whence we get an action of this group on
$C^{\blt}_{\fu^-}(\Phi(\CX_1)\otimes\ldots\otimes\Phi(\CX_n))$ respecting
the $X$-grading. Therefore we get a complex of local systems
over $\CO(K)$; let us denote it
$C^{\blt}_{\fu^-}(\Phi(\CX_1)\otimes\ldots\otimes\Phi(\CX_n))^{\heartsuit}$.

\subsection{Theorem} {\em There is a canonical isomorphism
in $\CCD(\CO(K))$
$$
R\eta_*\CX^{\alpha}(K)[-n]\iso
C^{\blt}_{\fu^-}(\Phi(\CX_1)\otimes\ldots\otimes\Phi(\CX_n))^{\heartsuit}_{\lambda
-\alpha}.
$$}

{\bf Proof} follows from ~\ref{chetyre phi braid} and  Theorem ~\ref{chetyre coh tens}.
$\Box$

\subsection{Corollary} {\em Set
$\lambda_{\infty}:=\alpha+2(l-1)\rho-\lambda.$
There is a canonical isomorphism in $\CCD(\CO(K))$
$$
R\eta_*\CX^{\alpha}(K)[-n]\iso
C^{\blt}_{\fu}(\Phi(\CX_1)\otimes\ldots\otimes\Phi(\CX_n)
\otimes DM(\lambda_{\infty})_{\zeta^{-1}})^{\heartsuit}_{0}.
$$}

{\bf Proof.} By Shapiro's lemma, we have a canonical morphism
of complexes which is a quasiisomorphism
$$
C^{\blt}_{\fu^-}(\Phi(\CX_1)\otimes\ldots\otimes\Phi(\CX_n))_{\lambda
-\alpha}\lra
C^{\blt}_{\fu}(\Phi(\CX_1)\otimes\ldots\otimes\Phi(\CX_n)
\otimes M^+(\alpha-\lambda))_{0}.
$$
By Lemma ~\ref{chetyre dm}, $M^+(\alpha-\lambda)=
DM(\lambda_{\infty})_{\zeta^{-1}}$. $\Box$

\section{Global sections over $\CP$}

\subsection{} Let $J$ be a finite set, $|J|=m$, and
$\{\CX_j\}$ a $J$-tuple of finite factorizable sheaves. Set
$\mu_j:=\lambda(\CX_j),\ \vmu=(\mu_j)\in X^J$. Let
$\alpha\in \BN[I]$ be such that $(\vmu,\alpha)$ is admissible,
cf. ~\ref{chetyre bal fun}. Let $\CX^{\alpha}_{\vmu}$ be the preverse
sheaf on $\CP^{\alpha}_{\vmu}$ obtained by gluing the sheaves $\CX_j$,
cf. ~\ref{chetyre fact p} and ~\ref{chetyre glu thm}.

Note that the group $\PGL_2(\BC)=\Aut(\BP^1)$ operates naturally on
$\CP^{\alpha}_{\vmu}$ and the sheaf $\CX^{\alpha}_{\vmu}$ is equivariant
with respect to this action.

Let
$$
\breta:\CP^{\alpha}_{\vmu}\lra T\CP^{oJ}
$$
denote the natural projection; we will denote this map also by
$\breta_J$ or $\breta_J^{\alpha}$. Note that $\breta$ commutes
with the natural action of $\PGL_2(\BC)$ on these spaces.
Therefore $R\breta_*\CX^{\alpha}_{\vmu}$ is a smooth
$\PGL_2(\BC)$-equivariant
complex on $T\CP^{oJ}$. Our aim in this section will be to compute
this complex algebraically.

Note that $R\breta_*\CX^{\alpha}_{\vmu}$ descends uniquely
to the quotient
$$
\ul{T\CP}^{oJ}:= T\CP^{oJ}/\PGL_2(\BC)
$$

\subsection{} Let us pick a bijection $J\iso [m]$.
Let $\ul{\CZ}$ be a contractible real submanifold
of $\ul{T\CP}^{oJ}$ defined in ~\cite{kl}II, 13.1
(under the name $\ul{\CV}_0$). Its points are configurations
$(z_1,\tau_1,\ldots,z_m,\tau_m)$ such that $z_j\in\BP^1(\BR)=S\subset
\BP^1(\BC)$; the points $z_j$ lie on $S$ in this cyclic order;
they orient $S$ in the same way as $(0,1,\infty)$ does;
the tangent vectors $\tau_j$ are real and compatible with this orientation.

\subsection{Definition}
\label{chetyre pos} An $m$-tuple of weights $\vmu\in X^m$
is called {\em positive} if
$$
\sum_{j=1}^m\ \mu_j+(1-l)2\rho\in\BN[I]\subset X.
$$
If this is so, we will denote
$$
\alpha(\vmu):= \sum_{j=1}^m\ \mu_j+(1-l)2\rho\ \Box
$$

\subsection{Theorem}
\label{chetyre stalk glob} {\em Let $\CX_1,\ldots,\CX_m\in\FS$.
Let $\vmu$ be a positive $m$-tuple of weights, $\mu_j\geq\lambda(\CX_j)$,
and let $\alpha=\alpha(\vmu)$. Let $\CX^{\alpha}_{\vmu}$ be the sheaf
on $\CP^{\alpha}_{\vmu}$ obtained by gluing the sheaves $\CX_j$. There is a
canonical isomorphism
$$
R^{\blt}\breta_*\CX^{\alpha}_{\vmu}[-2m]_{\CZ}\cong
\Tor^{\CC}_{\infh-\blt}(\sk,\Phi(\CX_1)\otimes\ldots\otimes\Phi(\CX_m)).
$$}

{\bf Proof} is sketched in the next few subsections.

\subsection{Two-sided \v{C}ech resolutions} The idea of the construction
below is inspired by ~\cite{b}, p. 40.

Let $P$ be a topological space, $\CU=\{U_i|\ i=1,\ldots, N\}$ an open
covering of $P$. Let $j_{i_0i_1\ldots i_a}$ denote the embedding
$$
U_{i_0}\cap\ldots\cap U_{i_a}\hra P.
$$

Given a sheaf $\CF$ on $P$, we have a canonical morphism
\begin{equation}
\label{chetyre right ch}
\CF\lra\Cch^{\blt}(\CU;\CF)
\end{equation}
where
$$
\Cch^a(\CU;\CF)=\oplus_{i_0<i_1<\ldots<i_a}\ j_{i_0i_1\ldots i_a*}
j_{i_0i_1\ldots i_a}^*\CF,
$$
the differential being the usual \v{C}ech one.

Dually, we define a morphism
\begin{equation}
\label{chetyre left ch}
\Cch_{\blt}(\CU;\CF)\lra\CF
\end{equation}
where
$$
\Cch_{\blt}(\CU;\CF):\
0\lra \Cch_{N}(\CU;\CF)\lra \Cch_{N-1}(\CU;\CF)\lra\ldots\lra
\Cch_{0}(\CU;\CF)\lra 0
$$
where
$$
\Cch_a(\CU;\CF)=\oplus_{i_0<i_1<\ldots<i_a}\ j_{i_0i_1\ldots i_a!}
j_{i_0i_1\ldots i_a}^*\CF.
$$
If $\CF$ is injective then the arrows ~(\ref{chetyre right ch}) and
{}~(\ref{chetyre left ch}) are quasiisomorphisms.

Suppose we have a second open covering of $P$, $\CV=\{V_j|\
j=1,\ldots, N\}$. Let us define sheaves
$$
\Cch^a_b(\CU,\CV;\CF):=
\Cch_b(\CV;\Cch^a(\CU;\CF));
$$
they form a bicomplex. Let us consider the associated simple complex
$\Cch^{\blt}(\CU,\CV;\CF)$, i.e.
$$
\Cch^{i}(\CU,\CV;\CF)=\oplus_{a-b=i}\ \Cch^a_b(\CU,\CV;\CF).
$$
It is a complex concentrated in degrees from $-N$ to $N$. We have
canonical morphisms
$$
\CF\lra\Cch^{\blt}(\CU;\CF)\lla\Cch^{\blt}(\CU,\CV;\CF)
$$
If $\CF$ is injective then both arrows are quasiisomorphisms, and
the above functors are exact on injective sheaves. Therefore,
they pass to derived categories, and we get a functor
$\CK\mapsto\Cch^{\blt}(\CU,\CV;\CK)$ from the bounded derived
category $\CCD^b(P)$ to the bounded filtered derived category
$\CCD F(P)$. This implies

\subsection{Lemma}
\label{chetyre bi cech} {\em Suppose that $\CK\in\CCD^b(P)$ is
such that $R^i\Gamma(P;\Cch^a_b(\CU,\CV;\CK))
=0$ for all $a,b$ and all $i\neq 0$. Then we have a canonical isomorphism
in $\CCD^b(P)$,
$$
R\Gamma(P;\CK)\iso R^0\Gamma(P;\Cch^{\blt}(\CU,\CV;\CK)).\ \Box
$$}

\subsection{} Returning to the assumptions of
theorem ~\ref{chetyre stalk glob}, let us pick a point
$\bz=(z_1,\tau_1,\ldots,z_m,\tau_m)\in T\CP^{oJ}$ such that $z_j$ are
real numbers
$z_1<\ldots <z_m$ and tangent vectors are directed to the right.

By definition, we have canonically
$$
R\breta_*(\CP^{\alpha}_{\vmu};\CX^{\alpha}_{\vmu})_{\ul{\CZ}}=
R\Gamma(\CP^{\alpha};\CK)
$$
where
$$
\CK:=\CX^{\alpha}_{\vmu}|_{\breta^{-1}(\bz)}[-2m].
$$
Let us pick $N\geq |\alpha|$ and reals $p_1,\ldots, p_N,q_1,\ldots,q_N$
such that
$$
p_1<\ldots <p_N<z_1<\ldots <z_m<q_N<\ldots<q_1.
$$
Let us define two open coverings $\CU=\{U_i|\ i=1,\ldots,N\}$ and
$\CV=\{V_i|\ i=1,\ldots, N\}$ of the space $\CP^{\alpha}$ where
$$
U_i=\CP^{\alpha}-\bigcup_k\ \{t_k=p_i\};\
V_i=\CP^{\alpha}-\bigcup_k\ \{t_k=q_i\},
$$
where $t_k$ denote the standard coordinates.

\subsection{Lemma} (i) {\em We have
$$
R^i\Gamma(\CP^{\alpha};\Cch^a_b(\CU,\CV;\CK))=0
$$
for all $a,b$ and all $i\neq 0$.}

(ii) {\em We have canonical isomorphism
$$
R^0\Gamma(\CP^{\alpha};\Cch^{\blt}(\CU,\CV;\CK))\cong
sDK^{\blt}_{N,\zeta^{-1}}\otimes_{\CC}
(K^{\blt}_N\otimes\Phi(\CX_1)\otimes\ldots\otimes\Phi(\CX_m)),
$$
in the notations of ~\ref{chetyre two side}.}

{\bf Proof} (sketch).  We should regard the computation
of $R\Gamma(\CP^{\alpha};\Cch^a_b(\CU,\CV;\CK))$ as the computation
of global sections over $\CP^{\alpha}$ of a sheaf obtained
by gluing $\CX_j$ into points $z_j$,
the Verma sheaves $\CM(0)$ or irreducibles $\CL(0)$ into the points
$p_j$, and dual sheaves $D\CM(0)_{\zeta^{-1}}$ or
$D\CL(0)_{\zeta^{-1}}$ into the points $q_j$.

Using $\PGL_2(\BR)$-invariance, we can
move one of the points $p_j$ to infinity. Then, the desired global sections
are reduced to global sections over an affine space $\CA^{\alpha}$, which are
calculated by means of Theorem ~\ref{chetyre coh tens}.

Note that in our situation all the sheaves $\Cch^a_b(\CU,\CV;\CK)$ actually
belong to the abelian category $\CM(\CP^{\alpha})$ of perverse
sheaves. So $\Cch^{\blt}(\CU,\CV;\CK)$ is a resolution of $\CK$ in
$\CM(\CP^{\alpha})$.
$\Box$

\subsection{} The conclusion of ~\ref{chetyre stalk glob} follow from
the previous lemma and Theorem ~\ref{chetyre two side}.
$\Box$

\subsection{} The group $\pi_1(\ul{T\CP}^{om},\CZ)$ operates on
the spaces
$\Tor^{\CC}_{\infh+\blt}(\sk,\Phi(\CX_1)\otimes\ldots\otimes\Phi(\CX_m))$
via its action on the object $\Phi(\CX_1)\otimes\ldots\otimes\Phi(\CX_m)$
induced by the braiding and balance in $\CC$. Let us denote by
$$
 \Tor^{\CC}_{\infh+\blt}(\sk,\Phi(\CX_1)\otimes\ldots\otimes\Phi(\CX_m))
^{\heartsuit}
$$
the corresponding local system on $\ul{T\CP}^{om}$.

\subsection{Theorem}
\label{chetyre global thm} {\em There is a canonical isomorphism of local
systems on $\ul{T\CP}^{om}$:
$$
R^{\blt-2m}\breta_*\CX^{\alpha}_{\vmu}\cong
\Tor^{\CC}_{\infh-\blt}(\sk,\Phi(\CX_1)\otimes\ldots\otimes\Phi(\CX_m))
^{\heartsuit}
$$}

{\bf Proof} follows immediately from ~\ref{chetyre phi braid} and
Theorem ~\ref{chetyre stalk glob}. $\Box$

\section{Application to conformal blocks}

\subsection{} In applications to conformal blocks we will encounter
the roots of unity $\zeta$ of not necessarily odd degree $l$. So we have
to generalize all the above considerations to the case of arbitrary $l$.

The definitions of the categories $\CC$ and $\FS$ do not change (for
the category $\CC$ the reader may consult ~\cite{ap}, \S3). The construction
of the functor $\Phi: \FS\lra\CC$ and the proof that $\Phi$ is an
equivalence repeats the one in III word for word.

Here we list the only minor changes (say, in the definition of the Steinberg
module) following ~\cite{l1} and ~\cite{ap}.

\subsubsection{}
So suppose $\zeta$ is a primitive root of unity of an {\em even} degree $l$.

We define $\ell:=\frac{l}{2}$. For the sake of unification of notations,
in case $l$ is {\em odd} we define $\ell:=l$.
For $i\in I$ we define $\ell_i:=\frac{\ell}
{(\ell,d_i)}$ where $(\ell,d_i)$ stands for the greatest common divisor
of $\ell$ and $d_i$.

For a coroot $\alpha\in\CR\in Y$ we can find an element $w$ of the Weyl
group $W$ and a simple coroot $i\in Y$ such that $w(i)=\alpha$ (notations
of ~\cite{l1}, 2.3). We define $\ell_\alpha:=\frac{\ell}{(\ell,d_i)}$, and
the result does not depend on a choice of $i$ and $w$.

We define $\gamma_0\in\CR$ to be the highest coroot, and $\beta_0\in\CR$
to be the coroot dual to the highest root. Note that $\gamma_0=\beta_0$ iff
our root datum is simply laced.

\subsubsection{}
We define
$$
Y_\ell:=\{\lambda\in X|\lambda\cdot\mu\in\ \ell\BZ\ \mbox{for any }\mu\in
X\}
$$

One should replace the congruence modulo $lY$ in the Definition ~\ref{chetyre admis}
and in ~\ref{chetyre bal fun} by the congruence modulo $Y_\ell$.

We define $\rho_\ell\in X$ as the unique element such that
$\langle i,\rho_\ell\rangle=\ell_i-1$ for any $i\in I$.

Then the {\em Steinberg module} $L(\rho_\ell)$ is irreducible projective
in $\CC$ (see ~\cite{ap} 3.14).

Note also that $\rho_\ell$ is the highest weight of $\fu^+$.

One has to replace
all the occurences of $(l-1)2\rho$ in the above sections
by $2\rho_\ell$.

In particular, the new formulations of the Definition ~\ref{chetyre pos} and
the Theorem ~\ref{chetyre stalk glob} force us to make the following changes
in ~\ref{chetyre bal fun} and ~\ref{chetyre change}.

In ~\ref{chetyre bal fun} we choose a balance function $n$ in the form
$$n(\mu)=\frac{1}{2}\mu\cdot\mu-\mu\cdot\rho_\ell$$
In other words, we set $\nu_0=-\rho_\ell$.
{\em This balance function does not necessarily have the property that
$n(-i')\equiv 0 \bmod l$. It is only true that $n(-i')\equiv 0\bmod\ell$.}

We say that a pair $(\vmu,\alpha)$ is admissible if $\sum_k\mu_k-\alpha
\equiv 2\rho_\ell\bmod Y_\ell$.


\subsubsection{}
\label{chetyre alcove}
The last change concerns the definition of the first alcove in ~\ref{chetyre alc}.

The corrected definition reads as follows:

if $\ell_i=\ell$ for any $i\in I$, then
$$
\Delta_l=\{\lambda\in X|\ \langle i,\lambda+\rho\rangle>0,\
\mbox{for all }i\in I;\ \langle\gamma_0,\lambda+\rho\rangle<\ell\};
$$
if not, then
$$
\Delta_l=\{\lambda\in X|\ \langle i,\lambda+\rho\rangle>0,\
\mbox{for all }i\in I;\ \langle\beta_0,\lambda+\rho\rangle<\ell_{\beta_0}\}
$$

\subsection{}
Let $\hfg$ denote the affine
Lie algebra associated with $\fg$:
$$
0\lra\BC\lra\hfg\lra\fg((\epsilon))\lra 0.
$$
Let $\tCO_\kappa$ be the category of integrable $\hfg$-modules with the central
charge $\kappa-h$ where $h$ stands for the dual
 Coxeter number of $\fg$.
It is a semisimple balanced braided rigid tensor category
(see e.g. ~\cite{ms} or ~\cite{f}).

Let $\CO_{-\kappa}$ be the category of $\fg$-integrable $\hfg$-modules
of finite length with the central charge $-\kappa-h$. It is a
balanced braided rigid tensor (bbrt) category (see ~\cite{kl}).
Let $\tCO_{-\kappa}$ be the semisimple subcategory of $\CO_{-\kappa}$
formed by direct sums of simple $\hfg$-modules with highest weights
in the alcove $\nabla_\kappa$:
$$
\nabla_\kappa:=\{\lambda\in X|\ \langle i,\lambda+\rho\rangle>0,\
\mbox{for all }i\in I;\ \langle\beta_0,\lambda+\rho\rangle<\kappa\}
$$
The bbrt structure
on $\CO_{-\kappa}$ induces the one on $\tCO_{-\kappa}$, and one can construct
an equivalence
$$
\tCO_\kappa\iso\tCO_{-\kappa}
$$
respecting bbrt structure (see ~\cite{f}).
D.Kazhdan and G.Lusztig have constructed an equivalence
$$
\CO_{-\kappa}\iso\ _R\CC_\zeta
$$
(notations of ~\ref{chetyre br}) respecting bbrt structure (see
{}~\cite{kl} and ~\cite{l3}).
Here $\zeta=\exp(\frac{\pi\sqrt{-1}}{d\kappa})$ where $d=\max_{i\in I}d_i$.
Thus $l=2d\kappa$, and $\ell=d\kappa$.

Note that the alcoves $\nabla_\kappa$ and $\Delta_l$ (see ~\ref{chetyre alcove})
coincide.

The Kazhdan-Lusztig equivalence induces an equivalence
$$
\tCO_{-\kappa}\iso\tCO_{\zeta}
$$
where $\tCO_{\zeta}$ is the semisimple subcategory of $_R\CC_\zeta$
formed by direct sums of simple $U_\sk$-modules $\tL(\lambda)$
with $\lambda\in\Delta$ (see ~\cite{an} and ~\cite{ap}).
 The bbrt structure on $_R\CC_\zeta$
induces the one on $\tCO_{\zeta}$, and the last equivalence respects bbrt
structure. We denote the composition
of the above equivalences by
$$
\phi:\tCO_\kappa\iso\tCO_{\zeta}.
$$

Given any bbrt category $\CB$ and objects $L_1,\ldots, L_m\in\CB$ we obtain
a local system $\Hom_{\CB}(\One,L_1\otimes\ldots\otimes L_m)^{\heartsuit}$ on
$T\CP^{om}$ with monodromies induced by the action of braiding and balance on
$L_1\otimes\ldots\otimes L_m$.

Here and below we write
a superscript $X^{\heartsuit}$ to denote a local system over $T\CP^{om}$
with the fiber at a standard real point
$z_1<\ldots<z_m$ with tangent vectors looking to the right, equal to $X$.

Thus, given $L_1,\ldots,L_m\in\tCO_\kappa$, the local system
$$
\Hom_{\tCO_\kappa}(\One,L_1\totimes\ldots\totimes L_m)^{\heartsuit}
$$
called {\em local system of conformal blocks} is isomorphic
to the local system
$\Hom_{\tCO_{\zeta}}(\One,\phi(L_1)\totimes\ldots\totimes \phi(L_m))
^{\heartsuit}$.
Here $\totimes$ will denote the tensor product in "tilded" categories.

To unburden the notations we leave out the subscript $\zeta$ in
$_R\CC_\zeta$ from now on.

For an object $X\in\ _R\CC$ let us define  a vector space
$\langle X\rangle_{_R\CC}$ in the same manner as in ~\ref{chetyre triv sub}, i.e.
as an image of the canonical map from the maximal trivial
subobject of $X$ to the maximal trivial quotient of $X$. Given $X_1,\ldots,
X_m\in\ _R\CC$, we denote
$$
\langle X_1,\ldots X_m\rangle :=
\langle X_1\otimes\ldots\otimes X_m\rangle_{_R\CC}.
$$

\subsubsection{} {\bf Lemma.} {\em We have an isomorphism of local systems
$$
\Hom_{\tCO_{\zeta}}(\One,\phi(L_1)\totimes\ldots\totimes\phi(L_m))
^{\heartsuit}\cong
\langle\phi(L_1),\ldots,\phi(L_m)\rangle_{_R\CC}^{\heartsuit}
$$}

{\bf Proof.} Follows from ~\cite{an}. $\Box$

\subsection{Lemma}
\label{chetyre compar conf} {\em The restriction functor
$\Upsilon:\ _R\CC\lra\CC$ (cf. ~\ref{chetyre ups}) induces isomorphism
$$
\langle\phi(L_1),\ldots,\phi(L_m)\rangle_{_R\CC}\iso
\langle\Upsilon\phi(L_1),\ldots,\Upsilon\phi(L_m)\rangle_{\CC}.
$$}

{\bf Proof.} We must prove that if $\lambda_1,\ldots,\lambda_m\in\Delta$,
$\wt{L(\lambda_1)},\ldots,\wt{L(\lambda_m)}$ are corresponding simples
in $_R\CC$, and $L(\lambda_i)=\Upsilon\wt{L(\lambda_i)}$ ---
the corresponding simples in $\CC$, then
the maximal trivial direct summand of
$\wt{L(\lambda_1)}\otimes\ldots\otimes\wt{L(\lambda_m)}$ in $_R\CC$
maps isomorphically to the maximal trivial direct summand of
$L(\lambda_1)\otimes\ldots\otimes L(\lambda_m)$ in $\CC$.

According to ~\cite{an}, ~\cite{ap},
$\wt{L(\lambda_1)}\otimes\ldots\otimes\wt{L(\lambda_m)}$ is a direct sum
of a module
$\wt{L(\lambda_1)}\totimes\ldots\totimes\wt{L(\lambda_m)}\in\CO_{\zeta}$ and
a negligible module $N\in\ _R\CC$. Here {\em negligible} means
that any endomorphism of $N$ has quantum trace zero (see {\em loc. cit.}).
Moreover, it is proven in {\em loc cit.} that $N$ is a direct summand
of $W\otimes M$ for some $M\in\ _R\CC$ where
$W=\oplus_{\omega\in\Omega}\ \wt{L(\omega)}$,
$$
\Omega=\{\omega\in X|\ \langle i,\omega+\rho\rangle>0\
\mbox{for all}\ i\in I;\ \langle\beta_0,\omega+\rho\rangle=\kappa\}
$$
being the affine wall of the first alcove.
By {\em loc. cit.}, $W$ is negligible. Since
$\Upsilon\wt{L(\omega)}=L(\omega)$, $\omega\in\Omega$ and since
$\Upsilon$ commutes with braiding, balance and rigidity, we see that
the modules $L(\omega)$ are negligible in $\CC$. Hence $\Upsilon W$ is
negligible, and $\Upsilon W\otimes\Upsilon M$ is negligible, and
finally $\Upsilon N$ is negligible. This implies that $\Upsilon N$ cannot have
trivial summands (since $L(0)$) is not negligible).

We conclude that
\begin{eqnarray}\nonumber
\langle\Upsilon(\wt{L(\lambda_1)}\otimes\ldots\otimes \wt{L(\lambda_m)})\rangle
_{\CC} =
\langle\Upsilon\wt{L(\lambda_1)}\totimes\ldots\totimes
\Upsilon\wt{L(\lambda_m)}\rangle_{\CC}=\\ \nonumber
\langle\wt{L(\lambda_1)}\totimes\ldots\totimes\wt{L(\lambda_m)}\rangle_{_R\CC}
=
\langle\wt{L(\lambda_1)}\otimes\ldots\otimes\wt{L(\lambda_m)}\rangle_{_R\CC}
\ \Box\nonumber
\end{eqnarray}

\subsection{} Corollary ~\ref{chetyre conf subq} implies that the local system
$$
\langle\Upsilon\phi(L_1),\ldots,\Upsilon\phi(L_m)\rangle_{\CC}^{\heartsuit}
$$
is canonically a subquotient of the local system
$$
\Tor^{\CC}_{\infh+0}(\sk,\Upsilon\phi(L_1)\otimes\ldots\otimes
\Upsilon\phi(L_m)\otimes L(2\rho_\ell)^{\heartsuit}
$$
(the action of monodromy being induced by braiding
and balance on the first $m$ factors).

\subsection{} Let us fix a point $\infty\in\BP^1$ and a nonzero tangent
vector $v\in T_{\infty}\BP^1$.
This defines an open subset
$$
T\CA^{om}\subset T\CP^{om}
$$
and the locally closed embedding
$$
\xi: T\CA^{om}\hra T\CP^{om+1}.
$$
Given $\lambda_1,\ldots,\lambda_m\in\Delta$, we consider the integrable
$\hfg$-modules $\hL(\lambda_1),\ldots,\hL(\lambda_m)$ of central charge
$\kappa-h$.

Suppose that
$$
\lambda_1+\ldots +\lambda_m=\alpha\in\BN[I]\subset X.
$$
We define $\lambda_{\infty}:=2\rho_\ell$, and
$\vlambda:=(\lambda_1,\ldots,\lambda_m,\lambda_{\infty})$.
Note that $\vlambda$ is positive and
$\alpha=\alpha(\vlambda)$, in the notations of ~\ref{chetyre pos}.

Denote by $\CX^{\alpha}_{\vlambda}$ the sheaf on $\CP_{\vlambda}^{\alpha}$
obtained by gluing $\CL(\lambda_1),\ldots,\CL(\lambda_m),
\CL(\lambda_{\infty})$. Note that
$$
\CX^{\alpha}_{\vlambda}=j_{!*}\CI^{\alpha}_{\vlambda}
$$
where $j:\CP^{o\alpha}_{\vlambda}\hra \CP^{\alpha}_{\vlambda}$.

Consider the local system of conformal blocks
$$
\Hom_{\tCO_\kappa}(\One,\hL(\lambda_1)\totimes\ldots\totimes\hL(\lambda_m))
^{\heartsuit}.
$$
If $\sum_{i=1}^m\lambda_i\not\in\BN[I]\subset X$ then it vanishes by the above
comparison with its "quantum group" incarnation.

\subsection{Theorem} {\em Suppose that $\sum_{i=1}^m\lambda_i=\alpha\in\BN[I]$.
Then the local system of conformal blocks restricted to
$T\CA^{om}$ is isomorphic to a canonical subquotient of a
"geometric" local system
$$
\xi^*R^{-2m-2}\breta^{\alpha}_{m+1*}j_{!*}\CI^{\alpha}_{\vlambda}.
$$}

{\bf Proof.} This follows from Theorem ~\ref{chetyre global thm} and the previous
discussion. $\Box$

\subsection{Corollary} {\em The above local system of conformal blocks
is semisimple. It is a direct summand of the geometric local system
above.}

{\bf Proof.} The geometric system is semisimple by Decomposition theorem,
{}~\cite{bbd}, Th\'{e}or\`{e}me 6.2.5. $\Box$

\subsection{} Example ~\ref{chetyre ex conf} shows that in general
a local system of conformal blocks is a {\em proper} direct summand
of the corresponding geometric system.


\newpage
\setcounter{section}{0}
\begin{center}{\large \bf Part V. MODULAR STRUCTURE}\end{center}

\begin{center}{\large\bf ON THE CATEGORY $\FS$}\end{center}

\section{Introduction}

\subsection{} 

Let $C\lra S$ be a smooth proper morphism of relative dimension 1.
Let $\vx=(x_1,\ldots,x_m)$ be an $m$-tuple of disjoint sections $x_k:\
S\lra C$. In this part we will show how to localize $\fu$-modules to the
sections $\vx$. To this end we will need a version of cohesive local system
on the space of (relative) configurations on $C$.

The main difference from the case $C=\BA^1$ is that the local systems are
in general no more abelian one-dimensional.
In fact, the monodromy in these local systems factors through the finite
Heisenberg group.
To stress the difference we will call them
the {\em Heisenberg local systems.}

These local systems are constructed in Chapter 1.

\subsection{}
Given a $K$-tuple of $\fu$-modules (or, equivalently, factorizable
sheaves) $\{\CX_k\}$
we study the sheaf $\og(\{\CX_k\})$ on $C$ obtained by
gluing the sheaves $\{\CX_k\}$. Namely, we study its behaviour when
the curve $C$ degenerates into a stable curve $\ul\oC$ with nodes.

It appears that the sheaf $\og(\{\CX_k\})$ degenerates into a sheaf
$\ul{\og}(\{\CX_k\},\{\CR_j\})$ obtained by gluing the sheaves
$\{\CX_k\}$ and a few copies of the sheaf $\CR$: one for each node of $\ul\oC$
(Theorem ~\ref{pjat' main}).

The sheaf $\CR$ is not an object of $\FS$, but rather of $\FS^{\otimes 2}$
(or, strictly speaking, of $\Ind\FS\otimes\Ind\FS$). It corresponds to
the regular $\fu$-bimodule $\bR$ under the equivalence $\Phi$.

The Theorem ~\ref{pjat' main} is the central result of this part. Its proof occupies
Chapters 2--5. We study the degeneration away from the nodes in Chapter 2.
We study the degeneration near the nodes in Chapter 4, after we
collect the necessary information about the regular bimodule $\bR$ in
Chapter 3.

As a byproduct of geometric construction of the regular bimodule we derive
the hermitian autoduality of $\bR$ and the adjoint representation $\ad$.

\subsection{}
In Chapter 5 we investigate the global sections of the sheaf
$\og(\{\CX_k\})$. They form a local system on the moduli space of
curves with $K$ marked points and nonzero tangent vectors at these points
(strictly speaking, the local system lives on the punctured determinant
line bundle over this space). The collection of all such local systems
equips the category $\FS$ with the {\em fusion}, or {\em modular}, structure
in the terminology of ~\cite{bfm} (Theorem ~\ref{pjat' modular functor}).

Historically, first examples of modular categories appeared in the
conformal field theory (WZW models), see e.g. ~\cite{ms} and ~\cite{tuy}.
Namely, the category $\tCO_\kappa$ of integrable $\hat{\fg}$-modules of
central charge $\kappa-h$ has a natural modular structure.

As far as we know, the category $\FS$ is the first example of
nonsemisimple modular category.

\subsection{}
In Chapter 6 we study the connection between modular categories
$\tCO_\kappa$ and $\CC_\zeta$ for $\zeta=\exp(\frac{\pi\sqrt{-1}}{d\kappa})$.

It appears that the modular structure on the former category can be
reconstructed in terms of the modular structure on the latter one.

As a corollary we get a description of local systems of conformal
blocks in WZW models in arbitrary genus as natural subquotients of
some semisimple local systems of geometric origin (Theorem
~\ref{pjat' integral representations}).

The geometric local systems are
equipped with natural hermitian nondegenerate fiberwise scalar product
(being direct images of perverse sheaves which are Verdier-autodual up to
the replacement $\zeta\mapsto\zeta^{-1}$)
which gives rise to a hermitian nondegenerate scalar product on conformal
blocks in WZW models.

\subsection{} Our work on this part began 5 years ago as
an attempt to understand the
remarkable paper ~\cite{cfw}. In fact, the key ingredients --- Heisenberg
local system, and adjoint representation --- were already present in this
paper.

V.Ginzburg has drawn our attention to this paper.
D.Kazhdan's interest to our work proved extremely stimulating.

During these years we benefited a lot from discussions with many people.
The idea of Chapter 4 is due to P.Deligne. The idea to study the degeneration
of Heisenberg local system is due to B.Feigin.
We are grateful to R.Hain, J.Harris and T.Pantev who took pain of answering
our numerous questions about various line bundles on the moduli spaces.
The second author is obliged to V.Ostrik for useful discussions of adjoint
representation.



\newpage
\begin{center}
{\bf Chapter 1. Heisenberg local system}
\end{center}
\vspace{.8cm}

\bigskip

\section{Notations and statement of the main result}

\subsection{}
\label{pjat' notations}
Let $\alpha\in\BN[X],\ \alpha=\sum a_\mu\mu$.
We denote by $\supp\alpha$ the subset of $X$ consisting of all $\mu$ s.t.
$a_\mu\not=0$.
Let  $\pi:\ J\lra X$ be an unfolding of $\alpha$, that is $\sharp\pi^{-1}(\mu)=
a_\mu$ for any $\mu\in X$. As always, $\Sigma_\pi$ denotes the group of
automorphisms of $J$ preserving the fibers of $\pi$.

\subsubsection{}
The fibered product $C\times_S\ldots\times_SC$ ($J$ times) will be denoted
by $C^J$. The group $\Sigma_\pi$ acts naturally on $C^J$, and the quotient
space $C^J/\Sigma_\pi$ will be denoted by $C^\alpha$.

$\Co^J$ (resp. $\Co^\alpha$) stands for the complement to diagonals in $C^J$
(resp. in $C^\alpha$).

\subsubsection{}
\label{pjat' TO}
$TC^J$ stands for the complement to the zero sections in the relative
(over $S$) tangent bundle. So $TC^J\lra C^J$ is a $(\BC^*)^J$-torsor.
We denote by $\TCo^J$ its restriction to $\Co^J$.
The group $\Sigma_\pi$ acts freely on $\TCo^J$, and we denote the quotient
$\TCo^J/\Sigma_\pi$ by $\TCo^\alpha$.

The natural projection $\TCo^J\lra\TCo^\alpha$ will be denoted by $\pi$,
or sometimes $\pi_J$.

\subsubsection{}
Given $j\in J$ we consider the
relative tangent bundle on $C^J$
along the $j$-th coordinate $T_j$. It is a line bundle on $C^J$.
For $\mu\in\supp\alpha$ we consider the line bundle
$\underset{j\in\pi^{-1}(\mu)}{\otimes}T_j$ on $C^J$. It has a natural
$\Sigma_\pi$-equivariant structure and can be descended to the line bundle
$T_\mu$ on $C^\alpha$.

\subsubsection{}
Given $\varepsilon>0$ we denote by $D_\varepsilon$ the standard disk
of radius $\varepsilon$. If there is no danger of confusion we will omit
$\varepsilon$ from our notations and will denote $D_\varepsilon$ simply
by $D$.

The definitions of $D^J, D^\alpha, \Do^J, \Do^\alpha, \TDo^J, \TDo^\alpha$
simply copy the above definitions, and we do not reproduce them.

\subsubsection{}
\label{pjat' tree}
Given a surjection $\tau:\ J\lra K$ we consider the map
$\pi_K:\ K\lra X,\ k\mapsto\sum_{j\in\tau^{-1}(k)}\pi(j)$.
We will use the notation $\alpha_K$ for $\sum_{\mu\in X}\sharp\pi_K^{-1}(\mu)
\mu\in\BN[X]$.

We consider the following (infinite dimensional) manifold
$$\TCo^\tau=\widetilde{\TCo^K}\times\prod_{k\in K}\TDo^{\tau^{-1}(k)}$$
where $\widetilde{\TCo^K}$ is the space of analytic open embeddings
$S_K\times D\hra C$ such that the restriction to $S_K\times 0$ is just
a $K$-tuple of sections $S\lra C$. Here $S_K$ denotes the disjoint union
of $K$ copies of $S$.

We have an evident projection $p_K:\ \widetilde{\TCo^K}\lra\TCo^K$ taking
the first jet. Note that $p_K$ is a homotopy equivalence.

We denote by $q_\tau$ the natural substitution map $\TCo^\tau\lra\TCo^J$.

\subsubsection{}
\label{pjat' d}
Recall (see IV.9.1) that $\ell$ stands for $\frac{l}{2}$ in case $l$ is even,
and for $l$ in case $l$ is odd. Recall (see {\em loc. cit.}) that
$$Y_\ell:=\{\lambda\in X\ |\ \lambda\cdot\mu\in\ell\BZ\ \forall\mu\in X\}$$
We denote by ${\od}$ the cardinality of $X/Y$, and by ${\od}_\ell$
the determinant of the form $\frac{1}{\ell}?\cdot ?$
restricted to the sublattice
$Y_\ell\subset X$
Note that if $\ell$ is divisible by $d:=\max_{i\in I}d_i$
then ${\od}_\ell=\sharp(X/Y_\ell)$.
This will be the case in our applications to
conformal blocks.

To handle the general case we need to introduce some new characters.
We define $$X_\ell:=\{\mu\in X\otimes\BQ\ |\ \mu\cdot Y_\ell\in\ell\BZ\}$$
Evidently, $Y_\ell\subset X\subset X_\ell$, and $X_\ell$ is generated by
$X$ and $\{\frac{\ell}{d_i}i',i\in I\}$. So if $d|\ell$ then $X_\ell=X$ but in
general this is not necessarily the case.

Note that ${\od}_\ell=\sharp(X_\ell/Y_\ell)$.

To study the modular properties of the Heisenberg local system and the
category $\FS$ (cf. especially the Theorem ~\ref{pjat' modular Heisenberg}(b))
we will have to modify slightly the definition of the latter one, and,
correspondingly, of the category $\CC$. We start with the category $\CC$.

Consider the subalgebra $\fu'\subset\fu$ (see II.12.3) generated by
$\theta_i,\epsilon_i,\tilde{K}_i^{\pm1},\ i\in I$
(notations of {\em loc. cit.}).

We define $\CC'$ to be a category of finite dimensional $X_\ell$-graded
vector spaces $V=\oplus V_\lambda$ equipped with a structure of a left
$\fu'$-module compatible with $X_\ell$-gradings and such that
$$\tilde{K}_ix=\zeta^{\langle d_ii,\lambda\rangle}x$$
for $x\in V_\lambda,\ i\in I$.

This is well defined since $\langle d_ii,X_\ell\rangle\in\BZ$ for any $i\in I$.

We have a natural inclusion of a full subcategory $\CC\hra\CC'$.

We define the category $\FS'$ exactly as in III.5.2, just replacing all
occurences of $X$ by $X_\ell$.

We have a natural inclusion of a full bbrt
subcategory $\FS\hra\FS'$, and the
equivalence $\Phi:\ \FS\lra\CC$ extends to the same named equivalence
$\Phi:\ \FS'\lra\CC'$. The proof is the same as in III; one only has to
replace $X$ by $X_\ell$ everywhere.

Recall that in case $d|\ell$ (the case of interest for applications to
conformal blocks)
we have $X_\ell=X,\FS'=\FS,\CC'=\CC$.

{\em From now on we will restrict ourselves to the study of the categories
$\FS',\CC'$. However, in order not to scare the reader away by a bunch of
new notations we will denote them by $\FS,\CC$.}

The interested reader will readily perform substitutions in the text below.

\subsubsection{}
\label{pjat' n}
Consider a function $n:\ {X_\ell}\lra\BZ[\frac{1}{2\od}]$
such that $n(\mu+\nu)=n(\mu)+n(\nu)+
\mu\cdot\nu$. We will choose $n$ of the following form:
$$n(\mu)=\frac{1}{2}\mu\cdot\mu+\mu\cdot\nu_0$$ for some $\nu_0\in
{X_\ell}$. From now on we fix such a function $n$ and the corresponding
$\nu_0$.

Let $g$ be the genus of our relative curve $C\lra S$.

\subsubsection{Definition}
\label{pjat' tilde}
(a) We will denote the unique homomorphism
$$\BN[{X_\ell}]\lra {X_\ell}$$ identical on ${X_\ell}$ by $\alpha\mapsto\alpha^\sim$;

(b) $\alpha\in\BN[{X_\ell}]$ is
called {\em $g$-admissible} if
$\alpha^\sim=(2g-2)\nu_0$.

\subsubsection{}
\label{pjat' k}
From now on we assume that $\zeta$ is a primitive root of unity of degree $l$.
Then $\zeta=\exp(2\pi\sqrt{-1}\frac{k}{l})$ for some integer $k$ prime to $l$.
We fix $k$, and
for a rational number $q$ we define $\zeta^q:=\exp(2\pi\sqrt{-1}q\frac{k}{l})$.

\subsubsection{} Given $\alpha=\sum a_\mu\mu\in\BN[Y]$ and its unfolding
$\pi:\ J\lra Y$, we consider the following one-dimensional local system
$\CI^J$ on $\TDo^J$:

by definition, its monodromies are as follows:

around diagonals: $\zeta^{2\pi(j_1)\cdot\pi(j_2)}$;

around zero sections of tangent bundle: $\zeta^{2n(\pi(j))}$.

We define the one-dimensional local system $\CI^\alpha$ on $\TDo^\alpha$
as $\CI^\alpha:=(\pi_{J*}\CI^J)^{\Sigma_\pi,-}$ (cf. III(46)).

\subsubsection{} For a line bundle $\CL$ we denote the corresponding
$\BC^*$-torsor by $\dot{\CL}$.

\subsection{Statement of the main result}

\subsubsection{Definition} The {\em Heisenberg local system $\CH$} is the
following collection of data:

1) A local system $\CH^\alpha$ on $\TCo^\alpha$ for each admissible
$\alpha\in\BN[X_\ell]$;

2) {\em  Factorization isomorphisms:} for each $\alpha\in\BN[X_\ell]$, unfolding
$\pi:\ J\lra X_\ell$, surjection $\tau:\ J\lra K$, the following isomorphisms
are given:
$$\phi_\tau:\  q_\tau^*\pi_J^*\CH^\alpha\iso p_K^*\pi_K^*\CH^{\alpha_K}
\boxtimes\fbox{$\times$}_{k\in K}\ \CI^{\tau^{-1}(k)}$$
satisfying the usual associativity constraints.

\subsubsection{Theorem} Let $\delta\lra S$ denote the {\em determinant}
line bundle of the family $C\lra S$ (see e.g. ~\cite{km}). Then
after the base change
$$
\begin{array}{ccc}
C_\delta&\lra&C\\
\downarrow&&\downarrow\\
\dot{\delta}&\lra&S
\end{array}
$$
there exists a Heisenberg Local System $\CH$.

The dimension of $\CH$ is equal to ${\od}_\ell^g$, and
the monodromy around the zero section of $\delta$ is equal to
$(-1)^{\rk {X_\ell}}\zeta^{12\nu_0\cdot\nu_0}$.

\subsubsection{Remark}
In the case $g=1$ the line bundle $\delta^{12}$ is known to be trivial.
It is easy to see that there exists a one-dimensional local system on
$\dot\delta$ with any given monodromy around the zero section.

We will construct the Heisenberg local system $\CH$ over $S$ (as opposed
to $\dot\delta$). Lifting it to $\dot\delta$ and twisting by the above
one-dimensional systems we can obtain a Heisenberg local system
with any given scalar monodromy around the zero section.

\medskip

 The construction of the desired Heisenberg System will be given in the
rest of this Chapter.

\section{The scheme of construction}

\subsection{} First note that it suffices to construct the desired
local system $\CH$ for $\zeta=\exp(\frac{\pi\sqrt{-1}}{\ell})$.
If $k$ is prime to $l$, and $\zeta'= \exp(2\pi\sqrt{-1}\frac{k}{l})$, then
$\CH_{\zeta'}$ is obtained from $\CH_\zeta$ just by application of a
Galois automorphism of our field $\sk$. So till the end of construction we
will assume that $\zeta=\exp(\frac{\pi\sqrt{-1}}{\ell})$.

\subsection{}
In what follows everything is relative over the base $S$. To unburden the
notations we will pretend though that $C$ is an absolute curve. Thus
$H^1(C)$ stands, say, for the local system of $\BZ$-modules of rank $2g$
over $S$.

For $\alpha\in\BN[{X_\ell}]$ we introduce the following divisor
$\CCD^\alpha$ on $C^\alpha$:
$$\CCD^\alpha=\frac{{\od}_\ell^g}{\ell}(\sum_{\mu\not=\nu}\mu\cdot\nu\Delta_{\mu\nu}+
\frac{1}{2}\sum_\mu\mu\cdot\mu\Delta_{\mu\mu})$$
where $\Delta_{\mu\nu},\ \mu,\nu\in\supp\alpha$, stands for the corresponding
diagonal in $C^\alpha$. Note that for $g\geq2$ all the coefficients of the
above sum are integers. To simplify the exposition we will assume that
$g\geq2$. For the case of $g=1$ see ~\ref{pjat' elliptic}.

Given an unfolding $\pi:\ J\lra {X_\ell}$
we denote the pullback of $\CCD^\alpha$ under $\pi_J$ by $\CCD^J$; this is a divisor
on $C^J$.

We consider the (relative) Picard scheme $\Pic(C)\otimes {X_\ell}$. The group of
its connected components is naturally isomorphic to ${X_\ell}$. Each component
carries a canonical polarization $\omega$ which we presently describe.
It is a skew-symmetric bilinear form on $H_1(\Pic^0(C)\otimes {X_\ell})=
H_1(C)\otimes {X_\ell}$ equal to the tensor product of the canonical cup-product
form on $H_1(C)$ and the symmetric bilinear form $\frac{{\od}_\ell^g}{\ell}
(?\cdot?)$ on ${X_\ell}$.
Note that $\frac{{\od}_\ell^g}{\ell}(?\cdot?)$ is positive definite,
so $\omega$ is (relatively) ample.

We denote by $\AJ_\alpha:\ C^\alpha\lra \Pic(C)\otimes {X_\ell}$ the Abel-Jacobi
map $\AJ_\alpha(\sum\mu x_\mu)=\sum(x_\mu)\otimes\mu$.

The admissibility condition implies that the Abel-Jacobi map lands into
the connected component $(\Pic(C)\otimes {X_\ell})_{(2g-2)\nu_0}$ to be denoted by $A$.
Note that the projection $A\lra S$ has a canonical section $\Omega
\otimes\nu_0$, so $A$ is a genuine abelian variety, not just a torsor over one.
Here $\Omega$ denotes the (relative) canonical line bundle on $C$.

\subsubsection{Definition}
\label{pjat' g>1}
We define the following line bundle $\CL_\alpha$
on $C^\alpha$:
$$\CL_\alpha=\otimes_{\mu\in\supp(\alpha)}T_\mu^{\otimes\frac{{\od}_\ell^g}{\ell}
n(\mu)}\otimes
    \CO(\CCD^\alpha)$$

The desired construction is an easy consequence of the following Propositions:

\subsection{Proposition}
\label{pjat' L}
There is a unique line bundle $\CL$ on $A$ such
that for any $g$-admissible $\alpha$ we have $\CL_\alpha=\AJ_\alpha^*(\CL)$.
The first Chern class $c_1(\CL)=-[\omega]$.

\subsection{Proposition}
\label{pjat' H}
There exists a local system $\fH$ on $\dot{\CL}\times_S\dot\delta$ such that
$\dim\fH={\od}_\ell^g$ (see ~\ref{pjat' d});
the monodromy around the zero section of $\CL$
is equal to $\zeta^{\frac{2\ell}{{\od}_\ell^g}}$;
and the monodromy around the zero
section of $\delta$ is equal to $(-1)^{\rk {X_\ell}}\zeta^{12\nu_0\cdot\nu_0}$.

\subsection{}
\label{pjat' construction}
In the remainder of this section we derive the desired construction from
the above  Propositions.

We fix a $g$-admissible $\alpha$.
We denote by $\dot{\AJ_\alpha}$ the natural map between
the total spaces of the corresponding $\BC^*$-torsors:
$$\dot{\AJ_\alpha}:\ (\AJ_\alpha^*(\CL))\dot{}
\lra\dot{\CL}$$

By the Proposition ~\ref{pjat' L} we have an isomorphism
\begin{equation}
\label{pjat' choice}
\AJ_\alpha^*(\CL)\iso\CL_\alpha
\end{equation}
It is clear from the definition of $\CL_\alpha$ that the pullback of
$\CL_\alpha$
to $TC^\alpha$ has the canonical meromorphic
section $s_\alpha$. The restriction of this section to
$\TCo^\alpha$ does not have poles nor zeros, and
hence it defines the same named section $s_\alpha$ of the  $\BC^*$-torsor
$(\AJ_\alpha^*(\CL))\dot{}$.

We change the base to $\dot\delta$, and preserve the notations
$\AJ_\alpha,s_\alpha$ for the base change of the same named morphism and
section. By the Proposition ~\ref{pjat' H} we have the local system $\fH$
on $\dot{\CL}\times_S\dot\delta$.

We define $\CH^\alpha$ to be $s_\alpha^*\dot{\AJ_\alpha}^*\fH$ twisted
by the one-dimensional sign local system.

The proof of the above Propositions and the construction of factorization
isomorphisms will be given in the following sections.

\subsection{}
\label{pjat' elliptic}
The above construction does not work as stated in the case of elliptic curves:
the line bundle $\CO(\CCD^\alpha)$ in the Definition ~\ref{pjat' g>1} does not make
sense since the coefficients $\frac{\od_\ell}{2\ell}\mu\cdot\mu$ of the divisor
$\CCD^\alpha$ apriori may be halfintegers.

We will indicate how to carry out the construction in this case.

For any $N\in\BZ$ such that $N\frac{\od_\ell}{2\ell}\mu\cdot\mu,N\frac{\od_\ell}{\ell}
n(\mu)\in\BZ\ \forall\mu\in {X_\ell}$ we define
$$\CL_\alpha^N:=\bigotimes_{\mu\in\supp\alpha}T^{N\frac{\od_\ell}{\ell}n(\mu)}
\otimes\CO(N\CCD^\alpha).$$

We will prove the following versions of the above Propositions.

\subsubsection{}
\label{pjat' L^N}
{\bf Proposition.} There is a unique line bundle $\CL^N$ on $A$ such that
for any $g$-admissible $\alpha\in\BN[{X_\ell}]$ we have $\CL^N_\alpha=\AJ^*_\alpha
\CL^N$. The first Chern class $c_1(\CL^N)=-N[\omega]$.
\subsubsection{}
\label{pjat' H^N}
{\bf Proposition.} There exists a local system $\fH^N$ on
 $(\CL^N)\dot{}$ such that $\dim\fH^N=\od_\ell$.
The monodromy of $\fH^N$ around the zero section of $\CL^N$ is equal to
$\zeta^{\frac{2}{N\od_\ell}}$.

Moreover, for any $N|N'$ we have $\fH^N=[\frac{N'}{N}]^*\fH^{N'}$ where
$[\frac{N'}{N}]:\ \CL^N\lra\CL^{N'}$ is the map of raising to the
$(\frac{N'}{N})$-the power.

\subsubsection{}
Now $\CH^\alpha$ is defined as $s_{\alpha,N}^*\dot{\AJ}^*\fH^N$ for any
$N$ as above exactly as in ~\ref{pjat' construction}. In fact, it is enough
to take $N=2$.

\bigskip

\section{The universal line bundle}
\label{pjat' proof L}

In this section we will give a proof of the Propositions ~\ref{pjat' L}
and ~\ref{pjat' L^N}.

\subsection{}
\label{pjat' Lambda}
First we formulate a certain generalization. Suppose given a free $\BZ$-module
$\Lambda$ of finite rank with an even symmetric bilinear pairing $(,):\
\Lambda\times\Lambda\lra\BZ$. We fix an element $\nu\in\Lambda$, and
a function $\ob:\ \Lambda\lra\BZ,\ \ob(\lambda):=\frac{1}{2}(\lambda,\lambda)
+(\lambda,\nu)$.

\subsubsection{}
For $\alpha\in\BN[\Lambda]$ we introduce the following divisor
$\CCD^\alpha$ on $C^\alpha$:
$$\CCD^\alpha=\sum_{\mu\not=\lambda}
(\mu,\lambda)\Delta_{\mu\lambda}+
\frac{1}{2}\sum_\mu(\mu,\mu)\Delta_{\mu\mu}$$
where $\Delta_{\mu\lambda},\ \mu,\lambda\in\supp\alpha$,
stands for the corresponding
diagonal in $C^\alpha$. Note that all the coefficients of the
above sum are integers.

Given an unfolding $\pi:\ J\lra \Lambda$
we denote the pullback of $\CCD^\alpha$ under $\pi_J$ by $\CCD^J$; this is a divisor
on $C^J$.

We consider the (relative) Picard scheme $\Pic(C)\otimes \Lambda$. The group of
its connected components is naturally isomorphic to $\Lambda$. Each component
carries a canonical polarization $\omega$ which we presently describe.
It is a skew-symmetric bilinear form on $H_1(\Pic^0(C)\otimes \Lambda)=
H_1(C)\otimes \Lambda$ equal to the tensor product of the canonical cup-product
form on $H_1(C)$ and the symmetric bilinear form $(,)$ on $\Lambda$.

We denote by $\AJ_\alpha:\ C^\alpha\lra \Pic(C)\otimes\Lambda$ the Abel-Jacobi
map $\AJ_\alpha(\sum\mu x_\mu)=\sum(x_\mu)\otimes\mu$.

The admissibility condition implies that the Abel-Jacobi map lands into
the connected component $(\Pic(C)\otimes \Lambda)_{(2g-2)\nu}$
to be denoted by $A_\Lambda$.

Note that the projection $A_\Lambda\lra S$ has a canonical section $\Omega
\otimes\nu$, so $A_\Lambda$ is a genuine abelian variety,
not just a torsor over one.
Here $\Omega$ denotes the (relative) canonical line bundle on $C$.

\subsubsection{Definition}
\label{pjat' L alpha}
We define the following line bundle $\CL_\alpha$
on $C^\alpha$:
$$\CL_\alpha=\otimes_{\mu\in\supp(\alpha)}T_\mu^{\otimes \ob(\mu)}\otimes
    \CO(\CCD^\alpha)$$

An element $\alpha=\sum a_\lambda\lambda\in\BN[\Lambda]$ (formal sum) is called
$g$-admissible if $\alpha^\sim=(2g-2)\nu$.

Now we are able to formulate the Proposition generalizing ~\ref{pjat' L}.

\subsection{Proposition}
\label{pjat' L(Lambda)}
There is a unique line bundle $\CL(\Lambda,\ (,),\ \nu)$ on $A_\Lambda$ such
that for any $g$-admissible $\alpha$ we have $\CL_\alpha=\AJ_\alpha^*
(\CL(\Lambda,\ (,),\ \nu))$.
The first Chern class $c_1(\CL(\Lambda,\ (,),\ \nu))=-[\omega]$.

\subsection{} We start the proof of the Proposition ~\ref{pjat' L(Lambda)}
with the following Lemma. Let $\pi_J:\ J\lra \Lambda$
be an unfolding of $\alpha$. Let $\tau:\ J\lra K$ be a surjection, and
$\pi_K:\ K\lra \Lambda$ be an unfolding of $\alpha_K$ as in ~\ref{pjat' tree}.
To simplify the notations we will denote $\alpha_K$ by $\beta$.

Let $\sigma_\tau$ denote the natural (``diagonal'') embedding
$C^K\hra C^J$.

\subsubsection{Lemma}
\label{pjat' compatibility}
There is a canonical isomorphism
$$\sigma_\tau^*\pi_J^*(\CL_\alpha)=\pi_K^*(\CL_\beta)$$

{\em Proof.} It is enough to prove the Lemma in the case $|J|=|K|+1$.
So we fix $i,j\in J$ such that $\tau(i)=\tau(j)=k$, and we denote
$\pi_J(i)$ by $\mu_i$, and $\pi_J(j)$ by $\mu_j$.

We have $\pi_J^*(\CCD^\alpha)=\mu_i\cdot\mu_j\Delta_{ij}+\CCD'$ where
$\Delta_{ij}\not\subset\supp(\CCD')$.

Moreover, it is clear that $\CCD'\cap\Delta_{ij}=\sigma_\tau(\pi_K^*(\CCD^\beta))$,
and hence $$\sigma_\tau^*\pi_\tau^*(\CO(\CCD'))=\pi_K^*(\CO(\CCD^\beta))$$

On the other hand, for any smooth divisor $\CCD$ we have a canonical isomorphism
$\CO(\CCD)|_\CCD= \CN_\CCD$ (the normal bundle).

In particular, we have $\sigma_\tau^*(\CO(\Delta_{ij}))=\sigma_\tau^*(T_i)=
\sigma_\tau^*(T_j)=T_k$.

Thus $\sigma_\tau^*(T_i^{\ob(\mu_i)}\otimes T_j^{\ob(\mu_j)}\otimes
\CO((\mu_i,\mu_j)\Delta_{ij}))=T_k^{\ob(\mu_i)+\ob(\mu_j)+(\mu_i,\mu_j)}=
T_k^{\ob(\mu_i+\mu_j)}=T_k^{\ob(\pi(k))}$.

Finally, if $m\not=i,j$, then evidently $\sigma_\tau^*(T_m)=T_{\tau(m)}$.

Putting all this together we obtain the statement of the Lemma. $\Box$

\subsection{}
To prove the Proposition we have to check a necessary condition: that the
first Chern class of $\CL_\alpha$ is a pullback of some cohomology class
on $A_\Lambda$ under the Abel-Jacobi map $\AJ_\alpha$.
This is the subject of the
following Lemma.

\subsubsection{Lemma}
\label{pjat' chern class}
If $\alpha$ is admissible, then
$c_1(\CL_\alpha)=\AJ_\alpha^*(-[\omega])$.

{\em Proof.} Let us choose an unfolding $\pi:\ J\lra \Lambda$.
We denote by $\pi_J$
the corresponding projection $C^J\lra C^\alpha$.

It is enough to prove that the pullback of the both sides to $C^J$ under
$\pi_J$ coincide.

We introduce the following family of 2-cycles in $C^J$. For 1-cycles $a,b$
in $C$, and $i\not=j\in J$, we denote by $a_i\times b_j$ the following
product cycle: the $i$-th coordinate runs along the cycle $a$, the
$j$-th coordinate runs along the cycle $b$, all the rest coordinates are
fixed. The homology class of $a_i\times b_j$ depends only on $i,j$ and
the classes of $a$ and $b$.

We denote by $f_i$ the following 2-cycle: the $i$-th coordinate runs along
the fundamental cycle of $C$, all the rest coordinates are fixed. The
homology class of $f_i$ depends on $i$ only.

It suffices to check that the pairings of both sides of ~\ref{pjat' chern class}
with this family of cycles coincide.

We have:
\begin{equation}
\langle-\omega,a_i\times b_j\rangle=(\pi(i),\pi(j))a\cap b=
(a_i\times b_j)\cap((\pi(i),\pi(j))\Delta_{ij})=
\langle c_1(\CL_\alpha),a_i\times b_j\rangle;
\end{equation}
\begin{equation}
\label{pjat' 1}
\langle-\omega,f_i\rangle=-g(\pi(i),\pi(j));
\end{equation}
\begin{equation}
\label{pjat' 2}
\langle c_1(\CL_\alpha),f_i\rangle=(2-2g)\ob(\pi(i))+f_i\cap\pi^*(\CCD^\alpha)=
(2-2g)\ob(\pi(i))+(\pi(i),(\sum_{j\in J}\pi(j)-\pi(i))).
\end{equation}
To assure the equality of (\ref{pjat' 1}) and (\ref{pjat' 2}) we must have
\begin{equation}
(1-g)(\pi(i),\pi(i))=(2-2g)\ob(\pi(i))+(\pi(i),\sum_{j\in J}\pi(j)),
\end{equation}
that is,
\begin{equation}
(1-g)((\pi(i),\pi(i))-2\ob(\pi(i)))=(\pi(i),\sum_{j\in J}\pi(j))
\end{equation}
which is precisely the admissibility condition. $\Box$

\subsection{}
Let us choose a basis $I$ of $\Lambda$.

Using the Lemma ~\ref{pjat' compatibility} we see that it suffices to prove
the Proposition for $\alpha$ of a particular kind, namely
$$\alpha=\sum_{i\in I}a_ii+\sum_{i\in I}a_{-i}(-i)$$
where all the positive integers $a_i,a_{-i}$ are big enough.

\subsubsection{}
We define $\alpha_+:=\sum_{i\in I}a_ii$,
and $\alpha_-:=\sum_{i\in I}a_{-i}(-i)$.

We consider the following Abel-Jacobi maps:
$$\AJ_+:\ C^{\alpha_+}\lra(\Pic(C)\otimes \Lambda)_{\alpha_+}=:A_+,$$
and
$$\AJ_-:\ C^{\alpha_-}\lra(\Pic(C)\otimes \Lambda)_{\alpha_-}=:A_-.$$
We have
$$\AJ_\alpha=m\circ(\AJ_+\times\AJ_-)$$
where $m:\ A_+\times_SA_-\lra A_\Lambda$ is the addition map.

If all the $a_{\pm i}$ are bigger than $g$, then the map
$$(\AJ_+\times\AJ_-)^*:\ \Pic^0(A_+\times_SA_-)\lra\Pic^0(C^{\alpha_+}
\times_SC^{\alpha_-})$$ is an isomorphism, and the induced map on the whole
Picard groups is an inclusion.

According to the Lemma ~\ref{pjat' chern class}, $c_1(\CL_\alpha)=\AJ_\alpha^*
(-[\omega])=(\AJ_+\times\AJ_-)^*(m^*(-[\omega]))$. So we deduce that there
exists a unique line bundle $\CL'$ on $A_+\times_SA_-$ such that
$\CL_\alpha=(\AJ_+\times\AJ_-)^*\CL'$.

It remains to show that $\CL'=m^*\CL$ for some line bundle $\CL$ on $A_\Lambda$
(necessarily uniquely defined). To this end it is enough to verify that
the restrictions of $\CL'$ to the fibers of $m$ are trivial line bundles.

\subsubsection{}
We choose an unfolding $J$ of $\alpha$. We choose a surjection $\tau:\
J\lra K$ with the following property: $K=K_0\sqcup K_1;\ \tau$ is
one-to-one over $K_1$, and for any $k\in K_0$ we have
$\tau^{-1}(k)=\{i,-i\}$ for some $i\in I$.

Moreover, we assume that $\alpha$ is big enough so that for each $i\in I$
both $i$ and $-i$ appear at least $g$ times in $\tau^{-1}$ of both $K_0$
and $K_1$.

Recall that $\sigma_\tau$ stands for the diagonal embedding
$C^K\hra C^J$.

It is clear that $C^K=C^{K_0}\times_SC^{K_1}$, and the Abel-Jacobi map
$\AJ_J\circ\sigma_\tau:\ C^K\lra A_\Lambda$
factors through the projection onto the second factor.

Fix a point $a\in A_\Lambda$. Let us choose $p\in C^{K_1}$ such that
$\AJ_J(\sigma_\tau(C^{K_0}\times p))=a$.
Then $(\AJ_+\times\AJ_-)\circ\sigma_\tau$
maps $C^{K_0}\times p$ to the fiber $m^{-1}(a)$.

It is easy to see that the induced maps on the Picard groups
$((\AJ_+\times\AJ_-)\circ\sigma_\tau)^*:\ \Pic(m^{-1}(a))\lra
\Pic(C^{K_0})$ is injective.

Thus we only have to check that
$$((\AJ_+\times\AJ_-)\circ\sigma_\tau)^*\CL'|_{C^{K_0}\times p}=
(\sigma_\tau\circ\pi_J)^*\CL_\alpha|_{C^{K_0}\times p}$$ is trivial.

According to the Lemma ~\ref{pjat' compatibility}, this line bundle is equal to
$\pi_K^*(\CL_{\alpha_K})$. It is clear from the definition that
$\CL_{\alpha_K}$ is lifted from the projection to the factors carrying
nonzero weights. In particular, the restriction of $\CL_{\alpha_K}$ to a
fiber of this projection is trivial.

This completes the proof of the Proposition ~\ref{pjat' L(Lambda)}.

\subsection{}
To prove the Proposition ~\ref{pjat' L} it suffices to apply Proposition
~\ref{pjat' L(Lambda)} to the case $\Lambda={X_\ell};\ (?,?)=\frac{{\od}_\ell^g}{\ell}?\cdot?;\
\nu=\nu_0;\ \ob(\lambda)=\frac{{\od}_\ell^g}{\ell}n(\lambda)$.

To prove the Proposition ~\ref{pjat' L^N} we take $\Lambda={X_\ell};\ (?,?)=
N\frac{{\od}_\ell^g}{\ell}?\cdot?;\
\nu=\nu_0;\ \ob(\lambda)=N\frac{{\od}_\ell^g}{\ell}n(\lambda).\ \Box$

\bigskip

\section{The universal local system}
\label{pjat' proof H}

In this section we will give a proof of the Propositions ~\ref{pjat' H}
and ~\ref{pjat' H^N}.

\subsection{}
\label{pjat' nu}
Recall the notations of ~\ref{pjat' Lambda}.
Let us take $\Lambda={X_\ell}\oplus Y_\ell;\
((x_1,y_1),(x_2,y_2))=
-\frac{1}{\ell}(y_1\cdot x_2+y_2\cdot x_1+y_1\cdot y_2);\
\nu=(\nu_0,0)$.

We denote $A_\Lambda$ by $A'$. We have an obvious projection $\pr_1:\
A'\lra A$. We denote the line bundle $\CL(\Lambda,\ (,),\ \nu)$ (notations
of ~\ref{pjat' L(Lambda)}) by $\CL'$.

\subsection{Theorem}
\label{pjat' L'}
(a) $\CL'$ is relatively ample with respect to $\pr_1$.

(b) The direct image
$\CE:=\pr_{1*}\CL'$ is a locally free sheaf of rank ${\od}_\ell^g$.

(c) In the situation of ~\ref{pjat' H} (that is, $g>1$)
we have $\det(\CE)=\CL\otimes(p^*\delta)^{{\od}_\ell^g(-\frac{1}{2}\rk {X_\ell}+6\frac
{\nu_0\cdot\nu_0}{\ell})}$, where $p$ stands for the projection $A\lra S$.

(d) Assume $g=1$. Then in the notations of ~\ref{pjat' elliptic} we have
for any $N$ as in {\em loc. cit.}
$(\det(\CE))^{\otimes N}=\CL^N\otimes p^*\delta^{\otimes iN}$ for some
$i\in\BZ$.

\subsubsection{}
\label{pjat' Beilinson Kazhdan}
Let us construct the desired local system assuming the Theorem
~\ref{pjat' L'}.
By the virtue of ~\cite{bk}, Corollary 3.4 and Corollary 4.2,
$\CE$ is naturaly equipped with the flat projective connection.
Hence its lifting to  $(\det(\CE))\dot{}$ carries
a flat connection with scalar monodromy around the zero section equal to
$\exp (\frac{2\pi\sqrt{-1}}{{\od}_\ell^g})$.
Isomorphism ~\ref{pjat' L'}.(c) yields the
map $m:\dot{\CL} \times_S\dot\delta\lra
 (\det(\CE))\dot{} $, where $m(\lambda , t)=\lambda
\otimes t^{{\od}_\ell^g(-\frac{1}{2}\rk {X_\ell}+6\frac{\nu_0\cdot\nu_0}{\ell})}$.
  It is clear that $m^*\CE$ is a locally free sheaf
with flat connection, whose monodromy around the 0-section of
$\CL$ is equal to $\exp (\frac{2\pi\sqrt{-1}}{{\od}_\ell^g})=
\zeta^{\frac{2\ell}{{\od}_\ell^g}}$,
and monodromy around the zero section of $\delta$ is equal to
$(-1)^{\rk {X_\ell}}\zeta^{12\nu_0\cdot\nu_0}$.
This proves the Proposition ~\ref{pjat' H}.

The proof of the Proposition ~\ref{pjat' H^N} is even simpler.
Note that for any $M$ the lifting of $\CE$ to $((\det(\CE))^{\otimes M})
\dot{}$ carries a flat connection with scalar monodromy around the zero
section equal to $\exp(\frac{2\pi\sqrt{-1}}{M\od_\ell})$.

Now ~\ref{pjat' L'}(d) implies the isomorphism $(\det(\CE))^{\otimes 12N}
=(\CL^N)^{\otimes 12}$ for any $N$ as in ~\ref{pjat' elliptic} (recall that
$\delta^{\otimes12}$ is trivial).

So for any $N$ we can define $\fH^N$ to be $m^*\CE$ where $m:\ \CL^N\lra
(\det(\CE))^{\otimes 12N}$ is
the composition of raising to the 12-th power and the above isomorphism.
$\Box$

\subsection{}
We now proceed with the proof of the Theorem ~\ref{pjat' L'}.
The statements (a) and (b) follow immediately from the Proposition
~\ref{pjat' L(Lambda)} and the Riemann-Roch formula for abelian varieties
(see e.g. ~\cite{m}).

To prove (c) we need two auxilliary Propositions.

As usually, it is more convenient to work in greater
generality. Thus suppose given $\Lambda,\ (,),\ \nu$ as in ~\ref{pjat' Lambda}.

We will denote the canonical section $\Omega\otimes\nu$ of
$A_\Lambda=(\Pic(C)\otimes \Lambda)_{(2g-2)\nu}$ by $s$.
This section identifies $A_\Lambda$ with $\Pic^0(C)\otimes\Lambda$, and
we will use this identification in what follows.

\subsection{Proposition}
\label{pjat' 6}
$s^*(\CL (\Lambda,\ (,),\ \nu)=\delta^{-6 (\nu,\nu)}$
(notations of ~\ref{pjat' L(Lambda)}).

\subsection{Corollary}
\label{pjat' 8 1/2}
$\CL(\Lambda,\ (,),\ \nu)=\CL(\Lambda,\ (,),\ 0)\otimes
p^*\delta^{-6(\nu,\nu)}$.

{\em Proof.}
Let first $g=1$. Then note that the condition of $g$-admissibility
reads $\alpha^\sim=0$ independently of $\nu$. To stress the dependence
of $\CL^\alpha$ on $\nu$ we will include $\nu$ as the subindex for a
moment.

We have $T_\lambda=p^*(\delta^{-1})$ for any $\lambda\in\Lambda$.

For any $1$-admissible   $\alpha$, and any $\nu$ we have
$$\CL^\alpha_\nu=\CL^\alpha_0\otimes\bigotimes_{\lambda\in\supp\alpha}
T_\lambda^{(\lambda,\nu)}=\CL^\alpha_0\otimes p^*(\delta^{-(\alpha^\sim,\nu)})=
\CL^\alpha_0$$

It follows that $\CL(\Lambda,\ (,),\ \nu)=\CL(\Lambda,\ (,),\ 0)$ for
any $\nu$. On the other hand $6(\nu,\nu)\in 12\BZ$, and $\delta^{\otimes 12}$
is known to be trivial.

This takes care of the case $g=1$.

Assume now $g>1$.

In view of the Proposition ~\ref{pjat' 6} it is enough to show that
$\CL(\Lambda,\ (,),\ \nu)=\CL(\Lambda,\ (,), 0)\otimes p^*\Xi$ for some
line bundle $\Xi$ on $S$. So we can assume that $S$ is a point.

Let us choose a nonzero holomorphic form $\sigma$ with zero divisor $(\sigma)$
on $C$. Then, for any
$\alpha\in\BN[\Lambda]$, we have an embedding
$$\iota_\sigma: C^\alpha\hra C^{\alpha+(2g-2)\nu};\
x\mapsto x+(\sigma)\otimes\nu.$$

Recall the Definition ~\ref{pjat' L alpha} of $g$-admissibility. Since $g$ is
fixed, while $\nu$ varies, till the end of the proof we will replace this
term by {\em $\nu$-admissibility}.

Note that if $\alpha$ is $0$-admissible, then $\alpha+(2g-2)\nu$ is
$\nu$-admissible.

It is clear that for $0$-admissible $\alpha$ we have
$$\iota_\sigma^*\CL^{\alpha+(2g-2)\nu}=
\bigotimes_{\mu\in\supp\alpha}T_\mu^{\ob_\nu(\mu)}\otimes\CO(\CCD^\alpha)
\otimes\bigotimes_{\mu\in\supp\alpha}T_\mu^{-(\mu,\nu)}=$$
$$=\bigotimes_{\mu\in\supp\alpha}T_\mu^{\frac{1}{2}(\mu,\mu)}
\otimes\CO(\CCD^\alpha)=
\bigotimes_{\mu\in\supp\alpha}T_\mu^{\ob_0(\mu)}\otimes\CO(\CCD^\alpha)=
\CL^\alpha$$
(notations of ~\ref{pjat' Lambda}).

Now the Corollary follows from the definition
(and uniqueness) of $\CL(\Lambda,\ (,),\ \nu)$ (Proposition ~\ref{pjat' L(Lambda)}).
$\Box$

\subsection{Proposition}
\label{pjat' 1/2}
Assume $g>1$. Then
$[Rp_{\Lambda*}(\CL (\Lambda,\ (,),\ 0))]^1
=\delta^{-\frac{1}{2}d^g\rk\Lambda}$
where $d$ stands for $\det(-(,)\ )$, while $p_\Lambda$ denotes the projection
$A_\Lambda\lra S$,
and $[Rp_{\Lambda*}(?)]^1$ denotes the
determinant of the complex $Rp_{\Lambda*}(?)$ (see e.g. ~\cite{km}).

Note that $d$ equals the Euler characteristic of
$Rp_{\Lambda*}(\CL (\Lambda,\ (,),\ 0))$.

\subsection{}
\label{pjat' 57}
Let us derive the Theorem ~\ref{pjat' L'}(c)
assuming the above Propositions.

Note that the scalar product $((\nu_0,0),(\nu_0,0))=0$, and hence the
Corollary ~\ref{pjat' 8 1/2} implies that $\CL'$ does not depend on $\nu_0$.
It follows from {\em loc. cit.} that the RHS of ~\ref{pjat' L'}(c) does not
depend on $\nu_0$ either.

We will assume that $\nu_0=0$ until the end of the proof.

For $M\in\BZ$ and an abelian group $G$ let $[M]:\ G\lra G$ denote the
multiplication by $M$. We choose $M$ such that $M{X_\ell}\subset Y_\ell$.

\subsubsection{Lemma}
\label{pjat' 571}
$[M]^*\det(\CE)=[M]^*
(\CL\otimes(p^*\delta)^{{\od}_\ell^g(-\frac{1}{2}\rk {X_\ell}+6\frac
{\nu_0\cdot\nu_0}{\ell})})$.

{\em Proof.} Let us define $+_M:\ {X_\ell}\oplus Y_\ell\lra Y_\ell$ as
$+_M(x,y)=Mx+y$.

We have an equality of quadratic forms $$([M]\times\id)^*(,)=(+_M)^*
(-\frac{1}{\ell}?\cdot?)+\pr_1^*(\frac{M^2}{\ell}?\cdot?)$$
(note that all the quadratic forms involved are integer valued and even).

It implies the equality of line bundles:
$$([M]\times\id)^*\CL'=(+_M\otimes\id)^*(\CL(Y_\ell,-\frac{1}{\ell}?\cdot?,0)
\otimes\pr_1^*\CL({X_\ell},\frac{M^2}{\ell}?\cdot?,0))$$
Note that $(\pr_1,+_M):\ {X_\ell}\oplus Y_\ell\lra {X_\ell}\oplus Y_\ell$ is an automorphism.
Hence

\begin{multline}
\label{pjat' mzvezdochka}
[M]^*\det(\CE)=\det(\pr_{1*}(\pr_2^*\CL(Y_\ell,-\frac{1}{\ell}?\cdot?,0)\otimes
\pr_1^*\CL({X_\ell},\frac{M^2}{\ell}?\cdot?,0))=\\
=\det(\pr_{1*}\pr_2^*\CL(Y_\ell,-\frac{1}{\ell}?\cdot?,0))\otimes
\CL({X_\ell},\frac{M^2}{\ell}?\cdot?,0)^{\od_\ell^g}
\end{multline}

by the projection formula (note that $\od_\ell^g$ is the Euler characteristic
of $\pr_{1*}\pr_2^*\CL(Y_\ell,-\frac{1}{\ell}?\cdot?,0)$).

The first tensor multiple is equal to
$$p_A^*\det((p_{\Pic^0(C)\otimes Y_\ell})_*\CL(Y_\ell,-\frac{1}{\ell}?\cdot?,0))=
p_A^*\delta^{-\frac{1}{2}\od_\ell^g\rk Y_\ell}$$
by the Proposition ~\ref{pjat' 1/2}.

The second tensor multiple is identified with
$$\CL({X_\ell},\frac{\od_\ell^gM^2}{\ell}?\cdot?,0)=
[M]^*\CL({X_\ell},\frac{\od_\ell^g}{\ell}?\cdot?,0)$$
This completes the proof of the Lemma. $\Box$

\subsubsection{}
\label{pjat' above}
Now the difference of the LHS and RHS of ~\ref{pjat' L'}(c)
is a line bundle $\Xi$ on $\Pic^0(C)\otimes {X_\ell}$ which is defined functorially
with respect to $S$. Moreover, we have seen that $[M]^*\Xi$ is trivial.
This implies that $\Xi$ is trivial itself.

In effect, it defines a section of $\Pic^0(C)_M\otimes {X_\ell}$ (where the
subscript $_M$ stands for $M$-torsion).

But $\Pic^0(C)_M\otimes {X_\ell}=R^1p_{C*}(\ul{\BZ/M\BZ})\otimes {X_\ell}$, and it is
well known that $R^1p_{C*}(\ul{\BZ/M\BZ})$ does not have nonzero
functorial sections. Hence the above section vanishes, so the difference
of the LHS and RHS is a line bundle lifted from the base $S$.

But we know that the restrictions of the LHS and RHS to the zero section
of $\Pic^0(C)\otimes {X_\ell}$ coincide.

This completes the proof of ~\ref{pjat' L'}(c).

\subsubsection{}
Now we will prove ~\ref{pjat' L'}(d).

According to the proof of
Corollary ~\ref{pjat' 8 1/2} both sides are independent of
$\nu$, so we put $\nu=0$.

The formula ~(\ref{pjat' mzvezdochka}) and the argument just above it applies
to the case of elliptic curve as well provided $N|M$.
The first tensor multiple of
~(\ref{pjat' mzvezdochka}) is obviously lifted from the base. Since the
Picard group of the moduli space $\CM_1$ is generated by $\delta$,
{\em loc. cit.} implies that
$$[M]^*\det(\CE)=p^*(\delta^i)\otimes\CL({X_\ell},\frac{\od_\ell M^2}{\ell}?\cdot?,0)$$
for some $i$.
Hence $$[M]^*(\det(\CE))^{\otimes N}
=p^*(\delta^{iN})\otimes\CL({X_\ell},\frac{N\od_\ell M^2}{\ell}?\cdot?,0)=
[M]^*(p^*(\delta^{iN})\otimes\CL^N)$$

Exactly as in ~\ref{pjat' above} this implies that
$$(\det(\CE))^{\otimes N}=p^*(\delta^{iN})\otimes\CL^N$$
This completes the proof of ~\ref{pjat' L'}(d). $\Box$

\subsection{}
We start with the proof of Proposition ~\ref{pjat' 6}.

If $g=1$ then both sides are trivial. Indeed, the RHS is trivial
since $6(\nu,\nu)\in 12\BZ$, and $\delta^{\otimes 12}$ is trivial.

To see that the LHS is trivial it suffices to consider the case $\alpha=
\{0\}$ (neutral element of $\Lambda$ with multiplicity one).

\medskip

>From now on we assume that $g>1$.

Put $P:= {\Bbb P}((p_C)_*(\Omega_{C/S}))$ where $p_C$ denotes the projection
$C\lra S$.
We have the natural inclusion $i: P \hookrightarrow C^{(2g-2)}$.
Let $\alpha \in {\Bbb N}[\Lambda]$ be the multiset consisting of
$(2g-2)$ copies of $\nu$. We will identify $C^{(2g-2)}$ with $C^\alpha$.
Of course $\AJ_\alpha \circ i$ maps $P$ to the image of $s$.
So it suffices to prove that

\begin{equation}
(\AJ_\alpha \circ i )^*(\CL (\Lambda,\ (,),\ \nu))= p_P ^*(\delta^
{-6(\nu,\nu)})
\end{equation}

where $p_P$ denotes the projection $P\lra S$.

Using the fact that $\ob(\nu)=\frac{3}{2}(\nu,\nu)$ we deduce

\begin{equation}
(\AJ_\alpha \circ i )^*(\CL (\Lambda,\ (,),\ \nu))= i^*(\CL _\alpha)=
i^* (S^{2g-2}(\Omega^{-\frac{3}{2}(\nu,\nu)}) (\frac{1}{2}(\nu,\nu)\Delta)).
\end{equation}

Here $\Delta \subset C^{(n)}$ is the diagonal divisor, and for an invertible
sheaf ${\CF}$ on $C$ we denote by  $S^n({\CF})$ the invertible sheaf
$\pi _*({\CF}^
{\boxtimes n})^{\Sigma _n}$ on $C^{(n)}$ (where $\pi: C^n \rightarrow
C^{(n)}$ is the projection).

Let us introduce some more notation. Let $I\subset C\times
C^{(n)}$ be the incidence divisor. We will denote
the projection of $C\times C^{(n)}$ to the $i$-th factor by $\pr_i$ till the
further notice. This will not cause any confusion with our previous use
of the notation $\pr_1$. For an invertible
sheaf ${\CF}$ on $C$ we define a rank $n$  locally free sheaf  ${\CF}^{(n)}$
 on $C^{(n)}$ by: ${\CF}^{(n)}:= \pr_{2*}(\CO_I \otimes\pr_1^*({\CF}))$.

\subsubsection{Lemma}
(a) The map ${\CF}\mapsto S^n({\CF})$ defines a homomorphism
from $\Pic(C)$ to $\Pic(C^{(n)})$;

(b) For any line bundle ${\CF}$ on $C$ we have
$((\pi _*({\CF}^{\boxtimes n}))^{\Sigma _n ,-})^{\otimes2} =
S^n({\CF}^{\otimes 2})(-\Delta)$;

(c) $\det ({\CF}^{(n)})= (\pi _*({\CF}^{\boxtimes n}))^{\Sigma _n, -}$.

{\em Proof.} (a) is clear.

(b) We obviously have a morphism from the LHS to $S^n({\CF}^2)$.
We have to check
that the cokernel is supported on $\Delta$, and its stalk at the generic
point of $\Delta $ is 1-dimensional. For this we can assume that
 $C={\Bbb A}^1$, and
${\CF}=\CO$, and check the statement directly.

(c) Let us first construct a morphism from the LHS to the RHS.
Let $\tilde I \subset C^{n+1}$ be the union of diagonals $x_1=x_i$;
let $\tilde {\pr}_2$ (resp. $\tilde {\pr}_1$) be the projection of $C^{n+1}$
onto the last $n$ (resp. first) coordinates. Then

\begin{equation}
\label{pjat' one}
\pi ^*(\det ({\CF}^{(n)}))=
\det (\tilde {\pr}_{2*}(\tilde {\pr}_1^*{\CF}\otimes
\CO_{\tilde I}))
\end{equation}

 because the square is Cartesian.

We also have the arrow: $\CO_{\tilde I}\rightarrow \oplus _{i=2,..,n+1}
\CO _{ \{x_1=x_i\} }$ which yields the morphism

\begin{equation}
\label{pjat' two}
\det (\tilde  {\pr}_{2*}(\tilde {\pr}_1^* {\CF} \otimes \CO_{\tilde I})
\lra \det (\tilde {\pr}_{2*} (\tilde {\pr}_1^*{\CF}\otimes(\oplus \CO
_{\{x_1=x_i\}})))= {\CF}^{\boxtimes n}
\end{equation}

It is clear that the morphism ~(\ref{pjat' two})
anticommutes with the action
of $\Sigma _n$, hence, by ~(\ref{pjat' one}),
it defines the desired arrow.

It is easy to see that the image of the monomorphism ~(\ref{pjat' two}) is
${\CF}^{\boxtimes n}(-D )\subset {\CF}^{\boxtimes n}$, where $D$
is the union of all diagonals. But the latter inclusion
induces isomorphism $$(\pi _* ({\CF}^{\boxtimes n}(-D)))^{\Sigma _n,-}=
(\pi _*({\CF}^{\boxtimes n}))^{\Sigma _n, -}$$ This completes the proof
of the Lemma. $\Box$

\subsubsection{}
We return to the proof of Proposition ~\ref{pjat' 6}.
It suffices to prove it
assuming that the line bundle $\Omega^{1/2}$ exists over $S$. Indeed,
let $\CM_g$ denote the moduli space of curves of genus $g$. Then it is
known that $\Pic(\CM_g)$ is torsion free (see e.g. ~\cite{m1}, Lemma 5.14).
Hence it injects into the Picard group
of the moduli space of curves with $\theta$-characteristics.

We have
$$i^*(\CL _\alpha)= i^*(S^{2g-2}(\Omega^{-\frac{3}{2}(\nu,\nu)})(\frac{1}{2}
(\nu,\nu)\Delta)) = i^*(S^{2g-2}(\Omega^3)(-\Delta))^{-\frac{1}{2}
(\nu,\nu)}=$$
$$i^*({\operatorname{det}}^{\otimes 2}
((\Omega^{\frac{3}{2}})^{(2g-2)}))^{-\frac{1}{2}(\nu,\nu)}
=i^*(\det((\Omega^{\frac{3}{2}})^{(2g-2)}))^{-(\nu,\nu)}$$

So the Proposition follows from

\subsubsection{Lemma}
$\det (i^* (\Omega^{\frac{3}{2}})^{(2g-2)})=p_P^*\delta^6$.

{\em Proof.} Consider the exact sequence of sheaves on $C\times C^{(n)}$:
$$0 \lra {\CF} \lra\pr_1^*\Omega^{\frac{3}{2}} \lra\pr_1^*\Omega^{\frac{3}{2}}
\otimes \CO_I\lra 0$$ where ${\CF}$ denotes the kernel.
Let now (until the end of the proof)
$\pr_i$ denote the projections of $C\times P$ to the respective factors.
We see that $$(id \times i)^*({\CF})=
\pr_1^*\Omega^{\frac{1}{2}}\otimes \pr_2^*(?)$$
We apply the functor $R\pr_{2*}\circ
(id\times i)^*$ to this exact sequence (note that $[R\pr_{2*}
(\Omega^{\frac{1}{2}}
\otimes \pr_2 ^*(?))]^1$ does not depend on $?$ since the Euler characteristic
of $\Omega^{\frac{1}{2}}$ is $0$).

Using Mumford's formula (see ~\cite{m1}, Theorem 5.10) we get
$$\det (i^* (\Omega^{\frac{3}{2}})^{(2g-2)})=
[R\pr_{2*} (\pr_1^*\Omega^{\frac{3}{2}}\otimes i^*(\CO _I))]^1=$$
$$=[R\pr_{2*} (\pr_1^*\Omega^{\frac{3}{2}})]^1 -
[R\pr_{2*} (\pr_1^*\Omega^{\frac{1}{2}})]^1=$$
$$=[6((\frac{3}{2})^2 -\frac{3}{2})+1 - (6((\frac{1}{2})^2 -\frac{1}{2})+1)]
p_P^*\delta = 6 p_P^*\delta$$

This completes the proof of the Lemma along with the Proposition ~\ref{pjat' 6}.

\subsection{}
We start the proof of the Proposition ~\ref{pjat' 1/2}.

\subsubsection{Lemma}
\label{pjat' isogeny}
Let ${\fA}\rightarrow S$ be a family of abelian varieties
over a smooth base $S$. Let $\fa:\ {\fA}' \lra {\fA}$ be an isogeny. Let
$\CL$ be a line bundle over ${\fA}$. Then $\deg (\fa) [Rp_{\fA*}\CL]^1-
[Rp_{\fA'*}\fa ^*\CL]^1$ lies in the torsion of $\Pic(S)$.

{\em Proof.}
By the relative Riemann-Roch theorem (see e.g. ~\cite{fu}
Theorem 15.2)
we have the equalities in the Chow ring $A^S_{\Bbb Q}$

$$ch [R(p_{\fA'*}\fa ^*\CL]= p_{\fA'*}(ch (\fa^* \CL) td_{{\fA}'/S})=
p_{\fA*}\fa_*\fa^*[ch ( \CL) td _{{\fA}/S}]=$$
$$=\deg(\fa) p_{\fA*}[ch (\CL)td _{{\fA}/S}]=
\deg (\fa) ch [Rp_{\fA*}(\CL)]$$

Taking the components of degree 1 we get the Lemma. $\Box$

\subsubsection{Remark}
We will apply the Lemma to the situation $\fA=\Pic^0(C)\otimes\Lambda$.
In this case the hypothesis of smoothness of $S$ is redundant. Indeed,
it is enough to take for $S$ the moduli space of genus $g$ curves $\CM_g$.
It is well known that there exists a smooth covering $\pi:\ \widetilde{\CM_g}
\lra\CM_g$ inducing injection on Picard groups. One can take for
$\widetilde{\CM_g}$ the moduli space of curves with a basis in
$H^1(C,\BZ/3\BZ)$.

\medskip

We continue the proof of the Proposition.
We define the integer function $m(\Lambda,\ (,)\ )$ by the requirement
$$[Rp_{\Lambda*}(\CL (\Lambda,\ (,),\ 0))]^1 =\delta^{m(\Lambda,\ (,)\ )}$$

We will need one more Lemma.

\subsubsection{Lemma}
\label{pjat' properties}
The function $m$ satisfies the following properties:

(a) $m(\Lambda _1 \oplus \Lambda _2,\ (,)_1+(,)_2)=
m(\Lambda _1,\ (,)_1) d_2^g+
m(\Lambda _2,\ (,)_2) d_1^g$ (for the notation $d$ see Proposition ~\ref{pjat' 1/2});

(b) Let $\iota:\ \Lambda ' \hookrightarrow \Lambda $ be an embedding
with finite cokernel. Then

 $m(\Lambda ', \iota^*(,)\ ) =\# (\Lambda / \Lambda ')^g
m(\Lambda, \sigma)$;

(c) If $\Lambda ={\Bbb Z}$ and $(1,1)=-16$ then
$m(\Lambda,\ (,)\ )=-\frac{16^g}{2}$.

\subsubsection{}
We will prove the Lemma below, and now we derive the Proposition from the
Lemma.

We call 2 lattices isogenic if they contain isomorphic
sublattices of maximal rank. From (a) it follows that the Proposition
holds for the sum of 2 lattices if it holds for each of them. From
(b) it follows that the Proposition holds
for a lattice if it holds for an isogenic
one. Thus it is enough to prove it for $(\Lambda,\ (,)\ )=({\Bbb Z},n)$.

Note that the statement for the lattice $({\Bbb Z},n)$ is equivalent to
the statement for the lattice $({\Bbb Z}\oplus {\Bbb Z},n \oplus (-n))$
since $m({\Bbb Z}\oplus {\Bbb Z},n \oplus (-n))= 2 (-n)^g
\cdot m({\Bbb Z},n)$.
But the lattices  $({\Bbb Z}\oplus {\Bbb Z},n \oplus (-n))$ are
isogenic  for various $n$, and (c) says that the Proposition is true for
$n=-16$.

The Proposition is proved. This completes the proof
of the Theorem ~\ref{pjat' L'}(c). $\Box$

\subsection{}
It remains to prove the Lemma ~\ref{pjat' properties}.

(a) is easy.

(b) follows from the Lemma ~\ref{pjat' isogeny}.

To prove (c) we will need the following Lemma.

Let $C\lra S$ be a family such that a square root $\Omega^{\frac{1}{2}}$
exists globally over $C$. It gives the section $\varrho:\ S\lra\Pic^{g-1}(C)$.
On the other hand, $\Pic^{g-1}(C)$ posesses the canonical theta line bundle
$\CO(\theta)$.

\subsubsection{Lemma}
\label{pjat' Harris}
$\varrho^*(\CO(\theta))=\delta^{\frac{1}{2}}$

{\em Proof.} (J.Harris and T.Pantev)
Consider a relative curve $C$ over base $B$, and a line bundle $L$
over $C$ of relative degree $g-1$. This gives the section from $B$
to $\Pic^{g-1}(C)$, and we can restrict $\CO(\theta)$ to $B$. Evidently,
the class of $\CO(\theta)|_B$ equals the class of divisor $D\subset B:\
b\in D$ iff $h^0(L_b)>0$.

On the other hand, one checks readily that $\CO(-D)=[Rp_{C*}L]^1$.

We apply this equality to the case $B=S,\ L=\Omega^{\frac{1}{2}}$.
Then $[Rp_{C*}L]^1$ is computed by Mumford formula (see ~\cite{m1},
Theorem 5.10). Namely, we get $[Rp_{C*}\Omega^{\frac{1}{2}}]^1=\delta^k$
where $k=6(\frac{1}{2})^2-6(\frac{1}{2})+1=-\frac{1}{2}$.

We conclude that
$\varrho^*(\CO(\theta))=\delta^{\frac{1}{2}}$. $\Box$

\subsubsection{} Now we are able to prove the Lemma ~\ref{pjat' properties}(c).

We can assume that $\Omega^{\frac{1}{2}}$ exists globally over $C$.
The corresponding section $\varrho:\ S\lra\Pic^{g-1}(C)$ identifies
$\Pic^{g-1}(C)$ with $\Pic^0(C)$. We will use this identification, and
we will denote the projection $\Pic^0(C)\lra S$ by $p$.

We have $Rp_*(\CO (\theta))=\CO$ and $\varrho^*(\CO (\theta ))=
\delta^{\frac{1}{2}}$ by the Lemma ~\ref{pjat' Harris}.
We define $\CL _1 := \CO (\theta)\otimes p^*\delta^{-\frac{1}{2}}$.

Then $(4)^*\CL _1= \CL ({\Bbb Z}, -16)$.
Indeed, the $\theta$-divisor,
and hence $\CL_1$ is $(-1)$-invariant. Thus $(2)^*\CL _1$ and
$\CL ({\Bbb Z},-4)$ are the line bundles with the same Chern class,
and both are $(-1)$-invariant. So fiberwise they can differ
only by a 2-torsion element. We conclude that $\CL({\Bbb Z}, -16)=
(2)^* \CL({\Bbb Z}, -4)= (4)^*\CL _1$ fiberwise. But the restriction
of both these line bundles to the 0-section is trivial. Hence they are
isomorphic.

Now (c) follows from the Lemma ~\ref{pjat' isogeny} since obviously
 $p_*\CL _1=\delta^{-\frac{1}{2}}$. $\Box$

\subsection{}
We finish this section with a proof of one simple property
of the Heisenberg local system $\CH^\alpha,\alpha\in\BN[X_\ell]$. It states,
roughly speaking, that $\CH^\alpha$ depends only on the class of $\alpha$
modulo $\BN[Y_\ell]$. More precisely, the following Lemma holds.

\subsubsection{Lemma}
\label{pjat' Independence}
Let $\alpha=\sum_{k\in K}a_k\mu_k,\alpha'=\sum_{k\in K}a_k\lambda_k
\in\BN[X_\ell]$. Suppose $\lambda_k-\mu_k\in Y_\ell$ for every $k\in K$.
Then the natural identification $C^\alpha=C^{\alpha'}$ lifts to a canonical
isomorphism $\CH^\alpha=\CH^{\alpha'}$.

{\em Proof.} Let $\gamma\in\BN[Y_\ell],\ \gamma^\sim=0$.
We define $\beta\in\BN[X_\ell]$ as $\beta=\alpha+\gamma$.
We have an obvious projection $\pr:\ C^\beta\lra C^\alpha$.
It is enough to construct a canonical isomorphism
$\CH^\beta\iso\CH^\alpha$.

\subsubsection{}
Consider the addition homomorphism of abelian varieties
$$a:\ (\Pic(C)\otimes X_\ell)_{(2g-2)\nu_0}\times(\Pic(C)\otimes Y_\ell)_0
\lra(\Pic(C)\otimes X_\ell)_{(2g-2)\nu_0}$$
We will compute the inverse image $a^*\CE$.

To this end we introduce the
following bilinear form $\sigma$ on $X_\ell\oplus Y_\ell\oplus Y_\ell$:
$$\sigma((x,y_1,y_2),(x',y_1',y_2')):=-\frac{1}{\ell}(x\cdot y_2'+x'\cdot y_2+
y_1\cdot y_2'+y_1'\cdot y_2+y_2\cdot y_2').$$
The automorphism $(x,y_1,y_2)\mapsto(x,y_1,y_1+y_2)$ of
$X_\ell\oplus Y_\ell\oplus Y_\ell$ carries this form into $\sigma'$, where
$$\sigma'((x,y_1,y_2),(x',y_1',y_2')):=-\frac{1}{\ell}(x\cdot y_2'+x'\cdot y_2+
y_2\cdot y_2'-y_1'\cdot y_1-x\cdot y_1'-x'\cdot y_1).$$
Note that $\sigma'=\sigma_1+\sigma_2$ where
$$\sigma_1((x,y_1,y_2),(x',y_1',y_2')):=-\frac{1}{\ell}(x\cdot y_2'+x'\cdot y_2+
y_2'\cdot y_2)$$
and
$$\sigma_2((x,y_1,y_2),(x',y_1',y_2')):=-\frac{1}{\ell}(-y_1'\cdot y_1-x\cdot y_1'-
x'\cdot y_1)$$

\subsubsection{}
We will denote by $\pr_{12}$ the projection of
$(\Pic(C)\otimes X_\ell)_{(2g-2)\nu_0}\times(\Pic(C)\otimes Y_\ell)_0\times
(\Pic(C)\otimes Y_\ell)_0$ to
$(\Pic(C)\otimes X_\ell)_{(2g-2)\nu_0}\times(\Pic(C)\otimes Y_\ell)_0$ (the
product of first two factors).

We also denote by $\pr_1$ the projection of
$(\Pic(C)\otimes X_\ell)_{(2g-2)\nu_0}\times(\Pic(C)\otimes Y_\ell)_0$
to $(\Pic(C)\otimes X_\ell)_{(2g-2)\nu_0}$.

Then we have
$$a^*\CE=\pr_{12*}\CL(X_\ell\oplus Y_\ell\oplus Y_\ell,\sigma,(\nu_0,0,0))=
\pr_{12*}\CL(X_\ell\oplus Y_\ell\oplus Y_\ell,\sigma',(\nu_0,0,0))=$$
$$=\pr_{12*}(\CL(X_\ell\oplus Y_\ell\oplus Y_\ell,\sigma_1,(\nu_0,0,0))\otimes
\CL(X_\ell\oplus Y_\ell\oplus Y_\ell,\sigma_2,(\nu_0,0,0)))=$$
$$=\pr_1^*\CE\otimes\CL(X_\ell\oplus Y_\ell,\sigma_2,(\nu_0,0))$$
In the last line we view $\sigma_2$ as a bilinear form on
$X_\ell\oplus Y_\ell$ since anyway it factors through the projection of
$X_\ell\oplus Y_\ell\oplus Y_\ell$ to the sum of the first two summands.
The last equality follows from the projection formula.

\subsubsection{}
Let us denote the line bundle   $\CL(X_\ell\oplus Y_\ell,\sigma_2,(\nu_0,0))$
by $\fL$. Then we have
$$\det(a^*\CE)=\det(\pr_1^*\CE)\otimes\fL^{\od_\ell^g}$$
Hence we can define a map
$$m:\ (\det(\pr_1^*\CE))\dot{}\times\dot{\fL}\lra(\det(a^*\CE))\dot{},\
(t,u)\mapsto t\otimes u^{\od_\ell^g}.$$

Recall that $\fH$ is the universal local system on $(\det\CE)\dot{}$.
We obviously have the equality
$$m^*a^*\fH=\pr_1^*\fH$$
(To be more precise, one could put the exterior tensor product of
$\pr_1^*\fH$ with the constant sheaf on $\dot{\fL}$ in the RHS).

\subsubsection{}
We are ready to finish the proof of the Lemma.

Recall that $\CH^\beta$ is defined as $s_\beta^*\AJ_\beta^*\fH$ (notations
of ~\ref{pjat' construction}).

We consider also the Abel-Jacobi map $\AJ_{\alpha,\gamma}:\
C^\beta\lra(\Pic(C)\otimes X_\ell)_{(2g-2)\nu_0}\times(\Pic(C)\otimes Y_\ell)_0$.
As in {\em loc. cit.} we have the canonical section
$\varsigma$ of $\AJ_{\alpha,\gamma}^*(\det(\pr_1^*\CE))$
and the canonical section
$s$ of $\AJ_{\alpha,\gamma}^*(\det(a^*\CE))$.

We have
$$\CH^\beta=s^*\AJ_{\alpha,\gamma}^*a^*\fH$$
and
$$\pr^*\CH^\alpha=\varsigma^*\AJ_{\alpha,\gamma}^*\pr_1^*\fH$$

It is easy to see that there exists a section $\fs$ of
$\AJ_{\alpha,\gamma}^*\fL$ such that
$$m(\varsigma,\fs)=s.$$

Hence $\CH^\beta=s^*\AJ_{\alpha,\gamma}^*a^*\fH=(\varsigma,\fs)^*
\AJ_{\alpha,\gamma}^*\pr_1^*\fH =
\varsigma^*\AJ_{\alpha,\gamma}^*\pr_1^*\fH = \pr^*\CH^\alpha$.

This completes the proof of the Lemma. $\Box$

\bigskip

\section{Factorization isomorphisms}

\subsection{} We need to introduce some more notation. Recall the setup
of ~\ref{pjat' notations}.

First, we will need the space $TC^\tau$ containing
$\TCo^\tau$ as an open subset. It is defined in the same way as $\TCo^\tau$,
only we do not throw away nor the diagonals, neither the zero sections.

The space $TC^\tau$ decomposes into the direct product
$$TC^\tau=\widetilde{TC^K}\times\prod_{k\in K}TD^{\tau^{-1}(k)}$$
The projection to the first (resp. second) factor will be denoted by
$pr_1$ (resp. $pr_2$).

We have the natural substitution map $TC^\tau\lra TC^J$. We will denote it
by $q_\tau$ when there is no risk of confusing it with the same named map
$\TCo^\tau\lra\TCo^J$.

The evident projection map $TC^J\lra C^J$ will be denoted by $pr$.
We will use the same notation for the projections $\TCo^J\lra\Co^J$, and
$\TCo^\alpha\lra\Co^\alpha$ when there is no risk of confusion.

The open embedding $\Co^J\hra C^J$ will be denoted by $j$. We will keep
the same notation for the open embedding $\Co^\alpha\hra C^\alpha$.

\subsection{} Recall that we made a choice in the definition of $\CH^\alpha$:
the isomorphism (\ref{pjat' choice}) was defined up to a scalar multiple.

Let us explain how to make a consistent family of choices.

Given a surjection $\tau:\ J\lra K$ we denote $\alpha_K$ by $\beta$ for short.

\subsubsection{Definition} The isomorphisms $\theta_\alpha:\ \CL_\alpha
\iso\AJ_\alpha^*(\CL)$, and $\theta_\beta:\ \CL_\beta\iso\AJ_\beta^*(\CL)$
are called {\em compatible} if the following diagram commutes:
$$
\begin{array}{ccc}
(\pi_J\circ\sigma_\tau)^*\CL_\alpha&
\overset{\pi_J^*(\theta_\alpha)}{\iso}&
(\AJ_\alpha\circ\sigma_\tau)^*\CL\\
\|&&\|\\
\pi_K^*\CL_\beta&
\overset{\pi_K^*(\theta_\beta)}{\iso}&
\AJ_\beta^*\CL
\end{array}
$$
Here the left vertical equality is just the Lemma ~\ref{pjat' compatibility},
and the right vertical equality is tautological since $\AJ_\beta=
\AJ_\alpha\circ\sigma_\tau$.

Clearly, we can choose the isomorphisms (\ref{pjat' choice}) for all $J$ in a
compatible way. We will assume such a choice made once and for all.

This gives rise to the local system $\tilde{\CH}^\alpha$ on $\dot{\CL}_\alpha$.

\subsection{} Now we are ready for the construction of factorization
isomorphisms.

\subsubsection{}
\label{pjat' disk}
Let us start with the case of disk.

Given an unfolding $\pi:\ J\lra {X_\ell}$ we define the meromorphic function
$F_J$ on $TD^J$ as follows. We fix a total order $<$ on $J$.
The standard coordinates on $TD^J$ are denoted by $(x_j,\xi_j;\ j\in J)$.
We define
\begin{equation}
F_J(x_j,\xi_j)= \prod_{\pi(i)\not=\pi(j),i<j}(x_i-x_j)^{2\pi(i)\cdot\pi(j)}
\times\prod_{\pi(i)=\pi(j)}(x_i-x_j)^{\pi(i)\cdot\pi(j)}
\times\prod_{j\in J}\xi_j^{2 n(\pi(j))}
\end{equation}

The function does not depend on the  order on $J$.

Recall the setup of ~\ref{pjat' tree} for the case $C=D$.
We have the following two functions on $TD^\tau$:
\begin{equation}
\label{pjat' facfun}
F_J(q_\tau),\   \operatorname{and}\
F_K(p_K)\times\prod_{k\in K}F_{\tau^{-1}(k)}
\end{equation} 
It is easy to see that the quotient of these two functions
$TD^\tau\lra\BC$ actually lands to $\BC^*$; moreover the quotient
function $TD^\tau\lra\BC^*$ is ``canonically'' homotopic to the constant map
$TD^\tau\lra 1$.

 Let $\CQ$ denote the one
dimensional local system on $\BC^*$ with monodromy $\zeta$ around the
origin and with fiber at 1 trivialized.

We define $\CI ^J := F_J^*(\CQ)$.

 Let $m:\ \BC^*\times\BC ^*\lra\BC ^*$ denote the multiplication
map. Then $m^*(\CQ)\cong \CQ \boxtimes \CQ$ canonically. 

Hence the homotopy between the functions ~\ref{pjat' facfun} yields
the factorization isomorphisms for the local systems $\CI ^J$.

The associativity of these factorization isomorphisms
follows from  associativity of the above mentioned homotopies between
functions ~\ref{pjat' facfun}.

\subsubsection{}
We return to the case of an arbitrary curve.
We denote the ``zero section'' $C^K\lra TC^\tau$ by $z$.

Consider the following line bundle on $TC^\tau$:
$$\CL_\tau=(q_\tau\circ pr)^*\CL_\alpha$$
$\dot{\CL}_\tau$ carries the local system $\CH^\tau:=(q_\tau\circ pr)^*
\tilde{\CH}^\alpha$.

Since $z\circ q_\tau\circ pr=\sigma_\tau$, we have $z^*\CL_\tau=\CL_\beta$.

The latter line bundle over $C^K$
has a meromorphic section $\pi_K^*s_\beta$ constructed in ~\ref{pjat' construction}.
We will extend this section to the whole $TC^\tau$.

We do have a meromorphic section $s'_\tau:\ TC^\tau\lra\CL_\tau,\
s'_\tau=q_\tau^*\pi_J^*(s_\alpha)$, where $s_\alpha:\ C^\alpha\lra\CL_\alpha$
was defined in ~\ref{pjat' construction}.

Recall that in ~\ref{pjat' disk} we have defined the meromorphic function
$F_M$ on $TD^M$ for any unfolding $M\lra {X_\ell}$.

We have $$\divv(s'_\tau)= pr_1^{-1}(\divv(p_K^*\pi_K^*s_\beta))+
pr_2^{-1}(\divv(\prod_{k\in K}F_{\tau^{-1}(k)}))$$

We define the meromorphic section
\begin{equation}
s_\tau=s'_\tau\prod_{k\in K}F_{\tau^{-1}(k)}^{-1}
\end{equation}
Then $\divv(s_\tau)= pr_1^{-1}(\divv(p_K^*\pi_K^*s_\beta))$, and it is easy
to see that $z^*s_\tau=\pi_K^*s_\beta$.

Since $pr_1$ is a projection with contractible
fibers  it follows  that restricting to $\TCo^\tau$ we have
\begin{equation}
s_\tau^*(\CH^\tau)\simeq pr_1^*p_K^*\pi_K^*s_\beta^*\tilde{\CH}^\beta=
pr_1^*p_K^*\pi_K^*\CH^\beta
\end{equation}

\subsubsection{}
 Let $m:\ \BC^*\times\CL_\tau\lra\CL_\tau$ denote the multiplication
map. Then $m^*\CH^\tau=\CQ\boxtimes\CH^\tau$.

So putting all the above together, we get
$$q_\tau^*\pi_J^*\CH^\alpha=(s'_\tau)^*\CH^\tau=
(s_\tau\prod_{k\in K}F_{\tau^{-1}(k)})^*\CH^\tau=$$
$$=(\prod_{k\in K}F_{\tau^{-1}(k)};s_\tau)^*m^*\CH^\tau=
(\prod_{k\in K}F_{\tau^{-1}(k)};s_\tau)^*(\CQ\boxtimes\CH^\tau)=$$
$$=(\prod_{k\in K}F_{\tau^{-1}(k)})^*\CQ\boxtimes s_\tau^*\CH^\tau=
(\prod_{k\in K}F_{\tau^{-1}(k)})^*\CQ\boxtimes s_\beta^*\CH^\beta=$$
$$=p_K^*\pi_K^*\CH^{\alpha_K}
\boxtimes\fbox{$\times$}_{k\in K}\ \CI^{\tau^{-1}(k)}$$
which completes the construction of the factorization isomorphisms.

\subsection{} The  compatibility of the constructed
factorization isomorphism for $\CH$ with the ones for $\CI$ is immediate. 

\newpage
\begin{center}
{\bf Chapter 2. The modular property of the Heisenberg local system}
\end{center}
\vspace{.8cm}

\bigskip

\section{Degeneration of curves: recollections and notations}
\label{pjat' degeneration of curves}

\subsection{}
\label{pjat' cases}
Let $\CM_g$ be the moduli space of curves of genus $g$. Let ${\fD}$ be a smooth
locus of an irreducible component of the divisor at infinity in the
compactification $\ol{\CM_g}$. It has one of the following types.

(a) Either
we have a decomposition $g=g_1+g_2,\ g_1,g_2>0$. Then ${\fD}$ is the
moduli space of the following objects: curves $C_1,C_2$ of genera
$g_1,g_2$ respectively; points $x_1\in C_1,x_2\in C_2$.
Thus, $C_g$ degenerates into $C_1\sqcup C_2/(x_1=x_2)$.

(b) Or we put $g_0=g-1$, and then ${\fD}$ is the
moduli space of the following objects: curve $C_0$ of genus $g_0$
with two distinct points $x_1\not=x_2\in C_0$. Thus, $C_g$ degenerates
into $C_0/(x_1=x_2)$.

We denote by $C_{\fD}$ the universal family over ${\fD}$.

\subsection{}
Let $\tilde{{\fD}}\longrightarrow {\fD}$ denote the moduli space
of objects as above
plus nonzero tangent vectors $v_1$ at $x_1$, and $v_2$ at $x_2$. Let
$\TB_{\fD}\CM_g$ denote the normal bundle to ${\fD}$ in $\ol{\CM_g}$
with the zero section
removed.  There is the canonical map $\wp:
\tilde{{\fD}}\rightarrow\TB_{\fD}\ol{\CM_g}$
(see e.g. ~\cite{bfm}, \S 4).

\subsection{}
Let $\ol C$ be the universal curve over $\ol{\CM_g}$.
Let $\alpha\in\BN[{X_\ell}]$.

Consider
$\TO\ol{C}^\alpha\rightarrow\ol{\CM_g}$ (see ~\ref{pjat' TO}).
Note that $\ol C$ has singularities (nodes) over $\ol{\CM_g}-\CM_g$. To
simplify the exposition
we just throw these nodal points away. That is,
$\TO\ol{C}^\alpha$ is formed by configurations of nonsingular points
with nonzero tangent vectors.

Let us describe the preimage of ${\fD}$ in $\TO\ol{C}^\alpha$, to be
denoted by ${\cal T}C_{\fD}^\alpha$.

In case ~\ref{pjat' cases}(a) ${\cal T}C_{\fD}^\alpha$ is a union of connected
components numbered by the decompositions $\alpha={\alpha_1}+{\alpha_2}$.
The connected component ${\cal T}C_{\fD}^{{\alpha_1}{\alpha_2}}$ is the product
$\tilde{C}_1^{\alpha_1}\times\tilde{C}_2^{\alpha_2}$.
Here $\tilde{C}_r^{\alpha_r}$
denotes the configurations of ${\alpha_r}$ distinct points with
nonzero tangent vectors
on $C_r-x_r,\ r=1,2$.

In case ~\ref{pjat' cases}(b) ${\cal T}C_{\fD}^\alpha$ is the space of configurations
of $\alpha$ distinct points with nonzero tangent vectors
on $C_0-\{x_1,x_2\}$.

\subsection{}
Thus we obtain the map
$\wp_\alpha:\ \tilde{{\fD}}\times_{\fD}{\cal T}C^\alpha_{\fD}\lra
\TB_{{\cal T}C_{\fD}^\alpha}\TO\ol{C}^\alpha$ --- the normal bundle to
${\cal T}C^\alpha_{\fD}$ in $\TO\ol{C}^\alpha$ with the zero section removed.

Note that if we prescribe a weight $\mu_1$ to the point $x_1$, and
a weight $\mu_2$ to the point $x_2$ then

in case ~\ref{pjat' cases}(a) we can identify $\tilde{{\fD}}\times_{\fD}
{\cal T}C_{\fD}^{{\alpha_1}{\alpha_2}}$ with $\TO C_1^{{\alpha_1}+\mu_1}\times
\TO C_2^{{\alpha_2}+\mu_2}$ (notations of ~\ref{pjat' TO});

in case ~\ref{pjat' cases}(b) we can identify $\tilde{{\fD}}\times_{\fD}
{\cal T}C_{\fD}^\alpha$ with $\TO C_0^{\alpha+\mu_1+\mu_2}$

(the exponents are formal sums, and we do not perform symmetrization
over $\mu_1$ and $\mu_2$ if they happen to be equal or contained in
supp ~$\alpha$).

\subsection{}
It is known that the determinant line bundle $\delta_g$ over $\CM_g$
extends to the line bundle $\ol{\delta_g}$ over $\ol{\CM_g}$.

In case ~\ref{pjat' cases}(a) let $\varpi_i, i=1,2$, denote the natural projection
from $\fD$ to $\CM_{g_i}$.

In case ~\ref{pjat' cases}(b) let $\varpi_0$ denote the natural projection
from $\fD$ to $\CM_{g_0}$.

Then $\ol{\delta_g}|_{\fD}=\varpi_1^*\delta_{g_1}\otimes\varpi_2^*\delta_{g_2}$
(resp. $\ol{\delta_g}|_{\fD}=\varpi_0^*\delta_{g_0}$) (see ~\cite{bfm}).

Summing up we obtain the map also denoted by $\wp_\alpha$:

in case ~\ref{pjat' cases}(a):

\begin{equation}
\label{pjat' map a}
(\TO C_1^{{\alpha_1}+\mu_1}\times_{\CM_{g_1}}\delta_{g_1})\times
(\TO C_2^{{\alpha_2}+\mu_2}\times_{\CM_{g_2}}\delta_{g_2})\lra
\TB_{{\cal T}C_{\fD}^\alpha\times_{\ol{\CM_g}}\ol{\delta_g}}(\TO\ol{C}^\alpha
\times_{\ol{\CM_g}}\ol{\delta_g});
\end{equation}

in case ~\ref{pjat' cases}(b):

\begin{equation}
\label{pjat' map b}
\TO C_0^{\alpha+\mu_1+\mu_2}\times_{\CM_{g_0}}\delta_{g_0}\lra
\TB_{{\cal T}C_{\fD}^\alpha\times_{\ol{\CM_g}}\ol{\delta_g}}(\TO\ol{C}^\alpha
\times_{\ol{\CM_g}}\ol{\delta_g}).
\end{equation}

\subsection{}
\label{pjat' modular Heisenberg}

Consider the specialization $\bSp{\cal H}_g^\alpha$ of the
Heisenberg local system along the boundary component
${\cal T}C_{\fD}^\alpha\times_{\ol{\CM_g}}\ol{\delta_g}$.
It is a sheaf on
$\TB_{{\cal T}C_{\fD}^\alpha\times_{\ol{\CM_g}}\ol{\delta_g}}(\TO\ol{C}^\alpha
 \times_{\ol{\CM_g}}\ol{\delta_g})$. We will describe its inverse image
under the map $\wp_\alpha$.

{\bf Theorem.}
(a) In case ~\ref{pjat' cases}(a) we have $$\wp_\alpha^*\bSp{\cal H}_g^\alpha|_
{(\TO C_1^{{\alpha_1}+\mu_1}\times_{\CM_{g_1}}\delta_{g_1})\times
(\TO C_2^{{\alpha_2}+\mu_2}\times_{\CM_{g_2}}\delta_{g_2})}=
{\cal H}_{g_1}^{{\alpha_1}+\mu_1}\boxtimes{\cal H}_{g_2}^{{\alpha_2}+\mu_2}$$
for $\mu_1=(2g_1-2)\nu_0-{\alpha_1}^\sim$, and $\mu_2=(2g_2-2)\nu_0-
{\alpha_2}^\sim$ (notations of ~\ref{pjat' tilde});

(b) In case ~\ref{pjat' cases}(b) we have $$\wp_\alpha^*\bSp{\cal H}_g^\alpha=
\oplus_{\mu_1+\mu_2=-2\nu_0}{\cal H}_{g_0}^{\alpha+\mu_1+\mu_2}$$
Here the sum is taken over the set ${X_\ell}/Y_\ell$ (see IV.9.1.2 and
~\ref{pjat' d}). The statement makes sense by the virtue of the Lemma
~\ref{pjat' Independence}.  Thus there are exactly $\od_\ell$ summands.

\medskip

The proof occupies the rest of this Chapter.

\newpage
\section{Proof of Theorem ~\ref{pjat' modular Heisenberg}(a)}

We are in the situation of ~\ref{pjat' cases}(a). Let us denote
$\CM_g\cup\fD$ by $\ol{\CM}^\circ\subset\ol{\CM_g}$. In this section
we will extend the
constructions of sections ~\ref{pjat' proof L} and ~\ref{pjat' proof H} to the
universal family of curves over $\ol{\CM}^\circ$.

\subsection{}
It is well known that $\Pic(C/{\CM_g})$ extends to the family of algebraic
groups over $\ol{\CM}^\circ$ which will be denoted by $\ol{\Pic}(\ol{C}/
\ol{\CM}^\circ)$.

Restricting to $\fD$ we get

\begin{equation}
\label{pjat' decomp Pic}
\ol{\Pic}(\ol{C}/\ol{\CM}^\circ)|_{\fD}=\Pic(C_1/\fD)\times_{\fD}\Pic(C_2/\fD)
/\langle(x_1)-(x_2)\rangle,
\end{equation}

and $\Pic(\ol{C}/\fD)=\Pic(C_1/\fD)\times_{\fD}\Pic(C_2/\fD)$.

In particular, $\ol{\Pic}(\ol{C}/\fD)$ is an extension of $\BZ$ by an
abelian scheme.

\subsubsection{}
Recall the notations of ~\ref{pjat' Lambda}. Given $\Lambda,\ (,),\ \nu$
as in ~\ref{pjat' Lambda} one constructs the linear bundle
$\CL(\Lambda,\ (,),\ \nu)$ over $A_\Lambda$ (see ~\ref{pjat' L(Lambda)}).
To unburden the notations we will denote this line bundle simply by
$\CL$ when it does not cause confusion.

$\ol{A}_\Lambda$ will denote $(\ol{\Pic}(\ol{C}/\ol{\CM}^\circ)
\otimes\Lambda)_{(2g-2)\nu}$, and $p$ will denote the projection
$\ol{A}_\Lambda\lra\ol{\CM}^\circ$.

We have an isomorphism $$\mho:\ A_\Lambda(C_1/\fD)\times_{\fD}
A_\Lambda(C_2/\fD)\iso\ol{A}_\Lambda|_{\fD}$$ where
$$\mho(a,b):= (a+x_1\otimes\nu,b+x_2\otimes\nu)\in
(\Pic(C_1/\fD)\otimes\Lambda)_{(2g_1-1)\nu}\times_{\fD}
(\Pic(C_2/\fD)\otimes\Lambda)_{(2g_2-1)\nu}=$$
$$=\ol{\Pic}(\ol{C}/\fD)\otimes\Lambda)_{(2g-2)\nu}$$
(the latter isomorphism is obtained from ~(\ref{pjat' decomp Pic})).

\subsubsection{Lemma}
\label{pjat' mho}
There exists a line bundle $\ol{\CL}$ on
$\ol{A}_\Lambda$ such that $\ol{\CL}|_{p^{-1}(\CM_g)}=\CL$, and
$$\mho^*(\ol{\CL}|_{\fD})=\CL_{C_1}^{\od_\ell^{g_2}}
\boxtimes\CL_{C_2}^{\od_\ell^{g_1}}.$$

{\em Proof.} Note that the line bundle $\CL_\alpha$ (see ~\ref{pjat' L alpha})
defined over $C^\alpha$ extends to $C^\alpha\cup(C_1-x_1)^{\alpha_1}\times
(C_2-x_2)^{\alpha_2}\subset\ol{C}^\alpha$ by the same formula ~\ref{pjat' L alpha}.
Moreover, $\CL_\alpha$ automatically extends to $C^\alpha\cup C_1^{\alpha_1}\times
C_2^{\alpha_2}$ since the complement to
$C^\alpha\cup(C_1-x_1)^{\alpha_1}\times
(C_2-x_2)^{\alpha_2}\subset\ol{C}^\alpha$
in $C^\alpha\cup C_1^{\alpha_1}\times C_2^{\alpha_2}$ has codimension two.

We define the following line bundle on
$C^\alpha\cup C_1^{\alpha_1}\times C_2^{\alpha_2}$:

\begin{equation}
\label{pjat' formula}
\ol{\CL}_\alpha:=\CL_\alpha(\ob(\mu_1)
(C_1^{\alpha_1}\times C_2^{\alpha_2}))
\end{equation}

Here $C_1^{\alpha_1}\times C_2^{\alpha_2}$ is viewed as a divisor on
$C^\alpha\cup C_1^{\alpha_1}\times C_2^{\alpha_2}$. Note also that
$\ob(\mu_1)=\ob(\mu_2)$ since $\mu_1+\mu_2=-2\nu$.

The line bundle $\CL$ on $A_\Lambda$ extends uniquely to a line bundle
$\ol{\CL}$ on  $\ol{A}_\Lambda$
in such a way that

\begin{equation}
\label{pjat' cunt}
\AJ^*\ol{\CL}=\ol{\CL}_\alpha.
\end{equation}

 In effect, any two
extensions of $\CL$ to $\ol{A}_\Lambda$ differ by the twist by some power of
$\CO(p^{-1}\fD)$ ---
the line bundle lifted from $\ol{\CM}^\circ$. But evidently,
$\AJ^*\CO(p^{-1}\fD)=\CO(C_1^{\alpha_1}\times C_2^{\alpha_2})$, so we can choose
our twist uniquely to satisfy ~(\ref{pjat' cunt}).

\subsubsection{}
Let us choose a basis $I$ of $\Lambda$. Let us fix
${\alpha_1},{\alpha_2}\in\BN[\Lambda];\
{\alpha_1}=\sum_{i\in I}b_ii+\sum_{i\in I}b_{-i}(-i),
{\alpha_2}=\sum_{i\in I}c_ii+\sum_{i\in I}c_{-i}(-i)$
such that $b_i,b_{-i}>g_1;\
c_i,c_{-i}>g_2$.

We have the closed embeddings $\sigma_1:\ C_1^{\alpha_1}\hra C_1^{{\alpha_1}+
\mu_1}$  (resp. $\sigma_2:\ C_2^{\alpha_2}\hra C_2^{{\alpha_2}+\mu_2}$)
adding to a configuration on $C_1$ (resp. $C_2$) $\mu_1$ copies of $x_1$
(resp. $\mu_2$ copies of $x_2$).

It is clear that $\AJ_1\circ\sigma_1$ (resp. $\AJ_2\circ\sigma_2$)
induces injection from $\Pic(A_{1\Lambda}/\fD)$ to $\Pic(C_1^{\alpha_1}/\fD)$
(resp. from $\Pic(A_{2\Lambda}/\fD)$ to $\Pic(C_2^{\alpha_2}/\fD)$).

Hence it is enough to check that

\begin{equation}
\label{pjat' bum}
(\mho\circ(\AJ_1\times\AJ_2)\circ(\sigma_1\times\sigma_2))^*\ol{\CL}|_
{p^{-1}\fD}=\sigma_1^*\CL_{{\alpha_1}+\mu_1}^{\od_\ell^{g_2}}\boxtimes
\sigma_2^*\CL_{{\alpha_2}+\mu_2}^{\od_\ell^{g_1}}
\end{equation}

Consider unfoldings $\pi_1: J_1\lra\Lambda$ of ${\alpha_1}$, and
$\pi_2: J_2\lra\Lambda$ of ${\alpha_2}$. It is enough to check that

\begin{equation}
\label{pjat' bum'}
(\pi_1\times\pi_2)^*
(\mho\circ(\AJ_1\times\AJ_2)\circ(\sigma_1\times\sigma_2))^*\ol{\CL}|_
{p^{-1}\fD}=(\pi_1\times\pi_2)^*
(\sigma_1^*\CL_{{\alpha_1}+\mu_1}^{\od_\ell^{g_2}}
\boxtimes\sigma_2^*\CL_{{\alpha_2}+\mu_2}^{\od_\ell^{g_1}})
\end{equation}

One checks readily that the isomorphism ~(\ref{pjat' bum'}) holds after restriction
to $(C_1-x_1)^{J_1}\times(C_2-x_2)^{J_2}$.

The necessary condition to extend the isomorphism across the complement
divisor is the equality of (relative) first Chern classes:
$$
c_1((\mho\circ(\AJ_1\times\AJ_2)\circ(\sigma_1\times\sigma_2))^*\ol{\CL}|_
  {p^{-1}\fD})=
c_1(\sigma_1^*\CL_{{\alpha_1}+\mu_1}^{\od_\ell^{g_2}}
\boxtimes\sigma_2^*\CL_{{\alpha_2}+\mu_2}^{\od_\ell^{g_1}})$$
This is immediate since $c_1(\mho^*(\ol{\CL}_{\fD}))=
c_1(\CL_{C_1}^{\od_\ell^{g_2}})
\boxtimes c_1(\CL_{C_2}^{\od_\ell^{g_1}})$.
If the complement divisor were irreducible, this condition would be
sufficient as well. If the divisor is reducible, one has to check that
the fundamental classes of different irreducible components are linearly
independent. All the irreducible components of our divisor are diagonals
$z_j=x_r,\ r=1,2;j\in J_r$. Their fundamental classes are linearly
independent in $H^2(C_1^{J_1}\times C_2^{J_2})$.

This completes the proof of the Lemma. $\Box$

\subsection{}
\label{pjat' proceedings}
We proceed with the proof of the Theorem ~\ref{pjat' modular Heisenberg}(a).

Let us consider the line bundle $\ol{\CL}:=
\ol{\CL}({X_\ell},\frac{{\od}_\ell^g}{\ell}?\cdot?,
\nu_0)$ on $\ol{A}:=\ol{A}_{X_\ell}$. On the other hand, let us consider the
triple $\Lambda,\ (,),\ \nu$ defined in ~\ref{pjat' nu}, and the corresponding
line bundle $\ol{\CL}':=\ol{\CL}(\Lambda,\ (,),\ \nu)$ on
$\ol{A}':=\ol{A}_{{X_\ell}\oplus Y_\ell}$. We denote by $\fp$ (resp. $\fp_r,\ r=1,2$)
the obvious projection
$\ol{A}'\lra\ol{A}$ (resp. $A'(C_r)\lra A(C_r)$). 

We define $\ol{\CE}:=\fp_*\ol{\CL}',\ \CE_{C_r}:=\fp_{r*}\CL'_{C_r}$.
It is obvious that

\begin{equation}
\label{pjat' emho}
\ol{\CE}|_{p^{-1}(\fD)}=(\mho^{-1})^*(\CE_{C_1}\boxtimes\CE_{C_2})
\end{equation}

and hence

\begin{multline}
\label{pjat' demho}
\det(\ol{\CE})|_{p^{-1}(\fD)}=
(\mho^{-1})^*((\det(\CE_{C_1}))^{\rk(\CE_{C_2})} \boxtimes
(\det(\CE_{C_2}))^{\rk(\CE_{C_1})})=\\
=(\mho^{-1})^*((\det(\CE_{C_1}))^{\od_\ell^{g_2}} \boxtimes
(\det(\CE_{C_2}))^{\od_\ell^{g_1}})
\end{multline}

Recall that the line bundle $\delta_g$ on $\CM_g$ extends to the line
bundle $\ol{\delta_g}$ on $\ol{\CM_g}$, and $\ol{\delta_g}|_{\fD}=
\delta_{g_1}\boxtimes\delta_{g_2}$.

{\bf Proposition.}
$\det(\ol{\CE})=\ol{\CL}\otimes(p^*\ol{\delta_g})^
{{\od}_\ell^g(-\frac{1}{2}\rk {X_\ell}+6k\frac{\nu_0\cdot\nu_0}{\ell})}$
where $p$ denotes
the projection $\ol{A}\lra\ol{\CM}^\circ$.

{\em Proof.} If we restrict both sides to $p^{-1}\CM_g\subset p^{-1}
\ol{\CM}^\circ$ then we just have the Theorem ~\ref{pjat' L'}.

On the other hand, it is easy to see that the normal line bundle
to $\fD$ in $\ol{\CM}^\circ$ is nontrivial. This means that two extensions
of a line bundle from $p^{-1}\CM_g$ to $p^{-1}\ol{\CM}^\circ$ coincide
iff their restrictions to $p^{-1}\fD$ coincide.

So it remains to compare the LHS and RHS restricted to $p^{-1}\fD$.

By the virtue of ~(\ref{pjat' demho}) above we have
$$\det(\ol{\CE})|_{p^{-1}(\fD)}=
(\mho^{-1})^*((\det(\CE_{C_1}))^{\od_\ell^{g_2}} \boxtimes
  (\det(\CE_{C_2}))^{\od_\ell^{g_1}})=$$
$$=(\mho^{-1})^*
((\CL_{C_1}^{\od_\ell^{g_2}}\boxtimes\CL_{C_2}^{\od_\ell^{g_1}})\otimes
(\delta_{g_1}^{\od_\ell^g(-\frac{1}{2}\rk {X_\ell}+6\frac{\nu_0\cdot\nu_0}{\ell})}\boxtimes
\delta_{g_2}^{\od_\ell^g(-\frac{1}{2}\rk {X_\ell}+6\frac{\nu_0\cdot\nu_0}{\ell})}))=$$
$$=\ol{\CL}|_{p^{-1}(\fD)}\otimes
p^*\ol{\delta}^{\od_\ell^g(-\frac{1}{2}\rk {X_\ell}+6\frac{\nu_0\cdot\nu_0}{\ell})}$$
Here the second equality follows from the Theorem ~\ref{pjat' L'}, and the third
one follows from the Lemma ~\ref{pjat' mho}.

This completes the proof of the Proposition.
$\Box$

\subsection{}
We are ready to finish the proof of the Theorem ~\ref{pjat' modular Heisenberg}(a).

\subsubsection{}
Let $\CN$ be the standard deformation of $\TO\ol{C}^\alpha\times_
{\ol{\CM_g}}\ol{\delta_g}$ to the normal cone of
$\CT C_{\fD}^\alpha\times_{\ol{\CM_g}}\ol{\delta_g}$ (see e.g. ~\cite{fu}).

So we have
$$T_{{\cal T}C_{\fD}^\alpha\times_{\ol{\CM_g}}\ol{\delta_g}}(\TO\ol{C}^\alpha
   \times_{\ol{\CM_g}}\ol{\delta_g})\hra\CN\hookleftarrow
\TO\ol{C}^\alpha\times_{\ol{\CM_g}}\ol{\delta_g}\times\BC^*$$

Let $t:\ \CN\lra\BC;\ \pr:\ \CN\lra
\TO\ol{C}^\alpha\times_{\ol{\CM_g}}\ol{\delta_g}$ denote the canonical
projections. We define an open subset $U\subset \CN$ as
$$U:=(\TO C^\alpha\times_{\CM_g}\dot{\delta})\times\BC^*\sqcup
\TO_{{\cal T}C_{\fD}^\alpha\times_{\ol{\CM_g}}\ol{\delta_g}}(\TO\ol{C}^\alpha
 \times_{\ol{\CM_g}}(\ol{\delta_g})\dot{})$$

To finish the proof it is enough to construct a local system $\CH_{\CN}$
on $U$ such that

\begin{equation}
\label{pjat' open}
\CH_{\CN}|_{(\TO C^\alpha\times_{\CM_g}\dot{\delta}_g)\times
{\Bbb R}^{>0}}=\CH^\alpha_g;
\end{equation}

\begin{equation}
\label{pjat' closed}
\wp_\alpha^*(\CH_{\CN}|_{(\TO C_1^{{\alpha_1}+\mu_1}
\times_{\fD}\dot{\delta}_{g_1})\times
(\TO C_2^{{\alpha_2}+\mu_2}\times_{\fD}\dot{\delta}_{g_2})})=
\CH^{{\alpha_1}+\mu_1}_{g_1}\boxtimes\CH^{{\alpha_2}+\mu_2}_{g_2}.
\end{equation}

Recall the notations of ~\ref{pjat' construction}. The canonical meromorphic
section $s_\alpha$ of $\AJ_\alpha^*(\det(\CE))$
extends to the same named meromorphic
section of $\AJ_\alpha^*(\det(\ol{\CE}))$
over $T\ol{C}^\alpha\times_{\ol{\CM_g}}(\ol{\delta_g})\dot{}$.

\subsubsection{Lemma}
\label{pjat' universal}
The section $s_\alpha':=t^{-\frac{\od_\ell^g}{\ell}n(\mu_1)}\pr^*s_\alpha$
is regular on $U$, and we have
$$\wp_\alpha^*(s_\alpha'|_{p^{-1}(\fD)})=s_{\alpha_1}^{\od_\ell^{g_2}}
\otimes s_{\alpha_2}^{\od_\ell^{g_1}}$$
(see ~(\ref{pjat' demho}))

{\em Proof.} Immediate from the Proposition ~\ref{pjat' proceedings} and
the construction of $\ol{\CL}$ in ~\ref{pjat' mho}.
Note also that $n(\mu_1)=n(\mu_2).\ \Box$

\subsubsection{Remark}
If one of $g_1,g_2$ is equal to $1$, then the following remark is in order.
Suppose, say, $g_1=1$. Then $N=\od_\ell^{g_2}$ satisfies the assumptions of
~\ref{pjat' elliptic}, and $s_{\alpha_1}^{\od_\ell^{g_2}}$ should be understood
as the canonical section of $\AJ_{\alpha_1}^*(\det(\CE_{C_1}))^{\otimes
\od_\ell^{g_2}}$ arising from the Theorem ~\ref{pjat' L'}(d).

\subsubsection{}
We construct the local
system $\ol{\fH}$ on $(\det(\CE))\dot{}$ as in ~\ref{pjat' Beilinson Kazhdan}.
Now we define $$\CH_{\CN}:=(s_\alpha')^*\ol\fH$$
The property ~(\ref{pjat' open}) is clear. Let us prove the property ~(\ref{pjat' closed}).

It is clear that the isomorphism ~(\ref{pjat' emho}) is the isomorphism of bundles
with flat projective connections.

Suppose $g_1,g_2>1$.
Let us define a map
$$\widetilde{\mho}:\ (\det(\CE_{C_1}))\dot{}\times_{\fD}(\det(\CE_{C_1}))\dot{}
\lra(\det(\ol{\CE}|_{\fD}))\dot{}$$
as a composition of
$$(\det(\CE_{C_1}))\dot{}\times_{\fD}(\det(\CE_{C_1}))\dot{}\lra
(\det(\CE_{C_1})^{\otimes\od_\ell^{g_2}})\dot{}\otimes
(\det(\CE_{C_1})^{\otimes\od_\ell^{g_1}})\dot{},\ (\lambda_1,\lambda_2)\mapsto
\lambda_1^{\od_\ell^{g_2}}\otimes\lambda_2^{\od_\ell^{g_1}},$$
the isomorphism ~(\ref{pjat' demho}), and the obvious map
$(\mho^*(\det(\ol{\CE}|_{\fD})))\dot{}\lra(\det(\ol{\CE}|_{\fD}))\dot{}$.

Then we see that $$\widetilde{\mho}^*(\ol{\fH}|_{\fD})=\fH_{C_1}\boxtimes
\fH_{C_2}$$

The property ~(\ref{pjat' closed}) follows.

This completes the proof of ~\ref{pjat' modular Heisenberg} in case $g_1,g_2>1$.

The minor changes in case one of $g_1,g_2$ equals $1$ are safely left
to the interested reader. $\Box$

\section{Proof of Theorem ~\ref{pjat' modular Heisenberg}(b)}

We start the proof with the following weaker statement.

\subsection{Proposition}
\label{pjat' weak}
$\wp_\alpha^*\bSp{\cal H}_g^\alpha=
\oplus_{\mu_1+\mu_2=-2\nu_0}{\cal H}_{g_0}^{\alpha+\mu_1+\mu_2}\otimes
\CP_{\mu_1}$

for some one-dimensional local systems $\CP_{\mu_1}$.

\subsection{}
To prove the proposition, let us introduce some notations.
We denote $\CM_g\cup\fD$ by $\ol{\CM}^\circ\subset\ol{\CM_g}$.
It is well known that $\Pic(C/{\CM_g})$ extends to the family of algebraic
groups $\Pic(\ol{C}/\ol{\CM}^\circ)$ over $\ol{\CM}^\circ$.

Restricting to $\fD$ we get the projection
$$\pi:\ \Pic(\ol{C}/\ol{\CM}^\circ)|_\fD\lra\Pic(C_0/\fD)$$
which identifies $\Pic(\ol{C}/\ol{\CM}^\circ)|_\fD$ with a $\BC^*$-torsor
over $\Pic(C_0/\fD)$.

The class of this torsor is equal to $(x_1)-(x_2)\in\Pic^0(C_0/\fD)=
\Pic^0(\Pic^n(C_0/\fD))$ for any $n\in\BZ$.

For any $\mu_1,\mu_2\in X_\ell$ such that $\mu_1+\mu_2=-2\nu_0$ we have
a subtraction map

\begin{multline}
\label{pjat' subtract}
r_{\mu_1,\mu_2}:\ (\Pic(C_0/\fD)\otimes X_\ell)_{(2g-2)\nu_0}
\lra(\Pic(C_0/\fD)\otimes X_\ell)_{(2g_0-2)\nu_0},\\
x\mapsto x-\mu_1\otimes x_1-\mu_2\otimes x_2
\end{multline}

Recall that in ~\ref{pjat' L'} we have defined for any curve $C$ of genus $g$
a vector bundle with flat projective connection $\CE$ on
$(\Pic(C)\otimes X_\ell)_{(2g-2)\nu_0}$. For the curve $C_0$ we will
denote this bundle by $\CE_{C_0}$.

We define the vector bundle $\CE_{\mu_1}:=r_{\mu_1,\mu_2}^*\CE_{C_0}$ on
$(\Pic(C_0/\fD)\otimes X_\ell)_{(2g-2)\nu_0}$. It carries a flat projective
connection.

If $\mu_1'-\mu_1\in Y_\ell$ then the vector bundles $\CE_{\mu_1}$ and
$\CE_{\mu_1'}$ differ by a twist by a line bundle, and the corresponding
flat projective connections are identified by the Lemma ~\ref{pjat' Independence}.
In particular, the locally constant sheaf of algebras $\End(\CE_{\mu_1})$
canonically depends only on the class of $\mu_1$ modulo $Y_\ell$.
We will use the notation $\End(\CE_\mu)$ for $\mu\in X_\ell/Y_\ell$ from now on.

\subsection{}
The plan of the proof is as follows.

Consider the sheaf of algebras $\End(\CE)$ on
$(\Pic(C)\otimes X_\ell)_{(2g-2)\nu_0}$. The data of a flat projective
connection on $\CE$ is equivalent to the data of a flat genuine connection
on $\End(\CE)$ preserving the algebra structure.

Let $\CU$ be a small (punctured) neighbourhood (in classical topology)
of $p^{-1}(\fD)$ in $(\Pic(C)\otimes X_\ell)_{(2g-2)\nu_0}$.
The monodromy around $p^{-1}(\fD)$ acts on $\fH$, or, equivalently, on $\CE$,
semisimply. It acts on the sheaf of algebras $\End(\CE)|_{\CU}$ by the
algebra automorphism
$$\fm:\ \End(\CE)\lra\End(\CE)$$
We will construct a locally constant sheaf of algebras
 $\fM\subset \End(\CE)$ with the following properties.

\subsubsection{}
\label{pjat' M}
(i) $\fM$ can be extended to a locally constant sheaf of algebras
$\ol{\fM}$ on $(\Pic(\ol{C}/\ol{\CM}^\circ)\otimes X_\ell)_{(2g-2)\nu_0}$.

(ii) The sheaf $\CE|_{\CU}$ decomposes into direct sum of $\od_\ell$ eigensheaves
of $\fm:\ \CE|_{\CU}=\oplus_\lambda\CE_\lambda$; and
$\fM=\oplus_\lambda\End(\CE_\lambda)$.

(iii) $\ol{\fM}|_{p^{-1}(\fD)}=\oplus_{\mu\in X_\ell/Y_\ell}
\End(\CE_\mu)$.

\subsubsection{}
The construction and the proof of the above properties will occupy the rest of this section.
Let us derive the Proposition ~\ref{pjat' weak} from these properties.

The decomposition (ii) of  $\CE|_{\CU}$ induces a decomposition:
 $\bSp\fH_g=\oplus_\mu\fH_\mu$.
It follows from (ii),(iii) that there is a bijection between the set of summands
and $X_\ell/Y_\ell$ such that
$\ol{\fM}|_{p^{-1}(\fD)}=\oplus_{\mu\in X_\ell/Y_\ell}\End(\fH_\mu)$.
Moreover, $\End(\fH_\mu)=\End(\CE_\mu)$.

Now for any $g$-admissible $\alpha\in X_\ell$ we have an isomorphism of
local systems
$$\bigoplus_{\mu_1+\mu_2=-2\nu_0}^{\mu_1,\mu_2\in X_\ell/Y_\ell}
\End(\CH_{g_0}^{\alpha+\mu_1+\mu_2})=
\bigoplus_{\mu_1+\mu_2=-2\nu_0}^{\mu_1,\mu_2\in X_\ell/Y_\ell}
s_{\alpha+\mu_1+\mu_2}^*\AJ_{\alpha+\mu_1+\mu_2}^*\End(\CE_{\mu_1})=$$
$$=\bigoplus_{\mu_1+\mu_2=-2\nu_0}^{\mu_1,\mu_2\in X_\ell/Y_\ell}
s_{\alpha+\mu_1+\mu_2}^*\AJ_{\alpha+\mu_1+\mu_2}^*\End(\fH_{\mu_1})$$
Hence
$s_{\alpha+\mu_1+\mu_2}^*\AJ_{\alpha+\mu_1+\mu_2}^*\fH_{\mu_1}=
\CH_{g_0}^{\alpha+\mu_1+\mu_2}\otimes\CP_{\mu_1}$ for some one-dimensional
local system $\CP_{\mu_1}$.

On the other hand, obviously,
$$\bigoplus_{\mu_1+\mu_2=-2\nu_0}^{\mu_1,\mu_2\in X_\ell/Y_\ell}
  s_{\alpha+\mu_1+\mu_2}^*\AJ_{\alpha+\mu_1+\mu_2}^*\fH_{\mu_1}=
s_{\alpha}^*\AJ_{\alpha}^*\bSp\fH_g=\bSp\CH_g^\alpha$$
This completes the proof of the Proposition ~\ref{pjat' weak}. $\Box$

\subsection{}
\label{pjat' useful}
Before we prove the properties ~\ref{pjat' M}(i)--(iii) we need to introduce the
following general construction.

Let $p_\fA:\ \fA\lra S$ be an abelian scheme over a base $S$. Let $L$ be
a relatively ample line bundle on $\fA$. It defines a morphism $\phi_L$
from $\fA$ to the dual abelian variety $\check{\fA}:\ a\mapsto T_aL\otimes L^{-1}$
where $T_a$ denotes the translation by $a$.
Let $i_G:\ G\hra \fA$ be the kernel of this homomorphism.
Then $p_G:\ G\lra S$ is an etale cover.

Let $L_\heartsuit$ denote the line bundle
$L\otimes p_\fA^*0^*L^{-1}$ where $0$ stands
for the zero section $S\lra \fA$. Recall that $\dot{L}$ denotes the $\BC^*$-torsor
corresponding to the line bundle $L$.

\subsubsection{}
The sheaf $p_{G*}i_G^*\dot{L_\heartsuit}$ has a natural structure of a group scheme
acting on the sheaf $p_{\fA*}L$.
There is an exact sequence
$$0\lra\CO^*\lra p_{G*}i_G^*\dot{L_\heartsuit}\overset{\fp}{\lra}G\lra 0$$
The subsheaf $\CO^*$ acts on $p_{\fA*}L$ by multiplications with functions.
The action of an element $h\in p_{G*}i_G^*\dot{L_\heartsuit}$ covers the translation
by $\fp(h)\in G\subset \fA$.

The group scheme $p_{G*}i_G^*\dot{L_\heartsuit}$ carries a unique flat connection
compatible with multiplication. In effect, the elements of finite order
form a Zariski dense union of etale group subschemes.

\subsubsection{}
The sheaf $p_{G*}i_G^*L_\heartsuit$ has a natural structure of algebra acting on
the sheaf $p_{\fA*}L$. The natural inclusion $p_{G*}i_G^*\dot{L_\heartsuit}\hra
p_{G*}i_G^*L_\heartsuit$ is multiplicative and compatible with the action on $p_{\fA*}L$.
The sheaf of algebras $p_{G*}i_G^*L_\heartsuit$ carries a unique flat connection
compatible with the one on $p_{G*}i_G^*\dot{L_\heartsuit}$.

The map

\begin{equation}
\label{pjat' gift}
p_{G*}i_G^*L_\heartsuit\iso\End(p_{\fA*}L)
\end{equation}

 is an isomorphism.

\subsubsection{Remark} Note that the above construction gives an alternative
way to define the flat projective connection on $p_{\fA*}L$ --- without
referring to ~\cite{bk}.

\subsection{}
\label{pjat' useless}
We apply the above construction to the following situation.
We take $S=(\Pic(C/\CU)\otimes X_\ell)_{(2g-2)\nu_0};\
\fA=(\Pic(C/\CU)\otimes(X_\ell\oplus Y_\ell))_{(2g-2)\nu_0,0};\
\pi=\pr_1;\ L=\CL'$ as in ~\ref{pjat' nu}.

Then $G=p_{(\Pic(C/\CU)\otimes X_\ell)_{(2g-2)\nu_0}}^*R^1p_{C*}
(\ul{X_\ell/Y_\ell})$.

Recall that by the Picard-Lefschetz theory
the monodromy $\fm$ acts on $R^1p_{C*}(\ul{\BZ})$ by the formula
$$y\mapsto y\pm\langle y,x\rangle x$$ where $x$ is the vanishing cycle.
In particular, the invariants of $\fm$ in $R^1p_{C*}(\ul{\BZ})$
form a sublattice of codimension 1, namely, the orthogonal complement of the
vanishing cycle.

We define $\fM\subset\End(\CE)=p_{G*}i_G^*\CL'_\heartsuit$ as follows:
$$\fM:=p_{G^\fm*}i_{G^\fm}^*\CL'_\heartsuit$$ where $G^\fm$ stands for the invariants
of $\fm$ on $G$.

\subsection{}
Let us  check the properties ~\ref{pjat' M}.

The sheaf $G^\fm=(X_\ell/Y_\ell)\otimes(R^1p_{C*}(\ul{\BZ}))^\fm$ extends to
a sheaf $\ol{G}$ on $\ol{\CM}^\circ$. Moreover, we have an embedding
$$i_{\ol{G}}:\ \ol{G}\hra(\Pic(\ol{C}/\ol{\CM}^\circ)\otimes Y_\ell)_0$$
Let us choose some extension $\ol{\CL'}$ of $\CL'$ to
$(\Pic(\ol{C}/\ol{\CM}^\circ)\otimes(X_\ell\oplus Y_\ell))_{(2g-2)\nu_0,0}$.
Recall that $\ol{\CL'}_\heartsuit$ denotes
$\ol{\CL'}\otimes p^*0^*\ol{\CL'}^{-1}$.
We define
$$\ol{\fM}:=p_{\ol{G}*}i_{\ol{G}}^*\ol{\CL'}_\heartsuit$$
Evidently, $\ol{\fM}|_{\CU}=\fM$.

It is easy to see that the algebra structure on $\fM$ extends uniquely
to $\ol {\fM}$. This proves (i).

Property (i) implies that ${\fM}\subset (\End (\CE ))^\fm $. On the other
hand one can check that there exists a unique decomposition of $\CE$
in the sum of $\od _\ell $ summands of equal dimension:
$\CE =\oplus \CE _\lambda$ such that $\fM =\oplus  \End (\CE _\lambda)
\subset \End (\CE )$. Since  $\fm$  commutes with $\fM$ it follows that
each $\CE _\lambda$ is an eigenspace of $\fm$. (ii) is proved.

\subsection{}
It remains to check ~\ref{pjat' M}(iii). Recall that for $\mu_1,\mu_2\in X_\ell$ such
that $\mu_1+\mu_2=-2\nu_0$ we have introduced the substraction map
$r_{\mu_1,\mu_2}$ in ~(\ref{pjat' subtract}).

Let us denote the line bundle introduced in ~\ref{pjat' nu} for the curve $C_0$
by $\CL'_{C_0}$.

\subsubsection{Lemma} For any $\mu_1,\mu_2\in X_\ell$ such that
$\mu_1+\mu_2=-2\nu_0$ we have a canonical isomorphism
$$\ol{\CL'}_\heartsuit|_{p^{-1}(\fD)}=
(r_{\mu_1,\mu_2}\times\id)^*(\CL'_{C_0})_\heartsuit$$

{\em Proof.} Just as in the proof of the Lemma ~\ref{pjat' mho} we define
the extension of the line bundle $\CL_\alpha$ on $C^\alpha$ to the
line bundle $\ol{\CL}_\alpha$ on $C^\alpha\cup (C_0-x_1-x_2)^\alpha$
for any $\alpha\in\BN[X_\ell\oplus Y_\ell]$.
The extension is given by the same formula ~(\ref{pjat' formula}).

We necessarily have
$$\AJ_\alpha^*\ol{\CL'}=\ol{\CL}_\alpha(n(C_0-x_1-x_2)^\alpha)$$
for some $n\in\BZ$ where $(C_0-x_1-x_2)^\alpha$ is viewed as a divisor in
$\ol{C}^\alpha$.

On the other hand,
$$\ol{\CL}_\alpha|_{p^{-1}(\fD)}=\AJ_\alpha^*\pi^*r_{\mu_1,\mu_2}^*
(\CL'_{C_0})\otimes\CN^m$$
for some $m\in\BZ$ where $\CN$ denotes the normal line bundle to
$(C_0-x_1-x_2)^\alpha$ in $C^\alpha\cup (C_0-x_1-x_2)^\alpha$.

The line bundle $\CN$ is lifted from the base $\fD$.
Hence $(L\otimes\CN)_\heartsuit=L_\heartsuit$ for any line bundle $L$ on
$(\Pic(C_0/\fD)\otimes(X_\ell\oplus Y_\ell))_{(2g-2)\nu_0,0}$.

Finally, two line bundles on
$(\Pic(C_0/\fD)\otimes(X_\ell\oplus Y_\ell))_{(2g-2)\nu_0,0}$
are isomorphic iff their inverse images with respect to $\AJ_\alpha$
are isomorphic for any $\alpha$.

This completes the proof of the Lemma. $\Box$

\subsection{}
Now we can use the above Lemma together with
~(\ref{pjat' gift}) to define for any $\mu\in X_\ell/Y_\ell$
the morphism of algebras
$$\varphi_\mu:\ \ol{\fM}|_{p^{-1}(\fD)}\lra\End(\CE_\mu)$$
It remains to show that the morphism
$$\varphi:=\sum_{\mu\in X_\ell/Y_\ell}\varphi_\mu:\
\ol{\fM}|_{p^{-1}(\fD)}\lra\bigoplus_{\mu\in X_\ell/Y_\ell}\End(\CE_\mu)$$
is an isomorphism. It is enough to prove that $\varphi$ is surjective.

Since each morphism $\varphi_\mu$ is obviously surjective, it suffices to
check that all the idempotents of $\bigoplus_{\mu\in X_\ell/Y_\ell}\End(\CE_\mu)$
lie in the image of $\varphi$.

Consider the subgroup $Z\subset\ol{G}|_{p^{-1}(\fD)}$,
$$Z:=\operatorname{Ker}\pi\cap\ol{G}|_{p^{-1}(\fD)}$$
We have $Z\simeq\ul{X_\ell/Y_\ell}$.

Clearly, the morphism $\varphi_\mu$ maps $p_{Z*}i_Z^*\ol{\CL'}_\heartsuit\subset
p_{\ol{G}*}i_{\ol{G}}^*\ol{\CL'}_\heartsuit$ to the one-dimensional subspace
of $\End(\CE_\mu)$ generated by $\id$. Let us denote the corresponding linear
functional on $p_{Z*}i_Z^*\ol{\CL'}_\heartsuit$ by $\aleph_\mu$.

It is enough to prove that the set of functionals $\{\aleph_\mu,\mu\in
X_\ell/Y_\ell\}$ is linearly independent.

Note that the quadratic form $?\cdot?$ on $X_\ell$ induces a perfect pairing
$(X_\ell/Y_\ell)\times(X_\ell/Y_\ell)\lra\BQ/l\BZ$, and thus identifies the
finite abelian group $X_\ell/Y_\ell$ with its dual $(X_\ell/Y_\ell)^\vee$.

Since $p_{Z*}i_Z^*\ol{\CL'}_\heartsuit$
is $X_\ell/Y_\ell$-graded, the group $(X_\ell/Y_\ell)^\vee$
acts on it by multiplication, and the corresponding module is isomorphic
to the dual regular representation $B[X_\ell/Y_\ell]$.

Consider $\mu=0\in X_\ell/Y_\ell$.
The functional $\aleph_0$ does not vanish on any graded component of
$p_{Z*}i_Z^*\ol{\CL'}_\heartsuit$. Hence 
$\aleph_0$ generates $(p_{Z*}i_Z^*\ol{\CL'}_\heartsuit)^\vee$ as
a $(X_\ell/Y_\ell)^\vee$-module. On the other hand, it is easy to see that
for any $\lambda\in(X_\ell/Y_\ell)^\vee$ we have $\lambda(\aleph_0)=
\aleph_\lambda$.

This completes the proof of the property ~\ref{pjat' M}(iii) along with
the statement the Proposition ~\ref{pjat' weak}. $\Box$

\subsection{}
Now we will derive the Theorem ~\ref{pjat' modular Heisenberg}(b) from the
Proposition ~\ref{pjat' weak}.

Given $\mu\in X_\ell/Y_\ell$ we will denote $-2\nu_0-\mu$ by $\mu'$.

We have proved that $\bSp\CE=\oplus_{\mu\in X_\ell/Y_\ell}\CE_\mu$, and
the flat projective connections on $\CE,\CE_\mu$ agree (note that the
sum of bundles with flat projective connections does not have a natural
projective connection).

Recall that in order to define the local system $\CH^{\alpha+\mu+\mu'}$
we have used the canonical section $s_{\alpha+\mu+\mu'}$ of the line
bundle $\AJ^*_{\alpha+\mu+\mu'}\det(\CE_\mu)$.

The Theorem follows immediately from the next two statements.

\subsection{Lemma}
The canonical section $s_\alpha$ of $\AJ_\alpha^*\det(\CE)$ can be
extended to a section $\ol{s}_\alpha$ over $\ol{\CM}^\circ$ in such way that
$$\ol{s}_\alpha|_{p^{-1}(\fD)}=\prod_{\mu\in X_\ell/Y_\ell}s_{\alpha+\mu+\mu'}$$

The proof is completely parallel to the one of ~\ref{pjat' universal}. $\Box$

\subsection{}
We consider $\oplus_{\mu\in X_\ell/Y_\ell}\CE_\mu$ as a projectively flat
bundle via the isomorphism with $\bSp\CE$. Then
$\oplus_{\mu\in X_\ell/Y_\ell}\det(\CE_\mu)$ is a direct summand of an
exterior power of $\oplus_{\mu\in X_\ell/Y_\ell}\CE_\mu$ and thus inherits
a flat projective connection.

{\bf Proposition.} The collection $(s_{\alpha+\mu+\mu'},\
\mu\in X_\ell/Y_\ell)$ forms a projectively flat section of
$\oplus_{\mu\in X_\ell/Y_\ell}\det(\CE_\mu)$.

\subsection{Proof}
The flat projective connection on
$\oplus_{\mu\in X_\ell/Y_\ell}\det(\CE_\mu)$ induces a flat genuine
connection on $\det(\CE_\mu)\otimes\det(\CE_\lambda)^{-1}$ for any
$\mu,\lambda\in X_\ell/Y_\ell$. We have to prove that $s_\mu s_\lambda^{-1}$
is a flat section of $\det(\CE_\mu)\otimes\det(\CE_\lambda)^{-1}$.

It suffices to prove that $s_\mu^N s_\lambda^{-N}$ is a flat section of
$(\det(\CE_\mu)\otimes\det(\CE_\lambda)^{-1})^N$ for some positive
integer $N$.

This follows immediately from the next two Lemmas.

\subsubsection{Lemma}
\label{pjat' trivial}
There exists $N>0$ such that the one-dimensional local system
$(\det(\CE_\mu)\otimes\det(\CE_\lambda)^{-1})^N$ is trivial for any
$\mu,\lambda\in X_\ell/Y_\ell$.

\subsubsection{Lemma}
\label{pjat' constant}
Any invertible function on $\TO C_0^{\alpha+\mu+\lambda}\times_{\CM_{g_0}}
\dot{\delta}$ is constant if $g_0>1$.

\subsubsection{}
We will prove the Lemmas in the following subsections. Let us derive the
Proposition now. First of all, if $g_0=0$ there is nothing to prove.
The case $g_0=1$ is left to the reader.
Otherwise, we conclude that $s_\mu^N s_\lambda^{-N}$ being invertible
over $\TO C_0^{\alpha+\mu+\lambda}\times_{\CM_{g_0}}\dot{\delta}$ is a
constant function, and thus a flat section of a trivial local system
$(\det(\CE_\mu)\otimes\det(\CE_\lambda)^{-1})^N.\ \Box$

\subsection{}
To prove the Lemma ~\ref{pjat' trivial} it is enough to show that the image
of the monodromy representation of the fundamental group
$\pi_1(\TB_{{\cal T}C_{\fD}^\alpha
\times_{\ol{\CM_g}}\ol{\delta_g}}(\TO\ol{C}^\alpha
 \times_{\ol{\CM_g}}\ol{\delta_g}))$ (see ~(\ref{pjat' map b}))
in $PGL(\oplus_{\mu\in X_\ell/Y_\ell}\CE_\mu)$ is finite.

To this end it suffices to show that the image of $\pi_1(\TO C^\alpha)$
in $PGL(\CE)$ is finite.

Recall the notations of ~\ref{pjat' useful}, ~\ref{pjat' useless}.
Let us choose a finite subgroup $\tilde{G}\subset
p_{G*}i_G^*\dot{\CL'_\heartsuit}$
projecting surjectively onto $G$ and such that the kernel of this projection
contains at least three elements.

Then it follows from the discussion in ~\ref{pjat' useful} that
the action of $\pi_1(\TO C^\alpha)$ on $\CE$ factors through the group
$\operatorname{Aut}(\tilde{G})$.

This proves the Lemma ~\ref{pjat' trivial}. $\Box$

\subsection{}
Let us prove the Lemma ~\ref{pjat' constant}.
To unburden the notation we will omit the subindex $_0$.
Assume there is a nonconstant invertible function.
We can lift it to the nonsymmetrized Cartesian power and obtain a nonconstant
invertible function on $(\TB C)^J\times_{\CM_g}\dot{\delta}$ for
some set $J$ of positive cardinality $n$.

Let us choose a subset $K=J-j$ of cardinality
$n-1$ and consider the projection
$$\pr:\ (\TB C)^J\times_{\CM_g}\dot{\delta}\lra
(\TB C)^K\times_{\CM_g}\dot{\delta}$$

\subsubsection{Claim}
$\pr_*\CO^*=\CO^*$.

{\em Proof.}
Consider a projection
$$\pr':\ (\TB C)^K\times_{\CM_g}\dot{\delta}\times_{\CM_g}C_j\lra
(\TB C)^K\times_{\CM_{g_0}}\dot{\delta}$$

Then we have
$$\pr_*\CO^*=\bigoplus_{i\in\BZ}\pr'_*(\dot{\Omega^i}(C_j/\CM_g))$$

But for generic $(C,z_k\in C,\ k\in K)$ the subgroup of $\Pic(C)$
generated by $\{\Omega,(z_k),\ k\in K\}$ is free. That is, for any $i\in\BZ$
there are  no meromorphic sections of $\Omega^i$ invertible on
$C-\{z_k,\ k\in K\}$.

This completes the proof of the Claim. $\Box$

\subsubsection{}
Using the Claim inductively we conclude that
$$\Gamma((\TB C)^J\times_{\CM_g}\dot{\delta},\CO^*)=\Gamma(\dot{\delta},\CO^*)=
\bigoplus_{i\in\BZ}\Gamma(\CM_g,\dot{\delta^i})$$

Since $g>1$, we know that for $i\not=0$ the line bundle $\delta^i$ is nontrivial,
and hence does not have any invertible sections.

For $i=0$ we know that that the only invertible functions on $\CM_g$
are constants (see e.g. ~\cite{bm} 2.3).

This completes the proof of Lemma ~\ref{pjat' constant} along with the
Theorem ~\ref{pjat' modular Heisenberg}. $\Box$

\newpage
\begin{center}
{\bf Chapter 3. Regular Representation}
\end{center}
\vspace{.8cm}

\bigskip

\section{A characterization of the Regular bimodule}

This section belongs really to the Chapter 2 of IV.

\subsection{}
\label{pjat' bR}
Consider the algebra $\tfu$ defined in ~\cite{ajs}, Remark 1.4.
It is an infinite dimensional $\sk$-algebra containing $\fu$ as a
subalgebra. We have $\tfu=\fu^-\otimes\tfu^0\otimes\fu^+$ where
$\tfu^0$ is formed by some divided power expressions in $K_i$.

We have ${\CC}=\tfu-\bmod$, and $\tfu$ itself viewed as a left regular
$\fu$-module is isomorphic to $\oplus_{\lambda\in Y}\ \fu\otimes L(l\lambda)$
where $\fu$ is considered as a left regular $\fu$-module.

Thus $\tfu$ does not belong to ${\CC}$ being infinite dimensional, but belongs
to $\Ind{\CC}$, see ~\cite{d4}, \S 4.
If we take into account the right regular action of $\fu$
an $\tfu$ as well, then we can regard $\tfu$ as an object of
$\Ind{\CC}\otimes\Ind{\CC}_r$, that is, an infinite dimensional
$\fu$-bimodule. Finally, composing the right action with the antipode,
we may view $\tfu$ as an object of $\Ind{\CC}\otimes\Ind{\CC}$. We will
denote this object by $\bR$.

\subsection{}
\label{pjat' ker fun}
Given $V,W\in{\CC}$ we define a functor
$F_{V\otimes W}:{\CC}\lra{\CC}$,
$$
F_{V\otimes W}(M)=(sM\otimes_{\CC}V)\otimes W.
$$
The functor ${\CC}\times{\CC}\lra\Fun({\CC},{\CC}),\
(V,W)\mapsto F_{V\otimes W}$, is exact in $W$ and right exact in $V$.
Therefore, by the universal property it extends to a functor
$$
{\CC}\otimes{\CC}\lra\Fun({\CC},{\CC}),\ \CO\mapsto F_{\CO}.
$$
This in turn extends to a functor
$$
\Ind{\CC}\otimes\Ind{\CC}\lra\Fun(\Ind{\CC},\Ind{\CC}),\ \CO\mapsto F_{\CO}.
$$
We have $F_\bR=\Id$.

In concrete terms, ${\CC}\otimes{\CC}$ is equivalent to ${X_\ell}\times {X_\ell}$-graded
$\fu$-bimodules. This equivalence sends $V\otimes W$ to
$V\otimes sW$ with an evident bimodule structure
(the latter tensor product is understood as a product of vector spaces).
For a bimodule $\CO$ and $M\in{\CC}$ we have
$$
F_{\CO}(M)=s(sM\otimes_{\CC}\CO).
$$
Here $sM\otimes_{\CC}\CO$ has a right $\fu$-module structure and ${X_\ell}$-grading
inherited from $\CO$. So, for the bimodule $\bR$ we have
$$
s(sM\otimes_{\CC}\bR)=M
$$
for all $M\in{\CC}$.

\subsection{}
\label{pjat' kernel functor}
Let ${\CC}^{\BZ}$ denote the category of $\BZ$-graded
objects of ${\CC}$. For an ${X_\ell}\times {X_\ell}$-graded $\fu$-bimodule
$\CO$ we define a functor
$F_{\CO;\infh+\blt}:{\CC}\lra{\CC}^{\BZ}$,
$$
F_{\CO;\infh+\blt}(M):=s\Tor^{\CC}_{\infh+\blt}(sM,\CO).
$$
Here
$\Tor^{\CC}_{\infh+\blt}(sM,\CO)$ has a right $\fu$-module structure
and ${X_\ell}$-grading inherited from $\CO$.

\subsection{} {\bf Lemma.} (i) {\em
Suppose $\CO$ is $\fu^+$-induced as a left $\fu$-module.
Then $F_{\CO;\infh+n}=0$ for $n<0$, and $F_{\CO;\infh+0}=F_{\CO}$.}

(ii) {\em Suppose $\CO$ is projective as a left $\fu$-module.
Then $F_{\CO;\infh+n}=0$ for $n\neq 0$, and $F_{\CO;\infh+0}=F_{\CO}$.}

{\bf Proof} follows immediately from the next

\subsubsection{} {\bf Lemma.} {\em Let $V\in{\CC}_r$, and assume that
$N\in{\CC}$ is $\fu^+$-induced. Then
$$
\Tor^{\CC}_{\infh+\blt}(V,N)=\Tor^{\CC}_{\blt}(V,N).
$$

In particular, if $N$ is projective then
$$
\Tor^{\CC}_{\infh+\blt}(V,N)=\left\{\begin{array}{ll}
                                      0 & \mbox{ if $n\neq 0$}\\
                                   V\otimes_{\CC}N & \mbox{ if $n=0$.}
                                    \end{array}
                              \right.
$$}

{\bf Proof.} According to Theorem IV.4.9 we may choose
any $\fu^+$-induced right convex resolution of $N$ for the
computation of $\Tor^{\CC}_{\infh+\blt}(V,N)$. So, let us choose $N$ itself.
Then we have
$$
H^{\blt}(R^{\blt}_{\swarrow}(V)\otimes_{\CC}N)=
H^{\blt}((R^{\blt}_{\swarrow}(V)\otimes sN)\otimes_{\CC}\sk).
$$
But $R^{\blt}_{\swarrow}(V)\otimes sN$ is a left $\fu$-projective
resolution of $V\otimes sN$, so
$$
H^{-\blt}((R^{\blt}_{\swarrow}(V)\otimes sN)\otimes_{\CC}\sk)=
\Tor^{\CC}_{\blt}((V\otimes sN),\sk)=
\Tor^{\CC}_{\blt}(V,N).\ \Box
$$

\subsection{Corollary} {\em
$$
F_{\bR;\infh+n}=\left\{\begin{array}{ll}
                            0 & \mbox{ if $n\neq 0$}\\
                            \Id & \mbox{ if $n=0$.}
                                    \end{array}
                              \right.
$$}

{\bf Proof.} The bimodule $\bR$ is projective as a left $\fu$-module. $\Box$

\subsection{Lemma.} {\em Let $M_1,\ldots, M_n,N\in{\CC}$. There is a canonical
perfect pairing between
$\Ext^{\infh+k}_{\CC}(M_1\otimes\ldots\otimes M_n,N)$ and
$\Tor^{\CC}_{\infh+k}(\sk,M_1\otimes\ldots\otimes M_n\otimes N^*)$
for every $k\in\BZ$.}

{\bf Proof.} This is just IV.4.8. $\Box$

\subsection{Converse Theorem}
\label{pjat' Converse Theorem}
Suppose $Q\in\operatorname{Ind}{\CC}\otimes
\operatorname{Ind}{\CC}$, and an isomorphism of functors
$\phi:\ F_{\bR,\frac{\infty}{2}+\bullet}\iso F_{Q,\frac{\infty}{2}+\bullet}$
is given. Then $\phi$ is induced by the isomorphism $\bar{\phi}:\ \bR\iso Q$.

{\em Proof.} We start with the following Lemma:

\subsubsection{Lemma} Let $M\in\operatorname{Ind}{\CC}$. Suppose
$\Ext_{\CC}^{\frac{\infty}{2}+n}(M,N)=0$ for $n\not=0$ and arbitrary $N\in{\CC}$.
Suppose further that $\Ext_{\CC}^{\frac{\infty}{2}+0}(M,N)$ is exact in $N$.
Then $M$ is projective, and  $\Ext_{\CC}^{\frac{\infty}{2}+0}(M,N)=\Hom_{\CC}(M,N)$.

{\em Proof.} By the Theorem IV.4.9, $\Ext_{\CC}^{\frac{\infty}{2}+\bullet}(M,N)=
H^\bullet(\Hom_CC(P^\bullet(M),R^\bullet(N)))$ where $P^\bullet(M)$ is some
$\fu^+$-induced convex right resolution of $M$, and $R^\bullet(N)$ is some
$\fu^-$-induced concave right resolution of $N$.

Since $\Ext_{\CC}^{\frac{\infty}{2}+n}(M,?)$ vanishes for $n\not=0$, and
is exact for $n=0$, we deduce that $\Ext_{\CC}^{\frac{\infty}{2}+\bullet}(M,N)=
H^\bullet(\Hom_{\CC}(P^\bullet(M),N))$. In particular,
$$\Ext_{\CC}^{\frac{\infty}{2}+0}(M,N)=\operatorname{Coker}
(\Hom_{\CC}(P^1(M),N)\overset{d_0^*}{\lra}\Hom_{\CC}(P^0(M),N))$$
Let us denote this latter functor in $N$ by $h(N)$ for brevity.
We denote by $\xi^*:\ h(N)\lra\Hom_{\CC}(M,N)$ the evident morphism of
functors. By the assumption, $h(N)$ is exact and hence representable by
some projective $P^0\in\operatorname{Ind}{\CC}$. We denote by $\xi:\ M\lra P^0$
the map inducing $\xi^*$ on the functors they represent.

We denote by $\theta:\ P^0\lra P^0(M)$ the map inducing the evident
projection $\Hom_{\CC}(P^0(M),N)\lra h(N)$ on the functors they represent.
Then $\theta\circ\xi=\varepsilon:\ M\lra P^0(M)$ where $\varepsilon$ stands
for the augmentation.

Now $\varepsilon$ is injective, hence $\xi:\ M\lra P^0$ is injective. We
may extend $\xi:\ M\lra P^0$ to a $\fu^+$-induced convex right resolution
$P^0\lra P^1\lra\ldots$ of $M$, and use this resolution instead of
$P^\bullet(M)$ for the calculation of
$\Ext_{\CC}^{\frac{\infty}{2}+\bullet}(M,N)$. Repeating the above argument
for $P^\bullet$ we conclude that
$$\operatorname{Coker}
(\Hom_{\CC}(P^1,N)\overset{d_0^*}{\lra}\Hom_{\CC}(P^0,N))=\Hom_{\CC}(P^0,N),$$
i.e. $d^*_0=0$, whence $0=d_0:\ P^0\lra P^1$, whence $M=P^0$. $\Box$

\subsubsection{Proof of the Theorem}
Applying the above Lemma together with the Lemma IV.4.8, we see that
$Q$ is projective as a left $\fu$-module, and
$s\Tor^{\CC}_{\frac{\infty}{2}+n}(sN,Q)=s\Tor^{\CC}_n(sN,Q)=
s(sN\otimes_{\CC} Q)$ if $n=0$, and zero otherwise for any $N\in{\CC}$.

So $\phi$ boils down to the isomorphism of functors in $N$:
$$s(sN\otimes_{\CC} \bR)\iso s(sN\otimes_{\CC} Q)$$
Applying the Lemma IV.4.8 again, we obtain the isomorphism of functors
$${\CC}\lra{\CC}:\ s(\Hom_{\CC}(\bR,?))\iso s(\Hom_{\CC}(Q,?)),$$
whence the isomorphism of representing objects. $\Box$

\section{The adjoint representation}

\subsection{}
\label{pjat' ad}
For $\mu\in Y_\ell\subset {X_\ell}$ (see IV.9.1.2)
we denote by $T_\mu$ an autoequivalence
$\CC\lra\CC$ given by twisting by $L(\mu):\ T_\mu(N):= N\otimes L(\mu)$.

For $\mu,\nu\in Y_\ell\subset {X_\ell}$ we will consider an autoequivalence
$T_\mu\otimes T_\nu$ of $\Ind\CC\otimes\Ind\CC$.

The objects $\bR$ and $T_\mu\otimes T_{-\mu}\bR$ give rise to the same
functor $F_{\bR}=F_{T_\mu\otimes T_{-\mu}\bR}=\Id:\ \Ind\CC\lra\Ind\CC$
(notations of ~\ref{pjat' ker fun}). The equality of functors is induced by
the isomorphism of objects $t_\mu:\ \bR\iso T_\mu\otimes T_{-\mu}\bR$.

Consider the tensor product functor $\otimes:\
\Ind\CC\otimes\Ind\CC\lra\Ind\CC$.

Let us denote by $\hat{\ad}$ the object $\otimes(\bR)\in\Ind\CC$.
For $\mu\in Y_\ell\subset {X_\ell}$ we have $\otimes(T_\mu\otimes T_{-\mu}\bR)=
\hat{\ad}$. Thus the isomorphism $t_\mu$ induces the same named
{\em automorphism} of $\hat{\ad}$. This way we obtain the action
of the group $Y_\ell$ on $\hat{\ad}$.

We define the {\em adjoint representation} $\ad$ as the quotient
$\hat{\ad}/Y_\ell$.

It is easy to see that $\ad\in\CC\subset\Ind\CC$, and in fact
$\hat{\ad}=\ad\otimes \sk[Y_\ell]$ with the trivial action of $\fu$
on the second multiple, and the trivial action of $Y_\ell$ on the first multiple.

In concrete terms, $\ad$ coincides with $\fu$ considered as a $\fu$-module
via the adjoint action. The ${X_\ell}$-grading on $\ad$ is compatible with
multiplication on $\fu$, and the generators have the following weights:
$$\theta_i\in(\ad)_{-i'};\ \epsilon_i\in(\ad)_{i'};\ K_i\in(\ad)_0.$$

Note that the highest weight of $\ad$ is equal to $2\rho_\ell$ (notations
of IV.9.1.2).

\subsection{}
\label{pjat' maximal summand}
Recall the notations of IV.9.1, 9.2. Recall that $\tCO_\zeta$ is a
semisimple subcategory of $\CC$ consisting of all the direct sums of simples
$L(\lambda),\lambda\in\Delta=\Delta_\ell$.

The natural embedding $\tCO_\zeta\hra\CC$ admits the right adjoint functor
$N\mapsto N^\Delta$ (the maximal $\tCO_\zeta$-subobject), as well as the
left adjoint functor $N\mapsto N_\Delta$ (the maximal $\tCO_\zeta$-quotient).
Finally, $$N\mapsto\langle N\rangle_\Delta$$ denotes the image functor of
the canonical morphism from the right adjoint to the left adjoint.

One checks easily that if $M\subset N$ is a direct summand belonging to
$\tCO_\zeta$ and not properly contained in any of the kind (i.e. maximal),
then we have a canonical isomorphism $\langle N\rangle_\Delta\iso M$.

For this reason we call $\langle N\rangle_\Delta$ {\em the maximal
$\tCO_\zeta$-direct summand of $N$}.

H.Andersen and J.Paradowski introduced the tensor structure $\tilde\otimes$
on the category $\tCO_\zeta$ (see ~\cite{ap}).
Arguing like in the proof of Lemma IV.9.3,
one can show easily that for $M,N\in\tCO_\zeta\subset\CC$ we have
$$M\tilde{\otimes}N=\langle M\otimes N\rangle_\Delta$$

\subsubsection{Definition}
\label{pjat' ss ad}
{\em The semisimple adjoint representation}
is the following object of $\tCO_\zeta$:
$$\tilde{\ad}:=
\bigoplus_{\lambda\in\Delta}L(\lambda)\tilde{\otimes}L(\lambda)^*$$
Here $*$ is the rigidity (see e.g. IV.4.6.2).

\subsection{Theorem}
\label{pjat' try and prove}
$\tilde{\ad}=\langle\ad\rangle_\Delta$.

{\em Proof.} According to ~\cite{lm}, $\ad\in\CC$ represents the functor
$V\mapsto{\operatorname{Nat}}(\Id,\Id\otimes V)$ sending $V\in\CC$ to the
vector space of natural transformations between endofunctors $\CC\lra\CC$.

Hence $\langle\ad\rangle_\Delta$ represents the functor sending
$V\in\tCO_\zeta$ to the vector space of natural transformations between
endofunctors $\Id$ and $\langle\Id\otimes V\rangle_\Delta:\ \tCO_\zeta\lra
\tCO_\zeta$.

Similarly, $\tilde{\ad}\in\tCO_\zeta$ represents the functor
$W\mapsto{\operatorname{Nat}}(\Id,\Id\tilde{\otimes}W)$.

It remains to recall that for $W\in\tCO_\zeta\subset\CC$ we have
$\Id\tilde{\otimes}W=\langle\Id\otimes W\rangle_\Delta.\ \Box$

\newpage
\begin{center}
{\bf Chapter 4. Quadratic degeneration in genus zero}
\end{center}
\vspace{.8cm}

\bigskip

This chapter belongs really to the part IV. Here we construct a sheaf
$\CR\in\Ind\FS\otimes\Ind\FS$ corresponding to the regular bimodule
$\bR\in\Ind\CC\otimes\Ind\CC$ under the equivalence $\Phi$.

\section{I-sheaves}

Since $\CR$ is an object of $\Ind\FS\otimes\Ind\FS$ we start with
developing a little machinery to describe certain Ind-sheaves.

\subsection{} We consider the root datum $Y^\diamondsuit:=Y\times Y;$

$X^\diamondsuit=X\times X;$

$I^\diamondsuit:=I\sqcup I;$

$\langle(y_1,y_2),(x_1,x_2)\rangle:=\langle y_1,x_1\rangle+\langle y_2,x_2
\rangle$;

$I^\diamondsuit\ni(i_1,i_2)\mapsto(i_1,i_2)\in Y^\diamondsuit;$

$I^\diamondsuit\ni(i_1,i_2)\mapsto(i_1',i_2')\in X^\diamondsuit.$

The corresponding category $\FS^\diamondsuit$ is equivalent to $\FS\otimes\FS$.
Recall that an object $\CX$ of $\FS^\diamondsuit$ is the following collection
of data (see III.4.2):

(a) a weight $\lambda(\CX)=(\lambda_1,\lambda_2)\in {X_\ell}\times {X_\ell}$;

(b) for each $(\alpha_1,\alpha_2)\in\BN[I]\times\BN[I]=\BN[I^\diamondsuit]$,
a sheaf
$\CX^{\alpha_1,\alpha_2}\in\CM(\CA^{\alpha_1,\alpha_2}_{\lambda_1,\lambda_2};
\CS)$;

Note that $(\CA^{\alpha_1,\alpha_2}_{\lambda_1,\lambda_2},\CS)=
(\CA^{\alpha_1}_{\lambda_1},\CS_1)\times
(\CA^{\alpha_2}_{\lambda_2},\CS_2)$;

we will denote by $\CX^{\valpha_1,\valpha_2}(\vd)$ perverse sheaves over
$\CA^{\valpha,\valpha_2}_{\lambda_1,\lambda_2}(\vd)$
obtained by taking the restrictions
with respect to the embeddings
$\CA^{\valpha_1,\valpha_2}_{\lambda_1,\lambda_2}(\vd)\hra
\CA^{\alpha_1,\alpha_2}_{\lambda_1,\lambda_2}$;

(c) for each $(\alpha_1,\alpha_2),(\beta_1,\beta_2)\in \BN[I]\times\BN[I],\
 d>0$, a {\em factorization isomorphism}
\begin{equation}
\psi^{\alpha_1,\alpha_2;\beta_1,\beta_2}(d):
\CX^{(\alpha_1,\alpha_2;\beta_1,\beta_2)}(d)
\iso
\CID^{(\alpha_1,\alpha_2;0,0)}_{\lambda_1-\beta_1',\lambda_2-\beta_2'}(d)
\boxtimes
\CX^{(0,0;\beta_1,\beta_2)}(d)
\end{equation}

satisfying the usual associativity conditions.

Note that $\CA^{\alpha_1,\alpha_2;0,0}(d)=\CA^{\alpha_1;0}(d)\times
\CA^{\alpha_2;0}(d)$, and

$\CID^{(\alpha_1,\alpha_2;0,0)}_{\lambda_1-\beta_1',\lambda_2-\beta_2'}(d)=
\CID^{(\alpha_1;0)}_{\lambda_1-\beta_1'}(d) \boxtimes
\CID^{(\alpha_2;0)}_{\lambda_2-\beta_2'}(d)$.

\subsection{Definition}
\label{pjat' I-sheaf}
An {\em I-sheaf} is the following collection
of data:

(a) a weight $\chi\in {X_\ell}$, to be denoted by $\chi({X_\ell})$;

(b) for each $(\lambda_1,\lambda_2)\in {X_\ell}\times {X_\ell}$ such that
$\lambda_1+\lambda_2-\chi=\gamma\in\BN[I]$, and
$\alpha_1,\alpha_2\in\BN[I]$, such that $\alpha_1,\alpha_2\leq\gamma$,
a sheaf
$\CX^{\alpha_1,\alpha_2}_{\lambda_1,\lambda_2}\in
\CM(\CA^{\alpha_1,\alpha_2}_{\lambda_1,\lambda_2};
\CS)$;

we will denote by $\CX^{\valpha_1,\valpha_2}(\vd)$ perverse sheaves over
$\CA^{\valpha,\valpha_2}_{\lambda_1,\lambda_2}(\vd)$
obtained by taking the restrictions
with respect to the embeddings
$\CA^{\valpha_1,\valpha_2}_{\lambda_1,\lambda_2}(\vd)\hra
\CA^{\alpha_1,\alpha_2}_{\lambda_1,\lambda_2}$;

(c) for each $\lambda_1,\lambda_2\in {X_\ell}$ with $\lambda_1+\lambda_2-\chi
=\gamma\in\BN[I]$, and
$\alpha_1,\alpha_2,\beta_1,\beta_2\in \BN[I]$ such that $\alpha_1+\beta_1,
\alpha_2+\beta_2\leq\gamma,\
 d>0$, a {\em factorization isomorphism}
\begin{equation}
\label{pjat' factor iso I-sh}
\psi^{\alpha_1,\alpha_2;\beta_1,\beta_2}(d):
\CX^{(\alpha_1,\alpha_2;\beta_1,\beta_2)}(d)
\iso
\CID^{(\alpha_1,\alpha_2;0,0)}_{\lambda_1-\beta_1',\lambda_2-\beta_2'}(d)
\boxtimes
\CX^{(0,0;\beta_1,\beta_2)}(d)
\end{equation}

satisfying the usual associativity conditions;

(d) for each $(\lambda_1,\lambda_2)\in {X_\ell}\times {X_\ell}$ such that
   $\lambda_1+\lambda_2-\chi=\gamma\in\BN[I]$, and
   $\alpha_1,\alpha_2\in\BN[I]$ with $\alpha_1,\alpha_2\leq\gamma$,
and $\beta_1,\beta_2\in\BN[I]$, an isomorphism

$\CX^{\alpha_1+\beta_1,\alpha_2+\beta_2}_{\lambda_1+\beta_1',\lambda_2+\beta_2'}
\iso\sigma_*\CX^{\alpha_1,\alpha_2}_{\lambda_1,\lambda_2}$

satisfying the usual compatibilities.

Here $\sigma$ stands for the closed embedding
$\CA^{\alpha_1,\alpha_2}_{\lambda_1,\lambda_2}\hra
\CA^{\alpha_1+\beta_1,\alpha_2+\beta_2}_{\lambda_1+\beta_1',\lambda_2+\beta_2'}$.

\subsubsection{Remark}
\label{pjat' part of data}
Note that an I-sheaf $\CX$ can be uniquely reconstructed from the partial set
of data (b),(c),(d), namely the data (b),(c),(d) given only for $\alpha_1=
\alpha_2,\beta_1=\beta_2$. In what follows we will describe I-sheaves by
these partial data.

\section{Degeneration of quadrics}

\subsection{}
\label{pjat' family}
We recall the construction of ~\cite{kl}II, 15.2, 15.3.
Consider $\BP^1\times\BP^1\times\BA^1$ with coordinates $p,q,t$. Consider
the subvariety $Q\subset\BP^1\times\BP^1\times\BA^1$ given by the equation
$pq=t$. The fiber $Q_t$ of $Q$ over $t\in\BA^1$ is a projective line if
$t\not=0$, and the union of two projective lines if $t=0$.

Each fiber $Q_t$ has two marked points: $(p=0,q=\infty)$, and $(p=\infty,q=0)$.
For $t=0$ each irreducible component of $Q_0$ has exactly one marked point,
and the marked points lie away from the singularity of $Q_0$. The irreducible
components of $Q_0$ are denoted by $Q_{0v}$ (for vertical), and $Q_{0h}$
(for horizontal).

There are two maps $f_1,f_2$ from the standard projective
line $\BP^1_{st}$ to $Q_t$, taking $0\in\BP^1_{st}$ to the first (resp.
second) marked point on $Q_t$, and mapping $\BP^1_{st}$ isomorphically onto
the irreducible component of $Q_t$ containing this marked point. If $z$
is the standard coordinate on $\BP^1_{st}$, then
$$f_1(z)=(\frac{tz}{z-1},\frac{z-1}{z});\
f_2(z)=(\frac{z-1}{z},\frac{tz}{z-1}).$$
Thus, restricting to the open subset $'Q\subset Q$ given by the inequality
$t\not=0$, we obtain a map $\BC^*\lra \tCP^2$.

The limit for $t\lra 0$ is a boundary point of $\tCP^2$ (notations of IV.2.3).
Composing with the ``1-jet at 0'' projection $\tCP^2\lra T\CP^{\circ 2}$
(notations of IV.2.7) we obtain a map $\BC^*\lra T\CP^{\circ 2}$; and the
limit for $t\lra 0$ is a boundary point of $T\CP^{\circ 2}$.

\subsection{}
For $\alpha\in\BN[I]$ let $Q^\alpha\overset{\pi^\alpha}{\lra}{\BA^1}$
denote the space of relative configurations: $(\pi^\alpha)^{-1}(t)=Q^\alpha_t$.

\subsubsection{} $Q^\alpha$
contains the open subset $'Q^\alpha:=(\pi^\alpha)^{-1}(\BC^*)$.

$Q^{\bullet\alpha}\subset Q^\alpha$ is the open subset of configurations
where none of the points equals the marked ones.

$Q^{\circ\alpha}\subset Q^{\bullet\alpha}$ is the open subset of configurations
where all the points are distinct.

$'Q^{\bullet\alpha}:=\ 'Q^\alpha\cap Q^{\bullet\alpha};\ \ \
'Q^{\circ\alpha}:=\ 'Q^\alpha\cap Q^{\circ\alpha}$.

The open inclusion $'Q^{\circ\alpha}\hra\ 'Q^{\bullet\alpha}$ is denoted by $j$.

The projection $Q^{\bullet\alpha}\lra\BA^1$ is denoted by $\pi^{\bullet\alpha}$.

\subsubsection{} Let $\mu_1,\mu_2\in {X_\ell}$ be a pair of weights such that the
triple $(\mu_1,\mu_2,\alpha)$ is admissible (see IV.3.1).

The above construction defines the map $\theta^\alpha:\ 'Q^{\circ\alpha}\lra
\CP^{\circ\alpha}_{\mu_1\mu_2}$ (notations of IV.3.2), and we define the
local system $\CJ^\alpha_{\mu_1\mu_2}$ on $'Q^{\circ\alpha}$ as
$\theta^{\alpha*}\CI^\alpha_{\mu_1\mu_2}$ (notations of IV.3.4).

We denote by $\CJ^{\bullet\alpha}_{\mu_1\mu_2}$ the perverse sheaf
$j_{!*}\CJ^\alpha_{\mu_1\mu_2}[\dim Q^\alpha]$ on $'Q^{\bullet\alpha}$.

\subsection{} The aim of this Chapter is to compute the nearby cycles
$\Psi_{\pi^{\bullet\alpha}}\CJ^{\bullet\alpha}_{\mu_1\mu_2}$ as a perverse
sheaf on $Q^{\bullet\alpha}_0$. The problem is not quite trivial, but it
appears that one can combine the sheaves
$\Psi_{\pi^{\bullet\alpha}}\CJ^{\bullet\alpha}_{\mu_1\mu_2}$
for different $\alpha,\mu_1,\mu_2$ into a single I-sheaf; and this I-sheaf
is already easy to compute.

\section{The I-sheaf $\CR$}

\subsection{}
\label{pjat' vh}
We denote by $Q^\bullet_0$ the open subset of the special fiber
$Q^\bullet_0:=Q_0-\{(0,\infty),(\infty,0)\}$ obtained by throwing away
the marked points. Evidently, $Q^\bullet_0$ is the union of two irreducible
components $\BA^1_v$ and $\BA^1_h$ given by the equations $p=0$ (vertical
component), and $q=0$ (horizontal component) respectively. They intersect
at the point $(0,0)\in Q^\bullet_0$.

The irreducible components of $Q^{\bullet\alpha}_0, \alpha\in\BN[I]$,
are numbered by decompositions $\alpha=\alpha_1+\alpha_2;\
\alpha_1,\alpha_2\in\BN[I]$. Namely,
$$Q^{\bullet\alpha}_0=\bigcup_{\alpha_1+\alpha_2=\alpha}\CA^{\alpha_1}_v\times
\CA^{\alpha_2}_h$$
Each irreducible component embeds into $\CA^\alpha_v\times\CA^\alpha_h$.

Namely, $(q_1,\ldots,q_{|\alpha_1|};p_1,\ldots,p_{|\alpha_2|})$ goes to
$(q_1,\ldots,q_{|\alpha_1|},\overbrace{0,\ldots,0}^{|\alpha_2|};
p_1,\ldots,p_{|\alpha_2|},\overbrace{0,\ldots,0}^{|\alpha_1|})$.

We see readily that these embeddings agree on the intersections of
irreducible components, whence they combine together to the closed
embedding

\begin{equation}
\varsigma^\alpha:\ Q^{\bullet\alpha}_0\hra\CA^\alpha_v\times\CA^\alpha_h
\end{equation}

\subsection{}
\label{pjat' CR}
We define an I-sheaf $\CR$ as follows.

(a) We set $\chi(\CR)=(l-1)2\rho$;

(b) For $\mu_1,\mu_2\in {X_\ell}$ such that $\mu_1+\mu_2-(l-1)2\rho=\alpha\in\BN[I]$,
we define

\begin{equation}
\CR^{\alpha,\alpha}_{\mu_1\mu_2}:=\varsigma^\alpha_*\Psi_{\pi^{\bullet\alpha}}
\CJ^{\bullet\alpha}_{\mu_2\mu_1}
\end{equation}

{\em (Note the reverse order of $\mu_1,\mu_2$!)}

$\CR^{\alpha,\alpha}_{\mu_1\mu_2}$ is a perverse sheaf on
$\CA^\alpha_{\mu_1}\times\CA^{\alpha}_{\mu_2}=\CA^{\alpha,\alpha}_{\mu_1\mu_2}$.

For $\beta\in\BN[I]$ we define $\CR^{\alpha+\beta,\alpha+\beta}_{\mu_1+\beta,
\mu_2+\beta}:=\sigma_*\CR^{\alpha,\alpha}_{\mu_1,\mu_2}$ where $\sigma$
stands for the closed embedding $\CA^{\alpha,\alpha}_{\mu_1\mu_2}\hra
\CA^{\alpha+\beta,\alpha+\beta}_{\mu_1+\beta,\mu_2+\beta}$.

Thus the isomorphisms ~\ref{pjat' I-sheaf}(d) are simultaneously defined.

(c) We will construct the factorization isomorphisms ~\ref{pjat' I-sheaf}(c)
for the case $\alpha_1=\alpha_2=\alpha;\beta_1=\beta_2=\beta$ (cf. Remark
~\ref{pjat' part of data}).

\subsubsection{}
For $d>0$ we introduce the following analytic open subsets
$Q_{<d,<d},Q_{>d,>d},Q_{<d,>d},Q_{>d,<d}$:

$Q_{<d,<d}:=\{(p,q,t)\bigl|\ |p|<d>|q|\};\
Q_{>d,>d}:=\{(p,q,t)\bigl|\ |p|>d<|q|\};$

$Q_{<d,>d}:=\{(p,q,t)\bigl|\ |p|<d<|q|\};\
Q_{>d,<d}:=\{(p,q,t)\bigl|\ |p|>d>|q|\}$.

For $\alpha,\beta,\gamma\in\BN[I]$ we introduce the following open subset
$Q^{\alpha,\beta,\gamma}(d)\subset Q^{\alpha+\beta+\gamma}$: it is formed
by the configurations such that exactly $\alpha$ points of configuration
lie in $Q_{<d,<d}$; exactly $\beta$ points lie in $Q_{<d,>d}$, and exactly
$\gamma$ points lie in $Q_{<d,>d}$.

We have evident decomposition
$$Q^{\alpha,\beta,\gamma}(d)=Q^{0,\beta,0}(d)\times Q^{0,0,\gamma}(d)\times
Q^{\alpha,0,0}(d).$$
Intersecting the above subsets with $'Q^{\alpha+\beta+\gamma},\
'Q^{\bullet\alpha+\beta+\gamma},\ 'Q^{\circ\alpha+\beta+\gamma},\
Q^{\bullet\alpha+\beta+\gamma}$, etc. we obtain the opens
$'Q^{\alpha,\beta,\gamma}(d)$, etc. with similar decompositions.

We denote by $\CJ^{\bullet\alpha,\beta,\gamma}_{\mu_1\mu_2}(d)$ the
restriction of $\CJ^{\bullet\alpha+\beta+\gamma}_{\mu_1\mu_2}$ to
$'Q^{\bullet\alpha,\beta,\gamma}(d)$.

\subsubsection{}
\label{pjat' no name}
We have the canonical isomorphisms

\begin{equation}
\label{pjat' degener factor iso}
\CJ^{\bullet\alpha,\beta,\gamma}_{\mu_1\mu_2}(d)\iso
\CID^{(0,\beta)}_{\mu_1}(d)\boxtimes\CID^{(0,\gamma)}_{\mu_2}\boxtimes
\CJ^{\bullet\alpha,0,0}_{\mu_1-\beta',\mu_2-\gamma'}(d)
\end{equation}

satisfying the standard associativity constraints. Here
$\CID^{(0,\beta)}_{\mu_1}(d),\CID^{(0,\gamma)}_{\mu_2}$ have the following
meaning.

The set $Q_{<d,>d}$ (resp. $Q_{>d,<d}$) projects to the vertical component
$Q_{0v}$ (resp. horizontal component $Q_{0h}$) of the special fiber $Q_0:\
(p,q,t)\mapsto(0,q,0)$ (resp. $(p,q,t)\mapsto(p,0,0)$).

This induces the projections $$v^\beta:\ Q^{0,\beta,0}(d)\lra Q^\beta_{0v};$$
$$h^\gamma:\ Q^{0,0,\gamma}(d)\lra Q^\gamma_{0h}.$$

Recall that $Q_{0v}$ (resp. $Q_{0h}$) has the marked point $(0,\infty)$
(resp. $(\infty,0)$) lying on it, and $v^\beta(Q^{0,\beta,0}(d))$
(resp. $h^\gamma(Q^{0,0,\gamma}(d))$) is a standard open subset of
$Q^\beta_{0v}$ (resp. $Q^\gamma_{0h}$).

Restricting our maps to $'Q^\bullet$ we get
$$v^\beta:\ 'Q^{\bullet0,\beta,0}(d)\lra Q^{\bullet\beta}_{0v};$$
$$h^\gamma:\ 'Q^{\bullet0,0,\gamma}(d)\lra Q^{\bullet\gamma}_{0h}.$$
Note that $v^\beta(\ 'Q^{\bullet0,\beta,0}(d))=Q^{\bullet0,\beta}_{0v}(d)$
--- the standard open formed by the configurations with all $\beta$ points
running in the annular neighbourhood of the marked point on $Q_{0v}$.

Similarly, $h^\gamma(\ 'Q^{\bullet0,0,\gamma}(d))=Q^{\bullet0,\gamma}_{0h}(d)$
--- the standard open formed by the configurations with all $\gamma$ points
running in the annular neighbourhood of the marked point on $Q_{0h}$.

The sheaf $\CID^{0,\beta}_{\mu_1}(d)$ (resp. $\CID^{0,\gamma}_{\mu_2}(d)$)
is defined in III.3.5. We keep the same notation for its inverse image
to $'Q^{\bullet0,\beta,0}(d)$ (resp. $'Q^{\bullet0,0,\gamma}(d)$) with
respect to $v^\beta$ (resp. $h^\gamma$).

\subsubsection{}
One final remark is in order. Note that $Q^\bullet_{0v}=\BA^1_v;\
Q^\bullet_{0h}=\BA^1_h$ (notations of ~\ref{pjat' vh}). Hence $Q^\bullet\beta_{0v}=
\CA^\beta_v;\ Q^\bullet\gamma_{0h}=\CA^\gamma_h$.

Under this identification the open $Q^{\bullet0,\beta}_{0v}(d)$ corresponds
to the standard open $\CA^{\beta,0}_v(d)$, and the open
$Q^{\bullet0,\gamma}_{0h}(d)$ corresponds to the standard open
$\CA^{\gamma,0}_h(d)$ (notations of III.2.2).

{\em Warning.}
{\em Note the opposite roles played by the marked points on $Q_{0v}$ and $\BA^1_v$:
the open $Q^{\bullet0,\beta}_{0v}(d)=\CA^{\beta,0}_v(d)$ is formed by
configurations where all points
run {\em near} the marked point of $Q_{0v}$,  or equivalently,
where all points run {\em far from} the marked point $0\in\BA^1$.}

The sheaf $\CID^{0,\beta}_{\mu_1}(d)$ on $Q^{\bullet0,\beta}_{0v}(d)$
corresponds to the sheaf $\CID^{\beta,0}_{(l-1)2\rho+\beta-\mu_1}(d)$
on $\CA^{\beta,0}_v(d)$ (notations of III.3.5).

The sheaf $\CID^{0,\gamma}_{\mu_2}(d)$ on $Q^{\bullet0,\gamma}_{0h}(d)$
corresponds to the sheaf $\CID^{\gamma,0}_{(l-1)2\rho+\gamma-\mu_2}(d)$
on $\CA^{\gamma,0}_h(d)$.

\subsubsection{}
Putting all the above together, the desired factorization isomorphisms
~\ref{pjat' I-sheaf}(c) for the sheaf $\CR$ are induced by the factorization
isomorphisms ~(\ref{pjat' degener factor iso}) for the sheaf $\CJ^\bullet$.

\section{Convolution}

\subsection{}
\label{pjat' proper base change}
Let us throw away the line $p=\infty,q=0$ from $Q$. What remains is the
degenerating family $Q_a$ of affine lines with the marked point $(0,\infty)$.
The notations $'Q_a,Q^\bullet_a$, etc. speak for themselves.

Given a finite factorizable sheaf $\CX\in\FS$ we can glue it into the
marked point $(0,\infty)$ and obtain the factorizable sheaf $\tCX$ over the
family $'Q_a$ of affine lines. For each $t\in\BC^*$ the restriction of
$\tCX$ to the fiber $(\CQ_a)_t=\CA$ is the FFC $\CX$ we started with.

More precisely, let $\lambda=\lambda(\CX),\ \alpha\in\BN[I]$. We have the
sheaf $\tCX^\alpha_\lambda$ over $'Q^\alpha_a$ such that for any $t\in\BC^*$
its restriction to each fiber $(\ 'Q^\alpha_a)_t=\CA^\alpha$ coincides with
$\CX^\alpha_\lambda$.

We have the projection
$$\pr:\ Q_a\lra\BA^1_h\times\BC,\ \pr(p,q,t)=(p,t).$$
The restriction of $\pr$ to $'Q_a$ is one-to-one onto $\BA^1_h\times\BC^*$,
so $\pr_*\tCX$ is the constant family of FFC on $\CA_h\times\BC^*$.

Recall that $\pi$ denotes the projection to the $t$-coordinate.
The square


               \begin{equation*}
                  \begin{CD}
                      Q_a  @>\pi>>  \BC  \\
                      @V{\operatorname{pr}}VV    @A{\id}AA    \\
                      \BA^1\times \BC @>\pi>> \BC           \\
                   \end{CD}
               \end{equation*}

commutes, and the map $\pr$ is proper.

By the proper base change for nearby cycles we have
$$\pr^\alpha_*\Psi_{\pi^\alpha}\tCX^\alpha_\lambda=
\Psi_{\pi^\alpha}(\pr^\alpha_*\tCX^\alpha_\lambda)=\CX^\alpha_\lambda$$
Here $\pi^\alpha,\pr^\alpha$ form the evident commutative diagram

               \begin{equation*}
                  \begin{CD}
                      \CQ_a^\alpha  @>\pi^\alpha>>  \BC  \\
                      @V{\operatorname{pr}^\alpha}VV @A{\id}AA       \\
                      \CA^\alpha_h\times \BC @>\pi^\alpha>> \BC          \\
                   \end{CD}
               \end{equation*}

\subsection{}
Let us compute $\pr^\alpha_*\Psi_{\pi^\alpha}\tCX^\alpha_\lambda$ in another
way. First we describe $\Psi_{\pi^\alpha}\tCX^\alpha_\lambda$. It is a sheaf
on $(\CQ^\alpha_a)_0$ which will be denoted by $Q^\alpha_{0a}$ for short.

Note that $Q_{0a}=\BP^1_v\cup\BA^1_h$ is the union of two irreducible
components intersecting at the point $(0,0)$. Hence $Q^\alpha_{0a}$ is the
union of irreducible components numbered by the decompositions
$\alpha=\alpha_1+\alpha_2,\ \alpha_1,\alpha_2\in\BN[I]$. Namely,
$$Q^\alpha_{0a}=\bigcup_{\alpha_1+\alpha_2=\alpha}\CP_v^{\alpha_1}
\times\CA^{\alpha_2}_h.$$
Here $\CP^{\alpha_1}_v$ stands for the space of configurations on $\BP^1_v$.

Each irreducible component embeds into $\CP^\alpha_v\times\CA^\alpha_h$ as
in ~\ref{pjat' vh}, and one checks readily that these embeddings agree on
intersections, whence they combine to the closed embedding
$$\varsigma^\alpha:\ Q^\alpha_{0a}\hra\CP^\alpha_v\times\CA^\alpha_h.$$
We will describe $\varsigma^\alpha_*\Psi_{\pi^\alpha}\tCX^\alpha_\lambda$.

\subsubsection{}
Note that $\BP^1_v$ is equipped with two marked points with tangent vectors.
The marked points are $\infty,0$. The tangent vector at $0$ is $\partial_q$,
and the tangent vector at $\infty$ was defined in ~\ref{pjat' family}.

Recall that for $\mu\in {X_\ell}$ such that the triple $(\lambda,\mu,\alpha)$
is admissible, the sheaf $\varsigma^\alpha_*\Psi_{\pi^{\bullet\alpha}}
\CJ^{\bullet\alpha}_{\lambda\mu}$ on $\CA^\alpha_v\times\CA^\alpha_h$ was
introduced in ~\ref{pjat' CR}(b).

According to the Theorem IV.3.6 we can glue the sheaves
$\varsigma^\alpha_*\Psi_{\pi^{\bullet\alpha}}\CJ^{\bullet\alpha}_{\lambda\mu}$
and $\CX^\alpha_\lambda$ to obtain the sheaf $\bar{\CX}^\alpha_{\lambda\mu}$
on the space $\CP^\alpha_v\times\CA^\alpha_h$.

Here $\CX^\alpha_\lambda$ is viewed as a sheaf on $Q^{0,\alpha}_{0v}(d)
\times\CA^\alpha_h$ (notations of ~\ref{pjat' family} and ~\ref{pjat' no name}) constant
along $\CA^\alpha_h$. So it is glued into $\infty\times\CA^\alpha_h$, while
$\varsigma^\alpha_*\Psi_{\pi^{\bullet\alpha}}\CJ^{\bullet\alpha}_{\lambda\mu}$
is glued into $0\times\CA^\alpha_h$.

\subsubsection{Claim}
\label{pjat' tilda = bar}
$$\varsigma^\alpha_*\Psi_{\pi^{\bullet\alpha}}\tCX^\alpha_\lambda=
\bar{\CX}^\alpha_{\lambda\mu}$$

{\em Proof.} Evident. $\Box$

\subsection{}
We proceed with the computation of $\pr^\alpha_*\Psi_{\pi^\alpha}
\tCX^\alpha_\lambda$.

Recall the sheaf $\CR$ introduced in ~\ref{pjat' CR}.

\subsubsection{}
We define the {\em convolution} $\CX\star\CR$ on $\CP_v\times\CA_h$
as the following collection of data:

For $\mu\in {X_\ell}$ such that $\lambda+\mu-(l-1)2\rho=\alpha\in\BN[I]$,
the sheaf
$$(\CX\star\CR)^\alpha_{\lambda\mu\lambda}:=\bar{\CX}^\alpha_{\lambda\mu}$$
on the space $\CP^\alpha_{\lambda\mu}\times\CA^\alpha_\lambda$.

The notation $\CP^\alpha_{\lambda\mu}$ suggests that the monodromy around
$\infty$ is $\lambda$, and the monodromy around $0$ is $\mu$.

We leave the formulation of factorization isomorphisms to the interested
reader.

\subsubsection{}
Denote by $\pr^\alpha$ the projection $\CP^\alpha\times\CA^\alpha\lra
\CA^\alpha$. Then, evidently,
$$\pr^\alpha_*(\CX\star\CR)^\alpha_{\lambda\mu\lambda}=\pr^\alpha_*\bar
{\CX}^\alpha_{\lambda\mu}=\pr^\alpha_*\varsigma^\alpha_*\Psi_{\pi^\alpha}
\tCX^\alpha_\lambda=\CX^\alpha_\lambda$$
Here the second equality holds by the Claim ~\ref{pjat' tilda = bar},
and the last equality was explained in ~\ref{pjat' proper base change}.

\subsection{}
Thus, for each $\CX\in\FS$, we have the natural isomorphism
$\pr_*(\CX\star\CR)=\CX$. That is, for $\lambda=\lambda(\CX),
\alpha\in\BN[I],\mu=\alpha+(l-1)2\rho-\lambda$, we have
$\pr_*(\CX\star\CR)^\alpha_{\lambda\mu\lambda}=\CX^\alpha_\lambda$.

\subsubsection{Lemma}
There is a natural isomorphism of functors in $\CX$:
$$\Phi\pr_*(\CX\star\CR)\iso F_{\Phi(\CR);\frac{\infty}{2}+\bullet}
(\Phi(\CX))$$
(notations of ~\ref{pjat' kernel functor}).

{\em Proof.} This is just the Theorem IV.8.4. $\Box$

\subsection{Theorem}
\label{pjat' CR vs bR}
There is a natural isomorphism in $\Ind\CC\otimes\Ind\CC$:
$$\Phi(\CR)\iso \bR.$$

{\em Proof.} This is just the above Lemma combined with the Converse
Theorem ~\ref{pjat' Converse Theorem}. $\Box$

\newpage
\begin{center}
{\bf Chapter 5. Modular functor}
\end{center}
\vspace{.8cm}

\bigskip

\section{Gluing over $C$}

This section is quite parallel to the Chapter 1 of IV.

We will change our notations slightly to make them closer to those of
{\em loc. cit.}

\subsection{}
Let again $C\lra S$ be a smooth proper relative curve of genus $g$.

Let us assume that $g>1$ for a moment. It is well known that the isomorphism
classes of complex structures on a surface of genus $g$ correspond bijectively
to the conformal equivalence classes of Riemann metrics on this surface.

Each conformal equivalence class contains a unique metric of constant curvature
$-1$. Thus we may assume that each fiber $C_s$ is equipped with such a metric.

Given $0<\varepsilon\leq1$ and a finite set $K$, let
$\tC^K_\varepsilon$ denote the space of $K$-tuples $(u_k)_{k\in K}$ of
analytic morphisms over $S:\ S\times D_\varepsilon\lra C$
inducing isometry on each fiber, such that the images
$u_k(S\times D_\varepsilon)$
do not intersect. Here we equip the disk $D_\varepsilon$ with the Poincar\'e
metric of constant curvature $-1$.

Given a $K$-tuple $\valpha=(\alpha_k)\in\BN[{X_\ell}]^K$ such that
$\alpha=\sum_k\ \alpha_k$, define a space
$$
C^{\valpha}_\varepsilon:=\tC^K_\varepsilon\times\prod_{k\in K}\ \TB
D^{\alpha_k}_\varepsilon
$$
and an open subspace
$$
C^{o\valpha}_\varepsilon:=\tC^K_\varepsilon\times\prod_{k\in K}\
TD^{o\alpha_k}_\varepsilon.
$$
We have an evident "substitution" map
$$
q_{\valpha}:C^{\valpha}_\varepsilon\lra\TB C^{\alpha}
$$
which restricts to $q_{\valpha}:C^{o\valpha}_\varepsilon\lra T C^{o\alpha}$.

\subsection{}
\label{pjat' bal fun} Let us define an element $\rho_\ell\in {X_\ell}$ by the condition
$\langle i,\rho_\ell\rangle=\ell_i-1$ for all $i\in I$ (see IV.9.1.2).
>From now on we choose a balance
function $n$ in the form
$$
n(\mu)=\frac{1}{2}\mu\cdot\mu - \mu\cdot\rho_\ell.
$$
In other words, we set
$$
\nu_0=-\rho_\ell.
$$

We pick a corresponding Heisenberg system $\CH$.

\subsubsection{Remark} The monodromy of $\CH$ around the zero section
of determinant line bundle is equal to
$(-1)^{\rk {X_\ell}}\zeta^{12\rho_\ell\cdot\rho_\ell}$.
According to the strange formula of Freudenthal--deVries
 (see e.g. ~\cite{fz})
we have $12\rho\cdot\rho=dh\dim\fg$. Here $\fg$ is the
Lie algebra with the root datum $X,Y,\ldots$, and $\check h$ is the dual
Coxeter number of $\fg$, while $d:=\max_{i\in I}d_i$.

Let $\zeta=\exp(\frac{\pi\sqrt{-1}}{d\kappa})$
for some positive integer $\kappa$. Then the above monodromy equals
$\exp(\pi\sqrt{-1}(\rk\fg+\frac{h\dim\fg}{\kappa}))$ which coincides
with the multiplicative central charge of the conformal field theory
associated with the affine Lie algebra $\hat{\fg}$ at level $\kappa$
(see e.g.~\cite{bfm}).

We will offer an explanation of this coincidence in Chapter 6.

\subsubsection{}
Given $\alpha=\sum a_ii\in \BN[I]$ and $\vmu=(\mu_k)\in {X_\ell}^K$, we define an
element
$$
\alpha_{\vmu}=\sum a_i\cdot (-i')+\sum_k\mu_k\in\BN[{X_\ell}]
$$
where the sum in the right hand side is a formal one. We say that a pair
$(\vmu,\alpha)$ is {\em $g$-admissible} if
$$
\sum_k\ \mu_k-\alpha\equiv (2-2g)\rho_\ell\bmod Y_\ell.
$$
Note that given $\vmu$, there exists $\alpha\in\BN[I]$ such that $(\vmu,\alpha)$
is admissible if and only if $\sum_k\ \mu_k\in Y$; if this holds true,
such elements $\alpha$ form an obvious countable set.

We will denote by
$$
e:\BN[I]\lra \BN[{X_\ell}]
$$
a unique homomorphism sending $i\in I$ to $-i'\in {X_\ell}$.

\subsection{}
\label{pjat' triples}
Let us consider the space $\TB C^K\times\TB C^{e(\alpha)}$;
its points are quadruples $((z_k),(\tau_k),(x_j),(\omega_j))$ where
$(z_k)\in C^K,\ \tau_k$ --- a non-zero (relative) tangent vector to $C$ at
$z_k$, $(x_j)\in C^{e(\alpha)}$, $\omega_j$ --- a non-zero tangent
vector at $x_j$. To a point $z_k$ is assigned a weight $\mu_k$,
and to $x_j$ --- a weight $-\pi(j)'$. Here $\pi:J\lra I$ is an unfolding
of $\alpha$ (implicit in the notation $(x_j)=(x_j)_{j\in J}$).

We will be interested in some open subspaces:
$$
\TB C^{\alpha}_{\vmu}:=T C^{oK}\times\TB C^{e(\alpha)}\subset
\TB C^{\alpha_{\vmu}}
$$
and
$$
T C^{o\alpha}_{\vmu}\subset
\TB C^{\alpha}_{\vmu}
$$
whose points are quadruples $((z_k),(\tau_k),(x_j),(\omega_j))\in
\TB C^{\alpha}_{\vmu}$  with all $z_k\neq x_j$. We have an obvious
symmetrization projection
$$
p^{\alpha}_{\vmu}: T C^{o\alpha}_{\vmu}\lra T C^{o\alpha_{\vmu}}.
$$
Define a space
$$
 C^{\alpha}_{\vmu}=T C^{oK}\times C^{e(\alpha)};
$$
its points are triples $((z_k),(\tau_k),(x_j))$ where
$(z_k)$, $(\tau_k)$ and $(x_j)$ are as above; and to $z_k$ and $x_j$
the weights as above are assigned. We have the canonical projection
$$
\TB C^{\alpha}_{\vmu}\lra  C^{\alpha}_{\vmu}.
$$
We define the open subspaces
$$
 C^{o\alpha}_{\vmu}\subset
 C^{\bullet\alpha}_{\vmu}\subset
 C^{\alpha}_{\vmu}.
$$
Here the $\bullet$-subspace (resp., $o$-subspace)
consists of all $((z_k),(\tau_k),(x_j))$ with
$z_k\neq x_j$ for all $k,j$ (resp., with all $z_k$ and $x_j$ distinct).

We define the {\em principal stratification} $\CS$ of $ C^{\alpha}_{\vmu}$
as the stratification generated by subspaces $z_k=x_j$ and $x_j=x_{j'}$
with $\pi(j)\neq\pi(j')$. Thus, $ C^{o\alpha}_{\vmu}$ is the open
stratum of $\CS$. As usually, we will denote by the same letter the
induced stratifications on subspaces.

The above projection restricts to
$$
T C^{o\alpha}_{\vmu}\lra  C^{o\alpha}_{\vmu}.
$$

\subsection{Factorization structure}
\subsubsection{} Suppose we are given $\valpha\in\BN[I]^K,\ \beta\in\BN[I]$;
set $\alpha:=\sum_k\ \alpha_k$. Define a space
$$
 C^{\valpha,\beta}_{\vmu,\varepsilon}\subset
\tC^K_\varepsilon\times\prod_k\ D^{\alpha_k}_\varepsilon\times C^{e(\beta)}
$$
consisting of all collections $((u_k),((x^{(k)}_j)_k),(y_j))$
where $(u_k)\in\tC^K_\varepsilon ,\ (x^{(k)}_j)_k\in D^{\alpha_k},\
(y_j)\in C^{e(\beta)}$, such that
$$
y_j\in C-\bigcup_{k\in K}\ \overline{u_k(S\times D_\varepsilon)}
$$
for all $j$ (the bar means closure).

We have canonical maps
$$
q_{\valpha,\beta}: C^{\valpha,\beta}_{\vmu,\varepsilon}\lra
C^{\alpha+\beta}_{\vmu},
$$
assigning to $((u_k),((x^{(k)}_j)_k),(y_j))$ a configuration
$(u_k(0)),(\overset{\bullet}{u}_k(\tau)),(u_k(x_j^{(k)})),(y_j))$,
where $\tau$ is the unit tangent vector to $D_\varepsilon$ at $0$,
and
$$
p_{\valpha,\beta}: C^{\valpha,\beta}_{\vmu,\varepsilon}\lra
\prod_k\ D^{\alpha_k}_\varepsilon \times C^{\bullet\beta}_{\vmu-\valpha}
$$
sending $((u_k),((x^{(k)}_j)_k),(y_j))$ to $((u_k(0)),
(\overset{\bullet}{u}_k(\tau)),(y_j))$.

\subsection{}
{\em From now on
we change the base to $\dot{\delta}\lra S$ in all the above spaces
and morphisms. We preserve the old notations.}

Given a $g$-admissible pair $(\vmu,\alpha)$ we can find another
$K$-tuple of weights $\vmu'\in {X_\ell}^K$ such that $\vmu'$ is congruent to
$\vmu$ modulo $Y_\ell$ (see IV.9.1.2),
and $\alpha_{\vmu'}$ is $g$-admissible in the sense
of ~\ref{pjat' tilde}.

Let us consider a local system
$p^{\alpha*}_{\vmu'}\CH^{\alpha_{\vmu'}}$
over $T C^{o\alpha}_{\vmu'}$.

Note that all the spaces $T C^{o\alpha}_{\vmu'}$ (for different choices
of $\vmu'$) are identified with $T C^{o\alpha}_{\vmu}$, and the local
system $p^{\alpha*}_{\vmu'}\CH^{\alpha_{\vmu'}}$ does not depend on the
choice of $\vmu'$ by the virtue of the Lemma ~\ref{pjat' Independence}.
Thus we will identify all these local systems and
call the result $p^{\alpha*}_{\vmu}\CH^{\alpha_{\vmu}}$.

By our choice of the balance function $n$, its monodromies with
respect to the rotating of tangent vectors $\omega_j$ at points $x_j$
corresponding to negative simple roots, are trivial.

Therefore
it descends to a unique local system over
$C^{o\alpha}_{\vmu}$, to be denoted by $\CH^{\alpha}_{\mu}$.

We define a perverse sheaf
$$
\CH^{\bullet\alpha}_{\vmu}:= j_{!*}\CH^{\alpha}_{\vmu}
[\dim  C^{\alpha}_{\vmu}]\in
\CM( C^{\bullet\alpha}_{\vmu};\CS).
$$

\subsection{Factorizable sheaves over $C$}
\label{pjat' gluing sheaves}
Suppose we are given
a $K$-tuple of FFS's $\{\CX_k\},\ \CX_k\in\FS_{c_k},\ k\in K,\ c_k\in {X_\ell}/Y$,
where $\sum_k\ c_k=0$.
Let us pick
$\vmu=(\vmu_k)\geq (\lambda(\CX_k))$.

Let us call a {\em factorizable sheaf over $C$ obtained by
gluing the sheaves $\CX_k$} the following collection of data which
we will denote by $\operatorname{g}(\{\CX_k\})$.

(i) For each $\alpha\in\BN[I]$ such that $(\vmu,\alpha)$ is $g$-admissible,
a sheaf $\CX^{\alpha}_{\vmu}\in\CM( C^{\alpha}_{\vmu};\CS)$.

(ii) For each $0<\varepsilon\leq1,
\valpha=(\alpha_k)\in\BN[I]^K,\ \beta\in \BN[I]$ such that
$(\vmu,\alpha+\beta)$ is $g$-admissible (where $\alpha=\sum\ \alpha_k$),
a {\em factorization isomorphism}
$$
\phi_{\valpha,\beta}: q^*_{\valpha,\beta}\CX^{\alpha+\beta}_{\vmu}\iso
p^*_{\valpha,\beta}((\Boxtimes_{k\in K}\ \CX^{\alpha_k}_{\mu_k})
\boxtimes\CH^{\bullet\beta}_{\vmu-\valpha}).
$$

These isomorphisms should satisfy the standard associativity property.

\subsection{Theorem}
\label{pjat' glu thm} {\em There exists a unique up to a canonical
isomorphism factorizable sheaf over $C$ obtained by gluing
the sheaves $\{\CX_k\}$.}

{\bf Proof} is similar to III.10.3. $\Box$

\subsection{}
\label{pjat' gluing I-sheaves}
We will need the following version of the above Theorem.

Suppose we are given a $K$-tuple $\{\CX_k\}$
of objects of $\FS$, and a $J$-tuple $\{\CY_j\}$
of I-sheaves (see ~\ref{pjat' I-sheaf}).

Let us pick $\vlambda=(\vlambda_k)\geq(\lambda(\CX_k))$, and
$(\vmu,\vnu)=(\vmu_j,\vnu_j)$ such that $(\vmu_j+\vnu_j)
\geq(\chi(\CY_j))$.

We will denote the concatenation of $\vlambda,\ \vmu$, and $\vnu$
by $\vlambda\vmu\vnu$.

Let us call a {\em factorizable sheaf over $C$ obtained by
gluing the sheaves $\CX_k$ and I-sheaves $\CY_j$}
the following collection of data which
we will denote by $\operatorname{g}(\{\CX_k\},\{\CY_j\})$.

(i) For each $\alpha\in\BN[I]$ such that $(\vlambda\vmu\vnu,
\alpha)$ is $g$-admissible,
a sheaf $\CX^{\alpha}_{\vlambda\vmu\vnu}=
\og(\{\CX_k\},\{\CY_j\})^{\alpha}_{\vlambda\vmu\vnu}\in
\CM( C^{\alpha}_{\vlambda\vmu\vnu};\CS)$.

(ii) For each $0<\varepsilon\leq1,
\valpha=(\alpha_k)\in\BN[I]^K,\
(\vbeta,\vgamma)=(\beta_j,\gamma_j)\in\BN[I]^{J\sqcup J},\
\xi\in \BN[I]$ such that
$\beta_j,\gamma_j\leq\vmu_j+\vnu_j-\chi(\CY_j)$ for any $j\in J$, and
$(\vlambda\vmu\vnu,\alpha+\beta+\gamma+\xi)$
is $g$-admissible (where $\alpha=\sum\alpha_k,\ \beta=\sum\beta_j,\ \gamma=
\sum\gamma_j$),
a {\em factorization isomorphism}
$$
\phi_{\valpha\vbeta\vgamma,\xi}:\
q^*_{\valpha\vbeta\vgamma,\xi}
\CX^{\alpha+\beta+\gamma+\xi}_{\vlambda\vmu\vnu}\iso
p^*_{\valpha\vbeta\vgamma,\xi}
((\Boxtimes_{k\in K}\ \CX^{\alpha_k}_{\lambda_k})
\boxtimes
(\Boxtimes_{j\in J}\ \CY^{\beta_j,\gamma_j}_{\mu_j,\nu_j})\boxtimes
\CH^{\bullet\xi}_{\vlambda\vmu\vnu-\valpha\vbeta\vgamma}).
$$

These isomorphisms should satisfy the standard associativity property.

{\bf Theorem.}
{\em There exists a unique up to a canonical
isomorphism factorizable sheaf over $C$ obtained by gluing
the sheaves $\{\CX_k\}$ and I-sheaves $\{\CY_j\}$.}

{\bf Proof} is similar to III.10.3. $\Box$

\subsection{}
\label{pjat' gluing CR}
We will apply the above Theorem exclusively to the case $\CY_j=\CR$
(see ~\ref{pjat' CR}).

\subsubsection{}
\label{pjat' nodes}
Recall (see ~\ref{pjat' triples}) that $C^\alpha_{\vlambda\vmu\vnu}=
TC^{0K\sqcup J\sqcup J}\times C^{e(\alpha)}$ is
the space of 7-tuples $$((z_k),(\tau_k);(x_j),
(\omega_j);(y_j),(\eta_j);(t_m)).$$
Here $(t_m)\in C^{e(\alpha)};\ (z_k)\in C^K;\ (x_j),(y_j)\in C^J$.
All the points $z_k,x_j,y_j$ are distinct, and $\tau_k$ is a nonzero tangent
vector at $z_k;\ \omega_j$ is a nonzero tangent vector at $x_j;\ \eta_j$
is a nonzero tangent vector at $y_j$.

We define
$$\ul{C}^\alpha_{\vlambda\vmu\vnu}:=
TC^{0K\sqcup J\sqcup J}\times (C/[(x_j)=(y_j)])^{e(\alpha)}$$

Note that the natural projection $\fp:\
C^\alpha_{\vlambda\vmu\vnu}\lra\ul{C}^\alpha_{\vlambda\vmu\vnu}$
is the normalization map.

\subsubsection{}
\label{pjat' gluing into nodes}
Suppose we are in the situation of ~\ref{pjat' gluing I-sheaves},
and $\{\CY_j\}$ is just a set of $J$ copies of the I-sheaf $\CR$ (see
~\ref{pjat' CR}). Thus, $\{\CY_j\}=\{\CR_j\};\ \CR_j=\CR$ for any $j\in J$.

It is immediate from the definition of $\CR$ that the sheaf
$\og(\{\CX_k\},\{\CR_j\})^{\alpha}_{\vlambda\vmu\vnu}$ on
$C^{\alpha}_{\vlambda\vmu\vnu}$ descends to a sheaf on
$\ul{C}^{\alpha}_{\vlambda\vmu\vnu}$ to be denoted by
$\ul{\og}(\{\CX_k\},\{\CR_j\})^{\alpha}_{\vlambda\vmu\vnu}$.

\subsubsection{Remark}
\label{pjat' disconnected}
All the above constructions generalize immediately
to the case where $C$ is a union of smooth connected components $C_n$ of
genera $g_n,\ \sum g_n=g$.

Then $C^\alpha$ is a union of connected components numbered by the
partitions $\alpha=\sum\alpha_n$. Each connected component is the
product $\prod C_n^{\alpha_n}$. The Heisenberg local system is just
the product of $\CH_{g_n}$. A factorizable sheaf $\og(\{\CX_k\})$
on each connected component
also decomposes into external product of the corresponding sheaves
on factors.
Note that each connected component gives rise to its own admissibility
condition.

On the other hand, a sheaf $\ul{\og}(\{\CX_k\},\{\CR_j\})$ does
not necessarily decompose into external product. In effect, if for
some $j$ the sections $x_j$ and $y_j$ lie on the different connected
components of $C$, then some connected components of
$C^\alpha_{\vlambda\vmu\vnu}$ will be glued together in
$\ul{C}^\alpha_{\vlambda\vmu\vnu}$.

\section{Degeneration of factorizable sheaves}

\subsection{Definition} An $m$-tuple of weights $\vmu\in {X_\ell}^m$ is
called {\em  $g$-positive} if
$$\sum^m_{j=1}\mu_j+2(g-1)\rho_\ell\in\BN[I]\subset {X_\ell}$$
(notations of IV.9.1.2).
If this is so, we will denote
$$\alpha_g(\vmu):=\sum^m_{j=1}\mu_j+2(g-1)\rho_\ell $$

\subsection{}
\label{pjat' preparation}
We will use freely all the notations
 of section ~\ref{pjat' degeneration of curves}. In particular, $C$ is the
universal curve of genus $g$.

Given a $K$-tuple of FFS's $\{\CX_k\}$ and a $g$-positive $K$-tuple
of weights $\vmu=(\vmu_k)\geq(\lambda(\CX_k))$ we define $\alpha:=
\alpha_g(\vmu)$.

We consider a sheaf $\CX^\alpha_\vmu={\og}^\alpha_\vmu(\{\CX_k\})$
on $C^\alpha_\vmu$ obtained by gluing the sheaves $\CX_k$ (see
~\ref{pjat' gluing sheaves}). We will study its specialization along the
boundary component $\fD$.

Let $\oC$ denote the following object:

in case ~\ref{pjat' cases}(a) $\oC$ is the universal curve over $\CM_{g_1}\times
\CM_{g_2}$. It consists of two connected components $\oC_1$ and $\oC_2$.

in case ~\ref{pjat' cases}(b) $\oC$ is the universal curve $\oC_0$ over $\CM_{g_0}$.

The maps ~(\ref{pjat' map a}) and ~(\ref{pjat' map b}) give rise to the following
morphisms also denoted by $\wp_\alpha$:

\begin{equation}
\label{pjat' morphism a}
(\oC_1)^{\alpha_1}_{\vmu_1,\mu_1}\times(\oC_2)^{\alpha_2}_{\vmu_2,\mu_2}\lra
\TB_{(\ol{C}|_{\fD})^\alpha_\vmu}
\ol{C}^\alpha_{\vmu};
\end{equation}

\begin{equation}
\label{pjat' morphism b}
(\oC_0)^\alpha_{\vmu,\mu_1,\mu_2}\lra
\TB_{(\ol{C}|_{\fD})^\alpha_\vmu}
\ol{C}^\alpha_{\vmu}
\end{equation}

(recall that by default all the bases are changed to the punctured
determinant line bundles). In ~(\ref{pjat' morphism a})
the sets $\vmu_1$ and $\vmu_2$ form an arbitrary partition of the set $\vmu$.
Finally, we can choose $\mu_1,\mu_2$ arbitrarily:
the spaces we consider do not depend on a choice of $\mu_r$.

Recall that $\oC^\alpha_{\vmu,\mu_1,\mu_2}$ is the set of tuples
$((z_k),(\tau_k);x,\omega;y,\eta;(t_m))$.
Here $(t_m)\in C^{e(\alpha)};\ (z_k)\in C^K;\ (x_j),(y_j)\in C^J$.
All the points $z_k,x,y$ are distinct, and $\tau_k$ is a nonzero tangent
vector at $z_k;\ \omega$ is a nonzero tangent vector at $x;\ \eta$
is a nonzero tangent vector at $y$ (see ~\ref{pjat' nodes}).

We defined the nodal curve $\ul{\oC}$ in {\em loc. cit.}.
One can see that the morphisms ~(\ref{pjat' morphism a}) for different connected
components (resp. the morphism ~(\ref{pjat' morphism b})) can be glued together into
(resp. factor through) the same named morphism $\wp_\alpha$:

\begin{equation}
\label{pjat' morphism}
\ul{\oC}^\alpha_{\vmu,\mu_1,\mu_2}\lra
\TB_{(\ol{C}|_{\fD})^\alpha_\vmu}
\ol{C}^\alpha_{\vmu}
\end{equation}

\subsection{Theorem}
\label{pjat' main}
There is a canonical isomorphism
$$\wp_\alpha^*\bSp\og^\alpha_\vmu(\{\CX_k\})=\ul{\og}^\alpha_{\vmu,\mu_1,\mu_2}
(\{\CX_k\},\CR)$$
(notations of ~\ref{pjat' gluing into nodes} and ~\ref{pjat' disconnected}).

We have to specify $\mu_1,\mu_2$ in the formulation of the Theorem.

\subsubsection{}
\label{pjat' mu_r}
In case ~(\ref{pjat' morphism a})
the choice of $\mu_1,\mu_2$ depends on the connected component (see
~\ref{pjat' disconnected}). Namely, we set $\mu_r=\alpha_r-\alpha_{g_r}(\vmu_r),\
r=1,2$.

In case ~(\ref{pjat' morphism b}) we sum up over all the choices of
$(\mu_1,\mu_2)\in ({X_\ell}\times {X_\ell})/Y_\ell$ (antidiagonal action) such that
$\mu_1+\mu_2=2\rho_\ell$.
Thus, if we identify all the spaces $\ul{\oC}^\alpha_{\vmu,\mu_1,\mu_2}$
with one abstract space $\ul{\oC}^\alpha_{\vmu,?,?}$, then
$\wp_\alpha^*\bSp\og^\alpha_\vmu(\{\CX_k\})$ is a direct sum of ${\od}_\ell$
summands. Or else we can view $\wp_\alpha^*\bSp\og^\alpha_\vmu(\{\CX_k\})$
as a collection of sheaves on different spaces
$\ul{\oC}^\alpha_{\vmu,\mu_1,\mu_2}$.

\subsubsection{Proof} By the definition of gluing (~\ref{pjat' gluing I-sheaves})
it suffices to construct the desired isomorphism ``near the points $x,y$''
and, separately, ``away from the points $x,y$''.

The situation ``near the points $x,y$'' is formally equivalent to (a direct
sum of a few copies of) the situation of ~\ref{pjat' CR}. Thus we obtain
a copy of the I-sheaf $\CR$ glued into the node $x=y$.

The situation ``away from the points $x,y$'' was the subject of the
Theorem ~\ref{pjat' modular Heisenberg}.

We leave it to the reader to work out the values of $\mu_1,\mu_2$. We just
note that these are exactly the values making the tuple $(\alpha,\vmu,\mu_1,
\mu_2)$ $g$-admissible. $\Box$

\section{Global sections over $C$}

\subsection{}
Given a $K$-tuple of weights $\vmu$ and $\alpha=\alpha_g(\vmu)$,
consider the projection $$\eta:\ C^\alpha_\vmu\lra TC^{\circ K}=\TO C^\vmu$$
(see ~\ref{pjat' triples}).

Given a $K$-tuple of FFS's $\{\CX_k\}$ with $(\vmu_k)\geq(\lambda(\CX_k))$
we consider the sheaf $\og(\{\CX_k\})^\alpha_\vmu$ on $C^\alpha_\vmu$
obtained by gluing the sheaves $\CX_k$.

In this section we will study the complex of sheaves with smooth cohomology
$R\eta_*\og(\{\CX_k\})^\alpha_\vmu$ on $\TO C^\vmu$.
Note that it does not depend on a choice of $(\vmu_k)\geq(\lambda(\CX_k))$
(if we identify the spaces $\TO C^\vmu$ for different choices of $\vmu$).

Namely, we want to compute its specialization along the boundary
component $\fD$. Recall the notations of ~\ref{pjat' degeneration of curves}.

In case ~\ref{pjat' cases}(a) for any partition $\vmu=\vmu_1\sqcup\vmu_2$
we have the map (cf. ~(\ref{pjat' map a})) $\wp_\vmu$:

\begin{equation}
\label{pjat' projection a}
\TO\oC_1^{\vmu_1,\mu_1}\times\TO\oC_2^{\vmu_2,\mu_2}\lra
\TB_{\TO\ol{C}^\vmu|_{\fD}}\TO\ol{C}^\vmu;
\end{equation}

in case ~\ref{pjat' cases}(b) we have the map (cf. ~(\ref{pjat' map b})) $\wp_\vmu$:

\begin{equation}
\label{pjat' projection b}
\TO\oC_0^{\vmu,\mu_1,\mu_2}\lra
\TB_{\TO\ol{C}^\vmu|_{\fD}}\TO\ol{C}^\vmu.
\end{equation}

On the other hand, consider the space $\ul{\oC}^\alpha_{\vmu,\mu_1,\mu_2}$
(see ~\ref{pjat' preparation}). In the case ~\ref{pjat' cases}(b)
it projects to $\TO\oC_0^{\vmu,\mu_1,\mu_2}$. We denote the projection
by $\ul\eta$.

In the case ~\ref{pjat' cases}(a)
$\ul{\oC}^\alpha_{\vmu,\mu_1,\mu_2}$ is disconnected. Its connected components
are numbered by the partitions $\vmu=\vmu_1\sqcup\vmu_2$. Each component
projects to the corresponding space
$\TO\oC_1^{\vmu_1,\mu_1}\times\TO\oC_2^{\vmu_2,\mu_2}$. Each one of those
projections will be denoted by $\ul\eta$, as well as their totality.

\subsection{Theorem}
\label{pjat' modular functor}
In the notations of ~\ref{pjat' main} and ~\ref{pjat' mu_r}
we have a canonical isomorphism
$$\wp_\vmu^*\bSp R\eta_*\og(\{\CX_k\})^\alpha_\vmu=
R\ul{\eta}_*\ul{\og}^\alpha_{\vmu,\mu_1,\mu_2}(\{\CX_k\},\CR)$$

{\em Proof.} This is just the Theorem ~\ref{pjat' main} plus the proper
base change for nearby cycles. $\Box$

One can apply Theorems ~\ref{pjat' main} and ~\ref{pjat' modular functor} iteratively
--- degenerating the curve $C$ more and more, and inserting more and more
copies of $\CR$ into new nodes.

In the terminology of ~\cite{bfm} 4.5 (see especially 4.5.6)
the Theorem ~\ref{pjat' modular functor} means that the category $\FS$ is
equipped with the structure of {\em fusion category}. The correspondence
associating to any $K$-tuple of FFS's $\{\CX_k\}$  the cohomology local systems
$R\eta_*\og(\{\CX_k\})^\alpha_\vmu$ on the space $TC^{\circ K}$, satisfying
the compatibility conditions ~\ref{pjat' modular functor}, is called
{\em fusion functor}, or {\em modular functor}.

\newpage
\begin{center}
{\bf Chapter 6. Integral representations of conformal blocks}
\end{center}
\vspace{.8cm}

\bigskip

\section{Conformal blocks in arbitrary genus}

\subsection{}
Recall the notations of IV.9. We fix a positive integer $\kappa$ and
consider the category $\tCO_\kappa$ of integrable $\hat{\fg}$-modules
of central charge $\kappa-h$. We set $\zeta=\exp(\frac{\pi\sqrt{-1}}
{d\kappa})$, so that $l=2d\kappa$, and $\ell=d\kappa$.

We have an equivalence of bbrt categories $\phi:\ \tCO_\kappa\iso\tCO_\zeta$
where $\tCO_\zeta$ is a semisimple subcategory of $\CC$ defined in
{\em loc. cit.} It is well known that $\tCO_\kappa$ has a structure of
{\em fusion category} (see e.g. ~\cite{bfm} or ~\cite{tuy}).
The category $\tCO_\zeta$ is also equipped with the structure of
fusion category via the equivalence $\phi$. From now on we will not
distinguish between the categories $\tCO_\kappa$ and $\tCO_\zeta$.

\subsection{}
Here is a more elementary way to think of fusion functors.
Recall that a fusion functor is a correspondence associating to any
$K$-tuple of objects $\{L_k\}$ of our category the local systems
(or, more generally, the complexes of sheaves with smooth cohomology)
$\langle\{L_k\}\rangle^K_C$ on the spaces $TC^{\circ K}$, satisfying a
certain compatibility conditions.

We can view those local systems as representations of Modular Teichm\"{u}ller
groups $T_{g,K}$ in the stalks at some base points.

We will choose the base points ``at infinity'' (see ~\cite{d4}).
Namely, we degenerate the
curve $C$ into the projective line $\BP^1$ with $g$ nodal points.
More precisely, we choose the real numbers $x_1<y_1<x_2<\ldots<x_g<y_g$
and degenerate $C$ into $\BP^1/(x_j=y_j,\ 1\leq j\leq g)$. The marked
points degenerate into the real numbers $z_1<\ldots<z_m$ ($m$ is the cardinality
of $K$) such that $z_m<x_1$. The tangent vectors at all the points are
real and directed to the right. We will denote this homotopy
base point at infinity by $\infty_{g,K}$.

\subsection{}
\label{pjat' usual action}
For $L_1,\ldots,L_m\in\tCO_\zeta$ the stalk
$(\langle\{L_k\}\rangle^K_C)_{\infty_{g,K}}$ is canonically isomorphic
to $\langle L_1\tilde{\otimes}\ldots\tilde{\otimes}L_m\tilde{\otimes}
\tilde{\ad}^{\tilde{\otimes}g}\rangle$ (see ~\cite{tuy}). Here
we use the notation $\langle ?\rangle$ introduced in IV.5.9 and IV.9
(maximal trivial summand), and the notation $\tilde{\ad}$ introduced
in ~\ref{pjat' ss ad} (semisimple adjoint representation).

The action of $T_{g,K}$ on
$\langle L_1\tilde{\otimes}\ldots\tilde{\otimes}L_m\tilde{\otimes}
\tilde{\ad}^{\tilde{\otimes}g}\rangle$ is induced by its action on
$L_1\tilde{\otimes}\ldots\tilde{\otimes}L_m\tilde{\otimes}
\tilde{\ad}^{\tilde{\otimes}g}$: generators act by braiding, balance,
and Fourier transforms (see ~\cite{lm})
of the adjoint representation.

\subsection{}
The Theorem ~\ref{pjat' modular functor} equips the category $\FS$ (and, via
the equivalence $\Phi$, the category $\CC$) with the fusion structure.
To unburden notations and to distinguish from the case of $\tCO_\zeta$
we will denote the corresponding local systems (or, rather, complexes
of sheaves with smooth cohomology) by $[L_1,\ldots,L_m]^K_C$.

{\bf Theorem.}
The stalk $([L_1,\ldots,L_m]^K_C)_{\infty_{g,K}}$ is canonically
isomorphic to $$\Tor^\CC_{\frac{\infty}{2}+\bullet}(\sk,L_1\otimes\ldots
\otimes L_m\otimes{\ad}^{\otimes g})$$ (notations of IV.4.7 and
~\ref{pjat' bR}).

{\em Proof.}
In effect, the stalk in question is the global sections
of some sheaf on the space $\oC^\alpha$ where
$\oC=\BP^1/(x_j=y_j,\ 1\leq j\leq g)$. The global sections will not
change if we merge the points $y_j$ into $x_j$ for any $j$ along the
real line. The specialization of our sheaf will live on the space
$(\BP^1)^\alpha$ (no nodal points anymore).

It will be the sheaf
obtained by gluing $\{\CL_k\},\CT_1,\ldots,\CT_g$ into the marked points
$z_1,\ldots,z_m,x_1,\ldots,x_g$. Here $\CL_k=\Phi^{-1}L_k$, and the sheaves
$\CT_1,\ldots,\CT_g$ are just $g$ copies of the same sheaf $\CT=
\Phi^{-1}\ad$. This follows immediately from the Theorem ~\ref{pjat' main},
the Definition of $\ad$ in ~\ref{pjat' ad}, and the Theorem ~\ref{pjat' CR vs bR}.

It remains to apply the Theorem IV.8.4. $\Box$

\subsection{}
\label{pjat' unusual action}
The action of $T_{g,K}$ on
$\Tor^\CC_{\frac{\infty}{2}+\bullet}(\sk,L_1\otimes\ldots
\otimes L_m\otimes\ad^{\otimes g})$
is induced by its action on
$L_1\otimes\ldots
\otimes L_m\otimes\ad^{\otimes g}$:
generators act by braiding, balance,
and Fourier transforms of the adjoint representation.

Now suppose $L_k\in\tCO_\zeta$ for all $k\in K$.
Then
$L_1\tilde{\otimes}\ldots\tilde{\otimes}L_m\tilde{\otimes}
\tilde{\ad}^{\tilde{\otimes}g}=\langle L_1\otimes\ldots\otimes L_m\otimes
\ad^{\otimes g}\rangle_\Delta$ (notations of ~\ref{pjat' maximal summand})
by the Theorem ~\ref{pjat' try and prove}. The RHS is a natural subquotient
of $L_1\otimes\ldots\otimes L_m\otimes\ad^{\otimes g}$, hence
the action of $T_{g,K}$ on the latter module induces some action of
$T_{g,K}$ on the RHS. This action coincides with the one of ~\ref{pjat' usual action}
by the uniqueness result of ~\cite{ms} which claims that any two extensions
of bbrt structure on a category to a fusion structure coincide.

\subsection{}
By the virtue of Corollary IV.5.11,
$\langle L_1\tilde{\otimes}\ldots\tilde{\otimes}L_m\tilde{\otimes}
\tilde{\ad}^{\tilde{\otimes}g}\rangle$ is canonically a subquotient
of $\Tor^\CC_{\frac{\infty}{2}+\bullet}(\sk,(L_1\otimes\ldots
   \otimes L_m\otimes\ad^{\otimes g})\otimes L(2\rho_\ell))$.

The above discussion implies that the action ~\ref{pjat' usual action}
of $T_{g,K}$ on the former space is induced by its action
~\ref{pjat' unusual action} on the latter space. More precisely,
$T_{g,K}$ acts on
$\Tor^\CC_{\frac{\infty}{2}+\bullet}(\sk,(L_1\otimes\ldots
\otimes L_m\otimes\ad^{\otimes g})\otimes L(2\rho_\ell))$ in the evident
way, leaving $L(2\rho_\ell)$ alone.

\subsection{}
Let us translate the above discussion back into
geometric language. Let us consider the set $J:=K\sqcup j$.
We denote by $\xi$ the natural projection  $TC^{\circ J}\lra TC^{\circ K}$.

Summarizing the above discussion, we obtain

{\bf Theorem.} Let $L_k\in\tCO_\zeta,\ k\in K$. Then
the local system $\xi^*\langle L_1,\ldots,L_m\rangle^K_C$ on $TC^{\circ J}$
is a natural subquotient of the local system
$[L_1,\ldots,L_m,L(2\rho_\ell)]^J_C$.

\subsection{}
\label{pjat' integral representations}
Let us reformulate the above Theorem in more concrete terms.

Recall that $\nabla_\kappa=\Delta_\ell$ is the first alcove parametrizing
the irreducibles in $\tCO_\kappa$. Let $\lambda_1,\ldots,\lambda_m\in
\nabla_\kappa$ be a $K$-tuple of weights. We consider
the corresponding irreducibles $L(\lambda_1),\ldots,L(\lambda_m)\in
\tCO_\kappa$, and the corresponding local system of conformal blocks
$\langle L(\lambda_1),\ldots,L(\lambda_m)\rangle^K_C$
on the space $TC^{\circ K}$.

On the other hand, consider $\alpha:=\lambda_1+\ldots+\lambda_m+2g\rho_\ell$.
We define $J:=K\sqcup j$, and $\lambda_j:=2\rho_\ell$.
Let $\vlambda$ denote $\lambda_1,\ldots,\lambda_m,\lambda_j$.

Denote by $\CX^\alpha_\vlambda$ the sheaf on $C^\alpha_\vlambda$
obtained by gluing $\CL(\lambda_1),\ldots,\CL(\lambda_m),\CL(\lambda_j)$.

Note that $\CX^\alpha_\vlambda=u_{!*}\CH^\alpha_\vlambda$
where $u$ stands for the open embedding $C^{\circ\alpha}_\vlambda\hra
C^\alpha_\vlambda$.

Recall that $\eta$ stands for the projection $C^\alpha_\vlambda\lra
TC^{\circ J}$.

{\bf Theorem.} The local system of conformal blocks
$\xi^*\langle L(\lambda_1),\ldots,L(\lambda_m)\rangle^K_C$ is
isomorphic to a canonical subquotient of a ``geometric'' local system
$R\eta_*u_{!*}\CH^\alpha_\vlambda$. $\Box$

\subsection{Corollary} {\em The above local system of conformal blocks
is semisimple. It is a direct summand of the geometric local system
above.}

{\bf Proof.} The geometric system is semisimple by Decomposition theorem,
~\cite{bbd}, Th\'{e}or\`{e}me 6.2.5. $\Box$

\newpage

\centerline{\bf Index of Notation}

$\CA(2)^\circ,$
0.5.6

$\CAO_{\nu,\emp},$
II.6.12

$\CA_\nu,$
II.6.12

$\CA^{0}(\emp),$
III.2.2

$\CA^{\valpha}(\vd)\subset\CA_{\alpha},$
III.2.2

$\CA_{\mu}^{\valpha}(\vd),$
III.2.3

$\CA^\alpha_\mu,$
III.2.3

$\CA_c,$
III.2.6

$\CAO^{\alpha}_{\mu},$
III.2.7

$\CAD^{\alpha}_{\mu},$
III.2.7

$\CA^\alpha(L),$
III.7.5

$\CA^\alpha(L;d),$
III.7.5

$\CA^{\valpha}(\tau),$
III.7.5

$\CA_{\vmu}^{\alpha}(L),$
III.7.7

$\CA_{\vmu}^{\alpha}(L;d),$
III.7.7

$\CA_{\vmu}^{\valpha}(\tau),$
III.7.7

$\CAO^{\alpha}_{\vmu}(K),$
III.7.9

$\CAD^{\alpha}_{\vmu}(K),$
III.7.9

$\CA^\alpha_{\mu;1},$
III.10.1

$\CA_{\alpha;1},$
III.11.6

$\widehat{\CA^{\alpha}(K)},$
III.11.8

$\CA^\alpha(K)_d,$
III.11.9

$\bA,$
0.14.13

$\BA_\BR,$
I.3.1

$\BAO_{\CH,\BR},$
I.3.1

$\BAO_{\CH},$
I.3.1

$\BA^\pm_{\BR,H},$
I.4.3

$^{\pi}\BAO^+_{\BR},$
II.8.1

$^{\pi}\BAO_{\emp,\BR},$
II.6.1

$^{\pi}\BAO_{\emp},$
II.6.1

$^\pi\BA,$
II.6.1

$^\pi\BA_\BR,$
II.6.1

$\BA^+_\BR,$
II.7.1

$\BA^{(n)},$
II.9.1

$^n\BA,$
II.9.1

$^n\BAO,$
II.9.1

$\BAO^{(n)},$
II.9.1

$^\bz\BA,$
II.9.1

$\BA_\rho,$
II.9.2

$\BA^1_{>c},$
III.2.1

$\BA^{J_1,J_2}(d),$
III.12.8

$A,$
V.3.2

$A',$
V.5.1

$A_\pm,$
V.4.5

$A_\Lambda,$
V.4.1

$A(r_1,r_2),$   0.3.2

$A^0,$
0.14.13

$A^1,$
0.14.13

$\AJ_\alpha:\ C^\alpha\lra \Pic(C)\otimes {X_\ell},$
V.3.2

$\dot{\AJ_\alpha}:\ (\AJ_\alpha^*(\CL))\dot{}\lra\dot\CL,$
V.3.5

$\AJ_\pm,$
V.4.5

$A\mapsto A',$
II.2.6

$\ad,$
0.12.5

$\ad_\MS,$
0.12.5

$\had,$
0.12.5

$\tilde\ad,$
V.11.2

$\ad_{j,\lambda},$
II.2.20

$\ad_{\theta_i,\lambda}:\ {\ '\ff}\lra{\ '\ff},$
II.2.19

$[a],$
II.1.4

$[a,b],$
II.1.4

$[a]_i,$ 0.2.2

$[a]_\zeta,$
II.1.4

$a:\ \CC\iso\CC_r,$
0.12.3

$\ya_{\alpha}: C^{\alpha}_{/S}\lra\Pic(C/S)\otimes X_{\ell},$
0.14.6

$^{\pi}a:\
{\ '\ff}_{\nu_{\pi}}\lra (^{\pi}{\ '\ff}_{\chi_J})^{\Sigma_{\pi}},$
II.4.2

$^{\pi}a^*:{\ '\ff}_{\nu_{\pi}}^*\iso
(^{\pi}{\ '\ff}_{\chi_J}^*)^{\Sigma_{\pi}},$
II.4.3

$^{\pi}a:\ _{\nu_{\pi}}C^{-r}_{{\ '\ff}}(V(\Lambda_0)\otimes\ldots\otimes
V(\Lambda_{n-1}))\lra\
_{\chi_J}C^{-r}_{^{\pi}{\ '\ff}}(V(^{\pi}\Lambda_0)\otimes\ldots\otimes
V(^{\pi}\Lambda_{n-1}))^{\Sigma_{\pi}},$
II.4.4

$^{\pi}a^*:\
_{\nu_{\pi}}C^{-r}_{{\ '\ff}^*}(V(\Lambda_0)^*\otimes\ldots\otimes
V(\Lambda_{n-1})^*)\lra
_{\chi_J}C^{-r}_{^{\pi}{\ '\ff}^*}(V(^{\pi}\Lambda_0)^*\otimes\ldots\otimes
V(^{\pi}\Lambda_{n-1})^*)^{\Sigma_{\pi}},$
II.4.6

$\CB,$ 0.2.12

$^n\BB,$
II.9.1

$\ob(\lambda),$
V.4.1

$b_M,$      0.2.18

$b(\CL_?)_{F<C},$
I.4.15

$b_{\tau,*},$
II.6.8

$b_{\tau,!},$
II.6.7

$b_{\varrho\leq\tau},$
II.8.9

$b_{\varrho\leq\tau,*},$
II.8.11

$\CC,$ 0.2.11, II.11.3

$\tilde\CC,$
II.11.4

$\CC',$
V.2.1

$\tCC_{\zeta^{-1}},$
III.13.2

$\CC^{\geq0},$ 0.2.16

$\CC^{\leq0},$ 0.2.16

$\CC_{\leq\lambda},$
III.14.1

$\CC_r,$
0.9.1

$_R\CC,$
IV.6.1

$_R\CC_\zeta,$
IV.9.2

$\bC_A,$
0.14.13

$\bC^{\nu}_{A,g},$
0.14.14

$\tC(J)_{/S},$
0.14.2

$\tC^K_\varepsilon,$
V.16.1

$C_s,$
0.14.1

$C^J_{/S},$
0.14.2

$C(J)_{/S},$
0.14.2

$C^{\nu}(K)^{\circ}_g,$
0.14.16

$C(J)_{g;\delta},$
0.14.16

$C(J)_g,$
0.14.16

$C(J;K)_g,$
0.14.16

$C^{o\valpha}_\varepsilon,$
V.16.1

$C^{\valpha}_\varepsilon,$
V.16.1

$C^{\bullet\alpha}_{\vmu},$
V.16.3

$C^{o\alpha}_{\vmu},$
V.16.3

$C^{\alpha}_{\vmu},$
V.16.3

$C^{\valpha,\beta}_{\vmu,\varepsilon},$
V.16.4

$C^{\vnu}_{\vA^0,g,\Gamma},$
0.14.15

$C^{\nu\circ}_{A,g},$
0.14.15

$C^{\nu}_{\vA^0,g,\Gamma},$
0.14.14

$C^{\nu\circ}_{A,g},$
0.14.14

$C^{\nu}_{A,g},$
0.14.14

$\Co^\alpha,$
V.2.1

$\Co^J,$
V.2.1

$C^J,$
V.2.1

$C_\fD,$
V.7.1

$\ol{C},$
V.7.3

$\ul{C}^\alpha_{\vlambda\vmu\vnu},$
V.16.9

$C^\bullet_\CH(\BA;\CM),$
I.3.13

$^+C^{\bullet}(\BA;\CM),$
II.7.11

$C^\bullet_A(M),$
II.3.1

$_{\bnu}C_{{\ '\ff}}^{-r}(M)=\ _{\nu_r,\ldots,\nu_0}C_{{\ '\ff}}^{-r}(M),$
II.3.2

$_{\nu}C_{{\ '\ff}}^{-r}(M),$
II.3.2

$_{\varrho}C^{-r}_{{\ '\ff}}(V(\Lambda_{0})\otimes\ldots\otimes
V(\Lambda_{n-1})),$
II.3.3

$_{\varrho}C^{-r}_{{\ '\ff}^*}(V(\Lambda_{0})^*\otimes\ldots\otimes
V(\Lambda_{n-1})^*),$
II.3.5

$C^{\bullet}_{\ff}(L(\Lambda_0)\otimes\ldots\otimes L(\Lambda_{n-1})),$
II.5.3

$C^{\blt}_{\fu^-}(\Phi(\CX_1)\otimes\ldots\otimes\Phi(\CX_n))^{\heartsuit},$
IV.7.4

$^{ji}C,$
II.6.2

$C_\tau,$
II.6.2

$\Ch(F),$
I.3.2

$\Ch,$
II.6.6

$\Ch^+(F_\varrho),$
II.8.5

$\Cch^i(\CU,\CV;\CF),$
IV.8.5

$\Cch^a_b(\CU,\CV;\CF),$
IV.8.5

$\Cch_a(\CU;\CF),$
IV.8.5

$\Cch^a(\CU;\CF),$
IV.8.5

$c\in X/Y,$ 0.4.3

$cl(\sigma;s),$
I.2.7

$c(\CL_*)_{F<C},$
I.4.8

$c(\CL_!)_{F<C},$
I.4.7

$c_{\tau,!},$
II.6.7

$c_{\Delta<C},$
II.6.7

$c_{\tau,*},$
II.6.8

$c_{\varrho\leq\tau},$
II.8.6

$c_{\varrho\leq\tau,*},$
II.8.11

${\tilde c}_{\varrho\leq\tau},$
II.8.11

$c_i:\ \CA^{\valpha}(\vd)\iso \CA^{\valpha_{\leq i}}(\vd_{\leq i}) \times
\CA^{\valpha_{\geq i}}(\vd_{\geq i}),$
III.2.2

$c_X:\ \widetilde{X\times X}\lra X\times X,$
III.11.7

$\fD,$
V.7.1

$\CCD\CM(\mu),$
0.4.7

$\CCD(X),$
I.2.1

$\CCD^*_{\BR^{*+}}(X),$
I.2.13

$\CCD(\BA;\CS_{\CH}),$
I.3.1

$\CCD^b(^\pi\BA,\Sigma_\pi),$
III.12.5

$\CCD^J,$
V.3.2

$\CCD^\alpha,$
V.3.2

$D,$            0.3.2

${\bar D}(r),$  0.3.2, III.2.1

$D(r),$     0.3.2, III.2.1

$D^{J_1,\ldots,J_n},$  0.3.3

$D(z;r),$
0.5.1, III.7.3

$D(c,c'),$
III.2.1

${\bar D}(z,r),$
0.5.1, III.7.3

$D^\nu(2),$
0.5.2

$D^{\nu,\nu'}(2),$
0.5.2

$D^{\nu_1,\nu_2,\nu}(2),$
0.5.2

$D^{o\alpha},$
IV.2.2

$\Do^\alpha,$
V.2.1

$D^{oJ},$
IV.2.2

$\Do^J,$
V.2.1

$D^{\alpha},$
IV.2.2, V.2.1

$D^J,$
IV.2.2, V.2.1

$D_{(\epsilon,1)},$
IV.2.2

$D_\epsilon,$
IV.2.2, V.2.1

$D^{\alpha,\beta},$
IV.3.3

$D:\ \tCC^{opp}\lra\tCC_{\zeta^{-1}},$
III.13.2

$D:\ \FS^{opp}\lra\FS_{\zeta^{-1}},$
III.13.3

$\DD,$          0.3.2

$\DD^\nu(2),$
0.5.2

$\DD^\nu(2)^\circ,$
0.5.2

$\DD^{\nu,\nu'}(2),$
0.5.3

$\DD^{\nu,\nu'}(2)^\circ,$
0.5.3

$\DD^{\nu_1,\nu_2,\nu}(2),$
0.5.3

$\DD^{\nu_1,\nu_2,\nu}(2)^\circ,$
0.5.3

$\DO_{F<C},$
I.3.2

$\DO_F,$
I.3.2

$\DO_k(\tau),$
III.7.3

$D^J_\BR,$
0.6.1

$D^{J+},$
0.6.1

$D_F,$
0.6.1, I.3.2

$D(J),$
0.7.2

${\tilde D}(J),$
0.7.2

$\DO^+_F,$
II.7.2

$D^+_F,$
II.7.2

$D^\pi,$
0.14.6

$D^\alpha,$
0.14.6

$D_X,$
I.2.5

$d,$ 0.2.1

$d_i,$ 0.2.1

$d_i:\ \CA^{\valpha}(\tau)\lra\CA^{\dpar_i\valpha}(\dpar_i\tau),$
III.7.5

$d_i:\ \CA^{\valpha}_{\vmu}\lra
\CA^{\dpar_i\valpha}_{\dpar_i\vmu}(\dpar_i\tau),$
III.7.7

$\dd_\ell,$
0.12.1, V.2.1

$\det(S_{\Lambda,\nu})(\bq;\br),$
II.2.25

$\det(S_\nu)(\bq),$
II.2.25

$\vd_{\geq i},$
III.2.2

$\vd_{\leq i},$
III.2.2

$\CE,$
0.14.9, V.5.2

$E_i,$     0.2.12, II.11.3, II.12.2

$E^{(p)}_i,$
II.12.2

$E_n,$         0.3.3

$\Ext^{\infh+n}_{\CC}(M,M'),$
0.9.7, IV.4.7

$e:\ \fu\lra \sk,$ 0.2.8

$e:\ \BN[I]\lra \BN[X],$
IV.3.1

$e:\ \BN[I]\lra\BN[X_\ell],$
V.16.2

$\FS,$ 0.4.6, III.5.2

$\FS',$
V.2.1

$\FS_c,$ 0.4.6, III.5.2

$\FS_{\zeta^{-1}},$
III.13.2

$\FS^\diamondsuit,$
V.12.1

$\FS_{\leq\lambda},$
III.14.1

$\FS^{(m)},$
0.15.8

$\FS^{\sharp},$
0.15.8

$^K\tFS,$
III.9.4

$^K\FS,$
III.9.6

$\FS^{\otimes K},$
III.10.2

$\tFS,$                   0.4.5

$\tFS_c,$                 0.4.5

$\Fac^1(E),$
I.3.11

$\bF,$
0.6.1, I.3.2

$F_i,$     0.2.12, II.11.3, II.12.2

$F^{(p)}_i,$
II.12.2

$F_j,$
III.12.10

$F_j^\perp(d),$
III.12.11

$F_{J'},$
III.12.15

$F_J(x_j,\xi_j),$
V.6.3

$F_\CO,$
0.12.3

$F_{\CO;\infh+n}:\ \Ind\CC\lra\Ind\CC,$
0.12.3

$F_\pm,$
I.4.3

$F(\bq^2),$
I.4.4

$F_\rho,$
II.6.5

$F_\varrho,$
II.8.4

$'\ff,$ 0.2.5, II.2.2

$'\ff_\nu,$ 0.2.5, II.2.2

${\ '\ff}^*_\nu,$
II.2.9

${\ '\ff}^*,$
II.2.9

${\ '\ff}_Q,$
II.2.20

$^\pi{\ '\ff}_\nu,$
II.4.1

$^\pi{\ '\ff},$
II.4.1

$\ff,$
II.5.1

$f:\ C\lra S,$
0.14.1

$f_*:\tT^{\nu}_{A,g}\lra\tT^{\nu}_{B,g},$
0.14.20

$f^*\langle\otimes_A\ \{\CX_a;\CR_n\}\rangle_g^{(m)},$
0.15.5

$f^*(\Boxtimes^{(\nu)}_{A,g}\ \{\CX_a;\CR_n\}),$
0.15.5

$f_H,$
I.4.3

$f^\valpha(\tau_d),$
III.10.3

$\operatorname{g}(\{\CX_k\}),$
V.16.6

$\operatorname{g}(\{\CX_k\},\{\CY_j\}),$
V.16.8

$\ul{\og}(\{\CX_k\},\{\CR_j\})^{\alpha}_{\vlambda\vmu\vnu},$
V.16.9

$\fg,$
0.10.6, III.18.2

${\hat\fg},$
0.10.6, IV.9.2

$g,$
0.14.1

$g_K:\ \FS^{\otimes K}\lra\ ^K\FS,$
III.10.3

$\fH,$
0.14.10, V.3.4

$\fH^N,$
V.3.6

$\CH^\alpha,$
V.2.2

${\tilde\CH}{}^\alpha,$
V.6.2

$\CH,$
V.2.2, V.16.2

$\CH^\tau,$
V.6.3

$\CH(\pi),$
0.7.7

$\CH^{\nu\bullet}_{\vmu},$
0.8.3

$\CH^{\nu\sharp}_{\vmu},$
0.11.4

$\CH(\pi),$
0.14.4

$\CH^{\nu}_{\vmu;A,g},$
0.14.18

$\CH_\chi,$
0.14.18

$\CH^{\nu}_{\vmu},$
0.14.16

$\tCH_{\vmu}^{\nu},$
0.14.16

$\CH^{\vnu}_{\vmu;\vA^0,g,\Gamma},$
0.14.15

$\CH^{\nu}_{\vmu;A,g},$
0.14.15

$\CH^{\bullet\alpha}_\vmu,$
V.16.5

$\CH_L,$
I.4.13

$\CH_\emptyset,$
II.6.1

$\Hom(\CM,\CN),$          0.4.5

$\Hom_{\tFS}(\CX,\CY),$
III.4.4

${\overline{\Hom}}_{^K\FS}(\CX,\CY),$
III.9.6

$\Hom_{^K\FS}(\CX,\CY),$
III.9.6

$\Hom_{^K\tFS}(\CX,\CY),$
III.9.4

$\Hom_{\tCO_\kappa}(\One,L_1\totimes\ldots\totimes L_m)^{\heartsuit},$
IV.9.2

$\Hom_{\CB}(\One,L_1\otimes\ldots\otimes L_m)^{\heartsuit},$
IV.9.2

$H^i(X,Y;\CK),$
I.2.4

$H_i(X,Y;\CF),$
I.2.7

$H_j(z),$
II.7.1

$H_w,$
II.9.2

$\HO_{\BR,\CK},$
I.3.1

$\HO_{\CK},$
I.3.1

$h,$
0.10.6, IV.9.2

$\hgt(\tau),$
III.7.2

$\CI,$
IV.3.1

$\CI^\alpha,$
V.2.1

$\CI^J,$
V.2.1, V.6.3

$\CI=\CJ\otimes\Sign,$ 0.3.10

$\CI_\mu^{\nu\bullet},$  0.4.2

$\CI^\nu_{\mu_1,\mu_2},$
0.5.6

$\CI^{\nu\bullet}_{\mu_1,\mu_2},$
0.5.7

$\CI_{\nu!*},$
II.6.13

$\CI_{\nu*},$
II.6.13

$\CI_{\nu!},$
II.6.13

$\CI_\nu,$
II.6.12

$\CI_\nu(\Lambda)_?,$
II.8.16

$\CI_\nu(\Lambda),$
II.8.16

$\CI^\alpha_\mu,$
III.3.1

$\CI^{\valpha}_{\mu}(\vd),$
III.3.2

$\CI^\alpha_D,$
IV.2.5

$\CID_{\mu}^{\valpha}(\vd),$
III.3.5

$\CI^{\bullet\alpha}_{\vmu},$
IV.3.4

$^\pi\CI_{!*},$
II.6.4

$\CI(^\pi\Lambda),$
II.8.2

$\CI(^{\pi}\Lambda)_{!*},$
II.8.3

$\CI(\Lambda_0,\ldots,\Lambda_{-n};\nu)_?,$
II.10.5

$\CI(\Lambda_0,\ldots,\Lambda_{-n};\nu),$
II.10.5, II.12.6

$\CI^\valpha_\vmu(\tau),$
III.8.1

$\CI^\alpha_\vmu(L),$
III.8.1

$\CI^\alpha_\vmu,$
III.8.1

$^\pi\CI_\vmu,$
III.8.1

$^K\CI,$
III.8.4

$\CID^{\valpha}_{\vmu}(\tau),$
III.8.5

$(I,\cdot),$ 0.2.1

$I^\diamondsuit,$
V.12.1

$\Ind\CC,$
0.12.3

$\Ind^\fu_{\fp_i},$
III.14.10

$\vec{I}_Q,$
II.2.20

$\ind_{\fu^{\geq0}}^\fu:\ \CC^{\geq0}\lra\CC,$ 0.2.16

$\ind_{\fu^{\leq0}}^\fu:\ \CC^{\leq0}\lra\CC,$   0.2.16

$i_{Y,X},$
I.2.2

$i',$
II.12.1

$\CJ,$ 0.3.10

$\CJ^\nu_{\mu_1,\mu_2},$
0.5.4

$\CJ(\pi),$
0.7.5

$\CJ^{\nu\bullet}_{\mu_1,\mu_2},$
0.13.2

$\CJ^{\bullet\alpha}_{\mu_1\mu_2},$
V.13.2

$\CJ^{\bullet\alpha,\beta,\gamma}_{\mu_1\mu_2}(d),$
V.14.2

$\Jac,$
0.14.5

$\vec{J}_\tau,$
II.6.9

$J_{\varrho\leq\tau},$
II.8.5

$^nJ,$
II.9.1

$j(Q),$
II.2.24

$(j'j''),$
III.12.9

$j_{\nu-i'\leq\nu!},\ j_{\nu-i'\leq\nu*}:\
\FS_{\leq\nu}\lra\FS_{\leq\nu},$
III.14.5

$\CK^*,$
I.2.3

$\CK_L,$
I.3.4

$\CK|_Y,$
I.2.2

$K_\nu,$ 0.2.7

${\tilde K}_\nu,$ 0.2.7

${\tilde K}_i,$
II.12.2

$K_i^{\pm1},$
II.11.3, II.12.2

$\vec{K},$
II.2.3

$\vec{K}_A,$
II.2.6

$K'(\tau),$
III.7.2

$K(\tau),$
III.7.2

$K_\tau,$
III.7.2

$K_a,$
III.7.2

$K^\bullet_n,$
IV.5.2

$K^\bullet,$
IV.5.2

$\sk,$ 0.2.2, I.2.1

$\fL(\lambda),$
0.10.6

$\CL,$
V.3.3

$\CL',$
V.5.1

$\CL(\Lambda,\ (,),\ \nu),$
V.4.2

$\CL(\mu),$
0.4.8, III.6.1

$\CL(\pi),$
0.14.6

$\CL(\alpha),$
0.14.6

$\CL_\alpha,$
V.3.2

$\CL_\tau,$
V.6.3

$\CL_\alpha^N,$
V.3.6

$\dCL_\delta,$
0.14.8

$\CL_\Lambda,$
0.14.8

$\CL(\bq),$
I.4.2

$\CL_{!*},$
I.4.5

$\dot\CL,$
V.2.1

$\Locsys(Y;\CC),$
III.11.1

${\tilde L}(\lambda),$
IV.6.4

$\hL(\lambda),$
IV.9.5

$L(\lambda),$ 0.2.17, II.5.2

$\langle L(\lambda_1),\ldots,L(\lambda_n)\rangle,$
0.10.3

$\langle L_1,\ldots,L_m\rangle^K_C,$
V.19.8

$[L_1,\ldots,L_m]^K_C,$
V.19.4

$\fl_i,$
III.14.10

$l,$ 0.2.2

$\ell,$ 0.2.2, IV.9.1

$\ell_\beta,$ 0.2.2, IV.9.1

$\ell_i,$
IV.9.1

$\CM(T,\CS),$ 0.2.4

$\CM(\mu),$
0.4.7

$\CM\boxtimes\CN,$
0.5.8

$\CM\otimes\CN,$
0.5.9

$\MS,$
0.10.6

$\CM_{\vA^0,g,\Gamma},$
0.14.14

$\CM_{A,g},$
0.14.14

$\CM_A,$
0.14.14

$\bCM_{g;\delta},$
0.14.14

$\bCM_g,$
0.14.14, V.7.1

$\CM_g,$
0.14.14, V.7.1

$\CM(\BA;\CS_{\CH}),$
I.3.1

$\CM_F,$
I.3.10

$\CM^+_F,$
II.7.9

$M^G,$
II.4.2

$M(\lambda),$ 0.2.17

$M^+(\lambda),$
IV.5.1

$M\mapsto M^*,$ 0.2.14

$M\mapsto M^\vee:\ \CC^\opp\iso\CC_r,$
0.9.1

$\langle M\rangle,$
0.10.1, IV.5.9

$\langle M\rangle_{\MS},$
0.10.7

$\langle M\rangle_\Delta,$
V.11.2

$M^{\Sigma_\pi,-},$
II.6.12

$m(J_1,\ldots,J_n):\ D^{J_1,\ldots,J_n}\lra D^J,$ 0.3.3

$m(\nu_1,\ldots,\nu_n):\
D^{\nu_1,\ldots,\nu_n}\lra D^{\nu_1+\ldots+\nu_n},$ 0.3.3

$m_a,$
II.8.5

$m_a(\nu_1,\ldots,\nu_n):\ D^{\nu_1,\ldots,\nu_n}\lra D^{\nu_1,\ldots,
\nu_{a-1},\nu_a+\nu_{a+1},\nu_{a+2},\ldots,\nu_n},$  0.3.4

$m_a(\nu_1,\ldots,\nu_n):\ \DD^{\nu_1,\ldots,\nu_n}\lra \DD^{\nu_1,\ldots,
\nu_{a-1},\nu_a+\nu_{a+1},\nu_{a+2},\ldots,\nu_n},$          0.3.4

$m(\rho):\ \prod_J{\tilde D}(K_j)\times
{\tilde D}(J)\lra{\tilde D}(K),$
0.7.2

$m_P(\rho):\ \prod_J \tD(K_j)\times \tP(J)\lra\tP(K),$
0.7.6

$\bar{m}_P(\rho):\ \prod_J D(K_j)\times\tP(J)\lra P(K),$
0.7.6

$m(\vnu;\nu'),$
0.8.1

$m_{C/S}(\rho):\ \prod_J \tD(K_j)\times\tC(J)_{/S}\lra \tC(J)_{/S},$
0.14.2

$\bar{m}_{C/S}(\rho):\ \prod_J D(K_j)\times\tC(J)_{/S}\lra C(K)_{/S},$
0.14.2

$m:\ \CL_!\lra\CL_*,$
I.4.5

$m:\ ^{\pi}\CI_!\lra\ ^{\pi}\CI_*,$
II.6.4

$m:\ \CI(^{\pi}\Lambda)_!\lra\CI(^{\pi}\Lambda)_*,$
II.8.3

$m:\ \CM(\Lambda)\lra D\CM(\Lambda)_{\zeta^{-1}},$
III.6.1

$N\otimes_\CC M,$
0.9.6

$N_A,$
0.14.13

$N_i,$
II.4.1

$\vec{N},$
II.2.20

$n(\lambda),$ 0.2.21, IV.2.4, IV.9.1, V.16.2

$n_\nu(\mu),$
0.14.8

$(n),$
II.9.1

$\CO(K),$
III.7.3

$\CO(K;d),$
III.7.3

$\CO(\tau),$
III.7.3

$\CO,$
III.18.2

$\CO_\fg,$
III.18.2

$\CO_{-\kappa},$
IV.9.2

$\tCO_\kappa,$
IV.9.2

$\tCO_{-\kappa},$
IV.9.2

$\Or_X,$
I.2.6

$\Ord(\varrho),$
II.3.4

$\wp_f^{(\nu\circ)}:C^{\nu\circ}_{A,g;\delta}\lra
\dT_{\dpar\overline{C^{\nu\circ}_{B,g;\delta}}}\ \overline
{C^{\nu\circ}_{B,g;\delta}},$
0.14.20

$\wp_f: \CM_{A,g;\delta}\lra\dT_{\dpar\bCM_{B,g;\delta}}\bCM_{B,g;\delta},$
0.14.20

$\wp:\ \tilde{{\fD}}\rightarrow\TB_{\fD}\ol{\CM_g},$
V.7.2

$\wp_\alpha:\ \tilde{{\fD}}\times_{\fD}{\cal T}C^\alpha_{\fD}\lra
\TB_{{\cal T}C_{\fD}^\alpha}\TO\ol{C}^\alpha,$
V.7.4

$\wp_\alpha:\ \ul{\oC}^\alpha_{\vmu,\mu_1,\mu_2}\lra
\TB_{(\ol{C}|_{\fD})^\alpha_\vmu}\ol{C}^\alpha_{\vmu},$
V.17.2

$\CP_r(I;n),$
II.3.3

$\CP_r(J),$
II.6.5

$\CP^{o\alpha},$
IV.2.1

$\CP^{oJ},$
IV.2.1

$\CP^\alpha,$
IV.2.1

$\CP^J,$
IV.2.1

$\CP^{o\valpha},$
IV.2.3

$\CP^{\valpha},$
IV.2.3

$\CP^{\bullet\alpha}_{\vmu},$
IV.3.2

$\CP^{o\alpha}_{\vmu},$
IV.3.2

$\CP^{\alpha}_{\vmu},$
IV.3.2

$\CP^{\valpha,\vbeta,\gamma}_{\vmu},$
IV.3.3

$\CP^{\valpha,\beta}_{\vmu},$
IV.3.3

$\CP^{o\valpha_L;\xi},$
IV.2.3

$\CP^{\valpha_L;\xi},$
IV.2.3

$\tCP^K,$
IV.2.3

$\BP^1_{st},$
V.13.1

$P=\BP^1(\BC),$
0.7.6

$P(J),$
0.7.6

$\tP(J),$
0.7.6

$P^\nu(K),$
0.8.1

$P^\nu(K)^\circ,$
0.8.1

$P^\nu(K)^\bullet,$
0.8.1

$P^{\vnu;\nu'},$
0.8.1

$\tP^{\nu}(K),$
0.8.1

$P^{\vnu^1;\vnu^2;\nu},$
0.8.1

$\tP^{\vnu;\nu'},$
0.8.1

$P(K,J),$
0.8.2

$P(n),$
0.10.3

${\tilde P}_\rho,$
II.9.2

${\tilde P},$
II.9.2

$P,$
II.9.2

$\Pic(C/S),$
0.14.5

$\fp_i,$
III.14.10

$p(\vnu;\nu'),$
0.8.1

$[p]^!_i,$
II.12.2

$p(J),$
0.7.2

$\pr_1,$
V.5.1

$pr,$
V.6.1

$pr_1,$
V.6.1

$pr_2,$
V.6.1

$\pr^\alpha,$
V.15.1

$\pr:\ Q_a\lra\BA^1_h\times\BC,$
V.15.1

$p:\ A\lra S,$
V.5.2

$p_K:\ \tCP^K\lra T\CP^{oK},$
IV.2.7

$p^{\alpha}_{\vmu}:\ T\CP^{o\alpha}_{\vmu}\lra T\CP^{o\alpha_{\vmu}},$
IV.3.2

$p_{\valpha,\beta}:\ \CP^{\valpha,\beta}_{\vmu}\lra
\prod_k\ D^{\alpha_k}\times\CP^{\bullet\beta}_{\vmu-\valpha},$
IV.3.3

$p_{\valpha,\vbeta,\gamma}:\ \CP^{\valpha,\vbeta,\gamma}_{\vmu}\lra
\prod_k\ D_{\frac{1}{2}}^{\alpha_k}
\times\prod_k\ D_{(\frac{1}{2},1)}^{\beta_k}\times
\CP^{\bullet\gamma}_{\vmu-\valpha-\vbeta},$
IV.3.3

$p_a(\nu_1,\ldots,\nu_n):\ \DD^{\nu_1,\ldots,\nu_n}\lra
\DD^{\nu_1,\ldots,\nu_a}\times \DD^{\nu_{a+1},\ldots,\nu_n},$ 0.3.4

$p_a(\nu_1,\ldots,\nu_n):\ D^{\nu_1,\ldots,\nu_n}\lra
D^{\nu_1,\ldots,\nu_a}\times \DD^{\nu_{a+1},\ldots,\nu_n},$  0.3.4

$p(\nu_1,\ldots,\nu_n):\ D^{\nu_1,\ldots,\nu_n}\lra D^{\nu_1}\times
\DD^{\nu_2}\times\ldots\times \DD^{\nu_n},$  0.3.3

$p(J_1,\ldots,J_n):\ D^{J_1,\ldots,J_n}\lra D^{J_1}\times \DD^{J_2}\times
\ldots\times \DD^{J_n},$  0.3.3

$p^{\beta}_{\lambda}:\
\FS_{c;\leq\lambda}\lra\CM(\CA^{\beta}_{\lambda};\CS),$
III.5.4

$p:\ \Jac_\Lambda\lra\Jac,$
0.14.8

$^np:\ ^n\BA\lra\BA^{(n)},$
II.9.1

$p_K:\ ^K\FS\lra\FS^{\otimes K},$
III.10.2

$p_K:\ \widetilde{\TCo^K}\lra\TCo^K,$
V.2.1

$p_{\valpha,\beta}:\ C^{\valpha,\beta}_{\vmu,\varepsilon}\lra
\prod_k\ D^{\alpha_k}_\varepsilon \times C^{\bullet\beta}_{\vmu-\valpha},$
V.16.4

$p^{\alpha}_{\vmu}:\ T C^{o\alpha}_{\vmu}\lra T C^{o\alpha_{\vmu}},$
V.16.3

$\CQ,$
V.6.3

$\CQ_2(J),$
II.9.2

$Q_v,$
0.13.1

$Q_{0v},$
V.13.1

$Q_{0h},$
V.13.1

$Q_u,$
0.13.1

$Q_t,$
0.13.1

$Q,$
0.13.1

$Q_{>d,<d},$
V.14.2

$Q_{<d,>d},$
V.14.2

$Q_{>d,>d},$
V.14.2

$Q_{<d,<d},$
V.14.2

$Q^{\alpha,\beta,\gamma}(d),$
V.14.2

$'Q,$
0.13.1

$'Q^{\nu\bullet},$
0.13.2

$'Q^{\bullet\alpha},$
V.13.2

$'Q^{\circ\alpha},$
V.13.2

$Q:\ \tilde\CC\lra\CC,$
II.11.7, II.12.5

$\bq:\ \CH\lra\BC^*,$
I.4.2

$\bq(C,C'),$
I.4.11

$q_{ij},$
II.2.25

$q_\lambda,$
III.14.1

$q^*_{\lambda\leq\mu},$
III.14.8

$q^!_{\lambda\leq\mu},$
III.14.8

$q^*_\lambda,$
III.14.8

$q^!_\lambda,$
III.14.8

$q_{\valpha,\vbeta,\gamma},$
IV.3.3

$q_\tau,$
V.2.1, V.6.1

$q_{\lambda\leq\mu}:\
\CC_{\leq\lambda}\hra\CC_{\leq\mu},$
III.14.8

$q_{\valpha_L;\xi}^2:\ \CP^{\valpha_L;\xi}\lra\CP^{\valpha_L},$
IV.2.3

$q_{\valpha_L;\xi}^1:\ \CP^{\valpha_L;\xi}\lra\CP^{\valpha},$
IV.2.3

$q_{\valpha}:\ \CP^{\valpha}\lra\DT\CP^{\alpha},$
IV.2.3

$q_{\valpha,\beta}:\ \CP^{\valpha,\beta}_{\vmu}\lra\CP^{\alpha+\beta}_{\vmu},$
IV.3.3

$q_{\alpha,\beta}:\ D^{\alpha,\beta}\lra D^{\alpha+\beta},$
IV.3.3

$q^1_{\valpha,\vbeta,\gamma}:\ \CP^{\valpha,\vbeta,\gamma}_{\vmu}\lra
\CP^{\valpha+\vbeta,\gamma}_{\vmu},$
IV.3.3

$q^2_{\valpha,\vbeta,\gamma}:\ \CP^{\valpha,\vbeta,\gamma}_{\vmu}\lra
\CP^{\valpha,\beta+\gamma}_{\vmu},$
IV.3.3

$q_{\valpha}:\ C^{\valpha}_\varepsilon\lra\TB C^{\alpha},$
V.16.1

$q_{\valpha,\beta}:\ C^{\valpha,\beta}_{\vmu,\varepsilon}\lra
C^{\alpha+\beta}_{\vmu},$
V.16.4

$\CR,$
0.13.3

$\CR^{\nu,\nu}_{\mu_1,\mu_2},$
0.13.3, V.14.2

$\CR_n,$
0.15.2

$\Rep^{\nu}_{c;A,g},$
0.14.21

$\tRep_{c;A},$
0.15.8

$\tRep_{c},$
0.15.8

$\bR_\MS,$
0.16.1

$\bR,$
0.12.3

$R\Gamma(X,Y;\CK),$
I.2.4

$R^\FS,$
IV.6.6

${\tilde R}_{V,W},$
IV.6.1

$R_{M,M'},$ 0.2.18

$R_{\CM,\CN}:\ \CM\otimes\CN\iso\CN\otimes\CM,$
0.5.10

$R:\ U\lra\fu,$
II.11.7, II.12.5

$r:\ '\ff\lra\ '\ff\otimes\ '\ff,$ 0.2.5

$r_i,$
II.2.25

$r_\rho,$
II.9.2

$r:\ R\Gamma(\BA,H_w;\CK)\lra\oplus_{\rho\in\CQ_2(J)}
R\Gamma(\BA_{\rho},\tP_{\rho};\CK),$
II.9.3

$r_{\pi}:\ \Locsys(\CO(K);\CC)\lra
\Locsys(\CO(L)\times\prod_l\ \CO(K_l);\CC),$
III.11.2

$\Sign,$ 0.3.9

$\Sew,$
0.14.13

$\Sh(X),$
I.2.1

$\Sh_{\BR^{*+}}(X),$
I.2.13

$\CS_\CH,$
I.3.1

$\CS_{\emptyset,\BR},$
II.6.1

$\CS_\emptyset,$
II.6.1

$\CS,$
II.7.1

$^\bz\CS,$
II.9.1

$\CS_\Delta,$
II.9.1

$\bSp\CH_g^\alpha,$
V.7.6

$S_F,$
0.6.1, I.3.2

$S^+_F,$
II.7.2

$S_\delta,$
0.14.1

$\Supp(\alpha),$
0.14.6

$S:\ {\ '\ff}\lra{\ '\ff}^*,$
II.2.13

$S(\ ,\ ):\ {\ '\ff}\otimes{\ '\ff}\lra \sk,$
II.2.10

$S_\Lambda:\ V(\Lambda)\otimes V(\Lambda)\lra \sk,$
II.2.16

$S_\Lambda:\ V(\Lambda)\lra
V(\Lambda)^*,$
II.2.21, II.12.4

$S_{m;\Lambda_0,\ldots,\Lambda_{n-1}},$
II.2.22

$S_{^{\pi}\Lambda},$
II.4.7

$S:\ \fu_{\beta}^{-*}\iso\fu_{\beta}^{-},$
III.12.19

$S':\ \fu_{\beta}^{-*}\iso\fu_{\beta}^+,$
III.12.19

$s_\alpha,$
V.3.5

$s_\tau,$
V.6.3

$s'_\tau,$
V.6.3

$s=\Omega\otimes\nu,$
V.5.3

$s:\ \CC\iso\CC_r,$
IV.4.6

$s:\ \fu\lra\fu^\opp,$ 0.2.9, IV.4.6

$s':\ \fu\lra\fu^\opp,$ 0.2.9

$\sgn(F,E),$
I.4.6

$\sgn(\tau_1,\tau_2),$
II.6.9

$\sgn(\varrho<\varrho'\leq\tau),$
II.8.8

$\sgn(\varrho),$
II.8.9

$\sgn(\tau,\eta),$
II.8.9

$\supp\alpha,$
IV.2.1

$\CT_j,$
0.14.6

$\tCT^{\nu}_{A,g},$
0.14.19

${\cal T}C_{\fD}^\alpha,$
V.7.3

$\tTeich_{A},$
0.15.8

$^CT_{ij}:\CL_C\lra\CL_{^{ji}C},$
II.6.2

$^CT_{i0},$
II.8.1

$T^{\nu o},$     0.3.1

$T^{J o},$       0.3.1

$T^\nu,$         0.3.1

$T^J,$           0.3.1

$T_\mu,$
V.2.1

$T_j,$
V.2.1

$T_{g,K},$
V.19.2

$T_{C/S},$
0.14.2

$TC^J,$
V.2.1

$\TCo^\alpha,$
V.2.1

$\TB C^{\alpha}_{\vmu},$
V.16.3

$T C^{o\alpha}_{\vmu},$
V.16.3

$\TO\ol{C}^\alpha,$
V.7.3

$\TCo^J,$
V.2.1

$TC^\tau,$
V.6.1

$TD^{o\alpha},$
IV.2.2

$TD^{oJ},$
IV.2.2

$\TDo^J,$
V.2.1

$TD^{o\valpha},$
IV.2.3

$\TDo^\alpha,$
V.2.1

$TD^{\valpha},$
IV.2.3

$T\CP^{oJ},$
IV.2.1

$T\CP^{o\alpha},$
IV.2.1

$T\CP^{o\alpha}_{\vmu},$
IV.3.2

$\ul{T\CP}^{oJ},$
IV.8.1

$\DT\CP^{\alpha}_{\vmu},$
IV.3.2

$\DT\CP^J\subset T\CP^J,$
IV.2.1

$\DT D^J,$
IV.2.2

$\Tor^{\CC}_{\infh+n}(N,M),$
0.9.7, IV.4.7

$\Tor^{\CC}_{\infh+\blt}(\sk,\Phi(\CX_1)\otimes\ldots\otimes\Phi(\CX_m))
^{\heartsuit},$
IV.8.10

$\Tor_r^A(\sk,M),$
II.3.1

$t_i:\ {\ '\ff}^+\otimes V(\Lambda)\lra{\ '\ff}^+\otimes V(\Lambda),$
II.2.18

$t_j,$
II.7.1

$U,$
II.11.3, II.12.2

$U(\beta),$
III.12.14

$U_\sk,$
IV.6.1

$\fu,$ 0.2.7, II.11.4, II.12.3

$\fu^0,$ 0.2.15

$\fu^-,$ 0.2.15

$\fu^+,$ 0.2.15

$\fu^{\leq0},$ 0.2.15

$\fu^{\geq0},$    0.2.15

$\fu_\lambda,$    0.2.15

${\tilde\fu},$
0.12.2

$\dfu,$
0.12.2

$\fu',$
0.12.2, V.2.1

$\bu_\CB,$ 0.2.12

$\bu_\sk,$   0.2.12

$u_E^F(\CK):\ \Phi^+_F(\CK)\lra\Phi^+_E(\CK),$
0.6.6

$^+u^F_E:\ \Phi^+_F\lra \Phi^+_E,$
II.7.10

$^+u^*(c_{\varrho\leq\tau}),$
II.8.8

$u^F_E(\CK),$
I.3.11

$\CV(\alpha),$
III.10.3

$\Ve(\bGamma),$
0.14.13

$V_X,$
I.2.2

$V(\Lambda)_\lambda,$
II.2.15

$V(\Lambda),$
II.2.15

$V(\Lambda)^*_\nu,$
II.2.21

$V(\Lambda)^*,$
II.2.21

$V(^\pi\Lambda),$
II.4.4

$_iV,$
III.14.11

$v^E_F(\CK),$
I.3.11

$v_F^E(\CK):\ \Phi^+_E(\CK)\lra\Phi^+_F(\CK),$
0.6.6

$^+v^E_F:\ \Phi^+_E\lra\Phi^+_F,$
II.7.10

$v_j:\ \Phi^+_{\{0\}}(\pi^*\CK)\rightleftharpoons\Phi^+_{F_j}(\pi^*\CK)\ : u_j,$
III.12.10

$v_{J'}:\ \Phi^+_{\{0\}}(\pi^*\CK)\rightleftharpoons
\Phi^+_{F_{J'}}(\pi^*\CK)\ :u_{J'},$
III.12.15

$v_{\beta}:\ \Phi^+_{\{0\}}(\pi^*\CK)\rightleftharpoons
\oplus_{J'\in U(\beta)}\ \Phi^+_{F_{J'}}(\pi^*\CK)\ :u_{\beta},$
III.12.15

$\bw,$
I.3.2

$^Fw,$
I.3.2

$w_\rho,$
II.6.5

$w_\varrho,$
II.8.4

$w\cdot\lambda,$
III.18.3

$\CX^{\valpha}_{\mu}(\vd),$
III.4.4

$\CX^{\alpha}_{\lambda},$
III.4.4

$\CX^\valpha_\vmu(\tau),$
III.9.3

$\CX^\alpha_\vmu(K),$
III.9.3

$\CX^{\valpha}(\tau_d),$
III.10.3

$\CX^\alpha_\vmu,$
IV.8.1

${\bar\CX}{}^\alpha_{\lambda\mu},$
V.15.2

$\CX\star\CR,$
V.15.3

$X,$ 0.2.1

$X_\ell,$
0.12.1, V.2.1

$X^\diamondsuit,$
V.12.1

$\langle X_1,\ldots X_m\rangle,$
IV.9.2

$\langle X\rangle_{_R\CC},$
IV.9.2

$|x|,$
II.2.2

$Y,$ 0.2.1

$Y_\ell,$ 0.2.3, IV.9.1

$Y^\diamondsuit,$
V.12.1

$\ul\CZ,$
IV.8.2

$z_Q,$
II.2.24

$^n\alpha_i:\ ^{n-1}\psi_i\circ\ ^n\psi_i\iso\
^{n-1}\psi_i\circ\ ^n\psi_{i+1},$
II.10.2

$\valpha_{\geq i},$
III.2.2, III.7.4

$\valpha_{\leq i},$
III.2.2, III.7.4

$\valpha_L,$
IV.2.3

$\valpha_k,$
IV.2.3

$\alpha^\sim,$
IV.2.6, V.2.1

$\alpha_\vmu,$
IV.3.1, V.16.2

$\alpha_g(\vmu),$
V.17.1

$\alpha(\vmu),$
IV.8.3

$\alpha_K,$
V.2.1

$\alpha_\pm,$
V.4.5

$\beta_0,$
0.10.2, IV.9.1

$\Gamma=(\bGamma,\bg),$
0.14.13

$\bGamma(\bC_A),$
0.14.13

$\gamma_0,$
0.10.2, IV.9.1

$^C\gamma_{ij},$
II.6.2

$^C\gamma_{i0},$
II.8.1

$\Delta,$
II.6.6

$\Delta_l,$
IV.5.10, IV.9.1

$\Delta:\ \fu\lra\fu\otimes\fu,$ 0.2.8, II.11.4, II.12.3

$\Delta:\ {\ '\ff}\lra{\ '\ff}\otimes{\ '\ff},$
II.2.4

$\Delta:\ U\lra U\otimes U,$
II.11.3, II.12.2

$\Delta_{{\ '\ff}\otimes{\ '\ff}}:\
{\ '\ff}\otimes{\ '\ff}\lra({\ '\ff}\otimes{\ '\ff})\otimes
({\ '\ff}\otimes{\ '\ff}),$
II.2.5

$\Delta^{(N)}(\theta_{\vec{K}})^+,$
II.2.8

$\Delta^{(N)}:\ {\ '\ff}\lra{\ '\ff}^{\otimes N},$
II.2.7

$\Delta_\bF,$
0.6.1

$\Delta_\ell,$
0.10.2

$\Delta_{\mu\nu},$
0.14.6

$\Delta_{ij},$
II.6.1

$\Delta_\Lambda:\ V(\Lambda)\lra{\ '\ff}\otimes V(\Lambda),$
II.2.18

$\Delta_{\Lambda}^*:\ {\ '\ff}\otimes V(\Lambda)^*\lra V(\Lambda)^*,$
II.2.28

$^\bF\Delta,$
I.3.2

$\delta_{C/S},$
0.14.1

$\delta,$
V.2.2

${\dot\delta},$
V.2.2

$\ol{\delta_g},$
0.14.14, V.7.5

$\delta:\ {\ '\ff}^*\lra{\ '\ff}^*\otimes{\ '\ff}^*,$
II.2.9

$\delta^{(N)}:\
{\ '\ff}^*\lra{\ '\ff}^{*\otimes N},$
II.2.9

$\delta_i:\ {\ '\ff}\lra{\ '\ff},$
II.2.11

$\epsilon_i,$ 0.2.7, II.11.4, II.12.3

$\epsilon_i:\ \Phi(\CM)_{\nu-i}\lra\Phi(\CM)_\nu,$
0.6.9

$\epsilon_i:\ \Phi_{\lambda}(\CX)\rightleftharpoons
\Phi_{\lambda+i'}(\CX)\ :\theta_i,$
III.12.11

$\epsilon_{\beta}:\ \Phi_{\lambda}(\CX)\rightleftharpoons
\fu^-_{\beta}\otimes\Phi_{\lambda+\beta'}(\CX)\ :\theta_{\beta},$
III.12.17

$\epsilon_i:\ V(\Lambda)\lra V(\Lambda),$
II.2.16

$\epsilon,$
II.2.2

$\zeta,$ 0.2.2

$\zeta_i,$ 0.2.2, II.12.1

$\zeta(\vec{K};\tau),$
II.2.3

$\zeta',$
II.10.4

$\breta^\alpha_J,$
IV.8.1

$\breta_J,$
IV.8.1

$\eta:\ P^\nu(K)\lra P(K),$
0.11.1

$\eta:\ C_{A,g}\lra\CM_{A,g},$
0.14.14

$\eta:\ \CA^\alpha(K)\lra\CO(K),$
IV.7.1

$\eta^{\vbeta;d}_{\vmu;\ '\vbeta;\varepsilon}:\
\Hom_{\CA^{(0,\vbeta)}_{\vmu}(\tau_d)}(\CX^{(0,\vbeta)}_{\vmu}(\tau_d),
\CY^{(0,\vbeta)}_{\vmu}(\tau_d))\lra
\Hom_{\CA^{(0,\ '\vbeta)}_{\vmu}(\tau_\varepsilon)}
(\CX^{(0,\ '\vbeta)}_{\vmu}(\tau_\varepsilon),
\CY^{(0,\ '\vbeta)}_{\vmu}(\tau_\varepsilon)),$
III.9.4

$\eta:\ C^\alpha_\vmu\lra TC^{\circ K}=\TO C^\vmu,$
V.18.1

$\breta:\ \CP^{\alpha}_{\vmu}\lra T\CP^{oJ},$
IV.8.1

$\baeta:\ \bC_g\lra\bCM_g,$
0.14.14

$\baeta:\ \bC_{A,g}\lra\CM_{A,g},$
0.14.14

$\ul\eta,$
V.18.1

$\theta_i,$ 0.2.7, II.11.4, II.12.3

$\theta_i:\ \Phi(\CM)_\nu\lra\Phi(\CM)_{\nu-i},$
0.6.9

$\theta_{i,DV}:\ (DV)_{\lambda}\rightleftharpoons(DV)_{\lambda-i'}\
:\epsilon_{i,DV},$
III.13.2

$\theta^*_i,$
II.2.13

$\theta_{\vec{K}},$
II.2.3

$\theta_{\vec{I},Q},$
II.2.20

$[\theta_{\vec{I},Q,\Lambda}],$
II.2.20

$\theta_{\varrho\leq\tau}^*,$
II.3.5

$\theta^*_{0\leq\tau},$
II.8.14

$\theta_{\varrho\leq\tau},$
II.3.4

$\theta_{\varrho\leq\tau;a},$
II.3.4

${\tilde \theta}_j,$
II.4.1

${\tilde\theta}_V,$
IV.6.3

$\theta^\FS,$
IV.6.6

$\theta^{\beta,\alpha}_{\mu}:\
\Hom_{\CA^{\alpha}_{\mu}}(\CX^{\alpha}_{\mu},\CY^{\alpha}_{\mu})\iso
\Hom_{\CA^{\alpha+\beta}_{\mu+\beta'}}
(\CX^{\alpha+\beta}_{\mu+\beta'},\CY^{\alpha+\beta}_{\mu+\beta'}),$
III.4.4

$\theta=\theta^{\alpha}_{\vmu;\vbeta}:\
\Hom_{\CA^{\alpha}_{\vmu}(K)}(\CX^{\alpha}_{\vmu}(K),\CY^{\alpha}_{\vmu}(K))
\iso\Hom_{\CA^{\alpha+\beta}_{\vmu+\vbeta'}(K)}
(\CX^{\alpha+\beta}_{\vmu+\vbeta'}(K),\CY^{\alpha+\beta}_{\vmu+\vbeta'}(K)),$
III.9.4

$\iota(\nu'):\ D^\nu\lra D^{\nu+\nu'},$ 0.4.4

$\Lambda,$
V.4.1

$^\pi\Lambda,$
II.4.4

$\lambda(\CM),$ 0.4.3

$\lambda_F,$
I.4.6

$\lambda_a,$
II.2.20

$\lambda_0,$
III.16

$\lambda_\infty,$
IV.7.6, IV.9.5

$\vlambda\geq\vmu,$
III.7.1

$\mu_Q,$
II.2.24

$\vmu_{\leq i},$
III.2.3, III.7.6

$\vmu_{\geq k},$
III.7.6

$\bnu,$
II.3.2

$\nu_1\cdot\nu_2,$ 0.2.1, II.4.1

$\nu_1\leq\nu_2,$ 0.2.1

$\nu_Q,$
II.2.20

$\nu_0,$
IV.2.4, IV.9.1, V.16.2

$\nu_a(\varrho),$
II.3.3

$\nu_\pi,$
II.4.1

$\xi:\ K^\bullet_n\lra K^\bullet_{n+1},$
IV.5.2

$\xi_{\vmu;\vbeta;d}:\
\Hom_{\CA^\beta_{\vmu}(K)}(\CX^\beta_{\vmu}(K),\CY^\beta_{\vmu}(K))\lra
\Hom_{\CA^{(0,\vbeta)}_{\vmu}(\tau_d)}(\CX^{(0,\vbeta)}_{\vmu}(\tau_d),
\CY^{(0,\vbeta)}_{\vmu}(\tau_d)),$
III.9.4

$\xi:\ T\CA^{om}\hra T\CP^{om+1},$
IV.9.5

$\pi_1(X;Y),$
I.4.1

$\pi_j,$
III.12.9

$\pi_K,$
V.2.1

$\pi^\alpha,$
V.13.2

$\pi^{\bullet\alpha},$
V.13.2

$\pi:\ J\lra I,$ 0.3.1

$\pi^F_E:\ \Ch(F)\lra\Ch(E),$
I.4.13

$\pi:\ ^{\pi}{\ '\ff}_{\chi_J}\lra {\ '\ff}_{\nu_{\pi}},$
II.4.2

$\pi:\ ^\pi\BA\lra\CA_\nu,$
II.6.12

$\pi^{\alpha}_{\vmu}:\ (D(0;1)^K\times D(0;1)^J)^{\circ}\lra
\CAO^{\alpha}_{\vmu}(K),$
III.8.1

$\varpi_{0,1,2},$
V.7.5

$\rho,$ 0.2.2

$\rho_\ell,$ 0.2.2, IV.9.1

$\rho_0,$
II.6.6

$\rho_{ab},$
III.7.2

$\varrho,$
II.3.3

$\varrho<\varrho',$
II.8.8

$\varrho_{J'},$
III.12.15

$\Sigma_\pi,$    0.3.1, II.4.1

$\sigma_k,$
IV.7.4

$\sigma_{\geq k},$
III.7.2

$\sigma_{\leq i},$
III.7.2

$\sigma_\lambda,$
III.14.1

$\sigma_{\lambda\leq\mu}^!,$
III.14.4

$\sigma_{\lambda\leq\mu}^*,$
III.14.4

$\sigma_{\lambda\leq\mu},$
III.14.4

$\sigma^!_{\lambda},\sigma^*_{\lambda}:\ \FS\lra\FS_{\leq\lambda},$
III.14.2

$\sigma_C:\ J\iso [N],$
II.8.1

$\sigma=\sigma^{\alpha,\beta}_{\mu}:\ \CA_{\mu}^{\alpha}\hra
\CA_{\mu+\beta'}^{\alpha+\beta},$
III.2.4

$\sigma=\sigma^{\alpha}_{\vmu;\vbeta}:\ \CA^{\alpha}_{\vmu}(K)\lra
\CA^{\alpha+\beta}_{\vmu+\vbeta'}(K),$
III.7.8

$\varsigma=\varsigma_{d_1,d_2},$
III.10.3

$\varsigma^\alpha:\
Q^{\bullet\alpha}_0\hra\CA^\alpha_v\times\CA^\alpha_h,$
V.14.1

$\varsigma^\alpha:\ Q^\alpha_{0a}\hra\CP^\alpha_v\times\CA^\alpha_h,$
V.15.2

$\tau_\lambda(\nu,\nu'),$ 0.4.5

$\tau^\pm_F,$
I.4.3

$\tau(\vec{K}),$
II.2.3

$\tau_A,$
II.2.6

$\tau_Q,$
II.2.20

$\tau_d,$
III.9.4

$\tau_{\geq k},$
III.7.2

$\tau_{\leq i},$
III.7.2

$\tau_C:J\iso [N],$
II.6.2

$\tau^{\alpha\beta}_{\mu}:\
\Hom_{\CA^{\alpha}_{\mu}}(\CX^{\alpha}_{\mu},\CY^{\alpha}_{\mu})\lra
\Hom_{\CA^{\beta}_{\mu}}(\CX^{\beta}_{\mu},\CY^{\beta}_{\mu}),$
III.4.4

$\tau_{\pi;\epsilon,\delta},$
III.11.2

$\tau_{\pi;d},$
III.11.9

$\Upsilon:\ _R\CC\lra\CC,$
IV.6.2

$\Phi_F(\CK),$
I.3.3

$\Phi_F:\ \CM(\BA,\CS_{\CH})\lra\Vect,$
I.3.10

$\Phi^+_F(\CK),$
0.6.2

$\Phi:\ \CM(D^\nu,\CS)\lra\Vect,$
0.6.7

$\Phi(\CM)_\lambda,$
0.6.8

$\Phi:\ \tFS\lra\Vect_f^X,$
III.4.5

$\Phi:\ \FS\lra\CC,$
0.6.12, III.12.27

$\Phi_f,$
I.2.11

$\Phi_{\sum t_j}(\CK),$
II.7.13, III.12.6

$\Phi^+_F(\CK),$
II.7.3

$\Phi_\nu:\ \CM(\CA_\nu,\CS)\lra\Vect,$
II.7.14

$\Phi_\Delta(\CI_{\nu?}),$
II.6.14

$\Phi_\lambda(\CX),$
III.4.5

$\tPhi:\ \FS\lra\tCC,$
III.12.12

$\phi_\mu(\nu_1,\nu_2),$ 0.3.6

$\phi_{\mu_1,\mu_2}(\nu,\nu'),$
0.5.4

$\phi_{\mu_1,\mu_2}(\nu_1,\nu_2,\nu),$
0.5.4

$\phi(\rho),$
0.7.3

$\phi_P(\rho),$
0.7.7

$\phi:\ \MS\lra\CC,$
0.10.7

$\phi^{\flat},$
0.10.7

$\phi^{\sharp},$
0.10.7

$\phi_C(\rho),$
0.14.4

$^{\pi}\phi_{\Delta,!}^{(\eta)}:\ \Phi_{\Delta}(^{\pi}\CI_!)\lra\
^{\pi}{\ '\ff}_{\chi_J},$
II.6.10

$^{\pi}\phi_{\Delta,*}^{(\eta)}:\ \Phi_{\Delta}(^{\pi}\CI_*)\lra\
^{\pi}{\ '\ff}^*_{\chi_J},$
II.6.10

$^{\pi}\phi_{\Delta,!*}^{(\eta)}:\ \Phi_{\Delta}(^{\pi}\CI_{!*})\iso\
^{\pi}\ff_{\chi_J},$
II.6.11

$^{\pi}\phi_{\varrho,!}^{(\eta)}:\ \Phi^+_{F_{\varrho}}(\CI(^{\pi}\Lambda)_!)
\iso\ _{\varrho}C^{-r}_{^{\pi}{\ '\ff}}(V(^{\pi}\Lambda)),$
II.8.9

$^{\pi}\phi_{!}^{(\eta)}:\
^+C^{\bullet}(^{\pi}\BA;\CI(^{\pi}\Lambda)_!)\iso\
_{\chi_J}C^{\bullet}_{^{\pi}{\ '\ff}}(V(^{\pi}\Lambda)),$
II.8.10

$^{\pi}\phi_{0,*}^{(\eta)}:\
\Phi^+_0(\CI(^{\pi}\Lambda)_*)\iso\ V(^{\pi}\Lambda)^*_{\chi_J},$
II.8.14

$\phi_{\nu,\Lambda,!*}^{(\eta)}:\
\Phi_{\nu}(\CI_{\nu}(\Lambda)_{!*})\iso L(\Lambda)_{\nu},$
II.8.18

$^{\pi}\phi_{\varrho,*}^{(\eta)}:\
\Phi^+_{F_{\varrho}}(\CI(^{\pi}\Lambda)_*)\iso\
_{\varrho}C^{-r}_{^{\pi}{\ '\ff}^*}(V(^{\pi}\Lambda)^*),$
II.8.20

$^{\pi}\phi_{*}^{(\eta)}:\
^+C^{\bullet}(^{\pi}\BA;\CI(^{\pi}\Lambda)_*)\iso\
_{\chi_J}C^{\bullet}_{^{\pi}{\ '\ff}^*}(V(^{\pi}\Lambda)^*),$
II.8.21

$_{\Lambda}\phi^{(\eta)}_{\nu,!*,0}:\
\CI_{\nu}(\Lambda)_{!*,0}\iso\ _{\nu}C^{\bullet}_{\ff}(L(\Lambda)),$
II.8.23

$\phi_{!*}^{(\eta)}:\
\Phi_{\nu}(^{\psi}\CI_{\nu}(\Lambda_0,\ldots,\Lambda_{-n})_{!*})\iso
(L(\Lambda_0)\otimes\ldots\otimes L(\Lambda_{-n}))_{\nu},$
II.10.6, II.12.7

$\phi_{!*,0}^{(\eta)}:\
^{\psi}\CI_{\nu}(\Lambda_0,\ldots,\Lambda_{-n})_{!*0}\iso\
_{\nu}C^{\bullet}_{\ff}(L(\Lambda_0)\otimes\ldots\otimes L(\Lambda_{-n})),$
II.10.8, II.12.8

$\phi_i=\phi^{\valpha}_{\mu;i}(\vd):\ \CI^{\valpha}_{\mu}(\vd)
\iso \CI^{\valpha_{\leq i}}_{\mu_{\leq i}}(\vd_{\leq i})
\boxtimes\CI^{\valpha_{\geq i}}_{\mu}(\vd_{\geq i}),$
III.3.2

$\phi_i=\phi^{\valpha}_{\mu;i}(\vd):\ \CID^{\valpha}_{\mu}(\vd)
\iso\CID^{\valpha_{\leq i}}_{\mu_{\leq i}}(\vd_{\leq i})
\boxtimes\CID^{\valpha_{\geq i}}_{\mu}(\vd_{\geq i}),$
III.3.5

$\phi_i=\phi_{i;\vmu}^{\valpha}(\tau):\
\CI_{\vmu}^{\valpha}(\tau)\cong
\CI_{\vmu_{\leq i}}^{\valpha_{\leq i}}(\tau_{\leq i})
\boxtimes\fbox{$\times$}_{k\in K_i}\CI_{\vmu_{\geq k}}^{\valpha_{\geq k}}
(\tau_{\geq k}),$
III.8.2

$\phi_i=\phi_{i;\vmu}^{\valpha}(\tau):\
\CID_{\vmu}^{\valpha}(\tau)\cong
\CID_{\vmu_{\leq i}}^{\valpha_{\leq i}}(\tau_{\leq i})
\boxtimes\fbox{$\times$}_{k\in K_i}\CID_{\vmu_{\geq k}}
^{\valpha_{\geq k}}(\tau_{\geq k}),$
III.8.5

$\phi_{d_1,d_2}^{\valpha_1,\valpha_2}:\
\CX^{\valpha_1}(\tau_{d_1})|
_{\CA^{\valpha_1}_{\vlambda}(\tau_{d_1})\cap
  \CA^{\valpha_2}_{\vlambda}(\tau_{d_2})}\iso
\CX^{\valpha_2}(\tau_{d_2})|
_{\CA^{\valpha_1}_{\vlambda}(\tau_{d_1})\cap
  \CA^{\valpha_2}_{\vlambda}(\tau_{d_2})},$
III.10.3

$\phi_{\pi}:\
\otimes_K\ X_k\iso\otimes_L\ (\otimes_{K_l}\ X_k),$
III.11.3

$\phi_{K\lra L}:\ r_{\pi}\circ\Psi_K\iso
\Psi_L\circ\Psi_{K\lra L},$
III.11.10

$\phi_{\valpha}:\ q^*_{\valpha}\CI^{\alpha}\iso
p_K^*\pi_K^*\CI^{\alpha_K}\boxtimes\Boxtimes_k\ \CI^{\alpha_k^{\sim}}_D,$
IV.2.8

$\phi_{\valpha,\beta}:\  q^*_{\valpha,\beta}\CX^{\alpha+\beta}_{\vmu}
\iso p^*_{\valpha,\beta}((\Boxtimes_{k\in K}\ \CX^{\alpha_k}_{\mu_k})
\boxtimes\CI^{\bullet\beta}_{\vmu-\valpha}),$
IV.3.5

$\phi:\ \tCO_\kappa\iso\tCO_{\zeta},$
IV.9.2

$\phi_\tau:\  q_\tau^*\pi_J^*\CH^\alpha\iso p_K^*\pi_K^*\CH^{\alpha_K}
\boxtimes\fbox{$\times$}_{k\in K}\ \CI^{\tau^{-1}(k)},$
V.2.2

$\phi_{\valpha,\beta}:\ q^*_{\valpha,\beta}\CX^{\alpha+\beta}_{\vmu}\iso
p^*_{\valpha,\beta}((\Boxtimes_{k\in K}\ \CX^{\alpha_k}_{\mu_k})
\boxtimes\CH^{\bullet\beta}_{\vmu-\valpha}),$
V.16.6

$\phi_{\valpha\vbeta\vgamma,\xi}:\
q^*_{\valpha\vbeta\vgamma,\xi}
\CX^{\alpha+\beta+\gamma+\xi}_{\vlambda\vmu\vnu}\iso
p^*_{\valpha\vbeta\vgamma,\xi}
((\Boxtimes_{k\in K}\ \CX^{\alpha_k}_{\lambda_k})
\boxtimes
(\Boxtimes_{j\in J}\ \CY^{\beta_j,\gamma_j}_{\mu_j,\nu_j})\boxtimes
\CH^{\bullet\xi}_{\vlambda\vmu\vnu-\valpha\vbeta\vgamma}),$
V.16.8

$\chi_K,$
II.3.3, II.4.1

$\Psi_f,$
I.2.11

$\Psi_K,$
III.11.8

$\Psi_{\pi;d}:\ \CM(\CA^{\alpha}(\tau_{\pi;d}))\lra
\CM(\CA^{\alpha}(L)_d\times\prod_{l\in L}\CO(K_l)),$
III.11.9

$\Psi_{K\lra L}:\ \CM(\CA^{\alpha}(K))\lra
\CM(\CA^{\alpha}(L)\times\prod_l\ \CO(K_l)),$
III.11.9

$\Psi_K:\ ^K\FS\lra\Locsys(\CO(K);\FS),$
III.11.11

$\psi(\nu_1,\nu_2),$ 0.4.3

$\psi(\nu_1,\nu_2,\nu),$
0.5.8

$\psi^{\CM\boxtimes\CN}(\nu,\nu'),$
0.5.8

$^\psi\CI(\Lambda_0,\ldots,\Lambda_{-n};\nu)_?,$
II.10.5

$^n\psi_i:\ \CCD(^n\BA;\CS_{\Delta})\lra\ \CCD(^{n-1}\BA;\CS_{\Delta}),$
II.10.1

$\psi:\ \CCD(^n\BA;\CS_{\Delta})\lra\CCD(\BA;\CS),$
II.10.3

$\psi^{\alpha,\beta}(d):\
\CX^{(\alpha,\beta)}(d)\iso\CID^{(\alpha,0)}_{\lambda-\beta'}(d)
\boxtimes\CX^{(0,\beta)}(d),$
III.4.2

$\psi^{\alpha,\beta}_{\mu}(d):\ \CX^{(\alpha,\beta)}_{\mu}(d)
\iso\CID^{(\alpha,0)}_{\mu-\beta'}(d)\boxtimes\CX^{(0,\beta)}_{\mu}(d),$
III.4.4

$\psi(\tau):\ \CX^{(\alpha,\vbeta)}(\tau)\cong
\CID_{\vlambda(\tau)_{\leq 1}}^{\alpha}(\tau_{\leq 1})
\boxtimes\CX^{(0,\vbeta)}(\tau),$
III.9.2

$\Omega,$
V.3.2

$\omega,$
0.14.5, V.3.2

$\omega_X,$
I.2.5

$\One,$  0.2.14

$\One_\FS,$
III.11.11

$(\cdot,\cdot)$ on $'\ff,$ 0.2.6

$^{\vee}:\ \CC^{opp}\iso\CC_r,$
IV.4.6

$^*:\ \CC\lra\CC^{opp},$
IV.4.4

$\otimes:\ \FS\otimes\FS\lra\FS,$
0.5.9

$\langle\otimes_A\ \{\CX_a; \CR_n\}\rangle^{(m)}_{\vA^0,g,\Gamma},$
0.15.4

$\otimes_A\ \{\CX_a; \CR_n\},$
0.15.2

$\otimes_K:\ \FS^K\lra\Locsys(\CO(K);\FS),$
III.11.11

$\langle\ \rangle^{(m)}_{\CC},$
0.16.1

$\langle\otimes_A\ \{L_a;\bR_{\MS n}\}\rangle_{\MS},$
0.16.1

$\overset{\bullet}{\otimes},$
III.11.4

$\totimes,$
IV.9.2

$\Boxtimes^{(\nu)}_K\ \CM_k,$
0.8.5

$\Boxtimes_{A,g}^{(\nu)}\ \{\CX_a; \CR_n\},$
0.15.3

$\dpar_i\valpha,$
III.2.2

$\dpar_i\vd,$
III.2.2

$\dpar_i\mu,$
III.2.3

$\dpar_i\vd,$
III.7.2

$\dpar_i\sigma,$
III.7.2

$\dpar_i\tau,$
III.7.2

$\dpar_i\valpha,$
III.7.4

$\dpar_i\vmu,$
III.7.6

$(\bullet)^{\Sigma_{\pi}}:\ \CCD^b(\CA_{\alpha},\Sigma_{\pi})\lra
\CCD^b(\CA_{\alpha}),$
III.12.5

$\nabla_\kappa,$
IV.9.2

$\infty_{g,K},$
V.19.2

\newpage

\centerline{\bf Index of Terminology}

$A$-curve, 0.14.13

adjacent facet, I.3.2

adjoint representation, 0.12.5, V.11.1

admissible element $\nu\in Y^+$, 0.8.4, (in $\BN[X_\ell]$), 0.14.6

admissible pair $(\vmu,\alpha)$, IV.3.1, IV.9.1

affine Lie algebra $\hat\fg$, 0.10.6, IV.9.2

antipode, 0.2.10

arrangement (real) 0.6.1, I.3.1

associativity condition, 0.3.6, 0.4.3, 0.7.3, 0.7.7

baby Verma module, 0.2.17, III.13.8

balance, 0.2.18

balance function, IV.2.4

base (of a tree), III.7.2

bbrt category, IV.9.2

braiding, 0.2.18

braiding local system, 0.3.5

branch, III.7.2

canonical map, 0.6.6, I.3.11, II.7.10

Cartan datum, 0.2.1

cartesian functor, 0.15.8

chamber, I.3.1

coaction, II.2.18

cochain complex of an arrangement, I.3.13

cocycle condition, 0.8.5

cohesive local system (CLS for short) $\CH$ of level $\mu$ (over $P$), 0.7.7,
IV.2.8, (over $C/S$), 0.14.4, (over $\CA$), III.8.4

compatible isomorphisms, V.6.2

comultiplication, 0.2.10

configuration space of colored points on the affine line, II.6.12

concave complex, 0.9.3, IV.4.2

convex complex, 0.9.3, IV.4.2

conformal blocks, 0.10.3, IV.5.10

conic complex, I.2.13

convolution, V.15.3

coorientation, I.4.6

cutting, III.2.2, III.7.2

deleting morphism, 0.14.20

depth, II.2.2

determinant line bundle, 0.14.1, V.2.2

diagonal stratification, II.6.1

dominant weight, III.18.2

dropping, III.2.2

dual cell, I.3.2

dualizing complex, I.2.5

edge (of an arrangement), 0.6.1, I.3.1

elementary tree, III.7.2

enhanced disk operad, III.7.5

enhanced graph, 0.14.13

enhanced tree, III.7.3

enhancement, III.7.3

face, I.3.1

facet, 0.6.1, I.3.1

factorizable sheaf (supported at $c$), 0.4.3, III.4.2, (finite, FFS for short),
0.4.6, III.5.1, (over $^K\CA,\ ^K\CI$), III.9.2

factorization isomorphisms, 0.3.5, 0.4.3, 0.7.3, 0.7.7, V.2.2, V.16.6

first alcove, 0.10.2

flag, 0.6.1, I.3.2

fusion functor, V.18.2

fusion structure, 0.15.1, (of multiplicative central charge $c$), 0.15.8,
V.18.2

$g$-admissible element of $\BN[X_\ell]$, V.2.1

$g$-admissible pair, V.16.2

$g$-positive $m$-tuple of weights, V.17.1

generalized vanishing cycles of $\CK$ at a facet $F$, I.3.3

generalized baby Verma module, III.14.10

gluing, 0.5.8, 0.8.6, III.10.3, IV.3.5, V.16.6

good $A_0$-tuple (w.r.t. $(g,\nu)$), 0.14.15, 0.14.22

good object (w.r.t. $\fu^+$ or $\fu^-$), 0.9.2

good resolution (right or left), 0.9.4

good stratification, I.2.9

good surjection, IV.4.3

half-monodromy, II.6.2, II.8.1

height (of tree), III.7.2

Heisenberg local system, 0.14.11, V.2.2

Hochschild complex, II.3.1

homogeneous element, II.2.2

$i$-cutting, III.2.3

I-sheaf, V.12.2

level, 0.7.7

local system of conformal blocks, 0.11.7, IV.9.2

marked disk operad, III.7.7

marked tree, III.7.6

marking, I.3.2

maximal trivial direct summand, 0.10.1, IV.5.9

modular functor, V.18.2

multiplicative central charge, 0.14.12

nearby cycles, I.2.11

neighbour (left), II.6.2, II.8.1

open $I$-colored configuration space, 0.3.5

operad of disks, III.7.3

operad of disks with tangent vectors $\CCD$, 0.7.2

orientation sheaf

perverse sheaf, 0.2.4, I.2.9

Poincar\'e groupoid, I.4.1

positive facet, 0.6.1, II.7.1

positive flag, II.7.1

positive $m$-tuple of weights, IV.8.3

principal stratification, II.7.1, II.7.14, III.7.9, IV.3.2, V.16.3

quantum commutator, II.2.20

quantum group with divided powers, IV.6.1

real arrangement, I.3.1

refinement, II.3.4

regular object, 0.16.1

regular representation, 0.12.3, V.10.1

regular sheaf, 0.13.3, V.14.2

relative singular $n$-cell, I.2.7

rhomb diagram, 0.3.5

ribbon category, 0.2.23

rigidity, IV.4.4

rigid object, III.15.3

root datum, 0.2.1

semiinfinite Ext, Tor, 0.9.7, IV.4.7

sewing morphism, 0.14.20

shape (of a tree), III.7.2

simple sewing, 0.14.20

sign local system, 0.3.9

skew-antipode, 0.2.10

skew-$\Sigma_\pi$-equivariant morphism, II.6.12

small quantum group $\fu$, 0.2.10, $\bu_k$, 0.2.12

standard local system, 0.3.10, ($X$-colored, over $\CCD$), 0.7.5

standard braiding local system (over the configuration space $\CA(2)^\circ$),
0.5.6

standard sheaves, I.4.5, III.6.1

Steinberg module, III.16

Steinberg sheaf, III.16

substitution isomorphism, III.7.3

tensor product of categories, III.10.2

tensor structure (rigid), 0.2.14

thickness, III.7.2

toric stratification, III.2.7, III.7.9

tree, III.7.2

two-sided \v{C}ech resolution, IV.8.5

unfolding, 0.3.1, II.6.12

universal Heisenberg local system, 0.14.18

vanishing cycles, I.2.11

vanishing cycles of $\CK$ across $F$, 0.6.2

vanishing cycles at the origin, 0.6.7

variation map, 0.6.6, I.3.11, II.7.10

Verma module, II.2.15

weight, II.2.1

Weyl group, 0.2.2, III.18.3

$X$-colored local system over $\CCD$, 0.7.3

young tree, III.7.2

\end{document}